\documentclass[11pt,leqno]{article}
\normalfont
\renewcommand{\baselinestretch}{1.15}

\usepackage{amsfonts}
\usepackage{mathrsfs}
\usepackage{graphicx}
\usepackage{hieroglf}
\usepackage{MnSymbol}
\usepackage{aurical}
\usepackage[T1]{fontenc}
\usepackage[dvipsnames]{xcolor}
\usepackage{extarrows}
\usepackage{footmisc}
\usepackage{xpicture}

\newcommand{\IntroSection}{1}
\newcommand{\RoutineSection}{2}
\newcommand{\SubstructureSection}{3}
\newcommand{\GStructureIntroSection}{4}
\newcommand{\FrameSection}{5}
\newcommand{\GridPosetSection}{6}
\newcommand{\TwoColorSturdySection}{7}
\newcommand{\WeylSection}{8}
\newcommand{\LFKSection}{9}
\newcommand{\FlowerSection}{10}
\newcommand{\QuasiExampleSection}{11}
\newcommand{\MinusculeExampleSection}{12}
\newcommand{\MinusculeLatticePosetSection}{13}

\newcommand{\ASevenFigure}{Figure 1.1}
\newcommand{\ASevenFigureTwo}{Figure 1.2}

\newcommand{\MountainValleyLemma}{Lemma 2.1}
\newcommand{\ModularLatticeTheorem}{Proposition 2.2}
\newcommand{\ModularLawLemma}{Lemma 2.3}
\newcommand{\DistributiveIsModularLemma}{Lemma 2.4}
\newcommand{\AbstractModularDistributiveProperties}{Lemmas 2.3 and 2.4.1}
\newcommand{\ColorsLemma}{Lemma 2.5}
\newcommand{\EncapsulateProposition}{Proposition 2.6}
\newcommand{\BooleanExample}{Example 2.7}
\newcommand{\FundamentalTheorem}{Theorem 2.8}
\newcommand{\FirstCorollary}{Corollary 2.9}
\newcommand{\JMCorollary}{Corollary 2.10}
\newcommand{\CompressionPosetLattice}{Lemma 2.11} 

\newcommand{\DCDLFigs}{Figures 2.1 and 2.2}
\newcommand{\PosetAndLatticeFig}{Figure 2.1}
\newcommand{\DualAndRecolorFig}{Figure 2.2}
\newcommand{\ComponentFig}{Figure 2.3}
\newcommand{\TopoBalancedRankedFigure}{Figure 2.4}
\newcommand{\OrderIdealFig}{Figure 2.5}
\newcommand{\NewOrderIdealFig}{Figure 2.6}

\newcommand{\MfiveNfiveTheorem}{Theorem 3.1}
\newcommand{\CancellationTheorem}{Corollary 3.2}
\newcommand{\ModularDistributiveSublattices}{Lemma 3.3}
\newcommand{\FullLengthLemma}{Lemma 3.4}
\newcommand{\FullLengthWithinProduct}{Proposition 3.5}
\newcommand{\FullLengthTheorem}{Theorem 3.6}
\newcommand{\IntervalProp}{Proposition 3.7}
\newcommand{\JCompResult}{Theorem 3.8}

\newcommand{\DilworthsTheorem}{Corollary 3.9}
\newcommand{\TwoChainLemma}{Lemma 3.10}
\newcommand{\JCompTheorem}{Theorem 3.11} 
\newcommand{\JCompStuff}{Theorems 3.8 and 3.11}
\newcommand{\NewSublatticeResults}{Theorems 3.6, 3.8, and 3.11}
\newcommand{\CFourFigure}{Figure 3.1}
\newcommand{\FFourFigure}{Figure 3.2}

\newcommand{\ShawtyCFourFigure}{Figure 3.1.2/3/4}
\newcommand{\ShawtyCFourDistributiveFigure}{Figure 3.1.2/3}
\newcommand{\ShawtyFFourFigure}{Figure 3.2.2/3}
\newcommand{\ShawtyFigureList}{Figures 3.1.2/3/4 and 3.2.2/3}
\newcommand{\IEGGraphFigure}{Figure 4.1}
\newcommand{\GCMTheorem}{Theorem 4.1}
\newcommand{\NewLFKTheorem}{Theorem 4.2}
\newcommand{\DistCoreConjecture}{Conjecture 4.3}
\newcommand{\ESixLatticeFigure}{Figure 4.2}

\newcommand{\ESixFigures}{Figures 4.2 and 4.3}
\newcommand{\DistributiveProp}{Proposition 5.1} 
\newcommand{\ChainProp}{Proposition 5.2} 
\newcommand{\FramingProps}{Propositions 5.1/5.2}
\newcommand{\StackDefsOne}{Definitions 5.3}
\newcommand{\StackEarlyLemma}{Lemma 5.4} 
\newcommand{\StackDefsTwo}{Definitions 5.5}
\newcommand{\StackDefsOneTwo}{Definitions 5.3 and 5.5}
\newcommand{\StackLemma}{Lemma 5.6}
\newcommand{\StackDefsThree}{Definitions 5.7}

\newcommand{\IsomLatticeTheorem}{Theorem 5.8}

\newcommand{\ESixStackFigure}{Figure 5.1}
\newcommand{\GridPosets}{Figure 6.1}
\newcommand{\GridPosetsII}{Figure 6.2}
\newcommand{\GridPosetsFigureList}{Figures 6.1 and 6.2}
\newcommand{\FundPosets}{Figure 6.3}
\newcommand{\GridPosetsTotalOrderList}{Figures 6.1, 6.2, and 6.3}
\newcommand{\FundLatticeIdealsFigure}{Figure 6.4}
\newcommand{\WeightFigure}{Figure 6.5}
\newcommand{\FundFigures}{Figures 6.4 and 6.5}
\newcommand{\WeightsLemma}{Lemma 6.1}
\newcommand{\TechnicalProposition}{Lemma 6.2}
\newcommand{\ClassTheorems}{Theorems 6.3 and 6.4}
\newcommand{\FundamentalClass}{Theorem 6.3}
\newcommand{\SemistandardClass}{Theorem 6.4}
\newcommand{\TechPropAalphaFigure}{Figure 6.6}
\newcommand{\TechPropAbetaFigure}{Figure 6.7}
\newcommand{\TechPropBalphaFigure}{Figure 6.8}
\newcommand{\TechPropBbetaFigure}{Figure 6.9}
\newcommand{\TechPropGalphaFigure}{Figure 6.10}
\newcommand{\TechPropGbetaFigure}{Figure 6.11}
\newcommand{\SubordinateGridDecomp}{Lemma 7.1}
\newcommand{\GStructureDCDLTheorem}{Theorem 7.2}
\newcommand{\GeneratorsAndRelations}{Proposition 8.1}
\newcommand{\OrbitProposition}{Theorem 8.2}
\newcommand{\FiniteWeylGroupResult}{Corollary 8.3}
\newcommand{\FiniteIndexCorollary}{Corollary 8.4}
\newcommand{\NonconstantTheorem}{Corollary 8.5}
\newcommand{\InnerProductTheorem}{Corollary 8.6}
\newcommand{\NGCorollary}{Corollary 8.7}
\newcommand{\KacMoodyCorollary}{Corollary 8.8}
\newcommand{\KMRepCorollary}{Corollary 8.9}
\newcommand{\KMCorollaryList}{Corollaries 8.8 and 8.9}
\newcommand{\OrbitCorollaries}{Corollaries 8.4 and 8.5}
\newcommand{\LFK}{Theorem 9.1}
\newcommand{\LongestLemma}{Lemma 10.1}
\newcommand{\MinimalWeightsLemma}{Lemma 10.2}
\newcommand{\SaturatedLemma}{Lemma 10.3}
\newcommand{\FTWSF}{Theorem 10.4}
\newcommand{\SpecializeTheorem}{Theorem 10.5}
\newcommand{\SplittingPosetTheorem}{Theorem 10.6}
\newcommand{\VertexColoringTheorem}{Theorem 10.7}
\newcommand{\WeylKacCharacter}{Theorem 10.8}
\newcommand{\CompletelyReducible}{Theorem 10.9}
\newcommand{\MainLieRepTheorem}{Theorem 10.10}
\newcommand{\TopoNotLatticeExercise}{Exercise 10.11}
\newcommand{\RootsFigure}{Figure 10.1}
\newcommand{\HighestShawtyLemma}{Lemma 11.1}
\newcommand{\ShortAdjointDCML}{Lemma 11.2}
\newcommand{\QMSplitting}{Theorem 11.3}
\newcommand{\MinusculeCharacterizationLemma}{Lemma 12.1}
\newcommand{\MinusculeStructureLemma}{Lemma 12.2}
\newcommand{\MinusculeWeylBialt}{Proposition 12.3}
\newcommand{\MeetJoinMinusculeProp}{Lemma 12.4}
\newcommand{\MinusculeMeetJoinFigure}{Figure 12.1}
\newcommand{\MinusculeMeetJoinFigureTwo}{Figure 12.2}
\newcommand{ \MinusculeMeetJoinFigureThree}{Figure 12.3}
\newcommand{\MeetJoinFollowUp}{Lemma 12.5}
\newcommand{\MinusculeDCDL}{Theorem 12.6}
\newcommand{\StarBowtieLemma}{Lemma 12.7}
\newcommand{\ijChainLemma}{Lemma 12.8}
\newcommand{\IdealMatchingProp}{Theorem 12.9}
\newcommand{\MinusculePosetIsLattice}{Lemma 12.10}
\newcommand{\MinusculePosetIsTopoBalanced}{Lemma 12.11}
\newcommand{\MinusculePosetIsDistributive}{Theorem 12.12}
\newcommand{\MinusculeTopTree}{Theorem 12.13}
\newcommand{\MinusculeDistributiveResults}{Theorems 12.6 and 12.12}
\newcommand{\NewMinusculeResults}{Theorems 12.9 and 12.13}
\newcommand{\DoppelLemma}{Lemma 13.1}
\newcommand{\GStructureDCDLCorollary}{Theorem 13.2}

\newfont{\mymathcal}{ecci0600 at 11pt}

\newfont{\myscbolditalics}{ecoc0500 at 11pt}

\newfont{\mysmallscbolditalics}{ecoc0500 at 8pt}

\newfont{\mybolditalics}{ecbi0500 at 11pt}

\newcommand{\myqk}{\mbox{\mybolditalics k}}
\newcommand{\myqx}{\mbox{\mybolditalics x}}
\newcommand{\myqy}{\mbox{\mybolditalics y}}

\newcommand{\myqX}{\mbox{\mybolditalics X}}
\newcommand{\myqY}{\mbox{\mybolditalics Y}}
\newcommand{\myqz}{\mbox{\mybolditalics z}}
\newcommand{\myqh}{\mbox{\mybolditalics h}}

\newcommand{\myqP}{\mbox{\mybolditalics P}}
\newcommand{\myqQ}{\mbox{\mybolditalics Q}}

\newfont{\eulercursive}{eurm10 at 11pt}
\newcommand{\myp}{\mbox{\eulercursive p}}

\newcommand{\myd}{\mbox{\eulercursive d}}
\newcommand{\mya}{\mbox{\eulercursive a}}
\newcommand{\mym}{\mbox{\eulercursive m}}

\newfont{\smalleulercursive}{eurm10 at 9pt}

\newcommand{\mysmallm}{\mbox{\smalleulercursive m}}
\newcommand{\mysmallp}{\mbox{\smalleulercursive p}}

\newfont{\smallereulercursive}{eurm10 at 7pt}

\newfont{\largeeulercursive}{eurm10 at 14pt}

\newcommand{\myr}{\mbox{\eulercursive r}}
\newcommand{\mysmallr}{\mbox{\smalleulercursive r}}

\newcommand{\myi}{\mbox{\eulercursive i}}

\newcommand{\QED}{\raisebox{0.5mm}{\fbox{\rule{0mm}{1.5mm}\ }}}

\newcounter{myfn}[page]
\renewcommand{\thefootnote}{\fnsymbol{footnote}}

\newcounter{rone}
\newcounter{rtwo}
\newcounter{rthree}
\newcounter{rfour}
\newcounter{rfive}
\newcounter{rsix}
\newcounter{rseven}
\newcommand{\myA}{\mbox{\sffamily A}}
\newcommand{\mysmallA}{\mbox{\footnotesize \sffamily A}}
\newcommand{\mytinyA}{\mbox{\tiny \sffamily A}}
\newcommand{\myB}{\mbox{\sffamily B}}

\newcommand{\myC}{\mbox{\sffamily C}}
\newcommand{\mysmallC}{\mbox{\footnotesize \sffamily C}}
\newcommand{\mytinyC}{\mbox{\tiny \sffamily C}}
\newcommand{\myD}{\mbox{\sffamily D}}

\newcommand{\myE}{\mbox{\sffamily E}}
\newcommand{\mysmallE}{\mbox{\footnotesize \sffamily E}}
\newcommand{\mytinyE}{\mbox{\tiny \sffamily E}}
\newcommand{\myF}{\mbox{\sffamily F}}
\newcommand{\mysmallF}{\mbox{\footnotesize \sffamily F}}
\newcommand{\mytinyF}{\mbox{\tiny \sffamily F}}
\newcommand{\myG}{\mbox{\sffamily G}}

\newcommand{\mytinyG}{\mbox{\tiny \sffamily G}}

\newcommand{\myX}{\mbox{\sffamily X}}

\newcommand{\mytinyX}{\mbox{\tiny \sffamily X}}
\newcommand{\myT}{\mbox{\sffamily T}}

\newcommand{\mytinyT}{\mbox{\tiny \sffamily T}}
\newcommand{\myK}{\mbox{\sffamily K}}

\newcommand{\myM}{\mbox{\sffamily M}}

\newcommand{\myN}{\mbox{\sffamily N}}

\newcommand{\mystackedP}{\mbox{\small \Fontauri P}}
\newcommand{\mytinystackedP}{\mbox{\scriptsize \Fontauri P}}
\newcommand{\mylargestackedP}{\mbox{\Huge \Fontauri P}}
\newcommand{\myLieop}{\mbox{\small \Fontlukas L}\, }
\newcommand{\scaffoldS}{\mbox{\footnotesize \Fontauri S}}
\newcommand{\bigscaffoldS}{\mbox{\huge \Fontauri S}}
\newcommand{\mysetF}{\mbox{\small \Fontauri F}}
\newcommand{\mysmallsetF}{\mbox{\scriptsize \Fontauri F}}
\newcommand{\myspecialf}{\mbox{\small \Fontauri f}}
\newcommand{\zeroweight}{\mbox{\scriptsize \Fontlukas Z}}
\newcommand{\tinyzeroweight}{\mbox{\tiny \Fontlukas Z}}
\newcommand{\myvarZ}{\mbox{\scriptsize \sffamily Z}}
\newcommand{\myfancyQ}{\mbox{\Fontlukas Q}}
\newcommand{\mybigfancyQ}{\mbox{\LARGE \Fontlukas Q}}
\newcommand{\mysmallfancyQ}{\mbox{\footnotesize \Fontlukas Q}}

\newenvironment{frcseries}{\fontfamily{frc}\selectfont}{}
\newcommand{\textfrc}[1]{{\frcseries#1}}
\newcommand{\mathfrc}[1]{\text{\large \bf \textfrc{#1}}}
\newcommand{\mygens}{\mathfrc{s}}
\newcommand{\mysmallgens}{\text{\footnotesize \bf \textfrc{s}}}
\newcommand{\mathfrcthin}[1]{\text{\large \textfrc{#1}}}
\newcommand{\mygreeku}{\mathfrcthin{v}}
\newcommand{\myquant}{\text{\small \bf \textfrc{q}}}

\newcommand{\myroot}{\mbox{\Large \Fontauri r}}

\newcommand{\mysmallestroot}{\mbox{\small \Fontauri r}}
\newcommand{\mysroot}{\mbox{\Large \Fontlukas s}}

\newcommand{\aelt}{\mathbf{a}} \newcommand{\belt}{\mathbf{b}}
 \newcommand{\delt}{\mathbf{d}}
\newcommand{\eelt}{\mathbf{e}} \newcommand{\felt}{\mathbf{f}}
\newcommand{\gelt}{\mathbf{g}} \newcommand{\helt}{\mathbf{h}}

\newcommand{\melt}{\mathbf{m}} \newcommand{\nelt}{\mathbf{n}}
 \newcommand{\pelt}{\mathbf{p}}
\newcommand{\qelt}{\mathbf{q}} \newcommand{\relt}{\mathbf{r}}
\newcommand{\selt}{\mathbf{s}} \newcommand{\telt}{\mathbf{t}}
\newcommand{\uelt}{\mathbf{u}} \newcommand{\velt}{\mathbf{v}}
\newcommand{\welt}{\mathbf{w}} \newcommand{\xelt}{\mathbf{x}}
\newcommand{\yelt}{\mathbf{y}} \newcommand{\zelt}{\mathbf{z}}

\newcommand{\ecolor}{\mathscr{E}}
\newcommand{\vcolor}{\mathscr{V}}
\newcommand{\EdgeSet}{\mbox{\sffamily Edges}}
\newcommand{\SmallEdgeSet}{\mbox{\scriptsize \sffamily Edges}}
\newcommand{\VertexSet}{\mbox{\sffamily Vertices}}

\newcommand{\Jcolor}{\mathbf{J}_{\mbox{\tiny \textnormal{color}}}}
\newcommand{\jcolor}{\mathbf{j}_{\mbox{\tiny \textnormal{color}}}}
\newcommand{\Mcolor}{\mathbf{M}_{\mbox{\tiny \textnormal{color}}}}
\newcommand{\mcolor}{\mathbf{m}_{\mbox{\tiny \textnormal{color}}}}
\newcommand{\Lcolor}{\mathbf{L}_{\mbox{\tiny \textnormal{color}}}}
\newcommand{\dist}{\mbox{\sffamily dist}}
\newcommand{\pathlength}{\mbox{\sffamily path{\_}length}}

\newcommand{\comp}{\mbox{\sffamily comp}}
\newcommand{\mychain}{\mbox{\sffamily chain}}
\newcommand{\mycolor}{\mbox{\sffamily color}}

\newcommand{\myabs}{\rule[-1.5mm]{0.2mm}{5mm}\, }

\newcommand{\dichromatic}{two-color }

\newcommand{\wt}{\mbox{\sffamily wt}}
\newcommand{\mychar}{\mbox{\sffamily char}}
\newcommand{\mymult}{\mbox{\sffamily mult}}

\newcommand{\lwt}{wt}

\newcommand{\LAone}{L_{\mytinyA_{1} \oplus \mytinyA_{1}}}

\newcommand{\LAtwo}{L_{\mytinyA_{2}}}

\newcommand{\LBtwo}{L_{\mytinyC_{2}}}

\newcommand{\LGtwo}{L_{\mytinyG_{2}}}

\newcommand{\Pba}{P_{\mathscr{G}}^{\beta\alpha}(\lambda)}
\newcommand{\Pab}{P_{\mathscr{G}}^{\alpha\beta}(\lambda)}

\newcommand{\WGF}{\mbox{\small \sffamily WGF}}
\newcommand{\RGF}{\mbox{\small \sffamily RGF}}
\newcommand{\CARD}{\mbox{\small \sffamily CARD}}
\newcommand{\posetlength}{\mbox{\small \sffamily LENGTH}}

\newcommand{\mymax}{\mbox{\sffamily max}}
\newcommand{\mymin}{\mbox{\sffamily min}}

\newcommand{\barx}{\overline{x}}
\newcommand{\mygridchain}{\mbox{\small \sffamily grid{\_}chain}}
\newcommand{\mygridcolor}{\mbox{\small \sffamily grid{\_}color}}
\newcommand{\mydepth}{\mbox{\small \sffamily depth}}

\newcommand{\mytier}{\mbox{\sffamily tier}}

\newcommand{\mysign}{\mbox{\sffamily sgn}}

\newcommand{\tbeta}{{\beta}}

\newcommand{\myarrow}[1]{\stackrel{#1}{\rightarrow}}
\newcommand{\mybackarrow}[1]{\stackrel{\ #1}{\leftarrow}}
\newcommand{\mylongarrow}[1]{\stackrel{#1}{\longrightarrow}}
\newcommand{\mylongbackarrow}[1]{\stackrel{\ #1}{\longleftarrow}}

\newcommand{\eqset}{\xlongequal{\text{set}}}%
\newcommand{\eqmulti}{\xlongequal{\text{multiset}}}%

\newcommand{\disjointunion}{\setlength{\unitlength}{0.14cm}
\ 
\begin{picture}(2,2) 
\put(0,0){\Large$\cup$}
\put(1.065,1.45){\circle*{0.6}}
\end{picture}\ }

\newcommand{\CircleInteger}[1]{
\setlength{\unitlength}{0.14cm}
\begin{picture}(3,2) 
\put(1,1){\circle{2}}
\put(0.55,0.6){{\tiny #1}}
\end{picture}
}

\newcommand{\CircleIntegerm}{
\setlength{\unitlength}{0.14cm}
\begin{picture}(3,2) 
\put(1,1){\circle{2}}
\put(0.3,0.6){{\tiny $m$}}
\end{picture}
}

\newcommand{\myposetstack}{
\setlength{\unitlength}{1cm}
\begin{picture}(0.25,0) 
\put(0,0){\mbox{\normalsize $\downupharpoons$}}
\end{picture}
}	
\newcommand{\mysmallposetstack}{
\setlength{\unitlength}{1cm}
\begin{picture}(0.25,0) 
\put(0,0){\mbox{\footnotesize $\downupharpoons$}}
\end{picture}
}	

\newcommand{\mybirdseye}{
{^{\mbox{\scriptsize\textpmhg{e}}\!\!}}
}
\newcommand{\mysmallbirdseye}{
{^{\mbox{\tiny\textpmhg{e}}\!\!}}
}
\newcommand{\mylargebirdseye}{
{^{\mbox{\large\textpmhg{e}}\!\!}}
}

\newcommand{\TypeEboxDot}[1]{
\setlength{\unitlength}{0.5cm}
\begin{picture}(1,1)
\thinlines
\put(0.5,0.5){\color{#1}\circle*{0.6}}
\put(0.5,0.5){\circle{0.6}}
\end{picture}
}

\newcommand{\TwoCitiesGraphWithLabels}{
\setlength{\unitlength}{0.75in}
\begin{picture}(1.65,0.25)
\put(0.25,0){\begin{picture}(1,0)
            \put(0,0.1){\circle*{0.05}}
            \put(-0.20,-0.05){\large $\gamma_{1}$}
            \put(1,0.1){\circle*{0.05}}
            \put(1.05,-0.05){\large $\gamma_{2}$}
            \put(0,0.1){\line(1,0){1}}
            \put(0.2,0.1){\vector(1,0){0.1}}
            \put(0.8,0.1){\vector(-1,0){0.1}}
            \put(0.225,-0.05){\footnotesize $p$}
            \put(0.71,-0.05){\footnotesize $q$}
            \end{picture}}
\end{picture}}

\newcommand{\DNGTwoNodeArrows}{
\setlength{\unitlength}{0.75in}
\begin{picture}(1.3,0.25)
\put(0,0){\begin{picture}(1,0)
            \put(0,0.1){\circle*{0.05}}
            \put(-0.2,-0.05){\large $\gamma_{i}$}
            \put(1,0.1){\circle*{0.05}}
            \put(1.05,-0.05){\large $\gamma_{j}$}
            \put(0,0.1){\line(1,0){1}}
            \put(0.2,0.1){\vector(1,0){0.1}}
            \put(0.3,0.1){\vector(1,0){0.1}}
            \put(0.6,0.1){\vector(-1,0){0.1}}
            \put(0.7,0.1){\vector(-1,0){0.1}}
            \put(0.8,0.1){\vector(-1,0){0.1}}
            \end{picture}}
\end{picture}}

\newcommand{\DNGGraphCircleInfinity}{\setlength{\unitlength}{0.75in}
\begin{picture}(1.45,0.35)
\put(0.15,0){\begin{picture}(1,0)
            \put(0,0.1){\circle*{0.05}}
            \put(-0.20,-0.05){\large $\gamma_{i}$}
            \put(1,0.1){\circle*{0.05}}
            \put(1.05,-0.05){\large $\gamma_{j}$}
            \put(0,0.1){\line(1,0){1}}
            \put(0.3,0.15){\CircleInteger{$\hspace*{-0.5mm}\infty$}}
            \end{picture}}
\end{picture}}

\newcommand{\ATwoGraphNoEdgeLabels}{\setlength{\unitlength}{0.75in}
\begin{picture}(1.4,0.25)
\put(0.15,0){\begin{picture}(1,0)
            \put(0,0.1){\circle*{0.05}}
            \put(-0.20,-0.05){\large $\gamma_{i}$}
            \put(1,0.1){\circle*{0.05}}
            \put(1.05,-0.05){\large $\gamma_{j}$}
            \put(0,0.1){\line(1,0){1}}
            \end{picture}}
\end{picture}}

\newcommand{\VertexForPosetsOrig}[2]{
\setlength{\unitlength}{1cm}
\begin{picture}(0,0)
\put(-0.38,0){
\begin{picture}(0,0)
\put(0,0){\circle*{0.15}} 
\put(-0.6,-0.1){\footnotesize $v_{#1}$}
\put(0.2,-0.1){\footnotesize $#2$}
\end{picture}
}
\end{picture}
}

\newcommand{\VertexForLatticeI}[1]{
\setlength{\unitlength}{1.5cm}
\begin{picture}(0,0)
\put(-0.25,0){
\begin{picture}(0,0)
\put(0,0){\circle*{0.1}} 
\put(-0.35,0.0){\footnotesize $\telt_{#1}$}
\end{picture}
}
\end{picture}
}

\newcommand{\NEEdgeLabelForLatticeI}[1]{
\setlength{\unitlength}{1.5cm}
\begin{picture}(0,0)
\put(-0.25,0){
\begin{picture}(0,0)
\put(0.4,0.4){\footnotesize #1} 
\end{picture}
}
\end{picture}
}

\newcommand{\NWEdgeLabelForLatticeI}[1]{
\setlength{\unitlength}{1.5cm}
\begin{picture}(0,0)
\put(-0.25,0){
\begin{picture}(0,0)
\put(-0.525,0.4){\footnotesize #1} 
\end{picture}
}
\end{picture}
}

\newcommand{\NEEdgeLabelForLatticeII}[1]{
\setlength{\unitlength}{1cm}
\begin{picture}(0,0)
\put(-0.25,0){
\begin{picture}(0,0)
\put(0.25,0.4){\footnotesize #1} 
\end{picture}
}
\end{picture}
}

\newcommand{\NWEdgeLabelForLatticeII}[1]{
\setlength{\unitlength}{1cm}
\begin{picture}(0,0)
\put(-0.25,0){
\begin{picture}(0,0)
\put(-0.75,0.4){\footnotesize #1} 
\end{picture}
}
\end{picture}
}

\newcommand{\VerticalEdgeLabelForLatticeI}[1]{
\setlength{\unitlength}{1.5cm}
\begin{picture}(0,0)
\put(-0.25,0){
\begin{picture}(0,0)
\put(-0.05,0.4){\footnotesize #1} 
\end{picture}
}
\end{picture}
}

\newcommand{\VertexIdeals}[2]{
\setlength{\unitlength}{1.5cm}
\begin{picture}(0,0)
\put(-0.25,0){
\begin{picture}(0,0)
\put(0,0){\circle*{0.1}} 
\put(#2,-0.05){\footnotesize $\langle #1 \rangle$}
\end{picture}
}
\end{picture}
}

\newcommand{\VertexIdealsEmpty}[1]{
\setlength{\unitlength}{1.5cm}
\begin{picture}(0,0)
\put(-0.25,0){
\begin{picture}(0,0)
\put(0,0){\circle*{0.1}} 
\put(#1,-0.05){\footnotesize $\emptyset$}
\end{picture}
}
\end{picture}
}

\newcommand{\RedNode}{
\setlength{\unitlength}{1cm}
\begin{picture}(0.5,0.5)
\put(0.01,0.225){\color{Red}\rule[-0.25mm]{.11mm}{.1mm}}
\put(0.02,0.225){\color{Red}\rule[-0.35mm]{.11mm}{.35mm}}
\put(0.03,0.225){\color{Red}\rule[-0.35mm]{.11mm}{.35mm}}
\put(0.04,0.225){\color{Red}\rule[-0.35mm]{.11mm}{.45mm}}
\put(0.05,0.225){\color{Red}\rule[-0.4mm]{.11mm}{.55mm}}
\put(0.06,0.225){\color{Red}\rule[-0.4mm]{.11mm}{.55mm}}
\put(0.07,0.225){\color{Red}\rule[-0.45mm]{.11mm}{.65mm}}
\put(0.08,0.225){\color{Red}\rule[-0.45mm]{.11mm}{.75mm}}
\put(0.09,0.225){\color{Red}\rule[-0.5mm]{.11mm}{.8mm}}
\put(0.10,0.225){\color{Red}\rule[-0.55mm]{.11mm}{.85mm}}
\put(0.11,0.225){\color{Red}\rule[-0.55mm]{.11mm}{.85mm}}
\put(0.12,0.225){\color{Red}\rule[-0.6mm]{.11mm}{.9mm}}
\put(0.13,0.225){\color{Red}\rule[-0.6mm]{.11mm}{.9mm}}
\put(0.14,0.225){\color{Red}\rule[-0.65mm]{.11mm}{.95mm}}
\put(0.15,0.225){\color{Red}\rule[-0.65mm]{.11mm}{.95mm}}
\put(0.16,0.225){\color{Red}\rule[-0.65mm]{.11mm}{.95mm}}
\put(0.17,0.225){\color{Red}\rule[-0.7mm]{.11mm}{1.05mm}}
\put(0.18,0.225){\color{Red}\rule[-0.7mm]{.11mm}{1.05mm}}
\put(0.19,0.225){\color{Red}\rule[-0.75mm]{.11mm}{1.1mm}}
\put(0.20,0.225){\color{Red}\rule[-0.75mm]{.11mm}{1.1mm}}
\put(0.21,0.225){\color{Red}\rule[-0.75mm]{.11mm}{1.1mm}}
\put(0.22,0.225){\color{Red}\rule[-0.8mm]{.11mm}{1.1mm}}
\put(0.23,0.225){\color{Red}\rule[-0.8mm]{.11mm}{1.1mm}}
\put(0.24,0.225){\color{Red}\rule[-0.8mm]{.11mm}{1.1mm}}
\put(0.25,0.225){\color{Red}\rule[-0.8mm]{.11mm}{1.1mm}}
\put(0.26,0.225){\color{Red}\rule[-0.8mm]{.11mm}{1.1mm}}
\put(0.27,0.225){\color{Red}\rule[-0.8mm]{.11mm}{1.1mm}}
\put(0.28,0.225){\color{Red}\rule[-0.85mm]{.11mm}{1.1mm}}
\put(0.29,0.225){\color{Red}\rule[-0.85mm]{.11mm}{1.1mm}}
\put(0.30,0.225){\color{Red}\rule[-0.85mm]{.11mm}{1.1mm}}
\put(0.31,0.225){\color{Red}\rule[-0.85mm]{.11mm}{1.1mm}}
\put(0.32,0.225){\color{Red}\rule[-0.85mm]{.11mm}{1.1mm}}
\put(0.33,0.225){\color{Red}\rule[-0.85mm]{.11mm}{1mm}}
\put(0.34,0.225){\color{Red}\rule[-0.85mm]{.11mm}{1mm}}
\put(0.35,0.225){\color{Red}\rule[-0.85mm]{.11mm}{0.95mm}}
\put(0.36,0.225){\color{Red}\rule[-0.85mm]{.11mm}{0.95mm}}
\put(0.37,0.225){\color{Red}\rule[-0.85mm]{.11mm}{0.95mm}}
\put(0.38,0.225){\color{Red}\rule[-0.85mm]{.11mm}{0.95mm}}
\put(0.39,0.225){\color{Red}\rule[-0.8mm]{.11mm}{0.85mm}}
\put(0.40,0.225){\color{Red}\rule[-0.8mm]{.11mm}{0.85mm}}
\put(0.41,0.225){\color{Red}\rule[-0.8mm]{.11mm}{0.8mm}}
\put(0.42,0.225){\color{Red}\rule[-0.7mm]{.11mm}{0.65mm}}
\put(0.43,0.225){\color{Red}\rule[-0.65mm]{.11mm}{0.6mm}}
\put(0.44,0.225){\color{Red}\rule[-0.6mm]{.11mm}{0.55mm}}
\put(0.45,0.225){\color{Red}\rule[-0.55mm]{.11mm}{0.4mm}}
\put(0.46,0.225){\color{Red}\rule[-0.5mm]{.11mm}{0.3mm}}
\put(0.47,0.225){\color{Red}\rule[-0.45mm]{.11mm}{0.25mm}}
\put(0.48,0.225){\color{Red}\rule[-0.4mm]{.11mm}{0.2mm}}
\put(0,0.2){\qbezier(0,0)(0.125,0.125)(0.5,0)}
\put(0,0.2){\qbezier(0,0)(0.375,-0.125)(0.5,0)}
\put(0,0.2){\circle*{0.0135}}
\put(0.5,0.2){\circle*{0.0135}}
\end{picture}}
\newcommand{\CyanNode}{
\setlength{\unitlength}{1cm}
\begin{picture}(0.5,0.5)
\put(0.01,0.225){\color{Cyan}\rule[-0.25mm]{.11mm}{.1mm}}
\put(0.02,0.225){\color{Cyan}\rule[-0.35mm]{.11mm}{.35mm}}
\put(0.03,0.225){\color{Cyan}\rule[-0.35mm]{.11mm}{.35mm}}
\put(0.04,0.225){\color{Cyan}\rule[-0.35mm]{.11mm}{.45mm}}
\put(0.05,0.225){\color{Cyan}\rule[-0.4mm]{.11mm}{.55mm}}
\put(0.06,0.225){\color{Cyan}\rule[-0.4mm]{.11mm}{.55mm}}
\put(0.07,0.225){\color{Cyan}\rule[-0.45mm]{.11mm}{.65mm}}
\put(0.08,0.225){\color{Cyan}\rule[-0.45mm]{.11mm}{.75mm}}
\put(0.09,0.225){\color{Cyan}\rule[-0.5mm]{.11mm}{.8mm}}
\put(0.10,0.225){\color{Cyan}\rule[-0.55mm]{.11mm}{.85mm}}
\put(0.11,0.225){\color{Cyan}\rule[-0.55mm]{.11mm}{.85mm}}
\put(0.12,0.225){\color{Cyan}\rule[-0.6mm]{.11mm}{.9mm}}
\put(0.13,0.225){\color{Cyan}\rule[-0.6mm]{.11mm}{.9mm}}
\put(0.14,0.225){\color{Cyan}\rule[-0.65mm]{.11mm}{.95mm}}
\put(0.15,0.225){\color{Cyan}\rule[-0.65mm]{.11mm}{.95mm}}
\put(0.16,0.225){\color{Cyan}\rule[-0.65mm]{.11mm}{.95mm}}
\put(0.17,0.225){\color{Cyan}\rule[-0.7mm]{.11mm}{1.05mm}}
\put(0.18,0.225){\color{Cyan}\rule[-0.7mm]{.11mm}{1.05mm}}
\put(0.19,0.225){\color{Cyan}\rule[-0.75mm]{.11mm}{1.1mm}}
\put(0.20,0.225){\color{Cyan}\rule[-0.75mm]{.11mm}{1.1mm}}
\put(0.21,0.225){\color{Cyan}\rule[-0.75mm]{.11mm}{1.1mm}}
\put(0.22,0.225){\color{Cyan}\rule[-0.8mm]{.11mm}{1.1mm}}
\put(0.23,0.225){\color{Cyan}\rule[-0.8mm]{.11mm}{1.1mm}}
\put(0.24,0.225){\color{Cyan}\rule[-0.8mm]{.11mm}{1.1mm}}
\put(0.25,0.225){\color{Cyan}\rule[-0.8mm]{.11mm}{1.1mm}}
\put(0.26,0.225){\color{Cyan}\rule[-0.8mm]{.11mm}{1.1mm}}
\put(0.27,0.225){\color{Cyan}\rule[-0.8mm]{.11mm}{1.1mm}}
\put(0.28,0.225){\color{Cyan}\rule[-0.85mm]{.11mm}{1.1mm}}
\put(0.29,0.225){\color{Cyan}\rule[-0.85mm]{.11mm}{1.1mm}}
\put(0.30,0.225){\color{Cyan}\rule[-0.85mm]{.11mm}{1.1mm}}
\put(0.31,0.225){\color{Cyan}\rule[-0.85mm]{.11mm}{1.1mm}}
\put(0.32,0.225){\color{Cyan}\rule[-0.85mm]{.11mm}{1.1mm}}
\put(0.33,0.225){\color{Cyan}\rule[-0.85mm]{.11mm}{1mm}}
\put(0.34,0.225){\color{Cyan}\rule[-0.85mm]{.11mm}{1mm}}
\put(0.35,0.225){\color{Cyan}\rule[-0.85mm]{.11mm}{0.95mm}}
\put(0.36,0.225){\color{Cyan}\rule[-0.85mm]{.11mm}{0.95mm}}
\put(0.37,0.225){\color{Cyan}\rule[-0.85mm]{.11mm}{0.95mm}}
\put(0.38,0.225){\color{Cyan}\rule[-0.85mm]{.11mm}{0.95mm}}
\put(0.39,0.225){\color{Cyan}\rule[-0.8mm]{.11mm}{0.85mm}}
\put(0.40,0.225){\color{Cyan}\rule[-0.8mm]{.11mm}{0.85mm}}
\put(0.41,0.225){\color{Cyan}\rule[-0.8mm]{.11mm}{0.8mm}}
\put(0.42,0.225){\color{Cyan}\rule[-0.7mm]{.11mm}{0.65mm}}
\put(0.43,0.225){\color{Cyan}\rule[-0.65mm]{.11mm}{0.6mm}}
\put(0.44,0.225){\color{Cyan}\rule[-0.6mm]{.11mm}{0.55mm}}
\put(0.45,0.225){\color{Cyan}\rule[-0.55mm]{.11mm}{0.4mm}}
\put(0.46,0.225){\color{Cyan}\rule[-0.5mm]{.11mm}{0.3mm}}
\put(0.47,0.225){\color{Cyan}\rule[-0.45mm]{.11mm}{0.25mm}}
\put(0.48,0.225){\color{Cyan}\rule[-0.4mm]{.11mm}{0.2mm}}
\put(0,0.2){\qbezier(0,0)(0.125,0.125)(0.5,0)}
\put(0,0.2){\qbezier(0,0)(0.375,-0.125)(0.5,0)}
\put(0,0.2){\circle*{0.0135}}
\put(0.5,0.2){\circle*{0.0135}}
\end{picture}}
\newcommand{\PurpleNode}{
\setlength{\unitlength}{1cm}
\begin{picture}(0.5,0.5)
\put(0.01,0.225){\color{Purple}\rule[-0.25mm]{.11mm}{.1mm}}
\put(0.02,0.225){\color{Purple}\rule[-0.35mm]{.11mm}{.35mm}}
\put(0.03,0.225){\color{Purple}\rule[-0.35mm]{.11mm}{.35mm}}
\put(0.04,0.225){\color{Purple}\rule[-0.35mm]{.11mm}{.45mm}}
\put(0.05,0.225){\color{Purple}\rule[-0.4mm]{.11mm}{.55mm}}
\put(0.06,0.225){\color{Purple}\rule[-0.4mm]{.11mm}{.55mm}}
\put(0.07,0.225){\color{Purple}\rule[-0.45mm]{.11mm}{.65mm}}
\put(0.08,0.225){\color{Purple}\rule[-0.45mm]{.11mm}{.75mm}}
\put(0.09,0.225){\color{Purple}\rule[-0.5mm]{.11mm}{.8mm}}
\put(0.10,0.225){\color{Purple}\rule[-0.55mm]{.11mm}{.85mm}}
\put(0.11,0.225){\color{Purple}\rule[-0.55mm]{.11mm}{.85mm}}
\put(0.12,0.225){\color{Purple}\rule[-0.6mm]{.11mm}{.9mm}}
\put(0.13,0.225){\color{Purple}\rule[-0.6mm]{.11mm}{.9mm}}
\put(0.14,0.225){\color{Purple}\rule[-0.65mm]{.11mm}{.95mm}}
\put(0.15,0.225){\color{Purple}\rule[-0.65mm]{.11mm}{.95mm}}
\put(0.16,0.225){\color{Purple}\rule[-0.65mm]{.11mm}{.95mm}}
\put(0.17,0.225){\color{Purple}\rule[-0.7mm]{.11mm}{1.05mm}}
\put(0.18,0.225){\color{Purple}\rule[-0.7mm]{.11mm}{1.05mm}}
\put(0.19,0.225){\color{Purple}\rule[-0.75mm]{.11mm}{1.1mm}}
\put(0.20,0.225){\color{Purple}\rule[-0.75mm]{.11mm}{1.1mm}}
\put(0.21,0.225){\color{Purple}\rule[-0.75mm]{.11mm}{1.1mm}}
\put(0.22,0.225){\color{Purple}\rule[-0.8mm]{.11mm}{1.1mm}}
\put(0.23,0.225){\color{Purple}\rule[-0.8mm]{.11mm}{1.1mm}}
\put(0.24,0.225){\color{Purple}\rule[-0.8mm]{.11mm}{1.1mm}}
\put(0.25,0.225){\color{Purple}\rule[-0.8mm]{.11mm}{1.1mm}}
\put(0.26,0.225){\color{Purple}\rule[-0.8mm]{.11mm}{1.1mm}}
\put(0.27,0.225){\color{Purple}\rule[-0.8mm]{.11mm}{1.1mm}}
\put(0.28,0.225){\color{Purple}\rule[-0.85mm]{.11mm}{1.1mm}}
\put(0.29,0.225){\color{Purple}\rule[-0.85mm]{.11mm}{1.1mm}}
\put(0.30,0.225){\color{Purple}\rule[-0.85mm]{.11mm}{1.1mm}}
\put(0.31,0.225){\color{Purple}\rule[-0.85mm]{.11mm}{1.1mm}}
\put(0.32,0.225){\color{Purple}\rule[-0.85mm]{.11mm}{1.1mm}}
\put(0.33,0.225){\color{Purple}\rule[-0.85mm]{.11mm}{1mm}}
\put(0.34,0.225){\color{Purple}\rule[-0.85mm]{.11mm}{1mm}}
\put(0.35,0.225){\color{Purple}\rule[-0.85mm]{.11mm}{0.95mm}}
\put(0.36,0.225){\color{Purple}\rule[-0.85mm]{.11mm}{0.95mm}}
\put(0.37,0.225){\color{Purple}\rule[-0.85mm]{.11mm}{0.95mm}}
\put(0.38,0.225){\color{Purple}\rule[-0.85mm]{.11mm}{0.95mm}}
\put(0.39,0.225){\color{Purple}\rule[-0.8mm]{.11mm}{0.85mm}}
\put(0.40,0.225){\color{Purple}\rule[-0.8mm]{.11mm}{0.85mm}}
\put(0.41,0.225){\color{Purple}\rule[-0.8mm]{.11mm}{0.8mm}}
\put(0.42,0.225){\color{Purple}\rule[-0.7mm]{.11mm}{0.65mm}}
\put(0.43,0.225){\color{Purple}\rule[-0.65mm]{.11mm}{0.6mm}}
\put(0.44,0.225){\color{Purple}\rule[-0.6mm]{.11mm}{0.55mm}}
\put(0.45,0.225){\color{Purple}\rule[-0.55mm]{.11mm}{0.4mm}}
\put(0.46,0.225){\color{Purple}\rule[-0.5mm]{.11mm}{0.3mm}}
\put(0.47,0.225){\color{Purple}\rule[-0.45mm]{.11mm}{0.25mm}}
\put(0.48,0.225){\color{Purple}\rule[-0.4mm]{.11mm}{0.2mm}}
\put(0,0.2){\qbezier(0,0)(0.125,0.125)(0.5,0)}
\put(0,0.2){\qbezier(0,0)(0.375,-0.125)(0.5,0)}
\put(0,0.2){\circle*{0.0135}}
\put(0.5,0.2){\circle*{0.0135}}
\end{picture}}

\newcommand{\OliveGreenNode}{
\setlength{\unitlength}{1cm}
\begin{picture}(0.5,0.5)
\put(0.01,0.225){\color{OliveGreen}\rule[-0.25mm]{.11mm}{.1mm}}
\put(0.02,0.225){\color{OliveGreen}\rule[-0.35mm]{.11mm}{.35mm}}
\put(0.03,0.225){\color{OliveGreen}\rule[-0.35mm]{.11mm}{.35mm}}
\put(0.04,0.225){\color{OliveGreen}\rule[-0.35mm]{.11mm}{.45mm}}
\put(0.05,0.225){\color{OliveGreen}\rule[-0.4mm]{.11mm}{.55mm}}
\put(0.06,0.225){\color{OliveGreen}\rule[-0.4mm]{.11mm}{.55mm}}
\put(0.07,0.225){\color{OliveGreen}\rule[-0.45mm]{.11mm}{.65mm}}
\put(0.08,0.225){\color{OliveGreen}\rule[-0.45mm]{.11mm}{.75mm}}
\put(0.09,0.225){\color{OliveGreen}\rule[-0.5mm]{.11mm}{.8mm}}
\put(0.10,0.225){\color{OliveGreen}\rule[-0.55mm]{.11mm}{.85mm}}
\put(0.11,0.225){\color{OliveGreen}\rule[-0.55mm]{.11mm}{.85mm}}
\put(0.12,0.225){\color{OliveGreen}\rule[-0.6mm]{.11mm}{.9mm}}
\put(0.13,0.225){\color{OliveGreen}\rule[-0.6mm]{.11mm}{.9mm}}
\put(0.14,0.225){\color{OliveGreen}\rule[-0.65mm]{.11mm}{.95mm}}
\put(0.15,0.225){\color{OliveGreen}\rule[-0.65mm]{.11mm}{.95mm}}
\put(0.16,0.225){\color{OliveGreen}\rule[-0.65mm]{.11mm}{.95mm}}
\put(0.17,0.225){\color{OliveGreen}\rule[-0.7mm]{.11mm}{1.05mm}}
\put(0.18,0.225){\color{OliveGreen}\rule[-0.7mm]{.11mm}{1.05mm}}
\put(0.19,0.225){\color{OliveGreen}\rule[-0.75mm]{.11mm}{1.1mm}}
\put(0.20,0.225){\color{OliveGreen}\rule[-0.75mm]{.11mm}{1.1mm}}
\put(0.21,0.225){\color{OliveGreen}\rule[-0.75mm]{.11mm}{1.1mm}}
\put(0.22,0.225){\color{OliveGreen}\rule[-0.8mm]{.11mm}{1.1mm}}
\put(0.23,0.225){\color{OliveGreen}\rule[-0.8mm]{.11mm}{1.1mm}}
\put(0.24,0.225){\color{OliveGreen}\rule[-0.8mm]{.11mm}{1.1mm}}
\put(0.25,0.225){\color{OliveGreen}\rule[-0.8mm]{.11mm}{1.1mm}}
\put(0.26,0.225){\color{OliveGreen}\rule[-0.8mm]{.11mm}{1.1mm}}
\put(0.27,0.225){\color{OliveGreen}\rule[-0.8mm]{.11mm}{1.1mm}}
\put(0.28,0.225){\color{OliveGreen}\rule[-0.85mm]{.11mm}{1.1mm}}
\put(0.29,0.225){\color{OliveGreen}\rule[-0.85mm]{.11mm}{1.1mm}}
\put(0.30,0.225){\color{OliveGreen}\rule[-0.85mm]{.11mm}{1.1mm}}
\put(0.31,0.225){\color{OliveGreen}\rule[-0.85mm]{.11mm}{1.1mm}}
\put(0.32,0.225){\color{OliveGreen}\rule[-0.85mm]{.11mm}{1.1mm}}
\put(0.33,0.225){\color{OliveGreen}\rule[-0.85mm]{.11mm}{1mm}}
\put(0.34,0.225){\color{OliveGreen}\rule[-0.85mm]{.11mm}{1mm}}
\put(0.35,0.225){\color{OliveGreen}\rule[-0.85mm]{.11mm}{0.95mm}}
\put(0.36,0.225){\color{OliveGreen}\rule[-0.85mm]{.11mm}{0.95mm}}
\put(0.37,0.225){\color{OliveGreen}\rule[-0.85mm]{.11mm}{0.95mm}}
\put(0.38,0.225){\color{OliveGreen}\rule[-0.85mm]{.11mm}{0.95mm}}
\put(0.39,0.225){\color{OliveGreen}\rule[-0.8mm]{.11mm}{0.85mm}}
\put(0.40,0.225){\color{OliveGreen}\rule[-0.8mm]{.11mm}{0.85mm}}
\put(0.41,0.225){\color{OliveGreen}\rule[-0.8mm]{.11mm}{0.8mm}}
\put(0.42,0.225){\color{OliveGreen}\rule[-0.7mm]{.11mm}{0.65mm}}
\put(0.43,0.225){\color{OliveGreen}\rule[-0.65mm]{.11mm}{0.6mm}}
\put(0.44,0.225){\color{OliveGreen}\rule[-0.6mm]{.11mm}{0.55mm}}
\put(0.45,0.225){\color{OliveGreen}\rule[-0.55mm]{.11mm}{0.4mm}}
\put(0.46,0.225){\color{OliveGreen}\rule[-0.5mm]{.11mm}{0.3mm}}
\put(0.47,0.225){\color{OliveGreen}\rule[-0.45mm]{.11mm}{0.25mm}}
\put(0.48,0.225){\color{OliveGreen}\rule[-0.4mm]{.11mm}{0.2mm}}
\put(0,0.2){\qbezier(0,0)(0.125,0.125)(0.5,0)}
\put(0,0.2){\qbezier(0,0)(0.375,-0.125)(0.5,0)}
\put(0,0.2){\circle*{0.0135}}
\put(0.5,0.2){\circle*{0.0135}}
\end{picture}}
\newcommand{\SpringGreenNode}{
\setlength{\unitlength}{1cm}
\begin{picture}(0.5,0.5)
\put(0.01,0.225){\color{SpringGreen}\rule[-0.25mm]{.11mm}{.1mm}}
\put(0.02,0.225){\color{SpringGreen}\rule[-0.35mm]{.11mm}{.35mm}}
\put(0.03,0.225){\color{SpringGreen}\rule[-0.35mm]{.11mm}{.35mm}}
\put(0.04,0.225){\color{SpringGreen}\rule[-0.35mm]{.11mm}{.45mm}}
\put(0.05,0.225){\color{SpringGreen}\rule[-0.4mm]{.11mm}{.55mm}}
\put(0.06,0.225){\color{SpringGreen}\rule[-0.4mm]{.11mm}{.55mm}}
\put(0.07,0.225){\color{SpringGreen}\rule[-0.45mm]{.11mm}{.65mm}}
\put(0.08,0.225){\color{SpringGreen}\rule[-0.45mm]{.11mm}{.75mm}}
\put(0.09,0.225){\color{SpringGreen}\rule[-0.5mm]{.11mm}{.8mm}}
\put(0.10,0.225){\color{SpringGreen}\rule[-0.55mm]{.11mm}{.85mm}}
\put(0.11,0.225){\color{SpringGreen}\rule[-0.55mm]{.11mm}{.85mm}}
\put(0.12,0.225){\color{SpringGreen}\rule[-0.6mm]{.11mm}{.9mm}}
\put(0.13,0.225){\color{SpringGreen}\rule[-0.6mm]{.11mm}{.9mm}}
\put(0.14,0.225){\color{SpringGreen}\rule[-0.65mm]{.11mm}{.95mm}}
\put(0.15,0.225){\color{SpringGreen}\rule[-0.65mm]{.11mm}{.95mm}}
\put(0.16,0.225){\color{SpringGreen}\rule[-0.65mm]{.11mm}{.95mm}}
\put(0.17,0.225){\color{SpringGreen}\rule[-0.7mm]{.11mm}{1.05mm}}
\put(0.18,0.225){\color{SpringGreen}\rule[-0.7mm]{.11mm}{1.05mm}}
\put(0.19,0.225){\color{SpringGreen}\rule[-0.75mm]{.11mm}{1.1mm}}
\put(0.20,0.225){\color{SpringGreen}\rule[-0.75mm]{.11mm}{1.1mm}}
\put(0.21,0.225){\color{SpringGreen}\rule[-0.75mm]{.11mm}{1.1mm}}
\put(0.22,0.225){\color{SpringGreen}\rule[-0.8mm]{.11mm}{1.1mm}}
\put(0.23,0.225){\color{SpringGreen}\rule[-0.8mm]{.11mm}{1.1mm}}
\put(0.24,0.225){\color{SpringGreen}\rule[-0.8mm]{.11mm}{1.1mm}}
\put(0.25,0.225){\color{SpringGreen}\rule[-0.8mm]{.11mm}{1.1mm}}
\put(0.26,0.225){\color{SpringGreen}\rule[-0.8mm]{.11mm}{1.1mm}}
\put(0.27,0.225){\color{SpringGreen}\rule[-0.8mm]{.11mm}{1.1mm}}
\put(0.28,0.225){\color{SpringGreen}\rule[-0.85mm]{.11mm}{1.1mm}}
\put(0.29,0.225){\color{SpringGreen}\rule[-0.85mm]{.11mm}{1.1mm}}
\put(0.30,0.225){\color{SpringGreen}\rule[-0.85mm]{.11mm}{1.1mm}}
\put(0.31,0.225){\color{SpringGreen}\rule[-0.85mm]{.11mm}{1.1mm}}
\put(0.32,0.225){\color{SpringGreen}\rule[-0.85mm]{.11mm}{1.1mm}}
\put(0.33,0.225){\color{SpringGreen}\rule[-0.85mm]{.11mm}{1mm}}
\put(0.34,0.225){\color{SpringGreen}\rule[-0.85mm]{.11mm}{1mm}}
\put(0.35,0.225){\color{SpringGreen}\rule[-0.85mm]{.11mm}{0.95mm}}
\put(0.36,0.225){\color{SpringGreen}\rule[-0.85mm]{.11mm}{0.95mm}}
\put(0.37,0.225){\color{SpringGreen}\rule[-0.85mm]{.11mm}{0.95mm}}
\put(0.38,0.225){\color{SpringGreen}\rule[-0.85mm]{.11mm}{0.95mm}}
\put(0.39,0.225){\color{SpringGreen}\rule[-0.8mm]{.11mm}{0.85mm}}
\put(0.40,0.225){\color{SpringGreen}\rule[-0.8mm]{.11mm}{0.85mm}}
\put(0.41,0.225){\color{SpringGreen}\rule[-0.8mm]{.11mm}{0.8mm}}
\put(0.42,0.225){\color{SpringGreen}\rule[-0.7mm]{.11mm}{0.65mm}}
\put(0.43,0.225){\color{SpringGreen}\rule[-0.65mm]{.11mm}{0.6mm}}
\put(0.44,0.225){\color{SpringGreen}\rule[-0.6mm]{.11mm}{0.55mm}}
\put(0.45,0.225){\color{SpringGreen}\rule[-0.55mm]{.11mm}{0.4mm}}
\put(0.46,0.225){\color{SpringGreen}\rule[-0.5mm]{.11mm}{0.3mm}}
\put(0.47,0.225){\color{SpringGreen}\rule[-0.45mm]{.11mm}{0.25mm}}
\put(0.48,0.225){\color{SpringGreen}\rule[-0.4mm]{.11mm}{0.2mm}}
\put(0,0.2){\qbezier(0,0)(0.125,0.125)(0.5,0)}
\put(0,0.2){\qbezier(0,0)(0.375,-0.125)(0.5,0)}
\put(0,0.2){\circle*{0.0135}}
\put(0.5,0.2){\circle*{0.0135}}
\end{picture}}
\newcommand{\BurntOrangeNode}{
\setlength{\unitlength}{1cm}
\begin{picture}(0.5,0.5)
\put(0.01,0.225){\color{BurntOrange}\rule[-0.25mm]{.11mm}{.1mm}}
\put(0.02,0.225){\color{BurntOrange}\rule[-0.35mm]{.11mm}{.35mm}}
\put(0.03,0.225){\color{BurntOrange}\rule[-0.35mm]{.11mm}{.35mm}}
\put(0.04,0.225){\color{BurntOrange}\rule[-0.35mm]{.11mm}{.45mm}}
\put(0.05,0.225){\color{BurntOrange}\rule[-0.4mm]{.11mm}{.55mm}}
\put(0.06,0.225){\color{BurntOrange}\rule[-0.4mm]{.11mm}{.55mm}}
\put(0.07,0.225){\color{BurntOrange}\rule[-0.45mm]{.11mm}{.65mm}}
\put(0.08,0.225){\color{BurntOrange}\rule[-0.45mm]{.11mm}{.75mm}}
\put(0.09,0.225){\color{BurntOrange}\rule[-0.5mm]{.11mm}{.8mm}}
\put(0.10,0.225){\color{BurntOrange}\rule[-0.55mm]{.11mm}{.85mm}}
\put(0.11,0.225){\color{BurntOrange}\rule[-0.55mm]{.11mm}{.85mm}}
\put(0.12,0.225){\color{BurntOrange}\rule[-0.6mm]{.11mm}{.9mm}}
\put(0.13,0.225){\color{BurntOrange}\rule[-0.6mm]{.11mm}{.9mm}}
\put(0.14,0.225){\color{BurntOrange}\rule[-0.65mm]{.11mm}{.95mm}}
\put(0.15,0.225){\color{BurntOrange}\rule[-0.65mm]{.11mm}{.95mm}}
\put(0.16,0.225){\color{BurntOrange}\rule[-0.65mm]{.11mm}{.95mm}}
\put(0.17,0.225){\color{BurntOrange}\rule[-0.7mm]{.11mm}{1.05mm}}
\put(0.18,0.225){\color{BurntOrange}\rule[-0.7mm]{.11mm}{1.05mm}}
\put(0.19,0.225){\color{BurntOrange}\rule[-0.75mm]{.11mm}{1.1mm}}
\put(0.20,0.225){\color{BurntOrange}\rule[-0.75mm]{.11mm}{1.1mm}}
\put(0.21,0.225){\color{BurntOrange}\rule[-0.75mm]{.11mm}{1.1mm}}
\put(0.22,0.225){\color{BurntOrange}\rule[-0.8mm]{.11mm}{1.1mm}}
\put(0.23,0.225){\color{BurntOrange}\rule[-0.8mm]{.11mm}{1.1mm}}
\put(0.24,0.225){\color{BurntOrange}\rule[-0.8mm]{.11mm}{1.1mm}}
\put(0.25,0.225){\color{BurntOrange}\rule[-0.8mm]{.11mm}{1.1mm}}
\put(0.26,0.225){\color{BurntOrange}\rule[-0.8mm]{.11mm}{1.1mm}}
\put(0.27,0.225){\color{BurntOrange}\rule[-0.8mm]{.11mm}{1.1mm}}
\put(0.28,0.225){\color{BurntOrange}\rule[-0.85mm]{.11mm}{1.1mm}}
\put(0.29,0.225){\color{BurntOrange}\rule[-0.85mm]{.11mm}{1.1mm}}
\put(0.30,0.225){\color{BurntOrange}\rule[-0.85mm]{.11mm}{1.1mm}}
\put(0.31,0.225){\color{BurntOrange}\rule[-0.85mm]{.11mm}{1.1mm}}
\put(0.32,0.225){\color{BurntOrange}\rule[-0.85mm]{.11mm}{1.1mm}}
\put(0.33,0.225){\color{BurntOrange}\rule[-0.85mm]{.11mm}{1mm}}
\put(0.34,0.225){\color{BurntOrange}\rule[-0.85mm]{.11mm}{1mm}}
\put(0.35,0.225){\color{BurntOrange}\rule[-0.85mm]{.11mm}{0.95mm}}
\put(0.36,0.225){\color{BurntOrange}\rule[-0.85mm]{.11mm}{0.95mm}}
\put(0.37,0.225){\color{BurntOrange}\rule[-0.85mm]{.11mm}{0.95mm}}
\put(0.38,0.225){\color{BurntOrange}\rule[-0.85mm]{.11mm}{0.95mm}}
\put(0.39,0.225){\color{BurntOrange}\rule[-0.8mm]{.11mm}{0.85mm}}
\put(0.40,0.225){\color{BurntOrange}\rule[-0.8mm]{.11mm}{0.85mm}}
\put(0.41,0.225){\color{BurntOrange}\rule[-0.8mm]{.11mm}{0.8mm}}
\put(0.42,0.225){\color{BurntOrange}\rule[-0.7mm]{.11mm}{0.65mm}}
\put(0.43,0.225){\color{BurntOrange}\rule[-0.65mm]{.11mm}{0.6mm}}
\put(0.44,0.225){\color{BurntOrange}\rule[-0.6mm]{.11mm}{0.55mm}}
\put(0.45,0.225){\color{BurntOrange}\rule[-0.55mm]{.11mm}{0.4mm}}
\put(0.46,0.225){\color{BurntOrange}\rule[-0.5mm]{.11mm}{0.3mm}}
\put(0.47,0.225){\color{BurntOrange}\rule[-0.45mm]{.11mm}{0.25mm}}
\put(0.48,0.225){\color{BurntOrange}\rule[-0.4mm]{.11mm}{0.2mm}}
\put(0,0.2){\qbezier(0,0)(0.125,0.125)(0.5,0)}
\put(0,0.2){\qbezier(0,0)(0.375,-0.125)(0.5,0)}
\put(0,0.2){\circle*{0.0135}}
\put(0.5,0.2){\circle*{0.0135}}
\end{picture}}

\begin{document}

\vspace*{-0.7in}
\hfill \parbox{2in}{\hfill {\scriptsize This version: May 7, 2022}}

\vspace*{-0.1in}
\hfill{\scriptsize 1st version: November 26, 2018}

\begin{center}
{\large {\bf Finite diamond-colored modular and distributive lattices}}\\
{\large {\bf with applications to combinatorial Lie representation theory}} 

\vspace*{0.05in}
\renewcommand{\thefootnote}{1}
Robert G.\ Donnelly\footnote{Email: {\tt rob.donnelly@murraystate.edu}}\\
Department of Mathematics and Statistics, Murray State University, Murray, KY 42071
\end{center} 

\vspace*{-0.25in}
\begin{abstract}
A finite modular or distributive lattice is `diamond-colored' if its covering digraph edges are colored in such a way that, within any diamond of edges, parallel edges have the same color. 
For us, a `diamond' is any four-element Boolean sublattice corresponding to a four-element subgraph of the given covering digraph. 
Such lattices arise naturally in combinatorial representation theory, particularly in the study of poset models for semisimple Lie algebra representations in finite dimensions and their companion Weyl group symmetric functions. 
Here we gather in one place some elementary but foundational combinatorial results concerning these lattice structures; this presentation includes some new objects/results as well as some new interpretations of classical objects/results. 
We extend these ideas in our development of a new approach to representing diamond-colored distributive lattices as partial orderings of certain integer arrays, which we call `ideal arrays', over vertex-colored posets. 

We also develop points of contact between diamond-colored modular $\!\!$/$\!\!$ distributive lattices and the aforementioned algebraic contexts.  
These points of contact are rendered here using a series of finiteness problems, all of whose answers are the Dynkin diagrams corresponding to the finite-dimensional simple Lie algebras over the complex numbers. 
We believe that some of these problems, and some of our approaches to solving these problems, are new. 
This work culminates in our statement of a diverse collection of equivalences regarding Dynkin diagrams that dates back to the seminal work of W.\ Killing and E.\ Cartan in which they classified the finite-dimensional simple Lie algebras over the complex numbers. 
We hereby nominate {\em La Florado Klasado} ({\em LFK})-- `The Flowering Classification' in Esperanto -- as the generic name for this multifaceted classification result. 

Within this environment, we classify those posets with a bi-coloring of vertices whose associated diamond-colored distributive lattices can serve, in a certain sense, as models for analogs of Schur functions for the rank two Weyl groups.  
In later sections, we apply our perspective to a new study of minuscule splitting distributive lattices (typically called minuscule lattices) and their associated minuscule compression posets. 
We also find all of the diamond-colored modular lattices that serve naturally as models for the quasi-minuscule representations $\!\!$/$\!\!$ Weyl symmetric functions. 
We close with a demonstration that diamond-colored distributive lattices of ideal arrays associated with a kind of stacking of minuscule compression posets enjoy a certain structuring property that is peculiar to {\em LFK}. 
\begin{center}

\ 
\vspace*{-0.1in}


\vspace*{0.1in}
\end{center}
\end{abstract}

\begin{center}
\underline{\hspace*{4in}}
\end{center}

\vfill

\newpage
\def\abstractname{Table of contents}
\begin{abstract}
\begin{center}
\parbox{5in}{
\hspace*{-0.1in}\IntroSection. Introduction\dotfill 4\\%
\hspace*{0.5in} \\ 
\hspace*{-0.25in}\fbox{Part I}\ \ A gentle introduction to order-theoretic aspects of diamond-coloring\dotfill 9\\
\hspace*{-0.1in}\RoutineSection. Rudiments\dotfill 9\\%
\hspace*{0.1in}\RoutineSection.1. Some useful directed graph manipulations\dotfill 9\\%
\hspace*{0.1in}\RoutineSection.2. Conventional aspects of partially ordered sets\dotfill 10\\%
\hspace*{0.1in}\RoutineSection.3. Lattice basics\dotfill 15\\%
\hspace*{0.1in}\RoutineSection.4. An initial consequence of diamond-coloring\dotfill 17\\%
\hspace*{0.1in}\RoutineSection.5. A first look at some distinguished features of DCDL's\dotfill 18\\%
\hspace*{0.1in}\RoutineSection.6. A fundamental theorem\dotfill 22\\%
\hspace*{-0.1in}\SubstructureSection. Substructures\dotfill 25\\%
\hspace*{0.1in}\SubstructureSection.1. Types of subposets\dotfill 25\\%
\hspace*{0.1in}\SubstructureSection.2. Sublattices\dotfill 25\\%
\hspace*{0.1in}\SubstructureSection.3. Full-length sublattices\dotfill 26\\%
\hspace*{0.1in}\SubstructureSection.4. Boolean sublattices\dotfill 29\\%
\hspace*{0.1in}\SubstructureSection.5. $J$-components of DCML's and DCDL's\dotfill 31\\%
\hspace*{0.5in} \\ 
\hspace*{-0.25in}\fbox{Part II}\ \  Diamond-coloring and some new objects, methods, and results\dotfill 45\\
\hspace*{-0.1in}\GStructureIntroSection. Integral embryophytes and sturdiness\dotfill 45\\%
\hspace*{0.1in}\GStructureIntroSection.1. Embryophytes\dotfill 45\\%
\hspace*{0.1in}\GStructureIntroSection.2. Working with embryophytes\dotfill 46\\%
\hspace*{0.1in}\GStructureIntroSection.3. Coxeter and Dynkin\dotfill 46\\%
\hspace*{0.1in}\GStructureIntroSection.4. Sturdy posets\dotfill 47\\%
\hspace*{0.1in}\GStructureIntroSection.5. Sturdy DCML's\dotfill 49\\%
\hspace*{0.1in}\GStructureIntroSection.6. Some open problems\dotfill 50\\%
\hspace*{-0.1in}\FrameSection. Scaffolds, skew-stacks, and DCDL's\dotfill 51\\%
\hspace*{0.1in}\FrameSection.1. Scaffolds\dotfill 52\\%
\hspace*{0.1in}\FrameSection.2. Skew-stacks\dotfill 53\\%
\hspace*{-0.1in}\GridPosetSection. Two-color grid posets\dotfill 62\\%
\hspace*{0.1in}\GridPosetSection.1. Grid posets\dotfill 62\\%
\hspace*{0.1in}\GridPosetSection.2. Two-color grid posets: definitions and examples\dotfill 63\\%
\hspace*{0.1in}\GridPosetSection.3. Decomposable/indecomposable two-color grid posets\dotfill 64\\%
\hspace*{0.1in}\GridPosetSection.4. Fundamental and semistandard posets in two colors\dotfill 70\\%
\hspace*{0.1in}\GridPosetSection.5. A combinatorial characterization of rank two fundamental posets\dotfill 75\\%
\hspace*{0.1in}\GridPosetSection.6. A combinatorial characterization of rank two semistandard posets\dotfill 76\\%
\hspace*{-0.1in}\TwoColorSturdySection.\ A sturdiness result from the theory of two-color grid posets\dotfill 78\\%
\hspace*{0.1in}\TwoColorSturdySection.1. Two-color subordinates of certain skew-stacks\dotfill 78\\%
\hspace*{0.1in}\TwoColorSturdySection.2. Sturdy DCDL's from certain scaffolds and skew-stacks\dotfill 80\\%
\hspace*{0.5in} \\ 
\hspace*{-0.25in}\fbox{Part III}\ \  A primer on algebraic contexts\dotfill 81\\ 
\hspace*{-0.1in}\WeylSection. Finitistic aspects of algebraic structures related to integral embryophytes\dotfill 81\\%
\hspace*{0.1in}\WeylSection.1. Basic set-up\dotfill 81\\%
\hspace*{0.1in}\WeylSection.2. The Weyl group\dotfill 81\\%
\hspace*{0.1in}\WeylSection.3. Coxeter-group-related notions\dotfill 82\\%
\hspace*{0.1in}\WeylSection.4. Weights\dotfill 83\\%
\hspace*{0.1in}\WeylSection.5. The Networked-numbers Game\dotfill 83}
\end{center}
\end{abstract}

\newpage
\def\abstractname{Table of contents, continued}
\begin{abstract}
\begin{center}
\parbox{5in}{
\hspace*{0.1in}\WeylSection.6. Finite orbits\dotfill 84\\%
\hspace*{0.1in}\WeylSection.7. Finite Weyl groups and finite co-root systems\dotfill 85\\%
\hspace*{0.1in}\WeylSection.8. Parabolic subgroups of finite index\dotfill 86\\%
\hspace*{0.1in}\WeylSection.9. Weyl symmetric functions\dotfill 87\\%
\hspace*{0.1in}\WeylSection.10. A Weyl-group-invariant inner product\dotfill 87\\%
\hspace*{0.1in}\WeylSection.11. Terminating Networked-numbers Games\dotfill 89\\%
\hspace*{0.1in}\WeylSection.12. Kac--Moody Lie algebras\dotfill 90\\%
\hspace*{0.1in}\WeylSection.13. Linear representations of Kac--Moody Lie algebras\dotfill 93\\%
\hspace*{-0.1in}\LFKSection. {\em La Florado Klasado} (The Flowering Classification)\dotfill 95\\%
\hspace*{0.1in}\LFKSection.1. Some reflections on the nature of mathematical classifications\dotfill 95\\%
\hspace*{0.1in}\LFKSection.2. Our statement of {\em LFK}, with proofs\dotfill 96\\%
\hspace*{0.1in}\LFKSection.3. Some further remarks on our version of {\em LFK}\dotfill 97\\%
\hspace*{-0.1in}\FlowerSection. Consequences of {\em LFK} for our algebraic-combinatorial environment\dotfill 98\\%
\hspace*{0.1in}\FlowerSection.1. A $\mathcal{W}$-invariant inner product space\dotfill 98\\%
\hspace*{0.1in}\FlowerSection.2. The root system and co-root system\dotfill 98\\%
\hspace*{0.1in}\FlowerSection.3. The longest Weyl group element\dotfill 98\\%
\hspace*{0.1in}\FlowerSection.4. Root lattice cosets as connected components of the weight lattice\dotfill 99\\%
\hspace*{0.1in}\FlowerSection.5. $\mathscr{G}$-structured posets and saturated sets of weights\dotfill 99\\%
\hspace*{0.1in}\FlowerSection.6. Monomial symmetric functions and Weyl bialternants\dotfill 100\\%
\hspace*{0.1in}\FlowerSection.7. Specializations\dotfill 101\\%
\hspace*{0.1in}\FlowerSection.8. Splitting posets\dotfill 102\\%
\hspace*{0.1in}\FlowerSection.9. Vertex coloring\dotfill 103\\%
\hspace*{0.1in}\FlowerSection.10. Kac--Moody Lie algebras and representations in finite dimensions\dotfill 104\\%
\hspace*{0.1in}\FlowerSection.11. DCML's as supporting graphs\dotfill 106\\%
\hspace*{0.5in} \\ 
\hspace*{-0.25in}\fbox{Part IV}\ \  Examples and applications\dotfill 114\\
\hspace*{-0.1in}\QuasiExampleSection. Examples of DCML's that are splitting posets: The quasi-minuscule case\dotfill 114\\%
\hspace*{0.1in}\QuasiExampleSection.1. Root lengths\dotfill 114\\%
\hspace*{0.1in}\QuasiExampleSection.2. Highest short root\dotfill 115\\%
\hspace*{0.1in}\QuasiExampleSection.3. DCML's associated with quasi-minuscule dominant weights\dotfill 116\\%
\hspace*{0.1in}\QuasiExampleSection.4. Comments and questions about quasi-minuscule DCML's\dotfill 118\\%
\hspace*{0.1in}\QuasiExampleSection.5. Quasi-minuscule DCML's are splitting posets and supporting graphs\dotfill 119\\%
\hspace*{-0.1in}\MinusculeExampleSection. Examples of DCDL's that are splitting posets: The minuscule case\dotfill 122\\%
\hspace*{0.1in}\MinusculeExampleSection.1. Some characterizations of minuscule dominant weights\dotfill 123\\%
\hspace*{0.1in}\MinusculeExampleSection.2. Classification of minuscule fundamental weights\dotfill 124\\%
\hspace*{0.1in}\MinusculeExampleSection.3. Some poset-structural aspects of the saturated set of weights $\Pi(\omega_{f})$\dotfill 125\\%
\hspace*{0.1in}\MinusculeExampleSection.4. Some properties and distinctions of minuscule Weyl bialternants\dotfill 126\\%
\hspace*{0.1in}\MinusculeExampleSection.5. Meets and joins in $\Pi(\omega_{f})$\dotfill 127\\%
\hspace*{0.1in}\MinusculeExampleSection.6. The saturated set of weights $\Pi(\omega_{f})$ as a DCDL\dotfill 133\\%
\hspace*{0.1in}\MinusculeExampleSection.7. Properties of certain subposets of minuscule compression posets\dotfill 134\\%
\hspace*{0.1in}\MinusculeExampleSection.8. Distributivity of minuscule compression posets\dotfill 138\\%
\hspace*{0.1in}\MinusculeExampleSection.9. Vascular graph doppelg\"{a}ngers in minuscule compression posets\dotfill 143\\%
\hspace*{-0.1in}\MinusculeLatticePosetSection. Skew-stacks of minuscule compression posets\dotfill 144\\%
\hspace*{0.5in} \\ 
\hspace*{-0.25in}References\dotfill 147}
\end{center}
\end{abstract}

\newpage
\noindent
{\bf \S 1. Introduction.} 
Modular lattices and distributive lattices are elegant ordered structures that, typically, are not explicitly featured in the general education of many mathematicians. 
Even so, these structures are commonly encountered in the undergraduate mathematics curriculum, although they are rarely named as such. 
The lattice of normal subgroups of a group is modular, as are the lattice of ideals of a ring and the lattice of subspaces of a vector space. 
The divisibility lattice of a positive integer is distributive, and the set of all $2^{n}$ subsets of an $n$-element set is a distributive lattice when the subsets are ordered by inclusion. 
The latter is an example of how some of the most basic data structures arising in elementary enumeration are naturally associated with distributive lattices. 
For another example of this phenomenon, regard each $k$-element subset of the set $\{1,2,\ldots,n+1\}$ as an increasing $k$-tuple, and then order these {\footnotesize $\left(\!\!\begin{array}{c}n+1\\ k\end{array}\!\!\right)$} $k$-tuples by `reverse component-wise comparison', i.e.\ $S \leq T$ for $S=(S_{1},\ldots,S_{k})$ and $T=(T_{1},\ldots,T_{k})$ if and only if $S_{i} \geq T_{i}$ for all $i \in \{1,2,\ldots,k\}$. 
The resulting partially ordered set of these $k$-tuples is a distributive lattice -- see \ASevenFigure\ for a depiction of this lattice when $n=7$ and $k=2$. 

For us, algebraic combinatorics is the venue for some of the most striking and natural occurrences of modular and distributive lattices. 
When such structures arise in algebraic combinatorics, they are generally equipped with additional information, such as colors supplied to the edges of graphical representations of the lattices (much like Cayley graphs in group theory). 
Indeed, the distributive lattice from the close of the preceding paragraph is one such example, and we will revisit this example later within one of its several algebraic manifestations. 
In this manuscript, our main interest is in finite modular $\!\!$/$\!\!$ distributive lattices whose covering digraphs (aka Hasse diagrams) have edges colored in such a way that on any `diamond' of edges\footnote{Within the covering digraph, a `diamond' is a four-element subgraph that corresponds to a length two distributive sublattice.}, opposite edges have the same color. 
These are the diamond-colored modular $\!\!$/$\!\!$ distributive lattices advertised in our title. 

Some of the purposes of earlier versions of this manuscript were to establish within the literature some foundational results concerning diamond-colored modular and distributive lattices, to provide a browsable tutorial on the rudiments of this topic, to present some apparently new results as well as new interpretations of classical results\footnote{Within the present manuscript, these include, in our view, some of the edge-colored aspects of \ColorsLemma, \EncapsulateProposition, and \FundamentalTheorem; \NewSublatticeResults; \NewLFKTheorem; most all of \S \FrameSection; \ClassTheorems; \GStructureDCDLTheorem; \OrbitProposition\ and aspects of \OrbitCorollaries; \QMSplitting; \NewMinusculeResults\ as well as our completely general proofs of \MinusculeDistributiveResults; and \GStructureDCDLCorollary.}, and to serve as a convenient reference for readers who are mainly interested in the contexts in which these lattice structures arise. 
For us, these contexts are algebraic-combinatorial and include specifically the theory of semisimple Lie algebra representations and their companion Weyl group symmetric functions. 
After a major revision, this manuscript now provides detailed overviews of these particular contexts. 

A central goal in our programmatic combinatorial study of these latter contexts is to find interesting poset models for families of such representations $\!\!$/$\!\!$ symmetric functions. 
For a tabular summary of case-wise results, see Table 1.1 of \cite{DonPosetModels}; those results provide some evidence that diamond-colored modular and distributive lattices arise naturally in the study of semisimple Lie algebra representations and Weyl group symmetric functions. 
Most of the diamond-colored lattices featured in that table are distributive, but we know that not all such representations $\!\!$/$\!\!$ symmetric functions can be modelled using distributive lattices. 
Based on our investigations so far, it seems that modular lattices will suffice for that eventual purpose.  
Thus, many results here are stated in terms of modular lattices.

A detailed outline of the manuscript's organization is provided in the table of contents. 
Here, we provide some additional qualitative comments.   
Part I of this manuscript develops, from a purely combinatorial point of view, the order-theoretic ideas we intend to focus on. 
Section \RoutineSection\ sets the overall environment. 
This includes definitions and notation, some of which are idiosyncratic, as well as rudimentary results. 
Perhaps the central result of this section is \FundamentalTheorem, which we call the Fundamental Theorem for Finite Diamond-colored Distributive Lattices. 
This is an edge-colored version of Birkhoff's famous representation theorem for distributive lattices. 
In \S \SubstructureSection, we develop some sublattice results -- some standard, some not so standard -- that will be useful in later sections. 

Part II extends the ideas of Part I by making connections with the classification of complex finite-dimensional simple Lie algebras and of finite irreducible Weyl groups with associated root systems. 
These latter objects are historically classified, irredundantly, by type -- $\myA_{n}$ ($n \geq 1$), $\myB_{n}$ ($n \geq 3$), $\myC_{n}$ ($n \geq 2$), $\myD_{n}$ ($n \geq 4$), $\myE_{6}$, $\myE_{7}$, $\myE_{8}$, $\myF_{4}$, $\myG_{2}$. 
That said, these algebraic structures are only referenced implicitly in Part II via the Dynkin diagrams of \IEGGraphFigure; otherwise, the development in Part II is still purely combinatorial. 
More specifically, in \S \GStructureIntroSection, we articulate a structure property of ranked posets that is naturally associated with Dynkin-type diagrams, which are here called `integral embryophytic graphs' or IEG's. 
\NewLFKTheorem, which concerns diamond-colored modular lattices that possess this structure property, provides a new instance of phenomena that are classified by the IEG's of \IEGGraphFigure, which we prefer to call `Coxeter--Dynkin flowers'. 
In \S \FrameSection, we introduce the concept of a `scaffold' which provides an alternative to \FundamentalTheorem\ for representing diamond-colored distributive lattices (see \FramingProps). 
We also show how, in some circumstances, a scaffold can be understood as a kind of stacking of smaller posets (\IsomLatticeTheorem). 
In \S \GridPosetSection, we recount some of this author's contributions to the theory of two-color grid posets from \cite{ADLMPPW} and \cite{ADLP}, and we develop a purely combinatorial characterization those two-color grid posets used in \cite{ADLMPPW} to construct diamond-colored distributive lattice models of symmetric functions associated with irreducible rank two Weyl groups (i.e.\ of type $\myA_{2}$, $\myC_{2}$, $\myG_{2}$). 
Among other special properties, these diamond-colored distributive lattices are `$\myX_{2}$-structured' for the appropriate $\myX \in \{\myA, \myC, \myG\}$. 
In \GStructureDCDLTheorem\ of \S \TwoColorSturdySection\ (and again in \GStructureDCDLCorollary\ of \S \MinusculeLatticePosetSection), we see how our structure property can, in some circumstances, be demonstrated for those diamond-colored distributive lattices obtained via certain scaffolds or by the corresponding stacking of posets. 

In Part III, our focus shifts to more algebraic combinatorial ideas. 
In \S \WeylSection, we consider a series of finiteness problems related to IEG's, many of which are algebraic in nature. 
These include well-known problems, such as classifying finite irreducible Weyl groups and finite-dimensional simple Lie algebras over $\mathbb{C}$, and lesser-known or new problems, such as identifying circumstances under which a Weyl group has a finite-index parabolic subgroup or a non-constant symmetric function. 
We demonstrate how each of the finiteness problems posed is answered by the collection of Coxeter--Dynkin flowers exhibited in \IEGGraphFigure. 
In \S \LFKSection, we gather the finiteness results of \S \WeylSection\ into one theorem -- \LFK\ -- which we call {\em La Florado Klasado} (in Esperanto), i.e.\ `The Flowering Classification'. 
We hereby propose `{\em La Florado Klasado}', or `{\em LFK}', as a name for the vast menagerie of classification results that are equivalent to some aspect of the version of {\em LFK} presented here. 
In \S \FlowerSection, we consider further implications of the finiteness hypothesis of {\em LFK} for the study of Weyl groups and simple Lie algebra representations, and in later parts of this section, we reconnect our discourse with diamond-colored modular and distributive lattices. 
With its detours into the theory of Kac--Moody Lie algebras and the Weyl--Kac character formula, this part of the manuscript will be harder going for the novice. 
For such a reader, a good goal for this part might be to achieve a familiarity with the equivalences stated in our version of {\em LFK} (\LFK). 

In Part IV, we develop some applications and consider some special cases. 
The main goal of \S \QuasiExampleSection\ is to develop for so-called quasi-minuscule representations the same kind of algebraic-combinatorial results obtained in \cite{DonAdjoint} for the adjoint representations. 
In particular, we establish in \QMSplitting\ that for a finite-dimensional simple Lie algebra over $\mathbb{C}$ whose root system has exactly $m$ short simple roots, there are exactly $m$ diamond-colored modular lattice `supporting graphs' for its quasi-minuscule representation, and these are also the only diamond-colored modular lattice `splitting posets' for the associated Weyl symmetric function. 
For any so-called minuscule representation, it is known that there is only one supporting graph and that this supporting graph coincides with the only splitting poset. 
We re-establish this result in \MinusculeWeylBialt\ of \S \MinusculeExampleSection. 
The main purpose of \S \MinusculeExampleSection\ is to offer a completely general proof that any such minuscule supporting graph $\!\!$/$\!\!$ splitting poset is a diamond-colored distributive lattice and that the companion poset of join irreducibles (which we call a `compression poset') is itself a distributive lattice that can be realized as a sublattice of a product of two nonempty chains. 
In the closing section (\S \MinusculeLatticePosetSection), we appeal to many of the results of \S \MinusculeExampleSection\ in order to demonstrate that the diamond-colored distributive lattices obtained from certain stackings of minuscule compression posets enjoy some of the structural properties of \S \GStructureIntroSection\ and \S \FrameSection. 

Here are some general notions that readers should keep in mind. 
The combinatorial objects of primary interest for us are partially ordered sets thought of as graphs when identified with their covering digraphs. 
Generally speaking, any graph we work with will be finite and directed, with no loops and at most one edge between any two vertices, i.e.\ simple directed graphs.  
Most often, such a graph will either have its edges or its vertices colored by elements from some index set (usually a set of positive integers).  
Such coloring can provide crucial information when we view these structures within the algebraic contexts that primarily motivate our interest.  
The conventions and notation we use here largely borrow from \cite{DonSupp}, \cite{ADLP}, \cite{ADLMPPW}, and \cite{StanleyText}. 
Many of the poset-related concepts are depicted in figures throughout the manuscript.  
Throughout, we alert the reader to what we believe are new results or new proofs of known results, and we attempt to find references for results that are non-routine or not very well known.  

Before we embark, we offer two more general thoughts. 
We regard many of our ideas here to be an outgrowth of R.\ A.\ Proctor's pioneering work in \cite{PrEur} with distributive lattices and compression posets related to minuscule representations of simple Lie algebras. 
Also, to the extent that lattice theory still exists on the periphery of mathematics as a tidy but dismissible area -- see \cite{Rota} for some historical reflections on such attitudes -- we hope our work here might be experienced by readers as an affirmative demonstration of its inherence, its usefulness, and its connectedness to some central areas of mathematics.

\newcommand{\ASevenFlower}{
\setlength{\unitlength}{0.5cm}
\begin{picture}(20,1)
\put(-2.5,-0.25){\LARGE $\myA_{7}$}
\put(-0.9,0){\vector(1,0){1.2}}
{\color{Cyan}\put(1,0){\circle*{0.7}}}
\put(1,0){\circle{0.7}}
{\color{Red}\put(4,0){\circle*{0.7}}}
\put(4,0){\circle{0.7}}
{\color{Purple}\put(7,0){\circle*{0.7}}}
\put(7,0){\circle{0.7}}
{\color{OliveGreen}\put(10,0){\circle*{0.7}}}
\put(10,0){\circle{0.7}}
{\color{BurntOrange}\put(13,0){\circle*{0.7}}}
\put(13,0){\circle{0.7}}
{\color{SpringGreen}\put(16,0){\circle*{0.7}}}
\put(16,0){\circle{0.7}}
{\color{Brown}\put(19,0){\circle*{0.7}}}
\put(19,0){\circle{0.7}}
\thicklines
\put(1.35,0){\line(1,0){2.3}}
\put(4.35,0){\line(1,0){2.3}}
\put(7.35,0){\line(1,0){2.3}}
\put(10.35,0){\line(1,0){2.3}}
\put(13.35,0){\line(1,0){2.3}}
\put(16.35,0){\line(1,0){2.3}}
\put(0.8,0.6){\textcolor{Cyan}{\em 1}}
\put(3.8,0.6){\textcolor{Red}{\em 2}}
\put(6.8,0.6){\textcolor{Purple}{\em 3}}
\put(9.8,0.6){\textcolor{OliveGreen}{\em 4}}
\put(12.8,0.6){\textcolor{BurntOrange}{\em 5}}
\put(15.8,0.6){\textcolor{SpringGreen}{\em 6}}
\put(18.8,0.6){\textcolor{Brown}{\em 7}}
\end{picture}
}

\begin{figure}
\begin{center}
\parbox{6in}{\small {\bf \ASevenFigure}\ \ A diamond-colored distributive lattice assembled from the $\left(\begin{array}{c}8\\ 2\end{array}\right) = 28$ two-element subsets of the eight-element set $\{\color{Cyan}{1},\color{Red}{2},\color{Purple}{3},\color{OliveGreen}{4},\color{BurntOrange}{5},\color{SpringGreen}{6},\color{Brown}{7},\color{Black}{8}\}$. 
Our two-element subsets are written as ordered pairs with the smallest element first. 
We say $(S_{1},S_{2}) \leq (T_{1},T_{2})$ if and only if $S_{1} \geq T_{1}$ and $S_{2} \geq T_{2}$. 
Moreover $(S_{1},S_{2}) \rightarrow (T_{1},T_{2})$ iff \underline{\sl either} $S_{1} = T_{1}+1$ with $S_{2}=T_{2}$ (in which case the covering digraph edge is given `color' $T_{1}$) \underline{\sl or} $S_{2} = T_{2}+1$ with $S_{1}=T_{1}$ (in which case the covering digraph edge is given `color' $T_{2}$).}

\vspace*{0.2in} 
\parbox{6in}{\footnotesize (In the language of \S \GStructureIntroSection, this DCDL is `$\myA_{7}$-structured', where $\myA_{7}$ is the graph from \IEGGraphFigure\ depicted immediately below as a path graph with seven colored nodes. 
In the language of \S \MinusculeExampleSection, this is the DCDL associated with the `minuscule fundamental weight $\omega_{\color{Red}2}$'. Thus our notation `$L_{\mytinyA_{7}}(\omega_{\color{Red}2})$'.)}

\vspace*{0.2in} 
\ASevenFlower

\vspace*{0.4in} 
\setlength{\unitlength}{1cm}
\begin{picture}(7,13)
\put(5,11.5){\LARGE $L_{\mysmallA_{7}}(\omega_{\mbox{\color{Red}{\small $2$}}})$}
\put(4.85,11.6){\vector(-2,-1){1.5}}
\put(0,0){\TypeEboxDot{Gray}}
\put(1,1){\TypeEboxDot{Gray}}
\put(0,2){\TypeEboxDot{Gray}}
\put(2,2){\TypeEboxDot{Gray}}
\put(1,3){\TypeEboxDot{Gray}}
\put(3,3){\TypeEboxDot{Gray}}
\put(0,4){\TypeEboxDot{Gray}}
\put(2,4){\TypeEboxDot{Gray}}
\put(4,4){\TypeEboxDot{Gray}}
\put(1,5){\TypeEboxDot{Gray}}
\put(3,5){\TypeEboxDot{Gray}}
\put(5,5){\TypeEboxDot{Gray}}
\put(0,6){\TypeEboxDot{Gray}}
\put(2,6){\TypeEboxDot{Gray}}
\put(4,6){\TypeEboxDot{Gray}}
\put(6,6){\TypeEboxDot{Gray}}
\put(1,7){\TypeEboxDot{Gray}}
\put(3,7){\TypeEboxDot{Gray}}
\put(5,7){\TypeEboxDot{Gray}}
\put(0,8){\TypeEboxDot{Gray}}
\put(2,8){\TypeEboxDot{Gray}}
\put(4,8){\TypeEboxDot{Gray}}
\put(1,9){\TypeEboxDot{Gray}}
\put(3,9){\TypeEboxDot{Gray}}
\put(0,10){\TypeEboxDot{Gray}}
\put(2,10){\TypeEboxDot{Gray}}
\put(1,11){\TypeEboxDot{Gray}}
\put(0,12){\TypeEboxDot{Gray}}
\thicklines
\put(0.625,0.375){\color{SpringGreen}\qbezier(0,0)(0.375,0.375)(0.75,0.75)}
\put(1.625,1.375){\color{BurntOrange}\qbezier(0,0)(0.375,0.375)(0.75,0.75)}
\put(2.625,2.375){\color{OliveGreen}\qbezier(0,0)(0.375,0.375)(0.75,0.75)}
\put(3.625,3.375){\color{Purple}\qbezier(0,0)(0.375,0.375)(0.75,0.75)}
\put(4.625,4.375){\color{Red}\qbezier(0,0)(0.375,0.375)(0.75,0.75)}
\put(5.625,5.375){\color{Cyan}\qbezier(0,0)(0.375,0.375)(0.75,0.75)}
\put(0.625,2.375){\color{BurntOrange}\qbezier(0,0)(0.375,0.375)(0.75,0.75)}
\put(1.625,3.375){\color{OliveGreen}\qbezier(0,0)(0.375,0.375)(0.75,0.75)}
\put(2.625,4.375){\color{Purple}\qbezier(0,0)(0.375,0.375)(0.75,0.75)}
\put(3.625,5.375){\color{Red}\qbezier(0,0)(0.375,0.375)(0.75,0.75)}
\put(4.625,6.375){\color{Cyan}\qbezier(0,0)(0.375,0.375)(0.75,0.75)}
\put(0.625,4.375){\color{OliveGreen}\qbezier(0,0)(0.375,0.375)(0.75,0.75)}
\put(1.625,5.375){\color{Purple}\qbezier(0,0)(0.375,0.375)(0.75,0.75)}
\put(2.625,6.375){\color{Red}\qbezier(0,0)(0.375,0.375)(0.75,0.75)}
\put(3.625,7.375){\color{Cyan}\qbezier(0,0)(0.375,0.375)(0.75,0.75)}
\put(0.625,6.375){\color{Purple}\qbezier(0,0)(0.375,0.375)(0.75,0.75)}
\put(1.625,7.375){\color{Red}\qbezier(0,0)(0.375,0.375)(0.75,0.75)}
\put(2.625,8.375){\color{Cyan}\qbezier(0,0)(0.375,0.375)(0.75,0.75)}
\put(0.625,8.375){\color{Red}\qbezier(0,0)(0.375,0.375)(0.75,0.75)}
\put(1.625,9.375){\color{Cyan}\qbezier(0,0)(0.375,0.375)(0.75,0.75)}
\put(0.625,10.375){\color{Cyan}\qbezier(0,0)(0.375,0.375)(0.75,0.75)}
\put(1.39,11.375){\color{Red}\qbezier(0,0)(-0.375,0.375)(-0.75,0.75)}
\put(2.39,10.375){\color{Purple}\qbezier(0,0)(-0.375,0.375)(-0.75,0.75)}
\put(3.39,9.375){\color{OliveGreen}\qbezier(0,0)(-0.375,0.375)(-0.75,0.75)}
\put(4.39,8.375){\color{BurntOrange}\qbezier(0,0)(-0.375,0.375)(-0.75,0.75)}
\put(5.39,7.375){\color{SpringGreen}\qbezier(0,0)(-0.375,0.375)(-0.75,0.75)}
\put(6.39,6.375){\color{Brown}\qbezier(0,0)(-0.375,0.375)(-0.75,0.75)}
\put(1.39,9.375){\color{Purple}\qbezier(0,0)(-0.375,0.375)(-0.75,0.75)}
\put(2.39,8.375){\color{OliveGreen}\qbezier(0,0)(-0.375,0.375)(-0.75,0.75)}
\put(3.39,7.375){\color{BurntOrange}\qbezier(0,0)(-0.375,0.375)(-0.75,0.75)}
\put(4.39,6.375){\color{SpringGreen}\qbezier(0,0)(-0.375,0.375)(-0.75,0.75)}
\put(5.39,5.375){\color{Brown}\qbezier(0,0)(-0.375,0.375)(-0.75,0.75)}
\put(1.39,7.375){\color{OliveGreen}\qbezier(0,0)(-0.375,0.375)(-0.75,0.75)}
\put(2.39,6.375){\color{BurntOrange}\qbezier(0,0)(-0.375,0.375)(-0.75,0.75)}
\put(3.39,5.375){\color{SpringGreen}\qbezier(0,0)(-0.375,0.375)(-0.75,0.75)}
\put(4.39,4.375){\color{Brown}\qbezier(0,0)(-0.375,0.375)(-0.75,0.75)}
\put(1.39,5.375){\color{BurntOrange}\qbezier(0,0)(-0.375,0.375)(-0.75,0.75)}
\put(2.39,4.375){\color{SpringGreen}\qbezier(0,0)(-0.375,0.375)(-0.75,0.75)}
\put(3.39,3.375){\color{Brown}\qbezier(0,0)(-0.375,0.375)(-0.75,0.75)}
\put(1.39,3.375){\color{SpringGreen}\qbezier(0,0)(-0.375,0.375)(-0.75,0.75)}
\put(2.39,2.375){\color{Brown}\qbezier(0,0)(-0.375,0.375)(-0.75,0.75)}
\put(1.39,1.375){\color{Brown}\qbezier(0,0)(-0.375,0.375)(-0.75,0.75)}
\put(0.9,0.65){\color{SpringGreen}{\em 6}}
\put(0.9,2.65){\color{BurntOrange}{\em 5}}
\put(1.9,1.65){\color{BurntOrange}{\em 5}}
\put(0.9,4.65){\color{OliveGreen}{\em 4}}
\put(1.9,3.65){\color{OliveGreen}{\em 4}}
\put(2.9,2.65){\color{OliveGreen}{\em 4}}
\put(0.9,6.65){\color{Purple}{\em 3}}
\put(1.9,5.65){\color{Purple}{\em 3}}
\put(2.9,4.65){\color{Purple}{\em 3}}
\put(3.9,3.65){\color{Purple}{\em 3}}
\put(0.9,8.65){\color{Red}{\em 2}}
\put(1.9,7.65){\color{Red}{\em 2}}
\put(2.9,6.65){\color{Red}{\em 2}}
\put(3.9,5.65){\color{Red}{\em 2}}
\put(4.9,4.65){\color{Red}{\em 2}}
\put(0.9,10.65){\color{Cyan}{\em 1}}
\put(1.9,9.65){\color{Cyan}{\em 1}}
\put(2.9,8.65){\color{Cyan}{\em 1}}
\put(3.9,7.65){\color{Cyan}{\em 1}}
\put(4.9,6.65){\color{Cyan}{\em 1}}
\put(5.9,5.65){\color{Cyan}{\em 1}}
\put(0.9,11.65){\color{Red}{\em 2}}
\put(0.9,9.65){\color{Purple}{\em 3}}
\put(1.9,10.65){\color{Purple}{\em 3}}
\put(0.9,7.65){\color{OliveGreen}{\em 4}}
\put(1.9,8.65){\color{OliveGreen}{\em 4}}
\put(2.9,9.65){\color{OliveGreen}{\em 4}}
\put(0.9,5.65){\color{BurntOrange}{\em 5}}
\put(1.9,6.65){\color{BurntOrange}{\em 5}}
\put(2.9,7.65){\color{BurntOrange}{\em 5}}
\put(3.9,8.65){\color{BurntOrange}{\em 5}}
\put(0.9,3.65){\color{SpringGreen}{\em 6}}
\put(1.9,4.65){\color{SpringGreen}{\em 6}}
\put(2.9,5.65){\color{SpringGreen}{\em 6}}
\put(3.9,6.65){\color{SpringGreen}{\em 6}}
\put(4.9,7.65){\color{SpringGreen}{\em 6}}
\put(0.9,1.65){\color{Brown}{\em 7}}
\put(1.9,2.65){\color{Brown}{\em 7}}
\put(2.9,3.65){\color{Brown}{\em 7}}
\put(3.9,4.65){\color{Brown}{\em 7}}
\put(4.9,5.65){\color{Brown}{\em 7}}
\put(5.9,6.65){\color{Brown}{\em 7}}
%
\put(-0.6,12.175){\small $( {\color{Cyan}1},{\color{Red}2} )$}
\put(0.4,11.175){\small $( {\color{Cyan}1},{\color{Purple}3} )$}
\put(-0.6,10.175){\small $( {\color{Red}2},{\color{Purple}3} )$}
\put(2.8,10.175){\small $( {\color{Cyan}1},{\color{OliveGreen}4} )$}
\put(0.4,9.175){\small $( {\color{Red}2},{\color{OliveGreen}4} )$}
\put(3.8,9.175){\small $( {\color{Cyan}1},{\color{BurntOrange}5} )$}
\put(-0.6,8.175){\small $( {\color{Purple}3},{\color{OliveGreen}4} )$}
\put(2.8,8.175){\small $( {\color{Red}2},{\color{BurntOrange}5} )$}
\put(4.8,8.175){\small $( {\color{Cyan}1},{\color{SpringGreen}6} )$}
\put(0.4,7.175){\small $( {\color{Purple}3},{\color{BurntOrange}5} )$}
\put(3.8,7.175){\small $( {\color{Red}2},{\color{SpringGreen}6} )$}
\put(5.8,7.175){\small $( {\color{Cyan}1},{\color{Brown}7} )$}
\put(-0.6,6.175){\small $( {\color{OliveGreen}4},{\color{BurntOrange}5} )$}
\put(2.8,6.175){\small $( {\color{Purple}3},{\color{SpringGreen}6} )$}
\put(4.8,6.175){\small $( {\color{Red}2},{\color{Brown}7} )$}
\put(6.8,6.175){\small $( {\color{Cyan}1},{\color{Black}8} )$}
\put(0.4,5.175){\small $( {\color{OliveGreen}4},{\color{SpringGreen}6} )$}
\put(3.8,5.175){\small $( {\color{Purple}3},{\color{Brown}7} )$}
\put(5.8,5.175){\small $( {\color{Red}2},{\color{Black}8} )$}
\put(-0.6,4.175){\small $( {\color{BurntOrange}5},{\color{SpringGreen}6} )$}
\put(2.8,4.175){\small $( {\color{OliveGreen}4},{\color{Brown}7} )$}
\put(4.8,4.175){\small $( {\color{Purple}3},{\color{Black}8} )$}
\put(0.4,3.175){\small $( {\color{BurntOrange}5},{\color{Brown}7} )$}
\put(3.8,3.175){\small $( {\color{OliveGreen}4},{\color{Black}8} )$}
\put(-0.6,2.175){\small $( {\color{SpringGreen}6},{\color{Brown}7} )$}
\put(2.8,2.175){\small $( {\color{BurntOrange}5},{\color{Black}8} )$}
\put(0.4,1.175){\small $( {\color{SpringGreen}6},{\color{Black}8} )$}
\put(-0.6,0.175){\small $( {\color{Brown}7},{\color{Black}8} )$}
\end{picture}
\end{center}
\end{figure}

\begin{figure}
\begin{center}
{\bf \ASevenFigureTwo}\ \ The diamond-colored distributive lattice $L_{\mytinyA_{7}}(\omega_{\color{Red}2})$ from \ASevenFigure\\ and its companion compression poset $P_{\mytinyA_{7}}(\omega_{\color{Red}2})$.

\vspace*{0.2in} 
\ASevenFlower

\vspace*{0.4in} 
\setlength{\unitlength}{1cm}
\begin{picture}(7,13)
\put(5,11.5){\LARGE $L_{\mysmallA_{7}}(\omega_{\mbox{\color{Red}{\small $2$}}})$}
\put(4.85,11.6){\vector(-2,-1){1.5}}
\put(0,0){\TypeEboxDot{Gray}}
\put(1,1){\TypeEboxDot{Gray}}
\put(0,2){\TypeEboxDot{Gray}}
\put(2,2){\TypeEboxDot{Gray}}
\put(1,3){\TypeEboxDot{Gray}}
\put(3,3){\TypeEboxDot{Gray}}
\put(0,4){\TypeEboxDot{Gray}}
\put(2,4){\TypeEboxDot{Gray}}
\put(4,4){\TypeEboxDot{Gray}}
\put(1,5){\TypeEboxDot{Gray}}
\put(3,5){\TypeEboxDot{Gray}}
\put(5,5){\TypeEboxDot{Gray}}
\put(0,6){\TypeEboxDot{Gray}}
\put(2,6){\TypeEboxDot{Gray}}
\put(4,6){\TypeEboxDot{Gray}}
\put(6,6){\TypeEboxDot{Gray}}
\put(1,7){\TypeEboxDot{Gray}}
\put(3,7){\TypeEboxDot{Gray}}
\put(5,7){\TypeEboxDot{Gray}}
\put(0,8){\TypeEboxDot{Gray}}
\put(2,8){\TypeEboxDot{Gray}}
\put(4,8){\TypeEboxDot{Gray}}
\put(1,9){\TypeEboxDot{Gray}}
\put(3,9){\TypeEboxDot{Gray}}
\put(0,10){\TypeEboxDot{Gray}}
\put(2,10){\TypeEboxDot{Gray}}
\put(1,11){\TypeEboxDot{Gray}}
\put(0,12){\TypeEboxDot{Gray}}
\thicklines
\put(0.625,0.375){\color{SpringGreen}\qbezier(0,0)(0.375,0.375)(0.75,0.75)}
\put(1.625,1.375){\color{BurntOrange}\qbezier(0,0)(0.375,0.375)(0.75,0.75)}
\put(2.625,2.375){\color{OliveGreen}\qbezier(0,0)(0.375,0.375)(0.75,0.75)}
\put(3.625,3.375){\color{Purple}\qbezier(0,0)(0.375,0.375)(0.75,0.75)}
\put(4.625,4.375){\color{Red}\qbezier(0,0)(0.375,0.375)(0.75,0.75)}
\put(5.625,5.375){\color{Cyan}\qbezier(0,0)(0.375,0.375)(0.75,0.75)}
\put(0.625,2.375){\color{BurntOrange}\qbezier(0,0)(0.375,0.375)(0.75,0.75)}
\put(1.625,3.375){\color{OliveGreen}\qbezier(0,0)(0.375,0.375)(0.75,0.75)}
\put(2.625,4.375){\color{Purple}\qbezier(0,0)(0.375,0.375)(0.75,0.75)}
\put(3.625,5.375){\color{Red}\qbezier(0,0)(0.375,0.375)(0.75,0.75)}
\put(4.625,6.375){\color{Cyan}\qbezier(0,0)(0.375,0.375)(0.75,0.75)}
\put(0.625,4.375){\color{OliveGreen}\qbezier(0,0)(0.375,0.375)(0.75,0.75)}
\put(1.625,5.375){\color{Purple}\qbezier(0,0)(0.375,0.375)(0.75,0.75)}
\put(2.625,6.375){\color{Red}\qbezier(0,0)(0.375,0.375)(0.75,0.75)}
\put(3.625,7.375){\color{Cyan}\qbezier(0,0)(0.375,0.375)(0.75,0.75)}
\put(0.625,6.375){\color{Purple}\qbezier(0,0)(0.375,0.375)(0.75,0.75)}
\put(1.625,7.375){\color{Red}\qbezier(0,0)(0.375,0.375)(0.75,0.75)}
\put(2.625,8.375){\color{Cyan}\qbezier(0,0)(0.375,0.375)(0.75,0.75)}
\put(0.625,8.375){\color{Red}\qbezier(0,0)(0.375,0.375)(0.75,0.75)}
\put(1.625,9.375){\color{Cyan}\qbezier(0,0)(0.375,0.375)(0.75,0.75)}
\put(0.625,10.375){\color{Cyan}\qbezier(0,0)(0.375,0.375)(0.75,0.75)}
\put(1.39,11.375){\color{Red}\qbezier(0,0)(-0.375,0.375)(-0.75,0.75)}
\put(2.39,10.375){\color{Purple}\qbezier(0,0)(-0.375,0.375)(-0.75,0.75)}
\put(3.39,9.375){\color{OliveGreen}\qbezier(0,0)(-0.375,0.375)(-0.75,0.75)}
\put(4.39,8.375){\color{BurntOrange}\qbezier(0,0)(-0.375,0.375)(-0.75,0.75)}
\put(5.39,7.375){\color{SpringGreen}\qbezier(0,0)(-0.375,0.375)(-0.75,0.75)}
\put(6.39,6.375){\color{Brown}\qbezier(0,0)(-0.375,0.375)(-0.75,0.75)}
\put(1.39,9.375){\color{Purple}\qbezier(0,0)(-0.375,0.375)(-0.75,0.75)}
\put(2.39,8.375){\color{OliveGreen}\qbezier(0,0)(-0.375,0.375)(-0.75,0.75)}
\put(3.39,7.375){\color{BurntOrange}\qbezier(0,0)(-0.375,0.375)(-0.75,0.75)}
\put(4.39,6.375){\color{SpringGreen}\qbezier(0,0)(-0.375,0.375)(-0.75,0.75)}
\put(5.39,5.375){\color{Brown}\qbezier(0,0)(-0.375,0.375)(-0.75,0.75)}
\put(1.39,7.375){\color{OliveGreen}\qbezier(0,0)(-0.375,0.375)(-0.75,0.75)}
\put(2.39,6.375){\color{BurntOrange}\qbezier(0,0)(-0.375,0.375)(-0.75,0.75)}
\put(3.39,5.375){\color{SpringGreen}\qbezier(0,0)(-0.375,0.375)(-0.75,0.75)}
\put(4.39,4.375){\color{Brown}\qbezier(0,0)(-0.375,0.375)(-0.75,0.75)}
\put(1.39,5.375){\color{BurntOrange}\qbezier(0,0)(-0.375,0.375)(-0.75,0.75)}
\put(2.39,4.375){\color{SpringGreen}\qbezier(0,0)(-0.375,0.375)(-0.75,0.75)}
\put(3.39,3.375){\color{Brown}\qbezier(0,0)(-0.375,0.375)(-0.75,0.75)}
\put(1.39,3.375){\color{SpringGreen}\qbezier(0,0)(-0.375,0.375)(-0.75,0.75)}
\put(2.39,2.375){\color{Brown}\qbezier(0,0)(-0.375,0.375)(-0.75,0.75)}
\put(1.39,1.375){\color{Brown}\qbezier(0,0)(-0.375,0.375)(-0.75,0.75)}
\put(0.9,0.65){\color{SpringGreen}{\em 6}}
\put(0.9,2.65){\color{BurntOrange}{\em 5}}
\put(1.9,1.65){\color{BurntOrange}{\em 5}}
\put(0.9,4.65){\color{OliveGreen}{\em 4}}
\put(1.9,3.65){\color{OliveGreen}{\em 4}}
\put(2.9,2.65){\color{OliveGreen}{\em 4}}
\put(0.9,6.65){\color{Purple}{\em 3}}
\put(1.9,5.65){\color{Purple}{\em 3}}
\put(2.9,4.65){\color{Purple}{\em 3}}
\put(3.9,3.65){\color{Purple}{\em 3}}
\put(0.9,8.65){\color{Red}{\em 2}}
\put(1.9,7.65){\color{Red}{\em 2}}
\put(2.9,6.65){\color{Red}{\em 2}}
\put(3.9,5.65){\color{Red}{\em 2}}
\put(4.9,4.65){\color{Red}{\em 2}}
\put(0.9,10.65){\color{Cyan}{\em 1}}
\put(1.9,9.65){\color{Cyan}{\em 1}}
\put(2.9,8.65){\color{Cyan}{\em 1}}
\put(3.9,7.65){\color{Cyan}{\em 1}}
\put(4.9,6.65){\color{Cyan}{\em 1}}
\put(5.9,5.65){\color{Cyan}{\em 1}}
\put(0.9,11.65){\color{Red}{\em 2}}
\put(0.9,9.65){\color{Purple}{\em 3}}
\put(1.9,10.65){\color{Purple}{\em 3}}
\put(0.9,7.65){\color{OliveGreen}{\em 4}}
\put(1.9,8.65){\color{OliveGreen}{\em 4}}
\put(2.9,9.65){\color{OliveGreen}{\em 4}}
\put(0.9,5.65){\color{BurntOrange}{\em 5}}
\put(1.9,6.65){\color{BurntOrange}{\em 5}}
\put(2.9,7.65){\color{BurntOrange}{\em 5}}
\put(3.9,8.65){\color{BurntOrange}{\em 5}}
\put(0.9,3.65){\color{SpringGreen}{\em 6}}
\put(1.9,4.65){\color{SpringGreen}{\em 6}}
\put(2.9,5.65){\color{SpringGreen}{\em 6}}
\put(3.9,6.65){\color{SpringGreen}{\em 6}}
\put(4.9,7.65){\color{SpringGreen}{\em 6}}
\put(0.9,1.65){\color{Brown}{\em 7}}
\put(1.9,2.65){\color{Brown}{\em 7}}
\put(2.9,3.65){\color{Brown}{\em 7}}
\put(3.9,4.65){\color{Brown}{\em 7}}
\put(4.9,5.65){\color{Brown}{\em 7}}
\put(5.9,6.65){\color{Brown}{\em 7}}
\end{picture}
\begin{picture}(7,13)
\put(5,11.5){\LARGE $P_{\mysmallA_{7}}(\omega_{\mbox{\color{Red}{\small $2$}}})$}
\put(4.85,11.6){\vector(-2,-1){1.5}}
\put(0,10){\TypeEboxDot{Cyan}}
\put(0.1,10){\color{Cyan}{\em \footnotesize 1}}
\put(1,11){\TypeEboxDot{Red}}
\put(1,9){\TypeEboxDot{Red}}
\put(1.7,11.3){\color{Red}{\em \footnotesize 2}}
\put(1.1,9){\color{Red}{\em \footnotesize 2}}
\put(2,10){\TypeEboxDot{Purple}}
\put(2,8){\TypeEboxDot{Purple}}
\put(2.7,10.3){\color{Purple}{\em \footnotesize 3}}
\put(2.1,8){\color{Purple}{\em \footnotesize 3}}
\put(3,9){\TypeEboxDot{OliveGreen}}
\put(3,7){\TypeEboxDot{OliveGreen}}
\put(3.7,9.3){\color{OliveGreen}{\em \footnotesize 4}}
\put(3.1,7){\color{OliveGreen}{\em \footnotesize 4}}
\put(4,8){\TypeEboxDot{BurntOrange}}
\put(4,6){\TypeEboxDot{BurntOrange}}
\put(4.7,8.3){\color{BurntOrange}{\em \footnotesize 5}}
\put(4.1,6){\color{BurntOrange}{\em \footnotesize 5}}
\put(5,7){\TypeEboxDot{SpringGreen}}
\put(5,5){\TypeEboxDot{SpringGreen}}
\put(5.7,7.3){\color{SpringGreen}{\em \footnotesize 6}}
\put(5.1,5){\color{SpringGreen}{\em \footnotesize 6}}
\put(6,6){\TypeEboxDot{Brown}}
\put(6.7,6.3){\color{Brown}{\em \footnotesize 7}}
\thicklines
\put(5.625,5.375){\color{Black}\qbezier(0,0)(0.375,0.375)(0.75,0.75)}
\put(4.625,6.375){\color{Black}\qbezier(0,0)(0.375,0.375)(0.75,0.75)}
\put(3.625,7.375){\color{Black}\qbezier(0,0)(0.375,0.375)(0.75,0.75)}
\put(2.625,8.375){\color{Black}\qbezier(0,0)(0.375,0.375)(0.75,0.75)}
\put(1.625,9.375){\color{Black}\qbezier(0,0)(0.375,0.375)(0.75,0.75)}
\put(0.625,10.375){\color{Black}\qbezier(0,0)(0.375,0.375)(0.75,0.75)}
\put(2.39,10.375){\color{Black}\qbezier(0,0)(-0.375,0.375)(-0.75,0.75)}
\put(3.39,9.375){\color{Black}\qbezier(0,0)(-0.375,0.375)(-0.75,0.75)}
\put(4.39,8.375){\color{Black}\qbezier(0,0)(-0.375,0.375)(-0.75,0.75)}
\put(5.39,7.375){\color{Black}\qbezier(0,0)(-0.375,0.375)(-0.75,0.75)}
\put(6.39,6.375){\color{Black}\qbezier(0,0)(-0.375,0.375)(-0.75,0.75)}
\put(1.39,9.375){\color{Black}\qbezier(0,0)(-0.375,0.375)(-0.75,0.75)}
\put(2.39,8.375){\color{Black}\qbezier(0,0)(-0.375,0.375)(-0.75,0.75)}
\put(3.39,7.375){\color{Black}\qbezier(0,0)(-0.375,0.375)(-0.75,0.75)}
\put(4.39,6.375){\color{Black}\qbezier(0,0)(-0.375,0.375)(-0.75,0.75)}
\put(5.39,5.375){\color{Black}\qbezier(0,0)(-0.375,0.375)(-0.75,0.75)}
\end{picture}
\end{center}
\end{figure}

\newpage
\begin{center}
\fbox{\Large \bf Part I}\\ 
\underline{\large \bf A gentle introduction to order-theoretic aspects of diamond-coloring}
\end{center}
Many elementary aspects of the structure of finite modular and distributive lattices generalize intuitively to our diamond-colored setting, as with the Fundamental Theorem of Finite Diamond-colored Distributive Lattices (\FundamentalTheorem) and its corollaries. 
However, some structural aspects of diamond-colored modular lattices (DCML's) and diamond-colored distributive lattices (DCDL's) have no `monochromatic' counterpart, as with \ColorsLemma, \JCompResult, and \JCompTheorem. 
To be sure we account for these distinctions, we mostly proceed from first principles. 
In addition to stating many rudimentary diamond-coloring results for the record, the next several sections are also intended partly to serve as a readable and -- for the beginner, perhaps -- informative tutorial.   
The proofs that are presented in this part of the manuscript are new or have some noteworthy features, such as comparing/contrasting multi-colored aspects with standard monochromatic approaches or perhaps highlighting idiosyncratic perspectives; otherwise references are given.

\vspace*{0.5cm} 
\noindent
{\bf \S \RoutineSection. Rudiments.} 
As we begin, we make note of some general conventions used in this manuscript. 
Unless stated otherwise, the graphs considered here -- directed and undirected graphs, covering digraphs for posets -- are assumed to be finite, although we sometimes re-iterate this finiteness hypothesis for emphasis. 
For integers $a$ and $b$, we take $[a,b]_{\mathbb{Z}}$ to be the empty set if $a > b$ and to be the set $\{a,a+1,\ldots,b-1,b\}$, totally ordered in the usual way, if $a \leq b$. 
We similarly define $[a,\infty)_{\mathbb{Z}}$, $(-\infty,b]_{\mathbb{Z}}$, and other open and clopen versions of these integer intervals. 

{\bf [\S \RoutineSection.1:\! Some useful directed graph manipulations.]} 
Suppose $R$ is a simple directed graph with vertex set $\VertexSet(R)$ and directed edge set $\EdgeSet(R)$.  
Most often, our vertex sets are collections of combinatorial objects, and in these cases it makes sense to use the phrase `$\relt$ is an element of $R$' (written `$\relt \in R$') as an occasional substitute for `$\relt$ is a vertex of $R$' (written `$\relt \in \VertexSet(R)$'). 
For $\selt, \telt \in R$, say $\selt$ and $\telt$ are {\em neighbors} if there is a directed edge between them. 
Let $I$ be a set, to be thought of as a collection of colors we sometimes refer to as our {\em palette}.  
If $R$ is accompanied by a not-necessarily-surjective function $\ecolor_{R}: \EdgeSet(R) \longrightarrow I$, then we say $R$ is {\em edge-colored by the set} $I$, or $I$-{\em edge-colored} (or simply $I$-{\em colored} if there is no potential ambiguity),  and we call $I$ the {\em color palette} of $R$. 
For any $i \in I$, let $\EdgeSet_{i}(R) := \ecolor_{R}^{-1}(i)$. 
If $J$ is a subset of $I$, remove all edges from $\EdgeSet(R)$ whose colors are not in $J$; the connected components of the resulting edge-colored directed graph are called the {\em J-components} of $R$. 
For any $\telt$ in $R$ and any $J \subseteq I$, we let $\comp_{J}(\telt)$ denote the $J$-component of $R$ containing $\telt$.  
The {\em dual} $R^{*}$ is the edge-colored directed graph whose vertex set $\VertexSet(R^{*})$ is the set of symbols $\{\telt^{*}\}_{\telt{\in}R}$ together with colored edges 
$\EdgeSet_{i}(R^{*}) := \{\telt^{*} \myarrow{i} \selt^{*}\, |\, \selt \myarrow{i} \telt \in \EdgeSet_{i}(R)\}$ 
for each $i \in I$. 
Let $Q$ be another $I$-colored directed graph. 
If $R$ and $Q$ have disjoint vertex sets, then the {\em disjoint sum} $R \oplus Q$ is the edge-colored directed graph whose vertex set is the disjoint union $\VertexSet(R) \disjointunion \VertexSet(Q)$ and whose colored edges are $\EdgeSet_{i}(R) \disjointunion \EdgeSet_{i}(Q)$ for each $i \in I$.  
If $\VertexSet(Q) \subseteq \VertexSet(R)$ and $\EdgeSet_{i}(Q) \subseteq \EdgeSet_{i}(R)$ for each $i \in I$, then $Q$ is an {\em edge-colored subgraph} of $R$. 
Let $R \times Q$ denote the edge-colored directed graph whose vertex set is the Cartesian product $\{(\selt,\telt)|\selt \in R,\telt \in Q\}$ and with colored edges $(\selt_{1},\telt_{1}) \myarrow{i} (\selt_{2},\telt_{2})$ if and only if $\selt_{1} = \selt_{2}$ in $R$ with $\telt_{1} \myarrow{i} \telt_{2}$ in $Q$ or $\selt_{1} \myarrow{i} \selt_{2}$ in $R$ with $\telt_{1} = \telt_{2}$ in $Q$.
Two edge-colored directed graphs are {\em isomorphic} if there is a bijection between their vertex sets that preserves edges and edge colors.  
If $R$ is an $I$-colored directed graph and $\sigma\, :\! I \longrightarrow I'$ is a mapping of sets, then we let $R^{\sigma}$ be the $I'$-colored directed graph with $\VertexSet(R^{\sigma}) := \VertexSet(R)$, $\EdgeSet(R^{\sigma}) := \EdgeSet(R)$, and $\ecolor_{R^{\sigma}} := \sigma \circ \ecolor_{R}$. 
We call $R^{\sigma}$ a {\em recoloring} of $R$. 
Observe that $(R^{*})^{\sigma} \cong (R^{\sigma})^{*}$. 
Similarly, we say a directed graph $R$ with a function $\vcolor_{R}\, :\! \VertexSet(R) \longrightarrow I$ that assigns colors to the vertices of $R$ is a {\em vertex-colored directed graph}, in which case we say $R$ is $I$-{\em vertex-colored}.  
In this latter context, we can speak of the {\em dual vertex-colored directed graph} $R^{*}$, the {\em disjoint sum} of two vertex-colored directed graphs with disjoint vertex sets, {\em isomorphism} of vertex-colored directed graphs, {\em recoloring}, etc.   

{\bf [\S \RoutineSection.2:\! Conventional aspects of partially ordered sets.]} 
A partially ordered set, or {\em poset}, is a pair $(R,\leq_{R})$ where $R$ is a set and $\leq_{R}$ is a reflexive, anti-symmetric, and transitive relation.  
We typically refer to a poset simply by invoking the name of the underlying set $R$, as long as the partial order $\leq_{R}$ (usually just written `$\leq$') is understood. 
Elements $\uelt$ and $\velt$ of $R$ are {\em comparable} if $\uelt \leq_{R} \velt$ or $\uelt \geq_{R} \velt$. 
If all elements of $R$ are comparable, the partial order can be called a {\em total order} and $R$ can be called a {\em totally ordered set}. 
A {\em chain} (respectively, {\em anti-chain}) in $R$ is set of pairwise comparable (respectively, incomparable) elements of $R$. 
A {\em covering relation} in $R$ is any ordered pair of elements $(\xelt,\yelt)$ such that $\xelt < \yelt$ and such that $\xelt = \zelt$ or $\zelt = \yelt$ whenever $\xelt \leq \zelt \leq \yelt$ for some $\zelt \in R$. 
Such a covering relation is notated as a directed edge $\xelt \rightarrow \yelt$, and we say that $\yelt$ {\em covers} $\xelt$ and that $\xelt$ is {\em covered by} $\yelt$. 
In this case, we say that $\xelt$ is {\em below} $\yelt$, is a {\em descendant} of $\yelt$, and is {\em below the edge} $\xelt \rightarrow \yelt$, and likewise we say that $\yelt$ is {\em above} $\xelt$, is an {\em ascendant} of $\xelt$, and is {\em above the edge} $\xelt \rightarrow \yelt$. 
The {\em covering digraph} of $R$ (traditionally called its `Hasse diagram') is the directed graph whose vertex set $\VertexSet(R)$ is $R$ and whose directed edge set $\EdgeSet(R)$ is the set of covering relations.  
For convenience, we sometimes identify a poset $R$ with its covering digraph. 
As such, operations on, properties of, and language about directed graphs from the previous paragraph also apply to posets. 
If a poset $R$ is edge-colored by a set $I$, then we use the notation $\xelt \myarrow{i} \yelt$ to indicate that the color $i$ is assigned to the covering relation $\xelt \rightarrow \yelt$. 
In figures, when such an edge is depicted without an arrowhead, the implied direction is ``up.'' 

For $\xelt, \yelt \in R$, the {\em interval} $[\xelt,\yelt] = [\xelt,\yelt]_{R}$ is the set of all elements $\zelt \in R$ with $\xelt \leq_{R} \zelt \leq_{R} \yelt$, with partial order induced\footnote{The induced partial order is obtained as follows: For $\selt$ and $\telt$ in the interval $\mathfrak{I} = [\xelt,\yelt]$, we take $\selt \leq_{\mathfrak{I}} \telt$ if and only if $\selt \leq_{R} \telt$.} by $R$. 
The interval $[\xelt,\yelt]$ is empty unless $\xelt \leq \yelt$. 
 A {\em lower order ideal} (or {\em down-set}) $\mathcal{I}$ from $R$ is a subset of $R$ with the property that, if $\yelt \in \mathcal{I}$ and $\xelt \leq \yelt$, then $\xelt \in \mathcal{I}$; endow $\mathcal{I}$ with the induced partial order from $R$. 
Similarly define an {\em upper order ideal} (or {\em up-set}) from $R$, which we often denote by an `$\mathcal{F}$' for the reason that, historically, such sets have sometimes been called `filters'. 
The lower order ideal $(-\infty,\yelt] = (-\infty,\yelt]_{R}$ is the down-set consisting of all $\zelt \in R$ with $\zelt \leq \yelt$, and we call $(-\infty,\yelt]$ the {lower order ideal} (or {down-set}) {\em generated by} $\yelt$. 
Similarly define an {upper order ideal} (or {up-set}) $[\xelt,\infty) = [\xelt,\infty)_{R}$ {\em generated by} $\xelt$. 
Observe that $[\xelt,\yelt]$ and $[\xelt,\infty) \cap (-\infty,\yelt]$ coincide as sets. 
The down-set $\langle \xelt_{1},\ldots,\xelt_{m} \rangle^{\downarrow}$ is the subset $\bigcup_{i=1}^{m}(-\infty,\xelt_{i}]$ with the induced partial order from $R$. 
Similarly define the up-set $\langle \xelt_{1},\ldots,\xelt_{m} \rangle^{\uparrow}$ {\em generated by} $\{\xelt_{1},\ldots,\xelt_{m}\} \subseteq R$. 
Note that $\selt \rightarrow \telt$ in an interval $[\xelt,\yelt]$, down-set $\mathcal{I}$, or up-set $\mathcal{F}$ if and only if $\selt \rightarrow \telt$ in the poset $R$ and $\selt$ and $\telt$ reside within the appropriate interval, down-set, or up-set. 
That is the covering digraph for $[\xelt,\yelt]$, $\mathcal{I}$, or $\mathcal{F}$ is just the induced subgraph of $R$ on the vertices of $[\xelt,\yelt]$, $\mathcal{I}$, or $\mathcal{F}$ respectively.  
See \S \SubstructureSection\ for further discussion of poset substructures. 

\begin{figure}[t]
\begin{center}
{\bf \PosetAndLatticeFig}\ \  A vertex-colored poset $P$ and an edge-colored distributive lattice $L$.

{\small (The set $\{${\em 1},{\em 2}$\}$ is both the set of vertex colors for $P$ and the set of edge colors for $L$.)}

\setlength{\unitlength}{1cm}
\begin{picture}(3,3.5)
\put(0,3){\begin{picture}(3,3.5)
\put(0,3.5){$P$}
\put(4,2){\circle*{0.15}} 
\put(3.4,1.9){\footnotesize $v_{6}$}
\put(4.2,1.9){\footnotesize {\em 2}} 
\put(1,1){\circle*{0.15}}
\put(0.4,0.9){\footnotesize $v_{5}$} 
\put(1.2,0.9){\footnotesize {\em 1}} 
\put(2,2){\circle*{0.15}} 
\put(1.4,1.9){\footnotesize $v_{4}$} 
\put(2.2,1.9){\footnotesize {\em 1}}
\put(3,3){\circle*{0.15}} 
\put(2.4,2.9){\footnotesize $v_{3}$}
\put(3.2,2.9){\footnotesize {\em 1}} 
\put(0,2){\circle*{0.15}}
\put(-0.6,1.9){\footnotesize $v_{2}$} 
\put(0.2,1.9){\footnotesize {\em 2}} 
\put(1,3){\circle*{0.15}} 
\put(0.4,2.9){\footnotesize $v_{1}$} 
\put(1.2,2.9){\footnotesize {\em 2}}
\put(1,1){\line(1,1){2}} \put(0,2){\line(1,1){1}}
\put(1,1){\line(-1,1){1}} \put(2,2){\line(-1,1){1}}
\put(3,3){\line(1,-1){1}}
\end{picture}
}
\end{picture}
\hspace*{1.5in}
\setlength{\unitlength}{1.5cm}
\begin{picture}(4,6.5)
\put(0.25,5.5){$L$}
\put(1,0){\line(-1,1){1}}
\put(1,0){\line(1,1){2}}
\put(0,1){\line(1,1){2}}
\put(2,1){\line(-1,1){1}}
\put(2,1){\line(0,1){1}}
\put(1,2){\line(0,1){1}}
\put(2,2){\line(-1,1){2}}
\put(2,2){\line(1,1){2}}
\put(3,2){\line(-1,1){1}}
\put(3,2){\line(0,1){1}}
\put(1,3){\line(1,1){2}}
\put(2,3){\line(0,1){1}}
\put(3,3){\line(-1,1){2}}
\put(0,4){\line(1,1){2}}
\put(4,4){\line(-1,1){2}}
\put(2,6){\VertexForLatticeI{0}}
\put(1,5){\VertexForLatticeI{1}}
\put(3,5){\VertexForLatticeI{2}}
\put(0,4){\VertexForLatticeI{3}}
\put(2,4){\VertexForLatticeI{4}}
\put(4,4){\VertexForLatticeI{5}}
\put(1,3){\VertexForLatticeI{6}}
\put(2,3){\VertexForLatticeI{7}}
\put(3,3){\VertexForLatticeI{8}}
\put(1,2){\VertexForLatticeI{9}}
\put(2,2){\VertexForLatticeI{10}}
\put(3,2){\VertexForLatticeI{11}}
\put(0,1){\VertexForLatticeI{12}}
\put(2,1){\VertexForLatticeI{13}}
\put(1,0){\VertexForLatticeI{14}}
\put(1,5){\NEEdgeLabelForLatticeI{{\em 2}}}
\put(3,5){\NWEdgeLabelForLatticeI{{\em 1}}}
\put(0,4){\NEEdgeLabelForLatticeI{{\em 2}}}
\put(2,4){\NWEdgeLabelForLatticeI{{\em 1}}}
\put(2,4){\NEEdgeLabelForLatticeI{{\em 2}}}
\put(4,4){\NWEdgeLabelForLatticeI{{\em 2}}}
\put(1,3){\NEEdgeLabelForLatticeI{{\em 2}}}
\put(1,3){\NWEdgeLabelForLatticeI{{\em 1}}}
\put(2,3){\VerticalEdgeLabelForLatticeI{{\em 1}}}
\put(3,3){\NWEdgeLabelForLatticeI{{\em 2}}}
\put(3,3){\NEEdgeLabelForLatticeI{{\em 2}}}
\put(1,2){\VerticalEdgeLabelForLatticeI{{\em 1}}}
\put(1.25,2.25){\NEEdgeLabelForLatticeI{{\em 2}}}
\put(2.2,1.8){\NWEdgeLabelForLatticeI{{\em 2}}}
\put(1.8,1.8){\NEEdgeLabelForLatticeI{{\em 2}}}
\put(3,2){\VerticalEdgeLabelForLatticeI{{\em 1}}}
\put(2.75,2.25){\NWEdgeLabelForLatticeI{{\em 2}}}
\put(0,1){\NEEdgeLabelForLatticeI{{\em 1}}}
\put(2,1){\VerticalEdgeLabelForLatticeI{{\em 1}}}
\put(2,1){\NWEdgeLabelForLatticeI{{\em 2}}}
\put(2,1){\NEEdgeLabelForLatticeI{{\em 2}}}
\put(1,0){\NWEdgeLabelForLatticeI{{\em 2}}}
\put(1,0){\NEEdgeLabelForLatticeI{{\em 1}}}
\end{picture}
\end{center}
\end{figure}

For the remainder of this subsection, regard $R$ to be an $I$-edge-colored poset. 
A {\em path} from $\selt$ to $\telt$ in $R$ is a sequence $\mathcal{P} = (\selt = \xelt_{0}, \xelt_{1}, \ldots, \xelt_{k} = \telt)$ such that for all $j \in [1,k]_{\mathbb{Z}}$ we have either $\xelt_{j-1} \myarrow{i_j} \xelt_{j}$ or $\xelt_{j} \myarrow{i_j} \xelt_{j-1}$, where $(i_{j})_{j=1}^{k}$ is a sequence of colors from $I$.  
This path has {\em length} $k$, written $\pathlength(\mathcal{P})=k$, and we allow paths to have length $0$.  
For any $i \in I$, we let $\mya_{i}(\mathcal{P}) := \myabs\{j \in [1,k]_{\mathbb{Z}}\, |\, \xelt_{j-1} \myarrow{i_j} \xelt_{j} \mbox{ in $\mathcal{P}$ and } i_j = i\}\myabs$ (a count of `ascents' of color $i$ in the path) and $\myd_{i}(\mathcal{P}) := \myabs\{j \in [1,k]_{\mathbb{Z}}\, |\, \xelt_{j} \myarrow{i_j} \xelt_{j-1} \mbox{ in $\mathcal{P}$ and } i_j = i\}\myabs$ (a count of `descents' of color $i$). 
Of course, $\pathlength(\mathcal{P}) = \sum_{i \in I}(\mya_{i}(\mathcal{P}) + \myd_{i}(\mathcal{P}))$. 
When $\sum_{i \in I}\myd_{i}(\mathcal{P}) = 0$ (respectively, $\sum_{i \in I}\mya_{i}(\mathcal{P}) = 0$), we say $\mathcal{P}$ is {\em nondecreasing} and is a path from $\selt$ {\em up to} $\telt$ (resp., {\em nonincreasing} and a path from $\selt$ {\em down to} $\telt$). 
Say $\mathcal{P}$ is {\em simple} if each vertex appearing in the path appears exactly once; if $\pathlength(\mathcal{P}) > 0$ and only the initial and terminal vertices coincide, we say $\mathcal{P}$ is a {\em simple cycle}. 
A {\em saturated chain} in $R$ has, as its elements, the vertices visited by some nondecreasing path together with the total order on these elements inherited from $R$. 
If $\selt$ and $\telt$ are within the same connected component of $R$, then the distance $\dist(\selt,\telt)$ between $\selt$ and $\telt$ is the minimum length achieved when all paths from $\selt$ to $\telt$ in $R$ are considered; any minimum-length-achieving path is {\em shortest}.  
Our poset $R$ is {\em diamond-colored} if, whenever \parbox{1.1cm}{\begin{center}
\setlength{\unitlength}{0.2cm}
\begin{picture}(4,3.5)
\put(2,0){\circle*{0.5}} \put(0,2){\circle*{0.5}}
\put(2,4){\circle*{0.5}} \put(4,2){\circle*{0.5}}
\put(0,2){\line(1,1){2}} \put(2,0){\line(-1,1){2}}
\put(4,2){\line(-1,1){2}} \put(2,0){\line(1,1){2}}
\put(0.75,0.55){\em \small k} \put(2.7,0.7){\em \small l}
\put(0.7,2.7){\em \small i} \put(2.75,2.55){\em \small j}
\end{picture} \end{center}} is an edge-colored subgraph of the covering digraph for $R$, then $i = l$ and $j = k$.  

\begin{figure}[t]
\begin{center}
{\bf \DualAndRecolorFig}\ \   $L^{*}$ and $(L^{*})^{\sigma}$ for the lattice $L$ from \PosetAndLatticeFig.

{\small (Here $\sigma(\mbox{\em 1}) = \alpha$ and $\sigma(\mbox{\em 2}) = \beta$.)}

\setlength{\unitlength}{1.5cm}
\begin{picture}(4,6.5)
\put(-0.15,3.35){$L^{*}$}
\put(1,6){\circle*{0.1}}
\put(0,5){\circle*{0.1}}
\put(2,5){\circle*{0.1}}
\put(1,4){\circle*{0.1}}
\put(2,4){\circle*{0.1}}
\put(3,4){\circle*{0.1}}
\put(1,3){\circle*{0.1}}
\put(2,3){\circle*{0.1}}
\put(3,3){\circle*{0.1}}
\put(0,2){\circle*{0.1}}
\put(2,2){\circle*{0.1}}
\put(4,2){\circle*{0.1}}
\put(1,1){\circle*{0.1}}
\put(3,1){\circle*{0.1}}
\put(2,0){\circle*{0.1}}
\put(2,0){\line(-1,1){1}}
\put(2,0){\NWEdgeLabelForLatticeI{{\em 2}}}
\put(2,0){\line(1,1){1}}
\put(2,0){\NEEdgeLabelForLatticeI{{\em 1}}}
\put(1,1){\line(-1,1){1}}
\put(1,1){\NWEdgeLabelForLatticeI{{\em 2}}}
\put(1,1){\line(1,1){1}}
\put(1,1){\NEEdgeLabelForLatticeI{{\em 1}}}
\put(3,1){\line(-1,1){1}}
\put(3,1){\NWEdgeLabelForLatticeI{{\em 2}}}
\put(3,1){\line(1,1){1}}
\put(3,1){\NEEdgeLabelForLatticeI{{\em 2}}}
\put(0,2){\line(1,1){1}}
\put(0,2){\NEEdgeLabelForLatticeI{{\em 1}}}
\put(2,2){\line(-1,1){1}}
\put(2,2){\NWEdgeLabelForLatticeI{{\em 2}}}
\put(2,2){\line(0,1){1}}
\put(2,2){\VerticalEdgeLabelForLatticeI{{\em 1}}}
\put(2,2){\line(1,1){1}}
\put(2,2){\NEEdgeLabelForLatticeI{{\em 2}}}
\put(4,2){\line(-1,1){1}}
\put(4,2){\NWEdgeLabelForLatticeI{{\em 2}}}
\put(1,3){\line(0,1){1}}
\put(1,3){\VerticalEdgeLabelForLatticeI{{\em 1}}}
\put(1,3){\line(1,1){1}}
\put(1.25,3.25){\NEEdgeLabelForLatticeI{{\em 2}}}
\put(2,3){\line(-1,1){1}}
\put(2.25,2.75){\NWEdgeLabelForLatticeI{{\em 2}}}
\put(2,3){\line(1,1){1}}
\put(1.75,2.75){\NEEdgeLabelForLatticeI{{\em 2}}}
\put(3,3){\line(0,1){1}}
\put(3,3){\VerticalEdgeLabelForLatticeI{{\em 1}}}
\put(3,3){\line(-1,1){1}}
\put(2.75,3.25){\NWEdgeLabelForLatticeI{{\em 2}}}
\put(1,4){\line(-1,1){1}}
\put(1,4){\NWEdgeLabelForLatticeI{{\em 1}}}
\put(1,4){\line(1,1){1}}
\put(1,4){\NEEdgeLabelForLatticeI{{\em 2}}}
\put(2,4){\line(0,1){1}}
\put(2,4){\VerticalEdgeLabelForLatticeI{{\em 1}}}
\put(3,4){\line(-1,1){1}}
\put(3,4){\NWEdgeLabelForLatticeI{{\em 2}}}
\put(0,5){\line(1,1){1}}
\put(0,5){\NEEdgeLabelForLatticeI{{\em 2}}}
\put(2,5){\line(-1,1){1}}
\put(2,5){\NWEdgeLabelForLatticeI{{\em 1}}}
\end{picture}
\hspace*{1in}
\setlength{\unitlength}{1.5cm}
\begin{picture}(4,6.5)
\put(-0.15,3.25){$(L^{*})^{\sigma}$}
\put(1,6){\circle*{0.1}}
\put(0,5){\circle*{0.1}}
\put(2,5){\circle*{0.1}}
\put(1,4){\circle*{0.1}}
\put(2,4){\circle*{0.1}}
\put(3,4){\circle*{0.1}}
\put(1,3){\circle*{0.1}}
\put(2,3){\circle*{0.1}}
\put(3,3){\circle*{0.1}}
\put(0,2){\circle*{0.1}}
\put(2,2){\circle*{0.1}}
\put(4,2){\circle*{0.1}}
\put(1,1){\circle*{0.1}}
\put(3,1){\circle*{0.1}}
\put(2,0){\circle*{0.1}}
\put(2,0){\line(-1,1){1}}
\put(2,0){\NWEdgeLabelForLatticeI{$\beta$}}
\put(2,0){\line(1,1){1}}
\put(2,0){\NEEdgeLabelForLatticeI{$\alpha$}}
\put(1,1){\line(-1,1){1}}
\put(1,1){\NWEdgeLabelForLatticeI{$\beta$}}
\put(1,1){\line(1,1){1}}
\put(1,1){\NEEdgeLabelForLatticeI{$\alpha$}}
\put(3,1){\line(-1,1){1}}
\put(3,1){\NWEdgeLabelForLatticeI{$\beta$}}
\put(3,1){\line(1,1){1}}
\put(3,1){\NEEdgeLabelForLatticeI{$\beta$}}
\put(0,2){\line(1,1){1}}
\put(0,2){\NEEdgeLabelForLatticeI{$\alpha$}}
\put(2,2){\line(-1,1){1}}
\put(2,2){\NWEdgeLabelForLatticeI{$\beta$}}
\put(2,2){\line(0,1){1}}
\put(2,2){\VerticalEdgeLabelForLatticeI{$\alpha$}}
\put(2,2){\line(1,1){1}}
\put(2,2){\NEEdgeLabelForLatticeI{$\beta$}}
\put(4,2){\line(-1,1){1}}
\put(4,2){\NWEdgeLabelForLatticeI{$\beta$}}
\put(1,3){\line(0,1){1}}
\put(1,3){\VerticalEdgeLabelForLatticeI{$\alpha$}}
\put(1,3){\line(1,1){1}}
\put(1.25,3.25){\NEEdgeLabelForLatticeI{$\beta$}}
\put(2,3){\line(-1,1){1}}
\put(2.25,2.75){\NWEdgeLabelForLatticeI{$\beta$}}
\put(2,3){\line(1,1){1}}
\put(1.75,2.75){\NEEdgeLabelForLatticeI{$\beta$}}
\put(3,3){\line(0,1){1}}
\put(3,3){\VerticalEdgeLabelForLatticeI{$\alpha$}}
\put(3,3){\line(-1,1){1}}
\put(2.75,3.25){\NWEdgeLabelForLatticeI{$\beta$}}
\put(1,4){\line(-1,1){1}}
\put(1,4){\NWEdgeLabelForLatticeI{$\alpha$}}
\put(1,4){\line(1,1){1}}
\put(1,4){\NEEdgeLabelForLatticeI{$\beta$}}
\put(2,4){\line(0,1){1}}
\put(2,4){\VerticalEdgeLabelForLatticeI{$\alpha$}}
\put(3,4){\line(-1,1){1}}
\put(3,4){\NWEdgeLabelForLatticeI{$\beta$}}
\put(0,5){\line(1,1){1}}
\put(0,5){\NEEdgeLabelForLatticeI{$\beta$}}
\put(2,5){\line(-1,1){1}}
\put(2,5){\NWEdgeLabelForLatticeI{$\alpha$}}
\end{picture}
\end{center}
\end{figure}

\begin{figure}[thb]
\begin{center}
{\bf \ComponentFig}\ \  The disjoint sum of the {\em 2}-components   

of the edge-colored lattice $L$ from \PosetAndLatticeFig. 

\setlength{\unitlength}{1cm}
\hspace*{0.2in}
\begin{picture}(14,3.5)
\put(2,2.5){\VertexForLatticeI{0}}
\put(1,1.5){\VertexForLatticeI{1}}
\put(1,1.5){\NEEdgeLabelForLatticeII{{\em 2}}}
\put(0,0.5){\VertexForLatticeI{3}}
\put(0,0.5){\NEEdgeLabelForLatticeII{{\em 2}}}
\put(6,3){\VertexForLatticeI{2}} 
\put(5,2){\VertexForLatticeI{4}}
\put(5,2){\NEEdgeLabelForLatticeII{{\em 2}}}
\put(4,1){\VertexForLatticeI{6}} 
\put(4,1){\NEEdgeLabelForLatticeII{{\em 2}}}
\put(7,2){\VertexForLatticeI{5}}
\put(7,2){\NWEdgeLabelForLatticeII{{\em 2}}}
\put(6,1){\VertexForLatticeI{8}} 
\put(6,1){\NWEdgeLabelForLatticeII{{\em 2}}}
\put(6,1){\NEEdgeLabelForLatticeII{{\em 2}}}
\put(5,0){\VertexForLatticeI{10}}
\put(5,0){\NWEdgeLabelForLatticeII{{\em 2}}}
\put(5,0){\NEEdgeLabelForLatticeII{{\em 2}}}
\put(10,2.5){\VertexForLatticeI{7}}
\put(9,1.5){\VertexForLatticeI{9}}
\put(9,1.5){\NEEdgeLabelForLatticeII{{\em 2}}}
\put(11,1.5){\VertexForLatticeI{11}}
\put(11,1.5){\NWEdgeLabelForLatticeII{{\em 2}}}
\put(10,0.5){\VertexForLatticeI{13}}
\put(10,0.5){\NEEdgeLabelForLatticeII{{\em 2}}}
\put(10,0.5){\NWEdgeLabelForLatticeII{{\em 2}}}
\put(13,2){\VertexForLatticeI{12}}
\put(14,1){\VertexForLatticeI{14}}
\put(14,1){\NWEdgeLabelForLatticeII{{\em 2}}}
\put(0,0.5){\line(1,1){2}}
\put(4,1){\line(1,1){2}} 
\put(5,0){\line(1,1){2}}
\put(5,0){\line(-1,1){1}} 
\put(6,1){\line(-1,1){1}}
\put(7,2){\line(-1,1){1}}
\put(9,1.5){\line(1,1){1}} 
\put(10,0.5){\line(1,1){1}}
\put(10,0.5){\line(-1,1){1}} 
\put(11,1.5){\line(-1,1){1}}
\put(14,1){\line(-1,1){1}}
\put(2.8,1.375){$\bigoplus$} \put(7.8,1.375){$\bigoplus$}
\put(11.8,1.375){$\bigoplus$}
\end{picture}
\end{center}
\end{figure}

A simple path $\mathcal{P}$ in $R$ is a {\em mountain path}  from $\selt$ to $\telt$ if for some $j \in [0,k]_{\mathbb{Z}}$ and $\uelt := \xelt_{j}$ we have $\selt=\xelt_{0} \rightarrow \xelt_{1} \rightarrow \cdots \rightarrow \uelt \leftarrow \cdots \leftarrow \xelt_{k} = \telt$ (with edge colors suppressed), in which case we call $\uelt$ the {\em apex} of the mountain path. 
Similarly, the simple path $\mathcal{P}$ is a {\em valley path} from $\selt$ to $\telt$ if for some $j \in [0,k]_{\mathbb{Z}}$ and $\velt := \xelt_{j}$ we have $\selt=\xelt_{0} \leftarrow \xelt_{1} \leftarrow \cdots \leftarrow \velt \rightarrow \cdots \rightarrow \xelt_{k} = \telt$ (with edge colors suppressed), in which case we call $\velt$ the {\em nadir} of the valley path. 
Our poset $R$ is {\em topographically balanced} if (1) whenever $\velt \rightarrow \selt$ and $\velt \rightarrow \telt$ for distinct $\selt$ and $\telt$ in $R$, then there exists a unique $\uelt$ in $R$ such that $\selt \rightarrow \uelt$ and $\telt \rightarrow \uelt$, and (2) whenever  $\selt \rightarrow \uelt$ and $\telt \rightarrow \uelt$ for distinct $\selt$ and $\telt$ in $R$, then there exists a unique $\velt$ in $R$ such that $\velt \rightarrow \selt$ and $\velt \rightarrow \telt$.  
Informally, this just says that any length two mountain path that is not a chain is uniquely balanced by a length two valley path that is not a chain, and vice-versa. 
See \TopoBalancedRankedFigure\ for an example.

A {\em rank function} on $R$ is a surjective function $\rho : R \longrightarrow [0,l]_{\mathbb{Z}}$ (where $l \in [0,\infty)_{\mathbb{Z}}$) with the property that if $\selt \rightarrow \telt$ in $R$, then $\rho(\selt) + 1 = \rho(\telt)$. 
If such a rank function $\rho$ exists, then $R$ is a {\em ranked} poset and $l = \posetlength(R)$ is the {\em length} of $R$ with respect to $\rho$. 
The companion {\em depth function} $\delta: R \longrightarrow [0,l]_{\mathbb{Z}}$ is given by $\delta(\relt) := l-\rho(\relt)$. 
A ranked poset that is connected has a unique rank function. 
In this setting, for any path $\mathcal{P}$ from $\selt$ to $\telt$, we have \[\rho(\telt) - \rho(\selt) =  \sum_{i \in I}\big(\mya_{i}(\mathcal{P}) - \myd_{i}(\mathcal{P})\big),\]  which serves as an expression for the rank of $\telt$ whenever $\rho(\selt) = 0$. 
We also see that there is a nondecreasing path $\mathcal{Q}$ from $\selt$ up to $\telt$ if and only if $\rho(\telt) - \rho(\selt) = \pathlength(\mathcal{Q})$ if and only if $\selt \leq \telt$ if and only if $\dist(\selt,\telt) = \rho(\telt) - \rho(\selt)$. 
The polynomial $\RGF(R;q) := \sum_{\relt \in R}q^{\rho(\relt)}$ in the variable $q$ is the {\em rank-generating function} for $R$. 

In a connected and topographically balanced poset $R$, there is a method by which any shortest path $\mathcal{P}$ between vertices can be converted to a shortest path that is a mountain (respectively, valley) path $\mathcal{P}_{\mbox{\em \scriptsize mountain}}$ (resp., $\mathcal{P}_{\mbox{\em \scriptsize valley}}$). 
Here is the {\bf mountain-path-construction} algorithm. 
Suppose $\mathcal{P} = (\selt = \xelt_{0}, \xelt_{1}, \ldots, \xelt_{k-1}, \xelt_{k} = \telt)$ is a shortest path between distinct elements $\selt$ and $\telt$; in particular, $\mathcal{P}$ is simple. 
We modify $\mathcal{P}$ as follows: 
\begin{enumerate}
\item If there is no $j \in [1,k-1]_{\mathbb{Z}}$ such that $\xelt_{j-1} \leftarrow \xelt_{j} \rightarrow \xelt_{j+1}$, then return $\mathcal{P}$ as $\mathcal{P}_{\mbox{\em \scriptsize mountain}}$. 
\item Otherwise, let $j \in [1,k-1]_{\mathbb{Z}}$ be least such that $\xelt_{j-1} \leftarrow \xelt_{j} \rightarrow \xelt_{j+1}$. 
Let $\xelt$ be the unique element from $R$ such that $\xelt_{j-1} \rightarrow \xelt \leftarrow \xelt_{j+1}$. 
If $\xelt = \xelt_{k}$ with $k \leq j-2$, then we could form a new path $\mathcal{P}'$ from $\selt$ to $\telt$ that is shorter than $\mathcal{P}$ by deleting $\xelt_{k+1}, \ldots, \xelt_{j}$ from $\mathcal{P}$ so that $\xelt_{j+1}$ now immediately succeeds $\xelt_{k}$. 
Similarly, if $\xelt = \xelt_{k}$ with $k \geq j+2$, then we could form a new path $\mathcal{P}'$ from $\selt$ to $\telt$ that is shorter than $\mathcal{P}$ by deleting $\xelt_{j}, \ldots, \xelt_{k-1}$ from $\mathcal{P}$ so that $\xelt_{k}$ now immediately succeeds $\xelt_{j-1}$. 
So, $\xelt =: \xelt_{j}'$ is not visited by $\mathcal{P}$. 
Form a new path $\mathcal{P}'$ from $\mathcal{P}$ by replacing $\xelt_{j}$ with $\xelt_{j}'$.
Observe that $\mathcal{P}'$ is also a shortest path from $\selt$ to $\telt$. 
\item Return to the first step of the process, using $\mathcal{P}'$ as $\mathcal{P}$. 
\end{enumerate} 
The ``mountain-ization'' $\mathcal{P}_{\mbox{\em \scriptsize mountain}}$ returned by {\bf mountain-path-construction} has the obvious properties that it is a shortest path and mountain path from $\selt$ to $\telt$. 
The ``valley-ization'' $\mathcal{P}_{\mbox{\em \scriptsize valley}}$, which is a shortest path and valley path from $\selt$ to $\telt$, can be obtained from $\mathcal{P}$ via a {\bf valley-path-construction} algorithm as follows: Apply {\bf mountain-path-construction} to the corresponding path $\mathcal{P}^{*}$ in $R^{*}$ to get $\mathcal{P\, }^{*}_{\mbox{\em $\!\!\!$\scriptsize mountain}}$, and then let $\mathcal{P}_{\mbox{\em \scriptsize valley}}$ be the path in $R$ corresponding to $\mathcal{P\, }^{*}_{\mbox{\em $\!\!\!$\scriptsize mountain}}$. 

Notice that if $\mathcal{P}$ is already a mountain (respectively, valley) path, then $\mathcal{P} = \mathcal{P}_{\mbox{\em \scriptsize mountain}}$ (resp., $\mathcal{P} = \mathcal{P}_{\mbox{\em \scriptsize valley}}$). 
Observe that $\selt < \telt$ if and only if $\mathcal{P}_{\mbox{\em \scriptsize mountain}}$ has apex $\telt$, 
and $\selt > \telt$ if and only if $\mathcal{P}_{\mbox{\em \scriptsize valley}}$ has nadir $\telt$. 

\noindent 
{\bf \MountainValleyLemma}\ \ {\sl A connected and topographically balanced poset $R$ is ranked.  
Moreover, if $\selt$ and $\telt$ are distinct elements of $R$ and if $\mathcal{P}$ is any shortest path from $\selt$ to $\telt$, then} $\dist(\selt,\telt) = 2\rho(\aelt)-\rho(\selt)- \rho(\telt) = \rho(\selt)+\rho(\telt) - 2\rho(\nelt)$, {\sl where $\aelt$ is the apex of} $\mathcal{P}_{\mbox{\em \scriptsize mountain}}$ {\sl and $\nelt$ is the nadir of} $\mathcal{P}_{\mbox{\em \scriptsize valley}}$. 

{\em Proof.} To demonstrate the conclusions of the second sentence of the result statement, suppose $R$ is ranked. 
Then $\dist(\selt,\aelt) = \rho(\aelt)-\rho(\selt)$ and $\dist(\aelt,\telt) = \rho(\aelt)-\rho(\telt)$, since $\selt \leq \aelt$ and $\telt \leq \aelt$; in this case $\dist(\selt,\telt) = \dist(\selt,\aelt) + \dist(\aelt,\telt) = 2\rho(\aelt)-\rho(\selt)- \rho(\telt)$. 
Similarly see that $\dist(\selt,\telt) = \rho(\selt)+\rho(\telt) - 2\rho(\nelt)$. 

Now we demonstrate that $R$ is ranked. 
We begin by saying why $R$ must have a unique maximal element and a unique minimal element. 
Well, a maximal element exists since $R$ is finite. 
Suppose there are two distinct, and therefore incomparable, maximal elements.  
Then any shortest path between them can be mountainized in order to produce an element that is above both maximal elements, contradicting their maximality. 
So, there exists a unique maximal element. 
Similarly see that there exists a unique minimal element, which, from here on, we call $\melt$. 

Next we show that for any $\xelt \in R$, any nondecreasing path from $\melt$ up to $\xelt$ has the same length as any shortest path from $\melt$ to $\xelt$. 
We argue this by induction on the distance of $\xelt$ from $\melt$. 
Of course, when $\dist(\melt,\xelt) = 0$, then $\xelt = \melt$ and any nondecreasing path from $\melt$ to $\xelt$ must have length zero, which is shortest. 
Suppose now that for some nonnegative integer $r-1$, any nondecreasing path from $\melt$ up to $\xelt$ has length equal to $\dist(\melt,\xelt) \leq r-1$. 
Let $\mathcal{P} = (\melt = \xelt_{0},\ldots,\xelt_{p}=\xelt)$ be a nondecreasing path from $\melt$ up to $\xelt$, where $\dist(\melt,\xelt) = r$. 
In particular, $p \geq r$. 
Say $\mathcal{S} = (\melt = \selt_{0},\ldots,\selt_{r}=\xelt)$ is a shortest path from $\melt$ to $\xelt$. 
Then $\mathcal{S}_{\mbox{\em \scriptsize valley}}$ is a length $r$ valley path from $\melt$ to $\xelt$, and, since $\melt$ is the unique minimal element, the nadir of $\mathcal{S}_{\mbox{\em \scriptsize valley}}$ must be $\melt$. 
In particular, $\mathcal{S}_{\mbox{\em \scriptsize valley}}$ is a length $r$ nondecreasing path from $\melt$ up to $\xelt$. 
Notice that for any $i \in [0,r]_{\mathbb{Z}}$, $(\melt = \selt_{0},\ldots,\selt_{i})$ is a shortest nondecreasing path from $\melt$ up to $\selt_{i}$. 
So, if $\xelt_{p-1} = \selt_{r-1}$, then our induction hypothesis applies to $\selt_{r-1}$, and therefore $p-1 = r-1$, i.e.\ $p=r$. 
That is, $\pathlength(\mathcal{P})=r$. 
Now suppose $\xelt_{p-1}$ and $\selt_{r-1}$ are distinct. 
Since $\xelt_{p-1} \rightarrow \xelt \leftarrow \selt_{r-1}$, then there exists a unique $\relt$ in $R$ such that $\xelt_{p-1} \leftarrow \relt \rightarrow \selt_{r-1}$. 
If we take any nondecreasing path $\mathcal{Q}$ from $\melt$ up to $\relt$, then we can append $\selt_{r-1}$ to get a nondecreasing path from $\melt$ up to $\selt_{r-1}$. 
The latter path has length $r-1$ by our induction hypothesis, and therefore the given nondecreasing path $\mathcal{Q}$ from $\melt$ up to $\relt$ has length $r-2$. 
So, if we append $\xelt_{p-1}$ to $\mathcal{Q}$, we get a path of length $r-1$ from $\melt$ up to $\xelt_{p-1}$. 
Again by our induction hypothesis, $\dist(\melt,\xelt_{p-1}) = r-1$, so we can conclude that $(\melt = \xelt_{0},\ldots,\xelt_{p-1})$ has length $r-1$, i.e.\ $p=r$. 
That is, $\mathcal{P}$ has length $r$. 
This completes the induction argument. 

Now define $\rho: R \longrightarrow [0,\infty)_{\mathbb{Z}}$ by $\rho(\xelt) := \dist(\melt,\xelt)$. 
We claim that, with a suitable restriction on the target set $[0,\infty)_{\mathbb{Z}}$, $\rho$ is a rank function. 
Suppose $\xelt \rightarrow \yelt$. 
Any nondecreasing path $\mathcal{P}$ from $\melt$ up to $\xelt$ has length $\dist(\melt,\xelt)$, as argued above. 
If we append $\yelt$ to the path $\mathcal{P}$, we get a nondecreasing path of length $\dist(\melt,\xelt)+1$ from $\melt$ up to $\yelt$. 
Therefore, $\dist(\melt,\yelt) = \dist(\melt,\xelt)+1$, which means $\rho(\xelt)+1=\rho(\yelt)$. 
So, $\rho$ is a rank function.\hfill\QED

\begin{figure}[ht]
\begin{center}
{\bf \TopoBalancedRankedFigure}\ \ A connected and (strongly) topographically balanced poset that is not a lattice.

\setlength{\unitlength}{0.4in} 
\begin{picture}(4,4.5)
\put(2,0){\circle*{0.1}} 
\put(1,1){\circle*{0.1}}
\put(2,1){\circle*{0.1}} 
\put(3,1){\circle*{0.1}} 
\put(0,2){\circle*{0.1}} 
\put(2,2){\circle*{0.1}} 
\put(4,2){\circle*{0.1}} 
\put(1,3){\circle*{0.1}} 
\put(2,3){\circle*{0.1}} 
\put(3,3){\circle*{0.1}} 
\put(2,4){\circle*{0.1}} 
\put(2,0){\line(-1,1){2}}
\put(2,0){\line(1,1){2}}
\put(0,2){\line(1,1){2}}
\put(4,2){\line(-1,1){2}}
\put(1,1){\line(1,1){2}}
\put(3,1){\line(-1,1){2}}
\put(2,0){\line(0,1){1}}
\put(2,4){\line(0,-1){1}}
\put(2,1){\line(-2,1){2}}
\put(2,1){\line(2,1){2}}
\put(2,3){\line(-2,-1){2}}
\put(2,3){\line(2,-1){2}}

\vspace*{-0.1in}
\end{picture} 
\end{center} 
\end{figure}

{\bf [\S \RoutineSection.3:\! Lattice basics.]} 
A {\em lattice} is a nonempty poset for which any two elements $\selt$ and $\telt$ have a unique least upper bound $\selt \vee \telt$ (the {\em join} of $\selt$ and $\telt$) and a unique greatest lower bound $\selt \wedge \telt$ (the {\em meet} of $\selt$ and $\telt$).  
That is, $(\selt \vee \telt) \leq \uelt$ whenever $\selt \leq \uelt$ and $\telt \leq \uelt$, and $\relt \leq (\selt \wedge \telt)$ whenever $\relt \leq \selt$ and $\relt \leq \telt$.  
A lattice $L$ is necessarily connected, and finiteness implies that there is a unique maximal element $\mymax(L)$ and a unique minimal element $\mymin(L)$. 
Associativity of the meet and join operations follow easily from transitivity and antisymmetry of the partial order on $L$.  
That is, for any $\relt, \selt, \telt \in L$, we have $\relt \wedge (\selt \wedge \telt) = (\relt \wedge \selt) \wedge \telt$ and $\relt \vee (\selt \vee \telt) = (\relt \vee \selt) \vee \telt$.  
Thus, for a nonempty subset $S$ of $L$, the meet $\bigwedge_{\mbox{\tiny $\selt \in S$}}(\selt)$ and the join $\bigvee_{\mbox{\tiny $\selt \in S$}}(\selt)$ are well-defined. 
One reason we require lattices to be nonempty is so that we can take $\bigwedge_{\mbox{\tiny $\selt \in S$}}(\selt) = \mymax(L)$ and $\bigvee_{\mbox{\tiny $\selt \in S$}}(\selt) = \mymin(L)$ if $S$ is empty.\footnote{Another reason is so that the Fundamental Theorem of Finite Diamond-colored Distributive Lattices (see \FundamentalTheorem\ below) affords a one-to-one correspondence between finite $I$-edge-colored diamond-colored distributive lattices and finite $I$-vertex-colored posets.}

A lattice $L$ is {\em modular} if it satisfies either of the equivalent conditions from the first sentence of the following theorem, a result which has precursors in Theorem 1 of \cite{Alvarez}, \S 2 of \cite{DuffusRival}, \S 3 of \cite{HG}, \S 3.3 of \cite{StanleyText}, etc.  
There are other ways to define modularity, but the definition we offer here seems to be the most natural for the algebraic contexts in which we work. 
See \TopoBalancedRankedFigure\ for an example of a connected and topographically balanced poset that is not a lattice. 

\noindent 
{\bf \ModularLatticeTheorem}\ \  {\sl Let $L$ be a lattice, possibly edge-colored. 
Then $L$ is topographically balanced if and only if $L$ is ranked with unique rank function $\rho$ satisfying} 
\[2\rho(\selt \vee \telt) - \rho(\selt) - \rho(\telt) = \rho(\selt) + \rho(\telt) - 2\rho(\selt \wedge \telt)\] 
{\sl for all $\selt, \telt \in L$. 
Assume now that these equivalent conditions hold. 
Let $l$ be the length of $L$ with respect to $\rho$, and let $\mathcal{P}$ be a shortest path in $L$ from an element $\selt$ to an element $\telt$. Then} 
\[\pathlength(\mathcal{P}) = \displaystyle \dist(\selt,\telt) = 2\rho(\selt \vee \telt) - \rho(\selt) - \rho(\telt) = \rho(\selt) + \rho(\telt) - 2\rho(\selt \wedge \telt).\] 
{\sl In particular, a mountain path from $\selt$ to $\telt$ is a shortest path if and only if its apex is $\selt \vee \telt$, and a  valley path from $\selt$ to $\telt$ is a shortest path if and only if its nadir is $\selt \wedge \telt$. 
We have} $\dist(\selt,\telt) \leq l$ {\sl for all $\selt$ and $\telt$ in $L$. 
Moreover,} $\dist(\selt,\telt) = l$ {\sl if and only if} $\selt \vee \telt = \mymax(L)$ {\sl and} $\selt \wedge \telt = \mymin(L)$. 

{\em Proof.} The equivalence asserted in the first sentence of the proposition follows from \S 3.3 of \cite{StanleyText} (see Proposition 3.3.2 of that text and the subsequent paragraphs). 
Assuming $L$ is modular, we now establish the second claim of the proposition statement. 
Suppose a path $\mathcal{P} = (\xelt_{0} = \selt, \xelt_{1}, \ldots, \xelt_{k-1}, \xelt_{k} = \telt)$ from $\selt$ to some distinct $\telt$ is shortest. 
Let $\aelt$ be the apex of the mountain-ization $\mathcal{P}_{\mbox{\em \scriptsize mountain}}$ of $\mathcal{P}$. 
Then $\pathlength(\mathcal{P}) = \dist(\selt,\telt) = 2\rho(\aelt)-\rho(\selt)-\rho(\telt)$, where the former equality is merely the definition of `shortest path' and the latter equality follows from \MountainValleyLemma. 
Now let $\mathcal{Q}$ be a path formed by a walk from $\selt$ up to $\selt \vee \telt$ followed by a walk from $\selt \vee \telt$ down to $\telt$. 
So, $\pathlength(\mathcal{Q}) = \left(\rho(\selt\vee \telt) - \rho(\selt)\right) + \left(\rho(\selt\vee \telt) - \rho(\telt)\right) = 2\rho(\selt \vee \telt) - \rho(\selt) - \rho(\telt)$. 
We see, then, that $2\rho(\aelt)-\rho(\selt)-\rho(\telt) = \dist(\selt,\telt) \leq \pathlength(\mathcal{Q}) = 2\rho(\selt \vee \telt) - \rho(\selt) - \rho(\telt)$, which means that $\rho(\aelt) \leq \rho(\selt \vee \telt)$. 
But since $\selt \leq \aelt$ and $\telt \leq \aelt$, then $\selt \vee \telt \leq \aelt$, and in particular $\rho(\selt \vee \telt) \leq \rho(\aelt)$. 
In any ranked poset, distinct elements with the same rank are incomparable. 
Therefore $(\selt \vee \telt) \leq \aelt$ together with $\rho(\selt \vee \telt) = \rho(\aelt)$ implies that $(\selt \vee \telt) = \aelt$. 
We conclude that $\pathlength(\mathcal{P}) = \dist(\selt,\telt) = 2\rho(\selt \vee \telt)-\rho(\selt)-\rho(\telt)$. 
Use the valley-ization $\mathcal{P}_{\mbox{\em \scriptsize mountain}}$ in a similar way to conclude that $\pathlength(\mathcal{P}) = \dist(\selt,\telt) = \rho(\selt) + \rho(\telt) - 2\rho(\selt \wedge \telt)$.

As part of our argument in the previous paragraph, we observed that any mountain path from $\selt$ to $\telt$ with apex $\selt \vee \telt$ is a shortest path from $\selt$ to $\telt$, and we saw that any mountain path from $\selt$ to $\telt$ that is shortest must have $\selt \vee \telt$ as its apex. 
Similarly see the corresponding statements for valley paths from $\selt$ to $\telt$. 

To prove that $\dist(\selt,\telt) \leq l$, it suffices to consider the following non-mutually-exclusive cases: 
(a) $l-\rho(\selt) - \rho(\telt) \leq 0$ and (b) $l-\rho(\selt) - \rho(\telt) \geq 0$. 
In case (a), we have $\dist(\selt,\telt) = 2\rho(\selt \vee \telt) - \rho(\selt) - \rho(\telt) = [\rho(\selt \vee \telt) - \rho(\selt)] + [\rho(\selt \vee \telt) - \rho(\telt)] \leq [l-\rho(\selt)] + [l-\rho(\telt)] = l + [l - \rho(\selt) - \rho(\telt)] \leq l$. 
In case (b), we have $\dist(\selt,\telt) = \rho(\selt) + \rho(\telt) - 2\rho(\selt \wedge \telt) = [\rho(\selt) - \rho(\selt \wedge \telt)] + [\rho(\telt) - \rho(\selt \wedge \telt)] \leq \rho(\selt) + \rho(\telt) \leq l$. 
This completes our case analysis proof that $\dist(\selt,\telt) \leq l$. 

Suppose $\dist(\selt,\telt) = l$. 
Then the facts that $2\rho(\selt \vee \telt) - \rho(\selt) - \rho(\telt) = l$ and that $\rho(\selt \vee \telt) \leq l$ imply that $\rho(\selt)+\rho(\telt) \leq l$. 
Similarly, the facts that $\rho(\selt) + \rho(\telt) - 2\rho(\selt \wedge \telt) = l$ and that $\rho(\selt \wedge \telt) \geq 0$ imply that $\rho(\selt)+\rho(\telt) \geq l$. 
So, $\rho(\selt)+\rho(\telt) = l$. 
Therefore both cases (a) and (b) of the previous paragraph apply here, with all of the inequalities of that paragraph becoming equalities. 
In particular, $\rho(\selt \vee \telt) = l$ and $\rho(\selt \wedge \telt) = 0$, and therefore $\selt \vee \telt = \mymax(L)$ and $\selt \wedge \telt = \mymin(L)$. 
Conversely, suppose $\selt \vee \telt = \mymax(L)$ and $\selt \wedge \telt = \mymin(L)$. 
Observe, then, that $2\rho(\selt \vee \telt) - \rho(\selt) - \rho(\telt) = \rho(\selt) + \rho(\telt) - 2\rho(\selt \wedge \telt)$ reduces to $\rho(\selt) + \rho(\telt) = l$. 
So, $\dist(\selt,\telt) = \rho(\selt) + \rho(\telt) - 2\rho(\selt \wedge \telt) = l - 2 \cdot 0 = l$.\hfill\QED

Modular lattices -- which we have defined as lattices meeting either of the equivalent conditions stated at the top of \ModularLatticeTheorem\ --  are often defined using the so-called `modular law', which we present in the next result as a characterizing identity. 

\noindent 
{\bf \ModularLawLemma}\ \ {\sl Let $L$ be a lattice.  Then $L$ is modular if and only if for all $\xelt, \yelt, \zelt \in L$ with $\xelt \leq \zelt$ we have $\xelt \vee (\yelt \wedge \zelt) = (\xelt \vee \yelt) \wedge \zelt$.} 

{\em Proof.} To begin, assume $L$ is modular, and let $\xelt, \yelt, \zelt \in L$ with $\xelt \leq \zelt$. 
Since $\xelt \leq (\xelt \vee \yelt)$ and $(\yelt \wedge \zelt) \leq (\xelt \vee \yelt)$, then $\xelt \vee (\yelt \wedge \zelt) \leq (\xelt \vee \yelt)$. 
By \ModularLatticeTheorem, a shortest path $\mathcal{P}_{1}$ from $\xelt$ to $\yelt$ is a mountain path that passes through $\xelt \vee (\yelt \wedge \zelt)$ as it goes from $\xelt$ up to $\xelt \vee \yelt$ and then from $\xelt \vee \yelt$ down to $\yelt$. 
Furthermore, we have $\xelt \leq \zelt$ (given) and $(\yelt \wedge \zelt) \leq \zelt$, and so $\xelt \vee (\yelt \wedge \zelt) \leq \zelt$. 
Again by \ModularLatticeTheorem, a shortest path $\mathcal{P}_{2}$ from $\zelt$ to $\yelt$ is a valley path that passes through $\xelt \vee (\yelt \wedge \zelt)$ as it goes from $\zelt$ down to $\yelt \wedge \zelt$ and then from $\yelt \wedge \zelt$ up to $\yelt$. 
In particular, we see that the part of $\mathcal{P}_{1}$ that goes from $\xelt \vee (\yelt \wedge \zelt)$ up to $\xelt \vee \yelt$ and then down to $\yelt$ is a shortest path from $\xelt \vee (\yelt \wedge \zelt)$ to $\yelt$, and similarly the part of $\mathcal{P}_{2}$ that goes from $\xelt \vee (\yelt \wedge \zelt)$ down to $\yelt \wedge \zelt$ and then up to $\yelt$ is also a shortest path from $\xelt \vee (\yelt \wedge \zelt)$ to $\yelt$. 
Therefore, a path $\mathcal{P}_{3}$ that goes from $\zelt$ down to $\xelt \vee (\yelt \wedge \zelt)$ and then up to $\xelt \vee \yelt$ and then down to $\yelt$ is a shortest path from $\zelt$ to $\yelt$, and hence the part of $\mathcal{P}_{3}$ that goes from $\zelt$ down to $\xelt \vee (\yelt \wedge \zelt)$ and then up to $\xelt \vee \yelt$ is valley path with $\xelt \vee (\yelt \wedge \zelt)$ as its nadir and is a shortest path from $\zelt$ to $\xelt \vee \yelt$. 
That is, $(\xelt \vee \yelt) \wedge \zelt = \xelt \vee (\yelt \wedge \zelt)$. 

For the converse, we suppose our given lattice $L$ has the property that $\xelt \vee (\yelt \wedge \zelt) = (\xelt \vee \yelt) \wedge \zelt$ for all $\xelt, \yelt, \zelt \in L$ with $\xelt \leq \zelt$. 
Suppose we are given $\selt, \telt, \uelt \in L$ with $\selt \rightarrow \uelt \leftarrow \telt$. 
Let $\relt := \selt \wedge \telt$. 
We claim that $\selt \leftarrow \relt \rightarrow \telt$. 
If so, then $\relt$ must be the unique element covered by both $\selt$ and $\telt$. 
(If $\selt \leftarrow \relt' \rightarrow \telt$, then $\relt' \leq (\selt \wedge \telt) = \relt$, but the facts that $\relt' \leq \relt < \telt$ and $\relt' \rightarrow \telt$ force $\relt' = \relt$.) 
Suppose that $\relt < \xelt \leq \selt$. 
Now, it is clear that $\xelt \vee \telt \leq \uelt$. 
If the preceding inequality is strict, then we will have $\telt \leq (\xelt \vee \telt) < \uelt$, and since $\telt \rightarrow \uelt$, it will follow that $\telt = \xelt \vee \telt$, i.e.\ $\xelt \leq \telt$. 
But then $\xelt \leq \selt$ and $\xelt \leq \telt$ together mean that $\xelt \leq (\selt \wedge \telt) = \relt$, contrary to our hypothesis that $\relt < \xelt$. 
So the inequality $\xelt \vee \telt \leq \uelt$ cannot be strict, hence $\xelt \vee \telt = \uelt$. 
Therefore $(\xelt \vee \telt) \wedge \selt = \uelt \wedge \selt = \selt$. 
We also have $\xelt \vee (\telt \wedge \selt) = \xelt \vee \relt = \xelt$.  
But since $\xelt \leq \selt$, then our hypothesized property tells us that $\xelt \vee (\telt \wedge \selt) = (\xelt \vee \telt) \wedge \selt$, and hence $\xelt = \selt$. 
That is, $\relt \rightarrow \selt$. 
An entirely similar argument shows that $\relt \rightarrow \telt$. 
This completes the proof.\hfill\QED

As with modular lattices, there are many ways to define `distributive' lattices. 
Here, we use the definition most appropriate to the `distributive' nomenclature. 
Say a lattice $L$ is {\em distributive} if and only if for all $\xelt, \yelt, \zelt \in L$ we have $\xelt \vee (\yelt \wedge \zelt) = (\xelt \vee \yelt) \wedge (\xelt \vee \zelt)$ or for all $\xelt, \yelt, \zelt \in L$ we have $\xelt \wedge (\yelt \vee \zelt) = (\xelt \wedge \yelt) \vee (\xelt \wedge \zelt)$. 
Each of the two latter identities is a {\em distributive law}. 
Further, say $L$ is {\em strongly topographically balanced} if it is topographically balanced and we have $\{\selt,\telt\} = \{\selt',\telt'\}$ whenever$\,$ \parbox{1.4cm}{\begin{center}
\setlength{\unitlength}{0.2cm}
\begin{picture}(6.5,3.5)
\put(3,0){\circle*{0.5}} \put(1,2){\circle*{0.5}}
\put(3,4){\circle*{0.5}} \put(5,2){\circle*{0.5}}
\put(1,2){\line(1,1){2}} \put(3,0){\line(-1,1){2}}
\put(5,2){\line(-1,1){2}} \put(3,0){\line(1,1){2}}
\put(3.5,-0.75){\footnotesize $\relt$} \put(5.75,1.75){\footnotesize $\telt$}
\put(3.5,4){\footnotesize $\uelt$} \put(-0.5,1.75){\footnotesize $\selt$}
\end{picture} \end{center}} and$\,$ \parbox{1.4cm}{\begin{center}
\setlength{\unitlength}{0.2cm}
\begin{picture}(6.5,3.5)
\put(3,0){\circle*{0.5}} \put(1,2){\circle*{0.5}}
\put(3,4){\circle*{0.5}} \put(5,2){\circle*{0.5}}
\put(1,2){\line(1,1){2}} \put(3,0){\line(-1,1){2}}
\put(5,2){\line(-1,1){2}} \put(3,0){\line(1,1){2}}
\put(3.5,-0.75){\footnotesize $\relt$} \put(5.75,1.75){\footnotesize $\telt'$}
\put(3.5,4){\footnotesize $\uelt$} \put(-0.5,1.75){\footnotesize $\selt'$}
\end{picture} \end{center}} are diamonds in $L$ with $\selt \not= \telt$ and $\selt' \not= \telt'$. 

\noindent 
{\bf \DistributiveIsModularLemma}\ \ {\sl (1) A lattice is distributive if and only if it satisfies both distributive laws. 
(2) Any distributive lattice is modular and strongly topographically balanced.} 

All parts of this result are well known, particularly if we regard ``strong topographical balance'' for distributive lattices to be a less general version of a Birkhoff result we consider in \S \SubstructureSection.2. 

{\em Proof of \DistributiveIsModularLemma.} 
For {\sl (1)}, assume that $\xelt \vee (\yelt \wedge \zelt) = (\xelt \vee \yelt) \wedge (\xelt \vee \zelt)$ for all $\xelt, \yelt, \zelt \in L$. 
Then, $(\xelt \wedge \yelt) \vee (\xelt \wedge \zelt) = [(\xelt \wedge \yelt) \vee \xelt] \wedge [(\xelt \wedge \yelt) \vee \zelt)] = 
(\xelt \vee \xelt) \wedge (\yelt \vee \xelt) \wedge (\xelt \vee \zelt) \wedge (\yelt \vee \zelt) = \big(\xelt \wedge (\yelt \vee \xelt) \wedge (\xelt \vee \zelt)\big) \wedge (\yelt \vee \zelt) = \xelt \wedge (\yelt \vee \zelt)$. 
That is, the defining distributive law implies the second of our distributive laws. 
The converse can be seen in a similar way. 
For {\sl (2)}, let $\xelt, \yelt, \zelt \in L$ with $\xelt \leq \zelt$. 
Then $\xelt \vee (\yelt \wedge \zelt) = (\xelt \vee \yelt) \wedge (\xelt \vee \zelt) = (\xelt \vee \yelt) \wedge \zelt$, so $L$ satisfies the modular law. 
Next suppose that for $\relt, \selt, \selt', \telt, \telt', \uelt$ in $L$ with $\selt \not= \telt$ and $\selt' \not= \telt'$ we have diamonds$\,$ \parbox{1.4cm}{\begin{center}
\setlength{\unitlength}{0.2cm}
\begin{picture}(6.5,3.5)
\put(3,0){\circle*{0.5}} \put(1,2){\circle*{0.5}}
\put(3,4){\circle*{0.5}} \put(5,2){\circle*{0.5}}
\put(1,2){\line(1,1){2}} \put(3,0){\line(-1,1){2}}
\put(5,2){\line(-1,1){2}} \put(3,0){\line(1,1){2}}
\put(3.5,-0.75){\footnotesize $\relt$} \put(5.75,1.75){\footnotesize $\telt$}
\put(3.5,4){\footnotesize $\uelt$} \put(-0.5,1.75){\footnotesize $\selt$}
\end{picture} \end{center}} and$\,$ \parbox{1.4cm}{\begin{center}
\setlength{\unitlength}{0.2cm}
\begin{picture}(6.5,3.5)
\put(3,0){\circle*{0.5}} \put(1,2){\circle*{0.5}}
\put(3,4){\circle*{0.5}} \put(5,2){\circle*{0.5}}
\put(1,2){\line(1,1){2}} \put(3,0){\line(-1,1){2}}
\put(5,2){\line(-1,1){2}} \put(3,0){\line(1,1){2}}
\put(3.5,-0.75){\footnotesize $\relt$} \put(5.75,1.75){\footnotesize $\telt'$}
\put(3.5,4){\footnotesize $\uelt$} \put(-0.5,1.75){\footnotesize $\selt'$}
\end{picture} \end{center}} in $L$. 
If $\{\selt,\telt\} \not= \{\selt',\telt'\}$, then we may assume, without loss of generality, that $\selt$, $\telt$, and $\selt'$ are all distinct. 
Observe that $\relt$ is the meet of any two of $\selt, \telt, \selt', \telt'$ and $\uelt$ is their join. 
So, $\selt \vee (\telt \wedge \selt') = \selt \vee \relt = \selt$ while $(\selt \vee \telt) \wedge (\selt \vee \selt') = \uelt \wedge \uelt = \uelt$. 
That is, $\selt \vee (\telt \wedge \selt') < (\selt \vee \telt) \wedge (\selt \vee \selt')$, contrary to our defining distributive law.\hfill\QED

{\bf [\S \RoutineSection.4:\! An initial consequence of diamond-coloring.]} 
The objects of primary interest for us are, of course, diamond-colored modular lattices (DCML's) and diamond-colored distributive lattices (DCDL's). 
Some DCDL's are depicted in \DCDLFigs. 
For a non-distributive DCML, see \CFourFigure.1. 
The next proposition shows how the modular lattice and diamond-coloring properties can interact. 
It is a utility that aids in proofs of other results (e.g.\ \FullLengthTheorem\ and \JCompResult) and helps in analyzing movement between vertices in diamond-colored modular lattices (e.g.\ computing the rank of an element). 

\noindent 
{\bf \ColorsLemma}\ \ {\sl 
Let $L$ be a modular lattice that is diamond-colored by a set $I$.  
Suppose $\selt \leq \telt$.  
Suppose $\mathcal{P} = (\selt = \relt_{0} \myarrow{i_{1}} \relt_{1} \myarrow{i_{2}} \relt_{2} \myarrow{i_{3}} \cdots \myarrow{i_{p-1}} \relt_{p-1} \myarrow{i_{p}} \relt_{p} = \telt)$ and  $\mathcal{Q} = (\selt = \relt'_{0} \myarrow{j_{1}} \relt'_{1} \myarrow{j_{2}} \relt'_{2} \myarrow{j_{3}} \cdots \myarrow{j_{q-1}} \relt'_{q-1} \myarrow{j_{q}} \relt'_{q} = \telt)$ are two paths from $\selt$ up to $\telt$.  
Then, $p = q$ and $a_{i}(\mathcal{P}) = a_{i}(\mathcal{Q})$ for all $i \in I$.  
Moreover, if $\relt_{1}$ and $\relt'_{p-1}$ are incomparable, then $i_{1} = j_{p}$.} 

{\em Proof.} Since $L$ is ranked, then $p = q$. 
We use induction on the length $p$ of the given paths to prove both claims of the lemma statement.  
If $p = 0$, then there is nothing to prove.  
For our induction hypothesis, we assume the theorem statement holds whenever $p \leq m$ for some nonnegative integer $m$.  
Suppose now that $p = m+1$.  
We consider two cases: (1) $\relt_{p-1} = \relt'_{p-1}$ and (2)  $\relt_{p-1} \not= \relt'_{p-1}$. 
In case (1), the induction hypothesis applies to the paths $\selt = \relt_{0} \myarrow{i_{1}} \relt_{1} \myarrow{i_{2}} \relt_{2} \myarrow{i_{3}} \cdots \myarrow{i_{p-1}} \relt_{p-1} = \relt'_{p-1}$ and $\selt = \relt'_{0} \myarrow{j_{1}} \relt'_{1} \myarrow{j_{2}} \relt'_{2} \myarrow{j_{3}} \cdots \myarrow{j_{p-1}} \relt'_{p-1} = \relt_{p-1}$.  
It follows that $\{i_{1}, i_{2}, \ldots, i_{p-1}\} \eqmulti \{j_{1}, j_{2}, \ldots, j_{p-1}\}$.  
Since in this case we have $i_{p} = j_{p}$, we conclude that $\{i_{1}, i_{2}, \ldots, i_{p-1}, i_{p}\} \eqmulti \{j_{1}, j_{2}, \ldots, j_{p-1}, j_{p}\}$.  
Also in this case, $\relt_{1} \leq \relt_{p-1} = \relt'_{p-1}$.  So $\relt_{1}$ and $\relt'_{p-1}$ are comparable.  

In case (2), $\relt_{p-1} \not= \relt'_{p-1}$. 
Let $\xelt := \relt_{p-1} \wedge \relt'_{p-1}$.  
Since $\selt \leq \relt_{p-1}$ and $\selt \leq \relt'_{p-1}$, it follows that $\selt \leq \xelt$.  
Consider a path $\selt = \relt''_{0} \myarrow{k_{1}} \relt''_{1} \myarrow{k_{2}} \relt''_{2} \myarrow{k_{3}} \cdots \myarrow{k_{p-3}} \relt''_{p-3} \myarrow{k_{p-2}} \relt''_{p-2} = \xelt$. 
Since $L$ is a diamond-colored modular lattice, we have $\xelt \myarrow{j_{p}} \relt_{p-1}$ and $\xelt \myarrow{i_{p}} \relt'_{p-1}$.  
Then by the induction hypothesis, we have $\{k_{1}, k_{2}, \ldots, k_{p-2}, j_{p}\} \eqmulti \{i_{1}, i_{2}, \ldots, i_{p-2}, i_{p-1}\}$ and $\{k_{1}, k_{2}, \ldots, k_{p-2}, i_{p}\} \eqmulti \{j_{1}, j_{2}, \ldots, j_{p-2}, j_{p-1}\}$.  
Then, \[\{i_{1}, i_{2}, \ldots, i_{p-2}, i_{p-1}, i_{p}\} \eqmulti \{k_{1}, k_{2}, \ldots, k_{p-2}, i_{p}, j_{p}\} \eqmulti \{j_{1}, j_{2}, \ldots, j_{p-2}, j_{p-1}, j_{p}\},\] as desired.  
Now suppose that $\relt_{1}$ and $\relt'_{p-1}$ are incomparable.  
Suppose $\relt_{1}$ and $\xelt$ are comparable.  
Then it must be the case that $\xelt < \relt_{1}$.  
(Else, $\relt_{1} \leq \xelt$ and $\xelt \leq \relt_{p-1}'$ means $\relt_{1}$ and $\relt_{p-1}'$ are comparable.)  
Since $\relt_{p-1} \not= \relt_{p-1}'$, then $p \geq 2$.  
If $p \geq 3$, then $\rho(\relt_{1}) \leq \rho(\telt) - 2$, while $\rho(\xelt) = \rho(\telt) - 2$.  
This contradicts the fact that $\xelt < \relt_{1}$, so when $\relt_{1}$ and $\xelt$ are comparable, it must be the case that $p = 2$.  
If $p = 2$, then $\xelt = \selt$ and we have the diamond \parbox{2.5cm}{\begin{center}
\setlength{\unitlength}{0.2cm}
\begin{picture}(6.5,3.5)
\put(1,0){\begin{picture}(4,3.5)
\put(2,0){\circle*{0.5}} \put(0,2){\circle*{0.5}}
\put(2,4){\circle*{0.5}} \put(4,2){\circle*{0.5}}
\put(0,2){\line(1,1){2}} \put(2,0){\line(-1,1){2}}
\put(4,2){\line(-1,1){2}} \put(2,0){\line(1,1){2}}
\put(0.75,0.55){\em \tiny $i_{1}$} \put(2.7,0.7){\em \tiny $j_{1}$}
\put(0.7,2.7){\em \tiny $i_{p}$} \put(2.7,2.55){\em \tiny $j_{p}$}
\put(2.5,-0.75){\footnotesize $\selt$} \put(4.75,1.75){\footnotesize 
$\relt_{p-1}'$}
\put(2.5,4){\footnotesize $\telt$} \put(-3.5,1.75){\footnotesize $\relt_{p-1}$}
\end{picture}} \end{picture} \end{center}} in $L$.  
From the diamond coloring property, we conclude that $i_{1} = j_{p}$.  
Suppose now that $\relt_{1}$ and $\xelt$ are incomparable.  
Then we can apply the induction hypothesis to the paths $\selt = \relt''_{0} \myarrow{k_{1}} \relt''_{1} \myarrow{k_{2}} \relt''_{2} \myarrow{k_{3}} \cdots \myarrow{k_{p-3}} \relt''_{p-3} \myarrow{k_{p-2}} \relt''_{p-2} = \xelt \myarrow{j_{p}} \relt_{p-1}$ and $\selt = \relt_{0} \myarrow{i_{1}} \relt_{1} \myarrow{i_{2}} \relt_{2} \myarrow{i_{3}} \cdots \myarrow{i_{p-1}} \relt_{p-1}$.  
From this, we see that $i_{1} = j_{p}$.  
This completes the induction step, and the proof.\hfill\QED  

{\bf [\S \RoutineSection.5:\! A first look at some distinguished features of DCDL's.]} 
The following discussion of diamond-colored distributive lattices and certain related vertex-colored posets encompasses the classical uncolored situation (see for example \S 3.4 of \cite{StanleyText}). 
These concepts have antecedents in the literature in work of Proctor, Stembridge, this author, and Green, among others (see e.g.\ \cite{PrEur}, \cite{PrGZ}, \cite{StemFC}, \cite{StemQuasi}, \cite{DonSupp}, \cite{ADLP}, \cite{ADLMPPW}, \cite{Green}), although there seems to be no standard treatment of these ideas. 

A diamond-colored distributive lattice can be constructed as follows.  
Let $P$ be a poset with vertices colored by a set $I$.  
Let $L$ be the set of lower order ideals (down-sets) from $P$.  
For $\xelt, \yelt \in L$, write $\xelt \leq \yelt$ if and only if $\xelt \subseteq \yelt$ (subset containment).  
With respect to this partial ordering, $L$ is a distributive lattice: $\xelt \vee \yelt = \xelt \cup \yelt$ (set union) and $\xelt \wedge \yelt = \xelt \cap \yelt$ (set intersection) for all $\xelt, \yelt \in L$.  
One can easily see that $\xelt \rightarrow \yelt$ in $L$ if and only if $\xelt \subset \yelt$ (proper containment) and $\yelt \setminus \xelt = \{v\}$ for some maximal element $v$ of $\yelt$ (with $\yelt$ thought of as a subposet of $P$ in the induced order).  
In this case, we declare that $\ecolor_{L}(\xelt \rightarrow \yelt) := \vcolor_{P}(v)$, thus making $L$ an edge-colored distributive lattice.  
One can easily check that $L$ is diamond-colored.  
The diamond-colored distributive lattice just constructed is given special notation: we write $L := \Jcolor(P)$.  
Note that if $P \cong Q$ as vertex-colored posets, then $\Jcolor(P) \cong \Jcolor(Q)$ as edge-colored posets.  
Moreover, $L$ is ranked with rank function given by $\rho(\telt) = |\telt|$, the number of elements in the subset $\telt$ from $P$.  
In particular, the length of $L$ is $|P|$. 

\begin{figure}[t]
\begin{center}
{\bf \OrderIdealFig}\ \   The lattice $L$ from \PosetAndLatticeFig\ recognized as $\Jcolor(P)$.

{\small (In this figure, each down-set from $P$ is identified by the indices of its maximal vertices.} 

\vspace*{-0.05in}
{\small 
For example, $\langle 2,3 \rangle$ in $L$ denotes the down-set $\{v_{2}, v_{3}, v_{4}, v_{5}, v_{6}\}$ from $P$.} 

\vspace*{-0.05in}
{\small 
A join irreducible in $L$ is a down-set $\langle k \rangle$ from $P$ whose only maximal element is $v_{k}$.)}

\setlength{\unitlength}{1cm}
\begin{picture}(3,3.5)
\put(0,3){\begin{picture}(3,3.5)
\put(-0.25,3.5){$P \cong \jcolor(L)$}
\put(4,2){\circle*{0.15}} 
\put(3.4,1.9){\footnotesize $v_{6}$}
\put(4.2,1.9){\footnotesize {\em 2}} 
\put(1,1){\circle*{0.15}}
\put(0.4,0.9){\footnotesize $v_{5}$} 
\put(1.2,0.9){\footnotesize {\em 1}} 
\put(2,2){\circle*{0.15}} 
\put(1.4,1.9){\footnotesize $v_{4}$} 
\put(2.2,1.9){\footnotesize {\em 1}}
\put(3,3){\circle*{0.15}} 
\put(2.4,2.9){\footnotesize $v_{3}$}
\put(3.2,2.9){\footnotesize {\em 1}} 
\put(0,2){\circle*{0.15}}
\put(-0.6,1.9){\footnotesize $v_{2}$} 
\put(0.2,1.9){\footnotesize {\em 2}} 
\put(1,3){\circle*{0.15}} 
\put(0.4,2.9){\footnotesize $v_{1}$} 
\put(1.2,2.9){\footnotesize {\em 2}}
\put(1,1){\line(1,1){2}} \put(0,2){\line(1,1){1}}
\put(1,1){\line(-1,1){1}} \put(2,2){\line(-1,1){1}}
\put(3,3){\line(1,-1){1}}
\end{picture}
}
\end{picture}
\hspace*{1.5in}
\setlength{\unitlength}{1.5cm}
\begin{picture}(4,6.5)
\put(-1,5.5){$L \cong \Jcolor(P)$}
\put(1,0){\line(-1,1){1}}
\put(1,0){\line(1,1){2}}
\put(0,1){\line(1,1){2}}
\put(2,1){\line(-1,1){1}}
\put(2,1){\line(0,1){1}}
\put(1,2){\line(0,1){1}}
\put(2,2){\line(-1,1){2}}
\put(2,2){\line(1,1){2}}
\put(3,2){\line(-1,1){1}}
\put(3,2){\line(0,1){1}}
\put(1,3){\line(1,1){2}}
\put(2,3){\line(0,1){1}}
\put(3,3){\line(-1,1){2}}
\put(0,4){\line(1,1){2}}
\put(4,4){\line(-1,1){2}}
\put(2,6){\VertexIdeals{1,3}{-0.65}}
\put(1,5){\VertexIdeals{2,3}{-0.65}}
\put(3,5){\VertexIdeals{1,6}{0.15}}
\put(0,4){\VertexIdeals{3}{-0.45}}
\put(2,4){\VertexIdeals{2,4,6}{0.15}}
\put(4,4){\VertexIdeals{1}{0.15}}
\put(1,3){\VertexIdeals{4,6}{-0.65}}
\put(2,3){\VertexIdeals{2,6}{0.15}}
\put(3,3){\VertexIdeals{2,4}{0.15}}
\put(1,2){\VertexIdeals{5,6}{-0.65}}
\put(2,2){\VertexIdeals{4}{0.15}}
\put(3,2){\VertexIdeals{2}{0.15}}
\put(0,1){\VertexIdeals{6}{-0.45}}
\put(2,1){\VertexIdeals{5}{0.15}}
\put(1,0){\VertexIdealsEmpty{0.15}}
\put(1,5){\NEEdgeLabelForLatticeI{{\em 2}}}
\put(3,5){\NWEdgeLabelForLatticeI{{\em 1}}}
\put(0,4){\NEEdgeLabelForLatticeI{{\em 2}}}
\put(2,4){\NWEdgeLabelForLatticeI{{\em 1}}}
\put(2,4){\NEEdgeLabelForLatticeI{{\em 2}}}
\put(4,4){\NWEdgeLabelForLatticeI{{\em 2}}}
\put(1,3){\NEEdgeLabelForLatticeI{{\em 2}}}
\put(1,3){\NWEdgeLabelForLatticeI{{\em 1}}}
\put(2,3){\VerticalEdgeLabelForLatticeI{{\em 1}}}
\put(3,3){\NWEdgeLabelForLatticeI{{\em 2}}}
\put(3,3){\NEEdgeLabelForLatticeI{{\em 2}}}
\put(1,2){\VerticalEdgeLabelForLatticeI{{\em 1}}}
\put(1.25,2.25){\NEEdgeLabelForLatticeI{{\em 2}}}
\put(2.2,1.8){\NWEdgeLabelForLatticeI{{\em 2}}}
\put(1.8,1.8){\NEEdgeLabelForLatticeI{{\em 2}}}
\put(3,2){\VerticalEdgeLabelForLatticeI{{\em 1}}}
\put(2.75,2.25){\NWEdgeLabelForLatticeI{{\em 2}}}
\put(0,1){\NEEdgeLabelForLatticeI{{\em 1}}}
\put(2,1){\VerticalEdgeLabelForLatticeI{{\em 1}}}
\put(2,1){\NWEdgeLabelForLatticeI{{\em 2}}}
\put(2,1){\NEEdgeLabelForLatticeI{{\em 2}}}
\put(1,0){\NWEdgeLabelForLatticeI{{\em 2}}}
\put(1,0){\NEEdgeLabelForLatticeI{{\em 1}}}
\end{picture}
\end{center}
\end{figure}

The process described in the previous paragraph can be reversed.  
Given a diamond-colored distributive lattice $L$, an element $\xelt$ is {\em join irreducible} if $\xelt \not= \mymin(L)$ and whenever $\xelt = \yelt \vee \zelt$ then $\xelt = \yelt$ or $\xelt = \zelt$.  
One can see that $\xelt$ is join irreducible if and only if $\xelt$ covers precisely one other vertex in $L$, i.e.\ $\myabs\{\xelt'\, |\, \xelt' \rightarrow \xelt\}\myabs = 1$.  
Let $P$ be the set of all join irreducible elements of $L$ together with the induced partial ordering.  
Color the vertices of the poset $P$ by the rule: $\vcolor_{P}(\xelt) := \ecolor_{L}(\xelt' \rightarrow \xelt)$.  
We call $P$ the {\em vertex-colored poset of join irreducibles} and denote it by $P := \jcolor(L)$.  
If $K \cong L$ is an isomorphism of diamond-colored lattices, then $\jcolor(K) \cong \jcolor(L)$ is an isomorphism of vertex-colored posets.  

What follows is a kind of dual to the above constructions of diamond-colored distributive lattices. 
Note that for $\xelt \subseteq P$, $\xelt$ is an upper order ideal (or up-set) if and only if the set complement $P \setminus \xelt$ is a lower order ideal.  
Partially order all up-sets from $P$ by reverse containment: $\xelt \leq \yelt$ if and only if $\xelt \supseteq \yelt$ for up-sets $\xelt, \yelt$ from $P$.  
The resulting partially ordered set $L$ is a distributive lattice.  
Color the edges of $L$ as in the case of down-sets.  
The result is a diamond-colored distributive lattice which we denote by $L  = \Mcolor(P)$.  
In the other direction, given a diamond-colored distributive lattice $L$, we say $\xelt \in L$ is {\em meet irreducible} if and only if $\xelt \not= \mymax(L)$ and whenever $\xelt = \yelt \wedge \zelt$ then $\xelt = \yelt$ or $\xelt = \zelt$.  
One can see that $\xelt$ is meet irreducible if and only if $\xelt$ is covered by exactly one other vertex in $L$.  
Now consider the set $P$ of meet irreducible elements in $L$ with the order induced from $L$.  
Color the vertices of $P$ in a manner analogous to the way we colored the vertices of the poset of join irreducibles.  
The vertex-colored poset $P$ is the {\em poset of meet irreducibles} for $L$.  
In this case, we write $P = \mcolor(L)$.  
We have $\Mcolor(P) \cong \Mcolor(Q)$ if $P$ and $Q$ are isomorphic  vertex-colored posets, and we have $\mcolor(L) \cong \mcolor(K)$ if $L$ and $K$ are isomorphic  diamond-colored distributive lattices.  

As a consequence of the forthcoming Fundamental Theorem of Finite DCDL's (\FundamentalTheorem, below), we show in \JMCorollary\ that if $L$ is a DCDL, then $\jcolor(L)$ and $\mcolor(L)$ are isomorphic as vertex-colored posets. 
Another consequence of \JMCorollary\ is that $\Jcolor(P) \cong \Mcolor(P)$ are isomorphic diamond-colored distributive lattices for any vertex-colored poset $P$. 
We refer to either of $\jcolor(L)$ and $\mcolor(L)$ as the {\em compression poset of} $L$, and we call either of $\Jcolor(P)$ or $\Mcolor(P)$ the {\em DCDL of order ideals from} $P$ (or just the {\em DCDL from} $P$). 

We modestly extend the discussion of the preceding paragraphs in the following proposition. 
This result recapitulates classical monochromatic aspects of finite distributive lattices to set them in place prior to our forthcoming Fundamental Theorem of Finite DCDL's (\FundamentalTheorem). 

\noindent
{\bf \EncapsulateProposition}\ \ {\sl (1) Let $L$ be a distributive lattice diamond-colored by a set $I$, and let $P := \jcolor(L)$  and $Q := \mcolor(L)$ be the associated $I$-vertex-colored posets of join irreducibles $\!\!$/$\!\!$ meet irreducibles respectively.  
If} $\posetlength(L)=l${\sl , then each of these compression posets has $l$ elements. 
For any $\xelt \in L$, let $\mathcal{I}_{\xelt} := \{\welt \in P\, |\, \welt \leq_{L} \xelt\}$, and for any $\yelt \in L$, let $\mathcal{F}_{\yelt} := \{\zelt \in Q\, |\, \yelt \leq_{L} \zelt\}$. 
Then, for a down-set $\mathcal{I}$ from $P$ and an $\xelt \in L$, we have $\xelt = \bigvee_{\welt \in \mathcal{I}}(\welt)$ if and only if $\mathcal{I} = \mathcal{I}_{\xelt}$. 
Similarly, for an up-set $\mathcal{F}$ from $Q$ and a $\yelt \in L$, we have $\yelt = \bigwedge_{\zelt \in \mathcal{F}}(\zelt)$ if and only if $\mathcal{F} = \mathcal{F}_{\yelt}$.}\\ 
{\sl (2) Let $P$ be an $I$-vertex-colored poset with $l$ elements. 
Then} $\Jcolor(P)$ {\sl and} $\Mcolor(P)$ {\sl are distributive lattices that are diamond-colored by the set $I$.  
Moreover, for any down-set $\xelt$ from $P$, the rank of $\xelt$ in $\Jcolor(P)$ is $|\xelt|$, that is, the size of $\xelt$ as a subset of $P$; similarly, for any up-set $\yelt$ from $P$, the rank of $\yelt$ in $\Mcolor(P)$ is $l-|\yelt|$. 
Finally, a down-set $\mathcal{I}$ from $P$ is a join irreducible in $\Jcolor(P)$ if and only if $\mathcal{I}=(-\infty,x]_{P}$ for some $x \in P$, and an up-set $\mathcal{F}$ from $P$ is a meet irreducible in $\Mcolor(P)$ if and only if $\mathcal{F} = [y,\infty)_{P}$ for some $y \in P$.} 

{\em Proof.} The first claim in part {\sl (1)} of the proposition is a straightforward consequence of the Fundamental Theorem of Finite Distributive Lattices (e.g.\ Theorem 3.4.1 of \cite{StanleyText}). 
Alternatively, this can be deduced directly from our set-up. 
Here we outline how one can show by induction on $l$ that $\myabs\jcolor(L)\myabs = l$ when $L$ is any length $l$ diamond-colored distributive lattice and $l \geq 1$. 
If $l=1$, then any diamond-colored distributive lattice $L$ is a two-element chain and $\mymax(L)$ is its only join irreducible element (i.e.\ $\myabs\jcolor(L)\myabs = 1$), thus establishing the basis step of the induction argument. 
Now suppose for some integer $l \geq 1$ and any integer $k$ with $1 \leq k \leq l$, we have $\myabs\jcolor(K)\myabs = k$ whenever $K$ is a DCDL whose length is $k$. 
Let $L$ be a DCDL with length $l+1$. 
Let $\mathcal{J} := \jcolor(L) = \{\xelt_{1}, \xelt_{2}, \ldots, \xelt_{p}\}$ be the set of all join irreducible elements in $L$ indexed in such a way that $\xelt_{p}$ is maximal when compared to the other join irreducibles. 
By induction on the rank of elements of $L$, it can be seen that for any $\telt \in L$, there is a subset $\mathcal{T}$ of $\mathcal{J}$ such that $\telt = \bigvee_{\mbox{\tiny $\yelt \in \mathcal{T}$}}\yelt$. 
Then, $\melt := \mymax(L) = \bigvee_{\mbox{\tiny $\yelt \in \mathcal{J}$}}\yelt$. 
Next, we make an observation about join irreducible elements of a distributive lattice: If $\xelt \in \mathcal{J}$ and if $\xelt \not\leq \selt$ and $\xelt \not\leq \telt$ for some given elements $\selt$ and $\telt$ in $L$, then $\xelt \not\leq \selt \vee \telt$. 
(To see this, note contrapositively that if $\xelt \leq \selt \vee \telt$, then the facts that $\xelt$ is a join irreducible and that $(\xelt \wedge \selt) \vee (\xelt \wedge \telt) = \xelt \wedge (\selt \vee \telt) = \xelt$ imply that $\xelt \wedge \selt = \xelt$ or $\xelt \wedge \telt = \xelt$, i.e.\ $\xelt \leq \selt$ or $\xelt \leq \telt$.)  
Let $\mathcal{J}' := \mathcal{J} \setminus \{\xelt_{p}\}$ and set $\melt' := \bigvee_{\mbox{\tiny $\yelt \in \mathcal{J}'$}}\yelt$. 
Based on the preceding observation, it is easy to argue that $\xelt_{p} \not\leq \melt'$, so $\melt' < \melt$. 
Suppose now that $\melt' < \telt \leq \melt$, where $\telt = \bigvee_{\mbox{\tiny $\yelt \in \mathcal{T}$}}\yelt$ for some subset $\mathcal{T}$ of $\mathcal{J}$. 
We must have $\xelt_{p} \in \mathcal{T}$, else $\telt \leq \melt'$; but the fact that $\xelt_{q} \leq \melt' < \telt$ for each $q \in \{1,\ldots,p-1\}$ implies that $\melt = \bigvee_{\mbox{\tiny $\yelt \in \mathcal{J}$}}\yelt \leq \telt$, so $\telt = \melt$. 
That is, $\melt' \rightarrow \melt$. 
Now note that $K := \{\yelt \in L\, |\, \yelt \leq \melt'\}$ can be viewed as a diamond-colored distributive lattice of length $l$ whose join irreducible elements are precisely those of the set $\mathcal{J}'$. 
Then $p-1=l$ by our inductive hypothesis. 
So $\myabs\jcolor(L)\myabs = p = l+1$, completing the induction argument. 
We conclude that  $\myabs\jcolor(L)\myabs = \posetlength(L)$. 
The proof that $\myabs\mcolor(L)\myabs = \posetlength(L)$ is entirely similar. 

Now let $\xelt \in L$. 
For the lower order ideal $\mathcal{I}_{\xelt}$ from $P$, it is clear that $\bigvee_{\yelt \in \mathcal{I}_{\xelt}}(\yelt) \leq_{L} \xelt$.  
We claim that $\xelt = \bigvee_{\yelt \in \mathcal{I}_{\xelt}}(\yelt)$.  
To see this we induct on the rank of $\xelt$.  
(That $L$ is ranked is a consequence of \DistributiveIsModularLemma.)  
If $\xelt = \mymin(L)$, then $\mathcal{I}_{\xelt} = \emptyset$, so the desired result follows.  
For our induction hypothesis, we suppose that $\zelt = \bigvee_{\yelt \in \mathcal{I}_{\zelt}}(\yelt)$ for all $\zelt$ with $\rho(\zelt) \leq k$ for some integer $k \geq 0$.  
Suppose now that $\xelt \in L$ with $\rho(\xelt) = k+1$.  
First, consider the case that $\xelt$ is join irreducible.  
Then $\xelt \in \mathcal{I}_{\xelt}$, so the result $\xelt = \bigvee_{\yelt \in \mathcal{I}_{\xelt}}(\yelt)$ follows immediately.  
Now suppose $\xelt$ is not join irreducible.  
Then we may write $\xelt = \selt \vee \telt$ for some $\selt <_{L} \xelt >_{L} \telt$.  
In particular, $\rho(\selt) \leq k$ and $\rho(\telt) \leq k$.  
So the induction hypothesis applies to $\selt$ and $\telt$.  
That is, $\selt = \bigvee_{\yelt \in \mathcal{I}_{\selt}}(\yelt)$ and $\telt = \bigvee_{\yelt \in \mathcal{I}_{\telt}}(\yelt)$.  
Note also that $(\mathcal{I}_{\selt} \cup \mathcal{I}_{\telt}) \subseteq \mathcal{I}_{\xelt}$.  
Then, 
\[\xelt = \selt \vee \telt  =  \bigvee_{\yelt \in (\mathcal{I}_{\selt} \cup \mathcal{I}_{\telt})}(\yelt)  \leq_{L}  \bigvee_{\yelt \in \mathcal{I}_{\xelt}}(\yelt) \leq_{L} \xelt,\]
so we have equality all the way through.  
That is, $\bigvee_{\yelt \in \mathcal{I}_{\xelt}}(\yelt) = \xelt$.  

Next we show that for any $\xelt \in L$, if $\xelt = \bigvee_{\yelt \in \mathcal{I}}(\yelt)$ for some down-set $\mathcal{I}$ from $P$, then $\mathcal{I} = \mathcal{I}_{\xelt}$ and $\myabs\mathcal{I}_{\xelt}\myabs = \rho(\xelt)$.  
To see this, we use induction on the rank of $\xelt$. 
When $\rho(\xelt) = 0$, then $\xelt = \mymin(L)$.  
In this case, if $\xelt = \bigvee_{\yelt \in \mathcal{I}}(\yelt)$ for some down-set $\mathcal{I}$ from $P$, then it must be the case that $\mathcal{I} = \emptyset$, hence $\mathcal{I} = \mathcal{I}_{\xelt}$ and $\myabs\mathcal{I}_{\xelt}\myabs = \rho(\xelt)$.  
For our induction hypothesis, suppose the claim holds for all elements of $L$ with rank no more than $k$ for some positive integer $k$.  
Next suppose that for some $\xelt \in L$ we have $\rho(\xelt) = k+1$ and $\xelt = \bigvee_{\yelt \in \mathcal{I}}(\yelt)$ for some down-set $\mathcal{I}$ from $P$.  
Choose a maximal element $\zelt$ in $\mathcal{I}$.  
Then let $\mathcal{J} := \mathcal{I} \setminus \{\zelt\}$. 
Clearly $\mathcal{J}$ is a down-set from $P$.  Let $\xelt' := \bigvee_{\yelt \in \mathcal{J}}(\yelt)$.  
Clearly $\xelt' \leq_{L} \xelt$.  
In order to apply the induction hypothesis to $\xelt'$, we need $\xelt' <_{L} \xelt$.  
Suppose otherwise, so $\xelt' = \xelt$.  
Then $\xelt' \not= \mymin(L)$, and hence $\mathcal{J} \not= \emptyset$.  
Further, $\xelt' = \xelt = \zelt \vee (\bigvee_{\yelt \in \mathcal{J}}(\yelt)) = \zelt \vee \xelt'$ implies that $\zelt \leq_{L} \xelt'$.  
So $\zelt \wedge \xelt' = \zelt$.  
But then $\zelt = \zelt \wedge \xelt' = \zelt \wedge (\bigvee_{\yelt \in \mathcal{J}}(\yelt)) = \bigvee_{\yelt \in \mathcal{J}}(\zelt \wedge \yelt)$.  
Since $\zelt$ is join irreducible, then $\zelt \wedge \yelt = \zelt$ for some $\yelt \in \mathcal{J}$.  
Then we have $\zelt \leq_{L} \yelt$.  
But $\zelt$ was chosen to be maximal in $\mathcal{I}$, and hence $\zelt \not\leq_{L} \welt$ for all $\welt \in \mathcal{J} = \mathcal{I} \setminus \{\zelt\}$.  
This is a contradiction, so we conclude that $\xelt' <_{L} \xelt$.  
Then $\rho(\xelt') < \rho(\xelt)$, so the induction hypothesis applies to $\xelt'$.  
We get $\mathcal{J} = \mathcal{I}_{\xelt'}$.  
In particular, $\myabs\mathcal{I}\myabs = \myabs\mathcal{I}_{\xelt'}\myabs+1$.  
Applying this reasoning to the particular down-set $\mathcal{I}_{\xelt}$, we conclude that $\myabs\mathcal{I}_{\xelt}\myabs = \myabs\mathcal{I}_{\xelt'}\myabs + 1$. 
Of course if $\telt \in \mathcal{I}$, then by definition $\telt \leq_{L} \xelt$. 
Hence $\telt \in \mathcal{I}_{\xelt}$. 
This shows that $\mathcal{I} \subseteq \mathcal{I}_{\xelt}$. 
Since $\myabs\mathcal{I}\myabs = \myabs\mathcal{I}_{\xelt}\myabs$, we conclude that $\mathcal{I} = \mathcal{I}_{\xelt}$.  
We hope to show that $\xelt$ covers $\xelt'$ in $L$. Suppose otherwise, and say $\xelt' <_{L} \xelt'' <_{L} \xelt$ for some $\xelt'' \in L$.  
Then $\mathcal{I}_{\xelt'} \subset \mathcal{I}_{\xelt''} \subset \mathcal{I}_{\xelt}$, both proper containments.  
(Otherwise, $\xelt' = \bigvee_{\yelt \in \mathcal{I}_{\xelt'}}(\yelt) = \bigvee_{\yelt \in \mathcal{I}_{\xelt''}}(\yelt) = \xelt''$, etc.)   
So $\myabs\mathcal{I}_{\xelt'}\myabs < \myabs\mathcal{I}_{\xelt''}\myabs < \myabs\mathcal{I}_{\xelt}\myabs$.  
But $\myabs\mathcal{I}_{\xelt'}\myabs +1 = \myabs\mathcal{I}_{\xelt}\myabs$, so both of the preceding inequalities cannot be strict.  
We conclude that there is no $\xelt'' \in L$ for which $\xelt' <_{L} \xelt'' <_{L} \xelt$.  
That is, $\xelt$ covers $\xelt'$. 
By the inductive hypothesis we have $\rho(\xelt') = \myabs\mathcal{I}_{\xelt'}\myabs$.  
So $\rho(\xelt) = \rho(\xelt') + 1 = \myabs\mathcal{I}_{\xelt'}\myabs+1 = \myabs\mathcal{I}_{\xelt}\myabs$. 
This completes the proof of our claim. 

Together, the two previous paragraphs have shown that for a down-set $\mathcal{I}$ from $P$ and an $\xelt \in L$, we have $\xelt = \bigvee_{\welt \in \mathcal{I}}(\welt)$ if and only if $\mathcal{I} = \mathcal{I}_{\xelt}$. 
Use analogous arguments to show that for an up-set $\mathcal{F}$ from $Q$ and a $\yelt \in L$, we have $\yelt = \bigwedge_{\zelt \in \mathcal{F}}(\zelt)$ if and only if $\mathcal{F} = \mathcal{F}_{\yelt}$.

For {\sl (2)}, we only consider the claims about $\Jcolor(P)$, since the claims about $\Mcolor(P)$ follow easily from analogous considerations. 
Let $L := \Jcolor(P)$. 
We have already observed in paragraphs prior to the proposition statement that $L$ is diamond-colored by the set $I$ and that its unique rank function is as described. 
Now, for any $x \in P$, the down-set $(-\infty,x]_{P}$ from $P$ has $x$ as its unique maximal element.  
It follows that for a down-set $\mathcal{I}'$ from $P$ we have $\mathcal{I}' \rightarrow (-\infty,x]_{P}$ in $L$ if and only if $\mathcal{I}' = (-\infty,x]_{P} \setminus \{x\}$.  
Hence, $(-\infty,x]_{P}$ is join irreducible in $L$.\hfill\QED

\noindent 
{\bf \BooleanExample}\ \ Let $P$ be an antichain whose elements all have the same color.  
Then the elements of $L := \Jcolor(P)$ are just the subsets of $P$.  
In particular, $|L| = 2^{|P|}$.  
Moreover, the rank $\rho_{L}(\telt)$ of a subset $\telt$ from $P$ is just $|\telt|$, so $|\rho_{L}^{-1}(k)| = \left(\begin{array}{c}|P|\\ k\end{array}\right)$ for any $k \in \{0,1,\ldots,|P|\}$.  
The join irreducible elements of $L$ are just the singleton subsets of $P$.  
Covering relations in $L$ are easy to describe: $\selt \rightarrow \telt$ if and only if $\telt$ is formed from $\selt$ by adding to $\selt$ exactly one  element from $P \setminus \selt$.  
Any such lattice $L$ is called a {\em Boolean lattice}.\hfill\QED 

{\bf [\S \RoutineSection.6:\! A fundamental theorem.]} 
The following theorem shows that the operations $\Jcolor$ and $\jcolor$ are inverses in a certain sense, as are the operations $\Mcolor$ and $\mcolor$.  
The statement of this result is a straightforward generalization of the Fundamental Theorem of Finite Distributive Lattices (e.g.\ Theorem 3.4.1 of \cite{StanleyText}; Theorem 9 from Ch.\ 2, \S 7 of \cite{Gratzer}; Theorem 2.5 of \cite{Aigner}). 
\EncapsulateProposition\ above provides much of what we need to construct the claimed isomorphisms of \FundamentalTheorem. 
Most of our work in the proof of \FundamentalTheorem\ concerns the preservation of edge or vertex colors. 
It is an open question whether the vertex- and edge-coloring notions for DCDL's and their companion compression posets encapsulated in the fundamental theorem below might be extended to DCML's via Markowsky's `poset of irreducibles' for modular lattices (see \cite{Markowsky}).

\noindent 
{\bf \FundamentalTheorem\ (Fundamental Theorem of Finite Diamond-colored Distributive Lattices)}\ \ 
{\sl (1) Let  $L$ be a diamond-colored distributive lattice. Then}   
\[L \cong \Jcolor(\jcolor(L)) \cong \Mcolor(\mcolor(L)),\]
{\sl with edge-preserving and edge-color-preserving isomorphisms $\Theta: L \longrightarrow \Jcolor(\jcolor(L))$ and $\Psi: L \longrightarrow \Mcolor(\mcolor(L))$ given by $\xelt \stackrel{\Theta}{\longmapsto} \mathcal{I}_{\xelt}$ and  $\xelt \stackrel{\Psi}{\longmapsto} \mathcal{F}_{\xelt}$ respectively.}\\ 
{\sl (2) Let $P$ be a  vertex-colored poset.  Then} 
\[P \cong \jcolor(\Jcolor(P)) \cong \mcolor(\Mcolor(P)),\] 
{\sl with edge-preserving and vertex-color-preserving isomorphisms $\theta: P \longrightarrow \jcolor(\Jcolor(P))$ and $\psi: P \longrightarrow \mcolor(\Mcolor(P))$ given by $x \stackrel{\theta}{\longmapsto} (-\infty,x]_{P}$ and  $x \stackrel{\psi}{\longmapsto} [x,\infty)_{P}$ respectively.}

{\em Proof.}  
For {\sl (1)}, set $P := \jcolor(L)$. 
Consider the function $\Theta: L \longrightarrow \Jcolor(P)$ defined by $\Theta(\xelt) := \mathcal{I}_{\xelt}$.  
We show that $\Theta$ is a bijection.  
If $\mathcal{I}_{\selt} = \mathcal{I}_{\telt}$, then, by \EncapsulateProposition.1, $\selt = \bigvee_{\yelt \in \mathcal{I}_{\selt}}(\yelt) = \bigvee_{\yelt \in \mathcal{I}_{\telt}}(\yelt) = \telt$.  
In particular, $\Theta$ is injective. 
Now suppose $\mathcal{I}$ is a down-set from $P$.  
Let $\xelt := \bigvee_{\yelt \in \mathcal{I}}(\yelt)$.  
By \EncapsulateProposition.1, $\mathcal{I} = \mathcal{I}_{\xelt}$.  
So, $\Theta$ is surjective. 

We wish to show that $\selt \myarrow{i} \telt$ in $L$ if and only if $\mathcal{I}_{\selt} \myarrow{i} \mathcal{I}_{\telt}$ in $\Jcolor(P)$.  
First, suppose $\selt \myarrow{i} \telt$ in $L$.  
It follows from the definitions that $\mathcal{I}_{s} \subseteq \mathcal{I}_{\telt}$.  
Now $\selt \not= \telt$ since $\telt$ covers $\selt$ in $L$.  
Since $\mathcal{I}_{\selt} = \Theta(\selt)$ and $\mathcal{I}_{\telt} = \Theta(\telt)$ and $\Theta$ is injective, then $\mathcal{I}_{\selt} \not= \mathcal{I}_{\telt}$.  
So $\mathcal{I}_{s} \subset \mathcal{I}_{\telt}$ is a proper containment.  
Suppose $\mathcal{I}_{\selt} \subseteq \mathcal{I} \subseteq \mathcal{I}_{\telt}$.  
Since $\Theta$ is surjective, then $\mathcal{I} = \mathcal{I}_{\xelt}$ for some $\xelt \in L$.  
But then $\bigvee_{\yelt \in \mathcal{I}_{\selt}}(\yelt) \leq_{L} \bigvee_{\yelt \in \mathcal{I}_{\xelt}}(\yelt) \leq_{L} \bigvee_{\yelt \in \mathcal{I}_{\telt}}(\yelt)$, and hence $\selt \leq_{L} \xelt \leq_{L} \telt$.  
Since $\telt$ covers $\selt$,  then $\selt = \xelt$ or $\xelt = \telt$.  
Hence $\mathcal{I}_{\selt} \rightarrow \mathcal{I}_{\telt}$ in $\Jcolor(P)$.  
In particular, there is some $\zelt \in P$ such that $\mathcal{I}_{\selt} = \mathcal{I}_{\telt} \setminus \{\zelt\}$.  
Moreover, by the definition of $\Jcolor$, $\mathcal{I}_{\selt} \myarrow{j} \mathcal{I}_{\telt}$ in $\Jcolor(P)$ where $j = \vcolor_{P}(\zelt)$.  
Now $j$ is just the color of the edge $\zelt' \myarrow{j} \zelt$ for the unique descendant $\zelt'$ of $\zelt$ in $L$.  
If $\zelt = \telt$, then necessarily $\zelt' = \selt$, and so $j = i$.  

So now suppose that $\zelt \not= \telt$.  So we have $\zelt <_{L} \telt$, and hence $\mathcal{I}_{\zelt} \subset \mathcal{I}_{\telt}$. 
We claim that $\zelt' \leq_{L} \selt$.  
To see this, apply the reasoning of the preceding paragraph to conclude that $\mathcal{I}_{\zelt'} \subset \mathcal{I}_{\zelt}$ with $\mathcal{I}_{\zelt} = \mathcal{I}_{\zelt'} \cup \{\zelt\}$.  
It follows that $\mathcal{I}_{\zelt'} \subset \mathcal{I}_{\telt}$.  
Since $\zelt \not\in \mathcal{I}_{\zelt'}$, $\mathcal{I}_{\zelt} \subset \mathcal{I}_{\telt}$, and $\mathcal{I}_{\selt} = \mathcal{I}_{\telt} \setminus \{\zelt\}$, we get $\mathcal{I}_{\zelt'} \subseteq \mathcal{I}_{\selt}$.  
Then $\zelt' = \bigvee_{\yelt \in \mathcal{I}_{\zelt'}}(\yelt) \leq_{L} \bigvee_{\yelt \in \mathcal{I}_{\selt}}(\yelt) = \selt$.  
Since $\zelt' \leq_{L} \selt$, there is a path $\zelt' = \zelt_{0}' \myarrow{i_{1}} \zelt_{1}' \myarrow{i_{2}} \cdots \myarrow{i_{p}} \zelt_{p}' = \selt$ in $L$ from $\zelt'$ up to $\selt$.  
Since $\zelt' \myarrow{j} \zelt$ and $\zelt' \myarrow{i_{1}} \zelt'_{1}$ and since $L$ is topographically balanced, then there is a unique $\zelt_{1}$ such that $\zelt \rightarrow \zelt_{1}$ and $\zelt'_{1} \rightarrow \zelt_{1}$.  
Since $L$ is diamond-colored, then $\zelt \myarrow{i_{1}} \zelt_{1}$ and $\zelt'_{1} \myarrow{j} \zelt_{1}$.  
Continue in this way, eventually obtaining a path $\zelt = \zelt_{0} \myarrow{i_{1}} \zelt_{1} \myarrow{i_{2}} \cdots \myarrow{i_{p}} \zelt_{p}$ with $\zelt'_{q} \myarrow{j} \zelt_{q}$ for $0 \leq q \leq p$. 
In particular, $\selt \leq_{L} \zelt_{p}$ and $\zelt \leq_{L} \zelt_{p}$, so $\selt \vee \zelt \leq_{L} \zelt_{p}$. 
We claim that $\zelt$ and $\selt$ are not comparable. 
Otherwise, $\selt \leq_{L} \zelt$ or $\zelt \leq_{L} \selt$.  
In the latter case, we would have $\zelt \in \mathcal{I}_{\selt}$, which is not true.  
In the former case, $\selt \leq_{L} \zelt <_{L} \telt$.  
Since $\telt$ covers $\selt$, then we must have $\selt = \zelt$.  
But then $\zelt \in \mathcal{I}_{\selt}$, which is not true.  
Since $\selt$ and $\zelt$ are not comparable, then $\selt < \selt \vee \zelt$.  
Since $\selt \vee \zelt \leq_{L} \zelt_{p}$ and $\selt \myarrow{j} \zelt_{p}$, it follows that $\zelt_{p} = \selt \vee \zelt$.  
But $\selt \vee \zelt = (\bigvee_{\yelt \in \mathcal{I}_{\selt}}(\yelt)) \vee \zelt = \telt$, and hence $\telt \leq \zelt_{p}$.  
Since $\selt$ is covered by both $\zelt_{p}$ and $\telt$, this can only mean that $\zelt_{p} = \telt$.  
Hence $j = i$.  
So $\mathcal{I}_{\selt} \myarrow{i} \mathcal{I}_{\telt}$ in $\Jcolor(P)$. 

On the other hand, suppose $\mathcal{I}_{\selt} \myarrow{i} \mathcal{I}_{\telt}$ in $\Jcolor(P)$.  
Then $\mathcal{I}_{\selt} = \mathcal{I}_{\telt} \setminus \{\zelt\}$ for some $\zelt \in P$, where $i = \vcolor(\zelt)$.  
That is, $\zelt' \myarrow{i} \zelt$ in $L$, where $\zelt'$ is the unique descendant of $\zelt$ in $L$.  
Then $\selt <_{L} \telt$.  
Since $\rho(\selt) = \myabs\mathcal{I}_{\selt}\myabs$ and $\rho(\telt) = \myabs\mathcal{I}_{\telt}\myabs$, then $\selt \rightarrow \telt$.  
Let $j$ be the color of this edge, so $\selt \myarrow{j} \telt$.  
The preceding two paragraphs showed that we must have $\mathcal{I}_{\selt} \myarrow{j} \mathcal{I}_{\telt}$.  
Then $i = j$. 

We conclude that $\Theta$ is an edge- and edge-color preserving bijection from $L$ to $\Jcolor(\jcolor(P))$.  
The argument that $\Psi: L \longrightarrow \Mcolor(\mcolor(P))$ is an isomorphism is entirely similar.  
This completes the proof of {\sl (1)}.

\begin{figure}[t]
\begin{center}
{\bf \NewOrderIdealFig}\ \   An illustration of the principles that  $\Jcolor(P_{1} \oplus P_{2}) \cong \Jcolor(P_{1}) \times \Jcolor(P_{2})$ 

\vspace*{-0.05in}
and $\jcolor(L_{1} \times L_{2}) \cong \jcolor(L_{1}) \oplus \jcolor(L_{2})$, cf.\ \JMCorollary.

{\small 
(As in \OrderIdealFig, here each down-set from $Q$ is identified by the indices of its maximal vertices.} 

\vspace*{-0.05in}
{\small 
A join irreducible in $K$ is a down-set $\langle k \rangle$ from $Q$ whose only maximal element is $v_{k}$.)}

\setlength{\unitlength}{1cm}
\begin{picture}(3,3.5)
\put(0,3){\begin{picture}(3,3.5)
\put(0.1,3.65){$Q \cong \jcolor(K)$}
\put(1.95,2){$\bigoplus$}
\put(4,1.75){\circle*{0.15}} 
\put(3.4,1.65){\footnotesize $v_{6}$}
\put(4.2,1.65){\footnotesize {\em 2}} 
\put(0,1){\circle*{0.15}}
\put(-0.6,0.9){\footnotesize $v_{5}$} 
\put(0.2,0.9){\footnotesize {\em 1}} 
\put(1,2){\circle*{0.15}} 
\put(0.4,1.9){\footnotesize $v_{4}$} 
\put(1.2,1.9){\footnotesize {\em 1}}
\put(3,2.75){\circle*{0.15}} 
\put(2.4,2.65){\footnotesize $v_{3}$}
\put(3.2,2.65){\footnotesize {\em 1}} 
\put(-1,2){\circle*{0.15}}
\put(-1.6,1.9){\footnotesize $v_{2}$} 
\put(-0.8,1.9){\footnotesize {\em 2}} 
\put(0,3){\circle*{0.15}} 
\put(-0.6,2.9){\footnotesize $v_{1}$} 
\put(0.2,2.9){\footnotesize {\em 2}}
\put(0,1){\line(1,1){1}} \put(-1,2){\line(1,1){1}}
\put(0,1){\line(-1,1){1}} \put(1,2){\line(-1,1){1}}
\put(3,2.75){\line(1,-1){1}}
\end{picture}
}
\end{picture}
\hspace*{1.5in}
\setlength{\unitlength}{1.5cm}
\begin{picture}(4,6.5)
\put(-1,5.5){$K \cong \Jcolor(Q)$}
\put(1,0){\line(-1,1){2}}
\put(1,0){\line(1,1){2}}
\put(0,1){\line(1,1){2}}
\put(2,1){\line(-1,1){2}}
\put(2,1){\line(0,1){1}}
\put(1,2){\line(0,1){1}}
\put(0,3){\line(0,1){1}}
\put(1,4){\line(0,1){1}}
\put(2,2){\line(-1,1){2}}
\put(2,2){\line(1,1){2}}
\put(3,2){\line(-1,1){2}}
\put(3,2){\line(0,1){1}}
\put(1,3){\line(1,1){2}}
\put(2,3){\line(0,1){1}}
\put(3,3){\line(-1,1){2}}
\put(0,4){\line(1,1){2}}
\put(4,4){\line(-1,1){2}}
\put(-1,2){\line(1,1){2}}
\put(2,6){\VertexIdeals{1,3}{-0.65}}
\put(1,5){\VertexIdeals{2,3,4}{-0.85}}
\put(1,4){\VertexIdeals{2,3}{-0.65}}
\put(3,5){\VertexIdeals{1,6}{0.15}}
\put(0,4){\VertexIdeals{3,4}{-0.65}}
\put(0,3){\VertexIdeals{3,5}{-0.65}}
\put(-1,2){\VertexIdeals{3}{-0.45}}
\put(2,4){\VertexIdeals{2,4,6}{0.15}}
\put(4,4){\VertexIdeals{1}{0.15}}
\put(1,3){\VertexIdeals{4,6}{-0.65}}
\put(2,3){\VertexIdeals{2,6}{0.15}}
\put(3,3){\VertexIdeals{2,4}{0.15}}
\put(1,2){\VertexIdeals{5,6}{-0.65}}
\put(2,2){\VertexIdeals{4}{0.15}}
\put(3,2){\VertexIdeals{2}{0.15}}
\put(0,1){\VertexIdeals{6}{-0.45}}
\put(2,1){\VertexIdeals{5}{0.15}}
\put(1,0){\VertexIdealsEmpty{0.15}}
\put(1,5){\NEEdgeLabelForLatticeI{{\em 2}}}
\put(3,5){\NWEdgeLabelForLatticeI{{\em 1}}}
\put(0,4){\NEEdgeLabelForLatticeI{{\em 2}}}
\put(2,4){\NWEdgeLabelForLatticeI{{\em 1}}}
\put(2,4){\NEEdgeLabelForLatticeI{{\em 2}}}
\put(4,4){\NWEdgeLabelForLatticeI{{\em 2}}}
\put(0.8,2.8){\NEEdgeLabelForLatticeI{{\em 2}}}
\put(0.3,3.3){\NEEdgeLabelForLatticeI{{\em 2}}}
\put(1.2,2.8){\NWEdgeLabelForLatticeI{{\em 1}}}
\put(1.7,3.3){\NWEdgeLabelForLatticeI{{\em 1}}}
\put(2,3){\VerticalEdgeLabelForLatticeI{{\em 1}}}
\put(3,3){\NWEdgeLabelForLatticeI{{\em 2}}}
\put(3,3){\NEEdgeLabelForLatticeI{{\em 2}}}
\put(1,2){\VerticalEdgeLabelForLatticeI{{\em 1}}}
\put(1.25,2.25){\NEEdgeLabelForLatticeI{{\em 2}}}
\put(2.2,1.8){\NWEdgeLabelForLatticeI{{\em 2}}}
\put(1.8,1.8){\NEEdgeLabelForLatticeI{{\em 2}}}
\put(3,2){\VerticalEdgeLabelForLatticeI{{\em 1}}}
\put(2.75,2.25){\NWEdgeLabelForLatticeI{{\em 2}}}
\put(0,1){\NEEdgeLabelForLatticeI{{\em 1}}}
\put(2,1){\VerticalEdgeLabelForLatticeI{{\em 1}}}
\put(2,1){\NWEdgeLabelForLatticeI{{\em 2}}}
\put(2,1){\NEEdgeLabelForLatticeI{{\em 2}}}
\put(1,0){\NWEdgeLabelForLatticeI{{\em 2}}}
\put(1,0){\NEEdgeLabelForLatticeI{{\em 1}}}
\put(-1,2){\NEEdgeLabelForLatticeI{{\em 1}}}
\put(0,1){\NWEdgeLabelForLatticeI{{\em 1}}}
\put(1,2){\NWEdgeLabelForLatticeI{{\em 1}}}
\put(0,3){\VerticalEdgeLabelForLatticeI{{\em 1}}}
\put(1,4){\VerticalEdgeLabelForLatticeI{{\em 1}}}
\end{picture}
\end{center}
\end{figure}

For {\sl (2)}, we only show $P \cong \jcolor(\Jcolor(P))$ since the argument that $P \cong  \mcolor(\Mcolor(P))$ is entirely similar. 
Let $L := \Jcolor(P)$, and let $Q := \jcolor(L)$, where we view the latter as the special collection of down-sets of $P$ that are join irreducibles in $L$.  
Consider the mapping $\theta: P \rightarrow Q$ by $\theta(x) := (-\infty,x]_{P}$. 
We know that $(-\infty,x]_{P}$ is a join irreducible element of $L$ by \EncapsulateProposition.2.
We claim that $\theta$ is a bijection.   
Indeed, if $\theta(x) = \theta(y)$ for $x, y \in P$, then $(-\infty,x]_{P} = (-\infty,y]_{P}$.  
But then $x \leq_{P} y$ and $y \leq_{P} x$.  
Therefore $x = y$, and hence $\theta$ is injective.  
On the other hand, if $\mathcal{I}$ is a down-set from $P$ that is join irreducible in $L$, then $\mathcal{I}$ must have a unique maximal element, say $x$.  
But then $\mathcal{I} = (-\infty,x]_{P} = \theta(x)$, so $\theta$ is injective. 

Finally we show $\theta$ preserves edges and vertex colors.  
If $u \rightarrow v$ in $P$, then $\langle u \rangle <_{Q} \langle v \rangle$.  
Now if $\langle u \rangle \leq_{Q} \langle z \rangle \leq_{Q} \langle v \rangle$, it follows that $u \leq_{P} z \leq_{P} v$.  
Since $v$ covers $u$ in $P$, then $u = z$ or $z = v$, and hence $\langle u \rangle = \langle z \rangle$ or $\langle z \rangle = \langle v \rangle$.  
That is, $u \rightarrow v$ in $P$ implies that $\theta(u) \rightarrow \theta(v)$ in $Q$.  
Conversely, if $\langle u \rangle \rightarrow \langle v \rangle$ in $Q$, then we must have $u <_{P} v$ in $P$.  
Suppose $u \leq_{P} z \leq_{P} v$.  
Then one easily sees that $\langle u \rangle \leq_{Q} \langle z \rangle \leq_{P} \langle v \rangle$, and hence $\langle u \rangle = \langle z \rangle$ or $\langle z \rangle = \langle v \rangle$. 
Then $u = z$ or $z = v$.  
That is, $\theta(u) \rightarrow \theta(v)$ in $Q$ implies that $u \rightarrow v$ in $P$.  
As for vertex colors, observe that $\vcolor_{P}(v) = i$ if and only if $\langle v \rangle \setminus \{v\} \myarrow{i} \langle v \rangle$ in $L$ if and only if $\vcolor_{Q}(\theta(v)) = i$. 
This completes the proof.\hfill\QED

\noindent 
{\bf \FirstCorollary}\ \ {\sl An edge-colored distributive lattice $L$ is isomorphic to $\Jcolor(P)$ or $\Mcolor(P)$ for some vertex-colored poset $P$ if and only if $L$ is diamond-colored.}   

{\em Proof.} The ``only if'' direction was observed in the  paragraphs preceding \EncapsulateProposition.  
For the ``if'' direction, we get $L \cong \Jcolor(P) \cong \Mcolor(Q)$ from \FundamentalTheorem, where $P := \jcolor(L)$ and $Q := \mcolor(L)$.\hfill\QED 

The next corollary states for the record how $\Jcolor$, $\jcolor$, $\Mcolor$, and $\mcolor$ interact with the vertex- and edge-colored poset operations  $*$ (dual), $\sigma$ (recoloring of vertices or edges), $\oplus$ (disjoint union), and $\times$ (Cartesian product). 

\noindent 
{\bf \JMCorollary}\ \ {\sl Let $P$ and $Q$ be posets with vertices colored by a set $I$, and let $K$ and $L$ be diamond-colored distributive lattices with edges colored by $I$.  
In what follows, $*$, $\sigma$, $\oplus$, $\times$, and $\cong$ account for colors on vertices/edges as appropriate. 
(1) If $K \cong L$, then $\jcolor(K) \cong \mcolor(K) \cong \mcolor(L) \cong \jcolor(L)$.  
If $P \cong Q$, then $\Jcolor(P) \cong \Mcolor(P) \cong \Mcolor(Q) \cong \Jcolor(Q)$.
(2) Also,  $\Jcolor(P^{*}) \cong  (\Jcolor(P))^{*}$, $\Jcolor(P^{\sigma}) \cong  (\Jcolor(P))^{\sigma}$ (recoloring), and $\Jcolor(P \oplus Q) \cong \Jcolor(P) \times \Jcolor(Q)$.  
Moreover, $\Mcolor(P^{*}) \cong  (\Mcolor(P))^{*}$, $\Mcolor(P^{\sigma}) \cong (\Mcolor(P))^{\sigma}$, and $\Mcolor(P \oplus Q) \cong \Mcolor(P) \times \Mcolor(Q)$. 
(3) Similarly, $\jcolor(L^{*}) \cong  (\jcolor(L))^{*}$, $\jcolor(L^{\sigma}) \cong  (\jcolor(L))^{\sigma}$, and $\jcolor(L \times K) \cong \jcolor(L) \oplus \jcolor(K)$. 
Moreover, $\mcolor(L^{*}) \cong  (\mcolor(L))^{*}$, $\mcolor(L^{\sigma}) \cong  (\mcolor(L))^{\sigma}$, and $\mcolor(L \times K) \cong \mcolor(L) \oplus \mcolor(K)$.} 

{\em Proof.}  Proofs of the claims in part {\sl (2)} are routine and therefore omitted.  
Apply \FundamentalTheorem\ to deduce {\sl (3)} from {\sl (2)}.  
For {\sl (1)}, when $P \cong Q$, the fact that $\Jcolor(P) \cong \Mcolor(Q)$ follows from the definitions.  
Suppose $K \cong L$.  From part {\sl (3)}, it follows that  $\jcolor(K) \cong (\jcolor(K^{*}))^{*}$. 
From the definitions, we get $(\jcolor(K^{*}))^{*} \cong \mcolor(K)$.  
Then $\jcolor(K) \cong \mcolor(L)$.\hfill\QED

To close this section, we record a routine result concerning compression posets. 
We will invoke this monochromatic result in later sections when we study minuscule splitting DCDL's and minuscule compression posets. 
The proof follows from first principles and is left as an exercise for the reader. 

\noindent
{\bf \CompressionPosetLattice}\ \ {\sl Let $P$ be a poset with lattice of order ideals $L := \mathbf{J}(P)$. 
Then $P$ is a lattice if and only if the join of any two meet irreducible elements of $L$ is a meet irreducible of $L$ and likewise the meet of any two join irreducible elements of $L$ is join irreducible.}

\begin{center}
\underline{\hspace*{4in}}
\end{center}

\vspace*{0.5cm} 
\noindent
{\bf \S \SubstructureSection. Substructures.} 
The results of this section concern certain notions of substructures of posets and lattices, especially DCML's and DCDL's. 

{\bf [\S \SubstructureSection.1:\! Types of subposets.]} 
Given a subset $Q$ of a poset $R$, let $Q$ inherit the partial ordering of $R$; call $Q$ a {\em subposet in the induced order}.  
For posets $(R,\leq_{R})$ and $(Q,\leq_{Q})$, suppose $Q \subseteq R$ and  $\selt \leq_{Q} \telt \Rightarrow \selt \leq_{R} \telt$ for all $\selt, \telt \in Q$.  
Then $Q$ is a {\em weak subposet} of $R$.  
If, in addition, $Q$ and $R$ are vertex-colored (respectively, edge-colored) by a set $I$ and $\vcolor_{Q}^{-1}(i) \subseteq \vcolor_{R}^{-1}(i)$ (resp.\ $\ecolor_{Q}^{-1}(i) \subseteq \ecolor_{R}^{-1}(i)$) for all $i \in I$, then $Q$ is a {\em vertex-colored} (resp.\ {\em edge-colored}) {\em weak subposet}. 

{\bf [\S \SubstructureSection.2:\! Sublattices.]} 
Let $L$ be a lattice with partial ordering $\leq_{L}$ and meet and join operations $\wedge_{L}$ and $\vee_{L}$ respectively.  
Let $K$ be a vertex subset of $L$.  
Suppose that $K$ has a lattice partial ordering $\leq_{K}$ of its own with meet and join operations $\wedge_{K}$ and $\vee_{K}$ respectively. 
We say $K$ is a {\em sublattice} of $L$ if for all $\xelt$ and $\yelt$ in $K$ we have $\xelt \wedge_{K} \yelt = \xelt \wedge_{L} \yelt$ and $\xelt \vee_{K} \yelt = \xelt \vee_{L} \yelt$.  
It is easy to see that if $K$ is a sublattice of $L$ then for all $\xelt$ and $\yelt$ in $K$ we have $\xelt \leq_{K} \yelt$ if and only if $\xelt \leq_{L} \yelt$.  
That is, $K$ is a subposet of $L$ in the induced order. 
If, in addition, $K$ and $L$ are edge-colored and $K$ is an edge-colored weak subposet of $L$, then call $K$ an {\em edge-colored sublattice} of $L$. 
As an example, suppose that a lattice $L$ is $I$-edge-colored and that $\xelt \leq \yelt$ in $L$. 
It is easy to see that the interval $[\xelt,\yelt]$ is an edge-colored sublattice of $L$. 

Some strictly monochromatic results are apropos here. 
Famously, R.\ Dedekind gave necessary and sufficient conditions for a lattice to be modular based on the absence of a certain sublattice. 
G.\ Birkhoff later gave similar necessary and sufficient conditions for a lattice to be distributive. 
To state these conditions, declare the `$\myN_{5}$' lattice to be \parbox[c]{1cm}{\begin{center}
\setlength{\unitlength}{0.2cm}
\begin{picture}(4.5,4.5)
\put(2.2,0){\circle*{0.6}} \put(0.2,2){\circle*{0.6}} \put(0.2,4){\circle*{0.6}}
\put(2.2,6){\circle*{0.6}} \put(4.2,3){\circle*{0.6}}
\put(0.2,2){\line(0,1){2}} \put(2.2,0){\line(-1,1){2}} \put(0.2,4){\line(1,1){2}}
\put(4.2,3){\line(-2,3){2}} \put(2.2,0){\line(2,3){2}}
\end{picture} \vspace*{-0.3in}  
\end{center}}, 
and define the `$\myM_{5}$' lattice to be \parbox[c]{1cm}{\begin{center}
\setlength{\unitlength}{0.2cm}
\begin{picture}(4.5,4.5)
\put(2.2,0){\circle*{0.6}} \put(0.2,3){\circle*{0.6}} \put(2.2,3){\circle*{0.6}}
\put(4.2,3){\circle*{0.6}} \put(2.2,6){\circle*{0.6}}
\put(0.2,3){\line(2,3){2}} \put(2.2,0){\line(-2,3){2}} \put(4.2,3){\line(-2,3){2}}
\put(2.2,0){\line(0,1){3}} \put(2.2,3){\line(0,1){3}} \put(2.2,0){\line(2,3){2}}
\end{picture} \vspace*{-0.3in} 
\end{center}}. 

\noindent 
{\bf \MfiveNfiveTheorem\ (Dedekind, Birkhoff)}\ \ {\sl A lattice is modular if and only if it has no sublattice isomorphic to} $\myN_{5}$.  {\sl A lattice is distributive if and only if it has no sublattice isomorphic either to} $\myN_{5}$ {\sl or} $\myM_{5}$. 

{\em Proof.} See \cite{Aigner}, Propositions 2.20 and 2.21.\hfill\QED

The preceding result has some further consequences for the structure of modular and distributive lattices. 
Say a lattice $L$ is {\em cancellative} if, for all $\xelt$, $\yelt$, $\zelt$ in $L$, we have $\xelt = \yelt$ whenever $\xelt \vee \zelt = \yelt \vee \zelt$ and $\xelt \wedge \zelt = \yelt \wedge \zelt$. 
Say $L$ is {\em semi-cancellative} if, for all $\xelt$, $\yelt$, $\zelt$ in $L$, we have $\xelt = \yelt$ whenever $\xelt \vee \zelt = \yelt \vee \zelt$, $\xelt \wedge \zelt = \yelt \wedge \zelt$, and $\xelt \leq \yelt$. 

\noindent 
{\bf \CancellationTheorem}\ \ {\sl A lattice is semi-cancellative if and only if it is modular, and it is cancellative if and only if it is distributive.}

{\em Proof.} See \cite{Aigner}, Corollary 2.22.\hfill\QED

The next (standard) result can be deduced from either of \MfiveNfiveTheorem\ or \CancellationTheorem, but in the proof below we appeal to \AbstractModularDistributiveProperties. 

\noindent 
{\bf \ModularDistributiveSublattices}\ \ {\sl (1) Any sublattice of a modular lattice is modular. (2) Any sublattice of a distributive lattice is distributive.}

{\em Proof.} Since a given sublattice of a modular (respectively, distributive) lattice inherits the properties expressed in \ModularLawLemma\ (resp.\ \DistributiveIsModularLemma.1), then the sublattice must itself be modular (resp.\ distributive).\hfill\QED

{\bf [\S \SubstructureSection.3:\! Full-length sublattices.]} 
Returning to our possibly-edge-colored context, if $K$ is a sublattice of $L$, if both $K$ and $L$ are ranked, and if both have the same length, then say $K$ is a {\em full-length sublattice} of $L$. 

\noindent 
{\bf \FullLengthLemma}\ \ {\sl Let $K$ be a full-length sublattice of $L$.  
Let $\rho^{(K)}$ and $\rho^{(L)}$ denote the rank functions of $K$ and $L$ respectively.  
Then $\rho^{(K)}(\xelt) = \rho^{(L)}(\xelt)$ for all $\xelt$ in $K$, and moreover for all $\xelt$ and $\yelt$ in $K$ we have $\xelt \rightarrow \yelt$ in $K$ if and only if $\xelt \rightarrow \yelt$ in $L$.} 

{\em Proof.} Let $l$ denote the common length of the ranked posets $K$ and $L$, and take any $\xelt$ in $K$.  
Then $\xelt = \xelt_{r}$ in some longest chain $\xelt_{0} \rightarrow \xelt_{1} \rightarrow \cdots \rightarrow \xelt_{l}$ in $K$.  
Now $\big(\rho^{(K)}(\xelt_{0}), \rho^{(K)}(\xelt_{1}), \ldots, \rho^{(K)}(\xelt_{l})\big) = (0, 1, \ldots, l)$.  
Since $\big(\rho^{(L)}(\xelt_{0}), \rho^{(L)}(\xelt_{1}), \ldots, \rho^{(L)}(\xelt_{l})\big)$ is an increasing sequence of integers bounded below by $0$ and above by $l$, then we get the equality $\big(\rho^{(L)}(\xelt_{0}), \rho^{(L)}(\xelt_{1}), \ldots, \rho^{(L)}(\xelt_{l})\big) = (0, 1, \ldots, l)$.  
Hence $\rho^{(K)}(\xelt) = \rho^{(K)}(\xelt_{r}) = \rho^{(L)}(\xelt_{r}) = \rho^{(L)}(\xelt)$. 

Finally, let $\xelt$ and $\yelt$ be elements of $K$.  
Assume $\xelt \rightarrow \yelt$ in $K$.  
Then $x <_{K} \yelt$ and $\rho^{(K)}(\xelt) + 1 = \rho^{(K)}(\yelt)$.  
So $x <_{L} \yelt$ in $L$ and $\rho^{(L)}(\xelt) + 1 = \rho^{(L)}(\yelt)$.  
Hence $\xelt \rightarrow \yelt$ is a covering relation in $L$ as well. 
Clearly this argument reverses to show that if $\xelt \rightarrow \yelt$ in $L$ then $\xelt \rightarrow \yelt$ in $K$.\hfill\QED 

The previous lemma gives us one way to know whether the edges of a sublattice are also edges of the `parent' lattice. 
Next is a situation in which a full-length sublattice can easily be discerned. 
This result is an elementary observation, and we claim no priority. 

\noindent
{\bf \FullLengthWithinProduct}\ \ {\sl Suppose $L_{1}, L_{2}, \ldots , L_{p}$ are all modular (respectively, distributive) lattices that are diamond-colored by a set $I$, with respective rank functions $\rho^{(1)}, \rho^{(2)}, \ldots , \rho^{(p)}$ and with respective lengths $l^{(1)}, l^{(2)}, \ldots , l^{(p)}$.  Let $L := L_{1} \times L_{2} \times \cdots \times L_{p}$, the Cartesian product poset (in the component-wise order on $p$-tuples) that is edge-colored by $I$ and with rank function $\rho$.\\  
(1) Then $L$ is also a modular (resp.\ distributive) lattice, is diamond-colored by $I$, and has length given by $\sum_{q=1}^{p}l^{(q)}$.  Moreover,} $\mymax(L)$ {\sl is the $p$-tuple} $(\mymax(L_{1}),\ldots,\mymax(L_{p}))$ {\sl while} $\mymin(L)$ {\sl is just} $(\mymin(L_{1}),\ldots,\mymin(L_{p}))$.  {\sl For any $\selt = (s_1, s_2, \ldots , s_p) \in L$ we have $\rho(\selt) = \sum_{q=1}^{p}\rho^{(q)}(s_{q})$.  For any other $\telt = (t_1, t_2, \ldots , t_p)$, we have $\selt \vee_{L} \telt = (s_{1} \vee t_{1}, \ldots , s_{p} \vee t_{p})$ -- a component-wise join -- while $\selt \wedge_{L} \telt = (s_{1} \wedge t_{1}, \ldots , s_{p} \wedge t_{p})$ -- a component-wise meet.\\ 
(2) Suppose $K$ is a vertex subset of $L$. 
Then $K$ is a full-length sublattice of $L$ if and only if the following two conditions hold: (}{\em i}{\sl ) There is a path in $L$ from} $\mymin(L)$ {\sl to} $\mymax(L)$ {\sl whose vertices are all from $K$, and (}{\em ii}{\sl ) $K$ is closed under component-wise joins and meets, i.e.\ for any $\selt, \telt \in K$ we have $\selt \vee_{L} \telt$ and $\selt \wedge_{L} \telt$ in $K$. 
In this case, $K$ is modular (resp.\ distributive) and is diamond-colored by the set $I$.} 

{\em Proof.} For {\sl (1)}, assume for the moment that each $L_{i}$ is a diamond-colored modular lattice. 
The facts that $L$ is diamond-colored; $\mymax(L)$ is the $p$-tuple $(\mymax(L_{1}),\ldots,\mymax(L_{p}))$;  $\mymin(L) = (\mymin(L_{1}),\ldots,\mymin(L_{p}))$; $\rho(\selt) = \sum_{q=1}^{p}\rho^{(q)}(s_{q})$ when  $\selt = (s_1, s_2, \ldots , s_p) \in L$; and $\posetlength(L) = \sum_{q=1}^{p}l^{(q)}$ are trivial consequences of definitions. 
Moreover, for any $\selt = (s_1, s_2, \ldots , s_p) \in L$ and $\telt = (t_1, t_2, \ldots , t_p)$ in $L$, it follows easily from definitions that the component-wise join $(s_{1} \vee t_{1}, \ldots , s_{p} \vee t_{p})$ is the unique least upper bound in $L$ of $\selt$ and $\telt$ and that the component-wise meet $(s_{1} \wedge t_{1}, \ldots , s_{p} \wedge t_{p})$ is their unique greatest lower bound. 
So $L$ is a lattice. 
Observe that:
\begin{eqnarray*}
2\rho(s_{1} \vee t_{1}, \ldots , s_{p} \vee t_{p}) - \rho(s_1, s_2, \ldots , s_p) - \rho(t_1, t_2, \ldots , t_p) &  & \\
& & \hspace*{-2in} =\ \sum_{i=1}^{p}\big(2\rho^{(i)}(s_{i} \vee \telt_{i}) - \rho^{(i)}(\selt_{i}) - \rho^{(i)}(\telt_{i})\big)\\
& & \hspace*{-2in} =\ \sum_{i=1}^{p}\big(\rho^{(i)}(\selt_{i}) + \rho^{(i)}(\telt_{i}) - 2\rho^{(i)}(s_{i} \wedge \telt_{i})\big)\\
& & \hspace*{-2in} =\ \rho(s_1, s_2, \ldots , s_p) + \rho(t_1, t_2, \ldots , t_p) - 2\rho(s_{1} \wedge t_{1}, \ldots , s_{p} \wedge t_{p}),
\end{eqnarray*}
so $L$ is modular. 
In a like manner, distributivity of $L$ follows from distributivity of each component, if each $L_{i}$ is distributive.  
For {\sl (2)}, the necessity of conditions ({\em i}) and  ({\em ii}) follows from the definitions. 
To see that these conditions are also sufficient, suppose ({\em i}) and  ({\em ii}) hold. 
It is enough to show that $K$ is ranked, because then the criteria of the definition of `full-length sublattice' will be met. 
Since a sublattice of a modular (respectively, distributive) lattice is also modular (resp.\ distributive) by \ModularDistributiveSublattices, then $K$ is modular (resp.\ distributive). 
In particular, $K$ is ranked.
Of course, $K$ inherits diamond-coloring from $L$, since by \FullLengthLemma, edges in $K$ are also edges in $L$.\hfill\QED

The following result, which is a diamond-colored version of Remark 2.1 of \cite{DonPeck}, can be applied to help find nice presentations of posets of join irreducibles for distributive lattices which arise as full-length sublattices of larger and more easily described distributive lattices (see e.g.\ \cite{DonPeck}, \cite{Gilliland}).

\noindent 
{\bf \FullLengthTheorem}\ \ {\sl (1) Let $P$ and $Q$ be vertex-colored posets with vertices colored by a set $I$.  
Suppose that for each $i \in I$, the vertices of color $i$ in $Q$ coincide with the vertices of color $i$ in $P$ (so in particular $P = Q$ as vertex sets).  
Further suppose that $Q$ is a weak subposet of $P$.  
Let $K := \Jcolor(P)$ and $L := \Jcolor(Q)$.  
Then $K$ is a full-length edge-colored sublattice of $L$.\\ 
(2) Conversely, suppose $L$ is a diamond-colored distributive lattice with edges colored by a set $I$. 
Suppose $K$ is a full-length edge-colored sublattice of $L$ (so $K$ is necessarily a diamond-colored distributive lattice).  
Let $Q := \jcolor(L)$ and $P := \jcolor(K)$. 
For any join irreducible $\xelt$ in $L$ (i.e.\ for any $\xelt \in Q$), the set $\{\yelt \in K\, |\, \xelt \leq_{L} \yelt\}$ has a unique minimal element $\welt_{\xelt}$, and $\welt_{\xelt}$ is a join irreducible in $K$ (so $\welt_{\xelt} \in P$). 
Moreover, the function $\phi:Q \longrightarrow P$ given by $\phi(\xelt) := \welt_{\xelt}$ is a vertex-color-preserving bijection, and if $\uelt \leq_{Q} \velt$ then $\phi(\uelt) \leq_{P} \phi(\velt)$. 
Now let $Q'$ be the set $P$ and declare that $\phi(\uelt) \leq_{Q'} \phi(\velt)$ if and only if $\uelt \leq_{Q} \velt$.  Then $Q'$ is a weak subposet of $P$ and $Q' \cong Q$ as vertex-colored posets.} 

{\em Proof.} The proof of {\sl (1)} is easy.  Let $\xelt$ be a down-set from $Q$.  
It follows from the definitions that $\xelt$ is also a down-set from $P$.  
So we get an inclusion $K = \Jcolor(Q) \subseteq \Jcolor(P) = L$.  
The length of $K$ (resp.\ $L$) is the cardinality of $Q$ (resp.\ $P$), and since $Q = P$ as vertex sets then $K$ and $L$ have the same length.  
Finally, note that for down-sets $\xelt$ and $\yelt$ from $Q$, $\xelt \vee_{K} \yelt = \xelt \cup \yelt = \xelt \vee_{L} \yelt$ and $\xelt \wedge_{K} \yelt = \xelt \cap \yelt = \xelt \wedge_{L} \yelt$.

For the proof of {\sl (2)}, we begin by choosing a join irreducible $\xelt$ in $L$.  
Let $\mathcal{F}_{\xelt} := \{\yelt \in K\, |\, \xelt \leq_{L} \yelt\}$.  
We claim that $\mathcal{F}_{\xelt}$ is an up-set from $K$ with a unique minimal element.  
First, if $\yelt \in \mathcal{F}_{\xelt}$ and $\yelt \leq_{K} \yelt'$ for some $\yelt' \in K$, then $\yelt \leq_{L} \yelt'$, and by transitivity of the partial order on $L$ it follows that $\xelt \leq_{L} \yelt'$.  
Hence $\yelt' \in \mathcal{F}_{\xelt}$. 
This shows that $\mathcal{F}_{\xelt}$ is an up-set from $K$. 
Second, if $\yelt$ and $\yelt'$ are both minimal elements of $\mathcal{F}_{\xelt}$, then whenever $\xelt \leq_{L} \yelt$ and $\xelt \leq_{L} \yelt'$ we will have $\xelt \leq_{L} (\yelt \wedge_{L} \yelt') = (\yelt \wedge_{K} \yelt')$.  
Hence $(\yelt \wedge_{K} \yelt') \in \mathcal{F}_{\xelt}$. 
Since $(\yelt \wedge_{K} \yelt') \leq_{K} \yelt$, $(\yelt \wedge_{K} \yelt') \leq_{K} \yelt'$, and $\yelt$ and $\yelt'$ are minimal elements of $\mathcal{F}_{\xelt}$, then we have $(\yelt \wedge_{K} \yelt') = \yelt$ and $(\yelt \wedge_{K} \yelt') = \yelt'$, i.e.\ $\yelt = \yelt'$.  
So $\mathcal{F}_{\xelt}$ has a unique minimal element. 

Let $\zelt$ be the unique minimal element of $\mathcal{F}_{\xelt}$, let $\mathcal{D}_{K}(\zelt) \subset K$ be the set of descendants of $\zelt$ in $K$, and let $\yelt$ be the unique descendant of $\xelt$ in $L$.  
We claim that for any $\zelt' \in \mathcal{D}_{K}(\zelt)$ we have $\xelt \vee_{L} \zelt' = \zelt$ and $\xelt \wedge_{L} \zelt' = \yelt$.  
To see this, note that when $\zelt' \rightarrow \zelt$ in $K$, we cannot have $\xelt \leq_{L} \zelt'$ or else $\zelt$ will not be minimal in $\mathcal{F}_{\xelt}$.  
So we cannot have $\zelt' = \xelt \vee_{L} \zelt'$.  
Therefore, $\zelt' <_{L} \xelt \vee_{L} \zelt'$.  
Then $\rho^{(L)}(\xelt \vee_{L} \zelt') \geq \rho^{(L)}(\zelt)$.  
But since $\xelt \leq_{L} \zelt$ and $\zelt' <_{L} \zelt$, we have $\rho^{(L)}(\xelt \vee_{L} \zelt') \leq \rho^{(L)}(\zelt)$.  
Hence $\rho^{(L)}(\xelt \vee_{L} \zelt') =  \rho^{(L)}(\zelt)$.  
It now follows that $\zelt = \xelt \vee_{L} \zelt'$.  
Next, since $\xelt \not\leq_{L} \zelt'$, then $\xelt \wedge_{L} \zelt' <_{L} \xelt$. 
But $\rho^{(L)}(\xelt \wedge_{L} \zelt') = \rho^{(L)}(\xelt) + \rho^{(L)}(\zelt') - \rho^{(L)}(\xelt \vee_{L} \zelt') =  \rho^{(L)}(\xelt) + \rho^{(L)}(\zelt) - 1 - \rho^{(L)}(\xelt \vee_{L} \zelt') = \rho^{(L)}(\xelt) - 1$.  
Thus $\xelt \wedge_{L} \zelt' \rightarrow \xelt$.  
But since $\yelt$ is the only element of $L$ covered by $\xelt$, then $\xelt \wedge_{L} \zelt' = \yelt$. 

Next we claim that $\zelt$ has exactly one descendant in $K$, i.e.\ $|\mathcal{D}_{K}(\zelt)| = 1$.  
Let $\zelt_{1}, \zelt_{2} \in \mathcal{D}_{K}(\zelt)$.  
Let $\zelt' := \zelt_{1} \wedge_{K} \zelt_{2}$.  
We will show that $\zelt' \vee_{L} \xelt = \zelt$ and $\zelt' \wedge_{L} \xelt = \yelt$.  
Since $\yelt \leq_{L} \zelt_{i}$ ($i = 1,2$)  by the previous paragraph, then $\yelt \leq_{L} \zelt_{1} \wedge_{L} \zelt_{2} = \zelt_{1} \wedge_{K} \zelt_{2} = \zelt'$.  
Since we also have $\yelt \leq_{L} \xelt$, then $\yelt \leq_{L} \zelt' \wedge_{L} \xelt$. 
Since $\zelt' \wedge_{L} \xelt \leq_{L} \xelt$ and $\yelt \rightarrow \xelt$, the only way to have $\yelt <_{L} \zelt' \wedge_{L} \xelt$ is if $\xelt = \zelt' \wedge_{L} \xelt$.  
But then we would have $\xelt \leq_{L} \zelt'$, which would mean $\zelt' \in \mathcal{F}_{\xelt}$.  
Then $\zelt \leq_{K} \zelt'$ by the minimality of $\zelt$.  
This contradicts the fact that $\zelt' \leq_{K} \zelt_{1} <_{K} \zelt$.  
So $\yelt = \zelt' \wedge_{L} \xelt$.  
Next, using a result from the previous paragraph we see that $\zelt' \vee_{L} \xelt = (\zelt_{1} \wedge_{K} \zelt_{2}) \vee_{L} \xelt = (\zelt_{1} \wedge_{L} \zelt_{2}) \vee_{L} \xelt = (\zelt_{1} \vee_{L} \xelt) \wedge_{L} (\zelt_{2} \vee_{L} \xelt) = \zelt \wedge_{L} \zelt = \zelt$. 
Now, $\rho^{(L)}(\zelt') + \rho^{(L)}(\xelt) = \rho^{(L)}(\zelt) + \rho^{(L)}(\yelt)$.  
Since $\rho^{(L)}(\yelt) = \rho^{(L)}(\xelt) - 1$, we have $\rho(\zelt') = \rho^{(L)}(\zelt) - 1$.  
Hence $\rho^{(L)}(\zelt') = \rho^{(L)}(\zelt_{i})$ for $i=1,2$.  
This can only happen if $\zelt_{1} = \zelt' = \zelt_{1} \wedge_{K} \zelt_{2} = \zelt_{2}$.  
Hence $\zelt_{1} = \zelt_{2}$.  
Next we argue that $\mathcal{D}_{K}(\zelt)$ is nonempty.  
We have that $\xelt \leq_{L} \zelt$ and (since $\xelt$ is join irreducible) $\rho^{(L)}(\xelt) > 0$.  
Therefore $\rho^{(K)}(\zelt) > 0$, so $\zelt$ is not the unique minimal element of $K$.  
In particular, $\mathcal{D}_{K}(\zelt)$ is nonempty. 
So $\zelt$ is join irreducible in $K$.  

With $P$ and $Q$ as in the theorem statement, we define a function $\phi: P \rightarrow Q$ by $\phi(\xelt) = \zelt$, where $\xelt$ and $\zelt$ are as in the preceding paragraphs.  
Next we show that $\phi$ is surjective.  
In doing so, we reset all variable names, so $\xelt, \yelt, \zelt,$ {\sl etc} are generic members of $K$ or $L$ until specified otherwise. 

Let $\zelt \in L$ be any join irreducible in $K$.  
Suppose $\zelt$ is also join irreducible in $L$. 
It follows that $\zelt$ is the unique minimal element of $\mathcal{F}_{\zelt}$.  
That is, $\zelt = \phi(\zelt)$.  
So now suppose $\zelt$ is not join irreducible in $L$.  
Let $\zelt'$ be the unique element of $K$ such that $\zelt' \rightarrow \zelt$.  
Define a set $\mathcal{S}_{\zelt} := \{\yelt \in L \setminus K\, |\, \yelt \vee_{L} \zelt' = \zelt\}$.  
Since $\zelt$ is not join irreducible in $L$, it follows that $\mathcal{S}_{\zelt}$ is nonempty.  
We claim $\mathcal{S}_{\zelt}$ has a unique minimal element.  
Indeed, suppose $\yelt$ and $\yelt'$ are minimal elements in $\mathcal{S}_{\zelt}$.  
Then $(\yelt \wedge_{L} \yelt') \vee_{L} \zelt' = (\yelt \vee_{L} \zelt') \wedge_{L} (\yelt' \vee_{L} \zelt') = \zelt \wedge_{L} \zelt = \zelt$.  
Since $\yelt <_{L} \zelt$ and $\yelt' <_{L} \zelt$, then $\yelt \wedge_{L} \yelt' <_{L} \zelt$.  
If $(\yelt \wedge_{L} \yelt') \in K$, then $(\yelt \wedge_{L} \yelt') <_{K} \zelt$. 
Then it must be the case that $(\yelt \wedge_{L} \yelt') \leq_{K} \zelt'$ since any path from $\yelt \wedge_{L} \yelt'$ up to $\zelt$ and that stays in $K$ must pass through $\zelt'$.  
But then we would have $(\yelt \wedge_{L} \yelt') \vee_{L} \zelt' = \zelt'$ instead of $(\yelt \wedge_{L} \yelt') \vee_{L} \zelt' = \zelt$.  
Then $(\yelt \wedge_{L} \yelt')$ is  in $L \setminus K$ and hence in $\mathcal{S}_{\zelt}$.  
Minimality of $\yelt$ and $\yelt'$ in $\mathcal{S}_{\zelt}$ then forces us to have $\yelt = (\yelt \wedge_{L} \yelt')  = \yelt'$.  
From here on, let $\xelt$ denote the unique minimal element of $\mathcal{S}_{\zelt}$.  

We have two claims: $\xelt$ is join irreducible in $L$, and $\zelt$ is the unique minimal element of $\mathcal{F}_{\xelt}$.  
Let $\xelt' := \xelt \wedge_{L} \zelt'$.  
Since $\rho(\xelt') = \rho(\xelt) + \rho(\zelt') - \rho(\zelt) = \rho(\xelt) - 1$, then $\xelt' \rightarrow \xelt$.  
Suppose $\xelt'' \rightarrow \xelt$ for some $\xelt'' \not= \xelt'$.  
It cannot be the case that $\xelt'' \leq_{L} \zelt'$, because otherwise $\xelt' \leq_{L} \zelt'$ and $\xelt'' \leq_{L} \zelt'$ means that $\xelt = (\xelt' \vee_{L} \xelt'') \leq_{L} \zelt'$, a contradiction.
Further, we have that $\xelt'' \in K$.  
Otherwise we would have $\xelt'' \in L \setminus K$, and since  $\xelt'' \not\leq_{L} \zelt'$ then $(\xelt'' \vee_{L} \zelt') = \zelt$.  
But then $\xelt''$ would be in $\mathcal{S}_{\zelt}$, violating minimality of $\xelt$.  
So $\xelt'' \in K$ and $\xelt'' \not\leq_{L} \zelt'$.  
Then there is a path from $\xelt''$ up to $\zelt$ that stays in $K$. 
But since $\zelt$ is join irreducible in $K$, then such a path must pass through $\zelt'$, implying that $\xelt'' \leq_{K} \zelt'$.  
But then $\xelt'' \leq_{L} \zelt'$, a contradiction.  
Therefore $\xelt'$ can be the only descendant of $\xelt$, hence $\xelt$ is join irreducible in $L$. 
Now if $\welt \in \mathcal{F}_{\xelt}$, then from the facts that $\xelt <_{L} \welt$ and $\xelt <_{L} \zelt$ we get $\xelt \leq_{L} (\welt \wedge_{L} \zelt)$.  
Since $(\welt \wedge_{L} \zelt) = (\welt \wedge_{K} \zelt)$, then $(\welt \wedge_{L} \zelt) \in K$, so we cannot have $\xelt = (\welt \wedge_{L} \zelt)$.  
Then $\xelt <_{L} (\welt \wedge_{L} \zelt)$.  
If $(\welt \wedge_{L} \zelt) <_{L} \zelt$, then we would have $\xelt \leq_{L} \zelt'$, which is not the case.  
So $(\welt \wedge_{L} \zelt) = \zelt$, and hence $\zelt \leq_{L} \welt$.  
So $\zelt$ is the unique minimal element of $\mathcal{F}_{\xelt}$.  
That is, $\zelt = \phi(\xelt)$.  

Our work in the preceding paragraphs shows that any join irreducible in $K$ is the image under $\phi$ of a join irreducible in $L$.  
That is, $\phi$ is surjective.
Since $|P| = |Q|$ ($K$ and $L$ have the same length), then $\phi$ is therefore a bijection.  
Suppose that $\zelt = \phi(\xelt) \not= \xelt$ for some $\xelt \in P$ and $\zelt \in Q$.  
Let $\xelt'$ be the unique descendant of $\xelt$ in $L$, with $\xelt' \myarrow{i} \xelt$ for some color $i$, and let $\zelt'$ be the unique descendant of $\zelt$ in $K$, with $\zelt' \myarrow{j} \zelt$ for some color $j$.  
Choose paths $\xelt' = \relt_{0} \mylongarrow{i=i_{1}} \xelt = \relt_{1} \myarrow{i_{2}} \relt_{2} \myarrow{i_{3}} \cdots \myarrow{i_{p-1}} \relt_{p-1} \myarrow{i_{p}} \relt_{p} = \zelt$ and $\xelt = \relt'_{0} \myarrow{j_{1}} \relt'_{1} \myarrow{j_{2}} \relt'_{2} \myarrow{j_{3}} \cdots \myarrow{j_{p-1}} \zelt' = \relt'_{p-1} \mylongarrow{j=j_{p}} \relt'_{p} = \zelt$ from $\xelt'$ up to $\zelt$.  
One path goes through $\xelt$ and the other through $\zelt'$.  
In particular, \ColorsLemma\ applies, so $i = i_{1} = j_{p} = j$. 
Since $\vcolor_{P}(\xelt) = i = j$ and $\vcolor_{Q}(\zelt) = j = i$, it follows that $\phi$ preserves vertex colors. 

To complete the proof of {\sl (2)}, we show that for $\uelt$ and $\velt$ in $P$, $\uelt \leq_{P} \velt$ implies that $\phi(\uelt) \leq_{Q} \phi(\velt)$.  
To see this, first note that $\uelt$ and $\velt$ are join irreducible elements of $L$ with $\uelt \leq_{L} \velt$.  
Consider $\mathcal{F}_{\uelt}$ and $\mathcal{F}_{\velt}$.  
If $\welt \in \mathcal{F}_{\velt}$, then $\welt \in K$ and $\velt \leq_{L} \welt$.  
Then $\uelt \leq_{L} \welt$ as well, so $\welt \in \mathcal{F}_{\uelt}$.  
So $\mathcal{F}_{\uelt} \supseteq \mathcal{F}_{\velt}$.  
Therefore $\phi(\uelt) \leq_{L} \phi(\velt)$.  
Since $\phi(\uelt)$ and $\phi(\velt)$ are both in $K$, then we have $\phi(\uelt) \leq_{K} \phi(\velt)$.  
Viewing $\phi(\uelt)$ and $\phi(\velt)$ as elements of $Q$, we then have $\phi(\uelt) \leq_{Q} \phi(\velt)$.\hfill\QED

{\bf [\S \SubstructureSection.4:\! Boolean sublattices.]} 
Our next result concerns the distributive lattice structure of certain intervals in diamond-colored distributive lattices. 
It is a diamond-colored version of comments from the penultimate paragraph of \S 3.4 of \cite{StanleyText}. 

\noindent 
{\bf \IntervalProp}\ \ {\sl Let $L$ be a diamond-colored distributive lattice.  
Let $\telt \in L$.  
Let $D$ be a subset of the descendants of $\telt$, to be notated as $D_{L}$ when viewed as an antichain in the induced order from $L$.  
For any $\selt \in D_{L}$,  let} $\vcolor_{D_{L}}(\selt) := \ecolor_{L}(\selt \rightarrow \telt)$, {\sl so $D_{L}$ is a vertex-colored poset.  
Let $\relt := \bigwedge_{\mbox{\tiny $\selt \in D$}}(\selt)$.  
Then $[\xelt, \telt]_{L} \cong \Mcolor(D_{L})$ and $D \subseteq [\xelt, \telt]_{L}$ for some $\xelt \in L$ if and only if $\xelt = \relt$, in which case $[\relt, \telt]_{L}$ is a Boolean lattice. 
Similarly, let $A$ be a subset of the ascendants of $\telt$, to be notated as $A_{L}$ when viewed as an antichain in the induced order from $L$. 
For any $\selt \in A_{L}$,  let} $\vcolor_{A}(\selt) := \ecolor_{L}(\telt \rightarrow \selt)$, {\sl so $A_{L}$ is a vertex-colored poset.   
Let $\uelt := \bigvee_{\mbox{\tiny $\selt \in A$}}(\selt)$.  
Then $[\telt, \xelt]_{L} \cong \Jcolor(A_{L})$ and $A \subseteq [\telt, \xelt]_{L}$ for some $\xelt \in L$ if and only if $\xelt = \uelt$, in which case $[\telt, \uelt]_{L}$ is a Boolean lattice.}

{\em Proof.}  
Of course, any two descendants of a given element of a poset are incomparable, as are any two ascendants.  
So, when the intervals $[\relt, \telt]_{L}$ and $[\telt, \uelt]_{L}$ meet the criteria described in the proposition statement, they must be Boolean lattices, cf.\ \BooleanExample. 
For the remainder of the proof, we only address the claim concerning the set $D$ since the proof for the claim concerning $A$ is entirely similar. 
For clarity, we write `$D_{\mbox{\tiny up-set}}$' when we view the set $D$ itself as an up-set taken from $D_{L}$. 
We reserve the notation `$D$' when $D$ is to be viewed simply as a subset of the descendants of $\telt$. 
Note that since $D_{L}$ is an antichain in the induced order from $L$, then every subset of $D$ is also an up-set from $D_{L}$. 

In the notation of the proposition statement, suppose $[\xelt, \telt]_{L} \cong \Mcolor(D_{L})$ and assume $D$ is a subset of $[\xelt,\telt]_{L}$.  
Let $\phi: \Mcolor(D_{L}) \rightarrow [\xelt, \telt]_{L}$ be the edge and edge-color preserving bijection. 
Since the unique maximal (resp.\ minimal) elements must correspond under the bijection $\phi$, then $\phi(\emptyset) = \telt$ (resp.\ $\phi(D_{\mbox{\tiny up-set}}) = \xelt$).  
For any $\selt \in D$ we have $\{\selt\} \rightarrow \emptyset$ in $\Mcolor(D_{L})$.  
Then $\phi(\{\selt\})$ must be covered by $\telt$ in $L$.  
So $\phi(D)$, thought of as a subset of elements of $[\xelt,\telt]_{L}$, is a subset of $D$, and since $\phi$ is a bijection we have that $\phi(D) \eqset D$.  
Then $\phi(D_{\mbox{\tiny up-set}}) = \phi(\cup_{\selt \in D}(\selt)) = \bigwedge_{\selt \in D}(\phi(\{\selt\}))$, where the meet is computed in $L$.  
But since $\phi(D) \eqset D$, then $\bigwedge_{\selt \in D}(\phi(\{\selt\})) = \bigwedge_{\selt \in D}(\selt)$, which is just $\relt$.  
That is, $\phi(D_{\mbox{\tiny up-set}}) = \relt$.  
Then $\xelt = \relt$.

For the converse, suppose that $\xelt = \relt$.  
Now, each $\selt \in D$ is a descendant of $\telt$, so $\selt \leq_{L} \telt$.  
By the definition of $\relt$ we have $\relt \leq_{L} \selt$.  
So $\selt \in [\relt,\telt]_{L}$.  
That is, $D \subseteq [\relt,\xelt]_{L}$.  
For any subset $S$ of $D$ (which is also an up-set from $D_{L}$), it follows from the definitions that $\relt \leq_{L} \bigwedge_{\selt \in S}(\selt) \leq_{L} \telt$, so $\bigwedge_{\selt \in S}(\selt) \in [\relt,\telt]_{L}$.  
Now define $\psi: \Mcolor(D_{L}) \rightarrow [\relt,\telt]_{L}$ by the rule $\psi(S) := \bigwedge_{\selt \in S}(\selt)$ for each subset $S \subseteq D$.  
Note that $\psi(\emptyset) = \telt$ and $\psi(D_{\mbox{\tiny up-set}}) = \relt$. 
We claim that if $S \myarrow{i} T$ in $\Mcolor(D_{L})$ then $\psi(S) \myarrow{i} \psi(T)$ in $[\relt,\telt]_{L}$.  
Now $S \myarrow{i} T$ in $\Mcolor(D)$ if and only if $|S| = |T| + 1$, $S = T \cup\{\selt\}$ for some $\selt \in D$, and $\vcolor_{D_{L}}(\selt) = i$. 
To establish our claim we induct on the size of $|S|$.  
If $|S| = 1$ then $S = \{\selt\}$ for some $\selt \in D$ and $T = \emptyset$. 
Then $\psi(S) = \selt$ and $\psi(T) = \telt$.  
Clearly $\selt \myarrow{i} \telt$ in this case.  
For our induction hypothesis we assume that $X \myarrow{i} Y$ in $\Mcolor(D_{L})$ implies  $\psi(X) \myarrow{i} \psi(Y)$ whenever $X$ has no more than $k$ elements, for some positive integer $k$.  
Now suppose $S \myarrow{i} T$ with $|S| = k+1$.  
So $|S| = |T| + 1$, $S = T \cup\{\selt\}$ for some $\selt \in D$, and $\vcolor_{D_{L}}(\selt) = i$.  
Let $Y := T \setminus \{\uelt\}$ for some $\uelt \in T$ with $\vcolor_{D_{L}}(\uelt) = j$, and let $X = S \setminus \{\uelt\}$.  
Then $Y = X \setminus \{\selt\}$. 
Then $T \myarrow{j} Y$, $X \myarrow{i} Y$, and $S \myarrow{j} X$.  
By the induction hypothesis, $\psi(T) \myarrow{j} \psi(Y)$ and $\psi(X) \myarrow{i} \psi(Y)$. 
Then $\psi(Y) = (\psi(X) \vee \psi(T))$.  
We claim that $\psi(S) = (\psi(X) \wedge \psi(T))$.  
Let $\yelt = \psi(Y)$.  
Then $\psi(X) = \yelt \wedge \selt$, $\psi(T) = \yelt \wedge \uelt$, and $\psi(S) = \yelt \wedge (\selt \wedge \uelt)$. 
So $\psi(S) = \yelt \wedge (\selt \wedge \uelt) = (\yelt \wedge \yelt) \wedge (\selt \wedge \uelt) = (\yelt \wedge \selt) \wedge (\yelt \wedge \uelt) = (\psi(X) \wedge \psi(T))$. 
Since $\psi(T) \myarrow{j} \psi(Y)$ and $\psi(X) \myarrow{i} \psi(Y)$ in our diamond-colored distributive lattice $L$, we must therefore have $\psi(S) \myarrow{j} \psi(X)$ and $\psi(S) \myarrow{i} \psi(T)$.  
This completes the induction step, and the proof of our claim. 

Let $d = |D|$.  
Let $D_{\mbox{\tiny up-set}} = S^{(0)} \rightarrow S^{(1)} \rightarrow S^{(2)} \rightarrow \cdots \rightarrow S^{(d-1)} \rightarrow S^{(d)} = \emptyset$ be a chain of maximal length in $\Mcolor(D_{L})$.  
Then $\relt = \psi(D_{\mbox{\tiny up-set}}) \rightarrow \psi(S^{(1)}) \rightarrow \cdots \rightarrow \psi(S^{(d-1)}) \rightarrow \psi(S^{(d)}) = \telt$, a chain of maximal length in $[\relt,\telt]_{L}$.  
In particular, the length of $[\relt,\telt]_{L}$ is $d$.  
In the paragraph preceding the proposition it was noted that intervals in diamond-colored modular lattices are edge-colored sublattices.  
Invoking the distributivity hypothesis for $L$,  we conclude that $[\relt,\telt]_{L}$ is a diamond-colored distributive lattice. 
Since $[\relt,\telt]_{L}$ has length $d$ as a ranked poset, it follows that $[\relt,\telt]_{L}$ must have precisely $d$ meet irreducibles.  
But each $\selt \in D$ is meet irreducible in $[\relt,\telt]_{L}$, so the set $D$ must account for all meet irreducibles in $[\relt,\telt]_{L}$.  
Therefore, $[\relt,\telt]_{L} \cong \Mcolor(D_{L})$ by \FundamentalTheorem.\hfill\QED  

{\bf [\S \SubstructureSection.5:\! $J$-components of DCML's and DCDL's.]} 
The next result concerns the structure of $J$-components of a DCML and will be applied in the development of \JCompTheorem. 
\JCompStuff\ have no direct analogs for monochromatic modular $\!$/$\!$ distributive lattices and are, to the best of our knowledge, new. 

\noindent 
{\bf \JCompResult}\ \ {\sl Let $L$ be a diamond-colored modular lattice with edge colors from a set $I$. 
If $\telt \in L$ and $J \subseteq I$, then} $\comp_{J}(\telt)$ {\sl is the covering digraph for a diamond-colored modular lattice. 
Moreover,} $\comp_{J}(\telt)$ {\sl is an edge-colored  sublattice of $L$. 
If $L$ is distributive, then so is} $\comp_{J}(\telt)$.  

{\em Proof.} 
The claim made in the last sentence of the theorem statement follows immediately once we prove the remaining assertions.  
To that end, let $K := \comp_{J}(\telt)$.  
The edge-colored and directed graph $K$ is a poset with partial order $\leq_{K}$ given as follows: 
For $\xelt$ and $\yelt$ in $K$, $\xelt \leq_{K} \yelt$ if and only if there is a path $(\xelt = \selt_{0},\selt_{1},\ldots,\selt_{p}=\yelt)$ of elements from $K$ such that $\xelt = \selt_{0} \mylongarrow{j_{1}} \selt_{1} \mylongarrow{j_{2}} \selt_{2} \mylongarrow{j_{3}} \cdots \mylongarrow{j_{p-1}} \selt_{p-1} \mylongarrow{j_{p}} \selt_{p} = \yelt$ is a sequence of edges in $L$ with each $j_{q} \in J$.  
It is easy to see that this defines a partial order on $K$ and that the edges $\pelt \myarrow{j} \qelt$ in $K$ are precisely the covering relations for this partial order.  

To complete the proof, it suffices to show the following: $\xelt \vee_{L} \yelt$ and $\xelt \wedge_{L} \yelt$ are in $K$ whenever $\xelt, \yelt \in K$.  
We actually prove the stronger claim that any shortest path from $\xelt$ to $\yelt$ in $K$ is also a shortest path in $L$ and moreover $\xelt \vee_{L} \yelt \in K$ and $\xelt \wedge_{L} \yelt \in K$.  
To do so, we induct on the distance $\dist_{K}(\xelt,\yelt)$ between vertices. 
If $\dist_{K}(\xelt,\yelt) = 0$, then $\xelt=\yelt$; the path $(\xelt)$, whose length is zero, is the only shortest path and is a shortest path in both $K$ and $L$; the common element $\xelt \vee_{L} \yelt = \xelt = \xelt \wedge_{L} \yelt$ is in $K$. 
As our inductive hypothesis, we assume that for some positive integer $p$, it is the case that any shortest path from $\xelt$ to $\yelt$ in $K$ is also a shortest path in $L$ and moreover $\xelt \vee_{L} \yelt \in K$ and $\xelt \wedge_{L} \yelt \in K$ whenever $\dist_{K}(\xelt,\yelt) \leq p-1$. 

Now suppose $\dist_{K}(\xelt,\yelt) = p$. 
Let $\mathcal{P} := (\xelt=\selt_{0},\ldots,\selt_{p}=\yelt)$ be a shortest path in $K$ from $\xelt$ to $\yelt$. 
Assume for the moment that $\selt_{p-1} \mylongbackarrow{j_{p}} \selt_{p}$, and let $\mathcal{P}' := (\xelt=\selt_{0},\ldots,\selt_{p-2},\selt_{p-1})$. 
Then $\mathcal{P}'$ from $\xelt$ to $\selt_{p-1}$ has length $p-1$ and must be a shortest path in $K$. 
By our inductive hypothesis, $\mathcal{P'}$ is also a shortest path in $L$. 
Since $L$ is topographically balanced and $\mathcal{P}'$ is a shortest path in $L$ from $\xelt$ to $\selt_{p-1}$, we can apply our {\bf mountain-path-construction} algorithm to $\mathcal{P}'$ to obtain a shortest mountain path $\mathcal{P}'_{\mbox{\em \scriptsize mountain}}$ from $\xelt$ to $\selt_{p-1}$ whose apex $\aelt = \xelt \vee_{L} \selt_{p-1}$. 
At each application of Step 2 during the execution of the {\bf mountain-path-construction} algorithm, diamond-coloring of $L$ ensures that the multiset of path colors $\{j_{1},\ldots,j_{p-1}\}$ is unchanged. 
So, the multiset of edge colors $\{j_{1},\ldots,j_{p-1}\}$ for $\mathcal{P}'$ is the same as the multiset of edge colors for $\mathcal{P}'_{\mbox{\em \scriptsize mountain}}$. 
Thus, $\mathcal{P}'_{\mbox{\em \scriptsize mountain}}$ is a path in $K$, and in particular $\aelt \in K$. 
We let $\aelt' := \xelt \vee_{L} \yelt$ and observe that $\aelt' \leq \aelt$ since $\xelt \leq \aelt$ and $\yelt < \selt_{p-1} \leq \aelt$. 
Take any path $\mathcal{Q}$ in $L$ from $\xelt$ up to $\aelt$ that goes through $\aelt'$. 
Another path $\mathcal{Q}'$ that starts at $\xelt$ and ascends to $\aelt$ is the initial part of $\mathcal{P}'_{\mbox{\em \scriptsize mountain}}$. 
Since $\mathcal{Q}'$ is a path in $K$, its edge colors are from the set $J$. 
By \ColorsLemma, the multiset of edge colors in the path $\mathcal{Q}$ coincides with the multiset of edge colors in the path $\mathcal{Q}'$. 
Therefore, $\mathcal{Q}$ is a path in $K$, and so $\aelt' \in K$. 
Let $\mathcal{R}'$ be the path that starts at $\aelt$ and descends to $\selt_{p-1}$ on edges from $\mathcal{P}'_{\mbox{\em \scriptsize mountain}}$ and then descends from $\selt_{p-1}$ to $\yelt$ along our given edge of color $j_{p}$. 
Further, let $\mathcal{R}$ be any path in $L$ from $\aelt$ down to $\yelt$ that goes through $\aelt'$. 
Again by \ColorsLemma, the multiset of edge colors in the path $\mathcal{R}$ coincides with the multiset of edge colors in the path $\mathcal{R}'$. 
Therefore $\mathcal{R}$ is a path in $K$. 
Consider the concatenation of $\mathcal{Q}$ and $\mathcal{R}$ but with any edges in the interval $[\aelt',\aelt]_{K}$ removed; the result is a mountain path $\mathcal{M}$ from $\xelt$ to $\yelt$ in $K$ whose apex is $\aelt' = \xelt \vee_{L} \yelt$. 
Moreover, this mountain path $\mathcal{M}$ is shortest in $L$ by \ModularLatticeTheorem, so $\dist_{L}(\xelt,\yelt)$ is the length of $\mathcal{M}$. 
So, $\mathcal{M}$ must also be shortest in $K$, and therefore $\mathcal{M}$ has length $p = \dist_{K}(\xelt,\yelt)$. 
That is, the path $\mathcal{P}$ of length $p$ from $\xelt$ to $\yelt$ in $K$ that we were originally given is a shortest path in $L$. 
Moreover, if we apply our {\bf valley-path-construction} algorithm to $\mathcal{M}$, the valley path $\mathcal{M}_{\mbox{\em \scriptsize valley}}$ we obtain is a shortest path in $K$ (since, due to earlier observations concerning Step 2 of the algorithm and diamond-coloring, its multiset of edge colors is the same as the multiset of edge colors for $\mathcal{M}$) and is shortest in $L$ (a feature of the algorithm is that it preserves shortest path lengths). 
Therefore by \ModularLatticeTheorem, its nadir is $\xelt \wedge_{L} \yelt$, which is therefore an element of $K$. 

To summarize our work so far on the inductive step of our proof, we have shown that if $\mathcal{P} := (\xelt=\selt_{0},\ldots,\selt_{p}=\yelt)$ is a shortest path in $K$ from $\xelt$ to $\yelt$ such that $\selt_{p-1} \mylongbackarrow{j_{p}} \selt_{p}=\yelt$, then $\mathcal{P}$ a shortest path in $L$ from $\xelt$ to $\yelt$, $\xelt \vee_{L} \yelt \in K$, and $\xelt \wedge_{L} \yelt_{K}$. 
An entirely similar argument can be given in the case that, for our given path $\mathcal{P}$, we have $\selt_{p-1} \mylongarrow{j_{p}} \selt_{p}=\yelt$. 
This completes the inductive step, and the proof.\hfill\QED

A famous result of R.\ P.\ Dilworth \cite{Dilworth} states that in any modular lattice and for any positive integer $k$, the number of elements with exactly $k$ descendants is the same as the number of elements with exactly $k$ ascendants (e.g.\ see the solution to Ch.\ 3 Exercise 101 (d) in \cite{StanleyText}). 
So, for example, the number of join irreducible elements is the same as the number of meet irreducible elements. 
Let $J$ be a subset of our color palette $I$. 
If $j \in J$ and $\xelt \myarrow{j} \yelt$ is an edge within some diamond-colored modular lattice $L$, then say that $\xelt$ is a $J${\em -descendant} of $\yelt$ and that $\yelt$ is a $J${\em -ascendant} of $\xelt$. 
Let $\#_{\mbox{\tiny $J$-asc}}(\xelt)$ denote the number of $J$-ascendants of $\xelt$, and let $\#_{\mbox{\tiny $J$-desc}}(\yelt)$ denote the number of $J$-descendants of $\yelt$. 

\noindent 
{\bf \DilworthsTheorem}\ \ {\sl Let $L$ be a DCML, let $J \subseteq I$, and let $k \in [1,\infty)_{\mathbb{Z}}$. 
Then the number of elements $\xelt \in L$ such that} $\#_{\mbox{\tiny $J$-asc}}(\xelt) = k$ {\sl equals the number of elements $\yelt \in L$ such that} $\#_{\mbox{\tiny $J$-desc}}(\yelt) = k$.

{\em Proof.} By \JCompResult, Dilworth's result generalizes to our multi-colored setting.\hfill\QED 

The next result considers full-length sublattices in a very specific monochromatic scenario. 
Our set-up and proof require the following notions. 
The {\em max cover number} of a modular lattice $L$ is $\max\{\#_{\mbox{\tiny $I$-desc}}(\xelt)\}_{\xelt \in L}$, which is the maximum number of elements any single element of $L$ covers. 
By \DilworthsTheorem, the max cover number is also $\max\{\#_{\mbox{\tiny $I$-asc}}(\xelt)\}_{\xelt \in L}$, which is the maximum number of elements above any single element of $L$.  
Define a {\em down-up walk} in $L$ to be a sequence of $2p+1$ terms $(\xelt_{0},\xelt_{1},\ldots,\xelt_{2p})$, where $p$ is a nonnegative integer, such that $\xelt_{2i} \leftarrow \xelt_{2i+1} \rightarrow \xelt_{2i+2}$ for $0 \leq i \leq p-1$; of course, this down-up walk is simple if the $2p+1$ terms are pairwise distinct vertices. 
Similarly define {\em up-down walk}. 
A reminder before we state the proposition: Viewing the partial order on a product $[a,b]_{\mathbb{Z}} \times [c,d]_{\mathbb{Z}}$ of nonempty chains as in \FullLengthWithinProduct.1, we see that if $(p,q) \rightarrow (r,s)$ is a covering relation within the product, then either $(p,q)+(1,0) = (r,s)$ or $(p,q)+(0,1) = (r,s)$.  

\noindent 
{\bf \TwoChainLemma}\ \ {\sl A modular lattice has max cover number not exceeding two if and only if it is isomorphic to a full-length sublattice of a product of two (nonempty) chains. 
In this case, the lattice is distributive.}  

{\em Proof.}  
Let $L$ be our given modular lattice with rank function $\rho_{L}$ and length $l := \posetlength(L)$. 
The last sentence of the proposition follows from the first: A modular lattice with max cover number not exceeding two cannot have a sublattice isomorphic to $\myM_{5}$. 
The `if' claim of the first sentence of the proposition statement is straightforward, so we just need a proof of the `only if' claim. 

So, assume $L$ has max cover number not exceeding two. 
We will algorithmically construct an isomorphism, to be called $\phi$, from $L$ to a full-length sublattice $K$ of a product of two chains. 
We do so by reading $L$ systematically across each rank, starting at $\mymin(L)$ and ending at $\mymax(L)$, and identifying each element of $L$ with some unique pair in $[0,l_{1}]_{\mathbb{Z}} \times [0,l_{2}]_{\mathbb{Z}}$, where $l_{1}$ and $l_{2}$ are yet-to-be-identified nonnegative integers with $l_{1}+l_{2} = l$. 
We use $K$ to denote the set of pairs in $[0,l_{1}]_{\mathbb{Z}} \times [0,l_{2}]_{\mathbb{Z}}$ realized by this process and regard $K$ to be a poset in the induced order. 
Throughout the process to follow, if some $\selt \in L$ has been assigned an integer pair via $\phi$, say $\phi(\selt) = (p,q)$, then we refer to $p =: \phi_{x}(\selt)$ as the `$x$-coordinate' of $\selt$ and $q =: \phi_{y}(\selt)$ as the `$y$-coordinate' of $\selt$. 
For each $r \in \{0,1,\ldots,l\}$, let $N(r) := |\rho_{L}^{-1}(r)|$, i.e.\ the number of elements of $L$ with rank equal to $r$.

\begin{description} 
\item[\em Step 0:] {To initialize our process of constructing $\phi$, we declare that $\phi(\mymin(L)) := (0,0)$, and temporarily set $l'_{1} := 0 =: l'_{2}$, $\mbox{\sffamily \small counter} := 0$, and $r:=\mbox{\sffamily \small counter}$. 
In the next step, we begin our process of reading ranks. 
In general, when reading the $k^{\mbox{\tiny th}}$ rank, we assume that a well-defined pair $\phi(\relt) = (p,q) \in [0,l'_{1}]_{\mathbb{Z}} \times [0,l'_{2}]_{\mathbb{Z}}$ has been assigned to each $\relt$ of the given rank; that, for any $\relt_{1}, \relt_{2} \in L$ with $\rho_{L}(\relt_{1}) \leq k$ and $\rho_{L}(\relt_{2}) \leq k$, we have $\phi(\relt_{1} \wedge \relt_{2}) = \phi(\relt_{1}) \wedge \phi(\relt_{2})$, where the meet on the right-hand side is the `component-wise min' $\big(\min(\phi_{x}(\relt_{1}),\phi_{x}(\relt_{2})),\min(\phi_{y}(\relt_{1}),\phi_{y}(\relt_{2}))\big)$; and that there is a simple down-up walk $(\selt_{0},\selt_{1},\ldots,\selt_{2N(k)-3},\selt_{2(N(k)-1)})$, where $\rho_{L}(\selt_{0}) = k = \phi_{x}(\selt_{0}) + \phi_{y}(\selt_{0})$ and where $\phi(\selt_{2i}) = \phi(\selt_{0})+(-i,i)$ when $i \in \{0,1,\ldots,N(k)-1\}$. 
Of course, the foregoing hypotheses hold when $k=0$.}

\item[\em Step 1:] {Let $\selt$ (resp.\ $\selt'$) be the element of $L$ with rank $r$ whose $x$-coordinate is largest (resp.\ smallest), so that, by the hypotheses expressed in Step 0, we have $\rho_{L}(\selt) = \phi_{x}(\selt)+\phi_{y}(\selt) = r$ and $\phi(\selt')=\phi(\selt)+(-N(r)+1,N(r)-1)$. 
For use later, we let $(p,q) := \phi(\selt)$ be the pair that $\phi$ has assigned to $\selt$. 
If $l > r$, we proceed to Step 2. 
Otherwise, we have $l = r$, so our construction of $\phi$ is complete, and to end the algorithm we take $l_{1} := l_{1}'$ and $l_{2} := l_{2}'$.}

\item[\em Step 2:] {Suppose that $N(r) = 1$. 
Let $\telt$ be an arbitrarily chosen ascendant of $\selt$ (so $\rho_{L}(\telt) = r+1$) and set $\phi(\telt) := (p+1,q)$; if $l'_{1} = p$, then replace $l'_{1}$ with $p+1$. 
If $\telt$ is the only ascendant of $\selt$, then skip to Step 4. 
Otherwise, let $\telt' \not= \telt$ denote the second ascendant of $\selt$ (so $\rho_{L}(\telt') = r+1$). Declare that $\phi(\telt') := (p,q+1)$, and if $l'_{2} = q$, then replace $l'_{2}$ with $q+1$. 
This accounts for all possible ascendants of $\selt$ and hence accounts for all elements of $L$ with rank $r+1$ and all of the edges below the rank $r+1$ elements. 
Note that for each element with rank $r+1$, the sum of the $x$- and $y$-coordinates of its assigned pair is $r+1$.  
Also observe that $(\telt,\selt,\telt')$ is a simple down-up walk with three distinct vertices and that $\phi(\telt') = \phi(\telt)+(-1,1)$. 
Finally, we verify that for all $\relt_{1}, \relt_{2} \in L$ with rank not exceeding $r+1$, we have $\phi(\relt_{1} \wedge \relt_{2}) = \phi(\relt_{1}) \wedge \phi(\relt_{2})$. 
This is true by hypothesis if $\max(\rho_{L}(\relt_{1}),\rho_{L}(\relt_{2})) \leq r$. 
If $\rho_{L}(\relt_{1})=r+1$ and $\rho_{L}(\relt_{2}) \leq r$, then $\relt_{2} \leq \selt$, hence $\relt_{2} \leq \relt_{1}$ and $\phi(\relt_{2}) \leq \phi(\relt_{1})$, and so $\phi(\relt_{1} \wedge \relt_{2}) = \phi(\relt_{2}) = \phi(\relt_{1}) \wedge \phi(\relt_{2})$. 
Similarly, $\phi(\relt_{1} \wedge \relt_{2}) = \phi(\relt_{1}) = \phi(\relt_{1}) \wedge \phi(\relt_{2})$ when $\rho_{L}(\relt_{1}) \leq r$ and $\rho_{L}(\relt_{2}) = r+1$
If $\rho_{L}(\relt_{1}) = \rho_{L}(\relt_{2})$, then $\relt_{1} \wedge \relt_{2} = \selt$ and hence $\phi(\relt_{1} \wedge \relt_{2}) = \phi(\selt) = \phi(\relt_{1}) \wedge \phi(\relt_{2})$. 
Now skip to Step 4.}

\item[\em Step 3:] {On the other hand, suppose $N(r) > 1$, so necessarily $r > 0$. 
We assume at this point in our process that $\phi(\relt_{1} \wedge \relt_{2}) = \phi(\relt_{1}) \wedge \phi(\relt_{2})$ when $\rho_{L}(\relt_{i}) \leq r$ (where $i \leq \{1,2\}$) and that we have a simple down-up walk $(\relt_{0},\relt_{1},\ldots,\relt_{2N(r)-3},\relt_{2(N(r)-1)})$ from $\relt_{0} = \selt$ to $\relt_{2(N(r)-1)} = \selt'$ wherein $\phi(\relt_{2i}) = \phi(\relt_{0})+(-i,i)$ for each $i \in \{0,1,\ldots,N(r)-1\}$ with $0 \leq \phi_{x}(\relt_{2i}) \leq l_{1}'$ and $0 \leq \phi_{y}(\relt_{2i}) \leq l_{2}'$. 
For convenience, set $\selt_{i} := \relt_{2i}$ for all $i \in \{0,1,\ldots,N(r)-1\}$. 
Now let $\telt_{i+1} := \selt_{i} \vee \selt_{i+1}$ for any $i \in \{0,1,\ldots,N(r)-2\}$. 
In particular, $\telt_{1}$ is an ascendant of $\selt$ and $\telt_{N(r)-1}$ is an ascendant of $\selt'$. 
Note that, for any $i,j \in \{1,2,\ldots,N(r)-1\}$, $\telt_{i} = \telt_{j}$ means that $\telt_{i}$ is above both $\selt_{i-1}$ and $\selt_{i}$ (which are distinct) and above both $\selt_{j-1}$ and $\selt_{j}$ (which are distinct); since the max cover number in $L$ cannot exceed two we are forced to have $\selt_{i-1} = \selt_{j-1}$ and $\selt_{i} = \selt_{j}$ and hence $i=j$. 
That is, the list $\telt_{1}, \telt_{2}, \ldots, \telt_{N(r)-1}$ consists of $N(r)-1$ distinct vertices. 
We now consider four cases that depend on the possibilities for other ascendants of $\selt$ and $\selt'$.}

\begin{description}
\item[\mbox{\textnormal{\underline{Case 1:}}}] Suppose neither $\selt$ nor $\selt'$ has any other ascendants besides those identified already. 
Then let $\phi(\telt_{i}) := \phi(\selt_{0})+(-i+1,i)$ for each $i \in \{1,2,\ldots,N(r)-1\}$. 
The $x$-coordinates of these $\telt_{i}$'s are between $p-N(r)+2 = \phi_{x}(\selt')$ and $p = \phi(\selt)$ inclusive while the $y$-coordinates are between $q+1 = \phi_{y}(\selt)+1$ and $q+N(r)-1 = \phi_{y}(\selt')$ inclusive. 
In particular, we have $\phi(\telt_{i+1}) \in [0,l'_{1}]_{\mathbb{Z}} \times [0,l'_{2}]_{\mathbb{Z}}$ and we have no need to update $l_{1}'$ or $l_{2}'$. 
Now, $\selt_{0}, \selt_{1}, \ldots, \selt_{N(r)-2}, \selt_{N(r)-1}$ are all of the elements in $L$ with rank $r$, and each of $\selt_{1}, \ldots, \selt_{N(r)-2}$ is covered by two elements with rank $r+1$. 
Since each of $\selt_{0} = \selt$ and $\selt_{N(r)-1} = \selt'$ has only one ascendant, then our list $\telt_{1}, \telt_{2}, \ldots, \telt_{N(r)-1}$ accounts for all elements with rank $r+1$ and the edges encountered in the foregoing discussion are all possible edges between these two ranks.  
Also, observe that $\selt_{i} \rightarrow \telt_{j}$ only if $j \in \{i,i+1\}$, in which case $\phi(\telt_{i}) = \phi(\selt_{i})+(1,0)$ and $\phi(\telt_{i+1}) = \phi(\selt_{i})+(0,1)$. 
So, in view of the sentence preceding the proposition statement, the rank $r$ elements $\phi(\selt_{0}), \ldots, \phi(\selt_{N(r)-1})$ and the rank $r+1$ elements $\phi(\telt_{1}), \ldots, \phi(\telt_{N(r)-1})$ together with the edges between them induced by the partial order on $[0,l'_{1}]_{\mathbb{Z}} \times [0,l'_{2}]_{\mathbb{Z}}$ perfectly mirror the edges between ranks $r$ and $r+1$ in $L$. 

\item[\mbox{\textnormal{\underline{Case 2:}}}] Suppose $\selt'$ has no other ascendants but $\selt$ has another ascendant, which we name $\telt_{0}$. 
Since, by hypothesis, $\telt_{1} \not= \telt_{0}$, and since each of $\telt_{2}, \ldots, \telt_{N(r)-1}$ already covers two elements from the list $\selt_{1}, \ldots, \selt_{N(r)-1}$, then $\telt_{0}$ is distinct from any of the vertices $\telt_{1}, \ldots, \telt_{N(r)-1}$. 
Then let $\phi(\telt_{i}) := \phi(\selt_{0})+(-i+1,i)$ for each for each $i \in \{0,1,\ldots,N(r)-1\}$; if $l'_{1} = p$, then replace $l'_{1}$ with $p+1$. 
Reasoning similar to Case 1 gives us $\phi(\telt_{i}) \in [0,l'_{1}]_{\mathbb{Z}} \times [0,l'_{2}]_{\mathbb{Z}}$ with no need to update $l_{2}'$. 
Also, we can see as in Case 1 the rank $r$ elements $\phi(\selt_{0}), \ldots, \phi(\selt_{N(r)-1})$ and the rank $r+1$ elements $\phi(\telt_{0}), \ldots, \phi(\telt_{N(r)-1})$ together with the edges between them induced by the partial order on $[0,l'_{1}]_{\mathbb{Z}} \times [0,l'_{2}]_{\mathbb{Z}}$ perfectly mirror the edges between ranks $r$ and $r+1$ in $L$.  

\item[\mbox{\textnormal{\underline{Case 3:}}}] Suppose $\selt$ has no other ascendants but $\selt'$ has another ascendant, which we name $\telt_{N(r)}$. 
Since, by hypothesis, $\telt_{N(r)} \not= \telt_{N(r)-1}$, and since each of $\telt_{1}, \ldots, \telt_{N(r)-2}$ already covers two elements from the list $\selt_{0}, \ldots, \selt_{N(r)-2}$, then $\telt_{N(r)}$ is distinct from any of the vertices $\telt_{1}, \ldots, \telt_{N(r)-1}$. 
Then let $\phi(\telt_{i}) := \phi(\selt_{0})+(-i+1,i)$ for each for each $i \in \{1,\ldots,N(r)\}$; if $l'_{2} = q+N(r)-1$, then replace $l'_{2}$ with $q+N(r)$. 
Reasoning similar to Case 1 gives us $\phi(\telt_{i}) \in [0,l'_{1}]_{\mathbb{Z}} \times [0,l'_{2}]_{\mathbb{Z}}$ with no need to update $l_{1}'$. 
Also, we can see as in Case 1 the rank $r$ elements $\phi(\selt_{0}), \ldots, \phi(\selt_{N(r)-1})$ and the rank $r+1$ elements $\phi(\telt_{1}), \ldots, \phi(\telt_{N(r)})$ together with the edges between them induced by the partial order on $[0,l'_{1}]_{\mathbb{Z}} \times [0,l'_{2}]_{\mathbb{Z}}$ perfectly mirror the edges between ranks $r$ and $r+1$ in $L$.  

\item[\mbox{\textnormal{\underline{Case 4:}}}] Last, suppose each of $\selt$ and $\selt'$ has another ascendant which we name $\telt_{0}$ and $\telt_{N(r)}$ respectively. 
Since each of $\telt_{2}, \ldots, \telt_{N(r)-2}$ already covers two elements from the list $\selt_{1}, \ldots, \selt_{N(r)-2}$, then $\{\telt_{0},\telt_{N(r)}\} \cap \{\telt_{1}, \ldots, \telt_{N(r)-1}\} = \emptyset$. 
We wish to show that the vertices $\telt_{0}$ and $\telt_{N(r)}$ are distinct.  
Suppose otherwise, so $\telt_{0} = \telt_{N(r)}$. 
Then the down-up walk $(\telt_{0}, \selt_{0}, \telt_{1},\ldots,\telt_{N(r)-1}, \selt_{N(r)-1}, \telt_{N(r)})$ is a cycle whose only coincident vertices are $\telt_{0} = \telt_{N(r)}$. 
For $i \in [0,N(r)-1]_{\mathbb{Z}}$, let $\uelt_{i} := \telt_{i} \vee \telt_{i+1}$, so $(\telt_{0}, \uelt_{0}, \telt_{1},\ldots,\telt_{N(r)-1}, \uelt_{N(r)-1}, \telt_{N(r)})$ is a cyclic up-down walk whose only coincident vertices are $\telt_{0} = \telt_{N(r)}$. 
We can view the preceding walk as a cyclic down-up walk by writing $(\uelt_{0}, \telt_{0}, \uelt_{1},\ldots,\uelt_{N(r)-1}, \telt_{N(r)-1}, \uelt_{N(r)})$, where $\uelt_{N(r)} := \uelt_{0}$ are the only coincident vertices. 
It is clear that we can indefinitely continue this process of constructing such cyclic down-up walks at higher ranks, which contradicts the finiteness of $L$. 
We conclude that we must have $\telt_{0} \not= \telt_{N(r)}$.

Now let $\phi(\telt_{i}) := \phi(\selt_{0})+(-i+1,i)$ for each $i \in \{0,1,\ldots,N(r)\}$; if $l'_{1} = p$, then replace $l'_{1}$ with $p+1$, and if $l'_{2} = q+N(r)-1$, then replace $l'_{2}$ with $q+N(r)$. 
Reasoning similar to Case 1 gives us $\phi(\telt_{i}) \in [0,l'_{1}]_{\mathbb{Z}} \times [0,l'_{2}]_{\mathbb{Z}}$. 
Also, we can see as in Case 1 the rank $r$ elements $\phi(\selt_{0}), \ldots, \phi(\selt_{N(r)-1})$ and the rank $r+1$ elements $\phi(\telt_{0}), \ldots, \phi(\telt_{N(r)})$ together with the edges between them induced by the partial order on $[0,l'_{1}]_{\mathbb{Z}} \times [0,l'_{2}]_{\mathbb{Z}}$ perfectly mirror the edges between ranks $r$ and $r+1$ in $L$.  
\end{description}

In each case, we see that, for all appropriate $i$, $\phi_{x}(\telt_{i})+\phi_{y}(\telt_{i}) = r+1$ and that there is a down-up walk that starts at the rank $r+1$ element with largest $x$-coordinate (either $\telt_{0}$ or $\telt_{1}$, depending on the case), ends at the element with the smallest $x$-coordinate (either $\telt_{N(r)-1}$ or $\telt_{N(r)}$, depending on the case), and with $\phi(\telt_{i+1}) = \phi(\telt_{i}) + (-1,1)$. 
Before concluding this step of the algorithm, we need to verify that $\phi(\relt_{1} \wedge \relt_{2}) = \phi(\relt_{1}) \wedge \phi(\relt_{2})$ when $\rho_{L}(\relt_{i}) \leq r+1$ for $i \in \{1,2\}$; that proof begins in the next paragraph. 
That will complete our verification that the assignment of integer pairs to the elements of rank $r+1$ satisfies the hypotheses expressed in Step 0, which are required in order to repeat the algorithm for this next rank. 

It suffices to show $\phi(\relt_{1} \wedge \relt_{2}) = \phi(\relt_{1}) \wedge \phi(\relt_{2})$ when at least one of $\rho_{L}(\relt_{1}) = r+1$ or $\rho_{L}(\relt_{2}) = r+1$. 
Say $\rho_{L}(\relt_{1}) = r+1$ and $\rho_{L}(\relt_{2}) \leq r$. 
If $\relt_{2} \leq \relt_{1}$, then $\phi(\relt_{2}) \leq \phi(\relt_{2})$ and so $\phi(\relt_{1} \wedge \relt_{2}) = \phi(\relt_{1}) = \phi(\relt_{1}) \wedge \phi(\relt_{2})$. 
So we will assume $\relt_{1}$ and $\relt_{2}$ are not comparable.  
Let $\mathcal{P}_{\mbox{\em \scriptsize valley}} = (\uelt_{0}=\relt_{1}, \uelt_{1}, \ldots, \uelt_{a}, \ldots, \uelt_{b}=\relt_{2})$ be a valley path from $\relt_{1}$ to $\relt_{2}$ whose nadir is $\uelt_{a} = \relt_{1} \wedge \relt_{2}$, so by \ModularLatticeTheorem\ $\mathcal{P}_{\mbox{\em \scriptsize valley}}$ is a shortest path. 
In particular, the valley path $\mathcal{P\ }'_{\mbox{\em $\!\!\!\!$\scriptsize valley}} = (\uelt_{1}, \ldots, \uelt_{a}, \ldots, \uelt_{b}=\relt_{2})$ is a shortest path from $\uelt_{1}$ to $\relt_{2}$ and its nadir $\uelt_{a}$ must be $\uelt_{1} \wedge \relt_{2}$. 
It is easy to reason that $\uelt_{1}$ and $\relt_{2}$ are not comparable. 
Now, $\phi(\uelt_{1} \wedge \relt_{2}) = \phi(\uelt_{1}) \wedge \phi(\relt_{2})$, since both elements have rank not exceeding $r$. 
Of course, $\phi(\relt_{1} \wedge \relt_{2}) \leq \phi(\relt_{1}) \wedge \phi(\relt_{2})$, which means that  
\begin{equation}\begin{array}{l}\mbox{$\phi(\uelt_{1}) \wedge \phi(\relt_{2}) \leq \phi(\relt_{1}) \wedge \phi(\relt_{2})$,}\, \mbox{or equivalently }\hspace*{1in}\\
\hspace*{1in}\mbox{$\min\{\phi_{z}(\uelt_{1}),\phi_{z}(\relt_{2})\} \leq \min\{\phi_{z}(\relt_{1}),\phi_{z}(\relt_{2})\}$}\, \mbox{for $z \in \{x,y\}$.} 
\end{array}
\end{equation}
We refer to the top-most inequality of (1) as a `poset inequality' and the bottom-most as a `numerical inequality'. 

For sake of argument, suppose $\phi_{x}(\uelt_{1})+1=\phi_{x}(\relt_{1})$, in which case $\phi_{y}(\uelt_{1})=\phi_{y}(\relt_{1})$. 
If $\phi_{x}(\relt_{2}) > \phi_{x}(\uelt_{1})$, then the left-hand side of the numerical inequality (1) becomes the quantity $\phi_{x}(\uelt_{1})$ while the right-hand side is $\phi_{x}(\relt_{1}) = \phi_{x}(\uelt_{1})+1$. 
In particular, the poset inequality (1) is strict. 
We must also have $\phi_{y}(\relt_{1}) = \phi_{y}(\uelt_{1}) > \phi_{y}(\relt_{2})$, because the opposite hypothesis necessitates $\phi(\relt_{2}) \geq \phi(\relt_{1})$, which is absurd. 
Now, 
$\phi(\uelt_{a+1}) = \big(\phi_{x}(\uelt_{1})+1,\phi_{y}(\relt_{2})\big)$. 
Then, we see that $\phi(\uelt_{a+1}) < \phi(\relt_{1})$, while we already know $\phi(\uelt_{a+1}) \leq \phi(\relt_{2})$. 
Over in $L$, we have $\relt_{1} \wedge \relt_{2} < \uelt_{a+1}$ by our set-up, but now $\uelt_{a+1} < \relt_{1}$ and $\uelt_{a+1} \leq \relt_{2}$ implies that $\uelt_{a+1} \leq \relt_{1} \wedge \relt_{2}$. 
From the contradiction in the previous sentence, we see that the hypothesis `$\phi_{x}(\relt_{2}) > \phi_{x}(\uelt_{1})$' is false, so that $\phi_{x}(\relt_{2}) \leq \phi_{x}(\uelt_{1})$. 
Therefore, $\phi_{x}(\relt_{2}) \leq \phi_{x}(\uelt_{1}) < \phi_{x}(\uelt_{1})+1 = \phi_{x}(\relt_{1})$, and hence $\phi_{y}(\relt_{2}) > \phi_{y}(\relt_{1})$ (else $\relt_{1}$ and $\relt_{2}$ are comparable).  
And so: $\phi(\relt_{1} \wedge \relt_{2}) = \phi(\uelt_{1} \wedge \relt_{2}) = \phi(\uelt_{1}) \wedge \phi(\relt_{2}) = (\phi_{x}(\relt_{2}),\phi_{y}(\uelt_{1})) = (\phi_{x}(\relt_{2}),\phi_{y}(\relt_{1})) = \phi(\relt_{1}) \wedge \phi(\relt_{2})$. 

We can similarly argue that $\phi(\relt_{1} \wedge \relt_{2}) = \phi(\relt_{1}) \wedge \phi(\relt_{2})$ when $\phi_{y}(\uelt_{1})+1=\phi_{y}(\relt_{1})$.  
Of course, the roles of $\relt_{1}$ and $\relt_{2}$ can be interchanged throughout the above arguments. 
We conclude that when exactly one of $\relt_{1}$ or $\relt_{2}$ has rank $r+1$, then $\phi(\relt_{1} \wedge \relt_{2}) = \phi(\relt_{1}) \wedge \phi(\relt_{2})$. 
Therefore, to conclude that the latter equality is an identity for all $\relt_{1}$ and $\relt_{2}$ with rank not exceeding $r+1$, we must establish the equality when $\rho_{L}(\relt_{1}) = r+1 = \rho_{L}(\relt_{2})$. 
To see this, let $\mathcal{P}_{\mbox{\em \scriptsize valley}} = (\uelt_{0}=\relt_{1}, \uelt_{1}, \ldots, \uelt_{a}, \ldots, \uelt_{b}=\relt_{2})$ be a valley path (and therefore a shortest path) from $\relt_{1}$ to $\relt_{2}$ whose nadir is $\uelt_{a} = \relt_{1} \wedge \relt_{2}$. 
Notice also that the valley paths $\mathcal{P\ }'_{\mbox{\em $\!\!\!\!$\scriptsize valley}} = (\uelt_{1}, \ldots, \uelt_{a}, \ldots, \uelt_{b}=\relt_{2})$ and $\mathcal{P\ }''_{\mbox{\em $\!\!\!\!$\scriptsize valley}} = (\relt_{1},\uelt_{1}, \ldots, \uelt_{a}, \ldots, \uelt_{b-1})$ are shortest with nadir $\uelt_{a} = \uelt_{1} \wedge \relt_{2} = \relt_{1} \wedge \uelt_{b-1}$. 
By the definition of `$\wedge$', we have $\phi(\relt_{1} \wedge \relt_{2}) \leq \phi(\uelt_{1}) \wedge \phi(\relt_{2}) = \phi(\relt_{1}) \wedge \phi(\uelt_{b-1})$. 
Without loss of generality, we can assume that $\phi_{x}(\relt_{1}) > \phi_{x}(\relt_{2})$ and that $\phi_{y}(\relt_{1}) < \phi_{y}(\relt_{2})$. 
So, $\phi_{x}(\uelt_{1}) \geq \phi_{x}(\relt_{2})$ and $\phi_{y}(\relt_{1}) \leq \phi_{y}(\uelt_{b-1})$. 
Then the $x$-coordinate of $\phi(\uelt_{1}) \wedge \phi(\relt_{2})$ is $\phi_{x}(\relt_{2})$ and the $y$-coordinate of $\phi(\relt_{1}) \wedge \phi(\uelt_{b-1})$ is $\phi_{y}(\relt_{1})$. 
Therefore, $\phi(\relt_{1} \wedge \relt_{2}) = \phi(\uelt_{a}) = (\phi_{x}(\relt_{2}),\phi_{y}(\relt_{1})) = \phi(\relt_{1}) \wedge \phi(\relt_{2})$. 

\item[\em Step 4:] {Add the value `$1$' to $\mbox{\sffamily \small counter}$, and let $r$ be this new value of $\mbox{\sffamily \small counter}$.  Then return to Step 1.}
\end{description} 

Let $K$ be the set of pairs $\{\phi(\relt)\}_{\relt \in L} \subseteq [0,l_{1}]_{\mathbb{Z}} \times [0,l_{2}]_{\mathbb{Z}}$, and give $K$ the order induced by the natural partial ordering on the latter product of chains. 
By the above construction, with respect to the bijection $\phi: L \longrightarrow K$, $\selt \rightarrow \telt$ in $L$ if and only if $\phi(\selt) \rightarrow \phi(\telt)$ in $K$, so $\phi$ is a poset isomorphism. 
Now, $\phi(\mymin(L)) = (0,0)$ and $\phi(\mymax(L)) = (l_{1},l_{2})$. 
Also, let $\selt_{0}^{(r)}$ denote the element in $L$ with rank $r$ and having the largest $x$-coordinate within its rank. 
Our algorithm shows that $\selt_{0}^{(0)} \rightarrow \selt_{0}^{(1)} \rightarrow \cdots \rightarrow \selt_{0}^{(l)}$ is a path in $L$ from $\mymin(L)$ up to $\mymax(L)$ and that $(0,0) = \phi(\selt_{0}^{(0)}) \rightarrow \phi(\selt_{0}^{(1)}) \rightarrow \cdots \rightarrow \phi(\selt_{0}^{(l)}) = (l_{1},l_{2})$ is a path in $K$ from the unique minimal element of $[0,l_{1}]_{\mathbb{Z}} \times [0,l_{2}]_{\mathbb{Z}}$ up to the unique maximal element. 

In order to apply \FullLengthWithinProduct.2 and complete the proof, we show that, for any $\phi(\selt)$ and $\phi(\telt)$ in $K$, their component-wise meet and join remain in $K$. 
The fact that $\phi(\relt_{1}) \wedge \phi(\relt_{2}) = \phi(\relt_{1} \wedge \relt_{2})$ when $\rho_{L}(\relt_{i}) \leq r$ (where $i \in \{1,2\}$) was verified for all $r$ as we algorithmically constructed $\phi$. 
Say the rank of any pair $(a,b) \in [0,l_{1}]_{\mathbb{Z}} \times [0,l_{2}]_{\mathbb{Z}}$ is given by $\rho(a,b) = a+b$. 
Now, clearly $\phi(\relt_{1}) \vee \phi(\relt_{2}) \leq \phi(\relt_{1} \vee \relt_{2})$ in $[0,l_{1}]_{\mathbb{Z}} \times [0,l_{2}]_{\mathbb{Z}}$. 
Also, $\rho\big(\phi(\relt_{1}) \vee \phi(\relt_{2})\big) = \rho\big(\phi(\relt_{1})\big) +\rho\big(\phi(\relt_{2})\big) -  \rho\big(\phi(\relt_{1}) \wedge \phi(\relt_{2})\big) = \rho\big(\phi(\relt_{1})\big) +\rho\big(\phi(\relt_{2})\big) -  \rho\big(\phi(\relt_{1} \wedge \relt_{2})\big) = \rho_{L}(\relt_{1})+\rho_{L}(\relt_{2}) - \rho_{L}(\relt_{1} \wedge \relt_{2}) = \rho_{L}(\relt_{1} \vee \relt_{2}) = \rho\big(\phi(\relt_{1} \vee \relt_{2})\big)$, i.e.\ $\rho\big(\phi(\relt_{1}) \vee \phi(\relt_{2})\big) = \rho\big(\phi(\relt_{1} \vee \relt_{2})\big)$. 
The latter fact, together with $\phi(\relt_{1}) \vee \phi(\relt_{2}) \leq \phi(\relt_{1} \vee \relt_{2})$, means that $\phi(\relt_{1}) \vee \phi(\relt_{2}) = \phi(\relt_{1} \vee \relt_{2})$.\hfill\QED

The next result (\JCompTheorem) provides a way to identify the compression posets of the $J$-components of DCDL as certain induced-order subposets of its vertex-colored compression poset.  
Suppose $P$ is an $I$-vertex-colored poset with vertex coloring function $\vcolor: P \longrightarrow I$. 
Let $L = \Jcolor(P)$, and suppose $J \subseteq I$.  
It follows from \JCompResult\ that for any $\telt \in L$, $\comp_{J}(\telt)$ is an edge-colored distributive sublattice of $L$.  
Regarding $\telt$ to be a lower order ideal from $P$, let $D_{J}(\telt)$ be the subset of $\telt$ such that the lower order ideal $\telt \setminus D_{J}(\telt)$ is the minimal element of $\comp_{J}(\telt)$.  
Notice that $\vcolor(D_{J}(\telt)) \subseteq J$. 
(``$D$'' is short for ``delete.'') 
Similarly let $A_{J}(\telt)$ be the subset of the upper order ideal $P \setminus \telt$ such that $\telt \cup A_{J}(\telt)$ is the maximal element of $\comp_{J}(\telt)$.  
Then $\vcolor(A_{J}(\telt)) \subseteq J$.  
(``$A$'' is short for ``add.'') 
The set $D_{J}(\telt)$ is ``largest'' in the following sense: If $D$ is any set of vertices in $P$ with colors from $J$ such that $D \subseteq \telt$ and $\telt \setminus D$ is a down-set from $P$, then $D \subseteq D_{J}(\telt)$. 
Similarly, $A_{J}(\telt)$ is ``largest'': If $A$ is any set of vertices in $P$ with colors from $J$ such that $A \cap \telt = \emptyset$ and $\telt \cup A$ is a down-set from $P$, then $A \subseteq A_{J}(\telt)$. 
Let $Q_{J}(\telt) := A_{J}(\telt) \cup D_{J}(\telt)$. View $Q_{J}(\telt)$ as a subposet in the induced order and with the inherited vertex coloring. 

Let $Q$ be a subposet of $P$ in the induced order and with the inherited vertex coloring.  
Call $Q$ a $J$-{\em subordinate} of $P$ if ({\em i}) $\vcolor(Q) \subseteq J$; ({\em ii}) there is a down-set $\relt$ from $P$ such that $\relt \cap Q = \emptyset$, $\relt \cup Q$ is a down-set from $P$, i.e.\ $Q = \selt \setminus \relt$ where for some down-set $\selt$ containing the down-set $\relt$; and ({\em iii}) $\vcolor(v) \in I \setminus J$ whenever $v$ is a maximal (respectively, minimal) element of $\relt$ (resp.\ $P \setminus (\relt \cup Q)$). 

\noindent 
{\bf \JCompTheorem}\ \ {\sl Given $I$, $J$, $P$, and $L$ as above.  
(1) For any $\telt \in L$, $Q_{J}(\telt)$ is a $J$-subordinate of $P$.  Each $J$-subordinate of $P$ is precisely $Q_{J}(\telt)$ for some $\telt \in L$. 
(2) Moreover, for each $\telt \in L$,} $\comp_{J}(\telt) \cong \Jcolor(Q_{J}(\telt))$ {\sl and equivalently} $\jcolor(\comp_{J}(\telt)) \cong Q_{J}(\telt)$. 

{\em Proof.} For (1), let $Q = Q_{J}(\telt)$.  
With $\relt := \telt \setminus D_{J}(\telt) = \mymin(\comp_{J}(\telt))$ and $\relt \cup Q = \mymax(\comp_{J}(\telt))$, it is easy to see that $Q$ meets the criteria for a $J$-subordinate of $P$. 
On the other hand, given some $J$-subordinate $Q$ of $P$, let $\telt$ be the down-set $\relt$ of part ({\em ii}) of the definition of $J$-subordinate.  
Then it is easy to see that $Q = Q_{J}(\telt)$. 
For (2), let $\phi: \Jcolor(Q_{J}(\telt)) \longrightarrow \comp_{J}(\telt)$ be given by $\phi(\xelt) = \xelt \cup \telt$. 
It is routine to check that $\phi$ is an edge and edge-color preserving mapping between diamond-colored distributive lattices.\hfill\QED

\newpage
\newcommand{\CFourFlower}{
\setlength{\unitlength}{0.5cm}
\begin{picture}(12,1)
\put(-2.5,-0.25){\LARGE $\myC_{4}$}
\put(-0.9,0){\vector(1,0){1.2}}
{\color{Cyan}\put(1,0){\circle*{0.7}}}
\put(1,0){\circle{0.7}}
{\color{Red}\put(4,0){\circle*{0.7}}}
\put(4,0){\circle{0.7}}
{\color{Purple}\put(7,0){\circle*{0.7}}}
\put(7,0){\circle{0.7}}
{\color{Green}\put(10,0){\circle*{0.7}}}
\put(10,0){\circle{0.7}}
\thicklines
\put(8,0){\vector(1,0){0.2}}
\put(8.65,0){\vector(-1,0){0.2}}
\put(9,0){\vector(-1,0){0.2}}
\put(1.35,0){\line(1,0){2.3}}
\put(4.35,0){\line(1,0){2.3}}
\put(7.35,0){\line(1,0){2.3}}
\put(0.8,0.6){\textcolor{Cyan}{\em 1}}
\put(3.8,0.6){\textcolor{Red}{\em 2}}
\put(6.8,0.6){\textcolor{Purple}{\em 3}}
\put(9.8,0.6){\textcolor{Green}{\em 4}}
\end{picture}
}

\begin{figure}
\begin{center}
{\bf \CFourFigure.1}\ \ The 28-element diamond-colored distributive lattice depicted below is denoted `$\widetilde{L_{\mytinyC_{4}}}(\omega_{\mbox{\color{Red}{\tiny $2$}}})$'. 
It can be viewed as a recoloring of the DCDL $L_{\mytinyA_{7}}(\omega_{\mbox{\color{Red}{\tiny $2$}}})$ from \ASevenFigureTwo.\\ 
{\footnotesize (The lattice depicted here is, in the language of \S \GStructureIntroSection, `$\mysmallC_{4}$-structured', where $\mysmallC_{4}$ is the graph}\\ 
{\footnotesize  immediately below that also appears in the list given in \IEGGraphFigure.)}

\vspace*{0.2in} 
\CFourFlower

\vspace*{0.4in} 
\setlength{\unitlength}{1cm}
\begin{picture}(7,13)
\put(3,12.5){\LARGE $\widetilde{L_{\mysmallC_{4}}}(\omega_{\mbox{\color{Red}{\small $2$}}})$}
\put(2.85,12.6){\vector(-2,-1){1.5}}
\put(3.4,11.7){\LARGE $\cong \mathbf{J}_{\mbox{\footnotesize color}}\big(\widetilde{P_{\mysmallC_{4}}}(\omega_{\mbox{\color{Red}{\small $2$}}})\big)$}
\put(0,0){\TypeEboxDot{Gray}}
\put(1,1){\TypeEboxDot{Gray}}
\put(0,2){\TypeEboxDot{Gray}}
\put(2,2){\TypeEboxDot{Gray}}
\put(1,3){\TypeEboxDot{Gray}}
\put(3,3){\TypeEboxDot{Gray}}
\put(0,4){\TypeEboxDot{Gray}}
\put(2,4){\TypeEboxDot{Gray}}
\put(4,4){\TypeEboxDot{Gray}}
\put(1,5){\TypeEboxDot{Gray}}
\put(3,5){\TypeEboxDot{Gray}}
\put(5,5){\TypeEboxDot{Gray}}
\put(0,6){\TypeEboxDot{Gray}}
\put(2,6){\TypeEboxDot{Gray}}
\put(4,6){\TypeEboxDot{Gray}}
\put(6,6){\TypeEboxDot{Gray}}
\put(1,7){\TypeEboxDot{Gray}}
\put(3,7){\TypeEboxDot{Gray}}
\put(5,7){\TypeEboxDot{Gray}}
\put(0,8){\TypeEboxDot{Gray}}
\put(2,8){\TypeEboxDot{Gray}}
\put(4,8){\TypeEboxDot{Gray}}
\put(1,9){\TypeEboxDot{Gray}}
\put(3,9){\TypeEboxDot{Gray}}
\put(0,10){\TypeEboxDot{Gray}}
\put(2,10){\TypeEboxDot{Gray}}
\put(1,11){\TypeEboxDot{Gray}}
\put(0,12){\TypeEboxDot{Gray}}
\thicklines
\put(0.625,0.375){\color{Red}\qbezier(0,0)(0.375,0.375)(0.75,0.75)}
\put(1.625,1.375){\color{Purple}\qbezier(0,0)(0.375,0.375)(0.75,0.75)}
\put(2.625,2.375){\color{Green}\qbezier(0,0)(0.375,0.375)(0.75,0.75)}
\put(3.625,3.375){\color{Purple}\qbezier(0,0)(0.375,0.375)(0.75,0.75)}
\put(4.625,4.375){\color{Red}\qbezier(0,0)(0.375,0.375)(0.75,0.75)}
\put(5.625,5.375){\color{Cyan}\qbezier(0,0)(0.375,0.375)(0.75,0.75)}
\put(0.625,2.375){\color{Purple}\qbezier(0,0)(0.375,0.375)(0.75,0.75)}
\put(1.625,3.375){\color{Green}\qbezier(0,0)(0.375,0.375)(0.75,0.75)}
\put(2.625,4.375){\color{Purple}\qbezier(0,0)(0.375,0.375)(0.75,0.75)}
\put(3.625,5.375){\color{Red}\qbezier(0,0)(0.375,0.375)(0.75,0.75)}
\put(4.625,6.375){\color{Cyan}\qbezier(0,0)(0.375,0.375)(0.75,0.75)}
\put(0.625,4.375){\color{Green}\qbezier(0,0)(0.375,0.375)(0.75,0.75)}
\put(1.625,5.375){\color{Purple}\qbezier(0,0)(0.375,0.375)(0.75,0.75)}
\put(2.625,6.375){\color{Red}\qbezier(0,0)(0.375,0.375)(0.75,0.75)}
\put(3.625,7.375){\color{Cyan}\qbezier(0,0)(0.375,0.375)(0.75,0.75)}
\put(0.625,6.375){\color{Purple}\qbezier(0,0)(0.375,0.375)(0.75,0.75)}
\put(1.625,7.375){\color{Red}\qbezier(0,0)(0.375,0.375)(0.75,0.75)}
\put(2.625,8.375){\color{Cyan}\qbezier(0,0)(0.375,0.375)(0.75,0.75)}
\put(0.625,8.375){\color{Red}\qbezier(0,0)(0.375,0.375)(0.75,0.75)}
\put(1.625,9.375){\color{Cyan}\qbezier(0,0)(0.375,0.375)(0.75,0.75)}
\put(0.625,10.375){\color{Cyan}\qbezier(0,0)(0.375,0.375)(0.75,0.75)}
\put(1.39,11.375){\color{Red}\qbezier(0,0)(-0.375,0.375)(-0.75,0.75)}
\put(2.39,10.375){\color{Purple}\qbezier(0,0)(-0.375,0.375)(-0.75,0.75)}
\put(3.39,9.375){\color{Green}\qbezier(0,0)(-0.375,0.375)(-0.75,0.75)}
\put(4.39,8.375){\color{Purple}\qbezier(0,0)(-0.375,0.375)(-0.75,0.75)}
\put(5.39,7.375){\color{Red}\qbezier(0,0)(-0.375,0.375)(-0.75,0.75)}
\put(6.39,6.375){\color{Cyan}\qbezier(0,0)(-0.375,0.375)(-0.75,0.75)}
\put(1.39,9.375){\color{Purple}\qbezier(0,0)(-0.375,0.375)(-0.75,0.75)}
\put(2.39,8.375){\color{Green}\qbezier(0,0)(-0.375,0.375)(-0.75,0.75)}
\put(3.39,7.375){\color{Purple}\qbezier(0,0)(-0.375,0.375)(-0.75,0.75)}
\put(4.39,6.375){\color{Red}\qbezier(0,0)(-0.375,0.375)(-0.75,0.75)}
\put(5.39,5.375){\color{Cyan}\qbezier(0,0)(-0.375,0.375)(-0.75,0.75)}
\put(1.39,7.375){\color{Green}\qbezier(0,0)(-0.375,0.375)(-0.75,0.75)}
\put(2.39,6.375){\color{Purple}\qbezier(0,0)(-0.375,0.375)(-0.75,0.75)}
\put(3.39,5.375){\color{Red}\qbezier(0,0)(-0.375,0.375)(-0.75,0.75)}
\put(4.39,4.375){\color{Cyan}\qbezier(0,0)(-0.375,0.375)(-0.75,0.75)}
\put(1.39,5.375){\color{Purple}\qbezier(0,0)(-0.375,0.375)(-0.75,0.75)}
\put(2.39,4.375){\color{Red}\qbezier(0,0)(-0.375,0.375)(-0.75,0.75)}
\put(3.39,3.375){\color{Cyan}\qbezier(0,0)(-0.375,0.375)(-0.75,0.75)}
\put(1.39,3.375){\color{Red}\qbezier(0,0)(-0.375,0.375)(-0.75,0.75)}
\put(2.39,2.375){\color{Cyan}\qbezier(0,0)(-0.375,0.375)(-0.75,0.75)}
\put(1.39,1.375){\color{Cyan}\qbezier(0,0)(-0.375,0.375)(-0.75,0.75)}
\put(0.9,0.65){\color{Red}{\em 2}}
\put(0.9,2.65){\color{Purple}{\em 3}}
\put(1.9,1.65){\color{Purple}{\em 3}}
\put(0.9,4.65){\color{Green}{\em 4}}
\put(1.9,3.65){\color{Green}{\em 4}}
\put(2.9,2.65){\color{Green}{\em 4}}
\put(0.9,6.65){\color{Purple}{\em 3}}
\put(1.9,5.65){\color{Purple}{\em 3}}
\put(2.9,4.65){\color{Purple}{\em 3}}
\put(3.9,3.65){\color{Purple}{\em 3}}
\put(0.9,8.65){\color{Red}{\em 2}}
\put(1.9,7.65){\color{Red}{\em 2}}
\put(2.9,6.65){\color{Red}{\em 2}}
\put(3.9,5.65){\color{Red}{\em 2}}
\put(4.9,4.65){\color{Red}{\em 2}}
\put(0.9,10.65){\color{Cyan}{\em 1}}
\put(1.9,9.65){\color{Cyan}{\em 1}}
\put(2.9,8.65){\color{Cyan}{\em 1}}
\put(3.9,7.65){\color{Cyan}{\em 1}}
\put(4.9,6.65){\color{Cyan}{\em 1}}
\put(5.9,5.65){\color{Cyan}{\em 1}}
\put(0.9,11.65){\color{Red}{\em 2}}
\put(0.9,9.65){\color{Purple}{\em 3}}
\put(1.9,10.65){\color{Purple}{\em 3}}
\put(0.9,7.65){\color{Green}{\em 4}}
\put(1.9,8.65){\color{Green}{\em 4}}
\put(2.9,9.65){\color{Green}{\em 4}}
\put(0.9,5.65){\color{Purple}{\em 3}}
\put(1.9,6.65){\color{Purple}{\em 3}}
\put(2.9,7.65){\color{Purple}{\em 3}}
\put(3.9,8.65){\color{Purple}{\em 3}}
\put(0.9,3.65){\color{Red}{\em 2}}
\put(1.9,4.65){\color{Red}{\em 2}}
\put(2.9,5.65){\color{Red}{\em 2}}
\put(3.9,6.65){\color{Red}{\em 2}}
\put(4.9,7.65){\color{Red}{\em 2}}
\put(0.9,1.65){\color{Cyan}{\em 1}}
\put(1.9,2.65){\color{Cyan}{\em 1}}
\put(2.9,3.65){\color{Cyan}{\em 1}}
\put(3.9,4.65){\color{Cyan}{\em 1}}
\put(4.9,5.65){\color{Cyan}{\em 1}}
\put(5.9,6.65){\color{Cyan}{\em 1}}
\end{picture}
\begin{picture}(7,13)
\put(3.5,11){\LARGE $\widetilde{P_{\mysmallC_{4}}}(\omega_{\mbox{\color{Red}{\small $2$}}})$}
\put(3.35,11){\vector(-1,-1){1}}
\put(4.1,10.2){\LARGE $\cong \mathbf{j}_{\mbox{\footnotesize color}}\big(\widetilde{L_{\mysmallC_{4}}}(\omega_{\mbox{\color{Red}{\small $2$}}})\big)$}
\put(0,9){\TypeEboxDot{Cyan}}
\put(0.1,9){\color{Cyan}{\em \footnotesize 1}}
\put(1,10){\TypeEboxDot{Red}}
\put(1,8){\TypeEboxDot{Red}}
\put(1.7,10.3){\color{Red}{\em \footnotesize 2}}
\put(1.1,8){\color{Red}{\em \footnotesize 2}}
\put(2,9){\TypeEboxDot{Purple}}
\put(2,7){\TypeEboxDot{Purple}}
\put(2.7,9.3){\color{Purple}{\em \footnotesize 3}}
\put(2.1,7){\color{Purple}{\em \footnotesize 3}}
\put(3,8){\TypeEboxDot{Green}}
\put(3,6){\TypeEboxDot{Green}}
\put(3.7,8.3){\color{Green}{\em \footnotesize 4}}
\put(3.1,6){\color{Green}{\em \footnotesize 4}}
\put(4,7){\TypeEboxDot{Purple}}
\put(4,5){\TypeEboxDot{Purple}}
\put(4.7,7.3){\color{Purple}{\em \footnotesize 3}}
\put(4.1,5){\color{Purple}{\em \footnotesize 3}}
\put(5,6){\TypeEboxDot{Red}}
\put(5,4){\TypeEboxDot{Red}}
\put(5.7,6.3){\color{Red}{\em \footnotesize 2}}
\put(5.1,4){\color{Red}{\em \footnotesize 2}}
\put(6,5){\TypeEboxDot{Cyan}}
\put(6.7,5.3){\color{Cyan}{\em \footnotesize 1}}
\thicklines
\put(5.625,4.375){\color{Black}\qbezier(0,0)(0.375,0.375)(0.75,0.75)}
\put(4.625,5.375){\color{Black}\qbezier(0,0)(0.375,0.375)(0.75,0.75)}
\put(3.625,6.375){\color{Black}\qbezier(0,0)(0.375,0.375)(0.75,0.75)}
\put(2.625,7.375){\color{Black}\qbezier(0,0)(0.375,0.375)(0.75,0.75)}
\put(1.625,8.375){\color{Black}\qbezier(0,0)(0.375,0.375)(0.75,0.75)}
\put(0.625,9.375){\color{Black}\qbezier(0,0)(0.375,0.375)(0.75,0.75)}
\put(2.39,9.375){\color{Black}\qbezier(0,0)(-0.375,0.375)(-0.75,0.75)}
\put(3.39,8.375){\color{Black}\qbezier(0,0)(-0.375,0.375)(-0.75,0.75)}
\put(4.39,7.375){\color{Black}\qbezier(0,0)(-0.375,0.375)(-0.75,0.75)}
\put(5.39,6.375){\color{Black}\qbezier(0,0)(-0.375,0.375)(-0.75,0.75)}
\put(6.39,5.375){\color{Black}\qbezier(0,0)(-0.375,0.375)(-0.75,0.75)}
\put(1.39,8.375){\color{Black}\qbezier(0,0)(-0.375,0.375)(-0.75,0.75)}
\put(2.39,7.375){\color{Black}\qbezier(0,0)(-0.375,0.375)(-0.75,0.75)}
\put(3.39,6.375){\color{Black}\qbezier(0,0)(-0.375,0.375)(-0.75,0.75)}
\put(4.39,5.375){\color{Black}\qbezier(0,0)(-0.375,0.375)(-0.75,0.75)}
\put(5.39,4.375){\color{Black}\qbezier(0,0)(-0.375,0.375)(-0.75,0.75)}
\end{picture}
\end{center}
\end{figure}

\newpage
\begin{figure}
\begin{center}
{\small {\bf \CFourFigure.2}\ \   A 27-element full-length sublattice of the DCDL $\widetilde{L_{\mytinyC_{4}}}(\omega_{\mbox{\color{Red}{\tiny $2$}}})$ from \CFourFigure.1. The compression poset for this sublattice can be discerned, via \FullLengthTheorem, from \CFourFigure.1's $\widetilde{P_{\mytinyC_{4}}}(\omega_{\mbox{\color{Red}{\tiny $2$}}})$.}\\  
{\footnotesize (This lattice, denoted $\mysmallfancyQ^{(\mbox{\color{Purple}{\tiny $3$}})}$, is ``$\mysmallC_{4}$-structured'', cf.\ \S \GStructureIntroSection, and first of three ``quasi-minuscule DCDL's'' for $\myC_{4}$, cf. \S \QuasiExampleSection.)}

\vspace*{0.1in} 
\CFourFlower

\vspace*{0.1in} 
\setlength{\unitlength}{1cm}
\begin{picture}(7,13)
\put(3.85,12){\LARGE $\mybigfancyQ^{(\mbox{\color{Purple}{\small $3$}})}$}
\put(3.7,12.1){\vector(-2,-1){1.5}}
\put(4.25,11.4){\LARGE $\cong \mathbf{J}_{\mbox{\footnotesize color}}\big(\mbox{\tt \huge Q}^{(\mbox{\color{Purple}{\small $3$}})}\big)$}
\put(0,0){\TypeEboxDot{Gray}}
\put(1,1){\TypeEboxDot{Gray}}
\put(0,2){\TypeEboxDot{Gray}}
\put(2,2){\TypeEboxDot{Gray}}
\put(1,3){\TypeEboxDot{Gray}}
\put(3,3){\TypeEboxDot{Gray}}
\put(0,4){\TypeEboxDot{Gray}}
\put(2,4){\TypeEboxDot{Gray}}
\put(4,4){\TypeEboxDot{Gray}}
\put(1,5){\TypeEboxDot{Gray}}
\put(3,5){\TypeEboxDot{Gray}}
\put(5,5){\TypeEboxDot{Gray}}
\put(0,6){\TypeEboxDot{Gray}}
\put(2,6){\TypeEboxDot{Gray}}
\put(4,6){\TypeEboxDot{Gray}}
\put(1,7){\TypeEboxDot{Gray}}
\put(3,7){\TypeEboxDot{Gray}}
\put(5,7){\TypeEboxDot{Gray}}
\put(0,8){\TypeEboxDot{Gray}}
\put(2,8){\TypeEboxDot{Gray}}
\put(4,8){\TypeEboxDot{Gray}}
\put(1,9){\TypeEboxDot{Gray}}
\put(3,9){\TypeEboxDot{Gray}}
\put(0,10){\TypeEboxDot{Gray}}
\put(2,10){\TypeEboxDot{Gray}}
\put(1,11){\TypeEboxDot{Gray}}
\put(0,12){\TypeEboxDot{Gray}}
\thicklines
\put(0.625,0.375){\color{Red}\qbezier(0,0)(0.375,0.375)(0.75,0.75)}
\put(1.625,1.375){\color{Purple}\qbezier(0,0)(0.375,0.375)(0.75,0.75)}
\put(2.625,2.375){\color{Green}\qbezier(0,0)(0.375,0.375)(0.75,0.75)}
\put(3.625,3.375){\color{Purple}\qbezier(0,0)(0.375,0.375)(0.75,0.75)}
\put(4.625,4.375){\color{Red}\qbezier(0,0)(0.375,0.375)(0.75,0.75)}
%
\put(0.625,2.375){\color{Purple}\qbezier(0,0)(0.375,0.375)(0.75,0.75)}
\put(1.625,3.375){\color{Green}\qbezier(0,0)(0.375,0.375)(0.75,0.75)}
\put(2.625,4.375){\color{Purple}\qbezier(0,0)(0.375,0.375)(0.75,0.75)}
\put(3.625,5.375){\color{Red}\qbezier(0,0)(0.375,0.375)(0.75,0.75)}
\put(4.625,6.375){\color{Cyan}\qbezier(0,0)(0.375,0.375)(0.75,0.75)}
\put(0.625,4.375){\color{Green}\qbezier(0,0)(0.375,0.375)(0.75,0.75)}
\put(1.625,5.375){\color{Purple}\qbezier(0,0)(0.375,0.375)(0.75,0.75)}
\put(2.625,6.375){\color{Red}\qbezier(0,0)(0.375,0.375)(0.75,0.75)}
\put(3.625,7.375){\color{Cyan}\qbezier(0,0)(0.375,0.375)(0.75,0.75)}
\put(0.625,6.375){\color{Purple}\qbezier(0,0)(0.375,0.375)(0.75,0.75)}
\put(1.625,7.375){\color{Red}\qbezier(0,0)(0.375,0.375)(0.75,0.75)}
\put(2.625,8.375){\color{Cyan}\qbezier(0,0)(0.375,0.375)(0.75,0.75)}
\put(0.625,8.375){\color{Red}\qbezier(0,0)(0.375,0.375)(0.75,0.75)}
\put(1.625,9.375){\color{Cyan}\qbezier(0,0)(0.375,0.375)(0.75,0.75)}
\put(0.625,10.375){\color{Cyan}\qbezier(0,0)(0.375,0.375)(0.75,0.75)}
\put(1.39,11.375){\color{Red}\qbezier(0,0)(-0.375,0.375)(-0.75,0.75)}
\put(2.39,10.375){\color{Purple}\qbezier(0,0)(-0.375,0.375)(-0.75,0.75)}
\put(3.39,9.375){\color{Green}\qbezier(0,0)(-0.375,0.375)(-0.75,0.75)}
\put(4.39,8.375){\color{Purple}\qbezier(0,0)(-0.375,0.375)(-0.75,0.75)}
\put(5.39,7.375){\color{Red}\qbezier(0,0)(-0.375,0.375)(-0.75,0.75)}
%
\put(1.39,9.375){\color{Purple}\qbezier(0,0)(-0.375,0.375)(-0.75,0.75)}
\put(2.39,8.375){\color{Green}\qbezier(0,0)(-0.375,0.375)(-0.75,0.75)}
\put(3.39,7.375){\color{Purple}\qbezier(0,0)(-0.375,0.375)(-0.75,0.75)}
\put(4.39,6.375){\color{Red}\qbezier(0,0)(-0.375,0.375)(-0.75,0.75)}
\put(5.39,5.375){\color{Cyan}\qbezier(0,0)(-0.375,0.375)(-0.75,0.75)}
\put(1.39,7.375){\color{Green}\qbezier(0,0)(-0.375,0.375)(-0.75,0.75)}
\put(2.39,6.375){\color{Purple}\qbezier(0,0)(-0.375,0.375)(-0.75,0.75)}
\put(3.39,5.375){\color{Red}\qbezier(0,0)(-0.375,0.375)(-0.75,0.75)}
\put(4.39,4.375){\color{Cyan}\qbezier(0,0)(-0.375,0.375)(-0.75,0.75)}
\put(1.39,5.375){\color{Purple}\qbezier(0,0)(-0.375,0.375)(-0.75,0.75)}
\put(2.39,4.375){\color{Red}\qbezier(0,0)(-0.375,0.375)(-0.75,0.75)}
\put(3.39,3.375){\color{Cyan}\qbezier(0,0)(-0.375,0.375)(-0.75,0.75)}
\put(1.39,3.375){\color{Red}\qbezier(0,0)(-0.375,0.375)(-0.75,0.75)}
\put(2.39,2.375){\color{Cyan}\qbezier(0,0)(-0.375,0.375)(-0.75,0.75)}
\put(1.39,1.375){\color{Cyan}\qbezier(0,0)(-0.375,0.375)(-0.75,0.75)}
\put(0.9,0.65){\color{Red}{\em 2}}
\put(0.9,2.65){\color{Purple}{\em 3}}
\put(1.9,1.65){\color{Purple}{\em 3}}
\put(0.9,4.65){\color{Green}{\em 4}}
\put(1.9,3.65){\color{Green}{\em 4}}
\put(2.9,2.65){\color{Green}{\em 4}}
\put(0.9,6.65){\color{Purple}{\em 3}}
\put(1.9,5.65){\color{Purple}{\em 3}}
\put(2.9,4.65){\color{Purple}{\em 3}}
\put(3.9,3.65){\color{Purple}{\em 3}}
\put(0.9,8.65){\color{Red}{\em 2}}
\put(1.9,7.65){\color{Red}{\em 2}}
\put(2.9,6.65){\color{Red}{\em 2}}
\put(3.9,5.65){\color{Red}{\em 2}}
\put(4.9,4.65){\color{Red}{\em 2}}
\put(0.9,10.65){\color{Cyan}{\em 1}}
\put(1.9,9.65){\color{Cyan}{\em 1}}
\put(2.9,8.65){\color{Cyan}{\em 1}}
\put(3.9,7.65){\color{Cyan}{\em 1}}
\put(4.9,6.65){\color{Cyan}{\em 1}}
%
\put(0.9,11.65){\color{Red}{\em 2}}
\put(0.9,9.65){\color{Purple}{\em 3}}
\put(1.9,10.65){\color{Purple}{\em 3}}
\put(0.9,7.65){\color{Green}{\em 4}}
\put(1.9,8.65){\color{Green}{\em 4}}
\put(2.9,9.65){\color{Green}{\em 4}}
\put(0.9,5.65){\color{Purple}{\em 3}}
\put(1.9,6.65){\color{Purple}{\em 3}}
\put(2.9,7.65){\color{Purple}{\em 3}}
\put(3.9,8.65){\color{Purple}{\em 3}}
\put(0.9,3.65){\color{Red}{\em 2}}
\put(1.9,4.65){\color{Red}{\em 2}}
\put(2.9,5.65){\color{Red}{\em 2}}
\put(3.9,6.65){\color{Red}{\em 2}}
\put(4.9,7.65){\color{Red}{\em 2}}
\put(0.9,1.65){\color{Cyan}{\em 1}}
\put(1.9,2.65){\color{Cyan}{\em 1}}
\put(2.9,3.65){\color{Cyan}{\em 1}}
\put(3.9,4.65){\color{Cyan}{\em 1}}
\put(4.9,5.65){\color{Cyan}{\em 1}}
\put(0.8,6.175){\small $\lsem {\color{Purple}3},{\color{Purple}3} \rsem$}
\put(2.8,6.175){\small $\lsem {\color{Red}2},{\color{Purple}3} \rsem$}
\put(4.8,6.175){\small $\lsem {\color{Cyan}1},{\color{Red}2} \rsem$}
\end{picture}
\begin{picture}(7,13)
\put(4.5,10.5){\LARGE $\mbox{\tt \huge Q}^{(\mbox{\color{Purple}{\small $3$}})}$}
\put(4.35,10.6){\vector(-2,-1){1.5}}
\put(4.9,9.9){\LARGE $\cong \mathbf{j}_{\mbox{\footnotesize color}}\big(\mybigfancyQ^{(\mbox{\color{Purple}{\small $3$}})}\big)$}
\put(0,9){\TypeEboxDot{Cyan}}
\put(0.1,9){\color{Cyan}{\em \footnotesize 1}}
\put(1,10){\TypeEboxDot{Red}}
\put(1,8){\TypeEboxDot{Red}}
\put(1.7,10.3){\color{Red}{\em \footnotesize 2}}
\put(1.1,8){\color{Red}{\em \footnotesize 2}}
\put(2,9){\TypeEboxDot{Purple}}
\put(2,7){\TypeEboxDot{Purple}}
\put(2.7,9.3){\color{Purple}{\em \footnotesize 3}}
\put(2.1,7){\color{Purple}{\em \footnotesize 3}}
\put(3,8){\TypeEboxDot{Green}}
\put(3,6){\TypeEboxDot{Green}}
\put(3.7,8.3){\color{Green}{\em \footnotesize 4}}
\put(3.1,6){\color{Green}{\em \footnotesize 4}}
\put(4,7){\TypeEboxDot{Purple}}
\put(4,5){\TypeEboxDot{Purple}}
\put(4.7,7.3){\color{Purple}{\em \footnotesize 3}}
\put(4.1,5){\color{Purple}{\em \footnotesize 3}}
\put(5,6){\TypeEboxDot{Red}}
\put(5,4){\TypeEboxDot{Red}}
\put(5.7,6.3){\color{Red}{\em \footnotesize 2}}
\put(5.1,4){\color{Red}{\em \footnotesize 2}}
\put(6,5){\TypeEboxDot{Cyan}}
\put(6.7,5.3){\color{Cyan}{\em \footnotesize 1}}
\thicklines
\put(5.625,4.375){\color{Black}\qbezier(0,0)(0.375,0.375)(0.75,0.75)}
\put(4.625,5.375){\color{Black}\qbezier(0,0)(0.375,0.375)(0.75,0.75)}
\put(3.625,6.375){\color{Black}\qbezier(0,0)(0.375,0.375)(0.75,0.75)}
\put(2.625,7.375){\color{Black}\qbezier(0,0)(0.375,0.375)(0.75,0.75)}
\put(1.625,8.375){\color{Black}\qbezier(0,0)(0.375,0.375)(0.75,0.75)}
\put(0.625,9.375){\color{Black}\qbezier(0,0)(0.375,0.375)(0.75,0.75)}
\put(2.39,9.375){\color{Black}\qbezier(0,0)(-0.375,0.375)(-0.75,0.75)}
\put(3.39,8.375){\color{Black}\qbezier(0,0)(-0.375,0.375)(-0.75,0.75)}
\put(4.39,7.375){\color{Black}\qbezier(0,0)(-0.375,0.375)(-0.75,0.75)}
\put(5.39,6.375){\color{Black}\qbezier(0,0)(-0.375,0.375)(-0.75,0.75)}
\put(6.39,5.375){\color{Black}\qbezier(0,0)(-0.375,0.375)(-0.75,0.75)}
\put(1.39,8.375){\color{Black}\qbezier(0,0)(-0.375,0.375)(-0.75,0.75)}
\put(2.39,7.375){\color{Black}\qbezier(0,0)(-0.375,0.375)(-0.75,0.75)}
\put(3.39,6.375){\color{Black}\qbezier(0,0)(-0.375,0.375)(-0.75,0.75)}
\put(4.39,5.375){\color{Black}\qbezier(0,0)(-0.375,0.375)(-0.75,0.75)}
\put(5.39,4.375){\color{Black}\qbezier(0,0)(-0.375,0.375)(-0.75,0.75)}
\put(6.325,5.25){\color{Black}\qbezier(0,0)(-1,0)(-5.675,3.925)}
\put(4.6275,6.075){\vector(-3,2){0.2}}
\end{picture}
\end{center}
\end{figure}

\newpage
\begin{figure}
\begin{center}
{\small {\bf \CFourFigure.3}\ \  Another 27-element full-length sublattice of the DCDL $\widetilde{L_{\mytinyC_{4}}}(\omega_{\mbox{\color{Red}{\tiny $2$}}})$ from \CFourFigure.1. The compression poset for this sublattice can be discerned, via \FullLengthTheorem, from \CFourFigure.1's $\widetilde{P_{\mytinyC_{4}}}(\omega_{\mbox{\color{Red}{\tiny $2$}}})$.}\\
{\footnotesize (This lattice, denoted $\mysmallfancyQ^{(\mbox{\color{Cyan}{\tiny $1$}})}$, is ``$\mysmallC_{4}$-structured'', cf.\ \S \GStructureIntroSection, and second of three ``quasi-minuscule DCML's'' for $\myC_{4}$, cf. \S \QuasiExampleSection.)}

\vspace*{0.1in} 
\CFourFlower

\vspace*{0.1in} 
\setlength{\unitlength}{1cm}
\begin{picture}(7,13)
\put(3.85,12){\LARGE $\mybigfancyQ^{(\mbox{\color{Cyan}{\small $1$}})}$}
\put(3.7,12.1){\vector(-2,-1){1.5}}
\put(4.25,11.4){\LARGE $\cong \mathbf{J}_{\mbox{\footnotesize color}}\big(\mbox{\tt \huge Q}^{(\mbox{\color{Cyan}{\small $1$}})}\big)$}
\put(0,0){\TypeEboxDot{Gray}}
\put(1,1){\TypeEboxDot{Gray}}
\put(0,2){\TypeEboxDot{Gray}}
\put(2,2){\TypeEboxDot{Gray}}
\put(1,3){\TypeEboxDot{Gray}}
\put(3,3){\TypeEboxDot{Gray}}
\put(0,4){\TypeEboxDot{Gray}}
\put(2,4){\TypeEboxDot{Gray}}
\put(4,4){\TypeEboxDot{Gray}}
\put(1,5){\TypeEboxDot{Gray}}
\put(3,5){\TypeEboxDot{Gray}}
\put(5,5){\TypeEboxDot{Gray}}
\put(2,6){\TypeEboxDot{Gray}}
\put(4,6){\TypeEboxDot{Gray}}
\put(6,6){\TypeEboxDot{Gray}}
\put(1,7){\TypeEboxDot{Gray}}
\put(3,7){\TypeEboxDot{Gray}}
\put(5,7){\TypeEboxDot{Gray}}
\put(0,8){\TypeEboxDot{Gray}}
\put(2,8){\TypeEboxDot{Gray}}
\put(4,8){\TypeEboxDot{Gray}}
\put(1,9){\TypeEboxDot{Gray}}
\put(3,9){\TypeEboxDot{Gray}}
\put(0,10){\TypeEboxDot{Gray}}
\put(2,10){\TypeEboxDot{Gray}}
\put(1,11){\TypeEboxDot{Gray}}
\put(0,12){\TypeEboxDot{Gray}}
\thicklines
\put(0.625,0.375){\color{Red}\qbezier(0,0)(0.375,0.375)(0.75,0.75)}
\put(1.625,1.375){\color{Purple}\qbezier(0,0)(0.375,0.375)(0.75,0.75)}
\put(2.625,2.375){\color{Green}\qbezier(0,0)(0.375,0.375)(0.75,0.75)}
\put(3.625,3.375){\color{Purple}\qbezier(0,0)(0.375,0.375)(0.75,0.75)}
\put(4.625,4.375){\color{Red}\qbezier(0,0)(0.375,0.375)(0.75,0.75)}
\put(5.625,5.375){\color{Cyan}\qbezier(0,0)(0.375,0.375)(0.75,0.75)}
\put(0.625,2.375){\color{Purple}\qbezier(0,0)(0.375,0.375)(0.75,0.75)}
\put(1.625,3.375){\color{Green}\qbezier(0,0)(0.375,0.375)(0.75,0.75)}
\put(2.625,4.375){\color{Purple}\qbezier(0,0)(0.375,0.375)(0.75,0.75)}
\put(3.625,5.375){\color{Red}\qbezier(0,0)(0.375,0.375)(0.75,0.75)}
\put(4.625,6.375){\color{Cyan}\qbezier(0,0)(0.375,0.375)(0.75,0.75)}
\put(0.625,4.375){\color{Green}\qbezier(0,0)(0.375,0.375)(0.75,0.75)}
\put(1.625,5.375){\color{Purple}\qbezier(0,0)(0.375,0.375)(0.75,0.75)}
\put(2.625,6.375){\color{Red}\qbezier(0,0)(0.375,0.375)(0.75,0.75)}
\put(3.625,7.375){\color{Cyan}\qbezier(0,0)(0.375,0.375)(0.75,0.75)}
%
\put(1.625,7.375){\color{Red}\qbezier(0,0)(0.375,0.375)(0.75,0.75)}
\put(2.625,8.375){\color{Cyan}\qbezier(0,0)(0.375,0.375)(0.75,0.75)}
\put(0.625,8.375){\color{Red}\qbezier(0,0)(0.375,0.375)(0.75,0.75)}
\put(1.625,9.375){\color{Cyan}\qbezier(0,0)(0.375,0.375)(0.75,0.75)}
\put(0.625,10.375){\color{Cyan}\qbezier(0,0)(0.375,0.375)(0.75,0.75)}
\put(1.39,11.375){\color{Red}\qbezier(0,0)(-0.375,0.375)(-0.75,0.75)}
\put(2.39,10.375){\color{Purple}\qbezier(0,0)(-0.375,0.375)(-0.75,0.75)}
\put(3.39,9.375){\color{Green}\qbezier(0,0)(-0.375,0.375)(-0.75,0.75)}
\put(4.39,8.375){\color{Purple}\qbezier(0,0)(-0.375,0.375)(-0.75,0.75)}
\put(5.39,7.375){\color{Red}\qbezier(0,0)(-0.375,0.375)(-0.75,0.75)}
\put(6.39,6.375){\color{Cyan}\qbezier(0,0)(-0.375,0.375)(-0.75,0.75)}
\put(1.39,9.375){\color{Purple}\qbezier(0,0)(-0.375,0.375)(-0.75,0.75)}
\put(2.39,8.375){\color{Green}\qbezier(0,0)(-0.375,0.375)(-0.75,0.75)}
\put(3.39,7.375){\color{Purple}\qbezier(0,0)(-0.375,0.375)(-0.75,0.75)}
\put(4.39,6.375){\color{Red}\qbezier(0,0)(-0.375,0.375)(-0.75,0.75)}
\put(5.39,5.375){\color{Cyan}\qbezier(0,0)(-0.375,0.375)(-0.75,0.75)}
\put(1.39,7.375){\color{Green}\qbezier(0,0)(-0.375,0.375)(-0.75,0.75)}
\put(2.39,6.375){\color{Purple}\qbezier(0,0)(-0.375,0.375)(-0.75,0.75)}
\put(3.39,5.375){\color{Red}\qbezier(0,0)(-0.375,0.375)(-0.75,0.75)}
\put(4.39,4.375){\color{Cyan}\qbezier(0,0)(-0.375,0.375)(-0.75,0.75)}
%
\put(2.39,4.375){\color{Red}\qbezier(0,0)(-0.375,0.375)(-0.75,0.75)}
\put(3.39,3.375){\color{Cyan}\qbezier(0,0)(-0.375,0.375)(-0.75,0.75)}
\put(1.39,3.375){\color{Red}\qbezier(0,0)(-0.375,0.375)(-0.75,0.75)}
\put(2.39,2.375){\color{Cyan}\qbezier(0,0)(-0.375,0.375)(-0.75,0.75)}
\put(1.39,1.375){\color{Cyan}\qbezier(0,0)(-0.375,0.375)(-0.75,0.75)}
\put(0.9,0.65){\color{Red}{\em 2}}
\put(0.9,2.65){\color{Purple}{\em 3}}
\put(1.9,1.65){\color{Purple}{\em 3}}
\put(0.9,4.65){\color{Green}{\em 4}}
\put(1.9,3.65){\color{Green}{\em 4}}
\put(2.9,2.65){\color{Green}{\em 4}}
%
\put(1.9,5.65){\color{Purple}{\em 3}}
\put(2.9,4.65){\color{Purple}{\em 3}}
\put(3.9,3.65){\color{Purple}{\em 3}}
\put(0.9,8.65){\color{Red}{\em 2}}
\put(1.9,7.65){\color{Red}{\em 2}}
\put(2.9,6.65){\color{Red}{\em 2}}
\put(3.9,5.65){\color{Red}{\em 2}}
\put(4.9,4.65){\color{Red}{\em 2}}
\put(0.9,10.65){\color{Cyan}{\em 1}}
\put(1.9,9.65){\color{Cyan}{\em 1}}
\put(2.9,8.65){\color{Cyan}{\em 1}}
\put(3.9,7.65){\color{Cyan}{\em 1}}
\put(4.9,6.65){\color{Cyan}{\em 1}}
\put(5.9,5.65){\color{Cyan}{\em 1}}
\put(0.9,11.65){\color{Red}{\em 2}}
\put(0.9,9.65){\color{Purple}{\em 3}}
\put(1.9,10.65){\color{Purple}{\em 3}}
\put(0.9,7.65){\color{Green}{\em 4}}
\put(1.9,8.65){\color{Green}{\em 4}}
\put(2.9,9.65){\color{Green}{\em 4}}
%
\put(1.9,6.65){\color{Purple}{\em 3}}
\put(2.9,7.65){\color{Purple}{\em 3}}
\put(3.9,8.65){\color{Purple}{\em 3}}
\put(0.9,3.65){\color{Red}{\em 2}}
\put(1.9,4.65){\color{Red}{\em 2}}
\put(2.9,5.65){\color{Red}{\em 2}}
\put(3.9,6.65){\color{Red}{\em 2}}
\put(4.9,7.65){\color{Red}{\em 2}}
\put(0.9,1.65){\color{Cyan}{\em 1}}
\put(1.9,2.65){\color{Cyan}{\em 1}}
\put(2.9,3.65){\color{Cyan}{\em 1}}
\put(3.9,4.65){\color{Cyan}{\em 1}}
\put(4.9,5.65){\color{Cyan}{\em 1}}
\put(5.9,6.65){\color{Cyan}{\em 1}}
\put(1.4,6.175){\small $\lsem {\color{Red}2},{\color{Purple}3} \rsem$}
\put(3.4,6.175){\small $\lsem {\color{Cyan}1},{\color{Red}2} \rsem$}
\put(5.4,6.175){\small $\lsem {\color{Cyan}1},{\color{Cyan}1} \rsem$}
\end{picture}
\begin{picture}(7,13)
\put(4.5,10.5){\LARGE $\mbox{\tt \huge Q}^{(\mbox{\color{Cyan}{\small $1$}})}$}
\put(4.35,10.6){\vector(-2,-1){1.5}}
\put(4.9,9.9){\LARGE $\cong \mathbf{j}_{\mbox{\footnotesize color}}\big(\mybigfancyQ^{(\mbox{\color{Cyan}{\small $1$}})}\big)$}
\put(0,9){\TypeEboxDot{Cyan}}
\put(0.1,9){\color{Cyan}{\em \footnotesize 1}}
\put(1,10){\TypeEboxDot{Red}}
\put(1,8){\TypeEboxDot{Red}}
\put(1.7,10.3){\color{Red}{\em \footnotesize 2}}
\put(1.1,8){\color{Red}{\em \footnotesize 2}}
\put(2,9){\TypeEboxDot{Purple}}
\put(2,7){\TypeEboxDot{Purple}}
\put(2.7,9.3){\color{Purple}{\em \footnotesize 3}}
\put(2.1,7){\color{Purple}{\em \footnotesize 3}}
\put(3,8){\TypeEboxDot{Green}}
\put(3,6){\TypeEboxDot{Green}}
\put(3.7,8.3){\color{Green}{\em \footnotesize 4}}
\put(3.1,6){\color{Green}{\em \footnotesize 4}}
\put(4,7){\TypeEboxDot{Purple}}
\put(4,5){\TypeEboxDot{Purple}}
\put(4.7,7.3){\color{Purple}{\em \footnotesize 3}}
\put(4.1,5){\color{Purple}{\em \footnotesize 3}}
\put(5,6){\TypeEboxDot{Red}}
\put(5,4){\TypeEboxDot{Red}}
\put(5.7,6.3){\color{Red}{\em \footnotesize 2}}
\put(5.1,4){\color{Red}{\em \footnotesize 2}}
\put(6,5){\TypeEboxDot{Cyan}}
\put(6.7,5.3){\color{Cyan}{\em \footnotesize 1}}
\thicklines
\put(5.625,4.375){\color{Black}\qbezier(0,0)(0.375,0.375)(0.75,0.75)}
\put(4.625,5.375){\color{Black}\qbezier(0,0)(0.375,0.375)(0.75,0.75)}
\put(1.625,8.375){\color{Black}\qbezier(0,0)(0.375,0.375)(0.75,0.75)}
\put(0.625,9.375){\color{Black}\qbezier(0,0)(0.375,0.375)(0.75,0.75)}
\put(2.39,9.375){\color{Black}\qbezier(0,0)(-0.375,0.375)(-0.75,0.75)}
\put(3.39,8.375){\color{Black}\qbezier(0,0)(-0.375,0.375)(-0.75,0.75)}
\put(4.39,7.375){\color{Black}\qbezier(0,0)(-0.375,0.375)(-0.75,0.75)}
\put(5.39,6.375){\color{Black}\qbezier(0,0)(-0.375,0.375)(-0.75,0.75)}
\put(6.39,5.375){\color{Black}\qbezier(0,0)(-0.375,0.375)(-0.75,0.75)}
\put(1.39,8.375){\color{Black}\qbezier(0,0)(-0.375,0.375)(-0.75,0.75)}
\put(2.39,7.375){\color{Black}\qbezier(0,0)(-0.375,0.375)(-0.75,0.75)}
\put(3.39,6.375){\color{Black}\qbezier(0,0)(-0.375,0.375)(-0.75,0.75)}
\put(4.39,5.375){\color{Black}\qbezier(0,0)(-0.375,0.375)(-0.75,0.75)}
\put(5.39,4.375){\color{Black}\qbezier(0,0)(-0.375,0.375)(-0.75,0.75)}
\put(4.375,7.125){\color{Black}\qbezier(0,0)(-0.875,-0.3)(-1.75,0)}
\put(3.5,6.975){\vector(1,0){0.2}}
\end{picture}
\end{center}
\end{figure}

\newpage
\begin{figure}
\begin{center}
{{\bf \CFourFigure.4}\ \  Below is a 27-element diamond-colored modular lattice\\ that is not a sublattice of $\widetilde{L_{\mytinyC_{4}}}(\omega_{\mbox{\color{Red}{\tiny $2$}}})$ from \CFourFigure.1.}\\ 
{\footnotesize (Since this modular lattice is not distributive, a companion compression poset of ideals is not available to us.\\ 
This lattice, denoted $\mysmallfancyQ^{(\mbox{\color{Red}{\tiny $2$}})}$, is ``$\mysmallC_{4}$-structured'', cf.\ \S \GStructureIntroSection, and third of three ``quasi-minuscule DCML's'' for $\myC_{4}$, cf. \S \QuasiExampleSection.\\  The other two are $\mysmallfancyQ^{(\mbox{\color{Cyan}{\tiny $1$}})}$ and $\mysmallfancyQ^{(\mbox{\color{Purple}{\tiny $3$}})}$ from \CFourFigure.2/3.)}

\vspace*{0.1in} 
\CFourFlower

\vspace*{0.1in} 
\setlength{\unitlength}{1cm}
\begin{picture}(7,13)
\put(5.5,11){\LARGE $\mybigfancyQ^{(\mbox{\color{Red}{\small $2$}})}$}
\put(5.35,11.1){\vector(-2,-1){1.5}}
\put(0,0){\TypeEboxDot{Gray}}
\put(1,1){\TypeEboxDot{Gray}}
\put(0,2){\TypeEboxDot{Gray}}
\put(2,2){\TypeEboxDot{Gray}}
\put(1,3){\TypeEboxDot{Gray}}
\put(3,3){\TypeEboxDot{Gray}}
\put(0,4){\TypeEboxDot{Gray}}
\put(2,4){\TypeEboxDot{Gray}}
\put(4,4){\TypeEboxDot{Gray}}
\put(1,5){\TypeEboxDot{Gray}}
\put(3,5){\TypeEboxDot{Gray}}
\put(5,5){\TypeEboxDot{Gray}}
\put(2,6){\TypeEboxDot{Gray}}
\put(4,6){\TypeEboxDot{Gray}}
\put(1,7){\TypeEboxDot{Gray}}
\put(3,7){\TypeEboxDot{Gray}}
\put(5,7){\TypeEboxDot{Gray}}
\put(0,8){\TypeEboxDot{Gray}}
\put(2,8){\TypeEboxDot{Gray}}
\put(4,8){\TypeEboxDot{Gray}}
\put(1,9){\TypeEboxDot{Gray}}
\put(3,9){\TypeEboxDot{Gray}}
\put(0,10){\TypeEboxDot{Gray}}
\put(2,10){\TypeEboxDot{Gray}}
\put(1,11){\TypeEboxDot{Gray}}
\put(0,12){\TypeEboxDot{Gray}}
\thicklines
\put(0.625,0.375){\color{Red}\qbezier(0,0)(0.375,0.375)(0.75,0.75)}
\put(1.625,1.375){\color{Purple}\qbezier(0,0)(0.375,0.375)(0.75,0.75)}
\put(2.625,2.375){\color{Green}\qbezier(0,0)(0.375,0.375)(0.75,0.75)}
\put(3.625,3.375){\color{Purple}\qbezier(0,0)(0.375,0.375)(0.75,0.75)}
\put(4.625,4.375){\color{Red}\qbezier(0,0)(0.375,0.375)(0.75,0.75)}
%
\put(0.625,2.375){\color{Purple}\qbezier(0,0)(0.375,0.375)(0.75,0.75)}
\put(1.625,3.375){\color{Green}\qbezier(0,0)(0.375,0.375)(0.75,0.75)}
\put(2.625,4.375){\color{Purple}\qbezier(0,0)(0.375,0.375)(0.75,0.75)}
\put(3.625,5.375){\color{Red}\qbezier(0,0)(0.375,0.375)(0.75,0.75)}
\put(4.625,6.375){\color{Cyan}\qbezier(0,0)(0.375,0.375)(0.75,0.75)}
\put(0.625,4.375){\color{Green}\qbezier(0,0)(0.375,0.375)(0.75,0.75)}
\put(1.625,5.375){\color{Purple}\qbezier(0,0)(0.375,0.375)(0.75,0.75)}
\put(2.625,6.375){\color{Red}\qbezier(0,0)(0.375,0.375)(0.75,0.75)}
\put(3.625,7.375){\color{Cyan}\qbezier(0,0)(0.375,0.375)(0.75,0.75)}
%
\put(1.625,7.375){\color{Red}\qbezier(0,0)(0.375,0.375)(0.75,0.75)}
\put(2.625,8.375){\color{Cyan}\qbezier(0,0)(0.375,0.375)(0.75,0.75)}
\put(0.625,8.375){\color{Red}\qbezier(0,0)(0.375,0.375)(0.75,0.75)}
\put(1.625,9.375){\color{Cyan}\qbezier(0,0)(0.375,0.375)(0.75,0.75)}
\put(0.625,10.375){\color{Cyan}\qbezier(0,0)(0.375,0.375)(0.75,0.75)}
\put(1.39,11.375){\color{Red}\qbezier(0,0)(-0.375,0.375)(-0.75,0.75)}
\put(2.39,10.375){\color{Purple}\qbezier(0,0)(-0.375,0.375)(-0.75,0.75)}
\put(3.39,9.375){\color{Green}\qbezier(0,0)(-0.375,0.375)(-0.75,0.75)}
\put(4.39,8.375){\color{Purple}\qbezier(0,0)(-0.375,0.375)(-0.75,0.75)}
\put(5.39,7.375){\color{Red}\qbezier(0,0)(-0.375,0.375)(-0.75,0.75)}
%
\put(1.39,9.375){\color{Purple}\qbezier(0,0)(-0.375,0.375)(-0.75,0.75)}
\put(2.39,8.375){\color{Green}\qbezier(0,0)(-0.375,0.375)(-0.75,0.75)}
\put(3.39,7.375){\color{Purple}\qbezier(0,0)(-0.375,0.375)(-0.75,0.75)}
\put(4.39,6.375){\color{Red}\qbezier(0,0)(-0.375,0.375)(-0.75,0.75)}
\put(5.39,5.375){\color{Cyan}\qbezier(0,0)(-0.375,0.375)(-0.75,0.75)}
\put(1.39,7.375){\color{Green}\qbezier(0,0)(-0.375,0.375)(-0.75,0.75)}
\put(2.39,6.375){\color{Purple}\qbezier(0,0)(-0.375,0.375)(-0.75,0.75)}
\put(3.39,5.375){\color{Red}\qbezier(0,0)(-0.375,0.375)(-0.75,0.75)}
\put(4.39,4.375){\color{Cyan}\qbezier(0,0)(-0.375,0.375)(-0.75,0.75)}
%
\put(2.39,4.375){\color{Red}\qbezier(0,0)(-0.375,0.375)(-0.75,0.75)}
\put(3.39,3.375){\color{Cyan}\qbezier(0,0)(-0.375,0.375)(-0.75,0.75)}
\put(1.39,3.375){\color{Red}\qbezier(0,0)(-0.375,0.375)(-0.75,0.75)}
\put(2.39,2.375){\color{Cyan}\qbezier(0,0)(-0.375,0.375)(-0.75,0.75)}
\put(1.39,1.375){\color{Cyan}\qbezier(0,0)(-0.375,0.375)(-0.75,0.75)}
\put(0.9,0.65){\color{Red}{\em 2}}
\put(0.9,2.65){\color{Purple}{\em 3}}
\put(1.9,1.65){\color{Purple}{\em 3}}
\put(0.9,4.65){\color{Green}{\em 4}}
\put(1.9,3.65){\color{Green}{\em 4}}
\put(2.9,2.65){\color{Green}{\em 4}}
%
\put(1.9,5.65){\color{Purple}{\em 3}}
\put(2.9,4.65){\color{Purple}{\em 3}}
\put(3.9,3.65){\color{Purple}{\em 3}}
\put(0.9,8.65){\color{Red}{\em 2}}
\put(1.9,7.65){\color{Red}{\em 2}}
\put(2.9,6.65){\color{Red}{\em 2}}
\put(3.9,5.65){\color{Red}{\em 2}}
\put(4.9,4.65){\color{Red}{\em 2}}
\put(0.9,10.65){\color{Cyan}{\em 1}}
\put(1.9,9.65){\color{Cyan}{\em 1}}
\put(2.9,8.65){\color{Cyan}{\em 1}}
\put(3.9,7.65){\color{Cyan}{\em 1}}
\put(4.9,6.65){\color{Cyan}{\em 1}}
%
\put(0.9,11.65){\color{Red}{\em 2}}
\put(0.9,9.65){\color{Purple}{\em 3}}
\put(1.9,10.65){\color{Purple}{\em 3}}
\put(0.9,7.65){\color{Green}{\em 4}}
\put(1.9,8.65){\color{Green}{\em 4}}
\put(2.9,9.65){\color{Green}{\em 4}}
%
\put(1.9,6.65){\color{Purple}{\em 3}}
\put(2.9,7.65){\color{Purple}{\em 3}}
\put(3.9,8.65){\color{Purple}{\em 3}}
\put(0.9,3.65){\color{Red}{\em 2}}
\put(1.9,4.65){\color{Red}{\em 2}}
\put(2.9,5.65){\color{Red}{\em 2}}
\put(3.9,6.65){\color{Red}{\em 2}}
\put(4.9,7.65){\color{Red}{\em 2}}
\put(0.9,1.65){\color{Cyan}{\em 1}}
\put(1.9,2.65){\color{Cyan}{\em 1}}
\put(2.9,3.65){\color{Cyan}{\em 1}}
\put(3.9,4.65){\color{Cyan}{\em 1}}
\put(4.9,5.65){\color{Cyan}{\em 1}}
\put(3,6){\TypeEboxDot{Gray}}
\put(3.4,6.65){\color{Red}{\em 2}}
\put(3.51,6.425){\color{Red}\qbezier(0,0)(0,0.325)(0,0.65)}
\put(3.4,5.65){\color{Red}{\em 2}}
\put(3.51,5.425){\color{Red}\qbezier(0,0)(0,0.325)(0,0.65)}
\put(1.4,6.175){\small $\lsem {\color{Red}2},{\color{Purple}3} \rsem$}
\put(4.8,6.175){\small $\lsem {\color{Cyan}1},{\color{Red}2} \rsem$}
\put(3.675,6.185){\tiny $\lsem {\color{Red}2},{\color{Red}2} \rsem$}
\end{picture}
\end{center}
\end{figure}

\newcommand{\FFourFlower}{
\setlength{\unitlength}{0.5cm}
\begin{picture}(12,1)
\put(-2.5,-0.25){\LARGE $\myF_{4}$}
\put(-0.9,0){\vector(1,0){1.2}}
{\color{Cyan}\put(1,0){\circle*{0.7}}}
\put(1,0){\circle{0.7}}
{\color{Red}\put(4,0){\circle*{0.7}}}
\put(4,0){\circle{0.7}}
{\color{Purple}\put(7,0){\circle*{0.7}}}
\put(7,0){\circle{0.7}}
{\color{Green}\put(10,0){\circle*{0.7}}}
\put(10,0){\circle{0.7}}
\thicklines
\put(5,0){\vector(1,0){0.2}}
\put(5.35,0){\vector(1,0){0.2}}
\put(6,0){\vector(-1,0){0.2}}
\put(1.35,0){\line(1,0){2.3}}
\put(4.35,0){\line(1,0){2.3}}
\put(7.35,0){\line(1,0){2.3}}
\put(0.8,0.6){\textcolor{Cyan}{\em 1}}
\put(3.8,0.6){\textcolor{Red}{\em 2}}
\put(6.8,0.6){\textcolor{Purple}{\em 3}}
\put(9.8,0.6){\textcolor{Green}{\em 4}}
\end{picture}
}

\begin{figure}
\begin{center}
{\bf \FFourFigure.1}\ \  The 27-element diamond-colored distributive lattice depicted below is denoted `$\widetilde{L_{\mytinyF_{4}}}(\omega_{\mbox{\color{Green}{\tiny $4$}}})$'. 
It can be viewed as a recoloring of the DCDL's $L_{\mytinyE_{6}}(\omega_{1})$ and $L_{\mytinyE_{6}}(\omega_{6})$ from \ESixLatticeFigure.\\ 
{\footnotesize (This lattice is, in the language of \S \GStructureIntroSection, `$\mysmallF_{4}$-structured', where $\mysmallF_{4}$ is the graph}\\ 
{\footnotesize  immediately below that also appears in the list given in \IEGGraphFigure.)}

\vspace*{0.1in} 
\FFourFlower

\vspace*{0.1in} 
\setlength{\unitlength}{1cm}
\begin{picture}(6,17)
\put(-1.9,12.2){\LARGE $\widetilde{L_{\mysmallF_{4}}}(\omega_{\mbox{\color{Green}{\small $4$}}})$}
\put(0.1,12.6){\vector(2,1){2}}
\put(-1.5,11.2){\LARGE $\cong \mathbf{J}_{\mbox{\footnotesize color}}\big(\widetilde{P_{\mysmallF_{4}}}(\omega_{\mbox{\color{Green}{\small $4$}}})\big)$}
\put(0,16){\TypeEboxDot{Gray}}
\put(0,8){\TypeEboxDot{Gray}}
\put(0,0){\TypeEboxDot{Gray}}
\put(1,15){\TypeEboxDot{Gray}}
\put(1,9){\TypeEboxDot{Gray}}
\put(1,7){\TypeEboxDot{Gray}}
\put(1,1){\TypeEboxDot{Gray}}
\put(2,14){\TypeEboxDot{Gray}}
\put(2,10){\TypeEboxDot{Gray}}
\put(2,8){\TypeEboxDot{Gray}}
\put(2,6){\TypeEboxDot{Gray}}
\put(2,2){\TypeEboxDot{Gray}}
\put(3,13){\TypeEboxDot{Gray}}
\put(3,11){\TypeEboxDot{Gray}}
\put(3,9){\TypeEboxDot{Gray}}
\put(3,7){\TypeEboxDot{Gray}}
\put(3,5){\TypeEboxDot{Gray}}
\put(3,3){\TypeEboxDot{Gray}}
\put(2,12){\TypeEboxDot{Gray}}
\put(2,4){\TypeEboxDot{Gray}}
\put(4,12){\TypeEboxDot{Gray}}
\put(4,10){\TypeEboxDot{Gray}}
\put(4,6){\TypeEboxDot{Gray}}
\put(4,4){\TypeEboxDot{Gray}}
\put(5,11){\TypeEboxDot{Gray}}
\put(5,5){\TypeEboxDot{Gray}}
\thicklines
\put(0.625,0.375){\color{Green}\qbezier(0,0)(0.375,0.375)(0.75,0.75)}
\put(0.9,0.65){\color{Green}{\em 4}}
\put(1.625,1.375){\color{Purple}\qbezier(0,0)(0.375,0.375)(0.75,0.75)}
\put(1.9,1.65){\color{Purple}{\em 3}}
\put(2.625,2.375){\color{Red}\qbezier(0,0)(0.375,0.375)(0.75,0.75)}
\put(2.9,2.65){\color{Red}{\em 2}}
\put(3.625,3.375){\color{Purple}\qbezier(0,0)(0.375,0.375)(0.75,0.75)}
\put(3.9,3.65){\color{Purple}{\em 3}}
\put(4.625,4.375){\color{Green}\qbezier(0,0)(0.375,0.375)(0.75,0.75)}
\put(4.9,4.65){\color{Green}{\em 4}}
\put(2.625,4.375){\color{Purple}\qbezier(0,0)(0.375,0.375)(0.75,0.75)}
\put(2.9,4.65){\color{Purple}{\em 3}}
\put(3.625,5.375){\color{Green}\qbezier(0,0)(0.375,0.375)(0.75,0.75)}
\put(3.9,5.65){\color{Green}{\em 4}}
\put(2.625,6.375){\color{Green}\qbezier(0,0)(0.375,0.375)(0.75,0.75)}
\put(2.9,6.65){\color{Green}{\em 4}}
\put(1.625,7.375){\color{Green}\qbezier(0,0)(0.375,0.375)(0.75,0.75)}
\put(1.9,7.65){\color{Green}{\em 4}}
\put(2.625,8.375){\color{Purple}\qbezier(0,0)(0.375,0.375)(0.75,0.75)}
\put(2.9,8.65){\color{Purple}{\em 3}}
\put(3.625,9.375){\color{Red}\qbezier(0,0)(0.375,0.375)(0.75,0.75)}
\put(3.9,9.65){\color{Red}{\em 2}}
\put(4.625,10.375){\color{Cyan}\qbezier(0,0)(0.375,0.375)(0.75,0.75)}
\put(4.9,10.65){\color{Cyan}{\em 1}}
\put(0.625,8.375){\color{Green}\qbezier(0,0)(0.375,0.375)(0.75,0.75)}
\put(0.9,8.65){\color{Green}{\em 4}}
\put(1.625,9.375){\color{Purple}\qbezier(0,0)(0.375,0.375)(0.75,0.75)}
\put(1.9,9.65){\color{Purple}{\em 3}}
\put(2.625,10.375){\color{Red}\qbezier(0,0)(0.375,0.375)(0.75,0.75)}
\put(2.9,10.65){\color{Red}{\em 2}}
\put(3.625,11.375){\color{Cyan}\qbezier(0,0)(0.375,0.375)(0.75,0.75)}
\put(3.9,11.65){\color{Cyan}{\em 1}}
\put(2.625,12.375){\color{Cyan}\qbezier(0,0)(0.375,0.375)(0.75,0.75)}
\put(2.9,12.65){\color{Cyan}{\em 1}}
\put(3.39,3.375){\color{Cyan}\qbezier(0,0)(-0.375,0.375)(-0.75,0.75)}
\put(2.9,3.65){\color{Cyan}{\em 1}}
\put(4.39,4.375){\color{Cyan}\qbezier(0,0)(-0.375,0.375)(-0.75,0.75)}
\put(3.9,4.65){\color{Cyan}{\em 1}}
\put(3.39,5.375){\color{Red}\qbezier(0,0)(-0.375,0.375)(-0.75,0.75)}
\put(2.9,5.65){\color{Red}{\em 2}}
\put(2.39,6.375){\color{Purple}\qbezier(0,0)(-0.375,0.375)(-0.75,0.75)}
\put(1.9,6.65){\color{Purple}{\em 3}}
\put(1.39,7.375){\color{Green}\qbezier(0,0)(-0.375,0.375)(-0.75,0.75)}
\put(0.9,7.65){\color{Green}{\em 4}}
\put(5.39,5.375){\color{Cyan}\qbezier(0,0)(-0.375,0.375)(-0.75,0.75)}
\put(4.9,5.65){\color{Cyan}{\em 1}}
\put(4.39,6.375){\color{Red}\qbezier(0,0)(-0.375,0.375)(-0.75,0.75)}
\put(3.9,6.65){\color{Red}{\em 2}}
\put(3.39,7.375){\color{Purple}\qbezier(0,0)(-0.375,0.375)(-0.75,0.75)}
\put(2.9,7.65){\color{Purple}{\em 3}}
\put(2.39,8.375){\color{Green}\qbezier(0,0)(-0.375,0.375)(-0.75,0.75)}
\put(1.9,8.65){\color{Green}{\em 4}}
\put(3.39,9.375){\color{Green}\qbezier(0,0)(-0.375,0.375)(-0.75,0.75)}
\put(2.9,9.65){\color{Green}{\em 4}}
\put(4.39,10.375){\color{Green}\qbezier(0,0)(-0.375,0.375)(-0.75,0.75)}
\put(3.9,10.65){\color{Green}{\em 4}}
\put(3.39,11.375){\color{Purple}\qbezier(0,0)(-0.375,0.375)(-0.75,0.75)}
\put(2.9,11.65){\color{Purple}{\em 3}}
\put(5.39,11.375){\color{Green}\qbezier(0,0)(-0.375,0.375)(-0.75,0.75)}
\put(4.9,11.65){\color{Green}{\em 4}}
\put(4.39,12.375){\color{Purple}\qbezier(0,0)(-0.375,0.375)(-0.75,0.75)}
\put(3.9,12.65){\color{Purple}{\em 3}}
\put(3.39,13.375){\color{Red}\qbezier(0,0)(-0.375,0.375)(-0.75,0.75)}
\put(2.9,13.65){\color{Red}{\em 2}}
\put(2.39,14.375){\color{Purple}\qbezier(0,0)(-0.375,0.375)(-0.75,0.75)}
\put(1.9,14.65){\color{Purple}{\em 3}}
\put(1.39,15.375){\color{Green}\qbezier(0,0)(-0.375,0.375)(-0.75,0.75)}
\put(0.9,15.65){\color{Green}{\em 4}}
\put(4,8){\TypeEboxDot{Gray}}
\put(4.39,8.375){\color{Purple}\qbezier(0,0)(-0.375,0.375)(-0.75,0.75)}
\put(3.9,8.65){\color{Purple}{\em 3}}
\put(3.625,7.375){\color{Purple}\qbezier(0,0)(0.375,0.375)(0.75,0.75)}
\put(3.9,7.65){\color{Purple}{\em 3}}
\end{picture}
\hspace*{1cm}
\begin{picture}(5,17)
\put(2.6,13){\LARGE $\widetilde{P_{\mysmallF_{4}}}(\omega_{\mbox{\color{Green}{\small $4$}}})$}
\put(2.5,13.2){\vector(-2,-1){0.9}}
\put(3,12){\LARGE $\cong \mathbf{j}_{\mbox{\footnotesize color}}\big(\widetilde{L_{\mysmallF_{4}}}(\omega_{\mbox{\color{Green}{\small $4$}}})\big)$}
\put(0,13){\TypeEboxDot{Green}}
\put(0,7){\TypeEboxDot{Green}}
\put(1,12){\TypeEboxDot{Purple}}
\put(1,8){\TypeEboxDot{Purple}}
\put(1,6){\TypeEboxDot{Purple}}
\put(2,11){\TypeEboxDot{Red}}
\put(2,9){\TypeEboxDot{Red}}
\put(2,7){\TypeEboxDot{Red}}
\put(2,5){\TypeEboxDot{Red}}
\put(1,10){\TypeEboxDot{Cyan}}
\put(3,6){\TypeEboxDot{Cyan}}
\put(3,10){\TypeEboxDot{Purple}}
\put(3,8){\TypeEboxDot{Purple}}
\put(3,4){\TypeEboxDot{Purple}}
\put(4,9){\TypeEboxDot{Green}}
\put(4,3){\TypeEboxDot{Green}}
\put(0.7,13.3){\color{OliveGreen}{\em \footnotesize 4}}
\put(0.1,7){\color{OliveGreen}{\em \footnotesize 4}}
\put(1.7,12.3){\color{Purple}{\em \footnotesize 3}}
\put(1.1,8.3){\color{Purple}{\em \footnotesize 3}}
\put(1.1,6){\color{Purple}{\em \footnotesize 3}}
\put(2.7,11.3){\color{Red}{\em \footnotesize 2}}
\put(2.1,9.15){\color{Red}{\em \footnotesize 2}}
\put(2.7,7.15){\color{Red}{\em \footnotesize 2}}
\put(2.1,5){\color{Red}{\em \footnotesize 2}}
\put(1.1,10){\color{Cyan}{\em \footnotesize 1}}
\put(3.7,6.3){\color{Cyan}{\em \footnotesize 1}}
\put(3.7,10.3){\color{Purple}{\em \footnotesize 3}}
\put(3.7,8){\color{Purple}{\em \footnotesize 3}}
\put(3.1,4){\color{Purple}{\em \footnotesize 3}}
\put(4.7,9.3){\color{Green}{\em \footnotesize 4}}
\put(4.1,3){\color{Green}{\em \footnotesize 4}}
\thicklines
\put(2.625,5.375){\color{Black}\qbezier(0,0)(0.375,0.375)(0.75,0.75)}
\put(1.625,6.375){\color{Black}\qbezier(0,0)(0.375,0.375)(0.75,0.75)}
\put(2.625,7.375){\color{Black}\qbezier(0,0)(0.375,0.375)(0.75,0.75)}
\put(3.625,8.375){\color{Black}\qbezier(0,0)(0.375,0.375)(0.75,0.75)}
\put(0.625,7.375){\color{Black}\qbezier(0,0)(0.375,0.375)(0.75,0.75)}
\put(1.625,8.375){\color{Black}\qbezier(0,0)(0.375,0.375)(0.75,0.75)}
\put(2.625,9.375){\color{Black}\qbezier(0,0)(0.375,0.375)(0.75,0.75)}
\put(1.625,10.375){\color{Black}\qbezier(0,0)(0.375,0.375)(0.75,0.75)}
\put(4.39,3.375){\color{Black}\qbezier(0,0)(-0.375,0.375)(-0.75,0.75)}
\put(3.39,4.375){\color{Black}\qbezier(0,0)(-0.375,0.375)(-0.75,0.75)}
\put(2.39,5.375){\color{Black}\qbezier(0,0)(-0.375,0.375)(-0.75,0.75)}
\put(1.39,6.375){\color{Black}\qbezier(0,0)(-0.375,0.375)(-0.75,0.75)}
\put(2.39,7.375){\color{Black}\qbezier(0,0)(-0.375,0.375)(-0.75,0.75)}
\put(3.39,6.375){\color{Black}\qbezier(0,0)(-0.375,0.375)(-0.75,0.75)}
\put(3.39,8.375){\color{Black}\qbezier(0,0)(-0.375,0.375)(-0.75,0.75)}
\put(2.39,9.375){\color{Black}\qbezier(0,0)(-0.375,0.375)(-0.75,0.75)}
\put(4.39,9.375){\color{Black}\qbezier(0,0)(-0.375,0.375)(-0.75,0.75)}
\put(3.39,10.375){\color{Black}\qbezier(0,0)(-0.375,0.375)(-0.75,0.75)}
\put(2.39,11.375){\color{Black}\qbezier(0,0)(-0.375,0.375)(-0.75,0.75)}
\put(1.39,12.375){\color{Black}\qbezier(0,0)(-0.375,0.375)(-0.75,0.75)}
\end{picture}
\end{center}
\end{figure}

\begin{figure}
\begin{center}
{\small {\bf \FFourFigure.2}\ \  A 26-element full-length sublattice of the DCDL $\widetilde{L_{\mytinyF_{4}}}(\omega_{\mbox{\color{Green}{\tiny $4$}}})$ from \FFourFigure.1. The compression poset for this sublattice can be discerned, via \FullLengthTheorem, from \FFourFigure.1's $\widetilde{P_{\mytinyF_{4}}}(\omega_{\mbox{\color{Green}{\tiny $4$}}})$.}\\ 
{\footnotesize (This lattice, denoted $\mysmallfancyQ^{(\mbox{\color{Green}{\tiny $4$}})}$, is ``$\mysmallF_{4}$-structured'', cf.\ \S \GStructureIntroSection, and first of two ``quasi-minuscule DCML's'' for $\myF_{4}$, cf. \S \QuasiExampleSection.)}

\vspace*{0.1in} 
\FFourFlower

\vspace*{0.1in} 
\setlength{\unitlength}{1cm}
\begin{picture}(6,17)
\put(-1.5,12.5){\LARGE $\mybigfancyQ^{(\mbox{\color{Green}{\small $4$}})}$}
\put(-0.3,12.7){\vector(2,1){2.1}}
\put(-1.1,11.8){\LARGE $\cong \mathbf{J}_{\mbox{\footnotesize color}}\big(\mbox{\tt \huge Q}^{(\mbox{\color{Green}{\small $4$}})}\big)$}
\put(0,9.5){$\lsem \mbox{\color{Green}{$4$}},\mbox{\color{Green}{$4$}} \rsem$}
\put(0.5,9.3){\vector(0,-1){0.6}}
\put(3.8,8.175){$\lsem \mbox{\color{Purple}{$3$}},\mbox{\color{Green}{$4$}} \rsem$}
\put(3.75,8.275){\vector(-1,0){0.8}}
\put(0,16){\TypeEboxDot{Gray}}
\put(0,8){\TypeEboxDot{Gray}}
\put(0,0){\TypeEboxDot{Gray}}
\put(1,15){\TypeEboxDot{Gray}}
\put(1,9){\TypeEboxDot{Gray}}
\put(1,7){\TypeEboxDot{Gray}}
\put(1,1){\TypeEboxDot{Gray}}
\put(2,14){\TypeEboxDot{Gray}}
\put(2,10){\TypeEboxDot{Gray}}
\put(2,8){\TypeEboxDot{Gray}}
\put(2,6){\TypeEboxDot{Gray}}
\put(2,2){\TypeEboxDot{Gray}}
\put(3,13){\TypeEboxDot{Gray}}
\put(3,11){\TypeEboxDot{Gray}}
\put(3,9){\TypeEboxDot{Gray}}
\put(3,7){\TypeEboxDot{Gray}}
\put(3,5){\TypeEboxDot{Gray}}
\put(3,3){\TypeEboxDot{Gray}}
\put(2,12){\TypeEboxDot{Gray}}
\put(2,4){\TypeEboxDot{Gray}}
\put(4,12){\TypeEboxDot{Gray}}
\put(4,10){\TypeEboxDot{Gray}}
\put(4,6){\TypeEboxDot{Gray}}
\put(4,4){\TypeEboxDot{Gray}}
\put(5,11){\TypeEboxDot{Gray}}
\put(5,5){\TypeEboxDot{Gray}}
\thicklines
\put(0.625,0.375){\color{Green}\qbezier(0,0)(0.375,0.375)(0.75,0.75)}
\put(1.625,1.375){\color{Purple}\qbezier(0,0)(0.375,0.375)(0.75,0.75)}
\put(2.625,2.375){\color{Red}\qbezier(0,0)(0.375,0.375)(0.75,0.75)}
\put(3.625,3.375){\color{Purple}\qbezier(0,0)(0.375,0.375)(0.75,0.75)}
\put(4.625,4.375){\color{Green}\qbezier(0,0)(0.375,0.375)(0.75,0.75)}
\put(2.625,4.375){\color{Purple}\qbezier(0,0)(0.375,0.375)(0.75,0.75)}
\put(3.625,5.375){\color{Green}\qbezier(0,0)(0.375,0.375)(0.75,0.75)}
\put(2.625,6.375){\color{Green}\qbezier(0,0)(0.375,0.375)(0.75,0.75)}
%
\put(1.625,7.375){\color{Green}\qbezier(0,0)(0.375,0.375)(0.75,0.75)}
\put(2.625,8.375){\color{Purple}\qbezier(0,0)(0.375,0.375)(0.75,0.75)}
\put(3.625,9.375){\color{Red}\qbezier(0,0)(0.375,0.375)(0.75,0.75)}
\put(4.625,10.375){\color{Cyan}\qbezier(0,0)(0.375,0.375)(0.75,0.75)}
\put(0.625,8.375){\color{Green}\qbezier(0,0)(0.375,0.375)(0.75,0.75)}
\put(1.625,9.375){\color{Purple}\qbezier(0,0)(0.375,0.375)(0.75,0.75)}
\put(2.625,10.375){\color{Red}\qbezier(0,0)(0.375,0.375)(0.75,0.75)}
\put(3.625,11.375){\color{Cyan}\qbezier(0,0)(0.375,0.375)(0.75,0.75)}
\put(2.625,12.375){\color{Cyan}\qbezier(0,0)(0.375,0.375)(0.75,0.75)}
\put(3.39,3.375){\color{Cyan}\qbezier(0,0)(-0.375,0.375)(-0.75,0.75)}
\put(4.39,4.375){\color{Cyan}\qbezier(0,0)(-0.375,0.375)(-0.75,0.75)}
\put(3.39,5.375){\color{Red}\qbezier(0,0)(-0.375,0.375)(-0.75,0.75)}
\put(2.39,6.375){\color{Purple}\qbezier(0,0)(-0.375,0.375)(-0.75,0.75)}
\put(1.39,7.375){\color{Green}\qbezier(0,0)(-0.375,0.375)(-0.75,0.75)}
\put(5.39,5.375){\color{Cyan}\qbezier(0,0)(-0.375,0.375)(-0.75,0.75)}
\put(4.39,6.375){\color{Red}\qbezier(0,0)(-0.375,0.375)(-0.75,0.75)}
\put(3.39,7.375){\color{Purple}\qbezier(0,0)(-0.375,0.375)(-0.75,0.75)}
\put(2.39,8.375){\color{Green}\qbezier(0,0)(-0.375,0.375)(-0.75,0.75)}
%
\put(3.39,9.375){\color{Green}\qbezier(0,0)(-0.375,0.375)(-0.75,0.75)}
\put(4.39,10.375){\color{Green}\qbezier(0,0)(-0.375,0.375)(-0.75,0.75)}
\put(3.39,11.375){\color{Purple}\qbezier(0,0)(-0.375,0.375)(-0.75,0.75)}
\put(5.39,11.375){\color{Green}\qbezier(0,0)(-0.375,0.375)(-0.75,0.75)}
\put(4.39,12.375){\color{Purple}\qbezier(0,0)(-0.375,0.375)(-0.75,0.75)}
\put(3.39,13.375){\color{Red}\qbezier(0,0)(-0.375,0.375)(-0.75,0.75)}
\put(2.39,14.375){\color{Purple}\qbezier(0,0)(-0.375,0.375)(-0.75,0.75)}
\put(1.39,15.375){\color{Green}\qbezier(0,0)(-0.375,0.375)(-0.75,0.75)}
\put(2.9,3.65){\color{Cyan}{\em 1}}
\put(3.9,4.65){\color{Cyan}{\em 1}}
\put(4.9,5.65){\color{Cyan}{\em 1}}
\put(4.9,10.65){\color{Cyan}{\em 1}}
\put(3.9,11.65){\color{Cyan}{\em 1}}
\put(2.9,12.65){\color{Cyan}{\em 1}}
\put(2.9,5.65){\color{Red}{\em 2}}
\put(3.9,6.65){\color{Red}{\em 2}}
\put(2.9,13.65){\color{Red}{\em 2}}
\put(3.9,9.65){\color{Red}{\em 2}}
\put(2.9,10.65){\color{Red}{\em 2}}
\put(2.9,2.65){\color{Red}{\em 2}}
\put(1.9,6.65){\color{Purple}{\em 3}}
\put(2.9,7.65){\color{Purple}{\em 3}}
%
\put(2.9,11.65){\color{Purple}{\em 3}}
\put(3.9,12.65){\color{Purple}{\em 3}}
\put(1.9,14.65){\color{Purple}{\em 3}}
\put(1.9,9.65){\color{Purple}{\em 3}}
\put(2.9,8.65){\color{Purple}{\em 3}}
%
\put(3.9,3.65){\color{Purple}{\em 3}}
\put(2.9,4.65){\color{Purple}{\em 3}}
\put(1.9,1.65){\color{Purple}{\em 3}}
\put(0.9,7.65){\color{Green}{\em 4}}
\put(1.9,8.65){\color{Green}{\em 4}}
\put(2.9,9.65){\color{Green}{\em 4}}
\put(3.9,10.65){\color{Green}{\em 4}}
\put(4.9,11.65){\color{Green}{\em 4}}
\put(0.9,15.65){\color{Green}{\em 4}}
\put(0.9,8.65){\color{Green}{\em 4}}
\put(1.9,7.65){\color{Green}{\em 4}}
\put(2.9,6.65){\color{Green}{\em 4}}
\put(3.9,5.65){\color{Green}{\em 4}}
\put(4.9,4.65){\color{Green}{\em 4}}
\put(0.9,0.65){\color{Green}{\em 4}}
\end{picture}
\hspace*{1cm}
\begin{picture}(5,17)
\put(2.6,13){\LARGE $\mbox{\tt \huge Q}^{(\mbox{\color{Green}{\small $4$}})}$}
\put(2.5,13.2){\vector(-2,-1){0.9}}
\put(3,12.3){\LARGE $\cong \mathbf{j}_{\mbox{\footnotesize color}}\big(\mybigfancyQ^{(\mbox{\color{Green}{\small $4$}})}\big)$}
\put(0,13){\TypeEboxDot{Green}}
\put(0,7){\TypeEboxDot{Green}}
\put(1,12){\TypeEboxDot{Purple}}
\put(1,8){\TypeEboxDot{Purple}}
\put(1,6){\TypeEboxDot{Purple}}
\put(2,11){\TypeEboxDot{Red}}
\put(2,9){\TypeEboxDot{Red}}
\put(2,7){\TypeEboxDot{Red}}
\put(2,5){\TypeEboxDot{Red}}
\put(1,10){\TypeEboxDot{Cyan}}
\put(3,6){\TypeEboxDot{Cyan}}
\put(3,10){\TypeEboxDot{Purple}}
\put(3,8){\TypeEboxDot{Purple}}
\put(3,4){\TypeEboxDot{Purple}}
\put(4,9){\TypeEboxDot{Green}}
\put(4,3){\TypeEboxDot{Green}}
\put(0.7,13.3){\color{OliveGreen}{\em \footnotesize 4}}
\put(0.1,7){\color{OliveGreen}{\em \footnotesize 4}}
\put(1.7,12.3){\color{Purple}{\em \footnotesize 3}}
\put(1.1,8.3){\color{Purple}{\em \footnotesize 3}}
\put(1.1,6){\color{Purple}{\em \footnotesize 3}}
\put(2.7,11.3){\color{Red}{\em \footnotesize 2}}
\put(2.1,9.15){\color{Red}{\em \footnotesize 2}}
\put(2.7,7.15){\color{Red}{\em \footnotesize 2}}
\put(2.1,5){\color{Red}{\em \footnotesize 2}}
\put(1.1,10){\color{Cyan}{\em \footnotesize 1}}
\put(3.7,6.3){\color{Cyan}{\em \footnotesize 1}}
\put(3.7,10.3){\color{Purple}{\em \footnotesize 3}}
\put(3.7,8){\color{Purple}{\em \footnotesize 3}}
\put(3.1,4){\color{Purple}{\em \footnotesize 3}}
\put(4.7,9.3){\color{Green}{\em \footnotesize 4}}
\put(4.1,3){\color{Green}{\em \footnotesize 4}}
\thicklines
\put(2.625,5.375){\color{Black}\qbezier(0,0)(0.375,0.375)(0.75,0.75)}
\put(1.625,6.375){\color{Black}\qbezier(0,0)(0.375,0.375)(0.75,0.75)}
\put(2.625,7.375){\color{Black}\qbezier(0,0)(0.375,0.375)(0.75,0.75)}
\put(3.625,8.375){\color{Black}\qbezier(0,0)(0.375,0.375)(0.75,0.75)}
\put(0.625,7.375){\color{Black}\qbezier(0,0)(0.375,0.375)(0.75,0.75)}
\put(1.625,8.375){\color{Black}\qbezier(0,0)(0.375,0.375)(0.75,0.75)}
\put(2.625,9.375){\color{Black}\qbezier(0,0)(0.375,0.375)(0.75,0.75)}
\put(1.625,10.375){\color{Black}\qbezier(0,0)(0.375,0.375)(0.75,0.75)}
\put(4.39,3.375){\color{Black}\qbezier(0,0)(-0.375,0.375)(-0.75,0.75)}
\put(3.39,4.375){\color{Black}\qbezier(0,0)(-0.375,0.375)(-0.75,0.75)}
\put(2.39,5.375){\color{Black}\qbezier(0,0)(-0.375,0.375)(-0.75,0.75)}
\put(1.39,6.375){\color{Black}\qbezier(0,0)(-0.375,0.375)(-0.75,0.75)}
%
\put(3.39,6.375){\color{Black}\qbezier(0,0)(-0.375,0.375)(-0.75,0.75)}
%
\put(2.39,9.375){\color{Black}\qbezier(0,0)(-0.375,0.375)(-0.75,0.75)}
\put(4.39,9.375){\color{Black}\qbezier(0,0)(-0.375,0.375)(-0.75,0.75)}
\put(3.39,10.375){\color{Black}\qbezier(0,0)(-0.375,0.375)(-0.75,0.75)}
\put(2.39,11.375){\color{Black}\qbezier(0,0)(-0.375,0.375)(-0.75,0.75)}
\put(1.39,12.375){\color{Black}\qbezier(0,0)(-0.375,0.375)(-0.75,0.75)}
\put(3.375,8.125){\color{Black}\qbezier(0,0)(-0.875,-0.3)(-1.75,0)}
\put(2.5,7.975){\vector(-1,0){0.2}}
\end{picture}
\end{center}
\end{figure}

\begin{figure}
\begin{center}
{\bf \FFourFigure.3}\ \  Another 26-element full-length sublattice of the DCDL $\widetilde{L_{\mytinyF_{4}}}(\omega_{\mbox{\color{Green}{\tiny $4$}}})$ from \FFourFigure.1. Its compression poset can be discerned, via \FullLengthTheorem, from \FFourFigure.1's $\widetilde{P_{\mytinyF_{4}}}(\omega_{\mbox{\color{Green}{\tiny $4$}}})$.\\
{\footnotesize (This lattice, denoted $\mysmallfancyQ^{(\mbox{\color{Purple}{\tiny $3$}})}$, is ``$\mysmallF_{4}$-structured'', cf.\ \S \GStructureIntroSection, and second of two ``quasi-minuscule DCML's'' for $\myF_{4}$, cf. \S \QuasiExampleSection.)}

\vspace*{0.1in} 
\FFourFlower

\vspace*{0.1in} 
\setlength{\unitlength}{1cm}
\begin{picture}(6,17)
\put(-1.5,12.5){\LARGE $\mybigfancyQ^{(\mbox{\color{Purple}{\small $3$}})}$}
\put(-0.3,12.7){\vector(2,1){2.1}}
\put(-1.1,11.8){\LARGE $\cong \mathbf{J}_{\mbox{\footnotesize color}}\big(\mbox{\tt \huge Q}^{(\mbox{\color{Purple}{\small $3$}})}\big)$}
\put(5.5,8.175){$\lsem \mbox{\color{Purple}{$3$}},\mbox{\color{Purple}{$3$}} \rsem$}
\put(5.4,8.275){\vector(-1,0){0.6}}
\put(0.3,8.175){$\lsem \mbox{\color{Purple}{$3$}},\mbox{\color{Green}{$4$}} \rsem$}
\put(1.3,8.275){\vector(1,0){0.8}}
\put(0,16){\TypeEboxDot{Gray}}
\put(0,0){\TypeEboxDot{Gray}}
\put(1,15){\TypeEboxDot{Gray}}
\put(1,9){\TypeEboxDot{Gray}}
\put(1,7){\TypeEboxDot{Gray}}
\put(1,1){\TypeEboxDot{Gray}}
\put(2,14){\TypeEboxDot{Gray}}
\put(2,10){\TypeEboxDot{Gray}}
\put(2,8){\TypeEboxDot{Gray}}
\put(2,6){\TypeEboxDot{Gray}}
\put(2,2){\TypeEboxDot{Gray}}
\put(3,13){\TypeEboxDot{Gray}}
\put(3,11){\TypeEboxDot{Gray}}
\put(3,9){\TypeEboxDot{Gray}}
\put(3,7){\TypeEboxDot{Gray}}
\put(3,5){\TypeEboxDot{Gray}}
\put(3,3){\TypeEboxDot{Gray}}
\put(2,12){\TypeEboxDot{Gray}}
\put(4,8){\TypeEboxDot{Gray}}
\put(2,4){\TypeEboxDot{Gray}}
\put(4,12){\TypeEboxDot{Gray}}
\put(4,10){\TypeEboxDot{Gray}}
\put(4,6){\TypeEboxDot{Gray}}
\put(4,4){\TypeEboxDot{Gray}}
\put(5,11){\TypeEboxDot{Gray}}
\put(5,5){\TypeEboxDot{Gray}}
\thicklines
\put(0.625,0.375){\color{Green}\qbezier(0,0)(0.375,0.375)(0.75,0.75)}
\put(1.625,1.375){\color{Purple}\qbezier(0,0)(0.375,0.375)(0.75,0.75)}
\put(2.625,2.375){\color{Red}\qbezier(0,0)(0.375,0.375)(0.75,0.75)}
\put(3.625,3.375){\color{Purple}\qbezier(0,0)(0.375,0.375)(0.75,0.75)}
\put(4.625,4.375){\color{Green}\qbezier(0,0)(0.375,0.375)(0.75,0.75)}
\put(2.625,4.375){\color{Purple}\qbezier(0,0)(0.375,0.375)(0.75,0.75)}
\put(3.625,5.375){\color{Green}\qbezier(0,0)(0.375,0.375)(0.75,0.75)}
\put(2.625,6.375){\color{Green}\qbezier(0,0)(0.375,0.375)(0.75,0.75)}
\put(3.625,7.375){\color{Purple}\qbezier(0,0)(0.375,0.375)(0.75,0.75)}
\put(1.625,7.375){\color{Green}\qbezier(0,0)(0.375,0.375)(0.75,0.75)}
\put(2.625,8.375){\color{Purple}\qbezier(0,0)(0.375,0.375)(0.75,0.75)}
\put(3.625,9.375){\color{Red}\qbezier(0,0)(0.375,0.375)(0.75,0.75)}
\put(4.625,10.375){\color{Cyan}\qbezier(0,0)(0.375,0.375)(0.75,0.75)}
%
\put(1.625,9.375){\color{Purple}\qbezier(0,0)(0.375,0.375)(0.75,0.75)}
\put(2.625,10.375){\color{Red}\qbezier(0,0)(0.375,0.375)(0.75,0.75)}
\put(3.625,11.375){\color{Cyan}\qbezier(0,0)(0.375,0.375)(0.75,0.75)}
\put(2.625,12.375){\color{Cyan}\qbezier(0,0)(0.375,0.375)(0.75,0.75)}
\put(3.39,3.375){\color{Cyan}\qbezier(0,0)(-0.375,0.375)(-0.75,0.75)}
\put(4.39,4.375){\color{Cyan}\qbezier(0,0)(-0.375,0.375)(-0.75,0.75)}
\put(3.39,5.375){\color{Red}\qbezier(0,0)(-0.375,0.375)(-0.75,0.75)}
\put(2.39,6.375){\color{Purple}\qbezier(0,0)(-0.375,0.375)(-0.75,0.75)}
%
\put(5.39,5.375){\color{Cyan}\qbezier(0,0)(-0.375,0.375)(-0.75,0.75)}
\put(4.39,6.375){\color{Red}\qbezier(0,0)(-0.375,0.375)(-0.75,0.75)}
\put(3.39,7.375){\color{Purple}\qbezier(0,0)(-0.375,0.375)(-0.75,0.75)}
\put(2.39,8.375){\color{Green}\qbezier(0,0)(-0.375,0.375)(-0.75,0.75)}
\put(4.39,8.375){\color{Purple}\qbezier(0,0)(-0.375,0.375)(-0.75,0.75)}
\put(3.39,9.375){\color{Green}\qbezier(0,0)(-0.375,0.375)(-0.75,0.75)}
\put(4.39,10.375){\color{Green}\qbezier(0,0)(-0.375,0.375)(-0.75,0.75)}
\put(3.39,11.375){\color{Purple}\qbezier(0,0)(-0.375,0.375)(-0.75,0.75)}
\put(5.39,11.375){\color{Green}\qbezier(0,0)(-0.375,0.375)(-0.75,0.75)}
\put(4.39,12.375){\color{Purple}\qbezier(0,0)(-0.375,0.375)(-0.75,0.75)}
\put(3.39,13.375){\color{Red}\qbezier(0,0)(-0.375,0.375)(-0.75,0.75)}
\put(2.39,14.375){\color{Purple}\qbezier(0,0)(-0.375,0.375)(-0.75,0.75)}
\put(1.39,15.375){\color{Green}\qbezier(0,0)(-0.375,0.375)(-0.75,0.75)}
\put(2.9,3.65){\color{Cyan}{\em 1}}
\put(3.9,4.65){\color{Cyan}{\em 1}}
\put(4.9,5.65){\color{Cyan}{\em 1}}
\put(4.9,10.65){\color{Cyan}{\em 1}}
\put(3.9,11.65){\color{Cyan}{\em 1}}
\put(2.9,12.65){\color{Cyan}{\em 1}}
\put(2.9,5.65){\color{Red}{\em 2}}
\put(3.9,6.65){\color{Red}{\em 2}}
\put(2.9,13.65){\color{Red}{\em 2}}
\put(3.9,9.65){\color{Red}{\em 2}}
\put(2.9,10.65){\color{Red}{\em 2}}
\put(2.9,2.65){\color{Red}{\em 2}}
\put(1.9,6.65){\color{Purple}{\em 3}}
\put(2.9,7.65){\color{Purple}{\em 3}}
\put(3.9,8.65){\color{Purple}{\em 3}}
\put(2.9,11.65){\color{Purple}{\em 3}}
\put(3.9,12.65){\color{Purple}{\em 3}}
\put(1.9,14.65){\color{Purple}{\em 3}}
\put(1.9,9.65){\color{Purple}{\em 3}}
\put(2.9,8.65){\color{Purple}{\em 3}}
\put(3.9,7.65){\color{Purple}{\em 3}}
\put(3.9,3.65){\color{Purple}{\em 3}}
\put(2.9,4.65){\color{Purple}{\em 3}}
\put(1.9,1.65){\color{Purple}{\em 3}}
\put(1.9,8.65){\color{Green}{\em 4}}
\put(2.9,9.65){\color{Green}{\em 4}}
\put(3.9,10.65){\color{Green}{\em 4}}
\put(4.9,11.65){\color{Green}{\em 4}}
\put(0.9,15.65){\color{Green}{\em 4}}
\put(1.9,7.65){\color{Green}{\em 4}}
\put(2.9,6.65){\color{Green}{\em 4}}
\put(3.9,5.65){\color{Green}{\em 4}}
\put(4.9,4.65){\color{Green}{\em 4}}
\put(0.9,0.65){\color{Green}{\em 4}}
\end{picture}
\hspace*{1cm}
\begin{picture}(5,17)
\put(2.6,13){\LARGE $\mbox{\tt \huge Q}^{(\mbox{\color{Purple}{\small $3$}})}$}
\put(2.5,13.2){\vector(-2,-1){0.9}}
\put(3,12.3){\LARGE $\cong \mathbf{j}_{\mbox{\footnotesize color}}\big(\mybigfancyQ^{(\mbox{\color{Purple}{\small $3$}})}\big)$}
\put(0,13){\TypeEboxDot{Green}}
\put(0,7){\TypeEboxDot{Green}}
\put(1,12){\TypeEboxDot{Purple}}
\put(1,8){\TypeEboxDot{Purple}}
\put(1,6){\TypeEboxDot{Purple}}
\put(2,11){\TypeEboxDot{Red}}
\put(2,9){\TypeEboxDot{Red}}
\put(2,7){\TypeEboxDot{Red}}
\put(2,5){\TypeEboxDot{Red}}
\put(1,10){\TypeEboxDot{Cyan}}
\put(3,6){\TypeEboxDot{Cyan}}
\put(3,10){\TypeEboxDot{Purple}}
\put(3,8){\TypeEboxDot{Purple}}
\put(3,4){\TypeEboxDot{Purple}}
\put(4,9){\TypeEboxDot{Green}}
\put(4,3){\TypeEboxDot{Green}}
\put(0.7,13.3){\color{OliveGreen}{\em \footnotesize 4}}
\put(0.1,7){\color{OliveGreen}{\em \footnotesize 4}}
\put(1.7,12.3){\color{Purple}{\em \footnotesize 3}}
\put(1.1,8.3){\color{Purple}{\em \footnotesize 3}}
\put(1.1,6){\color{Purple}{\em \footnotesize 3}}
\put(2.7,11.3){\color{Red}{\em \footnotesize 2}}
\put(2.1,9.15){\color{Red}{\em \footnotesize 2}}
\put(2.7,7.15){\color{Red}{\em \footnotesize 2}}
\put(2.1,5){\color{Red}{\em \footnotesize 2}}
\put(1.1,10){\color{Cyan}{\em \footnotesize 1}}
\put(3.7,6.3){\color{Cyan}{\em \footnotesize 1}}
\put(3.7,10.3){\color{Purple}{\em \footnotesize 3}}
\put(3.7,8){\color{Purple}{\em \footnotesize 3}}
\put(3.1,4){\color{Purple}{\em \footnotesize 3}}
\put(4.7,9.3){\color{Green}{\em \footnotesize 4}}
\put(4.1,3){\color{Green}{\em \footnotesize 4}}
\thicklines
\put(2.625,5.375){\color{Black}\qbezier(0,0)(0.375,0.375)(0.75,0.75)}
\put(1.625,6.375){\color{Black}\qbezier(0,0)(0.375,0.375)(0.75,0.75)}
\put(2.625,7.375){\color{Black}\qbezier(0,0)(0.375,0.375)(0.75,0.75)}
\put(3.625,8.375){\color{Black}\qbezier(0,0)(0.375,0.375)(0.75,0.75)}
\put(0.625,7.375){\color{Black}\qbezier(0,0)(0.375,0.375)(0.75,0.75)}
\put(1.625,8.375){\color{Black}\qbezier(0,0)(0.375,0.375)(0.75,0.75)}
\put(2.625,9.375){\color{Black}\qbezier(0,0)(0.375,0.375)(0.75,0.75)}
\put(1.625,10.375){\color{Black}\qbezier(0,0)(0.375,0.375)(0.75,0.75)}
\put(4.39,3.375){\color{Black}\qbezier(0,0)(-0.375,0.375)(-0.75,0.75)}
\put(3.39,4.375){\color{Black}\qbezier(0,0)(-0.375,0.375)(-0.75,0.75)}
\put(2.39,5.375){\color{Black}\qbezier(0,0)(-0.375,0.375)(-0.75,0.75)}
\put(1.39,6.375){\color{Black}\qbezier(0,0)(-0.375,0.375)(-0.75,0.75)}
\put(2.39,7.375){\color{Black}\qbezier(0,0)(-0.375,0.375)(-0.75,0.75)}
\put(3.39,6.375){\color{Black}\qbezier(0,0)(-0.375,0.375)(-0.75,0.75)}
\put(3.39,8.375){\color{Black}\qbezier(0,0)(-0.375,0.375)(-0.75,0.75)}
\put(2.39,9.375){\color{Black}\qbezier(0,0)(-0.375,0.375)(-0.75,0.75)}
\put(4.39,9.375){\color{Black}\qbezier(0,0)(-0.375,0.375)(-0.75,0.75)}
\put(3.39,10.375){\color{Black}\qbezier(0,0)(-0.375,0.375)(-0.75,0.75)}
\put(2.39,11.375){\color{Black}\qbezier(0,0)(-0.375,0.375)(-0.75,0.75)}
\put(1.39,12.375){\color{Black}\qbezier(0,0)(-0.375,0.375)(-0.75,0.75)}
\put(0.675,7.25){\color{Black}\qbezier(0,0)(1,0)(3.675,1.925)}
\put(2.7,8.075){\vector(3,2){0.2}}
\end{picture}
\end{center}
\end{figure}

\clearpage
\begin{center}
\fbox{\Large \bf Part II}\\ 
\underline{\large \bf Diamond-coloring and some new objects, methods, and results}
\end{center}
Within the algebraic-combinatorial context of our research, our work with diamond-colored modular and distributive lattices occasionally requires us to appeal to diamond-coloring results that are more general in nature. 
Some such results, which are mostly intuitive and natural, were recounted in Part I above.  
In the next several sections, and especially \S \GStructureIntroSection\ and \S \FrameSection, we present some general ideas -- for constructing DCDL's, for understanding certain structural limitations on DCML's and DCDL's, for re-conceptualizing combinatorial starting points -- that have been developed by the author and are largely new. 
We will also see how these ideas interact with some special bichromatic grid-like vertex-colored posets, called `two-color grid posets', which were first presented in \cite{ADLP} and \cite{ADLMPPW} and are re-considered here in \S \GridPosetSection\ and \S \TwoColorSturdySection.

\vspace*{0.5cm} 
\noindent 
{\bf \S \GStructureIntroSection.\ Integral embryophytes and sturdiness.} 
One striking feature of many naturally occurring DCDL's and DCML's is that the arrangement of edges is constrained by a certain arithmetic$\!$ /$\!$ combinatorial requirement.  
This structuring property has been considered in papers such as \cite{DonSupp}, \cite{DLP1}, \cite{DE}. 
It is, non-obviously, very closely connected to the classification of finite-dimensional simple Lie algebras over $\mathbb{C}$. 
Our purpose here is to further develop our understanding of this property within the context of DCDL's and DCML's and to make some additions to the growing collection of what are known as Dynkin diagram classification theorems. 

{\bf [\S \GStructureIntroSection.1:\! Embryophytes.]} To this end, let $I$ be a (finite, nonempty) palette of colors, and let $M = (M_{ij})_{i,j \in I}$ be an $I \times I$ matrix of real numbers. 
Let $\Gamma = \Gamma(M)$ be the simple graph\footnote{In particular, $\Gamma$ has no directed edges, no multiple edges, and no loops.} whose elements are $\{\gamma_{i}\}_{i \in I}$ and such that, for distinct $i, j \in I$, $\gamma_{i}$ and $\gamma_{j}$ are adjacent if and only if $M_{ij}M_{ji} \not= 0$. 
When the entries of $M$ satisfy 
\[M_{ij} \left\{\begin{array}{ccl} = & 2 & \hspace*{0.1in} \mbox{if $i$ = $j$}\\
= & 0 & \hspace*{0.1in} \mbox{if $\gamma_{i}$ is not adjacent to $\gamma_{j}$ in $\Gamma$}\\
< & 0 & \hspace*{0.1in} \mbox{if $\gamma_{i}$ is adjacent to $\gamma_{j}$ in $\Gamma$,}
\end{array}\right.\]
then we call $M$ a {\em pulsation matrix}, $\Gamma$ its {\em vascular graph}, and the pair $\mathscr{G} = (\Gamma,M)$ an {\em embryophytic graph} or just an {\em embryophyte}\footnote{This non-standard language is inspired by two considerations. The first is to make these primary objects more generic in the case that the matrix entries are integers; the standard name `generalized Cartan matrices' places them in, perhaps, a too specific context. 
The second is to support the botanical theme of our proposed name for the multifaceted classification result \LFK, which we call {\em La Florado Klasado}.}. 
To emphasize the role of $I$, we sometimes write $\mathscr{G} = (\Gamma_{I},M_{I \times I})$ when $I = I(\mathscr{G})$ is the color palette of the embryophyte $\mathscr{G}$. 
In such a setting, call the row vector $\myroot_{i} := (M_{ij})_{j \in I}$ the $i^{\mbox{\tiny th}}$ {\em simple root} for $\mathscr{G}$, and call the quantity $|M_{ij}|$ a {\em pulsation factor} when $M_{ij} < 0$. 
The number of colors in $I$ --- which is also the number of nodes of $\Gamma$ and the number of rows of $M$ --- is called the {\em rank} of $\mathscr{G}$. 

As we seek connections between combinatorics and algebra in the remainder of this manuscript, embryophytes will typically serve as our initial objects. 
For motivation, notice that if a pulsation matrix $M = M(\mathscr{G})$ of some embryophyte $\mathscr{G}$ has integer-valued pulsation factors, then $M$ meets the requirements for a {\em generalized Cartan matrix}, or `GCM';  
if, in addition, there is a diagonal matrix $D$ with positive main diagonal entries and a symmetric and positive definite matrix $S$ such that $DM=S$, then $M$ is a {\em Cartan matrix} (see, for example, \cite{Kac}). 
In our context, it is enough for now simply to think of $M$ as supplying additional information to the edges of the simple graph $\Gamma = \Gamma(\mathscr{G})$. 

{\bf [\S \GStructureIntroSection.2:\! Working with embryophytes.]} 
From here on, unless stated otherwise, assume that $\mathscr{G} = (\Gamma_{I},M_{I \times I})$ is a generic embryophyte. 
When $J$ is a subset of the color palette $I(\mathscr{G})$, we let $\Gamma_{J}$ be the subgraph of $\Gamma(\mathscr{G})$ induced by the nodes $\{\gamma_{j}\}_{j \in J}$ and $M_{J \times J}$ the $J \times J$ submatrix $\left(M_{jk}\right)_{j,k \in J}$ of the pulsation matrix $M(\mathscr{G})$. 
It is clear that the pair $\mathscr{H} := (\Gamma_{J},M_{J \times J})$ is an embryophyte with color palette $J = I(\mathscr{H})$, vascular graph $\Gamma_{J} = \Gamma(\mathscr{H})$, and pulsation matrix $M_{J \times J} = M(\mathscr{H})$. 
Any embryophyte $\mathscr{H}$ formed within $\mathscr{G}$ in this way is a {\em sub-embryophyte} of $\mathscr{G}$, and we say $\mathscr{H} =: \mathscr{G}_{J}$ is the sub-embryophyte {\em induced by the color sub-palette} $J$. 
When $\mathscr{G} = (\Gamma_{I},N_{I \times I})$ and $\mathscr{H} = (\Gamma'_{J},N'_{J \times J})$ are embryophytes such that $I \cap J = \emptyset$, then the {\em disjoint sum} $\mathscr{G} \oplus \mathscr{H}$ is the embryophyte whose color palette is $I \cup J$, vascular graph is the disjoint sum $\Gamma_{I} \oplus \Gamma'_{J}$, and pulsation matrix, denoted `$N_{I \times I} \oplus N'_{J \times J}$', is the matrix $M$ whose $(k,l)$-entry, for all $k$ and $l$ in $I \cup J$, is
\[M_{kl} = \left\{\begin{array}{cl} N_{kl} & \mbox{if } k,l \in I\\ N'_{kl} & \mbox{if } k,l \in J\\ 0 & \mbox{otherwise}\end{array}\right..\]
In this case, we can write $\mathscr{G} \oplus \mathscr{H} = \left(\Gamma(\mathscr{G}) \oplus \Gamma(\mathscr{H}), M(\mathscr{G}) \oplus M(\mathscr{H})\right)$. 
Say two embryophytes $\mathscr{G}$ and $\mathscr{H}$ are {\em isomorphic} and write $\mathscr{G} \cong \mathscr{H}$ if there is a bijection between their color palettes that preserves vascular graph vertices and edges and pulsation matrix entries. 
We say $\mathscr{G}$ is {\em connected} and its pulsation matrix is {\em indecomposable} if there is no partitioning of the color palette $I(\mathscr{G})$ into a disjoint union $J \cup J'$ such that $\mathscr{G} \cong \mathscr{G}_{J} \oplus \mathscr{G}_{J'}$. 
A {\em standard decomposition} of an embryophyte $\mathscr{G}$ is a partitioning of the color palette $I(\mathscr{G})$ into a disjoint union $J_{1} \cup \cdots \cup J_{p}$ such that each induced sub-embryophyte $\mathscr{G}_{J_{q}}$ is connected and $\mathscr{G} \cong \mathscr{G}_{J_{1}} \oplus \cdots \oplus \mathscr{G}_{J_{p}}$. 

An {\em integral embryophytic graph (IEG)}, or just an {\em integral embryophyte}, has all $M_{ij} \in \mathbb{Z}$. 
When all pulsation factors are unity, then we say $\mathscr{G} = (\Gamma,M)$ is a {\em unital embryophytic graph (UEG)}, or just a {\em unital embryophyte}. 
When depicting an IEG, we attach to each edge $\gamma_{i} \bullet\!\!\rule[0.75mm]{10mm}{0.25mm}\!\!\bullet \gamma_{j}$ in the vascular graph the pulsation factors $p=|M_{ij}|$ and $q=|M_{ji}|$ as follows: 

\vspace*{-0.075in}
\noindent
\begin{center}
\TwoCitiesGraphWithLabels
\end{center}

\vspace*{-0.075in}
\noindent
We sometimes use multiple arrows to represent the pulsation factors when these are small integers, say $|M_{ij}|=2$ and $|M_{ji}|=3$: \DNGTwoNodeArrows.
We use an unlabelled edge \ATwoGraphNoEdgeLabels\ as shorthand for an edge whose pulsation factors are unity, and we attach `$\infty$' if $M_{ij}M_{ji} \geq 4$, as in\  \DNGGraphCircleInfinity.  

{\bf [\S \GStructureIntroSection.3:\! Coxeter and Dynkin.]} Objects entirely similar to integral embryophytes were innovated by H.\ S.\ M.\ Coxeter in his seminal classification of finite reflection groups \cite{Coxeter} and independently by E.\ B.\ Dynkin in simplifying the classification of the finite-dimensional simple Lie algebras over the field of complex numbers (see \cite{DynkinClassify}, \cite{DynkinTranslation}, and \cite{KOV}). 
Indeed, embryophytes are typically called `Dynkin diagrams' when the pulsation matrix entries are integers.   
Amongst connected IEG's, those appearing in \IEGGraphFigure, which we call {\em integral Coxeter--Dynkin flowers}, impose many beautiful finiteness properties on related algebraic and combinatorial objects, some of which we will see below. 
We call an IEG an integral {\em Coxeter--Dynkin posy} if it is a disjoint sum of integral Coxeter--Dynkin flowers. 

\newcommand{\AnDNGGraph}{\setlength{\unitlength}{0.75in}
\begin{picture}(6,0.85)
\put(0,0){\begin{picture}(1,0)
            \put(0,0.1){\circle*{0.075}}
            \put(1,0.1){\circle*{0.075}}
            \put(2,0.1){\circle*{0.075}}
            \put(4,0.1){\circle*{0.075}}
            \put(5,0.1){\circle*{0.075}}
            \put(6,0.1){\circle*{0.075}}
            \put(0,0.1){\line(1,0){2}}
            \multiput(2,0.1)(0.4,0){5}{\line(1,0){0.2}}
            \put(4,0.1){\line(1,0){2}}
            \put(-0.05,-0.075){\tiny \em 1}
            \put(0.95,-0.075){\tiny \em 2}
            \put(1.95,-0.075){\tiny \em 3}
            \put(3.9,-0.075){\tiny \em n-2}
            \put(4.9,-0.075){\tiny \em n-1}
            \put(5.95,-0.075){\tiny \em n}
           \end{picture}}
\end{picture}}

\newcommand{\BnDNGGraph}{\setlength{\unitlength}{0.75in}
\begin{picture}(6,0.55)
\put(0,0){\begin{picture}(1,0)
            \put(0,0.1){\circle*{0.075}}
            \put(1,0.1){\circle*{0.075}}
            \put(2,0.1){\circle*{0.075}}
            \put(4,0.1){\circle*{0.075}}
            \put(5,0.1){\circle*{0.075}}
            \put(6,0.1){\circle*{0.075}}
            \put(0,0.1){\line(1,0){2}}
            \multiput(2,0.1)(0.4,0){5}{\line(1,0){0.2}}
            \put(4,0.1){\line(1,0){2}}
            \put(-0.05,-0.075){\tiny \em 1}
            \put(0.95,-0.075){\tiny \em 2}
            \put(1.95,-0.075){\tiny \em 3}
            \put(3.9,-0.075){\tiny \em n-2}
            \put(4.9,-0.075){\tiny \em n-1}
            \put(5.95,-0.075){\tiny \em n}
            \put(5.3,0.15){\CircleInteger{4}}
           \end{picture}}
\end{picture}}

\newcommand{\BnINGGraph}{\setlength{\unitlength}{0.75in}
\begin{picture}(6,0.55)
\put(0,0){\begin{picture}(1,0)
            \put(0,0.1){\circle*{0.075}}
            \put(1,0.1){\circle*{0.075}}
            \put(2,0.1){\circle*{0.075}}
            \put(4,0.1){\circle*{0.075}}
            \put(5,0.1){\circle*{0.075}}
            \put(6,0.1){\circle*{0.075}}
            \put(0,0.1){\line(1,0){2}}
            \multiput(2,0.1)(0.4,0){5}{\line(1,0){0.2}}
            \put(4,0.1){\line(1,0){2}}
            \put(-0.05,-0.075){\tiny \em 1}
            \put(0.95,-0.075){\tiny \em 2}
            \put(1.95,-0.075){\tiny \em 3}
            \put(3.9,-0.075){\tiny \em n-2}
            \put(4.9,-0.075){\tiny \em n-1}
            \put(5.95,-0.075){\tiny \em n}
            \put(5.2,0.1){\vector(1,0){0.1}}
            \put(5.3,0.1){\vector(1,0){0.1}}
            \put(5.8,0.1){\vector(-1,0){0.1}}
           \end{picture}}
\end{picture}}

\newcommand{\CnINGGraph}{\setlength{\unitlength}{0.75in}
\begin{picture}(6,0.55)
\put(0,0){\begin{picture}(1,0)
            \put(0,0.1){\circle*{0.075}}
            \put(1,0.1){\circle*{0.075}}
            \put(2,0.1){\circle*{0.075}}
            \put(4,0.1){\circle*{0.075}}
            \put(5,0.1){\circle*{0.075}}
            \put(6,0.1){\circle*{0.075}}
            \put(0,0.1){\line(1,0){2}}
            \multiput(2,0.1)(0.4,0){5}{\line(1,0){0.2}}
            \put(4,0.1){\line(1,0){2}}
            \put(-0.05,-0.075){\tiny \em 1}
            \put(0.95,-0.075){\tiny \em 2}
            \put(1.95,-0.075){\tiny \em 3}
            \put(3.9,-0.075){\tiny \em n-2}
            \put(4.9,-0.075){\tiny \em n-1}
            \put(5.95,-0.075){\tiny \em n}
            \put(5.2,0.1){\vector(1,0){0.1}}
            \put(5.7,0.1){\vector(-1,0){0.1}}
            \put(5.8,0.1){\vector(-1,0){0.1}}
           \end{picture}}
\end{picture}}

\newcommand{\DnDNGGraph}{\setlength{\unitlength}{0.75in}
\begin{picture}(6,0.75)
\put(0,-0.25){\begin{picture}(1,0)
            \put(0,0.35){\circle*{0.075}}
            \put(1,0.35){\circle*{0.075}}
            \put(2,0.35){\circle*{0.075}}
            \put(4,0.35){\circle*{0.075}}
            \put(5,0.35){\circle*{0.075}}
            \put(6,0.1){\circle*{0.075}}
            \put(6,0.6){\circle*{0.075}}
            \put(0,0.35){\line(1,0){2}}
            \multiput(2,0.35)(0.4,0){5}{\line(1,0){0.2}}
            \put(4,0.35){\line(1,0){1}}
            \put(5,0.35){\line(4,1){1}}
            \put(5,0.35){\line(4,-1){1}}
            \put(-0.05,0.175){\tiny \em 1}
            \put(0.95,0.175){\tiny \em 2}
            \put(1.95,0.175){\tiny \em 3}
            \put(3.9,0.175){\tiny \em n-3}
            \put(4.85,0.175){\tiny \em n-2}
            \put(5.9,0.425){\tiny \em n-1}
            \put(5.95,-0.075){\tiny \em n}
           \end{picture}}
\end{picture}}

\newcommand{\EEightDNGGraph}{\setlength{\unitlength}{0.75in}
\begin{picture}(6,0.75)
\put(0,-0.25){\begin{picture}(1,0)
            \put(0,0.1){\circle*{0.075}}
            \put(1,0.1){\circle*{0.075}}
            \put(2,0.1){\circle*{0.075}}
            \put(2,0.6){\circle*{0.075}}
            \put(3,0.1){\circle*{0.075}}
            \put(4,0.1){\circle*{0.075}}
            \put(5,0.1){\circle*{0.075}}
            \put(6,0.1){\circle*{0.075}}
            \put(0,0.1){\line(1,0){6}}
            \put(2,0.1){\line(0,1){0.5}}
            \put(-0.05,-0.075){\tiny \em 1}
            \put(2.05,0.55){\tiny \em 2}
            \put(0.95,-0.075){\tiny \em 3}
            \put(1.95,-0.075){\tiny \em 4}
            \put(2.95,-0.075){\tiny \em 5}
            \put(3.95,-0.075){\tiny \em 6}
            \put(4.95,-0.075){\tiny \em 7}
            \put(5.95,-0.075){\tiny \em 8}
           \end{picture}}
\end{picture}}

\newcommand{\ESevenDNGGraph}{\setlength{\unitlength}{0.75in}
\begin{picture}(6,0.75)
\put(0,-0.25){\begin{picture}(1,0)
            \put(0,0.1){\circle*{0.075}}
            \put(1,0.1){\circle*{0.075}}
            \put(2,0.1){\circle*{0.075}}
            \put(2,0.6){\circle*{0.075}}
            \put(3,0.1){\circle*{0.075}}
            \put(4,0.1){\circle*{0.075}}
            \put(5,0.1){\circle*{0.075}}
            \put(0,0.1){\line(1,0){5}}
            \put(2,0.1){\line(0,1){0.5}}
            \put(-0.05,-0.075){\tiny \em 1}
            \put(2.05,0.55){\tiny \em 2}
            \put(0.95,-0.075){\tiny \em 3}
            \put(1.95,-0.075){\tiny \em 4}
            \put(2.95,-0.075){\tiny \em 5}
            \put(3.95,-0.075){\tiny \em 6}
            \put(4.95,-0.075){\tiny \em 7}
           \end{picture}}
\end{picture}}

\newcommand{\ESixDNGGraph}{\setlength{\unitlength}{0.75in}
\begin{picture}(6,0.75)
\put(0,-0.25){\begin{picture}(1,0)
            \put(0,0.1){\circle*{0.075}}
            \put(1,0.1){\circle*{0.075}}
            \put(2,0.1){\circle*{0.075}}
            \put(2,0.6){\circle*{0.075}}
            \put(3,0.1){\circle*{0.075}}
            \put(4,0.1){\circle*{0.075}}
            \put(0,0.1){\line(1,0){4}}
            \put(2,0.1){\line(0,1){0.5}}
            \put(-0.05,-0.075){\tiny \em 1}
            \put(2.05,0.55){\tiny \em 2}
            \put(0.95,-0.075){\tiny \em 3}
            \put(1.95,-0.075){\tiny \em 4}
            \put(2.95,-0.075){\tiny \em 5}
            \put(3.95,-0.075){\tiny \em 6}
           \end{picture}}
\end{picture}}

\newcommand{\FFourDNGGraph}{\setlength{\unitlength}{0.75in}
\begin{picture}(1,0.85)
\put(0,0){\begin{picture}(1,0)
            \put(0,0.1){\circle*{0.075}}
            \put(1,0.1){\circle*{0.075}}
            \put(2,0.1){\circle*{0.075}}
            \put(3,0.1){\circle*{0.075}}
            \put(0,0.1){\line(1,0){3}}
            \put(-0.05,-0.075){\tiny \em 1}
            \put(0.95,-0.075){\tiny \em 2}
            \put(1.95,-0.075){\tiny \em 3}
            \put(2.95,-0.075){\tiny \em 4}
            \put(1.3,0.15){\CircleInteger{4}}
            \end{picture}}
\end{picture}}

\newcommand{\FFourINGGraph}{\setlength{\unitlength}{0.75in}
\begin{picture}(1,0.85)
\put(0,0){\begin{picture}(1,0)
            \put(0,0.1){\circle*{0.075}}
            \put(1,0.1){\circle*{0.075}}
            \put(2,0.1){\circle*{0.075}}
            \put(3,0.1){\circle*{0.075}}
            \put(0,0.1){\line(1,0){3}}
            \put(-0.05,-0.075){\tiny \em 1}
            \put(0.95,-0.075){\tiny \em 2}
            \put(1.95,-0.075){\tiny \em 3}
            \put(2.95,-0.075){\tiny \em 4}
            \put(1.2,0.1){\vector(1,0){0.1}}
            \put(1.3,0.1){\vector(1,0){0.1}}
            \put(1.8,0.1){\vector(-1,0){0.1}}
            \end{picture}}
\end{picture}}

\newcommand{\GTwoDNGGraph}{\setlength{\unitlength}{0.75in}
\begin{picture}(1,0.55)
\put(0,0){\begin{picture}(1,0)
            \put(0,0.1){\circle*{0.075}}
            \put(1,0.1){\circle*{0.075}}
            \put(0,0.1){\line(1,0){1}}
            \put(-0.05,-0.075){\tiny \em 1}
            \put(0.95,-0.075){\tiny \em 2}
            \put(0.3,0.15){\CircleInteger{6}}
            \end{picture}}
\end{picture}}

\newcommand{\GTwoINGGraph}{\setlength{\unitlength}{0.75in}
\begin{picture}(1,0.55)
\put(0,0){\begin{picture}(1,0)
            \put(0,0.1){\circle*{0.075}}
            \put(1,0.1){\circle*{0.075}}
            \put(0,0.1){\line(1,0){1}}
            \put(0.2,0.1){\vector(1,0){0.1}}
            \put(0.6,0.1){\vector(-1,0){0.1}}
            \put(0.7,0.1){\vector(-1,0){0.1}}
            \put(0.8,0.1){\vector(-1,0){0.1}}
            \put(-0.05,-0.075){\tiny \em 1}
            \put(0.95,-0.075){\tiny \em 2}
            \end{picture}}
\end{picture}}

\newcommand{\HFourDNGGraph}{\setlength{\unitlength}{0.75in}
\begin{picture}(1,0.55)
\put(0,0){\begin{picture}(1,0)
            \put(0,0.1){\circle*{0.075}}
            \put(1,0.1){\circle*{0.075}}
            \put(2,0.1){\circle*{0.075}}
            \put(3,0.1){\circle*{0.075}}
            \put(0,0.1){\line(1,0){3}}
            \put(-0.05,-0.075){\tiny \em 1}
            \put(0.95,-0.075){\tiny \em 2}
            \put(1.95,-0.075){\tiny \em 3}
            \put(2.95,-0.075){\tiny \em 4}
            \put(2.3,0.15){\CircleInteger{5}}
            \end{picture}}
\end{picture}}

\newcommand{\HThreeDNGGraph}{\setlength{\unitlength}{0.75in}
\begin{picture}(1,0.55)
\put(0,0){\begin{picture}(1,0)
            \put(0,0.1){\circle*{0.075}}
            \put(1,0.1){\circle*{0.075}}
            \put(2,0.1){\circle*{0.075}}
            \put(0,0.1){\line(1,0){2}}
            \put(-0.05,-0.075){\tiny \em 1}
            \put(0.95,-0.075){\tiny \em 2}
            \put(1.95,-0.075){\tiny \em 3}
            \put(1.3,0.15){\CircleInteger{5}}
            \end{picture}}
\end{picture}}

\newcommand{\HTwoDNGGraph}{\setlength{\unitlength}{0.75in}
\begin{picture}(1,0.55)
\put(0,0){\begin{picture}(1,0)
            \put(0,0.1){\circle*{0.075}}
            \put(1,0.1){\circle*{0.075}}
            \put(0,0.1){\line(1,0){1}}
            \put(-0.05,-0.075){\tiny \em 1}
            \put(0.95,-0.075){\tiny \em 2}
            \put(0.3,0.15){\CircleInteger{5}}
            \end{picture}}
\end{picture}}

\newcommand{\ITwoDNGGraph}{\setlength{\unitlength}{0.75in}
\begin{picture}(1,0.55)
\put(0,0){\begin{picture}(1,0)
            \put(0,0.1){\circle*{0.075}}
            \put(1,0.1){\circle*{0.075}}
            \put(0,0.1){\line(1,0){1}}
            \put(-0.05,-0.075){\tiny \em 1}
            \put(0.95,-0.075){\tiny \em 2}
            \put(0.3,0.15){\CircleIntegerm}
            \end{picture}}
\end{picture}}
\begin{figure}[t]
\begin{center}
{\bf \IEGGraphFigure}\ \  The integral Coxeter--Dynkin flowers of {\em La Florado Klasado}.\\ 
{\footnotesize (Nodes are numbered as in \cite{Hum}.)} 
\vspace*{0.04in}
\end{center}

\vspace*{-0.45in}
\hspace*{0.25in}
\begin{tabular}{cl}
$\myA_{n}$ ($n \geq 1$) & \AnDNGGraph\\

$\myB_{n}$ ($n \geq 3$) & \BnINGGraph\\

$\myC_{n}$ ($n \geq 2$) & \CnINGGraph\\

$\myD_{n}$ ($n \geq 4$) & \DnDNGGraph\\

$\myE_{6}$ & \ESixDNGGraph\\

$\myE_{7}$ & \ESevenDNGGraph\\

$\myE_{8}$ & \EEightDNGGraph\\

$\myF_{4}$ & \FFourINGGraph\\

$\myG_{2}$ & \GTwoINGGraph
\end{tabular}
\end{figure}

{\bf [\S \GStructureIntroSection.4:\! Sturdy posets.]} 
For our work here, embryophytes mainly occur to supply certain structure to the DCML's and DCDL's we consider in the remainder of this manuscript. 
That said, they also directly connect these lattices to the algebraic milieux (semisimple Lie algebra representations, Weyl symmetric functions) that is a large part of our motivation.  
To begin making these connections, let $\mathscr{G} = (\Gamma_{I},M_{I \times I})$ be an IEG, and let $R$ be an $I$-colored ranked poset. 
Given $\xelt \in R$, we let $\rho_{i}$ denote the rank function on $\comp_{i}(\xelt)$ and $\delta_{i}$ the depth function. 
So, the length of $\comp_{i}(\xelt)$ is $l_{i} = l_{i}(\xelt) = \rho_{i}(\xelt) + \delta_{i}(\xelt)$. 
We define the $i^{\mbox{\tiny th}}$ {\em centered coordinate} or {\em color $i$ weight of} $\xelt$ to be the quantity $\mym_{i}(\xelt) := \rho_{i}(\xelt) - \delta_{i}(\xelt)$. 
The {\em weight}\footnote{In \S \FlowerSection\ below, we express this vector in terms of a `fundamental' basis $\{\omega_{i}\}_{i \in I}$, so that $wt(\xelt) = \sum_{i \in I}\mysmallm_{i}(\xelt)\omega_{i}$.} of $\xelt$ is the vector $I$-tuple $wt(\xelt) := \left(\rule[-1.5mm]{-0.1mm}{4.75mm}\mym_{i}(\xelt)\right)_{i \in I}$. 
When $\mym_{i}(\xelt) \geq 0$ for all $i \in I$ (as happens, for example, when $\xelt$ is uniquely maximal in $R$), we say that $wt(\xelt)$ is {\em dominant}. 
\hfill Say $R$ is $\mathscr{G}$-structured if, whenever \ $\xelt \myarrow{i} \yelt$\  in $R$, we have $wt(\xelt) + \myroot_{i} = wt(\yelt)$,

\newcommand{\ESixFlower}{
\setlength{\unitlength}{0.5cm}
\begin{picture}(15,3.5)
\put(-2.5,1.8){\LARGE $\myE_{6}$}
\put(-1,2.1){\vector(3,-1){1.75}}
{\color{Cyan}\put(1,0){\circle*{0.7}}}
\put(1,0){\circle{0.7}}
{\color{Red}\put(4,0){\circle*{0.7}}}
\put(4,0){\circle{0.7}}
{\color{Purple}\put(7,0){\circle*{0.7}}}
\put(7,0){\circle{0.7}}
{\color{OliveGreen}\put(10,0){\circle*{0.7}}}
\put(10,0){\circle{0.7}}
{\color{BurntOrange}\put(13,0){\circle*{0.7}}}
\put(13,0){\circle{0.7}}
{\color{SpringGreen}\put(7,3){\circle*{0.7}}}
\put(7,3){\circle{0.7}}
\thicklines
\put(1.35,0){\line(1,0){2.3}}
\put(4.35,0){\line(1,0){2.3}}
\put(7.35,0){\line(1,0){2.3}}
\put(10.35,0){\line(1,0){2.3}}
\put(7,0.35){\line(0,1){2.3}}
\put(0.8,0.6){\textcolor{Cyan}{\em 1}}
\put(3.8,0.6){\textcolor{Red}{\em 3}}
\put(6.3,0.6){\textcolor{Purple}{\em 4}}
\put(9.8,0.6){\textcolor{OliveGreen}{\em 5}}
\put(12.8,0.6){\textcolor{BurntOrange}{\em 6}}
\put(6.3,2.1){\textcolor{SpringGreen}{\em 2}}
\end{picture}
}

\begin{figure}
\begin{center}
{\bf \ESixLatticeFigure.1}\ \  Two $\myE_{6}$-structured DCDL's.\\
{\small (These are the two ``$\mysmallE_{6}$-minuscule splitting DCDL's'', cf.\ \S \MinusculeExampleSection.)}

\vspace*{0.1in} 
\ESixFlower

\vspace*{0.1in} 
\setlength{\unitlength}{1cm}
\begin{picture}(6,17)
\put(-0.9,12.2){\LARGE $L_{\mysmallE_{6}}(\omega_{\mbox{\color{Cyan}{\small $1$}}})$}
\put(1.1,12.6){\vector(2,1){1}}
\put(0,16){\TypeEboxDot{Gray}}
\put(0,8){\TypeEboxDot{Gray}}
\put(0,0){\TypeEboxDot{Gray}}
\put(1,15){\TypeEboxDot{Gray}}
\put(1,9){\TypeEboxDot{Gray}}
\put(1,7){\TypeEboxDot{Gray}}
\put(1,1){\TypeEboxDot{Gray}}
\put(2,14){\TypeEboxDot{Gray}}
\put(2,10){\TypeEboxDot{Gray}}
\put(2,8){\TypeEboxDot{Gray}}
\put(2,6){\TypeEboxDot{Gray}}
\put(2,2){\TypeEboxDot{Gray}}
\put(3,13){\TypeEboxDot{Gray}}
\put(3,11){\TypeEboxDot{Gray}}
\put(3,9){\TypeEboxDot{Gray}}
\put(3,7){\TypeEboxDot{Gray}}
\put(3,5){\TypeEboxDot{Gray}}
\put(3,3){\TypeEboxDot{Gray}}
\put(2,12){\TypeEboxDot{Gray}}
\put(2,4){\TypeEboxDot{Gray}}
\put(4,12){\TypeEboxDot{Gray}}
\put(4,10){\TypeEboxDot{Gray}}
\put(4,6){\TypeEboxDot{Gray}}
\put(4,4){\TypeEboxDot{Gray}}
\put(5,11){\TypeEboxDot{Gray}}
\put(5,5){\TypeEboxDot{Gray}}
\thicklines
\put(0.625,0.375){\color{BurntOrange}\qbezier(0,0)(0.375,0.375)(0.75,0.75)}
\put(0.9,0.65){\color{BurntOrange}{\em 6}}
\put(1.625,1.375){\color{OliveGreen}\qbezier(0,0)(0.375,0.375)(0.75,0.75)}
\put(1.9,1.65){\color{OliveGreen}{\em 5}}
\put(2.625,2.375){\color{Purple}\qbezier(0,0)(0.375,0.375)(0.75,0.75)}
\put(2.9,2.65){\color{Purple}{\em 4}}
\put(3.625,3.375){\color{Red}\qbezier(0,0)(0.375,0.375)(0.75,0.75)}
\put(3.9,3.65){\color{Red}{\em 3}}
\put(4.625,4.375){\color{Cyan}\qbezier(0,0)(0.375,0.375)(0.75,0.75)}
\put(4.9,4.65){\color{Cyan}{\em 1}}
\put(2.625,4.375){\color{Red}\qbezier(0,0)(0.375,0.375)(0.75,0.75)}
\put(2.9,4.65){\color{Red}{\em 3}}
\put(3.625,5.375){\color{Cyan}\qbezier(0,0)(0.375,0.375)(0.75,0.75)}
\put(3.9,5.65){\color{Cyan}{\em 1}}
\put(2.625,6.375){\color{Cyan}\qbezier(0,0)(0.375,0.375)(0.75,0.75)}
\put(2.9,6.65){\color{Cyan}{\em 1}}
\put(1.625,7.375){\color{Cyan}\qbezier(0,0)(0.375,0.375)(0.75,0.75)}
\put(1.9,7.65){\color{Cyan}{\em 1}}
\put(2.625,8.375){\color{Red}\qbezier(0,0)(0.375,0.375)(0.75,0.75)}
\put(2.9,8.65){\color{Red}{\em 3}}
\put(3.625,9.375){\color{Purple}\qbezier(0,0)(0.375,0.375)(0.75,0.75)}
\put(3.9,9.65){\color{Purple}{\em 4}}
\put(4.625,10.375){\color{SpringGreen}\qbezier(0,0)(0.375,0.375)(0.75,0.75)}
\put(4.9,10.65){\color{SpringGreen}{\em 2}}
\put(0.625,8.375){\color{Cyan}\qbezier(0,0)(0.375,0.375)(0.75,0.75)}
\put(0.9,8.65){\color{Cyan}{\em 1}}
\put(1.625,9.375){\color{Red}\qbezier(0,0)(0.375,0.375)(0.75,0.75)}
\put(1.9,9.65){\color{Red}{\em 3}}
\put(2.625,10.375){\color{Purple}\qbezier(0,0)(0.375,0.375)(0.75,0.75)}
\put(2.9,10.65){\color{Purple}{\em 4}}
\put(3.625,11.375){\color{SpringGreen}\qbezier(0,0)(0.375,0.375)(0.75,0.75)}
\put(3.9,11.65){\color{SpringGreen}{\em 2}}
\put(2.625,12.375){\color{SpringGreen}\qbezier(0,0)(0.375,0.375)(0.75,0.75)}
\put(2.9,12.65){\color{SpringGreen}{\em 2}}
\put(3.39,3.375){\color{SpringGreen}\qbezier(0,0)(-0.375,0.375)(-0.75,0.75)}
\put(2.9,3.65){\color{SpringGreen}{\em 2}}
\put(4.39,4.375){\color{SpringGreen}\qbezier(0,0)(-0.375,0.375)(-0.75,0.75)}
\put(3.9,4.65){\color{SpringGreen}{\em 2}}
\put(3.39,5.375){\color{Purple}\qbezier(0,0)(-0.375,0.375)(-0.75,0.75)}
\put(2.9,5.65){\color{Purple}{\em 4}}
\put(2.39,6.375){\color{OliveGreen}\qbezier(0,0)(-0.375,0.375)(-0.75,0.75)}
\put(1.9,6.65){\color{OliveGreen}{\em 5}}
\put(1.39,7.375){\color{BurntOrange}\qbezier(0,0)(-0.375,0.375)(-0.75,0.75)}
\put(0.9,7.65){\color{BurntOrange}{\em 6}}
\put(5.39,5.375){\color{SpringGreen}\qbezier(0,0)(-0.375,0.375)(-0.75,0.75)}
\put(4.9,5.65){\color{SpringGreen}{\em 2}}
\put(4.39,6.375){\color{Purple}\qbezier(0,0)(-0.375,0.375)(-0.75,0.75)}
\put(3.9,6.65){\color{Purple}{\em 4}}
\put(3.39,7.375){\color{OliveGreen}\qbezier(0,0)(-0.375,0.375)(-0.75,0.75)}
\put(2.9,7.65){\color{OliveGreen}{\em 5}}
\put(2.39,8.375){\color{BurntOrange}\qbezier(0,0)(-0.375,0.375)(-0.75,0.75)}
\put(1.9,8.65){\color{BurntOrange}{\em 6}}
\put(3.39,9.375){\color{BurntOrange}\qbezier(0,0)(-0.375,0.375)(-0.75,0.75)}
\put(2.9,9.65){\color{BurntOrange}{\em 6}}
\put(4.39,10.375){\color{BurntOrange}\qbezier(0,0)(-0.375,0.375)(-0.75,0.75)}
\put(3.9,10.65){\color{BurntOrange}{\em 6}}
\put(3.39,11.375){\color{OliveGreen}\qbezier(0,0)(-0.375,0.375)(-0.75,0.75)}
\put(2.9,11.65){\color{OliveGreen}{\em 5}}
\put(5.39,11.375){\color{BurntOrange}\qbezier(0,0)(-0.375,0.375)(-0.75,0.75)}
\put(4.9,11.65){\color{BurntOrange}{\em 6}}
\put(4.39,12.375){\color{OliveGreen}\qbezier(0,0)(-0.375,0.375)(-0.75,0.75)}
\put(3.9,12.65){\color{OliveGreen}{\em 5}}
\put(3.39,13.375){\color{Purple}\qbezier(0,0)(-0.375,0.375)(-0.75,0.75)}
\put(2.9,13.65){\color{Purple}{\em 4}}
\put(2.39,14.375){\color{Red}\qbezier(0,0)(-0.375,0.375)(-0.75,0.75)}
\put(1.9,14.65){\color{Red}{\em 3}}
\put(1.39,15.375){\color{Cyan}\qbezier(0,0)(-0.375,0.375)(-0.75,0.75)}
\put(0.9,15.65){\color{Cyan}{\em 1}}
\put(4,8){\TypeEboxDot{Gray}}
\put(4.39,8.375){\color{OliveGreen}\qbezier(0,0)(-0.375,0.375)(-0.75,0.75)}
\put(3.9,8.65){\color{OliveGreen}{\em 5}}
\put(3.625,7.375){\color{Red}\qbezier(0,0)(0.375,0.375)(0.75,0.75)}
\put(3.9,7.65){\color{Red}{\em 3}}
\end{picture}
\begin{picture}(6,17)
\put(5,15){\LARGE $L_{\mysmallE_{6}}(\omega_{\mbox{\color{BurntOrange}{\small $6$}}})$}
\put(4.85,15.1){\vector(-2,-1){1.5}}
\put(5.5,14){\LARGE $\cong L_{\mysmallE_{6}}(\omega_{\mbox{\color{Cyan}{\small $1$}}})^{*}$}
\put(0,16){\TypeEboxDot{Gray}}
\put(0,8){\TypeEboxDot{Gray}}
\put(0,0){\TypeEboxDot{Gray}}
\put(1,15){\TypeEboxDot{Gray}}
\put(1,9){\TypeEboxDot{Gray}}
\put(1,7){\TypeEboxDot{Gray}}
\put(1,1){\TypeEboxDot{Gray}}
\put(2,14){\TypeEboxDot{Gray}}
\put(2,10){\TypeEboxDot{Gray}}
\put(2,8){\TypeEboxDot{Gray}}
\put(2,6){\TypeEboxDot{Gray}}
\put(2,2){\TypeEboxDot{Gray}}
\put(3,13){\TypeEboxDot{Gray}}
\put(3,11){\TypeEboxDot{Gray}}
\put(3,9){\TypeEboxDot{Gray}}
\put(3,7){\TypeEboxDot{Gray}}
\put(3,5){\TypeEboxDot{Gray}}
\put(3,3){\TypeEboxDot{Gray}}
\put(2,12){\TypeEboxDot{Gray}}
\put(2,4){\TypeEboxDot{Gray}}
\put(4,12){\TypeEboxDot{Gray}}
\put(4,10){\TypeEboxDot{Gray}}
\put(4,6){\TypeEboxDot{Gray}}
\put(4,4){\TypeEboxDot{Gray}}
\put(5,11){\TypeEboxDot{Gray}}
\put(5,5){\TypeEboxDot{Gray}}
\thicklines
\put(0.625,0.375){\color{Cyan}\qbezier(0,0)(0.375,0.375)(0.75,0.75)}
\put(0.9,0.65){\color{Cyan}{\em 1}}
\put(1.625,1.375){\color{Red}\qbezier(0,0)(0.375,0.375)(0.75,0.75)}
\put(1.9,1.65){\color{Red}{\em 3}}
\put(2.625,2.375){\color{Purple}\qbezier(0,0)(0.375,0.375)(0.75,0.75)}
\put(2.9,2.65){\color{Purple}{\em 4}}
\put(3.625,3.375){\color{OliveGreen}\qbezier(0,0)(0.375,0.375)(0.75,0.75)}
\put(3.9,3.65){\color{OliveGreen}{\em 5}}
\put(4.625,4.375){\color{BurntOrange}\qbezier(0,0)(0.375,0.375)(0.75,0.75)}
\put(4.9,4.65){\color{BurntOrange}{\em 6}}
\put(2.625,4.375){\color{OliveGreen}\qbezier(0,0)(0.375,0.375)(0.75,0.75)}
\put(2.9,4.65){\color{OliveGreen}{\em 5}}
\put(3.625,5.375){\color{BurntOrange}\qbezier(0,0)(0.375,0.375)(0.75,0.75)}
\put(3.9,5.65){\color{BurntOrange}{\em 6}}
\put(2.625,6.375){\color{BurntOrange}\qbezier(0,0)(0.375,0.375)(0.75,0.75)}
\put(2.9,6.65){\color{BurntOrange}{\em 6}}
\put(1.625,7.375){\color{BurntOrange}\qbezier(0,0)(0.375,0.375)(0.75,0.75)}
\put(1.9,7.65){\color{BurntOrange}{\em 6}}
\put(2.625,8.375){\color{OliveGreen}\qbezier(0,0)(0.375,0.375)(0.75,0.75)}
\put(2.9,8.65){\color{OliveGreen}{\em 5}}
\put(3.625,9.375){\color{Purple}\qbezier(0,0)(0.375,0.375)(0.75,0.75)}
\put(3.9,9.65){\color{Purple}{\em 4}}
\put(4.625,10.375){\color{SpringGreen}\qbezier(0,0)(0.375,0.375)(0.75,0.75)}
\put(4.9,10.65){\color{SpringGreen}{\em 2}}
\put(0.625,8.375){\color{BurntOrange}\qbezier(0,0)(0.375,0.375)(0.75,0.75)}
\put(0.9,8.65){\color{BurntOrange}{\em 6}}
\put(1.625,9.375){\color{OliveGreen}\qbezier(0,0)(0.375,0.375)(0.75,0.75)}
\put(1.9,9.65){\color{OliveGreen}{\em 5}}
\put(2.625,10.375){\color{Purple}\qbezier(0,0)(0.375,0.375)(0.75,0.75)}
\put(2.9,10.65){\color{Purple}{\em 4}}
\put(3.625,11.375){\color{SpringGreen}\qbezier(0,0)(0.375,0.375)(0.75,0.75)}
\put(3.9,11.65){\color{SpringGreen}{\em 2}}
\put(2.625,12.375){\color{SpringGreen}\qbezier(0,0)(0.375,0.375)(0.75,0.75)}
\put(2.9,12.65){\color{SpringGreen}{\em 2}}
\put(3.39,3.375){\color{SpringGreen}\qbezier(0,0)(-0.375,0.375)(-0.75,0.75)}
\put(2.9,3.65){\color{SpringGreen}{\em 2}}
\put(4.39,4.375){\color{SpringGreen}\qbezier(0,0)(-0.375,0.375)(-0.75,0.75)}
\put(3.9,4.65){\color{SpringGreen}{\em 2}}
\put(3.39,5.375){\color{Purple}\qbezier(0,0)(-0.375,0.375)(-0.75,0.75)}
\put(2.9,5.65){\color{Purple}{\em 4}}
\put(2.39,6.375){\color{Red}\qbezier(0,0)(-0.375,0.375)(-0.75,0.75)}
\put(1.9,6.65){\color{Red}{\em 3}}
\put(1.39,7.375){\color{Cyan}\qbezier(0,0)(-0.375,0.375)(-0.75,0.75)}
\put(0.9,7.65){\color{Cyan}{\em 1}}
\put(5.39,5.375){\color{SpringGreen}\qbezier(0,0)(-0.375,0.375)(-0.75,0.75)}
\put(4.9,5.65){\color{SpringGreen}{\em 2}}
\put(4.39,6.375){\color{Purple}\qbezier(0,0)(-0.375,0.375)(-0.75,0.75)}
\put(3.9,6.65){\color{Purple}{\em 4}}
\put(3.39,7.375){\color{Red}\qbezier(0,0)(-0.375,0.375)(-0.75,0.75)}
\put(2.9,7.65){\color{Red}{\em 3}}
\put(2.39,8.375){\color{Cyan}\qbezier(0,0)(-0.375,0.375)(-0.75,0.75)}
\put(1.9,8.65){\color{Cyan}{\em 1}}
\put(3.39,9.375){\color{Cyan}\qbezier(0,0)(-0.375,0.375)(-0.75,0.75)}
\put(2.9,9.65){\color{Cyan}{\em 1}}
\put(4.39,10.375){\color{Cyan}\qbezier(0,0)(-0.375,0.375)(-0.75,0.75)}
\put(3.9,10.65){\color{Cyan}{\em 1}}
\put(3.39,11.375){\color{Red}\qbezier(0,0)(-0.375,0.375)(-0.75,0.75)}
\put(2.9,11.65){\color{Red}{\em 3}}
\put(5.39,11.375){\color{Cyan}\qbezier(0,0)(-0.375,0.375)(-0.75,0.75)}
\put(4.9,11.65){\color{Cyan}{\em 1}}
\put(4.39,12.375){\color{Red}\qbezier(0,0)(-0.375,0.375)(-0.75,0.75)}
\put(3.9,12.65){\color{Red}{\em 3}}
\put(3.39,13.375){\color{Purple}\qbezier(0,0)(-0.375,0.375)(-0.75,0.75)}
\put(2.9,13.65){\color{Purple}{\em 4}}
\put(2.39,14.375){\color{OliveGreen}\qbezier(0,0)(-0.375,0.375)(-0.75,0.75)}
\put(1.9,14.65){\color{OliveGreen}{\em 5}}
\put(1.39,15.375){\color{BurntOrange}\qbezier(0,0)(-0.375,0.375)(-0.75,0.75)}
\put(0.9,15.65){\color{BurntOrange}{\em 6}}
\put(4,8){\TypeEboxDot{Gray}}
\put(4.39,8.375){\color{Red}\qbezier(0,0)(-0.375,0.375)(-0.75,0.75)}
\put(3.9,8.65){\color{Red}{\em 3}}
\put(3.625,7.375){\color{OliveGreen}\qbezier(0,0)(0.375,0.375)(0.75,0.75)}
\put(3.9,7.65){\color{OliveGreen}{\em 5}}
\end{picture}
\end{center}
\end{figure}

\clearpage 
\noindent 
which is equivalent to $\mym_{j}(\xelt) + M_{ij} = \mym_{j}(\yelt)$ for all $j \not= i$. 
If some $I$-edge-colored DCDL $\Jcolor(P)$ is $\mathscr{G}$-structured, then we say its $I$-vertex-colored compression poset $P$ is $\mathscr{G}$-structured also. 
(See \ESixFigures\ for some examples.) 
Say $R$ is {\em sturdy} if $R$ is $\mathscr{G}$-structured for some embryophyte $\mathscr{G}$. 
See \cite{DonPosetModels} for some other very basic results about sturdy posets. 

This author developed parts {\sl (1)} and {\sl (2)} of the following result in \cite{DE} (see Theorem 5.2); these constrain, under certain assumptions, the possible embryophytes that can actually arise amongst sturdy posets. 
We defer the proof of part {\sl (3)} to \S \FlowerSection.5, where it is a corollary of \SaturatedLemma. 

\noindent
{\bf \GCMTheorem\ (D.)}\ \ 
{\sl Let $\mathscr{G}_{I} = (\Gamma_{I},M_{I \times I})$ be an integral embryophytic graph. 
Suppose there exists a finite and $\mathscr{G}$-structured ranked poset $R$. 
(1) If} $\ecolor_{R}(\EdgeSet(R)) = I$, {\sl then $\mathscr{G}$ is an integral Coxeter--Dynkin posy. 
(2) If $\mathscr{G}$ is connected and if $R$ has at least one edge, then} $\ecolor_{R}(\EdgeSet(R)) = I$ {\sl and $\mathscr{G}$ is an integral Coxeter--Dynkin flower. 
(3) If $\mathscr{G}$ is an integral Coxeter--Dynkin flower and the vector $\lambda = (\lambda_{i})_{i \in I}$ is an $I$-tuple of nonnegative integers, then there exists a $\mathscr{G}$-structured ranked poset whose unique maximal element has weight $\lambda$.}

{\bf [\S \GStructureIntroSection.5:\! Sturdy DCML's.]}  As an application of \GCMTheorem, we obtain the following new Dynkin diagram classification result. 

\noindent 
{\bf \NewLFKTheorem}\ \ {\sl Let $M = (M_{ij})_{i,j\in I}$ be a matrix of real numbers with associated graph $\Gamma(M)$, and suppose $L$ is an $I$-colored DCML such that, for all $i,j \in I$ and $\xelt, \yelt \in L$, we have $\mym_{j}(\xelt) + M_{ij} = \mym_{j}(\yelt)$ whenever $\xelt \myarrow{i} \yelt$. 
(1) If} $\ecolor_{L}(\EdgeSet(L)) = I$, {\sl then $M$ is a pulsation matrix and the embryophyte $\mathscr{G} = (\Gamma(M),M)$ is an integral Coxeter--Dynkin posy. 
(2) If $\Gamma(M)$ is connected and if $L$ has at least one edge, then} $\ecolor_{L}(\EdgeSet(L)) = I$, {\sl $M$ is a pulsation matrix, and $\mathscr{G} = (\Gamma(M),M)$ is an integral Coxeter--Dynkin flower.
(3) If $\mathscr{H}$ is any integral Coxeter--Dynkin flower, then there exists an $\mathscr{H}$-structured diamond-colored modular lattice with at least one edge.}

{\em Proof.} For {\sl (1)}, we show that $M$ is a pulsation matrix of integers and invoke \GCMTheorem. 
Let $i, j \in I$ and pick an edge $\selt \myarrow{i} \telt$. 
Since $wt(\selt) + \myroot_{i} = wt(\telt)$, then $\mym_{j}(\selt)+M_{ij}=\mym_{j}(\telt)$. 
So, $M_{ij}$ is an integer. 
When $i=j$, we necessarily have $\mym_{i}(\telt)-\mym_{i}(\selt)=2$, hence $M_{ii}=2$. 
To confirm that $M$ is a pulsation matrix, it only remains to show that, when $i \not= j$, $M_{ij} \leq 0$ and that $M_{ij}=0$ if and only if $M_{ji}=0$. 

With $i \not= j$, suppose $\mathcal{C}$ is some $\{i,j\}$-component with at least one edge of color $i$ and at least one of color $j$. 
Since $\mathcal{C}$ is itself a modular lattice (\JCompResult), there is a unique minimal element (which we denote by $\xelt$) and a unique maximal element (denoted $\yelt$). 
Consider a path $\mathcal{P}$ from $\xelt$ up to $\yelt$ whose length is the same as the length of $\mathcal{C}$. 
By \ColorsLemma, the numbers $\mya_{i}:=\mya_{i}(\mathcal{P})$ and $\mya_{j}:=\mya_{j}(\mathcal{P})$ do not depend on the choice of path from $\xelt$ up to $\yelt$. 
Then $wt(\xelt) + \mya_{i}\alpha_{i} + \mya_{j}\alpha_{j} = wt(\yelt)$. In particular, $\mym_{i}(\xelt)+ \mya_{i} \cdot 2 + \mya_{j} \cdot M_{ji}  =  \mym_{i}(\yelt)$ and $\mym_{j}(\xelt)+ \mya_{i} \cdot M_{ij} + \mya_{j} \cdot 2 = \mym_{j}(\yelt)$. 
It follows that $\mya_{j} \cdot M_{ji} = (\mym_{i}(\yelt)- \mya_{i}) + (-\mym_{i}(\xelt) - \mya_{i})$ and $\mya_{i} \cdot M_{ij} = (\mym_{j}(\yelt) - \mya_{j}) + (-\mym_{j}(\xelt)- \mya_{j})$.  
Since $\yelt$ is maximal, then $\mym_{i}(\yelt) = l_{i}(\yelt) \leq \mya_{i}$, and similarly $\mym_{j}(\yelt) = l_{j}(\yelt) \leq \mya_{j}$, $-\mym_{i}(\xelt) = l_{i}(\xelt) \leq \mya_{i}$, and $-\mym_{j}(\xelt) = l_{j}(\xelt) \leq \mya_{j}$. 
That is, each of $\mym_{i}(\yelt) - \mya_{i}$, $-\mym_{i}(\xelt) - \mya_{i}$, $\mym_{j}(\yelt) - \mya_{j}$, and $-\mym_{j}(\xelt) - \mya_{j}$ is nonpositive. 
So, each of $\mya_{j} \cdot M_{ji}$ and $\mya_{i} \cdot M_{ij}$ is nonpositive,  and, since $\mya_{i} > 0$ and $\mya_{j} > 0$, then $M_{ji} \leq 0$ and $M_{ij} \leq 0$.

Continuing, suppose that $M_{ij} = 0$. 
Then $(\mym_{j}(\yelt) - \mya_{j}) + (-\mym_{j}(\xelt)- \mya_{j}) = 0$, and since each parenthetical quantity is nonpositive, we conclude that $\mym_{j}(\yelt) = \mya_{j} = l_{j}(\yelt)$ and that $-\mym(\xelt) = \mya_{j} = l_{j}(\xelt)$. 
Therefore, there is a path $\xelt \myarrow{j} \cdots \myarrow{j} \zelt$ from $\xelt$ up to $\zelt$ with exactly $\mya_{j}$ edges of color $j$. 
Now consider a path from $\zelt$ up to $\yelt$. 
No such path can use any more color $j$ edges, so we must have $\zelt \myarrow{i} \cdots \myarrow{i} \yelt$ with exactly $\mya_{i}$ edges of color $i$. 
Note that $\zelt$ must be minimal in the $i$-component $\comp_{i}(\yelt)$, so $\mym_{i}(\yelt) = \mya_{i} = l_{i}(\yelt)$.
Similarly, we can descend from $\yelt$ to some $\zelt'$ so that $\zelt' \myarrow{j} \cdots \myarrow{j} \yelt$ with $\mya_j$ color $j$ edges and then descend from $\zelt'$ to $\xelt$ via $\mya_i$ edges of color $i$. 
Then $\zelt'$ is maximal in $\comp_{i}(\xelt)$, hence $-\mym_{i}(\xelt) = \mya_{i} = l_{i}(\xelt)$. 
In particular, $(\mym_{i}(\yelt)- \mya_{i}) + (-\mym_{i}(\xelt) - \mya_{i}) = 0 = \mya_{j} \cdot M_{ji}$. 
Since we have at least one edge of color $j$, then $\mya_{j}>0$, which means we must have $M_{ji} = 0$. 
That is, from the hypothesis $M_{ij} = 0$, we can deduce $M_{ji}=0$. 
We can switch the roles of $i$ and $j$ from the beginning of this paragraph to see that $M_{ji} = 0$ implies $M_{ij}=0$. 

Now suppose that no $\{i,j\}$-component of $L$ has edges of both colors $i$ and $j$. 
Then, for an edge $\selt \myarrow{i} \telt$, we must have $\mym_{j}(\selt) = 0 = \mym_{j}(\telt)$, and hence $M_{ij} = 0$.
Similarly, if we take an edge $\relt \myarrow{j} \uelt$, we see that $\mym_{i}(\relt) = 0 = \mym_{i}(\uelt)$, and hence $M_{ji} = 0$. 

We conclude that $M$ is a pulsation matrix of integers. 
By \GCMTheorem.1, it follows that $\mathscr{G}$ is a Coxeter--Dynkin posy, which completes the proof of {\sl (1)}. 
For {\sl (2)}, apply \GCMTheorem.2 to see that $\ecolor(\EdgeSet(L)) = I$. 
The remaining claims of {\sl (2)} follow from part {\sl (1)}. 
For {\sl (3)}, the fact that $\mathscr{H}$ is an integral Coxeter--Dynkin flower means we can consider any splitting poset for a $W(\mathscr{H})$-symmetric function or any supporting graph for a representation of the related simple Lie algebra $\mathfrak{g}(\mathscr{H})$. 
In Theorem 1.2 of \cite{DonAdjoint}, we showed that there are exactly $n$ diamond-colored modular lattice supporting graphs for the adjoint representation of $\mathfrak{g}(\mathscr{H})$. 
By definition, a supporting graph is necessarily $\mathscr{H}$-structured, thereby establishing the existence claim of {\sl (3)}. 
Alternatively, {\sl (3)} follows from \ShortAdjointDCML\ below.\hfill\QED

{\bf [\S \GStructureIntroSection.6:\! Some open problems.]} 
Let $\mathscr{G}$ represent a generic Coxeter--Dynkin flower and $\lambda$ a dominant weight. 
A diamond-colored modular lattice $L$ is $(\mathscr{G},\lambda)$-structured if $L$ is $\mathscr{G}$-structured and its maximal element has weight $\lambda$. 
An interesting question is for which dominant weights $\lambda$ does there exist a $(\mathscr{G},\lambda)$-structured DCML? 
Our guess is that the answer is `all dominant weights'. 
If we rephrase this question by replacing `DCML' with `DCDL', then we already know the answer cannot be all dominant weights, as there exists no $\big(\myD_{4},(0,1,0,0)\big)$-structured DCDL, although $\big(\myD_{4},\lambda\big)$-structured DCDL's exist when $\lambda \in \{(1,0,0,0),(0,0,1,0),(0,0,0,1)\}$. 
At this time, we know of no $\myE_{8}$-structured DCDL's with more than one element, and our guess is that there is no $\myE_{8}$-structured DCDL whose maximal element has a nonzero dominant weight. 
Based on some evidence presented in \cite{Gilliland}, we believe there exists no $\big(\myF_{4},\omega_{2}\big)$-structured DCDL and no $\big(\myF_{4},\omega_{3}\big)$-structured DCDL. 

One of the main contributions of \cite{Gilliland} was the discovery of a $\big(\myF_{4},2\omega_{1}\big)$-structured DCDL and a $\big(\myF_{4},2\omega_{4}\big)$-structured DCDL. 
Computer verification was used to show that these two DCDL's are actually `splitting posets' for the associated `Weyl bialternants' $\chi_{_{2\omega_{1}}}^{\mytinyF_{4}}$ and $\chi_{_{2\omega_{4}}}^{\mytinyF_{4}}$ (for explication of these concepts, see \S \FlowerSection.6-\FlowerSection.8).  
The discovery of these DCDL's, and the demonstration that no $\big(\myF_{4},\omega_{2}\big)$- or $\big(\myF_{4},\omega_{3}\big)$-structured DCDL's exist, was aided by the development, in Chapter 3 of \cite{Gilliland}], of a general algorithm whose goal was to produce, when possible, a `smallest' $(\mathscr{G},\lambda)$-structured DCDL and to return the empty set if no $(\mathscr{G},\lambda)$-structured DCDL exists. 
Some experimental and some theoretical evidence suggested that whenever there exists a $(\mathscr{G},\lambda)$-structured DCDL, then the algorithm would yield a $(\mathscr{G},\lambda)$-structured DCDL that appears as a full-length sublattice within any other $(\mathscr{G},\lambda)$-structured DCDL. 
For this reason, the algorithm was called the {\bf distributive core} algorithm. 
In Chapter 6 of \cite{Gilliland}, the following conjecture was stated as Conjecture 6.1.  

\noindent 
{\bf \DistCoreConjecture\ (Gilliland and Donnelly)}\ \ {\sl If there exists a $(\mathscr{G},\lambda)$-structured DCDL $L$, then $L$ contains a full length edge-colored sublattice that is isomorphic to the output $K(\mathscr{G},\lambda)$ of the {\bf distributive core} algorithm from \cite{Gilliland}. 
In particular, if the {\bf distributive core} algorithm returns empty output for some input $\mathscr{G}$ and $\lambda$, then there does not exist a $(\mathscr{G},\lambda)$-structured DCDL.}

\begin{center}
\underline{\hspace*{4in}}
\end{center}

\begin{figure}
\begin{center}
{\bf \ESixLatticeFigure.2}\ \  Compression posets for the $\myE_{6}$-structured DCDL's of \ESixLatticeFigure.1.\\
{\small (These are the two ``$\mysmallE_{6}$-minuscule compression posets'', cf.\ \S \MinusculeExampleSection.)}

\vspace*{0.1in} 
\ESixFlower

\vspace*{0.1in} 
\setlength{\unitlength}{1cm}
\begin{picture}(5,12)
\put(-1.9,7.2){\LARGE $P_{\mysmallE_{6}}(\omega_{\mbox{\color{Cyan}{\small $1$}}})$}
\put(0.1,7.6){\vector(2,1){1}}
\put(0,10){\TypeEboxDot{Cyan}}
\put(0.7,10.3){\color{Cyan}{\em \footnotesize 1}}
\put(0,4){\TypeEboxDot{Cyan}}
\put(0.1,4){\color{Cyan}{\em \footnotesize 1}}
\put(1,9){\TypeEboxDot{Red}}
\put(1.7,9.3){\color{Red}{\em \footnotesize 3}}
\put(1,5){\TypeEboxDot{Red}}
\put(1.1,5.3){\color{Red}{\em \footnotesize 3}}
\put(1,3){\TypeEboxDot{Red}}
\put(1.1,3){\color{Red}{\em \footnotesize 3}}
\put(2,8){\TypeEboxDot{Purple}}
\put(2.7,8.3){\color{Purple}{\em \footnotesize 4}}
\put(2,6){\TypeEboxDot{Purple}}
\put(2.1,6.15){\color{Purple}{\em \footnotesize 4}}
\put(2,4){\TypeEboxDot{Purple}}
\put(2.7,4.15){\color{Purple}{\em \footnotesize 4}}
\put(2,2){\TypeEboxDot{Purple}}
\put(2.1,2){\color{Purple}{\em \footnotesize 4}}
\put(1,7){\TypeEboxDot{SpringGreen}}
\put(1.1,7){\color{SpringGreen}{\em \footnotesize 2}}
\put(3,3){\TypeEboxDot{SpringGreen}}
\put(3.7,3.3){\color{SpringGreen}{\em \footnotesize 2}}
\put(3,7){\TypeEboxDot{OliveGreen}}
\put(3.7,7.3){\color{OliveGreen}{\em \footnotesize 5}}
\put(3,5){\TypeEboxDot{OliveGreen}}
\put(3.7,5){\color{OliveGreen}{\em \footnotesize 5}}
\put(3,1){\TypeEboxDot{OliveGreen}}
\put(3.1,1){\color{OliveGreen}{\em \footnotesize 5}}
\put(4,6){\TypeEboxDot{BurntOrange}}
\put(4.7,6.3){\color{BurntOrange}{\em \footnotesize 6}}
\put(4,0){\TypeEboxDot{BurntOrange}}
\put(4.1,0){\color{BurntOrange}{\em \footnotesize 6}}
\thicklines
\put(2.625,2.375){\color{Black}\qbezier(0,0)(0.375,0.375)(0.75,0.75)}
\put(1.625,3.375){\color{Black}\qbezier(0,0)(0.375,0.375)(0.75,0.75)}
\put(2.625,4.375){\color{Black}\qbezier(0,0)(0.375,0.375)(0.75,0.75)}
\put(3.625,5.375){\color{Black}\qbezier(0,0)(0.375,0.375)(0.75,0.75)}
\put(0.625,4.375){\color{Black}\qbezier(0,0)(0.375,0.375)(0.75,0.75)}
\put(1.625,5.375){\color{Black}\qbezier(0,0)(0.375,0.375)(0.75,0.75)}
\put(2.625,6.375){\color{Black}\qbezier(0,0)(0.375,0.375)(0.75,0.75)}
\put(1.625,7.375){\color{Black}\qbezier(0,0)(0.375,0.375)(0.75,0.75)}
\put(4.39,0.375){\color{Black}\qbezier(0,0)(-0.375,0.375)(-0.75,0.75)}
\put(3.39,1.375){\color{Black}\qbezier(0,0)(-0.375,0.375)(-0.75,0.75)}
\put(2.39,2.375){\color{Black}\qbezier(0,0)(-0.375,0.375)(-0.75,0.75)}
\put(1.39,3.375){\color{Black}\qbezier(0,0)(-0.375,0.375)(-0.75,0.75)}
\put(2.39,4.375){\color{Black}\qbezier(0,0)(-0.375,0.375)(-0.75,0.75)}
\put(3.39,3.375){\color{Black}\qbezier(0,0)(-0.375,0.375)(-0.75,0.75)}
\put(3.39,5.375){\color{Black}\qbezier(0,0)(-0.375,0.375)(-0.75,0.75)}
\put(2.39,6.375){\color{Black}\qbezier(0,0)(-0.375,0.375)(-0.75,0.75)}
\put(4.39,6.375){\color{Black}\qbezier(0,0)(-0.375,0.375)(-0.75,0.75)}
\put(3.39,7.375){\color{Black}\qbezier(0,0)(-0.375,0.375)(-0.75,0.75)}
\put(2.39,8.375){\color{Black}\qbezier(0,0)(-0.375,0.375)(-0.75,0.75)}
\put(1.39,9.375){\color{Black}\qbezier(0,0)(-0.375,0.375)(-0.75,0.75)}
\end{picture}
\hspace*{1.5cm}
\begin{picture}(5,12)
\put(4.5,7){\LARGE $P_{\mysmallE_{6}}(\omega_{\mbox{\color{BurntOrange}{\small $6$}}})$}
\put(4.5,6.7){\vector(-2,-1){1}}
\put(5,6){\LARGE $\cong P_{\mysmallE_{6}}(\omega_{\mbox{\color{Cyan}{\small $1$}}})^{*}$}
\put(0,0){\TypeEboxDot{Cyan}}
\put(0.7,0){\color{Cyan}{\em \footnotesize 1}}
\put(0,6){\TypeEboxDot{Cyan}}
\put(0.1,6.3){\color{Cyan}{\em \footnotesize 1}}
\put(1,1){\TypeEboxDot{Red}}
\put(1.7,1){\color{Red}{\em \footnotesize 3}}
\put(1,5){\TypeEboxDot{Red}}
\put(1.1,5){\color{Red}{\em \footnotesize 3}}
\put(1,7){\TypeEboxDot{Red}}
\put(1.1,7.3){\color{Red}{\em \footnotesize 3}}
\put(2,2){\TypeEboxDot{Purple}}
\put(2.7,2){\color{Purple}{\em \footnotesize 4}}
\put(2,4){\TypeEboxDot{Purple}}
\put(2.1,4.15){\color{Purple}{\em \footnotesize 4}}
\put(2,6){\TypeEboxDot{Purple}}
\put(2.7,6.15){\color{Purple}{\em \footnotesize 4}}
\put(2,8){\TypeEboxDot{Purple}}
\put(2.1,8.3){\color{Purple}{\em \footnotesize 4}}
\put(1,3){\TypeEboxDot{SpringGreen}}
\put(1.1,3.3){\color{SpringGreen}{\em \footnotesize 2}}
\put(3,7){\TypeEboxDot{SpringGreen}}
\put(3.7,7){\color{SpringGreen}{\em \footnotesize 2}}
\put(3,3){\TypeEboxDot{OliveGreen}}
\put(3.7,3){\color{OliveGreen}{\em \footnotesize 5}}
\put(3,5){\TypeEboxDot{OliveGreen}}
\put(3.7,5.3){\color{OliveGreen}{\em \footnotesize 5}}
\put(3,9){\TypeEboxDot{OliveGreen}}
\put(3.1,9.3){\color{OliveGreen}{\em \footnotesize 5}}
\put(4,4){\TypeEboxDot{BurntOrange}}
\put(4.7,4){\color{BurntOrange}{\em \footnotesize 6}}
\put(4,10){\TypeEboxDot{BurntOrange}}
\put(4.1,10.3){\color{BurntOrange}{\em \footnotesize 6}}
\thicklines
\put(2.625,8.125){\color{Black}\qbezier(0,0)(0.375,-0.375)(0.75,-0.75)}
\put(1.625,7.125){\color{Black}\qbezier(0,0)(0.375,-0.375)(0.75,-0.75)}
\put(2.625,6.125){\color{Black}\qbezier(0,0)(0.375,-0.375)(0.75,-0.75)}
\put(3.625,5.125){\color{Black}\qbezier(0,0)(0.375,-0.375)(0.75,-0.75)}
\put(0.625,6.125){\color{Black}\qbezier(0,0)(0.375,-0.375)(0.75,-0.75)}
\put(1.625,5.125){\color{Black}\qbezier(0,0)(0.375,-0.375)(0.75,-0.75)}
\put(2.625,4.125){\color{Black}\qbezier(0,0)(0.375,-0.375)(0.75,-0.75)}
\put(1.625,3.125){\color{Black}\qbezier(0,0)(0.375,-0.375)(0.75,-0.75)}
\put(4.375,10.125){\color{Black}\qbezier(0,0)(-0.375,-0.375)(-0.75,-0.75)}
\put(3.375,9.125){\color{Black}\qbezier(0,0)(-0.375,-0.375)(-0.75,-0.75)}
\put(2.375,8.125){\color{Black}\qbezier(0,0)(-0.375,-0.375)(-0.75,-0.75)}
\put(1.375,7.125){\color{Black}\qbezier(0,0)(-0.375,-0.375)(-0.75,-0.75)}
\put(2.375,6.125){\color{Black}\qbezier(0,0)(-0.375,-0.375)(-0.75,-0.75)}
\put(3.375,7.125){\color{Black}\qbezier(0,0)(-0.375,-0.375)(-0.75,-0.75)}
\put(3.375,5.125){\color{Black}\qbezier(0,0)(-0.375,-0.375)(-0.75,-0.75)}
\put(2.375,4.125){\color{Black}\qbezier(0,0)(-0.375,-0.375)(-0.75,-0.75)}
\put(4.375,4.125){\color{Black}\qbezier(0,0)(-0.375,-0.375)(-0.75,-0.75)}
\put(3.375,3.125){\color{Black}\qbezier(0,0)(-0.375,-0.375)(-0.75,-0.75)}
\put(2.375,2.125){\color{Black}\qbezier(0,0)(-0.375,-0.375)(-0.75,-0.75)}
\put(1.375,1.125){\color{Black}\qbezier(0,0)(-0.375,-0.375)(-0.75,-0.75)}
\end{picture}
\end{center}
\end{figure}

\vspace*{0.5cm} 
\noindent 
{\bf \S \FrameSection. Scaffolds, skew-stacks, and DCDL's.} 
In our programmatic study of poset models for semisimple Lie algebra representations and their companion Weyl symmetric functions, we have encountered many posets that are diamond-colored distributive lattices and that arise naturally as partial orderings of certain integer arrays. 
Our purpose here is to develop this latter idea more systematically with the notion of `scaffolds'. 

{\bf [\S \FrameSection.1:\! Scaffolds.]} 
Let $P$ be a finite poset, to be identified with its covering digraph. 
Assign to each vertex $x$ in $P$ an `amplifier' $a_{x}$ and two `bounds' $(b_{x}^{(\mbox{\tiny lo})},b_{x}^{(\mbox{\tiny hi})})$, and require that the amplifier $a_x$ is a positive integer and that the bounds are integers satisfying $0 \leq b_{x}^{(\mbox{\tiny lo})} \leq b_{x}^{(\mbox{\tiny hi})}$.  
Let $\mathcal{A} := (a_{x})_{x \in P}$ denote the collection of amplifiers and $\mathcal{B} := \left((b_{x}^{(\mbox{\tiny lo})},b_{x}^{(\mbox{\tiny hi})})\right)_{x \in P}$ the collection of bounding pairs. 
Let $I$ be a color palette and $\vcolor\! :\! P \rightarrow I$ a vertex-coloring function. 
Then $\scaffoldS = (P,\mathcal{A},\mathcal{B},\vcolor)$ is called a {\em scaffold}, and the poset $P = P(\scaffoldS)$ is the scaffold's {\em footprint}. 

An array $\telt = (c_{x}(\telt))_{x \in P}$ is an {\em ideal array} over the scaffold $\scaffoldS$ if for all $x \in P$ we have 
\[b_{x}^{(\mbox{\tiny lo})} \leq c_{x}(\telt) \leq b_{x}^{(\mbox{\tiny hi})}\]
and for all $x \rightarrow y$ in $P$ we have 
\[a_{y}c_{x}(\telt) \geq a_{x}c_{y}(\telt).\]  
That is, $\telt$ satisfies all `bounding inequalities' and all `amplifier inequalities'. 
Let $\mathbf{L}(\scaffoldS)$ be the set of all ideal arrays over $\scaffoldS$.   
We intend to realize $\mathbf{L}(\scaffoldS)$ as a diamond-colored distributive lattice, but we will need the following concept:  
Say $\scaffoldS$ is a {\em proper} scaffold if $\mathbf{L}(\scaffoldS)$ is nonempty. 
For ideal arrays $\selt$ and $\telt$ in such an $\mathbf{L}(\scaffoldS)$, we say $\selt \leq \telt$ if and only if $c_{x}(\selt) \leq c_{x}(\telt)$ for all $x \in P$, i.e.\ if and only if $\selt$ is component-wise bounded by $\telt$. 
So, $\mathbf{L}(\scaffoldS)$ is the poset of ideal arrays over $\scaffoldS$ in the component-wise partial order. 

\noindent 
{\bf \DistributiveProp}\ \ {\sl Let $\scaffoldS$ be a proper scaffold. 
Then $L := \mathbf{L}(\scaffoldS)$ is a distributive lattice. 
Moreover, $\selt \rightarrow \telt$ in the covering digraph for $L$ if and only if there is some $y \in P$ such that $c_{y}(\selt)+1 = c_{y}(\telt)$ while $c_{x}(\selt) = c_{x}(\telt)$ for all $x \not= y$ in $P$.
} 

{\em Proof.} Let $\selt, \telt \in \mathbf{L}(\scaffoldS)$.  
Let $\relt := (\min\{c_{x}(\selt),c_{x}(\telt)\})_{x \in P}$ and $\uelt := (\max\{c_{x}(\selt),c_{x}(\telt)\})_{x \in P}$.  
It is then easy to check that $\relt$ and $\uelt$ are ideal arrays over $\scaffoldS$ and that $\relt = \selt \wedge \telt$ and $\uelt = \selt \vee \telt$.  
It is routine to see that $\vee$ distributes over $\wedge$, and that $\wedge$ distributes over $\vee$. 
So $\mathbf{L}(\scaffoldS)$ is distributive. 
If for some $y \in P$ we have $c_{y}(\selt)+1 = c_{y}(\telt)$ while $c_{x}(\selt) = c_{x}(\telt)$ for all $x \not= y$ in $P$, then it readily follows from the definition of the component-wise partial order on $\mathbf{L}(\scaffoldS)$ that $\selt \rightarrow \telt$. 

Now say $\selt \rightarrow \telt$ in $\mathbf{L}(\scaffoldS)$. 
Then $c_{x}(\selt) \leq c_{x}(\telt)$ for all $x \in P$, and this inequality is at least sometimes strict. 
Pick $y \in P$ so that $c_{y}(\selt) < c_{y}(\telt)$ while $c_{x}(\selt) = c_{x}(\telt)$ when $x < y$. 
Define a new array $\selt'$ wherein $c_{x}(\selt') = c_{x}(\selt)$ unless $x=y$, in which case $c_{y}(\selt') = c_{y}(\selt)+1$. 
It is easy to see that $\selt'$ satisfies all bounding inequalities. 
Then to see that $\selt'$ is an ideal array, we only need to check the amplifier inequalities. 
So suppose $u \rightarrow v$ in $P$. 
We must show that $a_{v}c_{u}(\selt') \geq a_{u}c_{v}(\selt')$. 
Of course this follows immediately from the associated amplifier inequality for $\selt$ unless $u = y$ or $v = y$. 
So suppose that $x \rightarrow y$ for some $x \in P$. Then $a_{y}c_{x}(\selt) = a_{y}c_{x}(\telt) \geq a_{x}c_{y}(\telt) \geq a_{x}(c_{y}(\selt)+1)$, so that $a_{y}c_{x}(\selt) \geq a_{x}c_{y}(\selt')$. 
Now suppose that $y \rightarrow z$ for some $z \in P$. 
Then $a_{z}c_{y}(\selt') = a_{z}(c_{y}(\selt)+1) > a_{y}c_{z}(\selt) = a_{y}c_{z}(\selt')$. 
Therefore $\selt'$ is an ideal array. 
Moreover, $\selt'$ was constructed so that $\selt < \selt' \leq \telt$. 
Since $\selt \rightarrow \telt$ in $\mathbf{L}(\scaffoldS)$, then we must have $\telt = \selt'$.
\hfill\QED

We note that any poset $P$ can be regarded as the footprint of some proper scaffold. 
For example, we might require for each $x \in P$ that the amplifier $a_{x}$ is $1$, that the lower bound $b_{x}^{(\mbox{\tiny lo})}$ is $0$, and that the upper bound $b_{x}^{(\mbox{\tiny hi})}$ is $k$, where $k$ is some fixed nonnegative integer not depending on $x$. 
For the associated scaffold $\scaffoldS$, it is easy to see that $\mathbf{L}(\scaffoldS)$ is just $\mathbf{J}(P \times [1,k]_{\mathbb{Z}})$, where the latter is the distributive lattice of order ideals from the poset product of $P$ with the chain $[1,k]_{\mathbb{Z}}$. 
This partly explains our use of the adjective ``ideal'' in reference to ``ideal arrays.''

If $\scaffoldS$ is a proper scaffold whose footprint $P = P(\scaffoldS)$ is $I$-vertex-colored, then, in view of \DistributiveProp, we can color the edges of $\mathbf{L}(\scaffoldS)$ by the rule $\selt \myarrow{i} \telt$ if $\selt \rightarrow \telt$ in $\mathbf{L}(\scaffoldS)$ and $\vcolor(y) = i$ where $y$ is the vertex from $P$ such that $c_{y}(\selt) + 1 = c_{y}(\telt)$. 
In this way, the vertex-coloring function $\vcolor\! :\! P \rightarrow I$ induces an edge-coloring function $\ecolor\! :\! \EdgeSet(\mathbf{L}(\scaffoldS)) \rightarrow I$. 
We use the notation $\mathbf{L}_{\mbox{\tiny color}}(\scaffoldS)$ to denote $\mathbf{L}(\scaffoldS)$ together with this assignment of edge colors. 
In view of the next result, we call $\mathbf{L}_{\mbox{\tiny color}}(\scaffoldS)$ the {\em diamond-colored distributive lattice of (ideal) arrays} over $\scaffoldS$. 

\noindent 
{\bf \ChainProp}\ \ {\sl Let $\scaffoldS$ be a proper scaffold whose footprint is $I$-vertex-colored, as above. 
Then} $\mathbf{L}_{\mbox{\tiny color}}(\scaffoldS)$ {\sl is a diamond-colored distributive lattice. 
Further, suppose that whenever $x \rightarrow y$ in $P$, we have $\vcolor(x) \not= \vcolor(y)$. 
Then for each $i \in I$, each $i$-component of} $\mathbf{L}_{\mbox{\tiny color}}(\scaffoldS)$ {\sl is isomorphic to a product of chains.} 

{\em Proof.} We must verify that on a diamond of edges$\,$ \parbox{1.4cm}{\begin{center}
\setlength{\unitlength}{0.2cm}
\begin{picture}(6.5,3.5)
\put(3,0){\circle*{0.5}} \put(1,2){\circle*{0.5}}
\put(3,4){\circle*{0.5}} \put(5,2){\circle*{0.5}}
\put(1,2){\line(1,1){2}} \put(3,0){\line(-1,1){2}}
\put(5,2){\line(-1,1){2}} \put(3,0){\line(1,1){2}}
\put(1.75,0.55){\em \small k} \put(3.7,0.7){\em \small l}
\put(1.7,2.7){\em \small i} \put(3.75,2.55){\em \small j}
\put(3.5,-0.75){\footnotesize $\relt$} \put(5.75,1.75){\footnotesize $\telt$}
\put(3.5,4){\footnotesize $\uelt$} \put(-0.5,1.75){\footnotesize $\selt$}
\end{picture} \end{center}}, we have $i=l$ and $k=j$. 
Now, to obtain $\uelt$ from $\selt$ requires that for some $y \in P$ with $\vcolor(y) = i$ we have $c_{x}(\uelt) = c_{x}(\selt)$ for all $x \not= y$ in $P$ and $c_{y}(\uelt)=c_{y}(\selt)+1$.  
To obtain the element $\telt$  (distinct from $\selt$) from $\uelt$ requires that for some $z \in P$ with $z \not= y$ and $\vcolor(\zelt) = j$ we have $c_{x}(\telt) = c_{x}(\uelt)$ for all $x \not= z$ in $P$ and $c_{z}(\telt)=c_{z}(\uelt)-1$. 
On the other hand, to obtain $\telt$ from $\selt$ via $\relt$ requires that we reduce an array entry by one and increase another array entry by one. 
Evidently, the decreased array entry is $c_{z}$, so $k = \vcolor(z) = j$, and the increased array entry is $c_{y}$, so $l = \vcolor(y) = i$. 
So $\mathbf{L}_{\mbox{\tiny color}}(\scaffoldS)$ is diamond-colored. 

For any $\telt \in \mathbf{L}_{\mbox{\tiny color}}(\scaffoldS)$ and $i \in I$, it is clear that $\comp_{i}(\telt)$ is comprised of precisely those ideal arrays $\selt$ such that $c_{x}(\selt) = c_{x}(\telt)$ for all $x \not\in \vcolor^{-1}(i)$. 
Now say $x,y \in \vcolor^{-1}(i)$. 
Since, by hypothesis, $x$ and $y$ are not adjacent in (the covering digraph for) $P$, then changes to the $c_{x}(\telt)$ entry of the array $\telt$ will not produce any violations of any bounding or amplifier inequality involving $c_{y}(\telt)$, and vice-versa. 
That is, changes to $\telt$ in such positions $x$ and $y$ are independent. 
So, $\displaystyle \comp_{i}(\telt) \cong \prod_{x \in \vcolor^{-1}(i)} \mathcal{C}_{x}$, where $\mathcal{C}_{x}$ is a chain whose length is the difference of the maximum value that can replace $c_{x}(\telt)$, without violating any amplifier or bounding inequalities, and the minimum such value for $c_{x}(\telt)$.\hfill\QED

{\bf [\S \FrameSection.2:\! Skew-stacks.]} 
Next, we describe a method of stacking posets that will give us a different perspective on certain types of scaffolds. 
To motivate the idea, we offer an example. 
Regard $[1,k]_{\mathbb{Z}}$ to be an uncolored chain with $m$ elements. 
For a given vertex-colored poset $P$, we can naturally view the product $P \times [1,k]_{\mathbb{Z}}$ as a vertex-colored poset by assigning to element $(v,j) \in P \times [1,k]_{\mathbb{Z}}$ the same color as $v \in P$.  
When $P$ is a `minuscule compression poset' (to be discussed in \S \MinusculeLatticePosetSection, below), then the diamond-colored distributive lattice $\Jcolor(P \times [1,k]_{\mathbb{Z}})$ is known to have many felicitous algebraic-combinatorial properties, see for example Section 9 of \cite{DonPosetModels}. 
The method of stacking posets that we offer here generalizes the notion of the product of a poset by a chain. 
See \ESixStackFigure\ for visual examples of some of the concepts presented below. 

\noindent 
{\bf \StackDefsOne}\ \   
As usual, $I$ is our palette of colors. Suppose we are given posets $P_{1},\ldots,P_{m}$ so that for each $r \in [1,m]_{\mathbb{Z}}$, $P_{r}$ has vertex-coloring function $\vcolor^{(r)}\! :\! P_{r} \rightarrow I$. 
For $r \in [1,m]_{\mathbb{Z}}$, suppose $\mathcal{F}_{r}$ is an up-set and $\mathcal{I}_{r}$ is a down-set from $P_{r}$, with $\mathcal{F}_{1} := \emptyset =: \mathcal{I}_{m}$. 
Further, suppose $\phi_{r}: \mathcal{I}_{r} \longrightarrow \mathcal{F}_{r+1}$ is a vertex-color preserving poset isomorphism when $1 \leq r \leq m-1$. 
Now take the covering digraph of the disjoint sum $P_{1} \oplus P_{2} \oplus \cdots \oplus P_{m}$ and add directed edges of the form $u \rightarrow v$ whenever $v = \phi_{r}(u)$ for some $r \in [1,m-1]_{\mathbb{Z}}$. 
Momentarily denote by $\mystackedP$ the resulting vertex-colored directed graph, wherein the set $\VertexSet(\mystackedP)$ of vertices is identified with the set $P_{1} \oplus \cdots \oplus P_{m}$ and colored accordingly, and the directed edge set $\EdgeSet(\mystackedP)$ is $\EdgeSet(P_{1} \oplus \cdots \oplus P_{m}) \disjointunion \left\{u \rightarrow v\, \rule[-2.75mm]{0.2mm}{7.5mm}\, v = \phi_{r}(u) \mbox{ for some } r \in [1,m-1]_{\mathbb{Z}}\right\}$. 
One can see that the digraph $\mystackedP$ is acyclic and is its own transitive reduction.  
So, $\mystackedP$ is the covering digraph of the poset whose vertices are $\VertexSet(\mystackedP)$ and whose partial order is the reachability relation. 
We call the vertex-colored poset $\mystackedP$ a {\em skew-stack} or {\em skew-stacked poset}, each $P_{r}$ is a {\em tier}, each $\phi_{r}$ is a {\em bond}, and these constituent parts are encapsulated in the notation $P_{1} \myposetstack_{\phi_{1}} P_{2} \myposetstack_{\phi_{2}} \cdots \myposetstack_{\phi_{m-1}} P_{m} =: \mystackedP$. 
Edges in $\mystackedP$ from $\EdgeSet(P_{1} \oplus \cdots \oplus P_{m})$ are {\em planks} and from $\left\{u \rightarrow v\, \rule[-2.75mm]{0.2mm}{7.5mm}\, v = \phi_{r}(u) \mbox{ for some } r \in [1,m-1]_{\mathbb{Z}}\right\}$ are {\em pillars}. 
Given the natural congruences $(P_{1} \myposetstack_{\phi_{1}} P_{2}) \myposetstack_{\phi_{2}} P_{3} \cong P_{1} \myposetstack_{\phi_{1}} P_{2} \myposetstack_{\phi_{2}} P_{3} \cong P_{1} \myposetstack_{\phi_{1}} (P_{2} \myposetstack_{\phi_{2}} P_{3})$ of vertex-colored posets, we regard ``$\myposetstack_{\phi}$'' to be an associative binary operation. 
For any $v \in \mystackedP$, let $\mytier(v) := t$ if $v$ is from the tier $P_{t}$. 
As a special case, observe that if, for some vertex-colored poset $Q$, we have $P_{r} \cong Q \cong \mathcal{I}_{r} \cong \mathcal{F}_{r+1} \cong P_{m}$ for isomorphisms $\phi_{r}$ when $r \in [1,m-1]_{\mathbb{Z}}$, then $P_{1} \myposetstack_{\phi_{1}} P_{2} \myposetstack_{\phi_{2}} \cdots \myposetstack_{\phi_{m-1}} P_{m} \cong Q \times [1,m]_{\mathbb{Z}}$. 
In this case, we give the skew-stack the notation $Q^{(\mysmallposetstack m)}$ and (sometimes) call it {\em prism}.  

\noindent
{\bf \StackEarlyLemma}\ \ 
{\sl In the above notation, let $u$ and $v$ be members of} $\mystackedP$, {\sl and assume that} $\mystackedP$ {\sl is connected.}\\ 
{\sl (1)  Suppose} $\mytier(u) = s$ {\sl and} $\mytier(v) = t$.  
{\sl Then $u < v$ in} $\mystackedP$ {\sl if and only if exactly one of the following three mutually exclusive conditions hold: ({\sl i}) $s=t$ and $u < v$ within the tier $P_{s}$; ({\sl ii}) $s < t$ and $v = (\phi_{t-1} \circ \cdots \circ \phi_{s})(u)$; or ({\sl iii}) $s < t$ and there exist $u' \in P_{t} \cap \overline{u}$ and $v' \in P_{s} \cap \overline{v}$ such that $u' < v$ and $u < v'$ within their respective tiers.}\\ 
{\sl (2) Suppose $u \leq v$ within some tier $P_{s}$. 
If, for some $t > s$ and $v' \in P_{t}$, we have $v' = (\phi_{t-1} \circ \cdots \phi_{s})(v)$, then there exists $u' \in P_{t}$ such that $u' = (\phi_{t-1} \circ \cdots \phi_{s})(u)$. 
Similarly, if, for some $r < s$ and $u'' \in P_{r}$, we have $u = (\phi_{s-1} \circ \cdots \phi_{r})(u'')$, then there exists $v'' \in P_{r}$ such that $v = (\phi_{s-1} \circ \cdots \phi_{r})(v'')$.}

{\em Proof.} For {\sl (1)}, the converse is evident, so assume $\mytier(u) = s$ and $\mytier(v) = t$ with $u < v$ in $\mystackedP$, and suppose neither {\sl (i)} nor {\sl (ii)} applies. 
In $\mystackedP$, since $u$ can only be smaller than elements from tiers $P_{r}$ where $s \leq r$, then it must be the case that $s \leq t$. 
If $s=t$, then the fact that $u < v$ in $\mystackedP$ means that $u < v$ in tier $P_{s}$, which coincides with {\sl (i)}; therefore $s < t$. 
Now say $u = u_{0} \rightarrow u_{1} \rightarrow \cdots \rightarrow u_{p} = v$ is a path from $u$ up to $v$ in $\mystackedP$, so each edge in this path is either a plank or a pillar. 
If all these edges are pillars, then {\sl (ii)} would apply, so there must be at least one plank. 
Our goal is to construct, from the given path, another path from $u$ up to $v$ such that $u = u_{0} \rightarrow u_{1}$ is a plank. 
So, let $i$ be smallest such that $u_{i} \rightarrow u_{i+1}$ is a plank. 
Therefore $u_{i} \in \overline{u}$. 
If $u_{i} = u$, then we are done. 
Otherwise, let $r := \mytier(u_{i})$, so $u_{i} = \phi_{r-1}(u_{i-1})$. 
In particular, $u_{i-1} \in \mathcal{I}_{r-1}$ and $u_{i} \in \mathcal{F}_{r}$. 
Since $u_{i} \rightarrow u_{i+1}$ in $P_{r}$, then $u_{i+1} \in \mathcal{F}_{r}$, and therefore $u_{i+1} = \phi_{r-1}(u_{i}')$ for some $u_{i}' \in \mathcal{I}_{r-1}$. 
Moreover, $u_{i-1} \rightarrow u_{i}'$ in $\mathcal{I}_{r-1}$ since $u_{i} \rightarrow u_{i+1}$ in $\mathcal{F}_{r}$. 
So we have a path $u = u_{0} \rightarrow u_{1} \rightarrow \cdots \rightarrow u_{i-1} \rightarrow u_{i}' \rightarrow u_{i+1} \rightarrow \cdots \rightarrow u_{p} = v$ from $u$ up to $v$ whose first plank $u_{i-1} \rightarrow u_{i}'$ appears earlier than than in the path we were given. 
Continue to apply this process until we obtain a path from $u$ up to $v$ whose first edge is a plank. 
If we apply this process now to the path $u_{1} \rightarrow \cdots \rightarrow u_{p} = v$ from $u_{1}$ up to $v$, we get either a path of pillars or else a path that starts with a plank. 

Repeating this process we finally get a path $u = v_{0} \rightarrow v_{1} \rightarrow \cdots \rightarrow v_{j-1} \rightarrow v_{j} \rightarrow v_{j+1} \rightarrow \cdots \rightarrow v_{p} = v$ from $u$ up to $v$ such that all edges $v_{i} \rightarrow v_{i+1}$ are planks for $i \in [0,j-1]_{\mathbb{Z}}$ and all edges $v_{i} \rightarrow v_{i+1}$ are pillars for $i \in [j,p-1]_{\mathbb{Z}}$. 
Then $v' := v_{j} \in \overline{v} \cap P_{s}$, and $u < v'$ in $P_{s}$. 
Similarly, we can obtain a path $u = w_{0} \rightarrow w_{1} \rightarrow \cdots \rightarrow w_{k-1} \rightarrow w_{k} \rightarrow w_{k+1} \rightarrow \cdots \rightarrow w_{p} = v$ from $u$ up to $v$ such that all edges $w_{i} \rightarrow w_{i+1}$ are pillars for $i \in [0,k-1]_{\mathbb{Z}}$ and all edges $w_{i} \rightarrow w_{i+1}$ are plankss for $i \in [k,p-1]_{\mathbb{Z}}$. 
Then $u' := w_{k} \in \overline{u} \cap P_{t}$, and $u' < v$ in $P_{t}$. 

For {\sl (2)}, we suppose $v' = (\phi_{t-1} \circ \cdots \phi_{s})(v)$ for some $t > s$ and $v' \in P_{t}$. 
In particular, $v \in \mathcal{I}_{s}$, or else the expression `$(\phi_{t-1} \circ \cdots \phi_{s})(v)$' is not defined. 
Since $u \leq v$ in $P_{s}$, it follows that $u$ is in the down-set $\mathcal{I}_{s}$ as well. 
In particular, $\phi_{s}(u)$ is defined. 
Continuing, we have $phi_{s}(v) \in \mathcal{I}_{s+1}$, or else the expression `$(\phi_{t-1} \circ \cdots \phi_{s})(v)$' is not defined. 
Since $u \leq v$ in $P_{s}$ and since $\phi_{s}$ preserves order relations between $\mathcal{I}_{s}$ and $\mathcal{F}_{s+1}$, then $\phi_{s}(u) \leq \phi_{s}(v)$. 
Therefore, $\phi_{s}(u)$ is in the down-set $\mathcal{I}_{s+1}$ as well. 
Continuing in this way, we see that $u' := (\phi_{t-1} \circ \cdots \phi_{s})(v)$ is a well-defined member of $P_{t}$ and that $u' \leq v'$ within the tier $P_{t}$. 
The conditional claim in the third sentence of the statement of part {\sl (2)} of the lemma follows by entirely similar reasoning.\hfill\QED

\noindent 
{\bf \StackDefsTwo}\ \ Continue with the prior notation. 
Let `$\sim$' $:=$ `$\sim_{\phi_{1},\phi_{2},\ldots,\phi_{m-1}}$' be the equivalence relation on $P_{1} \myposetstack_{\phi_{1}} P_{2} \myposetstack_{\phi_{2}} \cdots \myposetstack_{\phi_{m-1}} P_{m}$ given by $u \sim v$ if and only if $v = (\phi_{r_j}\circ\cdots\circ\phi_{r_2}\circ\phi_{r_1})(u)$ for some possibly empty subsequence $(r_1,r_2,\ldots,r_j)$ of consecutive and increasing numbers from the set $[1,m-1]_{\mathbb{Z}}$. 
For the moment, declare `$\mystackedP^{\mybirdseye}$' to be the set $(P_{1} \myposetstack_{\phi_{1}} P_{2} \myposetstack_{\phi_{2}} \cdots \myposetstack_{\phi_{m-1}} P_{m})/\!\! \sim$ of equivalence classes under `$\sim$', i.e.\ the set 
$\displaystyle \big\{\overline{v} \ \rule[-1.5mm]{0.2mm}{4.5mm}\ v \in P_{1} \myposetstack_{\phi_{1}} P_{2} \myposetstack_{\phi_{2}} \cdots \myposetstack_{\phi_{m-1}} P_{m}\big\}$, with the obvious inherited vertex-coloring function $\vcolor^{\mybirdseye}\! :\! \mystackedP^{\mybirdseye} \rightarrow I$ given by $\vcolor^{\mybirdseye}(\overline{x}) = \vcolor^{(r)}(x)$ when $\overline{x}$ is the equivalence class of some $x \in P_{r}$. 
For $\overline{u}, \overline{v} \in \mystackedP^{\mybirdseye}$, say $\overline{u} \leq_{\mysmallbirdseye} \overline{v}$ if for every $u' \in \overline{u}$ there is some $v' \in \overline{v}$ and a sequence $(u'=u_{0},u_{1},\ldots,u_{p}=v')$, possibly with $p=0$, wherein $u_{i+1}$ succeeds $u_{i}$ in the sequence if and only if $u_{i} \rightarrow u_{i+1}$ is a plank or $u_{i+1} \rightarrow u_{i}$ is a pillar. 
Say any such sequence is a {\em contra-bonding path} from $u'$ to $v'$. 
Let $u^{\uparrow}$ and $u^{\downarrow}$ denote the members of the equivalence class $\overline{u}$ such that $\mytier(u^{\downarrow}) \leq \mytier(u') \leq \mytier(u^{\uparrow})$ for all $u' \in \overline{u}$. 

\noindent 
{\bf \StackLemma}\ \ 
{\sl In the above notation, let $u$ and $v$ be members of} $\mystackedP$, {\sl and assume that} $\mystackedP$ {\sl is connected.}\\ 
{\sl (1) The relation} $\leq_{\mysmallbirdseye}$ {\sl is a partial order on} $\mystackedP^{\mybirdseye}$.\\ 
{\sl (2) The following statements (a)--(g) are equivalent: (a)} $\overline{u} \leq_{\mysmallbirdseye} \overline{v}${\sl ; (b) There exists a contra-bonding path from $u^{\uparrow}$ to some member $v'$ of $\overline{v}$; (c) There exists a contra-bonding path from $u^{\uparrow}$ to $v^{\uparrow}$; (d) For all $v' \in \overline{v}$, there exists a contra-bonding path from some $u' \in \overline{u}$ to $v'$; (e) There exists a contra-bonding path from some $u' \in \overline{u}$ to $v^{\downarrow}$; (f) There exists a contra-bonding path from $u^{\downarrow}$ to $v^{\downarrow}$; and (g) For all $u' \in \overline{u}$, there exists a contra-bonding path from $u'$ to some $v' \in \overline{v}$.}\\ 
{\sl (3) We have $\overline{u} \rightarrow \overline{v}$ in} $\mystackedP^{\mybirdseye}$ {\sl if and only if there exist $u' \in \overline{u}$ and $v' \in \overline{v}$ such that} $\mytier(u') = \mytier(v')$ {\sl and $u' \rightarrow v'$ is a plank within this tier.  In this case, whenever $u'' \in \overline{u}$ and $v'' \in \overline{v}$ such that} $\mytier(u'') = \mytier(v'')$, {\sl then $u'' \rightarrow v''$ is a plank within this tier.}\\
{\sl (4) Suppose each tier comprising} $\mystackedP$ {\sl has a unique maximal element and a unique minimal element.  Say $u$ is the minimal element of the `lowest' tier poset $P_{1}$ and $v$ is the maximal element of the `highest' tier poset $P_{m}$. Then $\overline{u}$ is the unique minimal element of} $\mystackedP^{\mybirdseye}$ {\sl and $\overline{v}$ is its unique maximal element.}\\ 
{\sl (5) Suppose each tier comprising} $\mystackedP$ {\sl is isomorphic to a full-length sublattice of some product of two (nonempty) chains. Then} $\mystackedP^{\mybirdseye}$ {\sl is isomorphic to a full-length sublattice of some product of two (nonempty) chains, and is therefore a distributive lattice.}

{\em Proof.} For {\sl (1)}, the fact that $\overline{u} \leq_{\mysmallbirdseye} \overline{u}$ follows since we can take our contra-bonding path to be the singleton $(u')$ for every $u' \in \overline{u}$. 
Now suppose $\overline{u} \leq_{\mysmallbirdseye} \overline{v}$ and $\overline{v} \leq_{\mysmallbirdseye} \overline{w}$. 
Let $u' \in \overline{u}$. 
There is a $v' \in \overline{v}$ and a contra-bonding path from $u'$ to $v'$, and there is a $w' \in \overline{w}$ and a contra-bonding path from $v'$ to $w'$. 
Concatenate these two contra-bonding paths to get a contra-bonding path from $u'$ to $w'$. 
Therefore $\overline{u} \leq_{\mysmallbirdseye} \overline{w}$. 

We have established that $\leq_{\mysmallbirdseye}$ is reflexive and transitive. 
Note that for any contra-bonding path $(w_{0}, w_{1}, \ldots, w_{p})$, we have $\mytier(w_{i}) \geq \mytier(w_{i+1})$ for all $i \in [0,p-1]_{\mathbb{Z}}$. 
Now suppose $\overline{u} \leq_{\mysmallbirdseye} \overline{v}$ and $\overline{v} \leq_{\mysmallbirdseye} \overline{u}$. 
So there is a contra-bonding path from $u_{0} := u$ to some $v_{1} \in \overline{v}$ and a contra-bonding path from $v_{1}$ to some $u_{1} \in \overline{u}$. 
Now, $\mytier(u_{1}) \leq \mytier(u_{0})$. 
If $\mytier(u_{1}) = \mytier(u_{0})$, then the edges in both of our contra-bonding paths are all planks, which forces $u_{0} = v_{1} = u_{1}$, and therefore $\overline{u} = \overline{v}$. 
Otherwise, $\mytier(u_{1}) < \mytier(u_{0})$.  
Then there is a contra-bonding path from $u_{1}$ to some $v_{2} \in \overline{v}$ and a contra-bonding path from $v_{2}$ to some $u_{2} \in \overline{u}$. 
As before, $\mytier(u_{2}) \leq \mytier(u_{1})$. 
If $\mytier(u_{2}) = \mytier(u_{1})$, then the edges in our two most recent contra-bonding paths are all planks, which forces $u_{1} = v_{2} = u_{2}$, and therefore $\overline{u} = \overline{v}$.
Otherwise, $\mytier(u_{2}) < \mytier(u_{1})$.  
Then there is a contra-bonding path from $u_{2}$ to some $v_{3} \in \overline{v}$ and a contra-bonding path from $v_{3}$ to some $u_{3} \in \overline{u}$. 
As before, $\mytier(u_{2}) \leq \mytier(u_{1})$. 
We continue this process, obtaining a sequence $\mytier(u_{0}) \geq \mytier(u_{1}) \geq \mytier(u_{2}) \geq \mytier(u_{3}) \geq \cdots$. 
This sequence cannot strictly decrease indefinitely, so for some $j \in [0,\infty)_{\mathbb{Z}}$ we get $\mytier(u_{j}) = \mytier(u_{j+1})$ which will yield $u_{j} = v_{j+1} = v_{j}$ and hence $\overline{u} = \overline{v}$. 

For {\sl (2)}, observe that {\sl (a)} is equivalent to {\sl (g)} by definition. 
Given this equivalence, then {\sl (b)} is a special case of {\sl (g)}. 
Now assume {\sl (b)}. 
We show by induction on $p$ that from any given contra-bonding path of length $p$ from $u^{\uparrow}$ to some $v' \in \overline{v}$, we can obtain a contra-bonding path from $u^{\uparrow}$ to $v^{\uparrow}$. 
If $p=0$, then $u^{\uparrow} = v'$ and therefore $u^{\uparrow} = v^{\uparrow}$; so, we get a length zero contra-bonding path from $u^{\uparrow}$ to $v^{\uparrow}$. 
Now say for some $p \in [0,\infty)_{\mathbb{Z}}$, suppose that from any length $q \leq p$ contra-bonding path from $u^{\uparrow}$ to some $v' \in \overline{v}$ we can obtain a contra-bonding path from $u^{\uparrow}$ to $v^{\uparrow}$. 
If we are given a a length $p+1$ contra-bonding path from $u^{\uparrow}$ to some $v' \in \overline{v}$, then, without loss of generality, we may assume that $v'$ is the largest member of $\overline{v}$ that is part of our contra-bonding path and that $u^{\uparrow} \not\in \overline{v}$. 
If $v' = v^{\uparrow}$, then we are done. 
So, suppose $v' < v^{\uparrow}$. 
Then, our contra-bonding path includes a plank $w \rightarrow v'$ within the tier $P_{s}$, where $\mytier(w) = \mytier(v') =: s$. 
We can apply the induction hypothesis to the contra-bonding path from $u^{\uparrow}$ to $w \in \overline{w}$ to get a contra-bonding path from $u^{\uparrow}$ to $w^{\uparrow}$. 
Now, $w < v^{\uparrow}$ in $\mystackedP$, so by \StackEarlyLemma.1, criteria {\sl (iii)} must apply to give us some $w' \in \overline{w}$ such that $\mytier(w') = \mytier(v^{\uparrow}) =: t$ and $w' < v^{\uparrow}$ in $P_{t}$. 
Since $w \rightarrow v'$ in $P_{s}$, then we must also have $w' \rightarrow v^{\uparrow}$ in $P_{t}$. 
So now take our contra-bonding path from $u^{\uparrow}$ to $w^{\uparrow}$, the follow pillar down to $w'$, and then follow the plank in $P_{t}$ to $v^{\uparrow}$. 
This completes the induction argument and our proof that {\sl (b)} implies {\sl (c)}. 
Taking $u' := u^{\uparrow}$, we see that {\sl (d)} follows from {\sl (c)}, as there is a sequence of pillars from $v^{\uparrow}$ down to any $v' \in \overline{v}$. 
Now {\sl (e)} is just a special case of {\sl (d)}. 
Our proof that {\sl (e)} $\Longrightarrow$ {\sl (f)} is entirely similar to our proof that {\sl (b)} $\Longrightarrow$ {\sl (c)}. 
Taking $v' := v^{\downarrow}$, we see that {\sl (g)} follows from {\sl (f)}, as there is a sequence of pillars from $u' \in \overline{u}$ down to any $u^{\downarrow}$. 
This completes our proof of {\sl (2)}. 

For {\sl (3)}, suppose $u \rightarrow v$ in some tier $P_{t}$. 
So, by {\sl (2)}, $\overline{u} <_{\mysmallbirdseye} \overline{v}$. 
Now suppose $\overline{u} \leq_{\mysmallbirdseye} \overline{w} \leq_{\mysmallbirdseye} \overline{v}$. 
Then there are contra-bonding paths from $u^{\uparrow}$ to $w^{\uparrow}$ and from $w^{\uparrow}$ to $v^{\uparrow}$, and there are contra-bonding paths from $u^{\downarrow}$ to $w^{\downarrow}$ and from $w^{\downarrow}$ to $v^{\downarrow}$. 
Say $t_{1} := \mytier(u^{\uparrow})$, $t_{2} := \mytier(w^{\uparrow})$, and $t_{3} := \mytier(v^{\uparrow})$, and let $s_{1} := \mytier(u^{\downarrow})$, $s_{2} := \mytier(w^{\downarrow})$, and $s_{3} := \mytier(v^{\downarrow})$. 
Then 
\[\begin{array}{ccccc}
t_{1} & \geq & t_{2} & \geq & t_{3}\\ 
\rotatebox{90}{$\leq$} & & \rotatebox{90}{$\leq$} & & \rotatebox{90}{$\leq$}\\
s_{1} & \geq & s_{2} & \geq & s_{3}
\end{array}\]
We also have $s_{1} \leq t \leq t_{1}$, and hence $s_{2} \leq t$. 
Also $s_{3} \leq t \leq t_{3}$, so  $t \leq t_{2}$. 
Then $\overline{w}$ has an element in $P_{t}$, which we name $w$ since, as yet, no tier has been declared for this representative of $\overline{w}$. 
If $u \leq w \leq v$, then, since $u \rightarrow v$, we get $u=w$ or $w=v$ and hence $\overline{u}=\overline{w}$ or $\overline{w}=\overline{v}$. 
So we are left to consider the two cases that in $P_{t}$, we have $u \not\leq w$ or $w \not\leq v$. 
We first consider the case that $w \not\leq v$. 
If $w^{\downarrow} \in P_{s_{3}}$, then the fact that we have a contra-bonding sequence from $w^{\downarrow}$ to $v^{\downarrow}$ will force $w^{\downarrow} \leq v^{\downarrow}$ in $P_{s_{3}}$.  
But, by applying $\phi_{t-1} \circ \cdots \circ \phi_{s_{3}}$, we see that $w^{\downarrow} \leq v^{\downarrow}$ in $P_{s_{3}}$ forces $w \leq v$ in $P_{t}$, contrary to our hypothesis that $w \not\leq v$. 
We conclude that $w^{\downarrow} \not\in P_{s_{3}}$ and hence $\mytier(w^{\downarrow}) > \mytier(v^{\downarrow})$. 
So, in our given contra-bonding path from $w^{\downarrow}$ to $v^{\downarrow}$, let $x_{1}$ be the last element in the same tier as $w^{\downarrow}$; since our contra-bonding path must have pillars (else $\mytier(w^{\downarrow}) = \mytier(v^{\downarrow})$), then within said path there exists a pillar below $x_{1}$ and hence $\mytier(x_{1}^{\downarrow}) < \mytier(x_{1}) = s_{2}$. 
Then, by part {\sl (2)}, we get $\overline{x_{1}} \leq_{\mysmallbirdseye} \overline{v}$, so there exists a contra-bonding path from $x_{1}^{\downarrow}$ to $v^{\downarrow}$. 
Let us say that $\mytier(x_{1}^{\downarrow}) = \mytier(v^{\downarrow})$; then we get $x_{1}^{\downarrow} \leq v^{\downarrow}$ in $P_{s_{3}}$. 
But then we get $w^{\downarrow} \leq x_{1} \leq v_{1}$ in $P_{s_{2}}$ (we know $\overline{v} \cap P_{s_{2}}$ includes some element $v_{1}$ since $s_{3} < s_{2} \leq t$ and $v \in P_{t}$), and therefore there must be some $x_{1}' \in \overline{x_{1}} \cap P_{t}$ with $w \leq x_{1}' \leq v$, contradicting $w \not\leq v$ in $P_{t}$. 
Therefore, $\mytier(x_{1}^{\downarrow}) > \mytier(v^{\downarrow})$. 
So, in our given contra-bonding path from $x_{1}^{\downarrow}$ to $v^{\downarrow}$, let $x_{2}$ be the last element in the same tier as $x_{1}^{\downarrow}$; since our contra-bonding path must have pillars (else $\mytier(x_{1}^{\downarrow}) = \mytier(v^{\downarrow})$), then within said path there exists a pillar below $x_{2}$ and hence $\mytier(x_{2}^{\downarrow}) < \mytier(x_{1})$. 
Then, by part {\sl (2)}, we get $\overline{x_{2}} \leq_{\mysmallbirdseye} \overline{v}$, so there exists a contra-bonding path from $x_{2}^{\downarrow}$ to $v^{\downarrow}$. 
But, as with $x_{1}^{\downarrow}$, we will see that $\mytier(x_{2}^{\downarrow}) > \mytier(v^{\downarrow})$, and, as we continue this process, $\mytier(x_{i}^{\downarrow}) = \mytier(v^{\downarrow})$ for each $i \in [1,\infty)_{\mathbb{Z}}$, contradicting the finiteness of our skew-stacked poset. 
Therefore it cannot be the case that $w \not\leq v$. 
We similarly rule out $u \not\leq w$. 
Therefore, $\overline{u} \leq_{\mysmallbirdseye} \overline{w} \leq_{\mysmallbirdseye} \overline{v}$ forces $\overline{u} = \overline{w}$ or $\overline{w}=\overline{v}$, then we conclude that $u \rightarrow v$ in $P_{t}$ implies that $\overline{u} \rightarrow \overline{v}$ in $\mystackedP^{\mybirdseye}$. 

To complete {\sl (3)}, suppose now that $\overline{u} \rightarrow \overline{v}$ in $\mystackedP^{\mybirdseye}$. 
Then, by {\sl (2)}, there is a contra-bonding path from $u$ to $v$. 
Let $u' \in \overline{u}$ be the last element from $\overline{u}$ that is part of said contra-bonding path, and similarly let $v' \in \overline{v}$ be the first such.   
Let $x'$ be the last element from our given contra-bonding path that shares a tier with $u'$, and similarly let $y'$ be the first element from said contra-bonding path that is in the same tier as $v'$. 
By {\sl (2)}, $\overline{u} \leq_{\mysmallbirdseye} \overline{x'} \leq_{\mysmallbirdseye} \overline{y'} \leq_{\mysmallbirdseye} \overline{v}$. 
Invoking the fact that $\overline{u} \rightarrow \overline{v}$ in $\mystackedP^{\mybirdseye}$, there are three possibilities to evaluate: $\overline{u} = \overline{x'} = \overline{y'}$, $\overline{x'} = \overline{y'} = \overline{v}$, or $\overline{u}=\overline{x'}$ and $\overline{y'}=\overline{v}$. 
With $\overline{u} = \overline{x'} = \overline{y'}$, then $x' = u'$ and $y'$ is an element of $\overline{u}$ in the same tier as $v'$. 
Of course $y' < v'$ is a strict inequality within this tier since $y'=v'$ would force the equalities $\overline{u} = \overline{x'} = \overline{y'}=\overline{v}$, contrary to our hypothesis that $\overline{v}$ covers $\overline{u}$.
If there is some $w$ in this tier such that $y' < w < v'$, then we would have $\overline{u} <_{\mysmallbirdseye} \overline{w} <_{\mysmallbirdseye} \overline{v}$, again contrary to $\overline{u} \rightarrow \overline{v}$. 
So, $y' \rightarrow v'$ with $y' \in \overline{u}$. 
A similar analysis of the case $\overline{x'} = \overline{y'} = \overline{v}$ leads us to conclude that $u' \rightarrow x'$ with $x' \in \overline{v}$.
Now suppose $\overline{u}=\overline{x'}$ and $\overline{y'}=\overline{v}$. 
In particular, we have $u'=x'$ and $y'=v'$. 
Suppose $\mytier(u^{\downarrow}) \leq \mytier(v^{\uparrow})$. 
Then in some tier $P_{t}$ with $\mytier(u^{\downarrow}) \leq t \leq \mytier(v^{\uparrow})$, there exist $u^{(t)} \in \overline{u}$ and $v^{(t)} \in \overline{v}$ such that $u^{(t)} < v^{(t)}$. 
Reasoning similar to our analysis of a previous case means that $u^{(t)} \rightarrow v^{(t)}$ in $P_{t}$. 
So now say $\mytier(u^{\downarrow}) > \mytier(v^{\uparrow})$. 
Any contra-bonding path from $u^{\downarrow}$ to $v^{\downarrow}$, which exists by {\sl (2)}, will require that the first step is a plank and not a pillar. 
Say $u^{\downarrow} \rightarrow x$ is this plank. 
But now $\overline{u} < \overline{x} \leq \overline{v}$, which means that $\overline{x} = \overline{v}$, and therefore $\mytier(u^{\downarrow}) = \mytier(x) \leq \mytier(v^{\uparrow})$, contradicting the fact that $\mytier(u^{\downarrow}) > \mytier(v^{\uparrow})$. 
This exhausts our analysis of cases, and we conclude that, when $\overline{u} \rightarrow \overline{v}$ in $\mystackedP^{\mybirdseye}$, there exist $u' \in \overline{u}$ and $v' \in \overline{v}$ such that $\mytier(u') = \mytier(v')$ and $u' \rightarrow v'$ within this tier. 

For {\sl (4)}, we note that if $u_{t}$ is the unique minimal element of $P_{t}$ and $u_{t+1}$ the unique minimal element of $P_{t+1}$, then, since $\mystackedP$ is connected, there must be a pillar from $u_{t}$ to some element $u_{t}' \in P_{t+1}$. Since $u_{t+1} \leq u_{t}'$ in $P_{t+1}$, we conclude that there is a path from $u_{t+1}$ up to $u_{t}'$. 
Follow this path with a step from $u_{t}'$ to $u_{t}$ along the pillar $u_{t}' \leftarrow u_{t}$ to get a contra-bonding path from $u_{t+1}$ to $u_{t}$. 
By concatenating such paths for all $t$, we can get a contra-bonding path from $u = u_{m}$ to $u_{t}$. 
For any $x \in P_{t}$, we can, therefore, obtain a contra-bonding path from $u$ to $x$. 
That is, for any $\overline{x} \in \mystackedP^{\mybirdseye}$, we have $\overline{u} \leq_{\mysmallbirdseye} \overline{x}$. 
So, $\overline{u}$ is the unique minimal element of $\mystackedP^{\mybirdseye}$. 
Similarly see that $\overline{v}$ is the unique maximal element of $\mystackedP^{\mybirdseye}$.

For {\sl (5)}, we make a preliminary observation. 
Suppose, for $i \in \{1,2\}$, that $Q_{i}$ is a full-length sublattice of $[0,x_{i}]_{\mathbb{Z}} \times [0,y_{i}]_{\mathbb{Z}}$ and that for a down-set $\mathcal{I}$ from $Q_{1}$ and an up-set $\mathcal{F}$ from $Q_{2}$ there is an isomorphism $\phi: \mathcal{I} \longrightarrow \mathcal{F}$. 
Let $\eelt := (1,0)$ and $\felt := (0,1)$, and set $\eelt' := \felt$ and $\felt' := \eelt$. 
Whenever $u \rightarrow v$ in $Q_{i}$, then it must be the case that $u + \eelt = v$ (in which case we write $u \mylongarrow{\eelt} v$) or $u + \felt = v$ (in which case we write $u \mylongarrow{\felt} v$). 
If $(0,0)$ is covered by an element from $\mathcal{I}$, then we have $(0,0) \mylongarrow{\eelt} (1,0)$ or $(0,0) \mylongarrow{\felt} (0,1)$ in $\mathcal{I}$. 
Say this edge in $\mathcal{I}$ is $(0,0) \mylongarrow{\gelt} \gelt$ in $\mathcal{I}$. 
Then $\phi(0,0) \rightarrow \phi(\gelt)$ in $\mathcal{F}$, so it must be the case that $\phi(0,0) \mylongarrow{\gelt} \phi(\gelt)$ or $\phi(0,0) \mylongarrow{\gelt'} \phi(\gelt)$. 
In the former case, one can now see that $u \mylongarrow{\helt} v$ in $\mathcal{I}$ with $\helt \in \{\eelt,\felt\}$ if and only if $\phi(u) \mylongarrow{\helt} \phi(v)$ in $\mathcal{F}$. 
In the latter case, let $\eta: [0,x_{i}]_{\mathbb{Z}} \times [0,y_{i}]_{\mathbb{Z}} \longrightarrow [0,y_{i}]_{\mathbb{Z}} \times [0,x_{i}]_{\mathbb{Z}}$ be given by $\eta(x,y) := (y,x)$, let $Q_{2}' := \eta(Q_{2})$, let $\mathcal{F}' := \eta(\mathcal{F})$, and let $\phi' := \eta \circ \phi$. 
One can now see that $u \mylongarrow{\helt} v$ in $\mathcal{I}$ with $\helt \in \{\eelt,\felt\}$ if and only if $\phi'(u) \mylongarrow{\gelt} \phi'(v)$ in $\mathcal{F}'$. 
That is, we can adjust our description of $Q_{2}$ as a full-length sublattice of a product of two (nonempty) chains so that whenever $u \mylongarrow{\gelt} v$ in $\mathcal{I}$ for $\gelt \in \{\eelt,\felt\}$, then we have $\phi'(u) \mylongarrow{\gelt} \phi'(v)$ in $\mathcal{F}'$. 
Call this process the {\em pairwise-orientation of} $Q_{2}$ {\em to} $Q_{1}$. 

In the setting of our skew-stack, let us suppose that each $P_{t}$ is identified as a full-length sublattice of $[0,x_{t}]_{\mathbb{Z}} \times [0,y_{t}]_{\mathbb{Z}}$ for nonnegative integers $x_{t}$ and $y_{t}$. 
In particular, by \FullLengthWithinProduct.2, each $P_{t}$ is a distributive lattice; for this same reason, once we show that $\mystackedP^{\mybirdseye}$ is a full-length sublattice of a product of two nonempty chains, then we can conclude that $\mystackedP^{\mybirdseye}$ is a distributive lattice. 
Let us assume that, for $t \in [1,m-1]_{\mathbb{Z}}$, $P_{t+1}$ is pairwise-oriented to the preceding $P_{t}$ so that, under the isomorphism $\phi_{t}: \mathcal{I}_{t} \longrightarrow \mathcal{F}_{t+1}$ we have $u \mylongarrow{\gelt} v$ in $\mathcal{I}_{t}$ for some $\gelt \in \{\eelt,\felt\}$ if and only if $\phi_{t}(u) \mylongarrow{\gelt} \phi_{t}(v)$ in $\mathcal{F}_{t+1}$. 
For such $t$, set $(p_{t+1},q_{t+1}) := \phi_{t}(0,0)$. 
Therefore, for any $u$ in $\mathcal{I}_{t}$, we have $\phi_{t}(u) = u+(p_{t+1},q_{t+1})$ in $\mathcal{F}_{t+1}$. 
Let 
\[\myfancyQ := \left\{(x,y) \in \mathbb{Z} \times \mathbb{Z}\ \rule[-3.75mm]{0.15mm}{9mm}\, \left(x - \sum_{s=t+1}^{m}p_{s},y - \sum_{s=t+1}^{m}q_{s}\right) \in P_{t} \mbox{ for some } t \in [1,m]_{\mathbb{Z}}\right\}\]
with partial order induced by the component-wise order on $\mathbb{Z} \times \mathbb{Z}$. 
Observe that $\myfancyQ$ is a subset of $[0,x_{1}+\sum_{s=2}^{m}p_{s}]_{\mathbb{Z}} \times [0,y_{1}+\sum_{s=2}^{m}q_{s}]_{\mathbb{Z}}$ and that both $(0,0)$ and $(x_{1}+\sum_{s=2}^{m}p_{s},y_{1}+\sum_{s=2}^{m}q_{s})$ are in $\myfancyQ$. 
Now, within each $P_{t}$, for $t \in [2,m]_{\mathbb{Z}}$, there is a path from $(0,0)$ up to $(p_{t},q_{t})$; by adding $(\sum_{s=t+1}^{m}p_{s},\sum_{s=t+1}^{m}q_{s})$ to each pair in this path, we get a path in $\myfancyQ$ from $(\sum_{s=t+1}^{m}p_{s},\sum_{s=t+1}^{m}q_{s})$ up to $(p_{t}+\sum_{s=t+1}^{m}p_{s},q_{t}+\sum_{s=t+1}^{m}q_{s}) = (\sum_{s=t}^{m}p_{s},\sum_{s=t}^{m}q_{s})$. 
By concatenation, we can therefore get a path in $\myfancyQ$ from $(0,0)$ up to $(p_{m},q_{m})$ up to {\sl etc} up to $(\sum_{s=t}^{m}p_{s},\sum_{s=t}^{m}q_{s})$ up to {\sl etc} up to $(x_{1}+\sum_{s=2}^{m}p_{s},y_{1}+\sum_{s=2}^{m}q_{s})$, and this is also a path within the chain product $[0,x_{1}+\sum_{s=2}^{m}p_{s}]_{\mathbb{Z}} \times [0,y_{1}+\sum_{s=2}^{m}q_{s}]_{\mathbb{Z}}$. 

We claim that $\myfancyQ$ is closed under component-wise joins and meets. 
That is, for $(x,y)$ and $(x',y')$ in $\myfancyQ$, we must show that the pairs $(\min(x,x'),\min(y,y'))$ and $(\max(x,x'),\max(y,y'))$ are both in $\myfancyQ$. 
We only do this for the component-wise join, as the argument for the component-wise meet is analogous. 
Let us say that $(a,b) := \big(x-\sum_{s=r+1}^{m}p_{s},y-\sum_{s=r+1}^{m}q_{s}\big) \in P_{r}$ and that $(a',b') := \big(x'-\sum_{s=t+1}^{m}p_{s},y'-\sum_{s=t+1}^{m}q_{s}\big) \in P_{t}$, with $r \leq t$. 
It suffices to consider the case that the pairs $(x,y)$ and $(x',y')$ are incomparable. 
Without loss of generality, assume that $x' \geq x$ and $y' \leq y$. 
From $x' \geq x$ it follows that $a'-\sum_{s=r+1}^{t}p_{s} \geq a$ and therefore $a' \geq \sum_{s=r+1}^{t}p_{s}$. 
Then, $(a',q_{t}) \in \mathcal{F}_{t}$; $\phi^{-1}_{t-1}(a',q_{t}) = (a'-p_{t},0) \in \mathcal{I}_{t-1}$; $(a'-p_{t},q_{t-1}) \in \mathcal{F}_{t-1}$; $\phi^{-1}_{t-2}(a'-p_{t},q_{t-1}) = (a'-p_{t}-p_{t-1},0) \in \mathcal{I}_{t-2}$; {\em et cetera}; $\big(a'-\sum_{s=r+2}^{t}p_{s},q_{r+1}\big) \in \mathcal{F}_{r+1}$; and $\phi^{-1}_{r}\big(a'-\sum_{s=r+2}^{t}p_{s},q_{r+1}\big) = \big(a'-\sum_{s=r+1}^{t}p_{s},0\big) \in \mathcal{I}_{r}$. 
Therefore $\big(a'-\sum_{s=r+1}^{t}p_{s},b\big)$, which is the component-wise join within $P_{r}$ of $(a,b)$ and $\big(a'-\sum_{s=r+1}^{t}p_{s},0\big)$, must also reside in $P_{r}$ by \FullLengthWithinProduct.2. 
So, $\big(a'-\sum_{s=r+1}^{t}p_{s}+\sum_{s=r+1}^{m}p_{s},b+\sum_{s=r+1}^{m}q_{s}\big) = \big(a'+\sum_{s=t+1}^{m}p_{s},b+\sum_{s=r+1}^{m}q_{s}\big) = (x',y) \in \myfancyQ\, $. 

Now define $\psi: \mystackedP \longrightarrow \myfancyQ$ by declaring that, for any $(x,y) \in P_{t}$, 
\[\psi(x,y) := \left(x + \sum_{s=t+1}^{m}p_{s},y + \sum_{s=t+1}^{m}q_{s}\right).\] 
Of course, if $(x,y) \in \myfancyQ\, $, then $(a,b) := \big(x - \sum_{s=t+1}^{m}p_{s},y - \sum_{s=t+1}^{m}q_{s}\big) \in P_{t}$ for some $t$, in which case $\psi(a,b) = (x,y)$;  that is, $\psi(\mystackedP) = \myfancyQ\, $. 
Observe that for any $(x,y) \in \mathcal{I}_{r}$ in some tier $P_{r}$, there is a path  $(0,0) \mylongarrow{\gelt_{1}} \cdots \mylongarrow{\gelt_{x+y}} (x,y)$ in $\mathcal{I}_{r}$ where each $\gelt_{i} \in \{\eelt,\felt\}$. 
Thanks to our pairwise-orientation of $P_{r+1}$ to $P_{r}$, we get $\phi_{r}(x,y) = \phi_{r}(0,0) + x(1,0) + y(0,1) = (x+p_{r+1},y+q_{r+1})$. 

We claim that $u \sim u'$ in $\mystackedP$ if and only if $\psi(u) = \psi(u')$. 
For specificity, we take $\mytier(u) := r \leq t =: \mytier(u')$ and write $u = (x,y)$ and $u' = (x',y')$.
First suppose $u \sim u'$. 
By our observation at the end of the last paragraph, we have $u' = (\phi_{t-1} \circ \cdots \circ \phi_{r})(x,y) = (x+\sum_{s=r+1}^{t}p_{s},y+\sum_{s=r+1}^{t}q_{s})$. 
So, $\psi(u') = \psi(x+\sum_{s=r+1}^{t}p_{s},y+\sum_{s=r+1}^{t}q_{s}) = (x+\sum_{s=r+1}^{t}p_{s}+\sum_{s=t+1}^{m}p_{s},y+\sum_{s=r+1}^{t}q_{s}+\sum_{s=t+1}^{m}q_{s}) = (x+\sum_{s=r+1}^{m}p_{s},y+\sum_{s=r+1}^{m}q_{s}) = \psi(u)$. 
On the other hand, say $\psi(u) = \psi(u')$, so $(x+\sum_{s=r+1}^{m}p_{s},y+\sum_{s=r+1}^{m}q_{s}) = (x'+\sum_{s=t+1}^{m}p_{s},y'+\sum_{s=t+1}^{m}q_{s})$ and therefore $u' = (x',y') = (x+\sum_{s=r+1}^{t}p_{s},y+\sum_{s=r+1}^{t}q_{s})$. 
But again, by our observation at the end of the last paragraph, we have $(x+\sum_{s=r+1}^{t}p_{s},y+\sum_{s=r+1}^{t}q_{s}) =   (\phi_{t-1} \circ \cdots \circ \phi_{r})(x,y)$, and therefore $u' = (\phi_{t-1} \circ \cdots \circ \phi_{r})(u)$. 
Then $u \sim u'$ in $\mystackedP$. 

Now let $\overline{\psi}: \mystackedP^{\mybirdseye} \longrightarrow \myfancyQ$ be given by $\overline{\psi}(\overline{u}) := \psi(u)$, which is well-defined and injective by the preceding paragraph and is surjective by the paragraph before that. 
We claim that $\overline{u} \rightarrow \overline{v}$ in $\mystackedP^{\mybirdseye}$ if and only if $\overline{\psi}(\overline{u}) \rightarrow \overline{\psi}(\overline{v})$ in $\myfancyQ\, $. 
Suppose $\overline{u} \rightarrow \overline{v}$ in $\mystackedP^{\mybirdseye}$. 
By part {\sl (3)}, we have $u' \rightarrow v'$ within some tier $P_{t}$ where $u' \in \overline{u}$ and $v' \in \overline{v}$. 
Taking $u' := (x',y')$, then $v' = u'+\gelt$ where $\gelt \in \{\eelt,\felt\}$. 
Write $\gelt := (\epsilon_{1},\epsilon_{2})$, so $v' = (x'+\epsilon_{1},y'+\epsilon_{2})$. 
Now, $\overline{\psi}(\overline{v}) = \psi(v') = \psi(x'+\epsilon_{1},y'+\epsilon_{2}) = \left(x'+\epsilon_{1} + \sum_{s=t+1}^{m}p_{s},y' +\epsilon_{2}+ \sum_{s=t+1}^{m}q_{s}\right) = \psi(x',y')+\gelt = \psi(u')+\gelt = \overline{\psi}(\overline{u})+\gelt$, and hence $\overline{\psi}(\overline{u}) \rightarrow \overline{\psi}(\overline{v})$ in $\myfancyQ\, $. 
For the converse, suppose $\overline{\psi}(\overline{u}) \rightarrow \overline{\psi}(\overline{v})$ in $\myfancyQ\, $. 
In particular, $\psi(u) + \gelt = \psi(v)$ for some $\gelt \in \{\eelt,\felt\}$, which we write as $\gelt := (\epsilon_{1},\epsilon_{2})$. 
We claim that there exist $u' \sim u$ and $v' \sim v$ in $\mystackedP$ such that $u', v' \in P_{t}$. 
If so, then we can see as follows that $u' \mylongarrow{\gelt} v'$ in $P_{t}$. 
Since $\mystackedP$ is connected, then so is $P_{t}$, and we may consider shortest paths from $u'$ to $v'$. 
\hfill As a full-length sublattice of the chain product $[0,x_{t}]_{\mathbb{Z}} \times [0,y_{t}]_{\mathbb{Z}}$, any shortest path in $P_{t}$ is

\newcommand{\ESixFlat}{\setlength{\unitlength}{1cm}
\begin{picture}(7.5,12)
\put(0,10){\CyanNode}
\put(0.7,10.3){\color{Cyan}{\em \footnotesize 1}}
\put(0,4){\CyanNode}
\put(-0.05,4){\color{Cyan}{\em \footnotesize 1}}
\put(1.5,9){\RedNode}
\put(2.2,9.3){\color{Red}{\em \footnotesize 3}}
\put(1.5,5){\RedNode}
\put(1.5,5.3){\color{Red}{\em \footnotesize 3}}
\put(1.5,3){\RedNode}
\put(1.45,3){\color{Red}{\em \footnotesize 3}}
\put(3,8){\PurpleNode}
\put(3.7,8.3){\color{Purple}{\em \footnotesize 4}}
\put(3,6){\PurpleNode}
\put(2.95,6.1){\color{Purple}{\em \footnotesize 4}}
\put(3,4){\PurpleNode}
\put(3.75,4.15){\color{Purple}{\em \footnotesize 4}}
\put(3,2){\PurpleNode}
\put(2.95,2){\color{Purple}{\em \footnotesize 4}}
\put(1.5,7){\SpringGreenNode}
\put(1.5,7){\color{SpringGreen}{\em \footnotesize 2}}
\put(4.5,3){\SpringGreenNode}
\put(5.2,2.9){\color{SpringGreen}{\em \footnotesize 2}}
\put(4.5,7){\OliveGreenNode}
\put(5.2,7.3){\color{OliveGreen}{\em \footnotesize 5}}
\put(4.5,5){\OliveGreenNode}
\put(5.25,5){\color{OliveGreen}{\em \footnotesize 5}}
\put(4.5,1){\OliveGreenNode}
\put(4.45,1){\color{OliveGreen}{\em \footnotesize 5}}
\put(6,6){\BurntOrangeNode}
\put(6.7,5.85){\color{BurntOrange}{\em \footnotesize 6}}
\put(6,0){\BurntOrangeNode}
\put(5.95,0){\color{BurntOrange}{\em \footnotesize 6}}
\thicklines
\put(3.625,2.275){\color{Black}\qbezier(0,0)(0.6,0.42)(1.2,0.84)}
\put(2.125,3.275){\color{Black}\qbezier(0,0)(0.6,0.42)(1.2,0.84)}
\put(3.625,4.275){\color{Black}\qbezier(0,0)(0.6,0.42)(1.2,0.84)}
\put(5.125,5.275){\color{Black}\qbezier(0,0)(0.6,0.42)(1.2,0.84)}
\put(0.625,4.275){\color{Black}\qbezier(0,0)(0.6,0.42)(1.2,0.84)}
\put(2.125,5.275){\color{Black}\qbezier(0,0)(0.6,0.42)(1.2,0.84)}
\put(3.625,6.275){\color{Black}\qbezier(0,0)(0.6,0.42)(1.2,0.84)}
\put(2.125,7.275){\color{Black}\qbezier(0,0)(0.6,0.42)(1.2,0.84)}
\put(6.325,0.275){\color{Black}\qbezier(0,0)(-0.6,0.42)(-1.2,0.84)}
\put(4.825,1.275){\color{Black}\qbezier(0,0)(-0.6,0.42)(-1.2,0.84)}
\put(3.325,2.275){\color{Black}\qbezier(0,0)(-0.6,0.42)(-1.2,0.84)}
\put(1.825,3.275){\color{Black}\qbezier(0,0)(-0.6,0.42)(-1.2,0.84)}
\put(3.325,4.275){\color{Black}\qbezier(0,0)(-0.6,0.42)(-1.2,0.84)}
\put(4.825,3.275){\color{Black}\qbezier(0,0)(-0.6,0.42)(-1.2,0.84)}
\put(4.825,5.275){\color{Black}\qbezier(0,0)(-0.6,0.42)(-1.2,0.84)}
\put(3.325,6.275){\color{Black}\qbezier(0,0)(-0.6,0.42)(-1.2,0.84)}
\put(6.325,6.275){\color{Black}\qbezier(0,0)(-0.6,0.42)(-1.2,0.84)}
\put(4.825,7.275){\color{Black}\qbezier(0,0)(-0.6,0.42)(-1.2,0.84)}
\put(3.325,8.275){\color{Black}\qbezier(0,0)(-0.6,0.42)(-1.2,0.84)}
\put(1.825,9.275){\color{Black}\qbezier(0,0)(-0.6,0.42)(-1.2,0.84)}
\end{picture}
}
\newcommand{\ESixFlatTwo}{\setlength{\unitlength}{1cm}
\begin{picture}(7.5,12)
\put(6,10){\BurntOrangeNode}
\put(6,10.1){\color{BurntOrange}{\em \footnotesize 6}}
\put(6,4){\BurntOrangeNode}
\put(6.75,4){\color{BurntOrange}{\em \footnotesize 6}}
\put(4.5,9){\OliveGreenNode}
\put(4.5,9.1){\color{OliveGreen}{\em \footnotesize 5}}
\put(4.5,5){\OliveGreenNode}
\put(5.25,5.3){\color{OliveGreen}{\em \footnotesize 5}}
\put(4.5,3){\OliveGreenNode}
\put(5.25,3){\color{OliveGreen}{\em \footnotesize 5}}
\put(3,8){\PurpleNode}
\put(3,8.1){\color{Purple}{\em \footnotesize 4}}
\put(3,6){\PurpleNode}
\put(3.85,6.1){\color{Purple}{\em \footnotesize 4}}
\put(3,4){\PurpleNode}
\put(3.85,4.1){\color{Purple}{\em \footnotesize 4}}
\put(3,2){\PurpleNode}
\put(3.75,2){\color{Purple}{\em \footnotesize 4}}
\put(4.5,7){\SpringGreenNode}
\put(5.25,7){\color{SpringGreen}{\em \footnotesize 2}}
\put(1.5,3){\SpringGreenNode}
\put(2.35,3.15){\color{SpringGreen}{\em \footnotesize 2}}
\put(1.5,7){\RedNode}
\put(1.5,7.1){\color{Red}{\em \footnotesize 3}}
\put(1.5,5){\RedNode}
\put(2.35,5.1){\color{Red}{\em \footnotesize 3}}
\put(1.5,1){\RedNode}
\put(2.25,1){\color{Red}{\em \footnotesize 3}}
\put(0,6){\CyanNode}
\put(0,6.1){\color{Cyan}{\em \footnotesize 1}}
\put(0,0){\CyanNode}
\put(0.75,0){\color{Cyan}{\em \footnotesize 1}}
\thicklines
\put(3.375,2.275){\color{Black}\qbezier(0,0)(-0.6,0.42)(-1.2,0.84)}
\put(4.875,3.275){\color{Black}\qbezier(0,0)(-0.6,0.42)(-1.2,0.84)}
\put(3.375,4.275){\color{Black}\qbezier(0,0)(-0.6,0.42)(-1.2,0.84)}
\put(1.875,5.275){\color{Black}\qbezier(0,0)(-0.6,0.42)(-1.2,0.84)}
\put(6.375,4.275){\color{Black}\qbezier(0,0)(-0.6,0.42)(-1.2,0.84)}
\put(4.875,5.275){\color{Black}\qbezier(0,0)(-0.6,0.42)(-1.2,0.84)}
\put(3.375,6.275){\color{Black}\qbezier(0,0)(-0.6,0.42)(-1.2,0.84)}
\put(4.875,7.275){\color{Black}\qbezier(0,0)(-0.6,0.42)(-1.2,0.84)}
\put(0.625,0.275){\color{Black}\qbezier(0,0)(0.6,0.42)(1.2,0.84)}
\put(2.175,1.275){\color{Black}\qbezier(0,0)(0.6,0.42)(1.2,0.84)}
\put(3.675,2.275){\color{Black}\qbezier(0,0)(0.6,0.42)(1.2,0.84)}
\put(5.175,3.275){\color{Black}\qbezier(0,0)(0.6,0.42)(1.2,0.84)}
\put(3.675,4.275){\color{Black}\qbezier(0,0)(0.6,0.42)(1.2,0.84)}
\put(2.175,3.275){\color{Black}\qbezier(0,0)(0.6,0.42)(1.2,0.84)}
\put(2.175,5.275){\color{Black}\qbezier(0,0)(0.6,0.42)(1.2,0.84)}
\put(3.675,6.275){\color{Black}\qbezier(0,0)(0.6,0.42)(1.2,0.84)}
\put(0.625,6.275){\color{Black}\qbezier(0,0)(0.6,0.42)(1.2,0.84)}
\put(2.175,7.275){\color{Black}\qbezier(0,0)(0.6,0.42)(1.2,0.84)}
\put(3.675,8.275){\color{Black}\qbezier(0,0)(0.6,0.42)(1.2,0.84)}
\put(5.175,9.275){\color{Black}\qbezier(0,0)(0.6,0.42)(1.2,0.84)}
\end{picture}
}

\newpage
\begin{figure}
\begin{center}
\parbox{6in}{{\bf \ESixStackFigure.1}\ \ A skew-stack $\mystackedP = P_{\mytinyE_{6}}(\omega_{\mbox{\color{Red}{\tiny $1$}}}) \myposetstack_{\phi_{1}} P_{\mytinyE_{6}}(\omega_{\mbox{\color{Red}{\tiny $1$}}}) \myposetstack_{\phi_{2}} P_{\mytinyE_{6}}(\omega_{\mbox{\color{BurntOrange}{\tiny $6$}}})$ of the $\myE_{6}$-structured posets from \ESixLatticeFigure.2. 
The bond $\phi_{1}$ is the identity mapping $P_{\mytinyE_{6}}(\omega_{\mbox{\color{Red}{\tiny $1$}}}) \longrightarrow P_{\mytinyE_{6}}(\omega_{\mbox{\color{Red}{\tiny $1$}}})$, and the bond $\phi_{2}$ identifies the lower order ideal of $P_{\mytinyE_{6}}(\omega_{\mbox{\color{Red}{\tiny $1$}}})$ generated by its color {\color{BurntOrange}{$6$}} vertices and the upper order ideal of $P_{\mytinyE_{6}}(\omega_{\mbox{\color{BurntOrange}{\tiny $6$}}})$ generated by its 
color {\color{BurntOrange}{$6$}} vertices.}

\vspace*{0.1in} 
\setlength{\unitlength}{1cm}
\begin{picture}(14,20)
\put(0.5,18.5){\LARGE $\mylargestackedP$}
\thicklines
\put(1.25,18.75){\vector(3,-1){2}}
\put(0,0){\ESixFlat}
\put(1.25,7){\ESixFlat}
\put(7.5,9){\ESixFlatTwo}
\thinlines
\put(6.75,0.275){\color{Gray}\qbezier(0,0)(0.55,3.42)(1.1,6.84)}
\put(5.25,1.275){\color{Gray}\qbezier(0,0)(0.55,3.42)(1.1,6.84)}
\put(3.75,2.275){\color{Gray}\qbezier(0,0)(0.55,3.42)(1.1,6.84)}
\put(2.25,3.275){\color{Gray}\qbezier(0,0)(0.55,3.42)(1.1,6.84)}
\put(0.75,4.275){\color{Gray}\qbezier(0,0)(0.55,3.42)(1.1,6.84)}
\put(5.25,3.275){\color{Gray}\qbezier(0,0)(0.55,3.42)(1.1,6.84)}
\put(3.75,4.275){\color{Gray}\qbezier(0,0)(0.55,3.42)(1.1,6.84)}
\put(2.25,5.275){\color{Gray}\qbezier(0,0)(0.55,3.42)(1.1,6.84)}
\put(5.25,5.275){\color{Gray}\qbezier(0,0)(0.55,3.42)(1.1,6.84)}
\put(3.75,6.275){\color{Gray}\qbezier(0,0)(0.55,3.42)(1.1,6.84)}
\put(2.25,7.275){\color{Gray}\qbezier(0,0)(0.55,3.42)(1.1,6.84)}
\put(6.75,6.275){\color{Gray}\qbezier(0,0)(0.55,3.42)(1.1,6.84)}
\put(5.25,7.275){\color{Gray}\qbezier(0,0)(0.55,3.42)(1.1,6.84)}
\put(3.75,8.275){\color{Gray}\qbezier(0,0)(0.55,3.42)(1.1,6.84)}
\put(2.25,9.275){\color{Gray}\qbezier(0,0)(0.55,3.42)(1.1,6.84)}
\put(0.75,10.275){\color{Gray}\qbezier(0,0)(0.55,3.42)(1.1,6.84)}
\put(7.9,13.275){\color{Gray}\qbezier(0,0)(3.1,2.92)(6.2,5.84)}
\put(6.4,12.275){\color{Gray}\qbezier(0,0)(3.1,2.92)(6.2,5.84)}
\put(4.9,11.275){\color{Gray}\qbezier(0,0)(3.1,2.92)(6.2,5.84)}
\put(3.4,10.275){\color{Gray}\qbezier(0,0)(3.1,2.92)(6.2,5.84)}
\put(6.4,10.275){\color{Gray}\qbezier(0,0)(3.1,2.92)(6.2,5.84)}
\put(4.9,9.275){\color{Gray}\qbezier(0,0)(3.1,2.92)(6.2,5.84)}
\put(6.4,8.275){\color{Gray}\qbezier(0,0)(3.1,2.92)(6.2,5.84)}
\put(7.9,7.275){\color{Gray}\qbezier(0,0)(3.1,2.92)(6.2,5.84)}
\end{picture}
\end{center}
\end{figure}

\begin{figure}
\begin{center}
\parbox{6in}{{\bf \ESixStackFigure.2}\ \ The $\myE_{6}$-structured birds-eye scaffold $\scaffoldS^{\mybirdseye}$ of the  skew-stack $\mystackedP$ from \ESixStackFigure.1. 
The vertex-colored footprint of our scaffold $\scaffoldS^{\mybirdseye}$ is the birds-eye poset $\mystackedP^{\mybirdseye}$ with vertex colors as depicted below (and ignoring just for a moment the ordered pair assigned to each vertex). 
The ordered pair at any vertex $\overline{x}$ in $\mystackedP^{\mybirdseye}$ records the bounding pair $(b_{\barx}^{(\mbox{\tiny lo})},b_{\barx}^{(\mbox{\tiny hi})})$ for that vertex, and the totality of these is the set $\mathcal{B}^{\mybirdseye} := \left((b_{\barx}^{(\mbox{\tiny lo})},b_{\barx}^{(\mbox{\tiny hi})})\right)_{\barx \in \mytinystackedP^{\mysmallbirdseye}}$ of bounding pairs for our scaffold. 
There are no non-unit amplifiers for this scaffold, so $\mathcal{A}^{\mybirdseye} := (1)_{\barx \in \mytinystackedP^{\mysmallbirdseye}}$. 
Our birds-eye scaffold is, then, $\scaffoldS^{\mybirdseye} := (\mystackedP^{\mybirdseye},\mathcal{A}^{\mybirdseye},\mathcal{B}^{\mybirdseye},\vcolor^{\mybirdseye})$, with each of $\mystackedP^{\mybirdseye}$, $\mathcal{A}^{\mybirdseye}$, $\mathcal{B}^{\mybirdseye}$, and $\vcolor^{\mybirdseye}$ as described here.}   

\vspace*{0.1in} 
\ESixFlower

\vspace*{-0.1in} 
\setlength{\unitlength}{1cm}
\begin{picture}(5,16)
\put(0,4){\begin{picture}(5,12)
\put(-2.5,7.2){\LARGE $\bigscaffoldS^{\mylargebirdseye}$}
\put(-1.4,7.6){\vector(2,1){1.5}}
\put(0,10){\TypeEboxDot{Cyan}}
\put(0.7,10.3){\color{Cyan}{\em \footnotesize 1}}
\put(0,9.8){{\tiny $(0,2)$}}
\put(0,4){\TypeEboxDot{Cyan}}
\put(0.1,4.3){\color{Cyan}{\em \footnotesize 1}}
\put(0,3.8){{\tiny $(0,2)$}}
\put(1,9){\TypeEboxDot{Red}}
\put(1.7,9.3){\color{Red}{\em \footnotesize 3}}
\put(1,8.8){{\tiny $(0,2)$}}
\put(1,5){\TypeEboxDot{Red}}
\put(1.1,5.3){\color{Red}{\em \footnotesize 3}}
\put(1.7,5.2){{\tiny $(0,2)$}}
\put(2,8){\TypeEboxDot{Purple}}
\put(2.7,8.3){\color{Purple}{\em \footnotesize 4}}
\put(1.7,8.2){{\tiny $(0,2)$}}
\put(2,6){\TypeEboxDot{Purple}}
\put(2.1,6.15){\color{Purple}{\em \footnotesize 4}}
\put(2.7,6.2){{\tiny $(0,2)$}}
\put(2,4){\TypeEboxDot{Purple}}
\put(2.7,4.15){\color{Purple}{\em \footnotesize 4}}
\put(1.7,4.2){{\tiny $(0,3)$}}
\put(1,7){\TypeEboxDot{SpringGreen}}
\put(1.1,7.15){\color{SpringGreen}{\em \footnotesize 2}}
\put(1.7,7.2){{\tiny $(0,2)$}}
\put(3,7){\TypeEboxDot{OliveGreen}}
\put(3.7,7.3){\color{OliveGreen}{\em \footnotesize 5}}
\put(2.7,7.2){{\tiny $(0,2)$}}
\put(3,5){\TypeEboxDot{OliveGreen}}
\put(3.7,5){\color{OliveGreen}{\em \footnotesize 5}}
\put(2.7,5.2){{\tiny $(0,3)$}}
\put(4,6){\TypeEboxDot{BurntOrange}}
\put(4.7,6.3){\color{BurntOrange}{\em \footnotesize 6}}
\put(3.7,6.2){{\tiny $(0,3)$}}
\thicklines
\put(2.625,2.375){\color{Black}\qbezier(0,0)(0.375,0.375)(0.75,0.75)}
\put(1.625,3.375){\color{Black}\qbezier(0,0)(0.375,0.375)(0.75,0.75)}
\put(2.625,4.375){\color{Black}\qbezier(0,0)(0.375,0.375)(0.75,0.75)}
\put(3.625,5.375){\color{Black}\qbezier(0,0)(0.375,0.375)(0.75,0.75)}
\put(0.625,4.375){\color{Black}\qbezier(0,0)(0.375,0.375)(0.75,0.75)}
\put(1.625,5.375){\color{Black}\qbezier(0,0)(0.375,0.375)(0.75,0.75)}
\put(2.625,6.375){\color{Black}\qbezier(0,0)(0.375,0.375)(0.75,0.75)}
\put(1.625,7.375){\color{Black}\qbezier(0,0)(0.375,0.375)(0.75,0.75)}
\put(4.39,0.375){\color{Black}\qbezier(0,0)(-0.375,0.375)(-0.75,0.75)}
\put(3.39,1.375){\color{Black}\qbezier(0,0)(-0.375,0.375)(-0.75,0.75)}
\put(2.39,2.375){\color{Black}\qbezier(0,0)(-0.375,0.375)(-0.75,0.75)}
\put(1.39,3.375){\color{Black}\qbezier(0,0)(-0.375,0.375)(-0.75,0.75)}
\put(2.39,4.375){\color{Black}\qbezier(0,0)(-0.375,0.375)(-0.75,0.75)}
\put(3.39,3.375){\color{Black}\qbezier(0,0)(-0.375,0.375)(-0.75,0.75)}
\put(3.39,5.375){\color{Black}\qbezier(0,0)(-0.375,0.375)(-0.75,0.75)}
\put(2.39,6.375){\color{Black}\qbezier(0,0)(-0.375,0.375)(-0.75,0.75)}
\put(4.39,6.375){\color{Black}\qbezier(0,0)(-0.375,0.375)(-0.75,0.75)}
\put(3.39,7.375){\color{Black}\qbezier(0,0)(-0.375,0.375)(-0.75,0.75)}
\put(2.39,8.375){\color{Black}\qbezier(0,0)(-0.375,0.375)(-0.75,0.75)}
\put(1.39,9.375){\color{Black}\qbezier(0,0)(-0.375,0.375)(-0.75,0.75)}
\end{picture}
}
\put(0,0){\TypeEboxDot{Cyan}}
\put(0.7,0){\color{Cyan}{\em \footnotesize 1}}
\put(0,-0.2){{\tiny $(2,3)$}}
\put(0,6){\TypeEboxDot{Cyan}}
\put(0.1,6.3){\color{Cyan}{\em \footnotesize 1}}
\put(0,5.8){{\tiny $(2,3)$}}
\put(1,1){\TypeEboxDot{Red}}
\put(1.7,1){\color{Red}{\em \footnotesize 3}}
\put(0.7,1.2){{\tiny $(2,3)$}}
\put(1,5){\TypeEboxDot{Red}}
\put(1.1,5){\color{Red}{\em \footnotesize 3}}
\put(1.7,5.2){{\tiny $(2,3)$}}
\put(1,7){\TypeEboxDot{Red}}
\put(1.1,7.15){\color{Red}{\em \footnotesize 3}}
\put(1.7,7.2){{\tiny $(0,3)$}}
\put(2,2){\TypeEboxDot{Purple}}
\put(2.7,2){\color{Purple}{\em \footnotesize 4}}
\put(1.7,2.2){{\tiny $(2,3)$}}
\put(2,4){\TypeEboxDot{Purple}}
\put(2.1,4.15){\color{Purple}{\em \footnotesize 4}}
\put(2.7,4.2){{\tiny $(2,3)$}}
\put(2,6){\TypeEboxDot{Purple}}
\put(2.7,6.15){\color{Purple}{\em \footnotesize 4}}
\put(1.7,6.2){{\tiny $(0,3)$}}
\put(1,3){\TypeEboxDot{SpringGreen}}
\put(1.1,3.15){\color{SpringGreen}{\em \footnotesize 2}}
\put(1.7,3.2){{\tiny $(2,3)$}}
\put(3,7){\TypeEboxDot{SpringGreen}}
\put(3.75,7.15){\color{SpringGreen}{\em \footnotesize 2}}
\put(2.7,7.2){{\tiny $(0,3)$}}
\put(3,3){\TypeEboxDot{OliveGreen}}
\put(3.7,3){\color{OliveGreen}{\em \footnotesize 5}}
\put(2.7,3.2){{\tiny $(2,3)$}}
\put(3,5){\TypeEboxDot{OliveGreen}}
\put(3.7,5.3){\color{OliveGreen}{\em \footnotesize 5}}
\put(2.7,5.2){{\tiny $(0,3)$}}
\put(4,4){\TypeEboxDot{BurntOrange}}
\put(3.7,4.2){{\tiny $(0,3)$}}
\put(4.7,4){\color{BurntOrange}{\em \footnotesize 6}}
\thicklines
\put(2.625,8.125){\color{Black}\qbezier(0,0)(0.375,-0.375)(0.75,-0.75)}
\put(1.625,7.125){\color{Black}\qbezier(0,0)(0.375,-0.375)(0.75,-0.75)}
\put(2.625,6.125){\color{Black}\qbezier(0,0)(0.375,-0.375)(0.75,-0.75)}
\put(3.625,5.125){\color{Black}\qbezier(0,0)(0.375,-0.375)(0.75,-0.75)}
\put(0.625,6.125){\color{Black}\qbezier(0,0)(0.375,-0.375)(0.75,-0.75)}
\put(1.625,5.125){\color{Black}\qbezier(0,0)(0.375,-0.375)(0.75,-0.75)}
\put(2.625,4.125){\color{Black}\qbezier(0,0)(0.375,-0.375)(0.75,-0.75)}
\put(1.625,3.125){\color{Black}\qbezier(0,0)(0.375,-0.375)(0.75,-0.75)}
\put(4.375,10.125){\color{Black}\qbezier(0,0)(-0.375,-0.375)(-0.75,-0.75)}
\put(3.375,9.125){\color{Black}\qbezier(0,0)(-0.375,-0.375)(-0.75,-0.75)}
\put(2.375,8.125){\color{Black}\qbezier(0,0)(-0.375,-0.375)(-0.75,-0.75)}
\put(1.375,7.125){\color{Black}\qbezier(0,0)(-0.375,-0.375)(-0.75,-0.75)}
\put(2.375,6.125){\color{Black}\qbezier(0,0)(-0.375,-0.375)(-0.75,-0.75)}
\put(3.375,7.125){\color{Black}\qbezier(0,0)(-0.375,-0.375)(-0.75,-0.75)}
\put(3.375,5.125){\color{Black}\qbezier(0,0)(-0.375,-0.375)(-0.75,-0.75)}
\put(2.375,4.125){\color{Black}\qbezier(0,0)(-0.375,-0.375)(-0.75,-0.75)}
\put(4.375,4.125){\color{Black}\qbezier(0,0)(-0.375,-0.375)(-0.75,-0.75)}
\put(3.375,3.125){\color{Black}\qbezier(0,0)(-0.375,-0.375)(-0.75,-0.75)}
\put(2.375,2.125){\color{Black}\qbezier(0,0)(-0.375,-0.375)(-0.75,-0.75)}
\put(1.375,1.125){\color{Black}\qbezier(0,0)(-0.375,-0.375)(-0.75,-0.75)}
\end{picture}
\end{center}
\end{figure}

\clearpage
\noindent 
also a shortest path in $[0,x_{t}]_{\mathbb{Z}} \times [0,y_{t}]_{\mathbb{Z}}$. 
For any shortest path from $u'$ to $v'$ in $P_{t}$, say $+\eelt$ is added $a_{\eelt}^{+}$ times and $-\eelt$ is added $a_{\eelt}^{-}$ times and that $+\felt$ is added $a_{\felt}^{+}$ times and $-\felt$ is added $a_{\felt}^{-}$. 
Now, $v' = (u'_{1}+a_{\eelt}^{+}-a_{\eelt}^{-},u'_{2}+a_{\felt}^{+}-a_{\felt}^{-})$. 
Then $\psi(u')+\gelt = \psi(v')$ if and only if $(u'_{1}+\sum_{s=t+1}^{m}p_{s}+\epsilon_{1},u'_{2}+\sum_{s=t+1}^{m}q_{s}+\epsilon_{2}) = (u'_{1}+a_{\eelt}^{+}-a_{\eelt}^{-}+\sum_{s=t+1}^{m}p_{s},u'_{2}+a_{\felt}^{+}-a_{\felt}^{-}+\sum_{s=t+1}^{m}q_{s})$ if and only if $\epsilon_{1} = a_{\eelt}^{+}-a_{\eelt}^{-}$ and $\epsilon_{2} = a_{\felt}^{+}-a_{\felt}^{-}$. 
Then, $u' + \gelt = v'$. 
So, $u' \mylongarrow{\gelt} v'$ is a shortest path in $[0,x_{t}]_{\mathbb{Z}} \times [0,y_{t}]_{\mathbb{Z}}$, and since this is the only shortest path possible, it must therefore be a shortest path in the full-length sublattice $P_{t}$. 
That is, $u' \rightarrow v'$ in $P_{t}$, so by \StackLemma.3, we have $\overline{u} \rightarrow \overline{v}$ in $\mystackedP^{\mybirdseye}$.\hfill\QED

\noindent 
{\bf \StackDefsThree}\ \ Continue with the notation from \StackDefsOneTwo. 
The pair $(\mystackedP^{\mybirdseye},\vcolor^{\mybirdseye})$ is, by \StackLemma.1, a vertex-colored poset we call the {\em birds-eye poset of} $\mystackedP = P_{1} \myposetstack_{\phi_{1}} P_{2} \myposetstack_{\phi_{2}} \cdots \myposetstack_{\phi_{m-1}} P_{m}$.  
Assume now that for each $r \in [1,m]_{\mathbb{Z}}$, we have $\vcolor^{(r)}(x) \not= \vcolor^{(r)}(y)$ whenever $x \rightarrow y$ in $P_{r}$. 
In this case, consider $\overline{x}$ in $\mystackedP^{\mybirdseye}$, where $\overline{x}$ is the equivalence class of an $x$ in some $P_{r}$. 
Let $b_{\overline{x}}^{(\mbox{\tiny lo})} := p$ if $p$ is smallest such that there exists $x' \in P_{p+1}$ with $x' \sim x$, and let $b_{\overline{x}}^{(\mbox{\tiny hi})} := q$ if $q$ is largest such that there exists $x' \in P_{q}$ with $x' \sim x$. 
Our hypothesis that $\vcolor^{(r)}(x) \not= \vcolor^{(r)}(y)$ whenever $x \rightarrow y$ in $P_{r}$ for each $r$ guarantees that the bounding pair $(b_{\barx}^{(\mbox{\tiny lo})},b_{\barx}^{(\mbox{\tiny hi})})$ depends only on $\overline{x}$. 
Let $\mathcal{B}^{\mybirdseye} := \left((b_{\barx}^{(\mbox{\tiny lo})},b_{\barx}^{(\mbox{\tiny hi})})\right)_{\barx \in \mytinystackedP^{\mysmallbirdseye}}$. 
Set $a_{\overline{x}} := 1$, so our amplifier set is $\mathcal{A}^{\mybirdseye} := (1)_{\barx \in \mytinystackedP^{\mysmallbirdseye}}$.  
Then set $\scaffoldS^{\mybirdseye} := (\mystackedP^{\mybirdseye},\mathcal{A}^{\mybirdseye},\mathcal{B}^{\mybirdseye},\vcolor^{\mybirdseye})$, a {\em birds-eye scaffold} with footprint $\mystackedP^{\mybirdseye}$. 
Define $\varphi: \Jcolor(\mystackedP) \longrightarrow \Lcolor(\scaffoldS^{\mybirdseye})$ as follows: For a down-set $\telt$ from $\mystackedP$, say $\varphi(\telt) := (c_{\barx}(\telt))_{\barx \in \mytinystackedP^{\mysmallbirdseye}}$, where $c_{\barx}(\telt) := r$ if $r$ is largest such that $x' \sim x$, $x' \in P_{r}$, and $x' \in \telt$ and $c_{\barx}(\telt) := b_{\overline{x}}^{(\mbox{\tiny lo})}$ if no such $x'$ exists. 
\ESixStackFigure\ depicts a skew-stack $\mystackedP$ and its associated birds-eye scaffold $\scaffoldS^{\mybirdseye}$. 

The next result follows by a straightforward walk through the definitions. 

\noindent
{\bf \IsomLatticeTheorem}\ \ {\sl Keep the hypotheses and notation of the preceding definitions.  Then the mapping} $\varphi: \Jcolor(\mystackedP) \longrightarrow \Lcolor(\scaffoldS^{\mybirdseye})$ {\sl of \StackDefsThree\ is an isomorphism of diamond-colored distributive lattices, so that} 
\[\Jcolor(\mystackedP) \cong \Lcolor(\scaffoldS^{\mybirdseye}).\]
{\sl That is, the DCDL of order ideals from the given skew-stack is naturally isomorphic to the DCDL of ideal arrays over the associated birds-eye scaffold.}

\begin{center}
\underline{\hspace*{4in}}
\end{center}

\noindent
{\bf \S \GridPosetSection. Two-color grid posets.} 
The `semistandard' vertex-colored posets depicted, by way of examples, in \GridPosetsFigureList\ were discovered through the collaboration of a number of researchers \cite{ADLP}, \cite{ADLMPPW} and are instances of the `two-colored grid posets' we discuss in this section. 
That said, the general theory of two-colored grid posets presented in the aforementioned papers and summarized again below was largely developed by this author. 

{\bf [\S \GridPosetSection.1:\! Grid posets.]} 
Given a finite poset $(P,\leq_{_{P}})$, a {\em chain function for} $P$ is a function $\mychain: P \longrightarrow [1,m]_{\mathbb{Z}}$ for some positive integer $m$ such that (1) $\mychain^{-1}(i)$ is a (possibly empty) chain in $P$ for $1 \leq i \leq m$, and (2) given any covering relation $u \rightarrow v$ in $P$, it is the case that either $\mychain(u) = \mychain(v)$ or $\mychain(u) = \mychain(v) + 1$.  
A {\em grid poset} is a finite poset $(P,\leq_{_{P}})$ together with a chain function $\mychain: P \longrightarrow [1,m]_{\mathbb{Z}}$ for some positive integer $m$.  
Depending on context, the notation $P$ can refer to the grid poset $(P,\leq_{_{P}}, \mychain: P \rightarrow [1,m]_{\mathbb{Z}})$ or the underlying poset $(P,\leq_{_{P}})$. 
In the covering digraph for $P$, each element covers at most two elements and is covered by at most two elements. 
Observe that if $i$ is the smallest (respectively largest) integer such that $\mychain^{-1}(i)$ is nonempty and  if $u$ is the maximal (respectively minimal) element of $\mychain^{-1}(i)$, then $u$ is a maximal (respectively minimal) element of the poset $P$.  
If $P$ is nonempty, the {\em rightmost maximal element} is the maximal element $z$ with the property that $\mychain(z) \geq \mychain(u)$ for all maximal elements $u$ in $P$. A {\em face vertex} of the grid poset $P$ is an element $u$ of $P$  such that whenever $u \rightarrow v$, then $\mychain(u) = \mychain(v)$. 
A grid poset $P$ is {\em connected} if and only if the covering digraph for the poset $P$ is connected.  
For $i = 1,2$, let $P_{i}$ be a grid poset with chain function $\mychain_{i}: P_{i} \longrightarrow [1,m_{i}]_{\mathbb{Z}}$ for some positive integer $m_{i}$.  
A one-to-one correspondence $\phi: P_{1} \longrightarrow P_{2}$ is an {\em isomorphism of grid posets} if we have $u \rightarrow v$ in $P_{1}$ with $\mychain_{1}(u) = \mychain_{1}(v)$ (respectively $\mychain_{1}(u) = \mychain_{1}(v) + 1$) if and only if  $\phi(u) \rightarrow \phi(v)$ in $P_{2}$ with $\mychain_{2}(\phi(u)) = \mychain_{2}(\phi(v))$ (respectively $\mychain_{2}(\phi(u)) = \mychain_{2}(\phi(v)) + 1$).  
In this case, say that $P_{1}$ and $P_{2}$ are {\em isomorphic grid posets}. 
For a nonempty grid poset $(P,\leq_{_{P}},\mychain: P \rightarrow [1,m]_{\mathbb{Z}})$, there exists a surjective chain function $\mychain':P \longrightarrow [1,m']_{\mathbb{Z}}$ such that the grid poset $P$ is isomorphic to $(P, {\leq_{_{P}},} \mychain': P \rightarrow [1,m']_{\mathbb{Z}})$; if $P$ is connected, then this surjective chain function $\mychain'$ is unique. 
We say $Q$ is a {\em grid subposet} of a given grid poset $P$ if (1) $Q$ is a subposet of $P$ in the induced order, and (2) whenever $u \rightarrow v$ is a covering relation in $Q$ then it is also a covering relation in $P$.  In this case, we regard $Q$ with the chain function $\mychain|_{Q}$ to be a grid poset on its own.  
For $1 \leq i \leq m$ we set $\mathcal{C}_{i} := \mychain^{-1}(i)$; when we depict grid posets, the chains $\mathcal{C}_{i}$ will be directed from SW to NE.  
The following are the six non-isomorphic connected grid posets with three elements; notice that for four of these, the underlying poset is a chain. 

\begin{figure}[ht]
\hspace*{0.7cm}
\setlength{\unitlength}{0.75cm}
\begin{picture}(3.5,4)
\put(0,0.75){
\begin{picture}(3,3)
\put(2,0){\circle*{0.15}} 
\put(1,1){\circle*{0.15}} 
\put(0,2){\circle*{0.15}} 
\put(2,0){\line(-1,1){2}} 
\put(2.2,0.2){\scriptsize $\mathcal{C}_{3}$}
\put(1.2,1.2){\scriptsize $\mathcal{C}_{2}$}
\put(0.2,2.2){\scriptsize $\mathcal{C}_{1}$}
\end{picture}
}
\end{picture}
\setlength{\unitlength}{0.75cm}
\begin{picture}(2.5,4)
\put(0,0.75){
\begin{picture}(2.5,3)
\put(0,0){\circle*{0.15}} 
\put(1,1){\circle*{0.15}} 
\put(0,2){\circle*{0.15}} 
\put(0,0){\line(1,1){1}}
\put(1,1){\line(-1,1){1}} 
\put(1.2,1.2){\scriptsize $\mathcal{C}_{2}$}
\put(0.2,2.2){\scriptsize $\mathcal{C}_{1}$}
\end{picture}
}
\end{picture}
\setlength{\unitlength}{0.75cm}
\begin{picture}(2.5,4)
\put(0,0.75){
\begin{picture}(2.5,3)
\put(1,0){\circle*{0.15}} 
\put(0,1){\circle*{0.15}} 
\put(1,2){\circle*{0.15}} 
\put(1,0){\line(-1,1){1}}
\put(0,1){\line(1,1){1}} 
\put(1.2,0.2){\scriptsize $\mathcal{C}_{2}$}
\put(1.2,2.2){\scriptsize $\mathcal{C}_{1}$}
\end{picture}
}
\end{picture}
\setlength{\unitlength}{0.75cm}
\begin{picture}(3,4)
\put(0,0.75){
\begin{picture}(3,3)
\put(0,0){\circle*{0.15}} 
\put(1,1){\circle*{0.15}} 
\put(2,2){\circle*{0.15}} 
\put(0,0){\line(1,1){2}}
\put(2.2,2.2){\scriptsize $\mathcal{C}_{1}$}
\end{picture}
}
\end{picture}
\setlength{\unitlength}{0.75cm}
\begin{picture}(3.5,4)
\put(0,1){
\begin{picture}(3,3)
\put(0,0){\circle*{0.15}} 
\put(1,1){\circle*{0.15}} 
\put(2,0){\circle*{0.15}} 
\put(0,0){\line(1,1){1}}
\put(1,1){\line(1,-1){1}}
\put(2.2,0.2){\scriptsize $\mathcal{C}_{2}$}
\put(1.2,1.2){\scriptsize $\mathcal{C}_{1}$}
\end{picture}
}
\end{picture}
\setlength{\unitlength}{0.75cm}
\begin{picture}(3.5,4)
\put(0,1){
\begin{picture}(3,3)
\put(0,1){\circle*{0.15}} 
\put(1,0){\circle*{0.15}} 
\put(2,1){\circle*{0.15}} 
\put(0,1){\line(1,-1){1}}
\put(1,0){\line(1,1){1}}
\put(2.2,1.2){\scriptsize $\mathcal{C}_{2}$}
\put(0.2,1.2){\scriptsize $\mathcal{C}_{1}$}
\end{picture}
}
\end{picture}
\end{figure}

\noindent
Any grid poset  $(P,\leq_{_{P}}, \mychain: P \rightarrow [1,m]_{\mathbb{Z}})$ has a natural and convenient total ordering $\mathcal{T}_{P}$ of its elements given by the following rule:  
For distinct $u$ and $v$ in $P$ write $u <_{_{\mathcal{T}_{P}}} v$ if and only if (1) $\mychain(u) < \mychain(v)$ or (2) $\mychain(u) = \mychain(v)$ with $v <_{_{P}} u$. 
Let $l := |P|$. Number the vertices of $P$ $v_{1}, v_{2},\ldots, v_{l}$ so that $v_{p} <_{_{\mathcal{T}_{P}}} v_{q}$ whenever $1 \leq p < q \leq l$. 
The numberings of the vertices of the grid posets $P$ depicted in \PosetAndLatticeFig\ and in  \GridPosetsTotalOrderList\ follow $\mathcal{T}_{P}$. 

{\bf [\S \GridPosetSection.2:\! Two-color grid posets: definitions and examples.]}
A {\em two-color function} for a grid poset $(P,\leq_{_{P}}, \mychain: P \rightarrow [1,m]_{\mathbb{Z}})$ is a function $\mycolor: P \longrightarrow \Delta$ such that (1) $|\Delta| = 2$, (2) $\mycolor(u) = \mycolor(v)$ if $\mychain(u) = \mychain(v)$, and (3) if $u$ and $v$ are in the same connected component of $P$ with $\mychain(u) = \mychain(v)+1$, then $\mycolor(u) \not= \mycolor(v)$. 
A {\em two-color grid poset} is a grid poset $(P, \leq_{_{P}}, \mychain: P \rightarrow [1,m]_{\mathbb{Z}})$ together with a two-color function $\mycolor: P \longrightarrow \Delta$; in some contexts we will use the notation $P$ to refer to the two-color grid poset $(P,\leq_{_{P}}, \mychain: P \rightarrow [1,m]_{\mathbb{Z}}, \mycolor: P \rightarrow \Delta)$. 
A two-color grid poset should be thought of as a certain kind of vertex-colored poset. 
So, associated to any two-color grid poset $P$ is the DCDL $L := \Jcolor(P)$.  
If $Q$ is a grid subposet of the two-color grid poset $P$, then $Q$ is itself a two-color grid poset with chain function $\mychain|_{Q}$ and two-color function $\mbox{\sffamily color}|_{Q}$; in this case we call $Q$ a {\em two-color grid subposet of} $P$.  
Two two-color grid posets $(P_{i},\leq_{_{P_{i}}}, \mychain_{i}: P_{i} \rightarrow [1,m_{i}]_{\mathbb{Z}}, \mycolor_{i}: P_{i} \rightarrow \Delta)$ for $i = 1,2$ are {\em isomorphic} if there is an isomorphism $\phi: P_{1}  \longrightarrow P_{2}$ of grid posets such that $\mycolor_{2}(\phi(u)) = \mycolor_{1}(u)$ for all $u$ in $P_{1}$.  
We will usually take $\Delta := \{\alpha,\beta\}$.  
When we {\em switch} (or {\em reverse}) the vertex colors of $P$ we replace the color function $\mycolor$ with the color function $\mycolor': P \longrightarrow \{\alpha, \beta\}$ given by: $\mycolor'(v) = \alpha$ if $\mycolor(v) = \beta$, and $\mycolor'(v) = \beta$ if $\mycolor(v) = \alpha$.  
Similarly, one can {\em switch} (or {\em reverse}) the edge colors of $L$. 

We say a two-color grid poset $P$  has the {\em max property} if $P$ is isomorphic to a two-color grid poset $(Q,\leq_{_{Q}}, \mychain: Q \rightarrow [1,m]_{\mathbb{Z}}, \mycolor: Q \rightarrow \Delta)$ with a surjective chain function such that (1) if $u$ is any maximal element in the poset $Q$, then $\mychain(u) \leq 2$, and (2) if $v \not= u$ is another maximal element in $Q$, then $\mycolor(u) \not= \mycolor(v)$. 
In \GridPosetsFigureList\ we depict eight two-color grid posets with the max property. 
The vertex-colored poset $P$ of \PosetAndLatticeFig\ is also a two-color grid poset with the max property; the lattice $L$ in that figure is $\Jcolor(P)$, and its corresponding down-sets are depicted in \OrderIdealFig. 

{\bf [\S \GridPosetSection.3:\! Decomposable/indecomposable two-color grid posets.]}
Let  $P$ be a grid poset with chain function $\mychain: P \longrightarrow [1,m]_{\mathbb{Z}}$.  
Suppose $P_{1}$ is a nonempty down-set and is a proper subset of $P$.  
Regard $P_{1}$ and $P_{2} := P \setminus P_{1}$ to be subposets of the poset $P$ in the induced order.  
Suppose that whenever $u$ is a maximal (respectively minimal) element of $P_{1}$ and $v$ is a maximal (respectively minimal) element of $P_{2}$, then $\mychain(u) \leq \mychain(v)$.  
Then we say that $P$ {\em decomposes into} $P_{1} \triangleleft P_{2}$, and we write $P = P_{1} \triangleleft P_{2}$. 
If no such down-set $P_{1}$ exists, then we say the grid poset $P$ is {\em indecomposable}. 
Note that if $P = P_{1} \triangleleft P_{2}$ and $u < v$ in $P$ with $u \in P_{2}$, then $v \in P_{2}$. 
Moreover, if $u \rightarrow v$ in $P$ with $u \in P_{1}$ and $v \in P_{2}$, then $\mychain(u) = \mychain(v)$. 
Also, if $u \rightarrow v$ is a covering relation in the poset $P_{i}$ for $i \in \{1,2\}$, note that $u \rightarrow v$ is also a covering relation in $P$; hence each $P_{i}$ is a grid subposet of $P$.  
If $P$ is a grid poset that decomposes into $P_{1} \triangleleft Q$, and if $Q$ decomposes into $P_{2} \triangleleft P_{3}$, then  $P = P_{1} \triangleleft (P_{2} \triangleleft P_{3})$.  
But now observe that $P = (P_{1} \triangleleft P_{2}) \triangleleft P_{3}$.  
So we may write $P = P_{1} \triangleleft P_{2} \triangleleft P_{3}$ unambiguously.  
In general, if $P = P_{1} \triangleleft P_{2} \triangleleft \cdots \triangleleft P_{k}$, then each $P_{i}$ with chain function $\mychain|_{P_{i}}$ is a grid subposet of $P$. 
Also,  a down-set $\selt$ taken from $P$ may be expressed as the disjoint union $(\selt \cap P_{1}) \cup (\selt \cap P_{2}) \cup \cdots \cup (\selt \cap P_{k})$, where each $\selt \cap P_{i}$ is a down-set taken from $P_{i}$. 
If in addition $P$ is a \dichromatic grid poset with two-color function $\mycolor$, then each $P_{i}$ with chain function $\mychain|_{P_{i}}$ and two-color function $\mycolor|_{P_{i}}$ is a \dichromatic grid subposet of $P$, and so $P_{1} \triangleleft P_{2} \triangleleft \cdots \triangleleft P_{k}$ is a decomposition of $P$ into \dichromatic grid posets.

Let $L := \Jcolor(P)$ for some given two-color grid poset $(P,\leq_{_{P}}, \mychain: P \rightarrow [1,m]_{\mathbb{Z}}, \mycolor: P \rightarrow \{\alpha,\beta\})$. 
Throughout what follows, we regard $\alpha$ to be the `1$^{\mbox{\tiny st}}$' color and $\beta$ to be the `2$^{\mbox{\tiny nd}}$' color. 
Then, for each $\xelt \in L$, our weight function is given by $wt(\xelt) = \big(\rho_{\alpha}(\xelt)-\delta_{\alpha}(\xelt),\rho_{\beta}(\xelt)-\delta_{\beta}(\xelt)\big)$, where, for $\gamma \in \{\alpha,\beta\}$, $\rho_{\gamma}$ and $\delta_{\gamma}$ are the rank and depth functions respectively for the color $\gamma$ components of $L$. 
If $P = P_{1} \triangleleft P_{2} \triangleleft \cdots \triangleleft P_{k}$, then for each $i$ we let $L_{i} := \Jcolor(P_{i})$ be the DCDL for the \dichromatic grid subposet $P_{i}$ of $P$. 
We let $wt_{i}$ denote the weight function for $L_{i}$; we let $\rho^{(i)}_{\alpha}$ and $\delta^{(i)}_{\alpha}$ (respectively $\rho^{(i)}_{\beta}$ and $\delta^{(i)}_{\beta}$) denote the rank and depth functions for color $\alpha$ (respectively color $\beta$).  
A version of the following appeared as Lemma 3.1 in \cite{ADLMPPW}. 

\noindent 
{\bf \WeightsLemma}\ \ {\sl Let $P$ be a \dichromatic grid poset as in the preceding paragraph, and suppose $P$ decomposes into $P = P_{1} \triangleleft P_{2} \triangleleft \cdots \triangleleft P_{k}$. 
Let $\xelt$ be an element of $L = \Jcolor(P)$, and let $\gamma \in \{\alpha,\beta\}$.  
Then: $\rho_{\gamma}(\xelt) = \sum_{i=1}^{k}\rho^{(i)}_{\gamma}(\xelt \cap P_{i})$,  $\delta_{\gamma}(\xelt) = \sum_{i=1}^{k}\delta^{(i)}_{\gamma}(\xelt \cap P_{i})$,  and $wt(\xelt) = \sum_{i=1}^{k}wt_{i}(\selt \cap P_{i})$.}

\newcommand{\ATwoFundamentalAA}{
\setlength{\unitlength}{1cm}
\begin{picture}(2,2)
\put(1,0){\circle*{0.15}} 
\put(0.4,-0.1){\footnotesize $u_{2}$}
\put(1.2,-0.1){\footnotesize $\beta$}
\put(0,1){\circle*{0.15}} 
\put(-0.6,0.9){\footnotesize $u_{1}$}
\put(0.2,0.9){\footnotesize $\alpha$}
\put(1,0){\line(-1,1){1}} 
\put(1.5,0.5){$\mathcal{C}_{2}$}
\put(0.5,1.5){$\mathcal{C}_{1}$}
\end{picture}
}

\newcommand{\ATwoFundamentalA}{
\setlength{\unitlength}{1cm}
\begin{picture}(2,3)
\put(0,2.75){\TechPropAalphaFigure}
\put(0,0.25){\ATwoFundamentalAA}
\end{picture} 
}

\newcommand{\ATwoFundamentalBB}{
\setlength{\unitlength}{1cm}
\begin{picture}(2,2)
\put(1,0){\circle*{0.15}} 
\put(0.4,-0.1){\footnotesize $u_{2}$}
\put(1.2,-0.1){\footnotesize $\alpha$}
\put(0,1){\circle*{0.15}} 
\put(-0.6,0.9){\footnotesize $u_{1}$}
\put(0.2,0.9){\footnotesize $\beta$}
\put(1,0){\line(-1,1){1}} 
\put(1.5,0.5){$\mathcal{C}_{2}$}
\put(0.5,1.5){$\mathcal{C}_{1}$}
\end{picture}
}

\newcommand{\ATwoFundamentalB}{
\setlength{\unitlength}{1cm}
\begin{picture}(2,3)
\put(0,2.75){\TechPropAbetaFigure}
\put(0,0.25){\ATwoFundamentalBB}
\end{picture} 
}

\newcommand{\BTwoFundamentalAA}{
\setlength{\unitlength}{1cm}
\begin{picture}(3,3)
\put(2,0){\circle*{0.15}} 
\put(1.4,-0.1){\footnotesize $u_{3}$}
\put(2.2,-0.1){\footnotesize $\alpha$}
\put(1,1){\circle*{0.15}} 
\put(0.4,0.9){\footnotesize $u_{2}$}
\put(1.2,0.9){\footnotesize $\beta$}
\put(0,2){\circle*{0.15}} 
\put(-0.6,1.9){\footnotesize $u_{1}$}
\put(0.2,1.9){\footnotesize $\alpha$}
\put(2,0){\line(-1,1){2}} 
\put(2.5,0.5){$\mathcal{C}_{3}$}
\put(1.5,1.5){$\mathcal{C}_{2}$}
\put(0.5,2.5){$\mathcal{C}_{1}$}
\end{picture}
}

\newcommand{\BTwoFundamentalA}{
\setlength{\unitlength}{1cm}
\begin{picture}(3,4)
\put(0,4){\TechPropBalphaFigure}
\put(0,0.5){\BTwoFundamentalAA}
\end{picture} 
}

\newcommand{\BTwoFundamentalBB}{
\setlength{\unitlength}{1cm}
\begin{picture}(3,4)
\put(1,0){\circle*{0.15}} 
\put(0.4,-0.1){\footnotesize $u_{4}$}
\put(1.2,-0.1){\footnotesize $\beta$}
\put(0,1){\circle*{0.15}} 
\put(-0.6,0.9){\footnotesize $u_{3}$}
\put(0.2,0.9){\footnotesize $\alpha$}
\put(1,2){\circle*{0.15}} 
\put(0.4,1.9){\footnotesize $u_{2}$}
\put(1.2,1.9){\footnotesize $\alpha$}
\put(0,3){\circle*{0.15}} 
\put(-0.6,2.9){\footnotesize $u_{1}$}
\put(0.2,2.9){\footnotesize $\beta$}
\put(0,1){\line(1,1){1}} 
\put(1,0){\line(-1,1){1}} 
\put(1,2){\line(-1,1){1}} 
\put(1.5,0.5){$\mathcal{C}_{3}$}
\put(1.5,2.5){$\mathcal{C}_{2}$}
\put(0.5,3.5){$\mathcal{C}_{1}$}
\end{picture}
}

\newcommand{\BTwoFundamentalB}{
\setlength{\unitlength}{1cm}
\begin{picture}(3,4.75)
\put(0,4.5){\TechPropBbetaFigure}
\put(0,0){\BTwoFundamentalBB}
\end{picture} 
}

\newcommand{\GTwoFundamentalAA}{
\setlength{\unitlength}{1cm}
\begin{picture}(4,6)
\put(3,0){\circle*{0.15}} 
\put(2.4,-0.1){\footnotesize $u_{6}$}
\put(3.2,-0.1){\footnotesize $\alpha$}
\put(2,1){\circle*{0.15}} 
\put(1.4,0.9){\footnotesize $u_{5}$}
\put(2.2,0.9){\footnotesize $\beta$}
\put(1,2){\circle*{0.15}} 
\put(0.4,1.9){\footnotesize $u_{4}$}
\put(1.2,1.9){\footnotesize $\alpha$}
\put(2,3){\circle*{0.15}} 
\put(1.4,2.9){\footnotesize $u_{3}$}
\put(2.2,2.9){\footnotesize $\alpha$}
\put(1,4){\circle*{0.15}} 
\put(0.4,3.9){\footnotesize $u_{2}$}
\put(1.2,3.9){\footnotesize $\beta$}
\put(0,5){\circle*{0.15}} 
\put(-0.6,4.9){\footnotesize $u_{1}$}
\put(0.2,4.9){\footnotesize $\alpha$}
\put(1,2){\line(1,1){1}} 
\put(3,0){\line(-1,1){2}} 
\put(2,3){\line(-1,1){2}} 
\put(3.5,0.5){$\mathcal{C}_{5}$}
\put(2.5,1.5){$\mathcal{C}_{4}$}
\put(2.5,3.5){$\mathcal{C}_{3}$}
\put(1.5,4.5){$\mathcal{C}_{2}$}
\put(0.5,5.5){$\mathcal{C}_{1}$}
\end{picture}
}

\newcommand{\GTwoFundamentalA}{
\setlength{\unitlength}{1cm}
\begin{picture}(4,7)
\put(0,7.25){\TechPropGalphaFigure}
\put(0,1){\GTwoFundamentalAA}
\end{picture} 
}

\newcommand{\GTwoFundamentalBB}{
\setlength{\unitlength}{1cm}
\begin{picture}(4,8)
\put(2,0){\circle*{0.15}} 
\put(1.4,-0.1){\footnotesize $u_{10}$}
\put(2.2,-0.15){\footnotesize $\beta$}
\put(1,1){\circle*{0.15}} 
\put(0.4,0.9){\footnotesize $u_{9}$}
\put(1.2,0.9){\footnotesize $\alpha$}
\put(2,2){\circle*{0.15}} 
\put(1.4,1.9){\footnotesize $u_{8}$}
\put(2.2,1.9){\footnotesize $\alpha$}
\put(1,3){\circle*{0.15}} 
\put(0.4,2.9){\footnotesize $u_{6}$}
\put(1.2,2.9){\footnotesize $\beta$}
\put(3,3){\circle*{0.15}} 
\put(2.4,2.9){\footnotesize $u_{7}$}
\put(3.2,2.9){\footnotesize $\alpha$}
\put(0,4){\circle*{0.15}} 
\put(-0.6,3.9){\footnotesize $u_{4}$}
\put(0.2,3.9){\footnotesize $\alpha$}
\put(2,4){\circle*{0.15}} 
\put(1.4,3.9){\footnotesize $u_{5}$}
\put(2.2,3.9){\footnotesize $\beta$}
\put(1,5){\circle*{0.15}} 
\put(0.4,4.9){\footnotesize $u_{3}$}
\put(1.2,4.9){\footnotesize $\alpha$}
\put(2,6){\circle*{0.15}} 
\put(1.4,5.9){\footnotesize $u_{2}$}
\put(2.2,5.9){\footnotesize $\alpha$}
\put(1,7){\circle*{0.15}} 
\put(0.4,6.9){\footnotesize $u_{1}$}
\put(1.2,6.9){\footnotesize $\beta$}
\put(1,1){\line(1,1){2}} 
\put(1,3){\line(1,1){1}} 
\put(0,4){\line(1,1){2}} 
\put(2,0){\line(-1,1){1}} 
\put(2,2){\line(-1,1){2}} 
\put(3,3){\line(-1,1){2}} 
\put(2,6){\line(-1,1){1}} 
\put(2.5,0.5){$\mathcal{C}_{5}$}
\put(3.5,3.5){$\mathcal{C}_{4}$}
\put(2.5,4.5){$\mathcal{C}_{3}$}
\put(2.5,6.5){$\mathcal{C}_{2}$}
\put(1.5,7.5){$\mathcal{C}_{1}$}
\end{picture} 
}

\newcommand{\GTwoFundamentalB}{
\setlength{\unitlength}{1cm}
\begin{picture}(4,8)
\put(-1,6){\TechPropGbetaFigure}
\put(0,0){\GTwoFundamentalBB}
\end{picture} 
}

\begin{figure}[ht]
\begin{center}
\GridPosets: Depicted below are four \dichromatic grid posets each possessing the max property.\\  
{\small (Each is the `$\mathscr{G}$-semistandard poset $P_{\mathscr{G}}^{\beta\alpha}(2,2)$' for the indicated integral Coxeter--Dynkin posy $\mathscr{G}$.)} 

\setlength{\unitlength}{1cm}
\begin{picture}(6,7)
\put(0,7){\fbox{$\mathscr{G} = \myA_{1} \oplus \myA_{1}$}}
\put(0.5,2.75){
\begin{picture}(3,3.5)
\put(3,2){\circle*{0.15}} 
\put(2.4,1.9){\footnotesize $v_{4}$}
\put(3.2,1.9){\footnotesize $\alpha$}
\put(4,3){\circle*{0.15}} 
\put(3.4,2.9){\footnotesize $v_{3}$}
\put(4.2,2.9){\footnotesize $\alpha$}
\put(0,2){\circle*{0.15}} 
\put(-0.6,1.9){\footnotesize $v_{2}$}
\put(0.2,1.9){\footnotesize $\beta$}
\put(1,3){\circle*{0.15}} 
\put(0.4,2.9){\footnotesize $v_{1}$}
\put(1.2,2.9){\footnotesize $\beta$}
\put(3,2){\line(1,1){1}}
\put(0,2){\line(1,1){1}} 
\put(1.5,3.5){$\mathcal{C}_{1}$}
\put(4.5,3.5){$\mathcal{C}_{2}$}
\put(1.8,2.375){$\bigoplus$}
\end{picture}
}
\put(0,2){\fbox{$\mathscr{G} = \myA_{2}$}}
\put(0,-2.5){
\begin{picture}(3,3.5)
\put(4,2){\circle*{0.15}} 
\put(3.4,1.9){\footnotesize $v_{8}$}
\put(4.2,1.9){\footnotesize $\beta$}
\put(5,3){\circle*{0.15}} 
\put(4.4,2.9){\footnotesize $v_{7}$}
\put(5.2,2.9){\footnotesize $\beta$}
\put(1,1){\circle*{0.15}} 
\put(0.4,0.9){\footnotesize $v_{6}$}
\put(1.2,0.9){\footnotesize $\alpha$}
\put(2,2){\circle*{0.15}} 
\put(1.4,1.9){\footnotesize $v_{5}$}
\put(2.2,1.9){\footnotesize $\alpha$}
\put(3,3){\circle*{0.15}} 
\put(2.4,2.9){\footnotesize $v_{4}$}
\put(3.2,2.9){\footnotesize $\alpha$}
\put(4,4){\circle*{0.15}} 
\put(3.4,3.9){\footnotesize $v_{3}$}
\put(4.2,3.9){\footnotesize $\alpha$}
\put(0,2){\circle*{0.15}} 
\put(-0.6,1.9){\footnotesize $v_{2}$}
\put(0.2,1.9){\footnotesize $\beta$}
\put(1,3){\circle*{0.15}} 
\put(0.4,2.9){\footnotesize $v_{1}$}
\put(1.2,2.9){\footnotesize $\beta$}
\put(1,1){\line(1,1){3}}
\put(4,2){\line(1,1){1}}
\put(0,2){\line(1,1){1}} 
\put(1,1){\line(-1,1){1}} 
\put(2,2){\line(-1,1){1}} 
\put(3,3){\line(1,-1){1}}
\put(4,4){\line(1,-1){1}}
\put(1.5,3.5){$\mathcal{C}_{1}$}
\put(4.5,4.5){$\mathcal{C}_{2}$}
\put(5.5,3.5){$\mathcal{C}_{3}$}
\end{picture}
}
\end{picture}
\hspace*{1.5cm}
\setlength{\unitlength}{1cm}
\begin{picture}(8,8)
\put(0,6){\fbox{$\mathscr{G} = \myC_{2}$}}
\put(0,0.25){
\begin{picture}(3,3.5)
\put(6,3){\circle*{0.15}} 
\put(5.4,2.9){\footnotesize $v_{14}$}
\put(6.2,2.9){\footnotesize $\alpha$}
\put(7,4){\circle*{0.15}} 
\put(6.4,3.9){\footnotesize $v_{13}$}
\put(7.2,3.9){\footnotesize $\alpha$}
\put(1,0){\circle*{0.15}} 
\put(0.4,-0.1){\footnotesize $v_{12}$}
\put(1.2,-0.1){\footnotesize $\beta$}
\put(3,2){\circle*{0.15}} 
\put(2.4,1.9){\footnotesize $v_{11}$}
\put(3.2,1.9){\footnotesize $\beta$}
\put(5,4){\circle*{0.15}} 
\put(4.4,3.9){\footnotesize $v_{10}$}
\put(5.2,3.9){\footnotesize $\beta$}
\put(6,5){\circle*{0.15}} 
\put(5.4,4.9){\footnotesize $v_{9}$}
\put(6.2,4.9){\footnotesize $\beta$}
\put(0,1){\circle*{0.15}} 
\put(-0.6,0.9){\footnotesize $v_{8}$}
\put(0.2,0.9){\footnotesize $\alpha$}
\put(1,2){\circle*{0.15}} 
\put(0.4,1.9){\footnotesize $v_{7}$}
\put(1.2,1.9){\footnotesize $\alpha$}
\put(2,3){\circle*{0.15}} 
\put(1.4,2.9){\footnotesize $v_{6}$}
\put(2.2,2.9){\footnotesize $\alpha$}
\put(3,4){\circle*{0.15}} 
\put(2.4,3.9){\footnotesize $v_{5}$}
\put(3.2,3.9){\footnotesize $\alpha$}
\put(4,5){\circle*{0.15}} 
\put(3.4,4.9){\footnotesize $v_{4}$}
\put(4.2,4.9){\footnotesize $\alpha$}
\put(5,6){\circle*{0.15}} 
\put(4.4,5.9){\footnotesize $v_{3}$}
\put(5.2,5.9){\footnotesize $\alpha$}
\put(0,3){\circle*{0.15}} 
\put(-0.6,2.9){\footnotesize $v_{2}$}
\put(0.2,2.9){\footnotesize $\beta$}
\put(2,5){\circle*{0.15}} 
\put(1.4,4.9){\footnotesize $v_{1}$}
\put(2.2,4.9){\footnotesize $\beta$}
\put(0,1){\line(1,1){5}}
\put(1,0){\line(1,1){5}}
\put(0,3){\line(1,1){2}} 
\put(6,3){\line(1,1){1}} 
\put(1,0){\line(-1,1){1}} 
\put(1,2){\line(-1,1){1}} 
\put(3,2){\line(-1,1){1}} 
\put(3,4){\line(-1,1){1}} 
\put(4,5){\line(1,-1){2}}
\put(5,6){\line(1,-1){2}}
\put(2.5,5.5){$\mathcal{C}_{1}$}
\put(5.5,6.5){$\mathcal{C}_{2}$}
\put(6.5,5.5){$\mathcal{C}_{3}$}
\put(7.5,4.5){$\mathcal{C}_{4}$}
\end{picture}
}
\end{picture}

\setlength{\unitlength}{1cm}
\begin{picture}(11,13)
\put(-0.5,9){\fbox{$\mathscr{G} = \myG_{2}$}}
\put(2,0){\circle*{0.15}} 
\put(1.4,-0.1){\footnotesize $v_{30}$}
\put(2.2,-0.15){\footnotesize $\beta$}
\put(1,1){\circle*{0.15}} 
\put(0.4,0.9){\footnotesize $v_{26}$}
\put(1.2,0.9){\footnotesize $\alpha$}
\put(2,2){\circle*{0.15}} 
\put(1.4,1.9){\footnotesize $v_{25}$}
\put(2.2,1.9){\footnotesize $\alpha$}
\put(1,3){\circle*{0.15}} 
\put(0.4,2.9){\footnotesize $v_{16}$}
\put(1.2,2.9){\footnotesize $\beta$}
\put(3,3){\circle*{0.15}} 
\put(2.4,2.9){\footnotesize $v_{24}$}
\put(3.2,2.9){\footnotesize $\alpha$}
\put(5,3){\circle*{0.15}} 
\put(4.4,2.9){\footnotesize $v_{29}$}
\put(5.2,2.85){\footnotesize $\beta$}
\put(0,4){\circle*{0.15}} 
\put(-0.6,3.9){\footnotesize $v_{10}$}
\put(0.2,3.9){\footnotesize $\alpha$}
\put(2,4){\circle*{0.15}} 
\put(1.4,3.9){\footnotesize $v_{15}$}
\put(2.2,3.9){\footnotesize $\beta$}
\put(4,4){\circle*{0.15}} 
\put(3.4,3.9){\footnotesize $v_{23}$}
\put(4.2,3.9){\footnotesize $\alpha$}
\put(1,5){\circle*{0.15}} 
\put(0.4,4.9){\footnotesize $v_{9}$}
\put(1.2,4.9){\footnotesize $\alpha$}
\put(5,5){\circle*{0.15}} 
\put(4.4,4.9){\footnotesize $v_{22}$}
\put(5.2,4.9){\footnotesize $\alpha$}
\put(9,5){\circle*{0.15}} 
\put(8.4,4.9){\footnotesize $v_{32}$}
\put(9.2,4.9){\footnotesize $\alpha$}
\put(2,6){\circle*{0.15}} 
\put(1.4,5.9){\footnotesize $v_{8}$}
\put(2.2,5.9){\footnotesize $\alpha$}
\put(4,6){\circle*{0.15}} 
\put(3.4,5.9){\footnotesize $v_{14}$}
\put(4.2,5.9){\footnotesize $\beta$}
\put(6,6){\circle*{0.15}} 
\put(5.4,5.9){\footnotesize $v_{21}$}
\put(6.2,5.9){\footnotesize $\alpha$}
\put(8,6){\circle*{0.15}} 
\put(7.4,5.9){\footnotesize $v_{28}$}
\put(8.2,5.9){\footnotesize $\beta$}
\put(1,7){\circle*{0.15}} 
\put(0.4,6.9){\footnotesize $v_{2}$}
\put(1.2,6.9){\footnotesize $\beta$}
\put(3,7){\circle*{0.15}} 
\put(2.4,6.9){\footnotesize $v_{7}$}
\put(3.2,6.9){\footnotesize $\alpha$}
\put(5,7){\circle*{0.15}} 
\put(4.4,6.9){\footnotesize $v_{13}$}
\put(5.2,6.9){\footnotesize $\beta$}
\put(7,7){\circle*{0.15}} 
\put(6.4,6.9){\footnotesize $v_{20}$}
\put(7.2,6.9){\footnotesize $\alpha$}
\put(11,7){\circle*{0.15}} 
\put(10.4,6.9){\footnotesize $v_{31}$}
\put(11.2,6.9){\footnotesize $\alpha$}
\put(4,8){\circle*{0.15}} 
\put(3.4,7.9){\footnotesize $v_{6}$}
\put(4.2,7.9){\footnotesize $\alpha$}
\put(8,8){\circle*{0.15}} 
\put(7.4,7.9){\footnotesize $v_{19}$}
\put(8.2,7.9){\footnotesize $\alpha$}
\put(10,8){\circle*{0.15}} 
\put(9.4,7.9){\footnotesize $v_{27}$}
\put(10.2,7.9){\footnotesize $\beta$}
\put(5,9){\circle*{0.15}} 
\put(4.4,8.9){\footnotesize $v_{5}$}
\put(5.2,8.9){\footnotesize $\alpha$}
\put(7,9){\circle*{0.15}} 
\put(6.4,8.9){\footnotesize $v_{12}$}
\put(7.2,8.9){\footnotesize $\beta$}
\put(9,9){\circle*{0.15}} 
\put(8.4,8.9){\footnotesize $v_{18}$}
\put(9.2,8.9){\footnotesize $\alpha$}
\put(4,10){\circle*{0.15}} 
\put(3.4,9.9){\footnotesize $v_{1}$}
\put(4.2,9.9){\footnotesize $\beta$}
\put(6,10){\circle*{0.15}} 
\put(5.4,9.9){\footnotesize $v_{4}$}
\put(6.2,9.9){\footnotesize $\alpha$}
\put(10,10){\circle*{0.15}} 
\put(9.4,9.9){\footnotesize $v_{17}$}
\put(10.2,9.9){\footnotesize $\alpha$}
\put(9,11){\circle*{0.15}} 
\put(8.4,10.9){\footnotesize $v_{11}$}
\put(9.2,10.9){\footnotesize $\beta$}
\put(8,12){\circle*{0.15}} 
\put(7.4,11.9){\footnotesize $v_{3}$}
\put(8.2,11.9){\footnotesize $\alpha$}
\put(9,5){\line(1,1){2}} 
\put(2,0){\line(1,1){8}} 
\put(1,1){\line(1,1){9}} 
\put(1,3){\line(1,1){8}} 
\put(0,4){\line(1,1){8}} 
\put(1,7){\line(1,1){3}}
\put(2,0){\line(-1,1){1}} 
\put(2,2){\line(-1,1){2}} 
\put(3,3){\line(-1,1){2}} 
\put(5,3){\line(-1,1){1}} 
\put(2,6){\line(-1,1){1}} 
\put(5,5){\line(-1,1){2}} 
\put(6,6){\line(-1,1){2}} 
\put(9,5){\line(-1,1){2}} 
\put(5,9){\line(-1,1){1}} 
\put(8,8){\line(-1,1){2}} 
\put(11,7){\line(-1,1){2}} 
\put(10,10){\line(-1,1){2}} 
\put(11.5,7.5){$\mathcal{C}_{6}$}
\put(10.5,8.5){$\mathcal{C}_{5}$}
\put(10.5,10.5){$\mathcal{C}_{4}$}
\put(9.5,11.5){$\mathcal{C}_{3}$}
\put(8.5,12.5){$\mathcal{C}_{2}$}
\put(4.5,10.5){$\mathcal{C}_{1}$}
\end{picture} 
\end{center}
\end{figure}

\begin{figure}[ht]
\begin{center}
\GridPosetsII: Depicted below are four \dichromatic grid posets each possessing the max property.\\  
{\small (Each is the `$\mathscr{G}$-semistandard poset $P_{\mathscr{G}}^{\alpha\beta}(2,2)$' for the indicated integral Coxeter--Dynkin posy $\mathscr{G}$.)} 

\setlength{\unitlength}{1cm}
\begin{picture}(6,7)
\put(0,7){\fbox{$\mathscr{G} = \myA_{1} \oplus \myA_{1}$}}
\put(0.5,2.75){
\begin{picture}(3,3.5)
\put(3,2){\circle*{0.15}} 
\put(2.4,1.9){\footnotesize $v_{4}$}
\put(3.2,1.9){\footnotesize $\beta$}
\put(4,3){\circle*{0.15}} 
\put(3.4,2.9){\footnotesize $v_{3}$}
\put(4.2,2.9){\footnotesize $\beta$}
\put(0,2){\circle*{0.15}} 
\put(-0.6,1.9){\footnotesize $v_{2}$}
\put(0.2,1.9){\footnotesize $\alpha$}
\put(1,3){\circle*{0.15}} 
\put(0.4,2.9){\footnotesize $v_{1}$}
\put(1.2,2.9){\footnotesize $\alpha$}
\put(3,2){\line(1,1){1}}
\put(0,2){\line(1,1){1}} 
\put(1.5,3.5){$\mathcal{C}_{1}$}
\put(4.5,3.5){$\mathcal{C}_{2}$}
\put(1.8,2.375){$\bigoplus$}
\end{picture}
}
\put(0,2){\fbox{$\mathscr{G} = \myA_{2}$}}
\put(0,-2.5){
\begin{picture}(3,3.5)
\put(4,2){\circle*{0.15}} 
\put(3.4,1.9){\footnotesize $v_{8}$}
\put(4.2,1.9){\footnotesize $\alpha$}
\put(5,3){\circle*{0.15}} 
\put(4.4,2.9){\footnotesize $v_{7}$}
\put(5.2,2.9){\footnotesize $\alpha$}
\put(1,1){\circle*{0.15}} 
\put(0.4,0.9){\footnotesize $v_{6}$}
\put(1.2,0.9){\footnotesize $\beta$}
\put(2,2){\circle*{0.15}} 
\put(1.4,1.9){\footnotesize $v_{5}$}
\put(2.2,1.9){\footnotesize $\beta$}
\put(3,3){\circle*{0.15}} 
\put(2.4,2.9){\footnotesize $v_{4}$}
\put(3.2,2.9){\footnotesize $\beta$}
\put(4,4){\circle*{0.15}} 
\put(3.4,3.9){\footnotesize $v_{3}$}
\put(4.2,3.9){\footnotesize $\beta$}
\put(0,2){\circle*{0.15}} 
\put(-0.6,1.9){\footnotesize $v_{2}$}
\put(0.2,1.9){\footnotesize $\alpha$}
\put(1,3){\circle*{0.15}} 
\put(0.4,2.9){\footnotesize $v_{1}$}
\put(1.2,2.9){\footnotesize $\alpha$}
\put(1,1){\line(1,1){3}}
\put(4,2){\line(1,1){1}}
\put(0,2){\line(1,1){1}} 
\put(1,1){\line(-1,1){1}} 
\put(2,2){\line(-1,1){1}} 
\put(3,3){\line(1,-1){1}}
\put(4,4){\line(1,-1){1}}
\put(1.5,3.5){$\mathcal{C}_{1}$}
\put(4.5,4.5){$\mathcal{C}_{2}$}
\put(5.5,3.5){$\mathcal{C}_{3}$}
\end{picture}
}
\end{picture}
\hspace*{1.5cm}
\setlength{\unitlength}{1cm}
\begin{picture}(8,8)
\put(0,6){\fbox{$\mathscr{G} = \myC_{2}$}}
\put(0,0.25){
\begin{picture}(3,3.5)
\put(1,3){\VertexForPosetsOrig{1}{\alpha}}
\put(0,2){\VertexForPosetsOrig{2}{\alpha}}
\put(6,6){\VertexForPosetsOrig{3}{\beta}}
\put(4,4){\VertexForPosetsOrig{4}{\beta}}
\put(2,2){\VertexForPosetsOrig{5}{\beta}}
\put(1,1){\VertexForPosetsOrig{6}{\beta}}
\put(7,5){\VertexForPosetsOrig{7}{\alpha}}
\put(6,4){\VertexForPosetsOrig{8}{\alpha}}
\put(5,3){\VertexForPosetsOrig{9}{\alpha}}
\put(4,2){\VertexForPosetsOrig{10}{\alpha}}
\put(3,1){\VertexForPosetsOrig{11}{\alpha}}
\put(2,0){\VertexForPosetsOrig{12}{\alpha}}
\put(7,3){\VertexForPosetsOrig{13}{\beta}}
\put(5,1){\VertexForPosetsOrig{14}{\beta}}
\put(0,2){\line(1,1){1}}
\put(1,1){\line(1,1){5}}
\put(2,0){\line(1,1){5}} 
\put(5,1){\line(1,1){2}} 
\put(2,0){\line(-1,1){2}} 
\put(3,1){\line(-1,1){2}} 
\put(5,1){\line(-1,1){1}} 
\put(5,3){\line(-1,1){1}} 
\put(7,3){\line(-1,1){1}}
\put(7,5){\line(-1,1){1}}
\put(1.5,3.5){$\mathcal{C}_{1}$}
\put(6.5,6.5){$\mathcal{C}_{2}$}
\put(7.5,5.5){$\mathcal{C}_{3}$}
\put(7.5,3.5){$\mathcal{C}_{4}$}
\end{picture}
}
\end{picture}

\setlength{\unitlength}{1cm}
\begin{picture}(11,13)
\put(-0.5,9){\fbox{$\mathscr{G} = \myG_{2}$}}
\put(2,7){\VertexForPosetsOrig{1}{\alpha}}
\put(0,5){\VertexForPosetsOrig{2}{\alpha}}
\put(9,12){\VertexForPosetsOrig{3}{\beta}}
\put(6,9){\VertexForPosetsOrig{4}{\beta}}
\put(3,6){\VertexForPosetsOrig{5}{\beta}}
\put(1,4){\VertexForPosetsOrig{6}{\beta}}
\put(10,11){\VertexForPosetsOrig{7}{\alpha}}
\put(9,10){\VertexForPosetsOrig{8}{\alpha}}
\put(8,9){\VertexForPosetsOrig{9}{\alpha}}
\put(7,8){\VertexForPosetsOrig{10}{\alpha}}
\put(6,7){\VertexForPosetsOrig{11}{\alpha}}
\put(5,6){\VertexForPosetsOrig{12}{\alpha}}
\put(4,5){\VertexForPosetsOrig{13}{\alpha}}
\put(3,4){\VertexForPosetsOrig{14}{\alpha}}
\put(2,3){\VertexForPosetsOrig{15}{\alpha}}
\put(1,2){\VertexForPosetsOrig{16}{\alpha}}
\put(10,9){\VertexForPosetsOrig{17}{\beta}}
\put(9,8){\VertexForPosetsOrig{18}{\beta}}
\put(7,6){\VertexForPosetsOrig{19}{\beta}}
\put(6,5){\VertexForPosetsOrig{20}{\beta}}
\put(4,3){\VertexForPosetsOrig{21}{\beta}}
\put(2,1){\VertexForPosetsOrig{22}{\beta}}
\put(11,8){\VertexForPosetsOrig{23}{\alpha}}
\put(10,7){\VertexForPosetsOrig{24}{\alpha}}
\put(9,6){\VertexForPosetsOrig{25}{\alpha}}
\put(8,5){\VertexForPosetsOrig{26}{\alpha}}
\put(7,4){\VertexForPosetsOrig{27}{\alpha}}
\put(6,3){\VertexForPosetsOrig{28}{\alpha}}
\put(5,2){\VertexForPosetsOrig{29}{\alpha}}
\put(3,0){\VertexForPosetsOrig{30}{\alpha}}
\put(10,5){\VertexForPosetsOrig{31}{\beta}}
\put(7,2){\VertexForPosetsOrig{32}{\beta}}
\put(7,2){\line(1,1){3}} 
\put(3,0){\line(1,1){8}} 
\put(2,1){\line(1,1){8}} 
\put(1,2){\line(1,1){9}} 
\put(1,4){\line(1,1){8}} 
\put(0,5){\line(1,1){2}}
\put(3,0){\line(-1,1){2}} 
\put(2,3){\line(-1,1){2}} 
\put(5,2){\line(-1,1){2}} 
\put(4,5){\line(-1,1){2}} 
\put(7,2){\line(-1,1){1}} 
\put(7,4){\line(-1,1){2}} 
\put(8,5){\line(-1,1){2}} 
\put(7,8){\line(-1,1){1}} 
\put(10,5){\line(-1,1){1}} 
\put(10,7){\line(-1,1){2}} 
\put(11,8){\line(-1,1){2}} 
\put(10,11){\line(-1,1){1}} 
\put(10.5,5.5){$\mathcal{C}_{6}$}
\put(11.5,8.5){$\mathcal{C}_{5}$}
\put(10.5,9.5){$\mathcal{C}_{4}$}
\put(10.5,11.5){$\mathcal{C}_{3}$}
\put(9.5,12.5){$\mathcal{C}_{2}$}
\put(2.5,7.5){$\mathcal{C}_{1}$}
\end{picture} 
\end{center}
\end{figure}

\newcommand{\AOneAOneAlpha}{
\setlength{\unitlength}{1cm}
\begin{picture}(0.5,0.75)
\put(-0.45,0.2){\footnotesize $v_{1}$}
\put(0,0.25){\circle*{0.15}} 
\put(0.2,0.25){\footnotesize $\alpha$}
\end{picture}
}
\newcommand{\AOneAOneBeta}{
\setlength{\unitlength}{1cm}
\begin{picture}(0.5,0.75)
\put(-0.45,0.2){\footnotesize $v_{1}$}
\put(0,0.25){\circle*{0.15}} 
\put(0.2,0.25){\footnotesize $\beta$}
\end{picture}
}
\newcommand{\ATwoAlpha}{
\setlength{\unitlength}{1cm}
\begin{picture}(2,2)
\put(0,0.25){
\begin{picture}(2,2)
\put(0.5,-0.1){\footnotesize $v_{2}$}
\put(1,0){\circle*{0.15}} 
\put(1.2,-0.1){\footnotesize $\beta$}
\put(-0.5,0.9){\footnotesize $v_{1}$}
\put(0,1){\circle*{0.15}} 
\put(0.2,0.9){\footnotesize $\alpha$}
\put(1,0){\line(-1,1){1}} 
\end{picture}
}
\end{picture}
}
\newcommand{\ATwoBeta}{
\setlength{\unitlength}{1cm}
\begin{picture}(2,2)
\put(0,0.25){
\begin{picture}(2,2)
\put(0.5,-0.1){\footnotesize $v_{2}$}
\put(1,0){\circle*{0.15}} 
\put(1.2,-0.1){\footnotesize $\alpha$}
\put(-0.5,0.9){\footnotesize $v_{1}$}
\put(0,1){\circle*{0.15}} 
\put(0.2,0.9){\footnotesize $\beta$}
\put(1,0){\line(-1,1){1}} 
\end{picture}
}
\end{picture}
}
\newcommand{\BTwoAAlpha}{
\setlength{\unitlength}{1cm}
\begin{picture}(3,4)
\put(0,0.75){
\begin{picture}(3,3)
\put(1.5,-0.1){\footnotesize $v_{3}$}
\put(2,0){\circle*{0.15}} 
\put(2.2,-0.1){\footnotesize $\alpha$}
\put(0.5,0.9){\footnotesize $v_{2}$}
\put(1,1){\circle*{0.15}} 
\put(1.2,0.9){\footnotesize $\beta$}
\put(-0.5,1.9){\footnotesize $v_{1}$}
\put(0,2){\circle*{0.15}} 
\put(0.2,1.9){\footnotesize $\alpha$}
\put(2,0){\line(-1,1){2}} 
\end{picture}
}
\end{picture}
}
\newcommand{\BTwoBBeta}{
\setlength{\unitlength}{1cm}
\begin{picture}(4,4)
\put(1,0.25){
\begin{picture}(3,3.5)
\put(0.5,-0.1){\footnotesize $v_{4}$}
\put(1,0){\circle*{0.15}} 
\put(1.2,-0.1){\footnotesize $\beta$}
\put(-0.5,0.9){\footnotesize $v_{3}$}
\put(0,1){\circle*{0.15}} 
\put(0.2,0.9){\footnotesize $\alpha$}
\put(0.5,1.9){\footnotesize $v_{2}$}
\put(1,2){\circle*{0.15}} 
\put(1.2,1.9){\footnotesize $\alpha$}
\put(-0.5,2.9){\footnotesize $v_{1}$}
\put(0,3){\circle*{0.15}} 
\put(0.2,2.9){\footnotesize $\beta$}
\put(0,1){\line(1,1){1}} 
\put(1,0){\line(-1,1){1}} 
\put(1,2){\line(-1,1){1}} 
\end{picture}
}
\end{picture}
}
\newcommand{\GTwoAAlpha}{
\setlength{\unitlength}{1cm}
\begin{picture}(4,6)
\put(0.25,1.25){
\begin{picture}(4,6)
\put(2.5,-0.1){\footnotesize $v_{6}$}
\put(3,0){\circle*{0.15}} 
\put(3.2,-0.1){\footnotesize $\alpha$}
\put(1.5,0.9){\footnotesize $v_{5}$}
\put(2,1){\circle*{0.15}} 
\put(2.2,0.9){\footnotesize $\beta$}
\put(0.5,1.9){\footnotesize $v_{4}$}
\put(1,2){\circle*{0.15}} 
\put(1.2,1.9){\footnotesize $\alpha$}
\put(1.5,2.9){\footnotesize $v_{3}$}
\put(2,3){\circle*{0.15}} 
\put(2.2,2.9){\footnotesize $\alpha$}
\put(0.5,3.9){\footnotesize $v_{2}$}
\put(1,4){\circle*{0.15}} 
\put(1.2,3.9){\footnotesize $\beta$}
\put(-0.5,4.9){\footnotesize $v_{1}$}
\put(0,5){\circle*{0.15}} 
\put(0.2,4.9){\footnotesize $\alpha$}
\put(1,2){\line(1,1){1}} 
\put(3,0){\line(-1,1){2}} 
\put(2,3){\line(-1,1){2}} 
\end{picture}
}
\end{picture}
}
\newcommand{\GTwoBBeta}{
\setlength{\unitlength}{1cm}
\begin{picture}(4,8)
\put(0.25,0.25){
\begin{picture}(4,8)
\put(1.5,-0.1){\footnotesize $v_{10}$}
\put(2,0){\circle*{0.15}} 
\put(2.2,-0.15){\footnotesize $\beta$}
\put(0.5,0.9){\footnotesize $v_{9}$}
\put(1,1){\circle*{0.15}} 
\put(1.2,0.9){\footnotesize $\alpha$}
\put(1.5,1.9){\footnotesize $v_{8}$}
\put(2,2){\circle*{0.15}} 
\put(2.2,1.9){\footnotesize $\alpha$}
\put(0.5,2.9){\footnotesize $v_{6}$}
\put(1,3){\circle*{0.15}} 
\put(1.2,2.9){\footnotesize $\beta$}
\put(2.5,2.9){\footnotesize $v_{7}$}
\put(3,3){\circle*{0.15}} 
\put(3.2,2.9){\footnotesize $\alpha$}
\put(-0.5,3.9){\footnotesize $v_{4}$}
\put(0,4){\circle*{0.15}} 
\put(0.2,3.9){\footnotesize $\alpha$}
\put(1.5,3.9){\footnotesize $v_{5}$}
\put(2,4){\circle*{0.15}} 
\put(2.2,3.9){\footnotesize $\beta$}
\put(0.5,4.9){\footnotesize $v_{3}$}
\put(1,5){\circle*{0.15}} 
\put(1.2,4.9){\footnotesize $\alpha$}
\put(1.5,5.9){\footnotesize $v_{2}$}
\put(2,6){\circle*{0.15}} 
\put(2.2,5.9){\footnotesize $\alpha$}
\put(0.5,6.9){\footnotesize $v_{1}$}
\put(1,7){\circle*{0.15}} 
\put(1.2,6.9){\footnotesize $\beta$}
\put(1,1){\line(1,1){2}} 
\put(1,3){\line(1,1){1}} 
\put(0,4){\line(1,1){2}} 
\put(2,0){\line(-1,1){1}} 
\put(2,2){\line(-1,1){2}} 
\put(3,3){\line(-1,1){2}} 
\put(2,6){\line(-1,1){1}} 
\end{picture} 
}
\end{picture}
}


\newcommand{\AoneAlphaIdeals}{
\setlength{\unitlength}{1.1cm}
\begin{picture}(3,3)
\put(0.25,2.5){$\LAone(\omega_{\alpha})$}
\put(1.5,1){\circle*{0.125}}
\put(1.5,2){\circle*{0.125}}
\put(1.7,0.95){\tiny $\emptyset$}
\put(1.6,1.95){\tiny $\langle 1 \rangle$}
\put(1.5,1){\line(0,1){1}}
\put(1.45,1.45){\footnotesize $\alpha$}
\end{picture}
}

\newcommand{\AoneAlphaTableaux}{
\setlength{\unitlength}{1.1cm}
\begin{picture}(3,3)
\put(0.25,2.5){$\LAone(\omega_{\alpha})$}
\put(1.5,1){\circle*{0.125}}
\put(1.5,2){\circle*{0.125}}
\put(1.6,0.95){\smallonebox{2}}
\put(1.6,1.95){\smallonebox{1}}
\put(1.5,1){\line(0,1){1}}
\put(1.45,1.45){\footnotesize $\alpha$}
\end{picture}
}

\newcommand{\AoneAlphaWeights}{
\setlength{\unitlength}{1.1cm}
\begin{picture}(3,3)
\put(0.25,2.5){$\LAone(\omega_{\alpha})$}
\put(1.5,1){\circle*{0.125}}
\put(1.5,2){\circle*{0.125}}
\put(1.6,0.95){(-1,0)}
\put(1.6,1.95){(1,0)}
\put(1.5,1){\line(0,1){1}}
\put(1.45,1.45){\footnotesize $\alpha$}
\end{picture}
}

\newcommand{\AoneBetaIdeals}{
\setlength{\unitlength}{1.1cm}
\begin{picture}(3,3)
\put(0.25,2.5){$\LAone(\omega_{\beta})$}
\put(1.5,1){\circle*{0.125}}
\put(1.5,2){\circle*{0.125}}
\put(1.7,0.95){\tiny $\emptyset$}
\put(1.6,1.95){\tiny $\langle 1 \rangle$}
\put(1.5,1){\line(0,1){1}}
\put(1.45,1.45){\footnotesize $\beta$}
\end{picture}
}

\newcommand{\AoneBetaTableaux}{
\setlength{\unitlength}{1.1cm}
\begin{picture}(3,3)
\put(0.25,2.5){$\LAone(\omega_{\beta})$}
\put(1.5,1){\circle*{0.125}}
\put(1.5,2){\circle*{0.125}}
\put(1.6,0.95){\smalltwobox{1}{3}}
\put(1.6,1.95){\smalltwobox{1}{2}}
\put(1.5,1){\line(0,1){1}}
\put(1.45,1.45){\footnotesize $\beta$}
\end{picture}
}

\newcommand{\AoneBetaWeights}{
\setlength{\unitlength}{1.1cm}
\begin{picture}(3,3)
\put(0.25,2.5){$\LAone(\omega_{\beta})$}
\put(1.5,1){\circle*{0.125}}
\put(1.5,2){\circle*{0.125}}
\put(1.6,0.95){(0,-1)}
\put(1.6,1.95){(0,1)}
\put(1.5,1){\line(0,1){1}}
\put(1.45,1.45){\footnotesize $\beta$}
\end{picture}
}

\newcommand{\AtwoAlphaIdeals}{
\setlength{\unitlength}{1.1cm}
\begin{picture}(3,4)
\put(0.25,3.5){$\LAtwo(\omega_{\alpha})$}
\put(1.5,1){\circle*{0.125}}
\put(1.5,2){\circle*{0.125}}
\put(1.5,3){\circle*{0.125}}
\put(1.7,0.95){\tiny $\emptyset$}
\put(1.6,1.95){\tiny $\langle 2 \rangle$}
\put(1.6,2.95){\tiny $\langle 1 \rangle$}
\put(1.5,1){\line(0,1){2}}
\put(1.45,1.45){\footnotesize $\beta$}
\put(1.45,2.45){\footnotesize $\alpha$}
\end{picture}
}

\newcommand{\AtwoAlphaTableaux}{
\setlength{\unitlength}{1.1cm}
\begin{picture}(3,4)
\put(0.25,3.5){$\LAtwo(\omega_{\alpha})$}
\put(1.5,1){\circle*{0.125}}
\put(1.5,2){\circle*{0.125}}
\put(1.5,3){\circle*{0.125}}
\put(1.6,0.95){\smallonebox{3}}
\put(1.6,1.95){\smallonebox{2}}
\put(1.6,2.95){\smallonebox{1}}
\put(1.5,1){\line(0,1){2}}
\put(1.45,1.45){\footnotesize $\beta$}
\put(1.45,2.45){\footnotesize $\alpha$}
\end{picture}
}

\newcommand{\AtwoAlphaWeights}{
\setlength{\unitlength}{1.1cm}
\begin{picture}(3,4)
\put(0.25,3.5){$\LAtwo(\omega_{\alpha})$}
\put(1.5,1){\circle*{0.125}}
\put(1.5,2){\circle*{0.125}}
\put(1.5,3){\circle*{0.125}}
\put(1.6,0.95){(0,-1)}
\put(1.6,1.95){(-1,1)}
\put(1.6,2.95){(1,0)}
\put(1.5,1){\line(0,1){2}}
\put(1.45,1.45){\footnotesize $\beta$}
\put(1.45,2.45){\footnotesize $\alpha$}
\end{picture}
}

\newcommand{\AtwoBetaIdeals}{
\setlength{\unitlength}{1.1cm}
\begin{picture}(3,4)
\put(0.25,3.5){$\LAtwo(\omega_{\beta})$}
\put(1.5,1){\circle*{0.125}}
\put(1.5,2){\circle*{0.125}}
\put(1.5,3){\circle*{0.125}}
\put(1.7,0.95){\tiny $\emptyset$}
\put(1.6,1.95){\tiny $\langle 2 \rangle$}
\put(1.6,2.95){\tiny $\langle 1 \rangle$}
\put(1.5,1){\line(0,1){2}}
\put(1.45,1.45){\footnotesize $\alpha$}
\put(1.45,2.45){\footnotesize $\beta$}
\end{picture}
}

\newcommand{\AtwoBetaTableaux}{
\setlength{\unitlength}{1.1cm}
\begin{picture}(3,4)
\put(0.25,3.5){$\LAtwo(\omega_{\beta})$}
\put(1.5,1){\circle*{0.125}}
\put(1.5,2){\circle*{0.125}}
\put(1.5,3){\circle*{0.125}}
\put(1.6,0.95){\smalltwobox{2}{3}}
\put(1.6,1.95){\smalltwobox{1}{3}}
\put(1.6,2.95){\smalltwobox{1}{2}}
\put(1.5,1){\line(0,1){2}}
\put(1.45,1.45){\footnotesize $\alpha$}
\put(1.45,2.45){\footnotesize $\beta$}
\end{picture}
}

\newcommand{\AtwoBetaWeights}{
\setlength{\unitlength}{1.1cm}
\begin{picture}(3,4)
\put(0.25,3.5){$\LAtwo(\omega_{\beta})$}
\put(1.5,1){\circle*{0.125}}
\put(1.5,2){\circle*{0.125}}
\put(1.5,3){\circle*{0.125}}
\put(1.6,0.95){(-1,0)}
\put(1.6,1.95){(1,-1)}
\put(1.6,2.95){(0,1)}
\put(1.5,1){\line(0,1){2}}
\put(1.45,1.45){\footnotesize $\alpha$}
\put(1.45,2.45){\footnotesize $\beta$}
\end{picture}
}

\newcommand{\BtwoAlphaIdeals}{
\setlength{\unitlength}{1.1cm}
\begin{picture}(3,5)
\put(0.25,4.5){$\LBtwo(\omega_{\alpha})$}
\put(1.5,1){\circle*{0.125}}
\put(1.5,2){\circle*{0.125}}
\put(1.5,3){\circle*{0.125}}
\put(1.5,4){\circle*{0.125}}
\put(1.7,0.95){\tiny $\emptyset$}
\put(1.6,1.95){\tiny $\langle 3 \rangle$}
\put(1.6,2.95){\tiny $\langle 2 \rangle$}
\put(1.6,3.95){\tiny $\langle 1 \rangle$}
\put(1.5,1){\line(0,1){3}}
\put(1.45,1.45){\footnotesize $\alpha$}
\put(1.45,2.45){\footnotesize $\beta$}
\put(1.45,3.45){\footnotesize $\alpha$}
\end{picture}
}

\newcommand{\BtwoAlphaTableaux}{
\setlength{\unitlength}{1.1cm}
\begin{picture}(3,5)
\put(0.25,4.5){$\LBtwo(\omega_{\alpha})$}
\put(1.5,1){\circle*{0.125}}
\put(1.5,2){\circle*{0.125}}
\put(1.5,3){\circle*{0.125}}
\put(1.5,4){\circle*{0.125}}
\put(1.6,0.95){\smallonebox{4}}
\put(1.6,1.95){\smallonebox{3}}
\put(1.6,2.95){\smallonebox{2}}
\put(1.6,3.95){\smallonebox{1}}
\put(1.5,1){\line(0,1){3}}
\put(1.45,1.45){\footnotesize $\alpha$}
\put(1.45,2.45){\footnotesize $\beta$}
\put(1.45,3.45){\footnotesize $\alpha$}
\end{picture}
}

\newcommand{\BtwoAlphaWeights}{
\setlength{\unitlength}{1.1cm}
\begin{picture}(3,5)
\put(0.25,4.5){$\LBtwo(\omega_{\alpha})$}
\put(1.5,1){\circle*{0.125}}
\put(1.5,2){\circle*{0.125}}
\put(1.5,3){\circle*{0.125}}
\put(1.5,4){\circle*{0.125}}
\put(1.6,0.95){(-1,0)}
\put(1.6,1.95){(1,-1)}
\put(1.6,2.95){(-1,1)}
\put(1.6,3.95){(1,0)}
\put(1.5,1){\line(0,1){3}}
\put(1.45,1.45){\footnotesize $\alpha$}
\put(1.45,2.45){\footnotesize $\beta$}
\put(1.45,3.45){\footnotesize $\alpha$}
\end{picture}
}

\newcommand{\BtwoBetaIdeals}{
\setlength{\unitlength}{1.1cm}
\begin{picture}(3,6)
\put(0.25,5.5){$\LBtwo(\omega_{\beta})$}
\put(1.5,1){\circle*{0.125}}
\put(1.5,2){\circle*{0.125}}
\put(1.5,3){\circle*{0.125}}
\put(1.5,4){\circle*{0.125}}
\put(1.5,5){\circle*{0.125}}
\put(1.7,0.95){\tiny $\emptyset$}
\put(1.6,1.95){\tiny $\langle 4 \rangle$}
\put(1.6,2.95){\tiny $\langle 3 \rangle$}
\put(1.6,3.95){\tiny $\langle 2 \rangle$}
\put(1.6,4.95){\tiny $\langle 1 \rangle$}
\put(1.5,1){\line(0,1){4}}
\put(1.45,1.45){\footnotesize $\beta$}
\put(1.45,2.45){\footnotesize $\alpha$}
\put(1.45,3.45){\footnotesize $\alpha$}
\put(1.45,4.45){\footnotesize $\beta$}
\end{picture}
}

\newcommand{\BtwoBetaTableaux}{
\setlength{\unitlength}{1.1cm}
\begin{picture}(3,6)
\put(0.25,5.5){$\LBtwo(\omega_{\beta})$}
\put(1.5,1){\circle*{0.125}}
\put(1.5,2){\circle*{0.125}}
\put(1.5,3){\circle*{0.125}}
\put(1.5,4){\circle*{0.125}}
\put(1.5,5){\circle*{0.125}}
\put(1.6,0.95){\smalltwobox{3}{4}}
\put(1.6,1.95){\smalltwobox{2}{4}}
\put(1.6,2.95){\smalltwobox{2}{3}}
\put(1.6,3.95){\smalltwobox{1}{3}}
\put(1.6,4.95){\smalltwobox{1}{2}}
\put(1.5,1){\line(0,1){4}}
\put(1.45,1.45){\footnotesize $\beta$}
\put(1.45,2.45){\footnotesize $\alpha$}
\put(1.45,3.45){\footnotesize $\alpha$}
\put(1.45,4.45){\footnotesize $\beta$}
\end{picture}
}

\newcommand{\BtwoBetaWeights}{
\setlength{\unitlength}{1.1cm}
\begin{picture}(3,6)
\put(0.25,5.5){$\LBtwo(\omega_{\beta})$}
\put(1.5,1){\circle*{0.125}}
\put(1.5,2){\circle*{0.125}}
\put(1.5,3){\circle*{0.125}}
\put(1.5,4){\circle*{0.125}}
\put(1.5,5){\circle*{0.125}}
\put(1.6,0.95){(0,-1)}
\put(1.6,1.95){(-2,1)}
\put(1.6,2.95){(0,0)}
\put(1.6,3.95){(2,-1)}
\put(1.6,4.95){(0,1)}
\put(1.5,1){\line(0,1){4}}
\put(1.45,1.45){\footnotesize $\beta$}
\put(1.45,2.45){\footnotesize $\alpha$}
\put(1.45,3.45){\footnotesize $\alpha$}
\put(1.45,4.45){\footnotesize $\beta$}
\end{picture}
}

\newcommand{\GtwoAlphaIdeals}{
\setlength{\unitlength}{1.1cm}
\begin{picture}(3,8)
\put(0.25,7.5){$\LGtwo(\omega_{\alpha})$}
\put(1.5,1){\circle*{0.125}}
\put(1.5,2){\circle*{0.125}}
\put(1.5,3){\circle*{0.125}}
\put(1.5,4){\circle*{0.125}}
\put(1.5,5){\circle*{0.125}}
\put(1.5,6){\circle*{0.125}}
\put(1.5,7){\circle*{0.125}}
\put(1.7,0.95){\tiny $\emptyset$}
\put(1.6,1.95){\tiny $\langle 6 \rangle$}
\put(1.6,2.95){\tiny $\langle 5 \rangle$}
\put(1.6,3.95){\tiny $\langle 4 \rangle$}
\put(1.6,4.95){\tiny $\langle 3 \rangle$}
\put(1.6,5.95){\tiny $\langle 2 \rangle$}
\put(1.6,6.95){\tiny $\langle 1 \rangle$}
\put(1.5,1){\line(0,1){6}}
\put(1.45,1.45){\footnotesize $\alpha$}
\put(1.45,2.45){\footnotesize $\beta$}
\put(1.45,3.45){\footnotesize $\alpha$}
\put(1.45,4.45){\footnotesize $\alpha$}
\put(1.45,5.45){\footnotesize $\beta$}
\put(1.45,6.45){\footnotesize $\alpha$}
\end{picture}
}

\newcommand{\GtwoAlphaTableaux}{
\setlength{\unitlength}{1.1cm}
\begin{picture}(3,8)
\put(0.25,7.5){$\LGtwo(\omega_{\alpha})$}
\put(1.5,1){\circle*{0.125}}
\put(1.5,2){\circle*{0.125}}
\put(1.5,3){\circle*{0.125}}
\put(1.5,4){\circle*{0.125}}
\put(1.5,5){\circle*{0.125}}
\put(1.5,6){\circle*{0.125}}
\put(1.5,7){\circle*{0.125}}
\put(1.6,0.95){\smallonebox{7}}
\put(1.6,1.95){\smallonebox{6}}
\put(1.6,2.95){\smallonebox{5}}
\put(1.6,3.95){\smallonebox{4}}
\put(1.6,4.95){\smallonebox{3}}
\put(1.6,5.95){\smallonebox{2}}
\put(1.6,6.95){\smallonebox{1}}
\put(1.5,1){\line(0,1){6}}
\put(1.45,1.45){\footnotesize $\alpha$}
\put(1.45,2.45){\footnotesize $\beta$}
\put(1.45,3.45){\footnotesize $\alpha$}
\put(1.45,4.45){\footnotesize $\alpha$}
\put(1.45,5.45){\footnotesize $\beta$}
\put(1.45,6.45){\footnotesize $\alpha$}
\end{picture}
}

\newcommand{\GtwoAlphaWeights}{
\setlength{\unitlength}{1.1cm}
\begin{picture}(3,8)
\put(0.25,7.5){$\LGtwo(\omega_{\alpha})$}
\put(1.5,1){\circle*{0.125}}
\put(1.5,2){\circle*{0.125}}
\put(1.5,3){\circle*{0.125}}
\put(1.5,4){\circle*{0.125}}
\put(1.5,5){\circle*{0.125}}
\put(1.5,6){\circle*{0.125}}
\put(1.5,7){\circle*{0.125}}
\put(1.6,0.95){(-1,0)}
\put(1.6,1.95){(1,-1)}
\put(1.6,2.95){(-2,1)}
\put(1.6,3.95){(0,0)}
\put(1.6,4.95){(2,-1)}
\put(1.6,5.95){(-1,1)}
\put(1.6,6.95){(1,0)}
\put(1.5,1){\line(0,1){6}}
\put(1.45,1.45){\footnotesize $\alpha$}
\put(1.45,2.45){\footnotesize $\beta$}
\put(1.45,3.45){\footnotesize $\alpha$}
\put(1.45,4.45){\footnotesize $\alpha$}
\put(1.45,5.45){\footnotesize $\beta$}
\put(1.45,6.45){\footnotesize $\alpha$}
\end{picture}
}

\newcommand{\GtwoBetaIdeals}{
\setlength{\unitlength}{1.1cm}
\begin{picture}(5,12)
\put(1.5,10){$\LGtwo(\omega_{\beta})$}
\put(5,6){\circle*{0.125}}
\put(2,5){\circle*{0.125}}
\put(2,7){\circle*{0.125}}
\put(3,1){\circle*{0.125}}
\put(3,2){\circle*{0.125}}
\put(3,3){\circle*{0.125}}
\put(3,4){\circle*{0.125}}
\put(3,6){\circle*{0.125}}
\put(3,8){\circle*{0.125}}
\put(3,9){\circle*{0.125}}
\put(3,10){\circle*{0.125}}
\put(3,11){\circle*{0.125}}
\put(4,5){\circle*{0.125}}
\put(4,7){\circle*{0.125}}
\put(3.15,10.95){\tiny $\langle 1 \rangle$}
\put(3.15,9.95){\tiny $\langle 2 \rangle$}
\put(3.15,8.95){\tiny $\langle 3 \rangle$}
\put(3.15,7.95){\tiny $\langle 4,5 \rangle$}
\put(4.15,7){\tiny $\langle 4,7 \rangle$}
\put(1.45,6.95){\tiny $\langle 5 \rangle$}
\put(2.35,5.95){\tiny $\langle 6,7 \rangle$}
\put(5.1,5.95){\tiny $\langle 4 \rangle$}
\put(1.45,4.95){\tiny $\langle 7 \rangle$}
\put(4.15,4.9){\tiny $\langle 6 \rangle$}
\put(3.15,3.95){\tiny $\langle 8 \rangle$}
\put(3.1,2.95){\tiny $\langle 9 \rangle$}
\put(3.1,1.95){\tiny $\langle 10 \rangle$}
\put(3.2,0.95){\tiny $\emptyset$}
\put(3,1){\line(0,1){3}}
\put(3,8){\line(0,1){3}}
\put(3,4){\line(-1,1){1}}
\put(3,4){\line(1,1){1}}
\put(2,5){\line(1,1){2}}
\put(4,5){\line(-1,1){2}}
\put(2,7){\line(1,1){1}}
\put(4,7){\line(-1,1){1}}
\put(4,5){\line(1,1){1}}
\put(5,6){\line(-1,1){1}}
\put(2.95,10.45){\footnotesize $\beta$}
\put(2.95,9.45){\footnotesize $\alpha$}
\put(2.95,8.45){\footnotesize $\alpha$}
\put(2.45,7.45){\footnotesize $\alpha$}
\put(3.45,7.45){\footnotesize $\beta$}
\put(2.45,6.45){\footnotesize $\beta$}
\put(3.45,6.45){\footnotesize $\alpha$}
\put(4.45,6.45){\footnotesize $\alpha$}
\put(2.45,5.45){\footnotesize $\beta$}
\put(3.45,5.45){\footnotesize $\alpha$}
\put(4.45,5.45){\footnotesize $\alpha$}
\put(2.45,4.45){\footnotesize $\alpha$}
\put(3.45,4.45){\footnotesize $\beta$}
\put(2.95,3.45){\footnotesize $\alpha$}
\put(2.95,2.45){\footnotesize $\alpha$}
\put(2.95,1.45){\footnotesize $\beta$}
\end{picture}
}

\newcommand{\GtwoBetaTableaux}{
\setlength{\unitlength}{1.1cm}
\begin{picture}(5,12)
\put(1.5,10){$\LGtwo(\omega_{\beta})$}
\put(5,6){\circle*{0.125}}
\put(2,5){\circle*{0.125}}
\put(2,7){\circle*{0.125}}
\put(3,1){\circle*{0.125}}
\put(3,2){\circle*{0.125}}
\put(3,3){\circle*{0.125}}
\put(3,4){\circle*{0.125}}
\put(3,6){\circle*{0.125}}
\put(3,8){\circle*{0.125}}
\put(3,9){\circle*{0.125}}
\put(3,10){\circle*{0.125}}
\put(3,11){\circle*{0.125}}
\put(4,5){\circle*{0.125}}
\put(4,7){\circle*{0.125}}
\put(3.15,10.95){\smalltwobox{1}{2}}
\put(3.15,9.95){\smalltwobox{1}{3}}
\put(3.15,8.95){\smalltwobox{1}{4}}
\put(3.15,7.95){\smalltwobox{1}{5}}
\put(4.15,7){\smalltwobox{1}{6}}
\put(1.45,6.95){\smalltwobox{2}{5}}
\put(2.35,5.95){\smalltwobox{2}{6}}
\put(5.1,5.95){\smalltwobox{1}{7}}
\put(1.45,4.95){\smalltwobox{3}{6}}
\put(4.15,4.9){\smalltwobox{2}{7}}
\put(3.15,3.95){\smalltwobox{3}{7}}
\put(3.1,2.95){\smalltwobox{4}{7}}
\put(3.1,1.95){\smalltwobox{5}{7}}
\put(3.1,0.95){\smalltwobox{6}{7}}
\put(3,1){\line(0,1){3}}
\put(3,8){\line(0,1){3}}
\put(3,4){\line(-1,1){1}}
\put(3,4){\line(1,1){1}}
\put(2,5){\line(1,1){2}}
\put(4,5){\line(-1,1){2}}
\put(2,7){\line(1,1){1}}
\put(4,7){\line(-1,1){1}}
\put(4,5){\line(1,1){1}}
\put(5,6){\line(-1,1){1}}
\put(2.95,10.45){\footnotesize $\beta$}
\put(2.95,9.45){\footnotesize $\alpha$}
\put(2.95,8.45){\footnotesize $\alpha$}
\put(2.45,7.45){\footnotesize $\alpha$}
\put(3.45,7.45){\footnotesize $\beta$}
\put(2.45,6.45){\footnotesize $\beta$}
\put(3.45,6.45){\footnotesize $\alpha$}
\put(4.45,6.45){\footnotesize $\alpha$}
\put(2.45,5.45){\footnotesize $\beta$}
\put(3.45,5.45){\footnotesize $\alpha$}
\put(4.45,5.45){\footnotesize $\alpha$}
\put(2.45,4.45){\footnotesize $\alpha$}
\put(3.45,4.45){\footnotesize $\beta$}
\put(2.95,3.45){\footnotesize $\alpha$}
\put(2.95,2.45){\footnotesize $\alpha$}
\put(2.95,1.45){\footnotesize $\beta$}
\end{picture}
}

\newcommand{\GtwoBetaWeights}{
\setlength{\unitlength}{1.1cm}
\begin{picture}(5,12)
\put(1.5,10){$\LGtwo(\omega_{\beta})$}
\put(5,6){\circle*{0.125}}
\put(2,5){\circle*{0.125}}
\put(2,7){\circle*{0.125}}
\put(3,1){\circle*{0.125}}
\put(3,2){\circle*{0.125}}
\put(3,3){\circle*{0.125}}
\put(3,4){\circle*{0.125}}
\put(3,6){\circle*{0.125}}
\put(3,8){\circle*{0.125}}
\put(3,9){\circle*{0.125}}
\put(3,10){\circle*{0.125}}
\put(3,11){\circle*{0.125}}
\put(4,5){\circle*{0.125}}
\put(4,7){\circle*{0.125}}
\put(3.15,10.95){(0,1)}
\put(3.15,9.95){(3,-1)}
\put(3.15,8.95){(1,0)}
\put(3.15,7.95){(-1,1)}
\put(4.15,7){(2,-1)}
\put(1.05,6.95){(-3,2)}
\put(2.,5.95){(0,0)}
\put(5,6.2){(0,0)}
\put(1.05,4.95){(3,-2)}
\put(4.15,4.9){(-2,1)}
\put(3.15,3.95){(1,-1)}
\put(3.1,2.95){(-1,0)}
\put(3.1,1.95){(-3,1)}
\put(3.1,0.95){(0,-1)}
\put(3,1){\line(0,1){3}}
\put(3,8){\line(0,1){3}}
\put(3,4){\line(-1,1){1}}
\put(3,4){\line(1,1){1}}
\put(2,5){\line(1,1){2}}
\put(4,5){\line(-1,1){2}}
\put(2,7){\line(1,1){1}}
\put(4,7){\line(-1,1){1}}
\put(4,5){\line(1,1){1}}
\put(5,6){\line(-1,1){1}}
\put(2.95,10.45){\footnotesize $\beta$}
\put(2.95,9.45){\footnotesize $\alpha$}
\put(2.95,8.45){\footnotesize $\alpha$}
\put(2.45,7.45){\footnotesize $\alpha$}
\put(3.45,7.45){\footnotesize $\beta$}
\put(2.45,6.45){\footnotesize $\beta$}
\put(3.45,6.45){\footnotesize $\alpha$}
\put(4.45,6.45){\footnotesize $\alpha$}
\put(2.45,5.45){\footnotesize $\beta$}
\put(3.45,5.45){\footnotesize $\alpha$}
\put(4.45,5.45){\footnotesize $\alpha$}
\put(2.45,4.45){\footnotesize $\alpha$}
\put(3.45,4.45){\footnotesize $\beta$}
\put(2.95,3.45){\footnotesize $\alpha$}
\put(2.95,2.45){\footnotesize $\alpha$}
\put(2.95,1.45){\footnotesize $\beta$}
\end{picture}
}

\begin{figure}[t]
\begin{center}
\FundPosets: The $\mathscr{G}$-fundamental compression posets. 

\ 

\begin{tabular}{|c||c|c|}
\hline
\parbox{1.5in}{\small \begin{center}Two-node\\ Coxeter--Dynkin posy $\mathscr{G}$\end{center}} & $P_{\mathscr{G}}(1,0)$ & $P_{\mathscr{G}}(0,1)$\\
\hline
\hline
\parbox[b]{0.5in}{\ 

$\myA_{1} \oplus \myA_{1}$

\ 

\vspace*{-0.1in}
} & \AOneAOneAlpha & \AOneAOneBeta\\
\hline
\parbox[b]{0.225in}{$\myA_{2}$
\vspace*{0.4in}} 
 & \ATwoAlpha & \ATwoBeta\\
\hline
\parbox[b]{0.225in}{$\myC_{2}$
\vspace*{1.1in}} 
 & \BTwoAAlpha & \BTwoBBeta\\
\hline
\parbox[b]{0.225in}{$\myG_{2}$
\vspace*{2.5in}} 
 & \GTwoAAlpha & \GTwoBBeta\\
\hline
\end{tabular}
\end{center}
\end{figure}

\newpage
\begin{figure}[ht]
\begin{center}
\setlength{\unitlength}{1.1cm}
\begin{picture}(14,18.5)
\put(0,18.1){\small \FundLatticeIdealsFigure: Elements of $\mathscr{G}$-fundamental DCDL's as down-sets from $\mathscr{G}$-fundamental compression posets.}
\put(1.25,17.5){\parbox{5in}{\scriptsize \begin{center}
(Each down-set from the fundamental poset is identified by the indices of its maximal vertices.\\
For example, $\langle 6,7 \rangle$ in $L_{\mytinyG_{2}}(\omega_{\beta})$ denotes the down-set $\{v_{6}, v_{7}, v_{8}, v_{9}, v_{10}\}$ from $P_{\mytinyG_{2}}(\omega_{\beta})$.) 
\end{center}
}}
\put(0,-0.5){
\begin{picture}(14,18.5)
\put(0,0){\line(0,1){17.5}}
\put(8,0){\line(0,1){17.5}}
\put(14,0){\line(0,1){17.5}}
\put(0,0){\line(1,0){14}}
\put(0,13){\line(1,0){8}}
\put(8,12){\line(1,0){6}}
\put(0,17.5){\line(1,0){14}}
\put(1,16.5){\fbox{\Large $\myA_{1} \oplus \myA_{1}$}}
\put(0.5,13){\AoneAlphaIdeals}
\put(3.5,13){\AoneBetaIdeals}
\put(8.5,16.5){\fbox{\Large $\myA_{2}$}}
\put(8,12){\AtwoAlphaIdeals}
\put(11,12){\AtwoBetaIdeals}
\put(8.5,10.5){\fbox{\Large $\myC_{2}$}}
\put(8,4.5){\BtwoAlphaIdeals}
\put(11,3.5){\BtwoBetaIdeals}
\put(1,11.5){\fbox{\Large $\myG_{2}$}}
\put(0,2){\GtwoAlphaIdeals}
\put(2,0){\GtwoBetaIdeals}
\end{picture}
}
\end{picture}
\end{center} 
\end{figure}

\begin{figure}[ht]
\begin{center}
\setlength{\unitlength}{1.1cm}
\begin{picture}(14,18.5)
\put(2.8,18.3){\small \WeightFigure: Weights of elements of $\mathscr{G}$-fundamental DCDL's.}
\put(1,17.9){\scriptsize {\large (}The weight of a given lattice element $\xelt$ is indicated by the integer pair $\big(\rho_{\alpha}(\xelt)-\delta_{\alpha}(\xelt),\rho_{\alpha}(\xelt)-\delta_{\alpha}(\xelt)\big)$.{\large )}}
\put(0,0){\line(0,1){17.5}}
\put(8,0){\line(0,1){17.5}}
\put(14,0){\line(0,1){17.5}}
\put(0,0){\line(1,0){14}}
\put(0,13){\line(1,0){8}}
\put(8,12){\line(1,0){6}}
\put(0,17.5){\line(1,0){14}}
\put(1,16.5){\fbox{\Large $\myA_{1} \oplus \myA_{1}$}}
\put(0.5,13){\AoneAlphaWeights}
\put(3.5,13){\AoneBetaWeights}
\put(8.5,16.5){\fbox{\Large $\myA_{2}$}}
\put(8,12){\AtwoAlphaWeights}
\put(11,12){\AtwoBetaWeights}
\put(8.5,10.5){\fbox{\Large $\myC_{2}$}}
\put(8,4.5){\BtwoAlphaWeights}
\put(11,3.5){\BtwoBetaWeights}
\put(1,11.5){\fbox{\Large $\myG_{2}$}}
\put(0,2){\GtwoAlphaWeights}
\put(2,0){\GtwoBetaWeights}
\end{picture}
\end{center} 
\end{figure}

\clearpage
{\bf [\S \GridPosetSection.4:\! $\mathscr{G}$-fundamental and $\mathscr{G}$-semistandard compression posets in two colors.]} 
Let $\mathscr{G}$ be a rank two Coxeter--Dynkin posy, so $\mathscr{G}$ is one of $\myA_{1} \oplus \myA_{1}$, $\myA_{2}$, $\myC_{2}$, or $\myG_{2}$, with the indices of the nodes $\gamma_{1}$ and $\gamma_{2}$ corresponding to $\alpha$ and $\beta$ respectively. 
We define the $\mathscr{G}$-{\em fundamental (compression) posets} $P_{\mathscr{G}}(1,0)$ and $P_{\mathscr{G}}(0,1)$ to be the  \dichromatic grid posets of \FundPosets. 
Now let $\lambda = (a,b)$ be a pair of nonnegative integers.  
There are exactly two possible ways that a \dichromatic grid poset $P$ with the max property can decompose as $P_{1} \triangleleft P_{2} \triangleleft \cdots \triangleleft P_{a+b}$ with $a$ of the $P_{i}$'s vertex-color isomorphic to  $P_{\mathscr{G}}(1,0)$ and the remaining $P_{i}$'s vertex-color isomorphic to $P_{\mathscr{G}}(0,1)$: we will either have $P_{i}$ isomorphic to $P_{\mathscr{G}}(0,1)$ for $1 \leq i \leq b$ and isomorphic to $P_{\mathscr{G}}(1,0)$ for $1+b \leq i \leq a+b$ (in which case we set $P_{\mathscr{G}}^{\beta\alpha}(\lambda) := P$), or we will have $P_{i}$ isomorphic to $P_{\mathscr{G}}(1,0)$ for $1 \leq i \leq a$ and isomorphic to $P_{\mathscr{G}}(0,1)$ for $a+1 \leq i \leq a+b$ (in which 
case we set $P_{\mathscr{G}}^{\alpha\beta}(\lambda) := P$).  
Note that $P_{\mathscr{G}}^{\beta\alpha}(1,0) = P_{\mathscr{G}}^{\alpha\beta}(1,0) = P_{\mathscr{G}}(1,0)$, and   $P_{\mathscr{G}}^{\beta\alpha}(0,1) = P_{\mathscr{G}}^{\alpha\beta}(0,1) = P_{\mathscr{G}}(0,1)$. When $a = b = 0$, then $P_{\mathscr{G}}^{\beta\alpha}(\lambda)$ and $P_{\mathscr{G}}^{\alpha\beta}(\lambda)$ are the empty set. 
We call $P_{\mathscr{G}}^{\beta\alpha}(\lambda)$ and $P_{\mathscr{G}}^{\alpha\beta}(\lambda)$ the $\mathscr{G}$-{\em semistandard (compression) posets} associated to $\lambda$. 
For each two-node Coxeter--Dynkin posy $\mathscr{G}$, the $\mathscr{G}$-semistandard posets $P_{\mathscr{G}}^{\beta\alpha}(2,2)$ and $P_{\mathscr{G}}^{\alpha\beta}(2,2)$ are depicted in \GridPosetsFigureList. 
One of the main goals of this section is to completely combinatorially characterize the $\mathscr{G}$-fundamental and $\mathscr{G}$-semistandard compression posets. 

Continuing with $\mathscr{G}$ as a two-node Coxeter--Dynkin posy, the $\mathscr{G}$-{\em fundamental diamond-colored distributive lattices} are the lattices $L_{\mathscr{G}}(1,0) := \Jcolor\big(P_{\mathscr{G}}(1,0)\big)$ and $L_{\mathscr{G}}(0,1) := \Jcolor\big(P_{\mathscr{G}}(0,1)\big)$. 
The $\mathscr{G}$-{\em semistandard DCDL's} are $L_{\mathscr{G}}^{\beta\alpha}(\lambda) := \Jcolor\big(P_{\mathscr{G}}^{\beta\alpha}(\lambda)\big)$ and $L_{\mathscr{G}}^{\alpha\beta}(\lambda) := \Jcolor\big(P_{\mathscr{G}}^{\alpha\beta}(\lambda)\big)$ for all possible nonnegative integer pairs $\lambda = (a,b)$. 
Momentarily set $L_{1} := L_{\mathscr{G}}^{\beta\alpha}(\lambda)$ and $L_{2} := L_{\mathscr{G}}^{\alpha\beta}(\lambda)$. 
Observe, then, that $\lambda$ is the weight of their maximal elements, i.e.\ $wt(\mymax(L_{1})) = \lambda = wt(\mymax(L_{2}))$. 
See \FundFigures\ for depictions of the $\mathscr{G}$-fundamental DCDL's and of the weights of all lattice elements. 

The following technical lemma is used in the proofs of \ClassTheorems, which give combinatorial characterizations of $\mathscr{G}$-fundamental and $\mathscr{G}$-semistandard posets respectively. 

\noindent 
{\bf \TechnicalProposition}\ \ {\sl Let $L = \Jcolor(P)$ be the diamond-colored distributive lattice of order ideals from a two-color grid poset} $(P,\leq_{_{P}}, \mychain: P \rightarrow [1,m]_{\mathbb{Z}}, \mycolor: P \rightarrow \{\alpha,\beta\})$, {\sl so in this context $I := \{\alpha,\beta\}$ is our color palette.} 
{\sl Suppose $L$ is $\mathscr{G}$-structured, where $\mathscr{G} = (\Gamma_{I},M_{I \times I})$ is one of the Coxeter--Dynkin flowers} $\myA_{2}$, $\myC_{2}$, {\sl or} $\myG_{2}$. 
{\sl In these three cases, when $M=M_{I \times I}$ is displayed as a $2 \times 2$ matrix, row $\alpha$ is the top row and row $\beta$ is the bottom row.} 
{\sl Let $z_{1}$ be a face vertex of $P$. 
Then there exists a unique subset $P_{2}$ of $P$ that contains $z_{1}$ and that, given the induced two-color grid poset structure from $P$, is isomorphic to the appropriate figure as identified via cases below. 
Moreover, if $z_{1}$ is the rightmost maximal element of $P$ and if $P_{1} := P \setminus P_{2}$ is nonempty, then $P = P_{1} \triangleleft P_{2}$.}

{\sl  
\noindent 
{\sc case 1:} Suppose $\displaystyle M = 
\left(\begin{array}{cc}2 & -1\\ -1 & 2\end{array}\right)$.  If} 
$\mycolor(z_{1}) = \alpha$ {\sl (respectively} $\mycolor(z_{1}) = 
\beta${\sl ), then there is a 
unique subset $P_{2} = \{z_{1},z_{2}\}$ of $P$ which, given the induced 
order from $P$, is a \dichromatic grid subposet isomorphic 
under the correspondence $u_{i} \leftrightarrow z_{i}$ to the 
\dichromatic grid poset of \TechPropAalphaFigure\ (respectively 
\TechPropAbetaFigure).} 
\begin{center}
\ATwoFundamentalA
\hspace*{1in}
\ATwoFundamentalB
\end{center}

{\sl 
\noindent 
{\sc case 2:} Suppose $\displaystyle M = 
\left(\begin{array}{cc}2 & -1\\ -2 & 2\end{array}\right)$.  If} 
$\mycolor(z_{1}) = \alpha$ {\sl (respectively} $\mycolor(z_{1}) = 
\beta${\sl ), then there is a 
unique subset $P_{2} = \{z_{1},z_{2},z_{3}\}$ (respectively 
$P_{2} = \{z_{1},z_{2},z_{3},z_{4}\}$) of $P$ which, given the induced 
order from $P$, is a \dichromatic grid subposet isomorphic 
under the correspondence $u_{i} \leftrightarrow z_{i}$ to the 
\dichromatic grid poset of \TechPropBalphaFigure\ (respectively 
\TechPropBbetaFigure).} 
\begin{center}
\BTwoFundamentalA
\hspace*{1in}
\BTwoFundamentalB
\end{center}

{\sl  
\noindent 
{\sc case 3:} Suppose $\displaystyle M = 
\left(\begin{array}{cc}2 & -1\\ -3 & 2\end{array}\right)$.  If} 
$\mycolor(z_{1}) = \alpha$ {\sl (respectively} $\mycolor(z_{1}) = 
\beta${\sl ), then there is a 
unique subset $P_{2} = \{z_{1},z_{2},z_{3},z_{4},z_{5},z_{6}\}$ (respectively 
$P_{2} = 
\{z_{1}, z_{2}, z_{3}, z_{4}, z_{5}, z_{6}, z_{7}, z_{8}, z_{9}, 
z_{10}\}$) of $P$ which, given the induced 
order from $P$, is a \dichromatic grid subposet isomorphic 
under the correspondence $u_{i} \leftrightarrow z_{i}$ to the 
\dichromatic grid poset of \TechPropGalphaFigure\ (respectively 
\TechPropGbetaFigure).} 
\begin{center}
\GTwoFundamentalA
\hspace*{1in}
\GTwoFundamentalB
\end{center}

{\em Proof.} 
Below is our proof for Case 2. 
Our proof for Case 1  
employs similar but shorter arguments. 
A similarly-structured argument for Case 3 is more involved.  
For Case 2, we 
use the following sequence of Claims 1 through 11 
to justify the conclusion. 
For $\telt \in L$ and $\gamma \in \{\alpha,\beta\}$, we make use the \S \SubstructureSection.5 notation `$D_{\gamma}(\telt)$' and `$A_{\gamma}(\telt)$' (see the paragraphs prior to \JCompTheorem). 
Let $r := \mychain(z_{1})$.  
Let $\uelt$ be the down-set from $P$ 
generated by $z_{1}$.  Let $\lwt(\uelt)$ be given by the pair $(a,b)$, where $a$ and $b$ are 
integers.\\ 
\underline{\em Claim 1:} If $\mycolor(z_{1}) = \alpha$, 
there exists a unique $z_{2}$ in $P$ with $\mychain(z_{2}) = r+1$
such that $z_{2} 
\rightarrow z_{1}$ is a covering relation in $P$. 
If $\mycolor(z_{1}) = \beta$, 
there exists a unique ordered pair $(z_{2},z_{3})$ of elements in 
$P$ with $\mychain(z_{2}) = \mychain(z_{3}) = r+1$
such that $z_{3} \rightarrow z_{2} \rightarrow z_{1}$ are 
covering relations in $P$. 
\underline{\em Proof.} Suppose $\mycolor(z_{1}) = \alpha$.  Let $\telt := 
\uelt \setminus \{z_{1}\}$. Then $\lwt(\telt) = \lwt(\uelt) - 
(2,-1) = (a-2,b+1)$.  Since $z_{1}$ is not covered by any color 
$\beta$ vertices in $P$, it follows that $A_{\beta}(\telt) = 
A_{\beta}(\uelt)$.  It follows that there exists a unique $z_{2}$ in 
$P$ such that $z_{2} \not\in D_{\beta}(\uelt)$ and $D_{\beta}(\telt) = 
D_{\beta}(\uelt) \cup \{z_{2}\}$.  Since $z_{2} \not\in 
D_{\beta}(\uelt)$, it must be the case that $z_{2} 
\rightarrow z_{1}$. Since $\mycolor(z_{2}) = \beta$, it follows 
that $\mychain(z_{2}) = r+1$.  Since $z_{1}$ can cover at most 
one element of 
$\mathcal{C}_{r+1}$, it follows that the element $z_{2}$ is the only 
element of $\mathcal{C}_{r+1}$ covered by $z_{1}$.  Use a similar 
argument if 
$\mycolor(z_{1}) = \beta$.  In this case 
$\lwt(\telt) = \lwt(\uelt) - (-2,2) = (a+2,b-2)$. 
As before observe that $A_{\alpha}(\telt) = 
A_{\alpha}(\uelt)$.  It follows that there exists a unique pair 
$\{z_{2},z_{3}\}$ in 
$P$ such that $z_{2}$ and $z_{3}$ are not in 
$D_{\alpha}(\uelt)$ and $D_{\alpha}(\telt) = 
D_{\alpha}(\uelt) \cup \{z_{2},z_{3}\}$.  Since $z_{2}$ and $z_{3}$ are not in 
$D_{\alpha}(\uelt)$, it must be the case that one of them is covered 
by $z_{1}$. Without loss of generality, call this element $z_{2}$.  
One can see then that $z_{3} \rightarrow z_{2}$. So 
$z_{3} \rightarrow z_{2} \rightarrow z_{1}$.  Since $z_{1}$ can cover at most 
one element of 
$\mathcal{C}_{r+1}$, it follows that the element $z_{2}$ is the only 
element of $\mathcal{C}_{r+1}$ covered by $z_{1}$.  And since 
$z_{2}$ can cover at most 
one element of 
$\mathcal{C}_{r+1}$, it follows that the element $z_{3}$ is the only 
element of $\mathcal{C}_{r+1}$ covered by $z_{2}$. \\  
\underline{\em Claim 2:} If $\mycolor(z_{1}) = \alpha$, then 
$z_{2}$ is a maximal element in $\uelt \setminus \{z_{1}\}$. 
If $\mycolor(z_{1}) = \beta$, then 
$z_{2}$ is a maximal element in $\uelt \setminus \{z_{1}\}$ and 
$z_{3}$ is a maximal element in $\uelt \setminus \{z_{1},z_{2}\}$. 
\underline{\em Proof.} Suppose 
$\mycolor(z_{1}) = \alpha$. If $z_{2}$ is not maximal in 
$\uelt \setminus \{z_{1}\}$, then $z_{2} \rightarrow z$ for some $z 
\in \uelt \setminus \{z_{1}\}$.  Since $z \not= 
z_{1}$, then $\mychain(z) = \mychain(z_{2})$.  However, 
$z \in \uelt$ implies that $z <_{_{P}} z_{1}$, so that $z_{2} <_{_{P}} 
z <_{_{P}} z_{1}$.  This contradicts the fact that $z_{2} \rightarrow 
z_{1}$.  Suppose 
$\mycolor(z_{1}) = \alpha$. The argument that $z_{2}$ is maximal 
in $\uelt \setminus \{z_{1}\}$ is the same as before.  If $z_{3}$ is 
not maximal in $\uelt \setminus \{z_{1},z_{2}\}$, then $z_{3} 
\rightarrow z$ is a covering relation in $P$ for some $z$ in 
$\uelt \setminus \{z_{1},z_{2}\}$. Since $z_{3} \rightarrow z_{2}$, 
it follows that $\mychain(z) = \mychain(z_{1})$.  But 
then $\telt \setminus D_{\beta}(\telt)$ cannot be a down-set, 
since $z$ is an element in 
$\telt \setminus D_{\beta}(\telt)$ but $z_{3}$ is not. \\
\underline{\em Claim 3:} If $\mycolor(z_{1}) = \alpha$, then 
there does not exist $u 
\in \mathcal{C}_{r+1}$ such that $u \rightarrow z_{2}$ and $u$ is 
maximal in $\uelt \setminus \{z_{1},z_{2}\}$. If 
$\mycolor(z_{1}) = \beta$, then 
there does not exist $u 
\in \mathcal{C}_{r+1}$ such that $u \rightarrow z_{3}$ and $u$ is 
maximal in $\uelt \setminus \{z_{1},z_{2},z_{3}\}$. 
\underline{\em Proof.} First consider the case 
$\mycolor(z_{1}) = \alpha$. 
If such an element $u$ exists, 
then $u \in D_{\beta}(\telt)$.  But $u \not\in D_{\beta}(\uelt)$ 
because $z_{2} \not\in D_{\beta}(\uelt)$.  Since $D_{\beta}(\telt) = 
D_{\beta}(\uelt) \cup \{z_{2}\}$, then it cannot be the case that 
$u \in D_{\beta}(\telt)$, so we have a contradiction.  Use a similar 
argument in the case that $\mycolor(z_{1}) = \beta$. \\
\underline{\em Claim 4:} If $\mycolor(z_{1}) = \alpha$, 
there exists a unique $z_{3}$ in $P$ with $\mychain(z_{3}) = r+2$
such that $z_{3} 
\rightarrow z_{2}$ is a covering relation in $P$. 
If $\mycolor(z_{1}) = \beta$, 
there exists a unique $z_{4}$ in 
$P$ with $\mychain(z_{4}) = r+2$
and $z_{4} \rightarrow z_{3}$ is a covering relation in $P$. 
\underline{\em Proof.} Suppose $\mycolor(z_{1}) = \alpha$.  
Let $\selt := 
\uelt \setminus \{z_{1},z_{2}\}$. By Claim 2, $\selt$ is a down-set from $P$.  Then $\lwt(\selt) = \lwt(\uelt) - 
(-2,2) = (a-2,b+1) - (-2,2) = (a,b-1)$.  Observe that $A_{\alpha}(\selt) = 
A_{\alpha}(\uelt)$.  It follows that there exists a unique $z_{3}$ in 
$P$ such that $z_{3} \not\in D_{\alpha}(\uelt)$ and $D_{\alpha}(\selt) = 
(D_{\alpha}(\uelt) \setminus \{z_{1}\}) \cup \{z_{3}\}$.  Since $z_{3} \not\in 
D_{\alpha}(\uelt)$ but $z_{3} \in \selt = \uelt \setminus \{z_{1},z_{2}\}$, 
it must be the case that $z_{3} 
\rightarrow z_{2}$. Since $\mycolor(z_{3}) = \alpha$, it follows 
that $\mychain(z_{3}) = r+2$.  
Since $z_{2}$ can cover at most one element of 
$\mathcal{C}_{r+2}$, it follows that the element $z_{3}$ is the only 
element of $\mathcal{C}_{r+2}$ covered by $z_{2}$. In the case that 
$\mycolor(z_{1}) = \beta$, use a similar argument to locate and 
characterize the element $z_{4}$. \\
\underline{\em Claim 5:} If $\mycolor(z_{1}) = \alpha$, then 
$z_{3}$ is maximal in $\uelt 
\setminus \{z_{1},z_{2}\}$. If $\mycolor(z_{1}) = \beta$, then 
$z_{4}$ is maximal in $\uelt 
\setminus \{z_{1},z_{2},z_{3}\}$. 
\underline{\em Proof.} 
For both cases ($\mycolor(z_{1}) = \alpha$ or 
$\mycolor(z_{1}) = \beta$), the argument is entirely similar to 
the $\mycolor(z_{1}) = \alpha$ case of Claim 2.\\
\underline{\em Claim 6:} If $\mycolor(z_{1}) = \alpha$, then 
there does not exist $u 
\in \mathcal{C}_{r+2}$ such that $u \rightarrow z_{3}$ and $u$ is 
maximal in $\uelt \setminus \{z_{1},z_{2},z_{3}\}$. 
If $\mycolor(z_{1}) = \beta$, then 
there does not exist $u 
\in \mathcal{C}_{r+2}$ such that $u \rightarrow z_{4}$ and $u$ is 
maximal in $\uelt \setminus \{z_{1},z_{2},z_{3},z_{4}\}$. 
\underline{\em Proof.}  Entirely similar to the proof of Claim 3. \\
\underline{\em Claim 7:} If $\mycolor(z_{1}) = \alpha$ (respectively 
$\mycolor(z_{1}) = \beta$), then there does 
not exist $z$ in $\mathcal{C}_{r+3}$ such that $z \rightarrow z_{3}$ 
(respectively $z \rightarrow z_{4}$). 
\underline{\em Proof.} Suppose 
$\mycolor(z_{1}) = \alpha$.  Let $\relt := \uelt \setminus 
\{z_{1},z_{2},z_{3}\}$.  Then $\relt$ is a down-set from $P$. 
Observe that $A_{\beta}(\relt) = A_{\beta}(\uelt)$, 
$D_{\beta}(\relt) \supseteq D_{\beta}(\uelt)$, and $\lwt(\relt) = 
\lwt(\selt) - (2,-1) = (a,b-1) - (2,-1) = (a-2,b)$, where $\selt$ 
is as in the proof of Claim 4. 
If such a $z$ exists, 
then $z \in D_{\beta}(\relt)$ and $z \not\in D_{\beta}(\uelt)$, in which 
case $D_{\beta}(\relt) \supset D_{\beta}(\uelt)$, hence 
$\rho_{\beta}(\relt) > \rho_{\beta}(\uelt)$, hence 
$\rho_{\beta}(\relt) - (l_{\beta}(\relt) - \rho_{\beta}(\relt)) 
> \rho_{\beta}(\uelt) - (l_{\beta}(\uelt) - \rho_{\beta}(\uelt)) = b$, 
which contradicts the fact that $\lwt(\relt) = (a-2,b)$. Use a 
similar argument in the case that 
$\mycolor(z_{1}) = \beta$. \\
\underline{\em Claim 8:} If $\mycolor(z_{1}) = \beta$, then for all $z \in 
\mathcal{C}_{i}$ for $r+2 < i \leq m$, $z \not<_{_{P}} z_{4}$. 
\underline{\em Proof.} Suppose not.  Then we may take $z$ to be the element of 
$\mathcal{C}_{r+3}$ such that if $w \in \mathcal{C}_{r+3}$ with $z 
\leq_{_{P}} w <_{_{P}} z_{4}$, then $z = w$.  
Note that $\mycolor(z) = \alpha$.  
Let $S := \{\xelt \in L | \xelt \subset \uelt \mbox{ and } 
z \not\in \xelt\}$.  Since $\xelt \cup 
\yelt$ and $\xelt \cap \yelt$ are in $S$ if $\xelt, \yelt 
\in S$, it follows that $S$ is a distributive sublattice of $L$.  Let 
$\welt$ denote the unique maximal element of $S$. One can see that 
$\zelt := \welt \cup \{z\}$ is a down-set from $P$, and hence an element of 
$L$. (Let $v \in \zelt$ and suppose $u <_{_{P}} v$.  If $v \in \welt$, 
then since $\welt$ is a down-set, we'll have $u \in \welt \subset 
\zelt$.  Now suppose $v = z$.  If $u \not\in \welt$, then consider 
$\welt \cup \{w \in P | w \leq_{_{P}} u\}$.  This is a down-set 
from $P$ (since it is the union of two down-sets) that does not 
contain $z$; $\{w \in P | w \leq_{_{P}} u\}$ is a subset of $\uelt$ 
since $u <_{_{P}} z <_{_{P}} z_{4} \in \uelt$.  
So $\welt \cup \{w \in P | w \leq_{_{P}} u\}$ is in 
$S$.  Moreover, $\welt <_{_{L}} \welt \cup \{w \in P | w \leq_{_{P}} 
u\}$, contradicting the maximality of $\welt$ in $S$.  Therefore $u 
\in \welt \subset \zelt$.  It now follows that $\zelt$ is a down-set from $P$.)  Then $\lwt(\welt) + (2,-1) = \lwt(\zelt)$.  Write 
$z \rightarrow y_{1} \rightarrow \cdots \rightarrow y_{q} = z_{4}$; 
here $q \geq 1$ and for $1 \leq i \leq q$, $y_{i} \in 
\mathcal{C}_{r+2}$. 
Note that if $y \rightarrow y_{1}$ in $\mathcal{C}_{r+2}$, then $y 
\in \welt$.  Now let $x$ be the element in $\mathcal{C}_{r+4}$ such that 
$x \leq_{_{P}} w$ for all $w 
\in \mathcal{C}_{r+4}$ satisfying: $w <_{_{P}} z$, and $w \not\leq_{_{P}} y$ 
for all $y$ in $\mathcal{C}_{r+3}$ with $y <_{_{P}} z$.  
Let $x_{1},\ldots,x_{p}$ be in 
$\mathcal{C}_{r+4}$ so that $x = x_{p} \rightarrow \cdots \rightarrow 
x_{1} \rightarrow z$; we have $p \geq 0$, with $p = 0$ if there is no 
such element $x$ in $\mathcal{C}_{r+4}$.  Note that for $1 \leq i \leq 
p$, $x_{i} \in \welt$.  Then $D_{\beta}(\zelt) = D_{\beta}(\welt) 
\setminus \{x_{1},\ldots,x_{p}\}$, and $A_{\beta}(\welt) = 
A_{\beta}(\zelt) \setminus \{y_{1},\ldots,y_{q}\}$.  So 
$\rho_{\beta}(\welt) = \rho_{\beta}(\zelt) + p$ and 
$l_{\beta}(\welt) - \rho_{\beta}(\welt) = l_{\beta}(\zelt) - 
\rho_{\beta}(\zelt) - q$. Then 
\begin{eqnarray*}
1 & = &2\rho_{\beta}(\welt) - l_{\beta}(\welt) - 
(2\rho_{\beta}(\zelt) - l_{\beta}(\zelt)\, )\\ 
 & = & \rho_{\beta}(\welt) - \rho_{\beta}(\zelt) - 
[l_{\beta}(\welt) - \rho_{\beta}(\welt)  
- (l_{\beta}(\zelt) - \rho_{\beta}(\zelt)\, )\, ]\\ 
 & = & p + q.
\end{eqnarray*} Therefore $p = 0$ and $q = 1$.  In particular, $z 
\rightarrow z_{4}$, which violates Claim 7.\\
\underline{\em Claim 9:} If $\mycolor(z_{1}) = \alpha$, then for all $z 
\in \mathcal{C}_{i}$ for $r+2 < i \leq m$, $z \not<_{_{P}} 
z_{3}$. 
\underline{\em Proof.} Totally order the face vertices of $P$ of color $\alpha$ 
as $f_{1}, f_{2},\ldots$ according to the following rule: for $i < j$, 
either $\mychain(f_{i}) < \mychain(f_{j})$ or 
$\mychain(f_{i}) = \mychain(f_{j})$ and $f_{i} <_{_{P}} 
f_{j}$.  We use 
induction on the face vertices of color $\alpha$.  Suppose the result 
is true if $z_{1} \in \{f_{1},\ldots,f_{k-1}\}$.  Now suppose $z_{1} = 
f_{k}$.  We use a contradiction argument similar to the proof 
of Claim 8.  We will take $z$ to be the element of 
$\mathcal{C}_{r+3}$ such that if $w \in \mathcal{C}_{r+3}$ with $z 
\leq_{_{P}} w <_{_{P}} z_{3}$, then $z = w$.  
Note that $\mycolor(z) = \beta$.  
Let $S := \{\xelt \in L | \xelt \subset \uelt \mbox{ and } 
z \not\in \xelt\}$.  Since $\xelt \cup 
\yelt$ and $\xelt \cap \yelt$ are in $S$ if $\xelt, \yelt 
\in S$, it follows that $S$ is a distributive sublattice of $L$.  Let 
$\welt$ denote the unique maximal element of $S$. As before, one can see that 
$\zelt := \welt \cup \{z\}$ is a down-set from $P$, and hence an element of 
$L$. Write 
$z \rightarrow y_{1} \rightarrow \cdots \rightarrow y_{q} = z_{3}$; 
here $q \geq 1$ and for $1 \leq i \leq q$, $y_{i} \in 
\mathcal{C}_{r+2}$. 
Note that if $y \rightarrow y_{1}$ in $\mathcal{C}_{r+2}$, then $y 
\in \welt$.  Now let $x$ be the element in $\mathcal{C}_{r+4}$ such that 
$x \leq_{_{P}} w$ for all $w 
\in \mathcal{C}_{r+4}$ satisfying: $w <_{_{P}} z$, and $w \not\leq_{_{P}} y$ 
for all $y$ in $\mathcal{C}_{r+3}$ with $y <_{_{P}} z$.  
Let $x_{1},\ldots,x_{p}$ be in 
$\mathcal{C}_{r+4}$ so that $x = x_{p} \rightarrow \cdots \rightarrow 
x_{1} \rightarrow z$; we have $p \geq 0$, with $p = 0$ if there is no 
such element $x$ in $\mathcal{C}_{r+4}$.  Note that for $1 \leq i \leq 
p$, $x_{i} \in \welt$.  Then $D_{\alpha}(\zelt) = D_{\alpha}(\welt) 
\setminus \{x_{1},\ldots,x_{p}\}$, and $A_{\alpha}(\welt) = 
A_{\alpha}(\zelt) \setminus \{y_{1},\ldots,y_{q}\}$.  So 
$\rho_{\alpha}(\welt) = \rho_{\alpha}(\zelt) + p$ and 
$l_{\alpha}(\welt) - \rho_{\alpha}(\welt) = l_{\alpha}(\zelt) - 
\rho_{\alpha}(\zelt) - q$. Since $\welt \myarrow{\beta} \zelt$, it 
follows that $\lwt(\welt) + (-2,2) = \lwt(\zelt)$.  Then
\begin{eqnarray*}
2 & = &2\rho_{\alpha}(\welt) - l_{\alpha}(\welt) - 
(2\rho_{\alpha}(\zelt) - l_{\alpha}(\zelt)\, )\\ 
 & = & \rho_{\alpha}(\welt) - \rho_{\alpha}(\zelt) - 
[l_{\alpha}(\welt) - \rho_{\alpha}(\welt)  
- (l_{\alpha}(\zelt) - \rho_{\alpha}(\zelt)\, )\, ]\\ 
 & = & p + q.
\end{eqnarray*} If $q = 1$, then $z \rightarrow z_{3}$, which 
violates Claim 7.  So $q = 2$ and $p = 0$.  Since $y_{1}$ is not 
maximal in $\uelt \setminus \{z_{1},z_{2},z_{3}\}$, then $y_{1} 
\rightarrow u$ for some $u <_{_{P}} z_{2}$ with $\mychain(u) = 
\mychain(z_{2}) = r+1$. So now let $y_{2} \in 
\mathcal{C}_{r+1}$ with $u \leq_{_{P}} y_{2} \rightarrow z_{2}$.  
Similarly find $y_{3} \rightarrow z_{1} = f_{k}$ in $\mathcal{C}_{r}$ with 
$y_{2} <_{_{P}} y_{3}$.  If $y_{3}$ is a face vertex, then $y_{3} = 
f_{k-1}$.  But then we may apply Claims 1 through 6 to the face 
vertex $z_{1}' := f_{k-1}$ to obtain $\{z_{1}' = f_{k-1}, z_{2}', 
z_{3}'\}$.  But then $z <_{_{P}} z_{3}'$, violating the inductive 
hypothesis.  So $y_{3}$ is not a face vertex.  Then 
$y_{3} \rightarrow v$ for some $v$ in $\mathcal{C}_{r-1}$.  Now $v 
\leq_{_{P}} f$ for some face vertex $f$ in $P$.  So  
$\mychain(f) < \mychain(z_{1})$.  If 
$\mycolor(f) = \alpha$, then $f = f_{i}$ for some $1 \leq i 
\leq k-1$, and we can obtain a contradiction in the 
same way as before.  So it must be the case 
that $f$ has color $\beta$.  But then we may apply Claims 1 through 6 
to construct $\{f = z_{1}'', z_{2}'', z_{3}'', z_{4}''\}$.  It follows 
that $z <_{_{P}} z_{4}''$, violating Claim 8.  So no such $z$ exists. \\
\underline{\em Claim 10:} Now suppose $z_{1}$ is the rightmost maximal element 
of $P$.  If $\mycolor(z_{1}) = \alpha$ (respectively 
$\mycolor(z_{1}) = \beta$), then let $P_{2} := 
\{z_{1},z_{2},z_{3}\}$ (respectively $P_{2} := 
\{z_{1},z_{2},z_{3},z_{4}\}$), and let $P_{1} := P 
\setminus P_{2}$.  If $u \rightarrow v$ is a covering relation in $P$ 
and $u \in P_{2}$, then $v \in P_{2}$. 
\underline{\em Proof:} Use Claims 2 and 
5.  For example, since $z_{3}$ is a maximal element in $\uelt 
\setminus \{z_{1},z_{2}\}$, then $z_{3}$ can only be covered by $z_{1}$ 
or $z_{2}$.\\
\underline{\em Claim 11:} Keep the hypothesis that $z_{1}$ 
is the rightmost maximal element 
of $P$. Regard $P_{1}$ as a subposet of $P$ in the induced 
order.  If $u \in P_{1}$ is a maximal element of $P_{1}$, 
then $\mychain(u) \leq r$. 
\underline{\em Proof.} Suppose 
$\mycolor(z_{1}) = \alpha$ (respectively 
$\mycolor(z_{1}) = \beta$).  If $u \in P_{1}$ is 
maximal in $P_{1}$ and $\mychain(u) > r$, then by Claim 9 
(respectively Claim 8), 
$u \in \mathcal{C}_{r+1} \cup \mathcal{C}_{r+2} \cup 
\mathcal{C}_{r+3}$. Then $u \rightarrow 
z_{2}$ or $u \rightarrow z_{3}$ (respectively $u \rightarrow 
z_{3}$ or $u \rightarrow z_{4}$) is a covering relation in $P$.  
But then $u$ is maximal in 
$\uelt \setminus \{z_{1},z_{2}\}$ or  
$\uelt \setminus \{z_{1},z_{2},z_{3}\}$  
(respectively 
$\uelt \setminus \{z_{1},z_{2},z_{3}\}$ or  
$\uelt \setminus \{z_{1},z_{2},z_{3},z_{4}\}$), 
which violates one of Claims 3 or 6.\\
\underline{\em Conclusion:} Suppose 
$\mycolor(z_{1}) = \alpha$ (respectively 
$\mycolor(z_{1}) = \beta$).  
Regard $P_{2}$ as a subposet of $P$ in the induced 
order.  By Claims 1 and 4, one can see that $P_{2}$ is an 
$\alpha\, \beta$ grid subposet of $P$ isomorphic to the $\alpha\, 
\beta$ grid poset of \TechPropBalphaFigure\ (respectively \TechPropBbetaFigure).  
Claims 1 and 4 also 
demonstrate that this is the only possible $\alpha\, \beta$ grid 
subposet of $P$ that contains $z_{1}$ and is isomorphic to the 
$\alpha\, 
\beta$ grid poset of \TechPropBalphaFigure\ (respectively \TechPropBbetaFigure).  
Now suppose that $z_{1}$ is the rightmost maximal element of $P$ and 
that $P_{1}$ is nonempty.  
If $u$ is any minimal element of the poset $P_{1}$, then 
$\mychain(u) \leq r+3$ by Claim 9 (respectively Claim 8); 
since $z_{3}$ (respectively $z_{4}$) is the unique 
minimal element of $P_{2}$ and $\mychain(z_{3}) = r+2$ 
(respectively $\mychain(z_{4}) = r+2$), 
it follows that $\mychain(u) \leq \mychain(z_{3})$ 
(respectively $\mychain(u) \leq \mychain(z_{4})$). This 
observation, together with Claims 10 and 11, shows that $P = P_{1} 
\triangleleft P_{2}$.\hfill\QED

{\bf [\S \GridPosetSection.5:\! A combinatorial characterization of rank two fundamental posets.]}

\newcommand{\AOneAOneAlphaMars}{
\setlength{\unitlength}{1cm}
\begin{picture}(0.5,0.5)
\put(0,0){\circle*{0.15}} 
\put(0.2,0){\footnotesize $\alpha$}
\end{picture}
}
\newcommand{\AOneAOneBetaMars}{
\setlength{\unitlength}{1cm}
\begin{picture}(0.5,0.5)
\put(0,0){\circle*{0.15}} 
\put(0.2,0){\footnotesize $\beta$}
\end{picture}
}
\newcommand{\ATwoAlphaMars}{
\setlength{\unitlength}{1cm}
\begin{picture}(2,2)
\put(1,0){\circle*{0.15}} 
\put(1.2,-0.1){\footnotesize $\beta$}
\put(0,1){\circle*{0.15}} 
\put(0.2,0.9){\footnotesize $\alpha$}
\put(1,0){\line(-1,1){1}} 
\end{picture}
}
\newcommand{\ATwoBetaMars}{
\setlength{\unitlength}{1cm}
\begin{picture}(2,2)
\put(1,0){\circle*{0.15}} 
\put(1.2,-0.1){\footnotesize $\alpha$}
\put(0,1){\circle*{0.15}} 
\put(0.2,0.9){\footnotesize $\beta$}
\put(1,0){\line(-1,1){1}} 
\end{picture}
}
\newcommand{\BTwoAAlphaMars}{
\setlength{\unitlength}{1cm}
\begin{picture}(3,3)
\put(2,0){\circle*{0.15}} 
\put(2.2,-0.1){\footnotesize $\alpha$}
\put(1,1){\circle*{0.15}} 
\put(1.2,0.9){\footnotesize $\beta$}
\put(0,2){\circle*{0.15}} 
\put(0.2,1.9){\footnotesize $\alpha$}
\put(2,0){\line(-1,1){2}} 
\end{picture}
}
\newcommand{\BTwoABetaMars}{
\setlength{\unitlength}{1cm}
\begin{picture}(3,3)
\put(2,0){\circle*{0.15}} 
\put(2.2,-0.1){\footnotesize $\beta$}
\put(1,1){\circle*{0.15}} 
\put(1.2,0.9){\footnotesize $\alpha$}
\put(0,2){\circle*{0.15}} 
\put(0.2,1.9){\footnotesize $\beta$}
\put(2,0){\line(-1,1){2}} 
\end{picture}
}
\newcommand{\BTwoBAlphaMars}{
\setlength{\unitlength}{1cm}
\begin{picture}(3,4)
\put(1,0){\circle*{0.15}} 
\put(1.2,-0.1){\footnotesize $\alpha$}
\put(0,1){\circle*{0.15}} 
\put(0.2,0.9){\footnotesize $\beta$}
\put(1,2){\circle*{0.15}} 
\put(1.2,1.9){\footnotesize $\beta$}
\put(0,3){\circle*{0.15}} 
\put(0.2,2.9){\footnotesize $\alpha$}
\put(0,1){\line(1,1){1}} 
\put(1,0){\line(-1,1){1}} 
\put(1,2){\line(-1,1){1}} 
\end{picture}
}
\newcommand{\BTwoBBetaMars}{
\setlength{\unitlength}{1cm}
\begin{picture}(3,4)
\put(1,0){\circle*{0.15}} 
\put(1.2,-0.1){\footnotesize $\beta$}
\put(0,1){\circle*{0.15}} 
\put(0.2,0.9){\footnotesize $\alpha$}
\put(1,2){\circle*{0.15}} 
\put(1.2,1.9){\footnotesize $\alpha$}
\put(0,3){\circle*{0.15}} 
\put(0.2,2.9){\footnotesize $\beta$}
\put(0,1){\line(1,1){1}} 
\put(1,0){\line(-1,1){1}} 
\put(1,2){\line(-1,1){1}} 
\end{picture}
}
\newcommand{\GTwoAAlphaMars}{
\setlength{\unitlength}{1cm}
\begin{picture}(4,6)
\put(3,0){\circle*{0.15}} 
\put(3.2,-0.1){\footnotesize $\alpha$}
\put(2,1){\circle*{0.15}} 
\put(2.2,0.9){\footnotesize $\beta$}
\put(1,2){\circle*{0.15}} 
\put(1.2,1.9){\footnotesize $\alpha$}
\put(2,3){\circle*{0.15}} 
\put(2.2,2.9){\footnotesize $\alpha$}
\put(1,4){\circle*{0.15}} 
\put(1.2,3.9){\footnotesize $\beta$}
\put(0,5){\circle*{0.15}} 
\put(0.2,4.9){\footnotesize $\alpha$}
\put(1,2){\line(1,1){1}} 
\put(3,0){\line(-1,1){2}} 
\put(2,3){\line(-1,1){2}} 
\end{picture}
}
\newcommand{\GTwoABetaMars}{
\setlength{\unitlength}{1cm}
\begin{picture}(4,6)
\put(3,0){\circle*{0.15}} 
\put(3.2,-0.1){\footnotesize $\beta$}
\put(2,1){\circle*{0.15}} 
\put(2.2,0.9){\footnotesize $\alpha$}
\put(1,2){\circle*{0.15}} 
\put(1.2,1.9){\footnotesize $\beta$}
\put(2,3){\circle*{0.15}} 
\put(2.2,2.9){\footnotesize $\beta$}
\put(1,4){\circle*{0.15}} 
\put(1.2,3.9){\footnotesize $\alpha$}
\put(0,5){\circle*{0.15}} 
\put(0.2,4.9){\footnotesize $\beta$}
\put(1,2){\line(1,1){1}} 
\put(3,0){\line(-1,1){2}} 
\put(2,3){\line(-1,1){2}} 
\end{picture}
}
\newcommand{\GTwoBAlphaMars}{
\setlength{\unitlength}{1cm}
\begin{picture}(4,8)
\put(2,0){\circle*{0.15}} 
\put(2.2,-0.15){\footnotesize $\alpha$}
\put(1,1){\circle*{0.15}} 
\put(1.2,0.9){\footnotesize $\beta$}
\put(2,2){\circle*{0.15}} 
\put(2.2,1.9){\footnotesize $\beta$}
\put(1,3){\circle*{0.15}} 
\put(1.2,2.9){\footnotesize $\alpha$}
\put(3,3){\circle*{0.15}} 
\put(3.2,2.9){\footnotesize $\beta$}
\put(0,4){\circle*{0.15}} 
\put(0.2,3.9){\footnotesize $\beta$}
\put(2,4){\circle*{0.15}} 
\put(2.2,3.9){\footnotesize $\alpha$}
\put(1,5){\circle*{0.15}} 
\put(1.2,4.9){\footnotesize $\beta$}
\put(2,6){\circle*{0.15}} 
\put(2.2,5.9){\footnotesize $\beta$}
\put(1,7){\circle*{0.15}} 
\put(1.2,6.9){\footnotesize $\alpha$}
\put(1,1){\line(1,1){2}} 
\put(1,3){\line(1,1){1}} 
\put(0,4){\line(1,1){2}} 
\put(2,0){\line(-1,1){1}} 
\put(2,2){\line(-1,1){2}} 
\put(3,3){\line(-1,1){2}} 
\put(2,6){\line(-1,1){1}} 
\end{picture} 
}
\newcommand{\GTwoBBetaMars}{
\setlength{\unitlength}{1cm}
\begin{picture}(4,8)
\put(2,0){\circle*{0.15}} 
\put(2.2,-0.15){\footnotesize $\beta$}
\put(1,1){\circle*{0.15}} 
\put(1.2,0.9){\footnotesize $\alpha$}
\put(2,2){\circle*{0.15}} 
\put(2.2,1.9){\footnotesize $\alpha$}
\put(1,3){\circle*{0.15}} 
\put(1.2,2.9){\footnotesize $\beta$}
\put(3,3){\circle*{0.15}} 
\put(3.2,2.9){\footnotesize $\alpha$}
\put(0,4){\circle*{0.15}} 
\put(0.2,3.9){\footnotesize $\alpha$}
\put(2,4){\circle*{0.15}} 
\put(2.2,3.9){\footnotesize $\beta$}
\put(1,5){\circle*{0.15}} 
\put(1.2,4.9){\footnotesize $\alpha$}
\put(2,6){\circle*{0.15}} 
\put(2.2,5.9){\footnotesize $\alpha$}
\put(1,7){\circle*{0.15}} 
\put(1.2,6.9){\footnotesize $\beta$}
\put(1,1){\line(1,1){2}} 
\put(1,3){\line(1,1){1}} 
\put(0,4){\line(1,1){2}} 
\put(2,0){\line(-1,1){1}} 
\put(2,2){\line(-1,1){2}} 
\put(3,3){\line(-1,1){2}} 
\put(2,6){\line(-1,1){1}} 
\end{picture} 
}
\noindent
{\bf \FundamentalClass}\ \ {\sl Let $P$ be a nonempty indecomposable \dichromatic grid poset with chain function} $\mychain: P \longrightarrow [1,m]_{\mathbb{Z}}$ {\sl and two-color function} $\mycolor: P \longrightarrow I$, {\sl where $I = \{\alpha,\beta\}$ is to be our color palette.   
Suppose the DCDL} $L = \Jcolor(P)$ {\sl is $\mathscr{G}$-structured by an IEG $\mathscr{G} = (M_{I \times I}, \Gamma_{I})$. 
Then $\mathscr{G}$ is a Coxeter--Dynkin posy, and $P$ is isomorphic to one of the following twelve connected \dichromatic grid posets:}
\begin{center}
\setlength{\unitlength}{1cm}
\begin{picture}(15,2.5)
\put(0,1){\AOneAOneAlphaMars}
\put(1.5,1){\AOneAOneBetaMars} 
\put(3,0.5){\ATwoAlphaMars} 
\put(6,0.5){\ATwoBetaMars}
\put(9,0){\BTwoAAlphaMars}
\put(13,0){\BTwoABetaMars}
\end{picture} 
\end{center}

\begin{center}
\setlength{\unitlength}{1cm}
\begin{picture}(15,5.5)
\put(0,1){\BTwoBAlphaMars}
\put(3,1){\BTwoBBetaMars} 
\put(6,0){\GTwoAAlphaMars} 
\put(11,0){\GTwoABetaMars}
\end{picture} 
\end{center}

\begin{center}
\setlength{\unitlength}{1cm}
\begin{picture}(15,7.5)
\put(2,0){\GTwoBAlphaMars}
\put(8,0){\GTwoBBetaMars} 
\end{picture} 
\end{center}

{\em Proof.} First note that $P$ must be connected (otherwise,  $P$ decomposes).  
Then without loss of generality, we may assume that the chain function $\mychain: P \longrightarrow [1,m]_{\mathbb{Z}}$ is surjective.  
If $\mycolor$ is not surjective, then $P$ is a chain whose vertices have the same color.  
But in this case, $P$ must be a one-element chain (otherwise, $P$ decomposes), and hence $P$ is isomorphic to one of the first two in our list of twelve $\alpha\, \beta$ grid posets. 
Since $L$ will be a two element chain, say \parbox{0.4cm}{\begin{center}
\setlength{\unitlength}{0.2cm}
\begin{picture}(1.5,1)
\put(0.75,-0.5){\circle*{0.5}} 
\put(0.75,2.5){\circle*{0.5}} 
\put(0.75,-0.5){\line(0,1){3}} 
\put(0.35,0.75){\footnotesize $\gamma$} 
\end{picture} \end{center}}, where $\gamma$ is one of $\alpha$ or $\beta$, then the $\gamma$ row of $M$ consists of a `$2$' and a `$0$', thereby forcing our pulsation matrix $M$ to be $\left(\begin{array}{cc}2 & 0\\ 0 & 2\end{array}\right)$. 

Now assume that the two-color function $\mycolor: P \longrightarrow \{\alpha,\beta\}$ is surjective.  
In particular, $P$ is not a disjoint union of chains.  
Then by \NewLFKTheorem, it is the case that (1) $M$ is one of the the matrices in the statement of \TechnicalProposition\ or (2) $M$ is the transpose of one of these matrices.  
In case (1), \TechnicalProposition\  --- together with the fact that $P$ is indecomposable --- implies that $P$ is isomorphic to one of the six \dichromatic grid posets depicted in the statement of \TechnicalProposition.  
In case (2), the same reasoning shows that $P$ is isomorphic to one of the six \dichromatic grid posets depicted in the statement of \TechnicalProposition, but with the vertex colors switched.\hfill\QED

{\bf [\S \GridPosetSection.6:\! A combinatorial characterization of rank two semistandard posets.]}

\noindent
{\bf \SemistandardClass}\ \ {\sl Let $P$ be a \dichromatic grid 
poset with chain function} $\mychain: P \longrightarrow [1,m]_{\mathbb{Z}}$ {\sl and with surjective two-color function} $\mycolor: P \longrightarrow I$, {\sl where $I = \{\alpha,\beta\}$ our color palette. 
Suppose $P$ has the max property.  
Further suppose the DCDL} $L = \Jcolor(P)$ {\sl is $\mathscr{G}$-structured by an IEG $\mathscr{G} = (M_{I \times I},  \Gamma_{I})$; switch the vertex colors of $P$ as needed so that $|M_{\beta,\alpha}| \geq |M_{\alpha,\beta}|$. 
Then $\mathscr{G}$ is a Coxeter--Dynkin posy.  
Let $a$ (respectively $b$) be the number of color $\alpha$ (respectively $\beta$) face vertices in $P$, and let $\lambda := (a,b)$. 
If the rightmost maximal vertex of $P$ has color $\alpha$ (respectively color $\beta$), then $P$ is isomorphic to the semistandard poset $\Pba$ (respectively $\Pab$).}

{\em Proof.} 
Set $\melt := \mymax(L)$, to be identified with the down-set consisting of all elements of $P$. 
Without loss of generality, we may assume that $\mychain$ is surjective.  
If $P$ is not connected, then $P$ has at least two maximal elements.  
Since $P$ has the max property, these maximal elements must be in $\mathcal{C}_{1} \cup \mathcal{C}_{2}$.  
It follows that chain $\mathcal{C}_{1}$ is in its own connected component.  
Let $v_{1}$ denote the maximal element of $\mathcal{C}_{1}$.  
Let $\selt := \melt \setminus \{v_{1}\}$; then $\selt$ is a down-set from $P$.  
Suppose $\mycolor(v_{1}) = \alpha$.  
Then $wt(\selt) + (M_{\alpha,\alpha},M_{\alpha,\beta}) = wt(\melt)$.  
Since $A_{\beta}(\selt) = \emptyset$, then $\selt$ is maximal in its color $\beta$ component.  
Also, one can see that $D_{\beta}(\selt) = D_{\beta}(\melt)$.  
Then $l_{\beta}(\melt) = \rho_{\beta}(\melt) = |D_{\beta}(\selt)| = \rho_{\beta}(\selt) = l_{\beta}(\selt)$. 
Then $M_{\alpha,\beta} = 0$.  
\NewLFKTheorem\ now implies that $M_{\beta,\alpha} = 0$ and that $P$ is a disjoint union of chains.  
So $\mathcal{C}_{2}$ is in its own connected component, and all vertices in $\mathcal{C}_{2}$ have color $\beta$.  
Similarly, if $\mycolor(v_{1}) = \beta$, then $M_{\beta,\alpha} = 0 = M_{\alpha,\beta}$ and $P$ is the disjoint union of the chains $\mathcal{C}_{1}$ and $\mathcal{C}_{2}$, where elements of $\mathcal{C}_{1}$ have color $\beta$ and elements of $\mathcal{C}_{2}$ have color $\alpha$.  
In either case, $P$ is isomorphic to the appropriate $\myA_{1} \oplus \myA_{1}$-semistandard poset and $M$ is the pulsation matrix $\left(\begin{array}{cc}2 & 0\\ 0 & 2\end{array}\right)$. 

Now assume $P$ is connected, with $|P| = N$. 
Then by \NewLFKTheorem, $M$ is one of the three matrices appearing in the statement of \TechnicalProposition.  
We proceed by induction. 
For the inductive hypothesis, assume that the theorem statement holds for all connected \dichromatic grid posets $Q$ meeting the hypotheses of the theorem statement and such that $|Q| < N$.  
Let $z_{1}$ be the rightmost maximal vertex of $P$ (in chain $\mathcal{C}_{2}$ of $P$ if there are two maximal vertices, in chain $\mathcal{C}_{1}$ if there is only one maximal vertex).  
Use \TechnicalProposition\ to locate a subset $P_{2} = \{z_{1},\ldots,z_{k}\}$ of $P$ (where $k \geq 1$) which, given the induced order from $P$ and viewed as a \dichromatic grid subposet of $P$, is isomorphic to one of the six \dichromatic grid posets appearing in the statement of \TechnicalProposition.  
If $P \setminus P_{2}$ is empty, then $P = P_{2}$ is isomorphic to one of the fundamental posets $P_{\mathscr{G}}(\omega_{\alpha})$ (if $\mycolor(z_{1}) = \alpha$) or $P_{\mathscr{G}}(\omega_{\beta})$ (if $\mycolor(z_{1}) = \beta$) for $\myA_{2}$, $\myC_{2}$, or $\myG_{2}$, and we are done.  

Now suppose $P_{1} := P \setminus P_{2}$ is nonempty.  
Suppose $\mychain(z_{1}) = 1$.  
Let $w$ be the rightmost maximal element of $P_{1}$.  
Then since $P = P_{1} \triangleleft P_{2}$, it follows that $\mychain(w) = 1$, and hence $P_{1}$ has $w$ as its unique maximal element. 
In particular $P_{1}$ must be connected.  
Now suppose $\mychain(z_{1}) = 2$.  
Let $w$ be the rightmost maximal element of $P_{1}$.  
Then since $P = P_{1} \triangleleft P_{2}$, it follows that $\mychain(w) \leq 2$.  
So $P_{1}$ has at most two maximal elements.  
If $\mychain(w) = 1$, then $P_{1}$ has a unique maximal element, and hence is connected.  
Then suppose $\mychain(w) = 2$.  
Since $P$ is connected, there exists some covering relation $u \rightarrow v$ in $P$ with $\mychain(u) = 2$ and $\mychain(v) = 1$.  
Since $u \not= z_{1}$ ($u$ is not  maximal in $P$), it follows that $u \in P_{1}$.  
For each $z \in P_{2}$ it is the case that $\mychain(z) \geq 2$; it follows that $v \in P_{1}$.  
Let $w'$ be the maximal element of $\mathcal{C}_{1}$ in $P$; then $w' \in P_{1}$.  
So we have $v \leq_{P_{1}} w'$, $u \rightarrow v$ in $P_{1}$, and $u \leq_{P_{1}} w$ in $P_{1}$; together imply that $w'$ and $w$ are in the same connected component of $P_{1}$.  
Since these are the only two maximal elements of $P_{1}$, it follows that $P_{1}$ is connected.  
In particular, $P_{1}$ must also have the max property. 

The previous paragraph established by an examination of cases that $P_{1}$ is connected and has the max property. 
Since $P_{1}$ is connected and $\mychain(v) = 1$ for some $v$ in $P_{1}$, then for some integer $m' \in [1,m]_{\mathbb{Z}}$, it is the case that $\mychain|_{P_{1}}: P_{1} \longrightarrow [1,m']_{\mathbb{Z}}$ is surjective. 
We claim that $m' > 1$. 
Suppose otherwise, so $m' = 1$.  
The connectedness of $P$  implies that $\mychain(z_{1}) = 1$.   
Suppose $\mycolor(z_{1}) = \alpha$. Let $w$ be the maximal element of $P_{1}$.  
Observe that $w \rightarrow z_{1}$ in $P$.  
Set $\telt := \melt \setminus \{z_{1}\}$, and let $\selt := \telt \setminus \{w\}$. 
One can see that $\selt$ and $\telt$ are down-sets from $P$ and that $\selt \mylongarrow{\alpha} \telt$ in $L$. 
Therefore $wt(\selt) + (M_{\alpha,\alpha},M_{\alpha,\beta}) = wt(\telt)$.     
However, the fact that $m' = 1$ implies that $\rho_{\beta}(\selt) = l_{\beta}(\selt) = \rho_{\beta}(\telt) = l_{\beta}(\telt) = 0$. 
Therefore $M_{\alpha,\beta} = 0$.  
This contradicts the fact that $M_{\alpha,\beta} \not= 0$.  
Arrive at a similar contradiction if $\mycolor(z_{1}) = \beta$.   
So $m' > 1$.  
It now follows that the color function $\mycolor|_{P_{1}}$ is also surjective. 
For each $i \in \{1,2\}$, let $L_{i}$ be the DCDL $\Jcolor(P_{i})$ with weight function $wt_{i}: L \longrightarrow \mathbb{Z} \times \mathbb{Z}$. 
We claim that $L_{1}$ is $\mathscr{G}$-structured. 
Suppose $\selt \mylongarrow{\alpha} \telt$ in $L_{1}$, thinking of $\selt$ and $\telt$ as down-sets from $P$ as well as from $P_{1}$.  
By \WeightsLemma, $wt(\selt) = wt_{1}(\selt) + wt_{2}(\emptyset)$ and $wt(\telt) = wt_{1}(\telt) + wt_{2}(\emptyset)$. 
Since $L$ is $\mathscr{G}$-structured, then $wt(\selt) + (M_{\alpha,\alpha},M_{\alpha,\beta}) = wt(\telt)$. 
Combine the preceding two sentences to conclude that $wt_{1}(\selt) + (M_{\alpha,\alpha},M_{\alpha,\beta}) = wt_{1}(\telt)$.  
A similar argument shows that if $\selt \mylongarrow{\beta} \telt$ in $L_{1}$, then $wt_{1}(\selt) + (M_{\beta,\alpha},M_{\beta,\beta}) = wt_{1}(\telt)$.  

So the inductive hypothesis applies to the connected \dichromatic grid poset $P_{1}$.  
Let $a'$ (respectively $b'$) be the number of face vertices of $P_{1}$ of color $\alpha$ (respectively color $\beta$).  
Let $\mu := (a',b')$. 
One can see that if $\mycolor(z_{1}) = \alpha$ (respectively $\mycolor(z_{1}) = \beta$), then $a' = a-1$ and $b' = b$ (respectively $a' = a$ and $b' = b-1$), in which case $P_{1}$ is isomorphic to the semistandard poset $P_{\mathscr{G}}^{\beta\alpha}(\mu)$ (respectively $P_{\mathscr{G}}^{\alpha\beta}(\mu)$) and $P_{2}$ is isomorphic to $P_{\mathscr{G}}(\omega_{\alpha})$ (respectively $P_{\mathscr{G}}(\omega_{\beta})$). 
Observe now that since $P$ has the max property, the only one way we can have $P = P_{1} \triangleleft P_{2}$ is if $P$ is isomorphic to $\Pba = P_{\mathscr{G}}^{\beta\alpha}(\mu) \triangleleft P_{\mathscr{G}}(\omega_{\alpha})$ (respectively $\Pab = P_{\mathscr{G}}^{\alpha\beta}(\mu) \triangleleft P_{\mathscr{G}}(\omega_{\beta})$).  
(Note that when $a' = a = 0$ or $b' = b = 0$, then $P_{\mathscr{G}}^{\beta\alpha}(\mu) = P_{\mathscr{G}}^{\alpha\beta}(\mu)$.)\hfill\QED

\begin{center}
\underline{\hspace*{4in}}
\end{center}

\vspace*{0.5cm} 
\noindent 
{\bf \S \TwoColorSturdySection.\ A sturdiness result from the theory of two-color grid posets.} 
The main result of this section (\GStructureDCDLTheorem) is the demonstration of sturdiness for DCDL's obtained from certain scaffolds/skew-stacks. 
To obtain this new result, we will utilize the theory of two-color grid posets. 

{\bf [\S \TwoColorSturdySection.1:\! Two-color subordinates of certain skew-stacks.]} 
\IsomLatticeTheorem\ established an isomorphism between the diamond-colored distributive lattice of order ideals from a skew-stack $\mystackedP := P_{1} \myposetstack_{\phi_{1}} P_{2} \myposetstack_{\phi_{2}} \cdots \myposetstack_{\phi_{m-1}} P_{m}$ and the DCDL of ideal arrays over the corresponding birds-eye scaffold $\scaffoldS^{\mybirdseye} = (\mystackedP^{\mybirdseye},\mathcal{A}^{\mybirdseye},\mathcal{B}^{\mybirdseye},\vcolor^{\mybirdseye})$. 
The next somewhat technical result invokes the theory of two-color grid posets and will be applied in our proof of \GStructureDCDLTheorem. 

\noindent 
{\bf \SubordinateGridDecomp}\ \ {\sl Take the hypotheses of \IsomLatticeTheorem. 
For each $r \in [1,m]_{\mathbb{Z}}$, assume that $\vcolor^{(r)}(x) \not= \vcolor^{(r)}(y)$ whenever $x \rightarrow y$ in $P_{r}$ and that $(\vcolor^{\mybirdseye})^{-1}(i)$ is a totally ordered subset of} $\mystackedP^{\mybirdseye}$ {\sl for each $i \in I_{n}$. 
Suppose further that for any pair $i,j \in I_{n}$ of distinct colors from $I_{n}$; any pair $r,s \in [1,m]_{\mathbb{Z}}$; and any pair of connected and nonempty $\{i,j\}$-subordinates $\mathcal{C}$ and $\mathcal{C}'$ from $P_{r}$ and $P_{s}$ respectively, it is the case that $\mathcal{C}$ and $\mathcal{C}'$ are saturated chains of the same length. 
Let $\mathcal{S}$ be a nonempty connected $\{i,j\}$-subordinate of} $\mystackedP$. 
{\sl Then $\mathcal{S}$ is a two-color grid poset that decomposes as}
\[\mathcal{S} = (\mathcal{S} \cap P_{1}) \triangleleft (\mathcal{S} \cap P_{2}) \triangleleft \cdots \triangleleft (\mathcal{S} \cap P_{m}).\]

{\em Proof.} . 
If $m=1$, then $\mystackedP = P_{1}$, in which case our connected $\{i,j\}$-subordinate $\mathcal{S}$ is a saturated chain. 
Now, $\vcolor^{(1)}(x) \not= \vcolor^{(1)}(y)$ whenever $x \rightarrow y$ in $P_{1}$, by hypothesis. 
Therefore, within our saturated chain, the vertex colors must alternate as we read vertices sequentially from top to bottom. 
So, $\mathcal{S} = \mathcal{S} \cap P_{1}$ is a two-color grid poset. 

Now suppose $m > 1$, and consider our generic skew-stacked poset $\mystackedP$ from the proposition statement together with a given connected $\{i,j\}$-subordinate $\mathcal{S}$. 
That is, $\mathcal{S} = \mathcal{K} \setminus \mathcal{J}$ is a set of color $i$ or $j$ elements of $\mystackedP$, where $\mathcal{K}$ and $\mathcal{J}$ are down-sets from $\mystackedP$ with $\mathcal{J} \subseteq \mathcal{K}$  such that $\mystackedP \setminus \mathcal{K}$ has no minimal element of color $i$ or $j$ and $\mathcal{J}$ has no maximal element of color $i$ or $j$. 
For each $r \in [1,m]_{\mathbb{Z}}$, let $\mathcal{K}_{r} := \mathcal{K} \cap P_{r}$, $\mathcal{J}_{r} := \mathcal{J} \cap P_{r}$, and $\mathcal{S}_{r} := \mathcal{S} \cap P_{r}$. 

Our first aim is to analyze the structure of each $\mathcal{S}_{r}$. 
Since $\mathcal{S}$ is connected, then there must be integers $r_{1}$ and $r_{2}$ such that $\mathcal{S}_{r}$ nonempty if and only if $r$ is an integer and $r_{1} \leq r \leq r_{2}$.  
In what follows, the empty $S_{r}$'s have no effect on our argument, so for simplicity we take $r_{1} = 1$ and $r_{2} = m$. 
It is routine to verify that $\mathcal{S}_{r}$ is an $\{i,j\}$-subordinate of $P_{r}$ with $\mathcal{S}_{r} = \mathcal{K}_{r} \setminus \mathcal{J}_{r}$ for the down-sets $\mathcal{K}_{r}$ and $\mathcal{J}_{r}$ from $P_{r}$. 
We claim that each $\mathcal{S}_{r}$ is connected and therefore a saturated chain with alternating vertex colors.  
Our approach is to take elements $u$ and $v$ in $\mathcal{S}_{r}$ and show they are path connected to one another by considering the cases (1) they have the same color and (2) they have different colors. 
In each case, we argue that $u$ and $v$ are contained in some connected union of intervals within $\mathcal{S}_{r}$. 
To this end, we analyze intervals $[x,y]_{P_{r}}$ in $P_{r}$ determined by certain elements $x$ and $y$ of $P_{r}$. 

To address case (1), suppose our given $u, v \in \mathcal{S}_{r}$ have the same color.  
Then these elements are comparable in $P_{r}$, and we observe, then, that $[u,v]_{P_{r}}$ is necessarily within the difference of down-sets $\mathcal{K}_{r} \setminus \mathcal{J}_{r} = \mathcal{S}_{r}$. 
For case (2), suppose $u$ and $v$ have different colors. 
If $u$ and $v$ are comparable, we see as in case (1) that $[u,v]_{P_{r}}$ is entirely within the difference of down-sets $\mathcal{K}_{r} \setminus \mathcal{J}_{r} = \mathcal{S}_{r}$. 
So, suppose $u$ and $v$ are incomparable. 
Now, the connected component of $u$ within $\mathcal{S}_{r}$ is a saturated chain whose elements alternate in color.  
If this chain has more than one element, then it must have an element $u'$ whose color agrees with $v$.  
As in case (1), $[u',v]_{P_{r}}$ will be inside the difference of down-sets $\mathcal{K}_{r} \setminus \mathcal{J}_{r} = \mathcal{S}_{r}$. 
But since $u'$ and $u$ are comparable, then $[u,u']_{P_{r}}$ is contained within $\mathcal{K}_{r} \setminus \mathcal{J}_{r} = \mathcal{S}_{r}$ as well. 
That is, $v$ must be within this connected component. 
Of course, we may interchange the roles of $u$ and $v$ in the preceding four sentences to see that if the connected component of $v$ within $\mathcal{S}_{r}$ has more than one element, then $u$ must be within this connected component. 
To exhaust all possibilities, we must lastly consider the case that that $\mathcal{S}_{r}$ consists of exactly of two incomparable (and differently colored) elements $u$ and $v$. 
Since the sets `$\{u\}$' and `$\{v\}$' are each, by themselves, $\{i,j\}$-subordinates, then all connected $\{i,j\}$-subordinates of each $P_{s}$ ($s \in [1,m]_{\mathbb{Z}}$) consist of exactly one element. 
This means $\mathcal{S}$ must consist of a saturated chain of vertices of color $i$ disconnected from a saturated chain of color $j$, contradicting the hypothesis that $\mathcal{S}$ is connected. 
Thus, in all cases, $\mathcal{S}_{r}$ is a connected $\{i,j\}$-subordinate and therefore a saturated chain whose elements alternate in color. 

To realize $\mathcal{S}$ as a two-color grid poset, we need a chain function $\mygridchain: \mathcal{S} \rightarrow \{1,2,\ldots,c\}$ (for some positive integer $c$) and a color function $\mygridcolor: \mathcal{S} \rightarrow \{i,j\}$.
Let $t_{r}$ be the maximal element of each $\mathcal{S}_{r}$ and $w_{r}$ its minimal element. 
We claim that for $r > 1$, there exists a unique element $u_{r-1} \in \mathcal{S}_{r-1}$ such that $u_{r-1} \rightarrow t_{r}$. 
Indeed, connectedness of $\mathcal{S}$ requires that $v \rightarrow \phi_{r-1}(v)$ for some $v \in \mathcal{S}_{r-1}$ and $\phi_{r-1}(v) \in \mathcal{S}_{r}$, where $\phi_{r-1}: \mathcal{I}_{r-1} \rightarrow \mathcal{F}_{r}$ is an isomorphism. 
Since $\phi_{r-1}(v)$ is in the up-set $\mathcal{F}_{r}$, then $t_{r}$ is also in $\mathcal{F}_{r}$. 
Since $t_{r}$ is in the image of $\phi_{r-1}$, then the preimage $\phi_{r-1}^{-1}(t_{r})$ is an element in the down-set $\mathcal{I}_{r-1}$ which we now name $u_{r-1}$. 
By definition, $t_{r} \in \mathcal{K} \setminus \mathcal{J}$, since $\mathcal{S}_{r} := (\mathcal{K} \cap P_{r}) \setminus (\mathcal{J} \cap P_{r}) = (\mathcal{K} \setminus \mathcal{J}) \cap P_{r}$. 
In particular, this means that $u_{r-1}$ must be in the down-set $\mathcal{K}$ and thus in $\mathcal{K}_{r-1}$. 
Suppose for the moment that $u_{r-1} \in \mathcal{J}$. 
Within the saturated chain $\phi_{r-1}^{-1}(\mathcal{S}_{r} \cap \mathcal{F}_{r})$, which includes $u_{r-1}$, let $u$ be largest with $u \in \mathcal{S}_{r-1}$. 
Now $\phi_{r-1}(v) \in \mathcal{S}_{r}$ and $\phi_{r-1}(v) \leq \phi_{r-1}(u_{r-1}) = t_{r}$ means that $v$ is in the saturated chain $\phi_{r-1}^{-1}(\mathcal{S}_{r} \cap \mathcal{F}_{r})$, and therefore $v \leq u_{r-1}$. 
But the latter fact forces $v$ to be in the down-set $\mathcal{J}$, which is not the case since $v \in \mathcal{S}_{r-1} = (\mathcal{K} \setminus \mathcal{J}) \cap P_{r-1}$. 
Thus, we have ruled out $u_{r-1} \in \mathcal{J}$, meaning $u_{r-1} \in (\mathcal{K} \setminus \mathcal{J}) \cap P_{r-1} = \mathcal{S}_{r-1}$, as desired. 
Similarly, we can show that if $r < m$, then there exists a unique element $v_{r+1} \in \mathcal{S}_{r+1}$ with $w_{r} \rightarrow v_{r+1}$. 

Within the context of the foregoing paragraph, for any $x \in \mathcal{S}$, let $\myr(x)$ be the number $r$ such that $x \in \mathcal{S}_{r}$, and let $\mydepth(x)$ represent the number of steps from $x$ up to $t_{\mysmallr(x)}$ within the saturated chain $\mathcal{S}_{r}$. 
Set $\mygridchain(t_{1}) := 1$ and $\mygridchain(x) := 1+\mydepth(x)$ if $\myr(x)=1$. 
Then for all $x$ with $\myr(x)>1$, we can recursively define $\mygridchain(t_{\mysmallr(x)}) := \mygridchain(u_{\mysmallr(x)-1})$ and set $\mygridchain(x) := \mygridchain(t_{\mysmallr(x)}) + \mydepth(x)$. 
Now suppose $u \rightarrow v$ for some $u, v \in \mathcal{S}$. 
If $\myr(u) = \myr(v)$, then $\mygridchain(u) = \mygridchain(t_{\mysmallr(u)}) + \mydepth(u) = \mygridchain(t_{\mysmallr(v)}) + \mydepth(v)+1 = \mygridchain(v)+1$. 
If $\myr(u) \not= \myr(v)$, then our poset stacking process forces $\myr(u)+1 = \myr(v) =: r$ and $v = \phi_{r-1}(u)$. 
In this case, 
\begin{eqnarray*}
\mygridchain(u) & = & \mygridchain(t_{r-1}) + \mydepth(u)\\ 
 & = & \mygridchain(t_{r-1}) - \mygridchain(u_{r-1}) + \mygridchain(u_{r-1})\\
 & &  \hspace*{0.25cm}+ \mydepth(v) - \mydepth(v) + \mydepth(u)\\ 
 & = & \mygridchain(t_{r}) + \mydepth(v) + \mygridchain(t_{r-1}) - \mygridchain(u_{r-1})\\ 
 & & \hspace*{0.25cm}+ \mydepth(u) - \mydepth(v)\\ 
 & = & \mygridchain(v) + \mydepth(u) - (\mygridchain(u_{r-1}) - \mygridchain(t_{r-1})) - \mydepth(v)\\
 & = & \mygridchain(v) + \mydepth(u) - \mydepth(u_{r-1}) + \mydepth(t_{r}) - \mydepth(v) 
\end{eqnarray*}

\newpage
\noindent
\begin{eqnarray*}
\hspace*{0.8in} & = & \mygridchain(v) + \mydepth(u) - \mydepth(u_{r-1}) + \mydepth(\phi_{r-1}(u_{r-1})) - \mydepth(\phi(u))\\ 
 & = & \mygridchain(v),   
\end{eqnarray*}
where the equality $\mydepth(u) - \mydepth(u_{r-1}) = \mydepth(\phi(u)) - \mydepth(\phi_{r-1}(u_{r-1}))$, which is invoked in the last step, follows from the fact that $\phi_{r-1}$ is an isomorphism. 
Thus, $\mygridchain$ is a chain function whose image is some interval of integers $[1,c]_{\mathbb{Z}}$. 
Moreover, if (say) $t_{1}$ has color $i$, then any given $v \in \mathcal{S}$ has color $i$ if $\mygridchain(v)$ is odd and color $j$ if $\mygridchain(v)$ is even. 
We can use the preceding rule to define the coloring function $\mygridcolor$.  
So, $(\mathcal{S},\leq_{\mathcal{S}},\mygridchain,\mygridcolor)$ is a two-color grid poset. 

To argue that the two-color grid poset $\mathcal{S}$ decomposes as $\mathcal{S}_{1} \triangleleft \mathcal{S}_{2} \triangleleft \cdots \triangleleft \mathcal{S}_{m}$, it only remains to show that $\mygridchain(w_{r}) \leq \mygridchain(w_{r+1})$ when $r \in [1,m-1]_{\mathbb{Z}}$ and that $\mygridchain(t_{r-1}) \leq \mygridchain(t_{r})$ when $r \in [2,m]_{\mathbb{Z}}$. 
For the latter, the inequality $\mygridchain(t_{r-1}) \leq \mygridchain(u_{r-1}) = \mygridchain(t_{r})$ follows from the definitions, and $\mygridchain(w_{r}) = \mygridchain(v_{r+1}) \leq \mygridchain(w_{r+1})$ follows similarly in the former case.\hfill\QED

{\bf [\S \TwoColorSturdySection.2:\! Sturdy DCDL's from certain scaffolds and skew-stacks.]} 
The next result applies \NewLFKTheorem\ as well as \IsomLatticeTheorem/\SubordinateGridDecomp\ and is primarily useful in establishing the sturdiness claim of \GStructureDCDLCorollary\ below for DCDL's obtained from skew-stacks of what are here called `minuscule compression posets'. 

\noindent 
{\bf \GStructureDCDLTheorem}\ \ {\sl Let $\mathscr{G} = (\Gamma_{I},M_{I \times I})$ be an integral embryophytic graph. 
Consider the skew-stacked poset} 
\[\mystackedP := P_{1} \myposetstack_{\phi_{1}} P_{2} \myposetstack_{\phi_{2}} \cdots \myposetstack_{\phi_{m-1}} P_{m}\] 
{\sl and the birds-eye scaffold} $\scaffoldS^{\mybirdseye} = (\mystackedP^{\mybirdseye},\mathcal{A}^{\mybirdseye},\mathcal{B}^{\mybirdseye},\vcolor^{\mybirdseye})$ {\sl corresponding to $\mystackedP$, where each of the tiers $P_{1}, P_{2}, \ldots, P_{m}$ is vertex-colored by $I$.
For each $r \in [1,m]_{\mathbb{Z}}$, assume that $\vcolor^{(r)}(x) \not= \vcolor^{(r)}(y)$ whenever $x \rightarrow y$ in $P_{r}$ and that $(\vcolor^{\mybirdseye})^{-1}(i)$ is a totally ordered subset of $\mystackedP^{\mybirdseye}$ for each $i \in I_{n}$. 
Suppose further that for any pair $i,j \in I_{n}$ of distinct colors from $I$; any pair $r,s \in \{1,2,\ldots,m\}$; and any pair of connected and nonempty $\{i,j\}$-subordinates $\mathcal{C}$ and $\mathcal{C}'$ from tiers $P_{r}$ and $P_{s}$ respectively, it is the case that $\mathcal{C}$ and $\mathcal{C}'$ are saturated chains of the same length. 
Finally, suppose that each} $\Jcolor(P_{r})$ {\sl is $\mathscr{G}$-structured and that the edge-coloring function} 
$\ecolor\! :\! \EdgeSet(\Jcolor(P_{r})) \rightarrow I$ {\sl is surjective. 
Then $\mathscr{G}$ is an integral Coxeter--Dynkin posy, and the diamond-colored distributive lattice} $\Jcolor(\mystackedP) \cong \Lcolor(\scaffoldS^{\mybirdseye})$ {\sl is $\mathscr{G}$-structured.}

{\em Proof.} 
Since $\Jcolor(P_{r})$ is $\mathscr{G}$-structured and $\mathscr{E}\! :\! \EdgeSet(\Jcolor(P_{r})) \longrightarrow I$ is surjective, it follows from  \NewLFKTheorem.1 that $\mathscr{G}$ is an integral Coxeter--Dynkin posy. 
By \IsomLatticeTheorem, $\Lcolor(\scaffoldS^{\mybirdseye})$ is $\mathscr{G}$-structured if $\Jcolor(\mystackedP)$ is. 
For any pair of distinct colors $i$ and $j$, let $\mathscr{G}_{i,j}$ be the corresponding sub-embryophyte of $\mathscr{G}$. 
To see that $\Jcolor(\mystackedP)$ is $\mathscr{G}$-structured, it suffices to check that $\Jcolor(\mathcal{S})$ is $\mathscr{G}_{i,j}$-structured when $\mathcal{S}$ is any nonempty connected $\{i,j\}$-subordinate of $\mystackedP$. 
By \SubordinateGridDecomp, Lemma 3.1 of \cite{ADLMPPW} applies to $\mathcal{S}$, which implies that $\Jcolor(\mathcal{S})$ is $\mathscr{G}_{i,j}$-structured.\hfill\QED

\newpage
\begin{center}
\fbox{\Large \bf Part III}\\ 
\underline{\large \bf A primer on algebraic contexts}
\end{center}
The integral embryophytic graphs (IEG's) introduced in \S \GStructureIntroSection\ above are, of course, a simple variation on objects known as Coxeter--Dynkin diagrams and Generalized Cartan Matrices. 
These latter objects are classical entry points into various Lie theoretic areas -- Weyl/Coxeter groups, root systems, semisimple Lie algebras and their representations, Kac--Moody algebras and their representations, the Weyl--Kac character formula and related Weyl symmetric functions, etc. 
In this part of our manuscript, we will develop some aspects of these Lie theoretic notions with three main purposes in mind. 
First, in \S \WeylSection\ we introduce and re-state/re-formulate some finitisitic aspects of this general algebraic setting that are crucial to our particular work. 
Second, in \S \LFKSection\ we highlight the resulting finitistic classifications as parts of a much larger framework of results we now call {\sl La Florado Klasado}, which is Esperantese for `The Flowering Classification' (\LFK). 
And third, \S \FlowerSection\ we record some basic facts (definitions, key results) about these Lie theoretic notions that will allow us later to convey precisely the specific contexts in which DCML's/DCDL's arise in our work. 
The depth of these topics and the brevity of our treatment here means that the level of difficulty for the reader is an order of magnitude greater than in earlier sections.

\vspace*{0.5cm} 
\noindent 
{\bf \S \WeylSection.\ Finitistic algebraic structures related to integral embryophytes.} 
Our main application of the foregoing \GStructureDCDLTheorem\ is the forthcoming \GStructureDCDLCorollary. 
To properly understand that latter result, we develop in this section a part of the algebraic context in which the key structures of the theorem -- minuscule compression posets and their companion minuscule splitting DCDL's -- most naturally arise. 
Along the way, we will record more Dynkin diagram classification results, two of which seem to be new (\OrbitProposition\ and \NonconstantTheorem). 

{\bf [\S \WeylSection.1:\! Basic set-up.]} 
For the whole of this section, we assume that $\mathscr{G} = (\Gamma_{I},M_{I \times I})$ is an integral embryophytic graph (IEG). 
Let $n := |I|$, the cardinality of our palette of colors, and let $\{\myroot_{i}\}_{i \in I}$ be a fixed basis for an $n$-dimensional real vector space $\mathfrak{E}$. 
Let $\mathscr{B} := \{\omega_{i}\}_{i \in I}$ be the dual basis -- which, from here on, we will call the {\em fundamental basis} --  for the dual space $\mathfrak{E}^{*}$. 
We use `$\zeroweight$' to denote the zero vector of $\mathfrak{E}^{*}$. 
Identify any simple root $\alpha_{i}$ as an element of $\mathfrak{E}^{*}$ by taking $\alpha_{i} := \sum_{j \in I}M_{ij}\omega_{j}$. 
Then, for all $i,j \in I$, we have $\alpha_{i}(\myroot_{j}) = \sum_{k \in I}M_{ik}\omega_{k}(\myroot_{j}) = M_{ij}$. 
For historical reasons and for reasons of compatibility with classical treatments of the ideas we develop next, we will initially focus on the dual space $\mathfrak{E}^{*}$ rather than on $\mathfrak{E}$. 

{\bf [\S \WeylSection.2:\! The Weyl group.]} 
We now define some linear transformations on $\mathfrak{E}^{*}$ which bear some similarities to reflections in hyperplanes orthogonal to the simple roots, although, at the moment, we are working without an inner product on $\mathfrak{E}^{*}$. 
For each $i \in I$, define $S^{*}_{i}: \mathfrak{E}^{*} \longrightarrow \mathfrak{E}^{*}$ by the rule $S^{*}_{i}(\mu) := \mu - \mu_{i}\alpha_{i}$, where $\mu = \sum_{j \in I}\mu_{j}\omega_{j}$. 
Observe that $(S^{*}_{i})^{2}$ is the identity, so each $S^{*}_{i}$ is an order two linear transformation.  
Moreover, each $S_{i}$ preserves the $\mathbb{Z}$-span $\Lambda$ of our fundamental basis, and we can view each $S^{*}_{i}$ as a member of the group $GL(\Lambda) := \{T: \Lambda \longrightarrow \Lambda\, |\, T \mbox{ is linear and invertible}\}$. 
With respect to $\mathscr{B}$, each $S^{*}_{i}$ is represented by a matrix $[S^{*}_{i}]_{\mathscr{B}}$ of integers whose inverse is also a matrix of integers, i.e.\ $[S^{*}_{i}]_{\mathscr{B}}$ resides within the group $GL_{n}(\mathbb{Z})$ of so-called unimodular matrices. 
The {\em Weyl group} $\mathcal{W}(\mathscr{G})$ associated with our given embryophyte is the subgroup of $GL(\Lambda)$ generated by the set $\{S^{*}_{i}\}_{i \in I}$ or, equivalently, the subgroup of $GL_{n}(\mathbb{Z})$ generated by $\{[S^{*}_{i}]_{\mathscr{B}}\}_{i \in I}$. 

The following result about Weyl groups is well known (see for example \cite{Kac}). 
For a proof from first principles, see \cite{DonRoots}. 

\noindent
{\bf \GeneratorsAndRelations}\ \ {\sl Via the correspondence $S^{*}_{i} \leftrightarrow \mygens_{i}$, the Weyl group $\mathcal{W}(\mathscr{G}) = \langle \{S^{*}_{i}\}_{i \in I} \rangle$ is isomorphic to the group with the following presentation by generators and relations:} $\displaystyle \left\langle \{\mygens_{i}\}_{i \in I}\, \rule[-2mm]{0.2mm}{6mm}\, (\mygens_{i}\mygens_{j})^{m_{ij}}=\varepsilon\right\rangle,$ {\sl where}
\[m_{ij} = \left\{\begin{array}{cl}
1 & \hspace*{0.1in}\mbox{if\ \ $i=j$}\\ 
2 & \hspace*{0.1in}\mbox{if\ \ $M_{ij}=M_{ji}=0$}\\ 
3 & \hspace*{0.1in}\mbox{if\ \ $M_{ij}M_{ji}=1$}\\ 
4 & \hspace*{0.1in}\mbox{if\ \ $M_{ij}M_{ji}=2$}\\ 
6 & \hspace*{0.1in}\mbox{if\ \ $M_{ij}M_{ji}=3$}\\ 
\infty & \hspace*{0.1in}\mbox{if\ \ $i \not= j$ and $M_{ij}M_{ji} \geq 4$} 
\end{array}\right.\] 
{\sl and where `$(\mygens_{i}\mygens_{j})^{\infty}=\varepsilon$' means that there is no relation between $\mygens_{i}$ and $\mygens_{j}$.}

From here on, we will view the Weyl group in terms of the concrete generators $S^{*}_{i}$ or in terms of the abstract generators $\mygens_{i}$, whichever is convenient. 
Notice that for $i \not= j$ and $0 \leq M_{ij}M_{ji} < 4$, the number $m_{ij}$ of the preceding result is the unique positive integer such that $M_{ij}M_{ji}=4\cos^{2}(\pi/m_{ij})$. 
Therefore $(\mathcal{W}(\mathscr{G}),\{\mygens_{i}\}_{i \in I})$ is a Coxeter group, cf.\ \cite{HumCoxeter}. 
Given a sub-embryophyte $\mathscr{H} = (\Gamma_{J},M_{J \times J})$, whose vascular graph node colors comprise a subset $J$ of $I$, the Weyl group $\mathcal{W}(\mathscr{H})$ is isomorphic to the subgroup of $\mathcal{W}(\mathscr{G})$ generated by $\{\mygens_{j}\}_{j \in J}$ and is itself a Coxeter group. 
In the theory of Coxeter groups, Coxeter subgroups obtained in this way are typically called {\em parabolic subgroups}. 
We will have occasion to invoke other aspects of the theory of Coxeter groups in some of the proofs that follow. 

{\bf [\S \WeylSection.3:\! Coxeter-group-related notions.]} 
There are several Coxeter-group-related notions that bear mentioning before we proceed further. 
First, for any $\sigma \in \mathcal{W}$, say the the {\em length} of $\sigma$, denoted $\ell(\sigma)$, is the number of `$\mygens_{i}$' generators appearing as factors in any shortest expression for $\sigma$ as a product of generators. 
We take $\ell(\varepsilon)$ to be zero. 
Say `$\mygens_{i_{p}}\cdots \mygens_{i_{2}}\mygens_{i_{1}}$' is a {\em reduced} expression for $\sigma$ if $\sigma = \mygens_{i_{p}}\cdots \mygens_{i_{2}}\mygens_{i_{1}}$ and $\ell(\sigma) = p$. 
Now, for all $i \in I$, note that $\det(S^{*}_{i}) = -1$. 
It is well known (see, for example, \cite{HumCoxeter}) that the function $\mysign: \mathcal{W} \longrightarrow \{\pm 1\}$ given by $\mysign(\sigma) = (-1)^{\ell(\sigma)}$ is a well-defined group homomorphism and that $\mysign(\sigma) = \prod_{q=1}^{p}\det(S^{*}_{i_{q}})$ if $\mygens_{i_{p}} \cdots \mygens_{i_{2}} \mygens_{i_{1}}$ is any expression for $\sigma$ as a product of group generators. 

Second, for each $i \in I$, let ${S}_{i}: \mathfrak{E} \longrightarrow \mathfrak{E}$ be given by ${S}_{i}(\myroot_{j}) := \myroot_{j} - M_{ij}\myroot_{i}$. 
It is easy to check that $({S}_{i})^{2}$ is the identity, so ${S}_{i} \in \mbox{\sffamily Aut}(\mathfrak{E})$. 
Moreoever, the set mapping $\varphi: \{\mygens_{i}\}_{i \in I} \longrightarrow \mbox{\sffamily Aut}(\mathfrak{E})$ given by $\varphi(\mygens_{i}) = {S}_{i}$ induces a faithful representation $\varphi: \mathcal{W} \longrightarrow \mbox{\sffamily Aut}(\mathfrak{E})$, cf.\ \cite{DonRoots}, \cite{BB}. 
Viewing this representation of $\mathcal{W}$ as an action of $\mathcal{W}$ on $\mathfrak{E}$, we have $\mygens_{i}.\myroot_{j} := \myroot_{j} - M_{ij}\myroot_{i}$ for all $i,j \in I$. 
Now consider the natural pairing $\llangle \cdot,\cdot \rrangle: \mathfrak{E}^{*} \times \mathfrak{E} \longrightarrow \mathbb{R}$, wherein $\llangle \mu,v \rrangle := \mu(v)$ for all $\mu \in \mathfrak{E}^{*}$ and $v \in \mathfrak{E}$. 
Of course this pairing is $\mathcal{W}$-invariant: 
$\llangle \mygens_{i}.\omega_{j},\myroot_{k} \rrangle = \llangle \omega_{j}-\delta_{ij}\alpha_{i},\myroot_{k} \rrangle = \delta_{jk}-\delta_{ij}\llangle \alpha_{i},\myroot_{k} \rrangle = \delta_{jk}-\delta_{ji}M_{ik} = \llangle \omega_{j},\myroot_{k} \rrangle - M_{ik}\llangle \omega_{j},\myroot_{i} \rrangle = \llangle \omega_{j},\mygens_{i}.\myroot_{k} \rrangle$. 
Moreover, for all $i \in I$, $v \in \mathfrak{E}$, and $\mu \in \mathfrak{E}^{*}$, we have $\mygens_{i}.v = v-\llangle \alpha_{i},v \rrangle\myroot_{i}$ and $\mygens_{i}.\mu = \mu - \llangle \mu,\myroot_{i} \rrangle\alpha_{i}$. 

Third, let ${\Phi}$ be the set $\{\sigma.\alpha_{i}\}_{i \in I, \sigma \in \mathcal{W}}$. 
We call ${\Phi} = {\Phi}(\mathscr{G})$ the {\em root system} associated with $\mathscr{G}$ and call members of ${\Phi}$ {\em roots}. 
It is known that for any $\alpha \in {\Phi}$, there exist integers $\{k_{i}\}_{i \in I}$ such that $\alpha = \sum_{i \in I}k_{i}\alpha_{i}$ and such that all the $k_{i}$'s are nonnegative (in which case we write $\alpha \in {\Phi}^{+}$ and call $\alpha$ a {\em positive root}) or all are nonpositive (in which case we write $\alpha \in {\Phi}^{-}$ and call $\alpha$ a {\em negative root}). 
Moreover, ${\Phi} = {\Phi}^{+} \disjointunion {\Phi}^{-}$, and for any $\alpha \in {\Phi}$ we have $k\alpha \in {\Phi}$ if and only if $k \in \{\pm 1\}$.  (See \cite{DonRoots} or Ch.\ 4 of \cite{BB} for details.) 
The length function and root system interact in at least the following ways: 
For any $\sigma \in \mathcal{W}$ and $i \in I$, we have $\ell(\sigma \mygens_{i}) > \ell(\sigma)$ if and only if $\sigma.\alpha_{i} \in {\Phi}^{+}$ and $\ell(\sigma \mygens_{i}) < \ell(\sigma)$ if and only if $\sigma.\alpha_{i} \in {\Phi}^{-}$; moreover, $\ell(\sigma) = \rule[-1.5mm]{0.2mm}{5mm}\{\alpha \in {\Phi}^{+}\, |\, \sigma.\alpha \in {\Phi}^{-}\}\rule[-1.5mm]{0.2mm}{5mm}$. 
Analogously, let the {\em toral root system} $\widetilde{\Phi}$ be the set $\{\sigma.\myroot_{i}\}_{i \in I, \sigma \in \mathcal{W}}$ wherein any {\em toral root} $\mysroot \in \widetilde{\Phi}$ is expressible as an integer linear combination of the set of {\em simple toral roots} $\{\myroot_{i}\}_{i \in I}$, and similarly define the {\em positive toral roots} $\widetilde{\Phi}^{+}$ and {\em negative toral roots} $\widetilde{\Phi}^{-}$. 
Our use of the adjective `toral' here anticipates eventual connections with maximal toral subalgebras of semisimple Lie algebras. 

Fourth, for any subset $J \subseteq I$ of our color palette, elements of the set $\mathcal{W}^{J} := \{\sigma \in \mathcal{W}\, |\, \ell(\sigma \mygens_{j}) > \ell(\sigma)\ \forall j \in J\}$ are called {\em minimal coset representatives}. 
Now, for any $\sigma\in\mathcal{W}$, it can be seen there exist unique elements $\sigma_{J} \in \mathcal{W}_{J}$ and $\sigma^{J} \in \mathcal{W}^{J}$ such that $\sigma = \sigma^{J}\sigma_{J}$ with $\ell(\sigma) = \ell(\sigma^{J})+\ell(\sigma_{J})$. 
Therefore, the cosets $\sigma \mathcal{W}_{J}$ and $\tau \mathcal{W}_{J}$ of the parabolic subgroup $\mathcal{W}_{J}$ coincide if and only if $\sigma^{J} = \tau^{J}$, in which case $\ell(\sigma^{J}) \leq \ell(\tau)$. 
This latter inequality demonstrates that elements of $\mathcal{W}^{J}$ are minimum-length representatives of the cosets of $\mathcal{W}_{J}$. 

{\bf [\S \WeylSection.4:\! Weights.]} 
Elements of the $\mathbb{Z}$-span of our fundamental basis are traditionally called {\em weights}, and the elements of $\mathscr{B}$ are {\em fundamental weights}. 
Let us now consider the orbits of weights under the action of the Weyl group, an understanding of which is a crucial part of the theory of certain ``$\mathcal{W}(\mathscr{G})$-symmetric functions''. 
As an example, let $\mathscr{G}$ be the three-cycle \parbox{1.7cm}{\begin{center}
\setlength{\unitlength}{0.2cm}
\begin{picture}(6.5,1) 
\put(1,0){\circle*{0.6}}
\put(3,2){\circle*{0.6}} 
\put(5,0){\circle*{0.6}}
\put(1,0){\qbezier(0,0)(1,1)(2,2)} 
\put(5,0){\qbezier(0,0)(-1,1)(-2,2)} 
\put(1,0){\qbezier(0,0)(2,0)(4,0)} 
\put(-1,-0.25){\footnotesize $\gamma_{2}$} 
\put(5.75,-0.25){\footnotesize $\gamma_{3}$} 
\put(3.75,2){\footnotesize $\gamma_{1}$} 
\end{picture} \end{center}}. 
The Weyl group $\mathcal{W}(\mathscr{G})$ has generators $\mygens_{1}, \mygens_{2}, \mygens_{3}$. 
Note that the orbit of $\mu = \omega_{1} - \omega_{3}$ is the set $\{\omega_{1} - \omega_{3}, -\omega_{1} + \omega_{2}, -\omega_{2} + \omega_{3}\}$. 
Of particular interest to us is the potential finiteness of orbits of certain special weights: Say a weight $\lambda = \sum_{i \in I}\lambda_{i}\omega_{i}$ is {\em dominant} if each $\lambda_{i} \geq 0$ and denote by $\Lambda^{+}$ the cone of dominant weights. 
This cone plays an important role in the study of algebraic objects related to IEG's,\footnote{For example, $\Lambda^{+}$ is a fundamental domain for the action of the related Weyl group cf.\ \cite{HumCoxeter}, and the dominant weights are in one-to-one correspondence with the irreducible representations of the associated semisimple Lie algebra when $\mathscr{G}$ is an integral Coxeter--Dynkin posy.} and it will be crucially featured in what follows. 

For a dominant weight $\lambda$, let $J^{c} := \{j \in I\, |\, \llangle \lambda,\myroot_{j} \rrangle > 0\}$. 
That is, the set $J = I \setminus J^{c}$ identifies the fundamental basis vectors that are \underline{not} utilized in the $\mathbb{Z}$-linear combination comprising $\lambda$. 
We say that $\lambda$ is $J^{c}$-{\em dominant}.  
It is a fact that the parabolic subgroup $\mathcal{W}_{J}$ is the stabilizer of any $J^{c}$-dominant weight $\lambda$ and that for all weights $\mu$ in the $\mathcal{W}$-orbit of $\lambda$ there exists a unique $\sigma \in \mathcal{W}^{J}$ such that $\mu = \sigma.\lambda$. 
(For reasons why these facts hold within our integral embryophytic setting, see the first paragraph of \S 4 of \cite{DonRoots}.) 
That is, we have a natural one-to-one correspondence between $\mathcal{W}\lambda$ and $\mathcal{W}^{J}$.  

{\bf [\S \WeylSection.5:\! The Networked-numbers Game.]} 
K.\ Eriksson was one of the first to explore how the Weyl group action on weights can be viewed as a one-player game, which we call here the {\em Networked-numbers Game} or {\em NG} (see \cite{ErikssonThesis}, \cite{ErikssonDiscrete}, \cite{ErikssonEur}). 
For us, the NG will help frame some of our discussion and, mainly, facilitate a number of our proofs. 
Given any weight $\mu = \sum_{j \in I}\mu_{j}\omega_{i}$ in $\Lambda$, say the action of some generator $\mygens_{i}$ on $\mu$ is {\em legal} if $\mu_{i} > 0$.  
Regardless of legality, we call the application of $\mygens_{i}$ to $\mu$ the $i^{\mbox{\tiny th}}$ {\em firing move}. 
We imagine each $\mu_{j}$ as a number assigned to node $\gamma_{j}$ of the vascular graph $\Gamma(\mathscr{G})$, so the $i^{\mbox{\tiny th}}$ firing move transforms these numbers by the rule $\mu_{j} \mapsto \mu_{j}-M_{ij}\mu_{i}$.  
More generally, say a sequence of nodes $(\gamma_{i_{1}},\gamma_{i_{2}},\ldots,\gamma_{i_{p}})$ is a {\em legal firing sequence from} $\mu$ if, for each $1 \leq q \leq p$, the $i_{q}^{\mbox{\tiny th}}$ firing move is legal when applied to $\mygens_{i_{q-1}}\cdots \mygens_{i_{2}}\mygens_{i_{1}}.\mu$. 
Notice that the latter condition is equivalent to $0 < \llangle \mygens_{i_{q-1}}\cdots \mygens_{i_{2}}\mygens_{i_{1}}.\mu,\myroot_{i_{q}} \rrangle = \llangle \mu,\mygens_{i_{1}}\mygens_{i_{2}}\cdots \mygens_{i_{q-1}}.\myroot_{i_{q}} \rrangle$ by $\mathcal{W}$-invariance of the pairing $\llangle \cdot,\cdot \rrangle$. 
Suppose now that $\mu$ is dominant.  
Let $\mysroot_{q} := \mygens_{i_{1}}\mygens_{i_{2}}\cdots \mygens_{i_{q-1}}.\myroot_{i_{q}}$ for each $1 \leq q \leq p$. 
Then the firing sequence $(\gamma_{i_{1}},\gamma_{i_{2}},\ldots,\gamma_{i_{p}})$ from $\mu$ is legal if and only if $0 < \llangle \mu,\mysroot_{q} \rrangle$ for each $1 \leq q \leq p$. 
Clearly, then, each $\mysroot_{q}$ is in $\widetilde{\Phi}^{+}$, and it follows from our discussion above about lengths of Weyl group elements that $\mygens_{i_{q}} \cdots \mygens_{i_{2}}\mygens_{i_{1}}$ is a reduced expression.

To play the Networked-numbers Game, a hypothetical player chooses a weight $\mu$ as an initial position and then applies legal firing moves of their choosing until the numbers at all of the nodes are nonpositive. 
The complete (possibly infinite and possibly empty) sequence $(\gamma_{i_{1}}, \gamma_{i_{2}}, \gamma_{i_{3}}, \ldots)$ of fired nodes is the {\em game sequence}. 
If the sequence is finite with $p$ node firings, then the game has {\em length} $p$, and $\mygens_{i_{p}}\cdots \mygens_{i_{2}}\mygens_{i_{1}}.\mu$ is the {\em terminal position} for the game. 
In this case, the weight $-\mygens_{i_{p}}\cdots \mygens_{i_{2}}\mygens_{i_{1}}.\mu$ is necessarily dominant. 
Of course, an empty game sequence has {\em length zero}, and a non-terminating game sequence has {\em infinite length}. 
For a nonnegative integer $q$ and a given {\em prefix} $(\gamma_{i_{1}}, \gamma_{i_{2}}, \ldots, \gamma_{i_{q}})$ of a game sequence $(\gamma_{i_{1}}, \gamma_{i_{2}}, \gamma_{i_{3}}, \ldots)$ played from initial position $\mu$, we call $\mygens_{i_{q}}\cdots \mygens_{i_{2}}\mygens_{i_{1}}.\mu$ an {\em intermediate position}. 

The following remarkable result from \cite{ErikssonThesis} (cf.\ \cite{ErikssonEur}) nicely constrains the possible lengths of games played on a connected IEG and will be invoked as needed in proofs of the succeeding results.  

\noindent 
{\bf Eriksson's Strong Convergence Theorem}\ \ {\sl Let $\mathscr{G}$ be a connected IEG. 
In Networked-numbers Game play, any two game sequences played from the same initial position must have the same length, and if the two game sequences are finite then they have the same terminal position.}  

{\bf [\S \WeylSection.6:\! Finite orbits.]} The next result is a Dynkin-diagram classification that does not seem to have appeared in a similarly explicit form elsewhere. 

\noindent 
{\bf \OrbitProposition}\ \ {\sl Let $\mathscr{G}$ be a connected IEG. 
The $\mathcal{W}(\mathscr{G})$-orbit of some nonzero dominant weight is finite if and only if $\mathscr{G}$ is an integral Coxeter--Dynkin flower.} 

{\em Proof.} For the `if' direction, we assume $\mathscr{G}$ is an integral Coxeter--Dynkin flower. 
In Proposition 2.3 of \cite{DonConvergeDiverge}, we demonstrated that for any such $\mathscr{G}$ (there called a `connected Dynkin diagram of finite type') and from any given dominant initial position $\lambda$, there is a finite-length game sequence. 
The argument used there is straightforward and entirely combinatorial, but a bit tedious. 
By Eriksson's Strong Convergence Theorem, all game sequences from $\lambda$ have the same length and the same terminal position. 
It follows that the $\mathcal{W}$-orbit of \underline{any} nonzero dominant $\lambda$ is finite. 

For the `only if' part of the theorem statement, we argue the contrapositive. 
So, let $\lambda$ be any nonzero dominant weight and suppose $\mathscr{G}$ is not a Coxeter--Dynkin flower. 
Now, it is easy to see that if $\mathscr{H}$ is any sub-embryophyte of $\mathscr{G}$, then via some legal firing sequence from initial position $\lambda$ we can produce nonnegative numbers on each of the nodes of $\mathscr{H}$ with at least one such number being positive. 
Thus, if $M_{ij}M_{ji} \geq 4$ for some pair of adjacent nodes $\gamma_{i}$ and $\gamma_{j}$, then game play will eventually result  in nonnegative numbers at nodes $\gamma_{i}$ and $\gamma_{j}$ with at least one of these numbers being positive. 
Without loss of generality, say node $\gamma_{i}$ has a positive number at this point in game play. 
It is clear that any finite prefix of the infinite `alternating' sequence of nodes $(\gamma_{i},\gamma_{j},\gamma_{i},\gamma_{j},\ldots)$ played from this initial position on the two-node embryophyte determined by nodes $\gamma_{i}$ and $\gamma_{j}$ is legal.  
Moreover, we see that for any positive integers $p$ and $q$ with $q < p$, the length $p$ prefix results in an intermediate position on $\mathscr{G}$ that is distinct from the intermediate position produced by the length $q$ sequence. 
That is, our game play has produced infinitely many intermediate positions and therefore infinitely many weights in the $\mathscr{W}$-orbit of $\lambda$. 

So now suppose $0 \leq M_{ij}M_{ji} < 4$ for any two distinct nodes $\gamma_{i}$ and $\gamma_{j}$. 
Since $\mathscr{G}$ is not, by hypothesis, a Coxeter--Dynkin flower, then $\mathscr{G}$ must have an embryophytic subgraph $\mathscr{H}$ that is isomorphic to one of the IEG's from Figure 3.1 of \cite{DonConvergeDiverge}. 
Let $H$ be the subset of $I$ corresponding to the nodes of $\mathscr{H}$. 
As in the preceding paragraph, game play from the initial position $\lambda$ will eventually result in nonnegative numbers at the nodes of $\mathscr{H}$ with at least one of these numbers being positive. 
Say a dominant $\lambda' = \sum_{i \in I}\lambda'_{i}\omega_{i}$ is $H$-dominant if $\lambda'_{i}=0$ for all $i \in I \setminus H$. 
We can regard any such $H$-dominant weight as an initial position on $\mathscr{H}$ for playing the NG. 
To show that the $\mathcal{W}(\mathscr{G})$-orbit of $\lambda$ is infinite, it now suffices to obtain, for each nonzero and $H$-dominant initial position on $\mathscr{H}$, a non-terminating game sequence that passes through infinitely many distinct intermediate positions. 

To further refine our strategy, we claim it suffices to find a non-terminating game sequence that passes through infinitely many distinct intermediate positions when the initial position is any one of the fundamental weights $\{\omega_{h}\}_{h \in H}$. 
To demonstrate this claim, we offer reasons pertaining to $\mathscr{G}$ that will subsequently apply to $\mathscr{H}$ as well. 
Suppose $\mu=\sum_{i \in I}\mu_{i}\omega_{i}$ and $\nu=\sum_{i \in I}\nu_{i}\omega_{i}$ are dominant weights in $\Lambda$ with $\nu_{i}\leq\mu_{i}$ for all $i \in I$. 
Further suppose that some firing sequence $(\gamma_{i_{1}},\gamma_{i_{2}},\ldots,\gamma_{i_{p}})$ is legal from initial position $\nu$. 
We aim to show that this firing sequence is legal from $\mu$ as well, an observation that is called `Eriksson's Comparison Theorem' in \cite{DonConvergeDiverge}. 
Now, for all $1 \leq q \leq p$, we have $0 < \llangle \mygens_{i_{q-1}}\cdots \mygens_{i_{2}}\mygens_{i_{1}}.\nu,\myroot_{i_{q}} \rrangle = \llangle \nu,\mygens_{i_{1}}\mygens_{i_{2}}\cdots \mygens_{i_{q-1}}.\myroot_{i_{q}} \rrangle$. 
Let $\mysroot_{q} := \mygens_{i_{1}}\mygens_{i_{2}}\cdots \mygens_{i_{q-1}}.\myroot_{i_{q}} = \sum_{j \in I}k^{(q)}_{j}\myroot_{j}$, so $0 < \llangle \nu,\mysroot_{q} \rrangle = \sum_{i,j}\nu_{i}k^{(q)}_{j}\llangle \omega_{i},\myroot_{j} \rrangle = \sum_{i} \nu_{i}k^{(q)}_{i} \leq \sum_{i} \mu_{i}k^{(q)}_{i} = \sum_{i,j}\mu_{i}k^{(q)}_{j}\llangle \omega_{i},\myroot_{j} \rrangle = \llangle \mu,\mysroot_{q} \rrangle = \llangle \mu,\mygens_{i_{1}}\mygens_{i_{2}}\cdots \mygens_{i_{q-1}}.\myroot_{i_{q}} \rrangle = \llangle \mygens_{i_{q-1}}\cdots \mygens_{i_{2}}\mygens_{i_{1}}.\mu,\myroot_{i_{q}} \rrangle$. 
So, the firing sequence $(\gamma_{i_{1}},\gamma_{i_{2}},\ldots,\gamma_{i_{p}})$ is legal from initial position $\mu$ as well. 
Now suppose that for some infinite game sequence $(\gamma_{i_{1}},\gamma_{i_{2}},\ldots)$ from $\nu$, it is the case that for some fixed $j \in I$ the sequence $(\llangle \mygens_{i_{p-1}}\cdots \mygens_{i_{2}}\mygens_{i_{1}}.\nu,\myroot_{j} \rrangle)_{p \geq 1,i_{p}=j}$, which consists of the positive numbers at node $\gamma_{j}$ when $\gamma_{j}$ is fired, is nondecreasing and unbounded. 
In this case, we say the pairing of the given game sequence together with the initial position $\nu$ has property $(\star)$. 
It follows from our prior work in this paragraph that $(\gamma_{i_{1}},\gamma_{i_{2}},\ldots)$ is also an infinite game sequence from $\mu$ and that the sequence $(\llangle \mygens_{i_{p-1}}\cdots \mygens_{i_{2}}\mygens_{i_{1}}.\mu,\myroot_{j} \rrangle)_{p \geq 1,i_{p}=j}$ is nondecreasing and unbounded, i.e.\ our given game sequence with initial position $\mu$ has property $(\star)$. 
In particular, the number of intermediate positions for the game sequence is infinite when played from either $\nu$ or $\mu$. 

The reasoning of the preceding paragraph also applies to $\mathscr{H}$. 
So, to see that $\mathcal{W}(\mathscr{G})\lambda$ is infinite, it is enough to find, for each $h \in H$, an infinite game sequence played from initial position $\omega_{h}$ that has property $(\star)$. 
This is exactly what we did in the proof of Proposition 3.1 of \cite{DonConvergeDiverge}, which suffices to complete this proof.\hfill\QED 

{\bf [\S \WeylSection.7:\! Finite Weyl groups and finite root systems.]} An easy corollary is the following well-known classification, although our take on what we call the `toral root system' is non-standard. 

\noindent 
{\bf \FiniteWeylGroupResult}\ \ {\sl Let $\mathscr{G}$ be a connected IEG. 
The Weyl group $\mathcal{W}(\mathscr{G})$ is finite if and only if the root system ${\Phi}(\mathscr{G})$ is finite if and only if the toral root system $\widetilde{\Phi}(\mathscr{G})$ is finite if and only if $\mathscr{G}$ is an integral Coxeter--Dynkin flower.} 

{\em Proof.} We precede the proof with two general observations that both hold for \underline{any} IEG. 
Let $\varrho := \sum_{i \in I}\omega_{i}$, and let $\mygens_{i_{p}} \cdots \mygens_{i_{12}}\mygens_{i_{1}}$ be a reduced expression for some member of $\mathcal{W} = \mathcal{W}(\mathscr{G})$. 
We aim to show that $(\gamma_{i_{1}},\gamma_{i_{2}},\ldots,\gamma_{i_{p}})$ is a legal NG firing sequence from the initial position $\varrho$. 
Since $\mygens_{i_{1}} \mygens_{i_{2}} \cdots \mygens_{i_{p}}$ is also reduced, then so is any sub-expression $\mygens_{i_{1}} \mygens_{i_{2}} \cdots \mygens_{i_{q}}$ when $1 \leq q \leq p$. 
In particular, $q-1 = \ell(\mygens_{i_{1}} \mygens_{i_{2}} \cdots \mygens_{i_{q-1}}) < \ell(\mygens_{i_{1}} \mygens_{i_{2}} \cdots \mygens_{i_{q-1}} \mygens_{i_{q}}) = q$, and hence $\mygens_{i_{1}} \mygens_{i_{2}} \cdots \mygens_{i_{q-1}} \mygens_{i_{q}}.\myroot_{i_{q}} \in \widetilde{\Phi}^{+}$. 
Therefore, $\llangle \mygens_{i_{q-1}} \cdots \mygens_{i_{2}}\mygens_{i_{1}}.\varrho,\myroot_{i_{q}} \rrangle = \llangle \varrho,\mygens_{i_{1}} \mygens_{i_{2}} \cdots \mygens_{i_{q-1}}.\myroot_{i_{q}} \rrangle = \llangle \sum_{i \in I}\omega_{i},\sum_{j \in I}k_{j}\myroot_{j} \rrangle = \sum_{i \in I}k_{i} > 0$ since the positive root $\mygens_{i_{1}} \mygens_{i_{2}} \cdots \mygens_{i_{q-1}}.\myroot_{i_{q}} = \sum_{j \in I}k_{j}\myroot_{j}$ has each $k_{j} \geq 0$ for all $j \in I$ and $k_{j} > 0$ for some $j \in I$. 
That is, for each $1 \leq q \leq p$, the firing sequence $(\gamma_{i_{1}},\gamma_{i_{2}},\ldots,\gamma_{i_{q}})$ from initial position $\varrho$ is legal. 

Our second general observation is that for given colors $i,j \in I$ and Weyl group element $\sigma \in \mathcal{W}$, we have, by Theorem 3.2 of \cite{DonEur}, $\sigma.\alpha_{i} = \alpha_{j}$ if and only if $\sigma.\myroot_{i} = \myroot_{j}$. 
Put another way, for given colors $i,j \in I$ and Weyl group elements $\sigma, \tau \in \mathcal{W}$, we have $\sigma.\alpha_{i} = \tau.\alpha_{j}$ if and only if $\sigma.\myroot_{i} = \tau.\myroot_{j}$. 
In particular, the root system ${\Phi} = {\Phi}(\mathscr{G}) = \{\sigma.\alpha_{i}\}_{\sigma \in \mathcal{W}, i \in I}$ is finite if and only if the toral root system $\widetilde{\Phi} = \widetilde{\Phi}(\mathscr{G}) = \{\sigma.\myroot_{i}\}_{\sigma \in \mathcal{W}, i \in I}$ is finite. 

To formally begin our proof of \FiniteWeylGroupResult, suppose $\mathscr{G}$ is an integral Coxeter--Dynkin flower. 
From Proposition 2.3 of \cite{DonConvergeDiverge}, we know that there exists a finite game sequence for NG play from initial position $\varrho$. 
By Eriksson's Strong Convergence Theorem, all game sequences played from $\varrho$ must have the same length and the same terminal position. 
By the previous paragraph, the factors of any reduced expression for an element of $\mathcal{W}$ correspond to a legal firing sequence from $\lambda$. 
So, there can be no reduced expression in $\mathcal{W}$ whose length is greater than the length of the game sequences for $\varrho$, which forces $\mathcal{W}$ to be finite.
Of course, $\mathcal{W}$ finite $\Longrightarrow$ the root system ${\Phi}$ is finite. 

So, now suppose the toral root system $\widetilde{\Phi}$ is finite. 
Let $\mu = \sum_{i \in I}\mu_{i}\omega_{i}$ be any member of the $\mathcal{W}$-orbit of $\varrho$, so $\mu = \mygens_{i_{p}} \cdots \mygens_{i_{2}} \mygens_{i_{1}}.\varrho$ for some reduced expression $\mygens_{i_{p}} \cdots \mygens_{i_{2}} \mygens_{i_{1}}$. 
Then for each $i \in I$, we have $\mu_{i} = \llangle \mu,\myroot_{i} \rrangle = \llangle \varrho,\mygens_{i_{1}}\mygens_{i_{2}}\cdots \mygens_{i_{p}}.\myroot_{i} \rrangle$, where $\mygens_{i_{1}}\mygens_{i_{2}}\cdots \mygens_{i_{p}}.\myroot_{i}$ is a toral root in $\widetilde{\Phi}$. 
Since $\widetilde{\Phi}$ is finite, then the number of possible $\mu_{i}$'s is finite for each $i \in I$, and therefore the orbit of $\lambda$ must be finite. 
By \OrbitProposition, $\mathscr{G}$ must be a Coxeter--Dynkin flower.\hfill\QED

{\bf [\S \WeylSection.8:\! Parabolic subgroups of finite index.]} Another easy corollary is the following result, first proved by Hosaka in \cite{Hosaka} using ideas from geometric group theory. 
We report this corollary, and its new proof, for the record but will not require it for any subsequent work. 

\noindent
{\bf \FiniteIndexCorollary}\ \ {\sl Let $\mathscr{G}$ be a connected IEG. 
There exists a proper embryophytic subgraph $\mathscr{H}$ of $\mathscr{G}$ such that $\mathcal{W}(\mathscr{H})$ has finite index as a subgroup of $\mathcal{W}(\mathscr{G})$ if and only if $\mathscr{G}$ is an integral Coxeter--Dynkin flower.}

{\em Proof.} The `if' part of the statement follows from the fact that for any Coxeter--Dynkin flower $\mathscr{G}$, the Weyl group $\mathcal{W}(\mathscr{G})$ is finite, cf.\ \FiniteWeylGroupResult. 
For the `only if' part, let $H$ denote the proper subset of our color palette $I$ corresponding to the nodes of the vascular graph of $\mathscr{H}$. 
Let $W := \mathcal{W}(\mathscr{G})$ and $W_{H} := \mathcal{W}(\mathscr{H})$.  
Let $\lambda = \sum_{i \in I}\lambda_{i}\omega_{i}$ be `$H^{c}$-dominant' in the sense that $\lambda_{i} > 0$ for all $i \in I \setminus H$ and $\lambda_{h} = 0$ for all $h \in H$. 
Clearly $W_{H}$ stabilizes $\lambda$.  
Now let $W^{H}$ denote the set of minimal coset representatives of $W_{H}$, cf.\ \cite{HumCoxeter}. 
This set consists of those $\sigma \in W$ such that $\ell(\sigma \mygens_{h}) > \ell(\sigma)$ for all $h \in H$ and is in one-to-one correspondence with the cosets of $W_{H}$. 
For each $\sigma \in W$, there exist unique $\sigma^{H} \in W^{H}$ and $\sigma_{H} \in W_{H}$ such that $\sigma = \sigma^{H}\sigma_{H}$. 
Of course $\sigma.\lambda = \sigma^{H}.\lambda$. 
Let $\sigma_{1}, \sigma_{2} \in W$. 
Partition $W^{H}$ by the equivalence relation $\sigma_{1}^{H} \sim \sigma_{2}^{H}$ if and only if $\sigma_{1}^{H}.\lambda = \sigma_{2}^{H}.\lambda$. 
Evidently, the $W$-orbit $W\lambda$ of $\lambda$ is in one-to-one correspondence with the set $W^{H}/\sim$ of equivalence classes under `$\sim$'. 
So, $\rule[-1.5mm]{0.2mm}{5mm}W\lambda\rule[-1.5mm]{0.2mm}{5mm} = \rule[-1.5mm]{0.2mm}{5mm}W^{H}/\sim\rule[-1.5mm]{0.2mm}{5mm} \leq \rule[-1.5mm]{0.2mm}{5mm}W^{H}\rule[-1.5mm]{0.2mm}{5mm} = \rule[-1.5mm]{0.2mm}{5mm}W/W_{H}\rule[-1.5mm]{0.2mm}{5mm} < \infty$, where the latter inequality follows from our hypothesis that $W_{H}$ has finite index in $W$. 
We may therefore invoke \OrbitProposition\ to conclude that $\mathscr{G}$ is a Coxeter--Dynkin flower.\hfill\QED

{\bf [\S \WeylSection.9:\! Weyl symmetric functions.]} Associate to any $\mu = \sum_{i \in I}\mu_{i}\omega_{i}$ in $\Lambda$ a Laurent monomial $\myvarZ^{\mu} := \prod_{i \in I}z_{i}^{\mu_{i}}$ in some indeterminates $\{z_{i}\}_{i \in I}$. 
Observe that each Weyl group generator  acts on these monomials by declaring that $\mygens_{i}.\myvarZ^{\mu} := \myvarZ^{\mysmallgens_{i}.\mu}$ and then extending in the natural way to an action of $\mathcal{W}(\mathscr{G})$ on the ring `$\mathbb{Z}[\Lambda]$' of all Laurent polynomials $\sum_{\mu \in \Lambda}c_{\mu}\myvarZ^{\mu}$ where at most finitely many of the integer $c_{\mu}$'s are nonzero. 
Of particular interest are certain Laurent polynomials that are invariant under the action of the Weyl group. 
We restrict our attention to $\mathcal{W}(\mathscr{G})$-invariant Laurent polynomials $\chi = \sum_{\mu \in \Lambda}c_{\mu}\myvarZ^{\mu}$ such that, if $\chi$ is not constant, there exists a nonzero dominant weight $\lambda$ with $c_{\lambda} \not= 0$. 
Call any such $\chi$ a $\mathcal{W}(\mathscr{G})${\em -symmetric function}. 
The requirement that a $\mathcal{W}(\mathscr{G})$-symmetric function should have a monomial term with exponents from a nonzero dominant weight may seem {\em ad hoc} but, as we asserted in \S \WeylSection.4, is natural within our context, since the cone $\Lambda^{+}$ of dominant weights plays such a large role in related theories.   
The following result -- yet another Dynkin diagram classification -- is not surprising but does not seem to have been stated formally elsewhere. 
 
\noindent
{\bf \NonconstantTheorem}\ \ {\sl Let $\mathscr{G}$ be a connected IEG.  
There exists a non-constant $\mathcal{W}(\mathscr{G})$-symmetric function if and only if $\mathscr{G}$ is an integral Coxeter--Dynkin flower.} 

{\em Proof.} First we consider the `only if' part of the statement. 
For our non-constant $\mathcal{W}(\mathscr{G})$-symmetric function $\sum_{\mu \in \Lambda}c_{\mu}\myvarZ^{\mu}$, there exists a nonzero dominant weight $\lambda$ with $c_{\lambda} \not= 0$. 
Then, for each $\sigma \in \mathcal{W}(\mathscr{G})$, the coefficient $c_{\sigma.\lambda} = c_{\lambda}$ is also nonzero. 
Finiteness of the symmetric function forces the $\mathcal{W}(\mathscr{G})$-orbit of $\lambda$ to be finite. 
It follows by \OrbitProposition\ that $\mathscr{G}$ is a Coxeter--Dynkin flower. 
For the `if' part of the statement, when $\mathscr{G}$ is a Coxeter--Dynkin flower, then $\mathcal{W} = \mathcal{W}(\mathscr{G})$ is finite by \FiniteWeylGroupResult. 
Let $\lambda$ be any nonzero dominant weight.  
Then the finite sum $\sum_{\sigma \in \mathcal{W}}\myvarZ^{\sigma.\lambda}$ is clearly non-constant and $\mathcal{W}$-invariant and therefore a $\mathcal{W}$-symmetric function.\hfill\QED 

{\bf [\S \WeylSection.10:\! A Weyl-group-invariant inner product.]} We now wish to endow $\mathfrak{E}^{*}$ with an inner product that is natural, at least in some sense, to our set-up. 
Here is a classical approach: 
We wish to view each Weyl group generator $S^{*}_{j}$ as a reflection with respect to $\alpha_{j}$. 
In general, for a nonzero vector $v$ from a Euclidean space with inner product $B(\cdot,\cdot)$, the corresponding reflection $\sigma_{v}$ is given by $\sigma_{v}(w) = w - \left(\frac{2B(w,v)}{B(v,v)}\right)v$. 
That is, $\sigma_{v}(w) = w - B(w,v^{\vee}) v$ with $v^{\vee} := \frac{2}{B(v,v)} v$. 
Therefore, in the  embryophytic setting of $\mathfrak{E}^{*}$, our desired inner product must satisfy $B(\omega_{i},\alpha_{j}^{\vee}) = \delta_{ij}$ for all $i,j \in I$, since we require that $S^{*}_{j}(\omega_{i}) = \omega_{i} - B(\omega_{i},\alpha_{j}^{\vee})\alpha_{j}$. 
For similar reasons, for all $i,j \in I$, we must also have $B(\alpha_{i},\alpha_{j}^{\vee}) = M_{ij}$ and hence $M_{ij}B(\alpha_{j},\alpha_{j}) = M_{ji}B(\alpha_{i},\alpha_{i})$. 
Since the $\alpha_{j}^{\vee}$'s will be dual to the fundamental basis $\mathscr{B} = \{\omega_{i}\}_{i \in I}$, then the set $\{\alpha_{j}^{\vee}\}_{j \in I}$ must be linearly independent, so the set of simple roots $\{\alpha_{i}\}_{i \in I}$ must be as well. 
From this, it is clear that, with respect to the fundamental basis, the desired inner product must have positive definite and symmetric representing matrix $S = DM$, where $D$ is the diagonal matrix with $D_{ii} = \mbox{\small $2/B(\alpha_{i},\alpha_{i})$}$ for each $i \in I$. 

Therefore, the existence of a positive definite diagonal matrix $D$ such that $DM$ is symmetric and positive definite is a necessary condition for the inner product we seek. 
It is also sufficient, which we can see as follows. 
Let us regard the simple roots $\alpha_{i}$ (each of which is necessarily nonzero), weights $\mu$, and elements of $\mathfrak{E}^{*}$ to be row vectors with respect to the fundamental basis $\mathscr{B} = \{\omega_{i}\}_{i \in I}$. 
That is, if $\mu = \sum_{i \in I}\mu_{i}\omega_{i}$, then $[\mu]_{\mathscr{B}} = [\mu_{i}]_{i \in I}$. 
If $S = DM$ is positive definite and symmetric for a diagonal matrix $D$ with positive diagonal entries, then we let our inner product $B(\cdot,\cdot)$ be defined by the rule $B(v,w) := [v]_{\mathscr{B}}S^{-1}[w]_{\mathscr{B}}^{\mytinyT}$ for all $v = \sum_{i \in I}v_{i}\omega_{i}$ and $w = \sum_{i \in I}w_{i}\omega_{i}$ in $\mathfrak{E}^{*}$. 
Note that $B(\omega_{i},\alpha_{j}) = B(\alpha_{j},\omega_{i}) = [\alpha_{j}]_{\mathscr{B}}M^{-1}D^{-1}[\omega_{i}]_{\mathscr{B}}^{\mytinyT} = D_{ii}^{-1}[\alpha_{j}]_{\mathscr{B}}M^{-1}[\omega_{i}]_{\mathscr{B}}^{\mytinyT} = D^{-1}_{ii}\delta_{ij}$, in which case $B(\alpha_{i},\alpha_{i}) = B(\sum_{k \in I}M_{ik}\omega_{k},\alpha_{i}) = D_{ii}^{-1}M_{ii} = 2D_{ii}^{-1}$ or simply $D_{ii} = 2/B(\alpha_{i},\alpha_{i})$. 
Then $B(\mygens_{i}.v,w) = B(v-v_{i}\alpha_{i},w) = B(v,w) - v_{i}B(\alpha_{i},\sum_{k}w_{k}\omega_{k}) = B(v,w)-v_{i}w_{i}D_{ii}^{-1} = B(v,w) - w_{i}B(\sum_{k \in I}v_{k}\omega_{k},\alpha_{i}) = B(v,w) - w_{i}B(v,\alpha_{i}) = B(v,\mygens_{i}.w)$. 
In particular, $B$ is $\mathcal{W}$-invariant. 
Moreover, $\mygens_{i}.v = v - v_{i}\alpha_{i} = v - D_{ii}B(v_{i}\omega_{i},\alpha_{i}) \alpha_{i} = v - D_{ii}B(\sum_{k \in I}v_{k}\omega_{k},\alpha_{i}) = v - \frac{2B(v,\alpha_{i})}{B(\alpha_{i},\alpha_{i})}\alpha_{i}$.
That is, each $\mygens_{i}$ acts as a reflection through the hyperplane that is orthogonal to $\alpha_{i}$, as desired. 

\noindent 
{\bf \InnerProductTheorem}\ \ {\sl Let $\mathscr{G}$ be a connected IEG. There exists a positive definite diagonal matrix $D$ such that $DM$ is positive definite and symmetric if and only if $\mathscr{G}$ is an integral Coxeter--Dynkin flower. 
In this case, the bilinear form} $\langle v,w \rangle := [v]_{\mathscr{B}}(M^{-1}D^{-1})[w]_{\mathscr{B}}^{\mytinyT}$, {\sl defined for all $v, w \in \mathfrak{E}^{*}$, is a $\mathcal{W}$-invariant inner product on $\mathfrak{E}^{*}$ with respect to which, for each $i \in I$, the Weyl group generator $\mygens_{i}$ acts on $\mathfrak{E}^{*}$ as a reflection in the hyperplane orthogonal to the simple root $\alpha_{i}$.} 

{\em Proof.} For brevity, the notation `$[\mu]$' in the proof means `$[\mu]_{\mathscr{B}}$'. 
The claims made in the second sentence of the corollary statement follow directly from the paragraph preceding the corollary. 
Suppose, now, that there exists a positive definite diagonal matrix $D$ such that $DM$ is positive definite and symmetric; the inverse $M^{-1}D^{-1}$ is positive definite and symmetric as well. 
From our work in the paragraph before the corollary, we see that the quadratic form $\myspecialf:V \longrightarrow \mathbb{R}$ given by $\myspecialf(\mu) := [\mu] M^{-1}D^{-1} [\mu]^{\mytinyT}$ is $\mathcal{W}$-invariant in the sense that $\myspecialf(\sigma.\mu) = \myspecialf(\mu)$ for all $\sigma \in \mathcal{W}$. 
Identify $\mathfrak{E}^{*}$ with $\mathbb{R}^{n}$ by identifying any $\mu \in \mathfrak{E}^{*}$ with the row vector $[\mu]$, and the topologize $\mathfrak{E}^{*}$ accordingly. 
Now fix a nonzero dominant $\lambda \in \Lambda$. 
Since $\mygens_{i}.\lambda = \lambda - \lambda_{i}\alpha_{i}$ for each $i \in I$, it follows that $\sigma.\lambda$ is in the coset $\lambda + \mathbb{Z}\{\alpha_{i}\}_{i \in I}$ of the root lattice $\mathbb{Z}\{\alpha_{i}\}_{i \in I}$. 
This coset is a discrete subspace of $\mathfrak{E}^{*}$ in that each singleton set $\{\lambda + \sum_{i \in I}k_{i}\alpha_{i}\}$, where each $k_{i}$ is some fixed integer, is an open subset in the subspace topology. 
Since our quadratic form $\myspecialf$ is defined by a symmetric and, in particular, positive definite matrix $M^{-1}D^{-1}$, it follows that, for any fixed nonnegative real number $c$, the set $\mathfrak{E}^{*}_{c} := \{v \in \mathfrak{E}^{*}\, |\, \myspecialf(v)=c\}$ is compact in $\mathfrak{E}^{*}$. 
Now let $c = \myspecialf(\lambda)$. 
As the intersection of a compact set and a discrete subspace, $\mathfrak{E}^{*}_{c} \cap (\lambda + \mathbb{Z}\{\alpha_{i}\}_{i \in I})$ must be finite. 
Since the $\mathcal{W}$-orbit of $\lambda$ resides within this intersection, it follows that $\mathcal{W}\lambda$ is finite. 
By \OrbitProposition, it follows that $\mathscr{G}$ is an integral Coxeter--Dynkin flower. 

Now suppose $\mathscr{G}$ is an integral Coxeter--Dynkin flower. 
In the $\myA-\myD-\myE$ cases, the pulsation matrix is $M=M(\mathscr{G})$ is symmetric. 
Positive definiteness in the cases that $\mathscr{G} \in \{\myE_{6},\myE_{7},\myE_{8}\}$ can be checked by hand. 
In the $\myA_{n}$ and $\myD_{n}$ cases, use induction. 
In the $\myF_{4}$ case, $M = \left(\begin{array}{cccc}2 & -1 & 0 & 0\\ -1 & 2 & -2 & 0\\ 0 & -1 & 2 & -1\\ 0 & 0 & -1 & 2\end{array}\right)$. 
Our investigation in the paragraph preceding the corollary implies that, in the presence of an inner product, the simple root lengths must satisfy $\langle \alpha_{1},\alpha_{1} \rangle = \langle \alpha_{2},\alpha_{2} \rangle = 2\langle \alpha_{3},\alpha_{3} \rangle = 2\langle \alpha_{4},\alpha_{4} \rangle$. 
Motivated by the idea that $\langle \alpha_{3},\alpha_{3} \rangle = \langle \alpha_{4},\alpha_{4} \rangle = 1$ and $\langle \alpha_{1},\alpha_{1} \rangle = \langle \alpha_{2},\alpha_{2} \rangle = 2$, we take $D = \left(\begin{array}{cccc}1 & 0 & 0 & 0\\ 0 & 1 & 0 & 0\\ 0 & 0 & 2 & 0\\ 0 & 0 & 0 & 2\end{array}\right)$, so $DM = \left(\begin{array}{cccc}2 & -1 & 0 & 0\\ -1 & 2 & -2 & 0\\ 0 & -2 & 4 & -2\\ 0 & 0 & -2 & 4\end{array}\right)$. 
By inspection, the symmetric matrix $DM$ is positive definite. 
Use a similar idea with $\mathscr{G} = \myB_{n}$: Take $D = \left(\begin{array}{ccccc}1 & 0 & \cdots & 0 & 0\\ 0 & 1 & \cdots & 0 & 0\\ \vdots & \vdots & \ddots & \vdots & \vdots\\ 0 & 0 & \cdots  & 1 & 0\\ 0 & 0 & \cdots  & 0 & 2\end{array}\right)$ and see that $DM = \left(\begin{array}{ccccc}2 & -1 & \cdots & 0 & 0\\ -1 & 2 & \cdots & 0 & 0\\ \vdots & \vdots & \ddots & \vdots & \vdots\\ 0 & 0 & \cdots  & 2 & -2\\ 0 & 0 & \cdots  & -2 & 4\end{array}\right)$ is positive definite and symmetric. 
When $\mathscr{G} = \myC_{n}$, take $D = \left(\begin{array}{ccccc}2 & 0 & \cdots & 0 & 0\\ 0 & 2 & \cdots & 0 & 0\\ \vdots & \vdots & \ddots & \vdots & \vdots\\ 0 & 0 & \cdots  & 2 & 0\\ 0 & 0 & \cdots  & 0 & 1\end{array}\right)$ and see that $DM = \left(\begin{array}{ccccc}4 & -2 & \cdots & 0 & 0\\ -2 & 4 & \cdots & 0 & 0\\ \vdots & \vdots & \ddots & \vdots & \vdots\\ 0 & 0 & \cdots  & 4 & -2\\ 0 & 0 & \cdots  & -2 & 2\end{array}\right)$ is positive definite and symmetric. 
When $\mathscr{G} = \myG_{2}$, take $D = \left(\begin{array}{cc}3 & 0\\ 0 & 1\end{array}\right)$ and see that $DM = \left(\begin{array}{cc}6 & -3\\ -3 & 2\end{array}\right)$ is positive definite and symmetric.\hfill\QED

Of course, integral pulsation matrices satisfying the criteria of \InnerProductTheorem\ are called `Cartan matrices'. 
The classification of indecomposable Cartan matrices that is proffered in the above corollary is well known, but we believe our proof of the `only-if' part of \InnerProductTheorem\ is new, in that it invokes \OrbitProposition. 
For a classical proof, see \cite{Kac}, which draws on work of Vinberg \cite{Vin}. 

{\bf [\S \WeylSection.11:\! Terminating Networked-numbers Games.]} Another consequence of \OrbitProposition\ is \NGCorollary. 
This classification is really a corollary of the {\sl proof} of \OrbitProposition\ and first appeared\footnote{In the YouTube video of an August 9, 2018 Univ.\ of New South Wales lecture on ``$\myA\myD\myE$ graphs and lattices, weights and quarks'', N.\ J.\ Wildberger attributes to K.\ Eriksson a version of this result for ``nonzero positive weights''.  (Go to time 34:45 of the video.)  
If ``nonzero positive weight'' means ``nonzero dominant weight'', then the result presented by Prof.\ Wildberger is a consequence of \NGCorollary. Also, one-half of the equivalence asserted in our \NGCorollary\ is an existential statement (as in, ``there exists a nonzero dominant initial position...'') rather than a universal statement (as in, ``for every nonzero dominant initial position...'').  The existential formulation of \NGCorollary\ seems to be original to \cite{DonEur}.} in \cite{DonEur}.

\noindent 
{\bf \NGCorollary}\ \ {\sl Let $\mathscr{G}$ be a connected IEG. 
There exists a terminating Networked-numbers Game on $\mathscr{G}$ played from a nonzero dominant initial position if and only if $\mathscr{G}$ is a Coxeter--Dynkin flower.} 

{\em Proof.} For the `if' direction, assume $\mathscr{G}$ is a Coxeter--Dynkin flower. 
Let $\mu := \sum_{i \in I}\omega_{i}$. 
In Proposition 2.3 of \cite{DonConvergeDiverge}, we exhibited, for each Coxeter--Dynkin flower, a finite game sequence for NG play from initial position $\mu$. 
For the converse, we demonstrate the contrapositive. 
So suppose our connected IEG $\mathscr{G}$ is not a Coxeter--Dynkin flower. 
The proof of the `only if' part of \OrbitProposition\ showed that for NG play from any nonzero dominant initial position, there is a non-terminating game. 
By Eriksson's Strong Convergence Theorem, all games played from such a position must therefore be non-terminating. 
That is, if $\mathscr{G}$ is not a Coxeter--Dynkin flower, then there is no nonzero dominant initial position from which there is a terminating game.\hfill\QED

{\bf [\S \WeylSection.12:\! Kac--Moody Lie algebras.]} 
The last two classification results we pursue in this section (\KMCorollaryList) require the most overhead and may be safely skipped if one is not interested in any Lie algebraic aspects of DCML's and DCDL's in later sections. 
We give a very distilled version of the set-up ideas here and refer the interested reader to \cite{Hum}, \cite{Kac}, and \cite{Kumar} for details. 
The following discussion is intended to facilitate the statement of \KacMoodyCorollary, which is that the only finite-dimensional Kac--Moody Lie algebras obtained from connected IEG's are those obtained from Coxeter--Dynkin flowers. 
This claim is demonstrated here as a consequence of \FiniteWeylGroupResult.

Until we declare otherwise, let $\mathbb{F}$ be any field. 
A {\em Lie algebra over} $\mathbb{F}$ is an $\mathbb{F}$-vector space $\mathfrak{g}$ together with a bilinear operation $[\cdot,\cdot]: \mathfrak{g} \times \mathfrak{g} \longrightarrow \mathfrak{g}$, called the {\em Lie bracket}, that is {\em one-dimensionally commutative} (i.e.\ $[x,x]=0$ for all $x \in \mathfrak{g}$) and {\em Jacobi-associative} (i.e.\ $[x,[y,z]]+[y,[z,x]]+[z,[x,y]]=0$ for all $x,y,z \in \mathfrak{g}$). 
Notice that one-dimensional commutativity implies {\em anti-commutativity} (i.e.\ $[x,y] = -[y,x]$ for all $x,y \in \mathfrak{g}$); anti-commutativity implies one-dimensional commutativity when $\mbox{\sffamily char}(\mathbb{F}) \not= 2$. 
A {\em Lie subalgebra} of $\mathfrak{g}$ is a subspace that is closed under the Lie bracket; a {\em Lie ideal} of $\mathfrak{g}$ is a subspace $\mathfrak{i}$ that is `i-o closed' under the bracket, i.e.\ $[x,y] \in \mathfrak{i}$ for all $x \in \mathfrak{i}$ and $y \in \mathfrak{g}$.  
By anti-commutativity, a Lie ideal of a Lie algebra is necessarily two-sided.  

The peculiar non-associative property of the Lie bracket is actually natural in certain contexts. 
For example, the real vector space $\mathbb{R}^{3}$ together with the cross product is a real Lie algebra. 
For a less singular example, any $\mathbb{F}$-algebra\footnote{We assume that an $\mathbb{F}$-algebra is unital and associative.} $\mathscr{A}$ becomes a Lie algebra $\myLieop(\mathscr{A})$ when we define its Lie bracket to be the commutator $[x,y] := xy-yx$ for all $x,y \in \mathscr{A}$. 
Let $V$ be a $d$-dimensional $\mathbb{F}$-vector space, and let $\mbox{\sffamily \small END}(V)$ be the $d^{2}$-dimensional $\mathbb{F}$-vector space of linear transformations $V \longrightarrow V$ (i.e.\ endomorphisms of $V$). 
Since $\mbox{\sffamily \small END}(V)$ is an $\mathbb{F}$-algebra, then it becomes a Lie algebra $\mathfrak{gl}(V) := \myLieop\left(\mbox{\sffamily \small END}(V)\right)$ when we take its Lie bracket to be the commutator $[S,T] := ST-TS$, where, for linear transformations $S$ and $T$, the quantity `$ST$' means $S \circ T$.  
We call $\mathfrak{gl}(V)$ the {\em general linear Lie algebra over} $V$. 
If we fix a basis for $V$ and identify each linear transformation of $V$ with its $d \times d$ representing matrix in this basis, then we refer to this concretization of the general linear Lie algebra over $\mathbb{F}$ as `$\mathfrak{gl}(d,\mathbb{F})$'. 
Now, the $(d^{2}-1)$-dimensional subspace $\mathfrak{sl}(V)$ of $\mathfrak{gl}(V)$ consisting of the traceless linear transformations \underline{is not} closed under composition. 
However, $\mathfrak{sl}(V)$ \underline{is} closed under the Lie bracket, so it is a Lie subalgebra (actually, a Lie ideal) of $\mathfrak{gl}(V)$, which we call the {\em special linear Lie algebra over} $V$. 
With respect to the concretization $\mathfrak{gl}(d,\mathbb{F})$, we view the special linear Lie algebra as the Lie subalgebra $\mathfrak{sl}(d,\mathbb{F})$ of matrices with trace equal to zero. 

Say $\mathfrak{g}$ and $\mathfrak{h}$ are Lie algebras over a field $\mathbb{F}$. 
A {\em Lie algebra homomorphism} $\phi: \mathfrak{g} \longrightarrow \mathfrak{h}$ from $\mathfrak{g}$ to $\mathfrak{h}$ is an $\mathbb{F}$-linear transformation that preserves Lie brackets. 
The kernel of $\phi$ is a Lie ideal of $\mathfrak{g}$, and the image of $\phi$ is a Lie subalgebra of $\mathfrak{h}$.  
The usual homomorphism theorems hold for Lie algebras. 
A Lie algebra is {\em simple} if it has no Lie ideals besides itself and the trivial `zero suspace' $\{0\}$, so it is isomorphic to its image under any nontrivial homomorphism. 
A {\em representation} of $\mathfrak{g}$ is a Lie algebra homomorphism $\phi$ from $\mathfrak{g}$ to the general linear Lie algebra $\mathfrak{gl}(V)$ over an $\mathbb{F}$-vector space $V$. 
We think of $\mathfrak{g}$ as acting on $V$ by defining $x.v := \phi(x)(v)$ for all $x \in \mathfrak{g}$ and $v \in V$, in which case we call $V$ a $\mathfrak{g}$-{\em module}. 
Suppose a given action `$x.v$' of elements $x \in \mathfrak{g}$ on vectors $v \in V$ is linear and satisfies $[y,z].w = y.(z.w) - z.(y.w)$ for all $y,z \in \mathfrak{g}$ and $w \in V$. 
Then $V$ is a $\mathfrak{g}$-module associated with the representation $\phi: \mathfrak{g} \longrightarrow \mathfrak{gl}(V)$ given by $\phi(x)(v) = x.v$ for all $x \in \mathfrak{g}$ and $v \in V$. 
For fixed $z \in \mathfrak{g}$, we define the linear transformation $\mbox{\sffamily ad}z: \mathfrak{g} \longrightarrow \mathfrak{g}$ by the rule $(\mbox{\sffamily ad}z)(x) := [z,x]$. 

There is a categorical object associated with any Lie algebra $\mathfrak{g}$ over a field $\mathbb{F}$. 
The {\em universal enveloping algebra} of $\mathfrak{g}$ is a pair $(\mathscr{U}(\mathfrak{g}),\myi)$, where $\mathscr{U}(\mathfrak{g})$ is an algebra over $\mathbb{F}$ (and therefore a Lie algebra $\myLieop\left(\mathscr{U}(\mathfrak{g})\right)$ with the commutator as its Lie bracket), $\myi: \mathfrak{g} \longrightarrow \myLieop\left(\mathscr{U}(\mathfrak{g})\right)$ is a Lie algebra homomorphism, and the pair satisfies the universal property depicted below:  
\begin{center}
\setlength{\unitlength}{1cm}
\begin{picture}(4,2)
\thinlines
\put(0,1.5){\Large $\mathfrak{g}$} 
\put(0.75,1.65){$\myi$}
\put(0.375,1.575){\vector(1,0){1}}
\put(1.5,1.5){\Large $\mbox{\large \Fontlukas L}\left(\mathscr{U}(\mathfrak{g})\right)$}
\put(2.3,1.3){\vector(0,-1){0.7}}
\put(2.4,0.85){$\exists!\psi$}
\put(1.8,0){\Large $\mbox{\large \Fontlukas L}\left(\mathscr{A}\right)$}
\put(0.3,1.3){\vector(3,-2){1.4}}
\put(0.4,0.625){$\forall\phi$}
\end{picture}
\end{center}
\noindent 
That is, any Lie algebra homomorphism $\phi: \mathfrak{g} \longrightarrow \myLieop\left(\mathscr{A}\right)$ from $\mathfrak{g}$ to the Lie algebra of an $\mathbb{F}$-algebra $\mathscr{A}$ extends uniquely to a Lie algebra homomorphism $\psi: \myLieop\left(\mathscr{U}(\mathfrak{g})\right) \longrightarrow \myLieop\left(\mathscr{A}\right)$ so that $\phi = \psi \circ \myi$.  
As with many such objects that are characterized by universal mapping properties, the pair $(\mathscr{U}(\mathfrak{g}),\myi)$ is unique, in the obvious sense, if it exists. 
The standard construction of $\mathscr{U}(\mathfrak{g})$ is as the quotient $\mathscr{T}(\mathfrak{g})/\mbox{\Fontauri C}$ of the tensor algebra $\mathscr{T}(\mathfrak{g})$ by the two-sided ideal $\mbox{\Fontauri C}$ generated by elements of the form $x \otimes y - y \otimes x - [x,y]$. 
A non-obvious aspect of this construction is that the function $\myi: \mathfrak{g} \longrightarrow \myLieop\left(\mathscr{U}(\mathfrak{g})\right)$ is injective. 
This latter fact is a corollary of the famous Poincar\'{e}--Birkhoff--Witt Basis Theorem. 

We can use the universal enveloping algebra as follows to define a {\em free Lie algebra} $\mathfrak{f}(\mathcal{S})$ on some given set $\mathcal{S}$. 
Let $V := V_{\mathbb{F}}(\mathcal{S})$ be the $\mathbb{F}$-vector space freely generated by $\mathcal{S}$. 
Note that $V$ is naturally embedded in the tensor algebra $\mathscr{T}(V)$, and therefore in the Lie algebra $\myLieop\left(\mathscr{T}(V)\right)$, via an injective set mapping we name $\epsilon: V \longrightarrow \myLieop\left(\mathscr{T}(V)\right)$. 
Declare $\mathfrak{f}(\mathcal{S})$ to be the smallest Lie subalgebra of $\myLieop\left(\mathscr{T}(V)\right)$ that contains $\mathcal{S}$, so $\epsilon(V) \subseteq \mathfrak{f}(\mathcal{S})$. 
We would like any set mapping from $\mathcal{S}$ to a given Lie algebra $\mathfrak{g}$ to extend uniquely to a Lie algebra homomorphism from $\mathfrak{f}(\mathcal{S})$ to $\mathfrak{g}$. 
To see this, we consider the following depiction of our set-up so far: 
\begin{center}
\hspace*{-0.35in}
\setlength{\unitlength}{1cm}
\begin{picture}(8,2)
\thinlines
\put(0,1.5){\Large $\mathcal{S}$}
\put(0.375,1.625){\vector(1,0){1.5}}
\put(0.6,1.7){\tiny inclusion}
\put(0.375,1.4){\vector(3,-2){1.5}}
\put(0.6,0.625){$\forall\phi$}
\put(2,1.5){\Large $V$} 
\put(2.375,1.625){\vector(1,0){1.7}}
\put(3.1,1.75){$\epsilon$}
\put(2.2,1.4){\vector(0,-1){0.9}}
\put(1.45,0.9){\small $\exists!T_{\phi}$}
\put(4.2,1.5){\Large $\mathfrak{f}(\mathcal{S})$}
\put(5.225,1.625){\vector(1,0){1.5}}
\put(5.45,1.7){\tiny inclusion}
\put(4.6,1.3){\vector(0,-1){0.8}}
\put(4.7,0.85){$\hat{\psi}|_{\mathfrak{f}(\mathcal{S})}$}
\put(4.075,1.375){\vector(-3,-2){1.55}}
\multiput(4.06,1.38)(-0.15,-0.1){10}{\color{white}{\rule{3mm}{0.4mm}}}
\put(3.4,0.6){$\psi${\scriptsize ?}}
\put(6.8,1.5){\Large $\mbox{\large \Fontlukas L}\left(\mathscr{T}(V)\right)$}
\put(7.475,1.35){\vector(-3,-2){1.55}}
\put(6.9,0.525){$\exists!\hat{\psi}$}
\put(2.075,0.1){\Large $\mathfrak{g}$}
\put(2.45,0.125){\vector(1,0){1.325}}
\put(2.95,-0.25){$\myi$}
\put(3.9,0){\Large $\mbox{\large \Fontlukas L}\left(\mathscr{U}(\mathfrak{g})\right)$}
\end{picture}
\end{center}
Here, the set mapping $\phi: \mathcal{S} \longrightarrow \mathfrak{g}$ extends uniquely to a linear transformation $T_{\phi}: V \longrightarrow \mathfrak{g}$. 
Moreover, by the universal mapping property enjoyed by the tensor algebra, there is a unique $\mathbb{F}$-algebra homomorphism $\hat{\psi}: \mathscr{T}(V) \longrightarrow \mathscr{U}(\mathfrak{g})$ -- which is also a Lie algebra homomorphism $\myLieop\left(\mathscr{T}(V)\right) \longrightarrow \myLieop\left(\mathscr{U}(\mathfrak{g})\right)$ -- such that $\myi \circ T_{\phi} = \hat{\psi} \circ \epsilon$. 
The restriction $\hat{\psi}|_{\mathfrak{f}(\mathcal{S})}$ of $\hat{\psi}$ to $\mathfrak{f}(\mathcal{S})$ is the unique Lie algebra homomorphism $\mathfrak{f}(\mathcal{S}) \longrightarrow \mathfrak{g}$ such that $\myi \circ T_{\phi} = \hat{\psi}|_{\mathfrak{f}(\mathcal{S})} \circ \epsilon$. 
In particular, $\hat{\psi}|_{\mathfrak{f}(\mathcal{S})}(\epsilon(V)) \subseteq \myi(\mathfrak{g})$. 
So, the smallest Lie subalgebra $\mathfrak{s}$ of $\myLieop\left(\mathscr{U}(\mathfrak{g})\right)$ containing $\hat{\psi}|_{\mathfrak{f}(\mathcal{S})}(\epsilon(V))$ must reside within $\myi(\mathfrak{g})$. 
Since $\hat{\psi}^{-1}(\mathfrak{s})$ is a Lie subalgebra of $\myLieop\left(\mathscr{T}(V)\right)$ that contains $V$, then it must be the case that $\mathfrak{f}(\mathcal{S}) \subseteq \hat{\psi}^{-1}(\mathfrak{s})$. 
Therefore $\hat{\psi}|_{\mathfrak{f}(\mathcal{S})}(\mathfrak{f}(\mathcal{S})) = \hat{\psi}(\mathfrak{f}(\mathcal{S})) \subseteq \mathfrak{s} \subseteq \myi(\mathfrak{g})$. 
To complete the demonstration of the desired universal mapping property, we now use the fact that $\myi$ is injective and take $\psi := \myi^{-1} \circ \hat{\psi}|_{\mathfrak{f}(\mathcal{S})}$. 

The above well-defined notion of a free Lie algebra is needed to make precise the concept of constructing a Lie algebra by generators and relations. 
Let $\mathfrak{f}(\mathcal{S})$ be the free Lie algebra over a set $\mathcal{S}$. 
Let $\mathcal{R}$ be a subset of $\mathfrak{f}(\mathcal{S})$, with $\mbox{\Fontauri I}(\mathcal{R})$ as the smallest Lie ideal in $\mathfrak{f}(\mathcal{S})$ containing $\mathcal{R}$. 
The Lie algebra {\em generated by} $\mathcal{S}$ {\em subject to relations} $\mathcal{R}$, denoted by $\langle \mathcal{S}\, |\, \mathcal{R} \rangle$, is the quotient $\mathfrak{f}(\mathcal{S})/\mbox{\Fontauri I}(\mathcal{R})$. 
Let $\myp: \mathfrak{f}(\mathcal{S}) \longrightarrow \langle \mathcal{S}\, |\, \mathcal{R} \rangle$ denote the natural projection of $\mathfrak{f}(\mathcal{S})$ onto the quotient $\mathfrak{f}(\mathcal{S})/\mbox{\Fontauri I}(\mathcal{R}) = \langle \mathcal{S}\, |\, \mathcal{R} \rangle$. 
Then the pair $(\langle \mathcal{S}\, |\, \mathcal{R} \rangle,\myp)$ uniquely satisfies the universal property depicted below:
\begin{center}
\hspace*{-0.35in}
\setlength{\unitlength}{1cm}
\begin{picture}(4.5,2)
\thinlines
\put(0,1.5){\Large $\mathcal{S}$}
\put(0.375,1.625){\vector(1,0){1.4}}
\put(0.55,1.7){\tiny inclusion}
\put(0.375,1.4){\vector(3,-2){1.5}}
\put(0.6,0.625){$\forall\phi$}
\put(1.9,1.5){\Large $\mathfrak{f}(\mathcal{S})$} 
\put(2.9,1.625){\vector(1,0){1.4}}
\put(3.4,1.75){$\mysmallp$}
\put(2.2,1.3){\vector(0,-1){0.8}}
\put(1.65,0.95){\footnotesize $\exists!\hat{\phi}$}
\put(2.3,0.95){\tiny $\mathcal{R} \subseteq \mbox{\sffamily \tiny ker}(\hat{\phi})$}
\put(4.4,1.5){\Large $\langle \mathcal{S}\, |\, \mathcal{R} \rangle$}
\put(5.1,1.3){\vector(-3,-1){2.65}}
\put(3.85,0.6){$\exists!\psi$}
\put(2.075,0.1){\Large $\mathfrak{g}$}
\end{picture}
\end{center}
\noindent 
That is, for any set mapping $\phi: \mathcal{S} \longrightarrow \mathfrak{g}$ whose induced Lie algebra homomorphism $\hat{\phi}: \mathfrak{f}(\mathcal{S}) \longrightarrow \mathfrak{g}$ has $\mathcal{R} \subseteq \mbox{\sffamily ker}(\hat{\phi})$, there exists a unique Lie algebra homomorphism $\psi: \langle \mathcal{S}\, |\, \mathcal{R} \rangle \longrightarrow \mathfrak{g}$ such that $\hat{\phi} = \psi \circ \myp$.

For the remainder of this section, our ground field for Lie algebras is $\mathbb{C}$. 
The following definitions closely follow \cite{Kumar}, except that we work with the transpose of the generalized Cartan matrix. 
Given an integral embryophyte $\mathscr{G}$, let $\myr := n - \mbox{\sffamily rank}(M(\mathscr{G}))$ be the co-rank of the pulsation matrix $M(\mathscr{G})$. 
Let $\mathfrak{h}_{\mathscr{G}}$ be a complex vector space of dimension $n+\myr$.  
The natural pairing $\llangle \cdot,\cdot \rrangle: \mathfrak{h}_{\mathscr{G}}^{*} \times \mathfrak{h}_{\mathscr{G}} \longrightarrow \mathbb{C}$ is taken to be $\llangle \tbeta,v \rrangle := \tbeta(v)$ for all $\tbeta \in \mathfrak{h}^{*}_{\mathscr{G}}$ and $v \in \mathfrak{h}_{\mathscr{G}}$. 
Choose linearly independent sets ${\Delta}_{\mathscr{G}} := \{\tbeta_{i}\}_{i \in I} \subseteq \mathfrak{h}^{*}_{\mathscr{G}}$ and $\widetilde{\Delta}_{\mathscr{G}} := \{{\mysroot}_{j}\}_{j \in I} \subseteq \mathfrak{h}_{\mathscr{G}}$ and  such that $M_{ij} = \tbeta_{i}(\mysroot_{j}) =: \llangle \tbeta_{i},\mysroot_{j} \rrangle$. 
The triple $(\mathfrak{h}_{\mathscr{G}},\Delta_{\mathscr{G}},\widetilde{\Delta}_{\mathscr{G}})$ is called a {\em realization} of $M(\mathscr{G})$. 
Suppose we are given another triple $(\mathfrak{h}',\Delta',\widetilde{\Delta}')$ such that $\dim_{\mathbb{C}}(\mathfrak{h}')=n+\myr$, ${\Delta}' = \{\alpha_{i}'\}_{i \in I} \subseteq (\mathfrak{h}')^{*}$, $\widetilde{\Delta}' = \{\myroot_{j}'\}_{j \in I} \subseteq \mathfrak{h}'$, and $\alpha_{i}'(\myroot_{j}') = M_{ij}$. 
Then there exists an isomorphism of vector spaces $\theta: \mathfrak{h}' \longrightarrow \mathfrak{h}_{\mathscr{G}}$ such that $\theta(\myroot_{i}') = \mysroot_{i}$ and $\theta^{*}(\alpha_{j}') = \tbeta_{j}$, where $\theta^{*}: (\mathfrak{h}')^{*} \longrightarrow \mathfrak{h}^{*}_{\mathscr{G}}$ is induced by $\theta$. 
In this sense, the realization $(\mathfrak{h}_{\mathscr{G}},\Delta_{\mathscr{G}},\widetilde{\Delta}_{\mathscr{G}})$ is unique to $M(\mathscr{G})$. 

We are now ready to define the Lie algebra associated with our IEG $\mathscr{G}$. 
Let $\{\myqx_{i},\myqy_{i}\}_{i \in I(\mathscr{G})}$ be a set of symbols. 
Given our realization $(\mathfrak{h}_{\mathscr{G}},\Delta_{\mathscr{G}},\widetilde{\Delta}_{\mathscr{G}})$ of $M(\mathscr{G})$, let $\mathcal{S}_{\mathscr{G}}$ be the set $\mathfrak{h}_{\mathscr{G}}  \disjointunion \{\myqx_{i},\myqy_{i}\}_{i \in I(\mathscr{G})}$. 
Next we identify some families of relations $\mathcal{R}_{1}$ through $\mathcal{R}_{4}^{\pm}$:
\begin{description}
\item[($\mathcal{R}_{1}$)] $[\myqh, \myqk] = 0$ for all $\myqh, \myqk \in \mathfrak{h}_{\mathscr{G}}$ 
\item[($\mathcal{R}_{2}$)] $[\myqh, \myqx_{j}] = \tbeta_{j}(\myqh)\myqx_{j}$, $[\myqh, \myqy_{j}] = -\tbeta_{j}(\myqh)\myqy_{j}$ for all $\myqh \in \mathfrak{h}_{\mathscr{G}}$ and $j \in I$
\item[($\mathcal{R}_{3}$)] $[\myqx_{i}, \myqy_{j}] = \delta_{ij}\mysroot_{i}$ for all $i, j \in I$
\item[($\mathcal{R}_{4}^{\pm}$)] For all $i, j \in I$ with $i \not= j$ we have ($+$) $(\mbox{\sffamily ad}\myqx_{i})^{1-M_{ji}}(\myqx_{j}) = 0$ and ($-$) $(\mbox{\sffamily ad}\myqy_{i})^{1-M_{ji}}(\myqy_{j}) = 0$ 
\end{description}
Let $\mathcal{R}_{\mathscr{G}}$ be the collection of all relations in the above families $\mathcal{R}_{1}$ through $\mathcal{R}_{4}^{\pm}$. 
The {\em Kac--Moody Lie algebra} $\mathfrak{g} = \mathfrak{g}(\mathscr{G})$ is the Lie algebra $\langle \mathcal{S}_{\mathscr{G}}\, |\, \mathcal{R}_{\mathscr{G}} \rangle$, i.e.\ the Lie algebra generated by $\mathfrak{h}_{\mathscr{G}}  \disjointunion \{\myqx_{i},\myqy_{i}\}_{i \in I(\mathscr{G})}$ subject to the above relations.  

Let $\mathcal{Q} := \mathbb{Z}{\Delta}$ be the $\mathbb{Z}$-lattice of integer linear combinations of the elements of ${\Delta}$, and let $\mathcal{Q}^{+} := \{\sum_{i \in I}a_{i}\tbeta_{i} \in \mathcal{Q}\, |\, a_{i} \geq 0 \mbox{ for all } i \in I\}$. 
For $\tbeta \in \mathcal{Q}$, set $\mathfrak{g}_{\tbeta} := \{\myqz \in \mathfrak{g}\, |\, [\myqh,\myqz] = \tbeta(\myqh)\myqz\  \forall \myqh \in \mathfrak{h}\}$. 
The root-space decomposition of $\mathfrak{g}$ from Theorem 1.2.1 of \cite{Kumar} is the equality $\mathfrak{g} = \mathfrak{h} \bigoplus_{\tbeta \in \mathcal{Q}^{+} \setminus \{0\}}\left(\mathfrak{g}_{\tbeta} \oplus \mathfrak{g}_{-\tbeta}\right)$ where, for any $\tbeta \in \mathcal{Q}^{+} \setminus \{0\}$, we have $\dim(\mathfrak{g}_{\tbeta}) = \dim(\mathfrak{g}_{-\tbeta}) < \infty$. 
Let $\mbox{\Fontauri \small R}(\mathscr{G}) := \{\tbeta \in \mathcal{Q} \setminus \{0\}\, |\, \dim(\mathfrak{g}_{\tbeta}) \not= 0\}$ be the set of {\em roots}.  
If we let $\mbox{\Fontauri \small R}(\mathscr{G})^{+} := \mbox{\Fontauri \small R}(\mathscr{G}) \cap \mathcal{Q}^{+}$ and $\mbox{\Fontauri \small R}(\mathscr{G})^{-} := \{-\tbeta\, |\, \tbeta \in \mbox{\Fontauri \small R}(\mathscr{G})^{+}\}$, then $\mbox{\Fontauri \small R}(\mathscr{G}) = \mbox{\Fontauri \small R}(\mathscr{G})^{+} \disjointunion \mbox{\Fontauri \small R}(\mathscr{G})^{-}$. 

Following \S 1.3 of \cite{Kumar}, we define, for each $i \in I$, an involutive linear transformation $T^{*}_{i}$ on $\mathfrak{h}_{\mathscr{G}}^{*}$ by $T^{*}_{i}(\mu) := \mu - \llangle \mu,\mysroot_{i} \rrangle \tbeta_{i}$ for all $\mu \in \mathfrak{h}_{\mathscr{G}}^{*}$. 
Comparing \GeneratorsAndRelations\ above to Proposition 1.3.21 of \cite{Kumar}, we see that the subgroup of $\mbox{\sffamily Aut}(\mathfrak{h}_{\mathscr{G}}^{*})$ generated by $\{T^{*}_{i}\}_{i \in I}$ is isomorphic to our Weyl group $\mathcal{W}(\mathscr{G})$. 
This connection effects the application of \FiniteWeylGroupResult\ in the proof of the next result. 
From here on, we view $\mathcal{W}(\mathscr{G})$ as acting on $\mathfrak{h}_{\mathscr{G}}^{*}$ via $\mygens_{i}.\mu := T^{*}_{i}(\mu)$ for all $\mu \in \mathfrak{h}_{\mathscr{G}}^{*}$ and $i \in I$. 

\noindent 
{\bf \KacMoodyCorollary}\ \ {\sl Let $\mathscr{G}$ be a connected IEG. 
The Kac--Moody Lie algebra $\mathfrak{g}(\mathscr{G})$ is finite-dimensional if and only if $\mathscr{G}$ is a Coxeter--Dynkin flower, in which case $\mathfrak{g}(\mathscr{G})$ is simple.} 

{\em Proof.} Considering the root-space decomposition of $\mathfrak{g} = \mathfrak{g}(\mathscr{G})$, we see that $\mathfrak{g}$ is finite-dimensional if and only if the set $\mbox{\Fontauri \small R}(\mathscr{G})$ is finite. 
By Proposition 1.4.2.e of \cite{Kumar}, $\mbox{\Fontauri \small R}(\mathscr{G})$ is finite if and only if the Weyl group $\mathcal{W}(\mathscr{G})$ is finite.  
So, the claimed equivalence now follows from \FiniteWeylGroupResult. 
By Proposition 1.7 of \cite{Kac}, $\mathfrak{g}$ is simple if and only if $\det(M(\mathscr{G})) \not= 0$ and $\Gamma(\mathscr{G})$ is connected. 
By hypothesis, $\Gamma(\mathscr{G})$ is connected; that $\det(M(\mathscr{G})) \not= 0$ follows from \InnerProductTheorem.\hfill\QED

{\bf [\S \WeylSection.13:\! Linear representations of Kac--Moody Lie algebras.]} 
Above we defined the Kac--Moody Lie algebra $\mathfrak{g}(\mathscr{G})$ associated with any given integral embryophytic graph $\mathscr{G}$. 
Our aim here is to investigate the possible finite-dimensionality of $\mathfrak{g}(\mathscr{G})$-modules within a class of certain well-behaved representations of $\mathfrak{g}(\mathscr{G})$. 
After some set-up, the classification result we arrive at in \KMRepCorollary\ is effected by the first classification result of this section, namely \OrbitProposition.

Key to the construction of $\mathfrak{g} = \mathfrak{g}(\mathscr{G})$ is the so-called realization $(\mathfrak{h}_{\mathscr{G}},\Delta_{\mathscr{G}},\widetilde{\Delta}_{\mathscr{G}})$, where $\mathfrak{h}_{\mathscr{G}}$ is an $(n+\myr)$-dimensional complex vector space ($\myr$ is the co-rank of $M = M(\mathscr{G})$) and where the sets ${\Delta}_{\mathscr{G}} = \{\tbeta_{i} \in \mathfrak{h}_{\mathscr{G}}^{*}\}_{i \in I}$ and $\widetilde{\Delta}_{\mathscr{G}} = \{\mysroot_{j} \in \mathfrak{h}_{\mathscr{G}}\}_{j \in I}$ satisfy $\llangle \tbeta_{i},\mysroot_{j} \rrangle = M_{ij}$. 
According to the root-space decomposition of $\mathfrak{g}$ discussed above, we have $\mathfrak{g} = \mathfrak{h}_{\mathscr{G}} \bigoplus_{\tbeta \in \mbox{\Fontauri \scriptsize R}(\mathscr{G})^{+}}\left(\mathfrak{g}_{\tbeta} \oplus \mathfrak{g}_{-\tbeta}\right)$. 
We introduce a partial order on $\mathfrak{h}_{\mathscr{G}}^{*}$ where, for any $\mu, \lambda \in \mathfrak{h}_{\mathscr{G}}^{*}$, we have $\mu \leq \lambda$ if and only if $\lambda - \mu \in \mathcal{Q}^{+}$. 

We can simultaneously think of a representation of $\mathfrak{g}$ as a $\mathfrak{g}$-module. 
Say a $\mathfrak{g}$-module $V$ is {\em nontrivial} if there exists $\myqz \in \mathfrak{g}$ and $v \in V$ with $\myqz.v \not= 0$. 
For any $\mathfrak{g}$-module $V$ and $\mu \in \mathfrak{h}_{\mathscr{G}}^{*}$, let $V_{\mu} := \{v \in V\, |\, \myqh.v = \mu(\myqh)v,  \forall \myqh \in \mathfrak{h}_{\mathscr{G}}\}$; this subspace of $V$ is called a {\em weight space} and its nonzero elements are {\em weight vectors with weight} $\mu$. 
Call $V$ a {\em weight module} if $V = \bigoplus_{\mu \in \mathfrak{h}_{\mathscr{G}}^{*}}V_{\mu}$; in this case, we call $\mathscr{W}(V) := \{\mu \in \mathfrak{h}_{\mathscr{G}}^{*}\, |\, \dim V_{\mu} > 0\}$ the {\em set of weights of} $V$. 
We note that for $\mu \in \mathscr{W}(V)$ and $v \in V_{\mu}$, then $\myqx_{i}.v \in V_{\mu+\tbeta_{i}}$ and $\myqy_{i}.v \in V_{\mu-\tbeta_{i}}$. 
One can see this for $\myqx_{i}.v$ as follows: $\myqh.(\myqx_{i}.v) = \myqx_{i}.(\myqh.v) + \myqh.(\myqx_{i}.v) - \myqx_{i}.(\myqh.v) = \mu(\myqh)\myqx_{i}.v+ [\myqh, \myqx_{i}].v = \mu(\myqh)\myqx_{i}.v+ \tbeta_{i}(\myqh)\myqx_{i}.v = (\mu+\tbeta_{i})(\myqh)\myqx_{i}.v$. 
Similarly see that $\myqh.(\myqy_{i}.v) = (\mu-\tbeta_{i})(\myqh)\myqy_{i}.v$. 
Any basis for the weight module $V$ that respects the weight space decomposition is called a {\em weight basis}. 

Say $V$ is a {\em highest weight module with highest weight $\lambda$} if, for some $\lambda \in \mathfrak{h}_{\mathscr{G}}^{*}$ there exists a weight vector $v_{\lambda}$ with weight $\lambda$ such that ({\em a}) $\myqx.v_{\lambda} = 0$ for any $\myqx \in \mathfrak{g}_{\tbeta}$ for all $\tbeta \in \mbox{\Fontauri \small R}(\mathscr{G})^{+}$ and ({\em b}) $\mathfrak{g}.v_{\lambda} = V$. 
Say a weight module $V$ is an {\em integrable module} if, under the Lie algebra homomorphism $\phi: \mathfrak{g} \longrightarrow \myLieop(\mbox{\sffamily End}(V))$ associated with the action of $\mathfrak{g}$ on $V$, each of $\phi(\myqx_{i}), \phi(\myqy_{i}) \in \mbox{\sffamily End}(V)$ is locally nilpotent. 
Now, for each $\lambda \in \mathfrak{h}_{\mathscr{G}}^{*}$, there exists a highest weight module $\mathcal{V}(\lambda)$ with highest weight $\lambda$, called a {\em Verma module}, which enjoys the following properties: (0) $\dim \mathcal{V}(\lambda)_{\lambda} = 1$, (1) $\mathcal{V}(\lambda)$ has finite-dimensional weight spaces, (2) $\mu \leq \lambda$ in $\mathfrak{h}_{\mathscr{G}}^{*}$ whenever $\mu \in \mathscr{W}(\mathcal{V}(\lambda))$, and (3) any highest weight module $W$ with highest weight $\lambda$ is a quotient of $\mathcal{V}(\lambda)$. 
From here on, $(-\infty,\lambda] := \{\mu \in \mathfrak{h}_{\mathscr{G}}^{*}\, |\, \mu \leq \lambda\}$. 
The Verma module $\mathcal{V}(\lambda)$ has a unique maximal proper $\mathfrak{g}$-submodule we denote $\mathcal{V}'(\lambda)$. 
Let $V(\lambda)$ be the quotient $\mathcal{V}(\lambda)/\mathcal{V}'(\lambda)$, which is clearly irreducible. 

Verma modules and their quotients are among the most important members of the `category $\mathcal{O}$ modules'. 
A module $V$ is {\em category} $\mathcal{O}$ if it is a weight module with finite-dimensional weight spaces and if there exists a finite list $\lambda_{1},\ldots,\lambda_{m}$ of elements from $\mathfrak{h}_{\mathscr{G}}^{*}$ such that $\mathscr{W}(V) \subset \bigcup_{i=1}^{m}(-\infty,\lambda_{i}]$. 
The {\em formal character} $\mychar(V)$ of a $\mathfrak{g}$-module $V$ in category $\mathcal{O}$ is the formal Laurent power series $\mychar(V) := \sum_{\mu \in \mathscr{W}(V)}(\dim V_{\mu})\myvarZ^{\mu}$, cf.\ \S \WeylSection.9 above. 

The irreducible module $V(\lambda)$ is category $\mathcal{O}$ for all $\lambda \in \mathfrak{h}_{\mathscr{G}}^{*}$, and any irreducible module in category $\mathcal{O}$ must be isomorphic to some such $V(\lambda)$. 
Let $\mathcal{P} := \{\mu \in \mathfrak{h}_{\mathscr{G}}^{*}\, |\, \llangle \mu,\mysroot_{i} \rrangle \in \mathbb{Z}, \forall i \in I\}$, and let $\mathcal{P}^{+} := \{\mu \in \mathcal{P}\, |\, \llangle \mu,\mysroot_{i} \rrangle \geq 0, \forall i \in I\}$. 
The irreducible category-$\mathcal{O}$-module $V(\lambda)$ is integrable if and only if $\lambda \in \mathcal{P}^{+}$. 
A final fact we need is that for any nontrivial integrable module $V$ in category $\mathcal{O}$, there exists a nonzero $\lambda \in \mathcal{P}^{+}$ such that irreducible and integrable category $\mathcal{O}$ module $V(\lambda)$ is isomorphic to some submodule $W$ of $V$.  

\noindent 
{\bf \KMRepCorollary}\ \ {\sl Let $\mathscr{G}$ be a connected IEG. 
The Kac--Moody Lie algebra $\mathfrak{g}(\mathscr{G})$ has a nontrivial and finite-dimensional weight module if and only if $\mathscr{G}$ is a Coxeter--Dynkin flower.}

{\em Proof.} 
Suppose a $\mathfrak{g}$-module $V$ is a nontrivial and finite-dimensional weight module. 
Finite-dimensionality implies that $V$ is integrable and in category $\mathcal{O}$. 
Let $\lambda$ be a nonzero weight in $\mathcal{P}^{+}$ such that the irreducible and integrable category $\mathcal{O}$ module $V(\lambda)$ is isomorphic to some submodule $W$ of $V$. 
By Lemma 1.3.5.a of \cite{Kumar}, $\dim W_{\lambda} = \dim W_{\sigma.\lambda}$ for all $\sigma \in \mathcal{W}(\mathscr{G})$. 
Since $W$ is finite-dimensional, then $\lambda$ must have a finite orbit under the action of $\mathcal{W} = \mathcal{W}(\mathscr{G})$. 

Let $\lambda' \in \Lambda \subset \mathfrak{E}^{*}$ be the dominant weight $\sum_{i \in I}\llangle \lambda,\mysroot_{i} \rrangle\omega_{i}$. 
Then for any $j \in I$, we have $\llangle \lambda',\myroot_{j} \rrangle = \sum_{k \in I} \llangle \lambda,\mysroot_{k} \rrangle \llangle \omega_{k},\myroot_{j} \rrangle = \llangle \lambda,\mysroot_{j} \rrangle$. 
Let $q$ be a positive integer, and let $\mathcal{S} = \{s_{1},s_{2},\ldots,s_{r}\} \subseteq [1,q]_{\mathbb{Z}}$ be a nonempty subset whose elements are listed in increasing order: $s_{1} < s_{2} < \cdots < s_{r}$. 
Let $\myquant(\mathcal{S}) := \llangle \lambda,\mysroot_{i_{s_{1}}} \rrangle\prod_{t=1}^{r-1}M_{s_{t},s_{t+1}}$, and let $\myquant'(\mathcal{S}) := \llangle \lambda',\myroot_{i_{s_{1}}} \rrangle\prod_{t=1}^{r-1}M_{s_{t},s_{t+1}}$. 
It is easy to see by induction that for any positive integer $p$, we have 
\[\mygens_{i_{p}}\cdots\mygens_{i_{1}}.\lambda = \lambda - \sum_{q=1}^{p}\left(\llangle \lambda,\mysroot_{i_{q}} \rrangle + \sum_{r=1}^{q-1}(-1)^{r}\sum_{\mathcal{S} \subseteq [1,q-1]_{\mathbb{Z}}, |S|=r}\myquant(\mathcal{S})\right)\beta_{i_{q}}\]
and that 
\[\mygens_{i_{p}}\cdots\mygens_{i_{1}}.\lambda' = \lambda' - \sum_{q=1}^{p}\left(\llangle \lambda',\myroot_{i_{q}} \rrangle + \sum_{r=1}^{q-1}(-1)^{r}\sum_{\mathcal{S} \subseteq [1,q-1]_{\mathbb{Z}}, |S|=r}\myquant'(\mathcal{S})\right)\alpha_{i_{q}}.\]
Equality of $\llangle \lambda,\mysroot_{j} \rrangle$ and $\llangle \lambda',\myroot_{j} \rrangle$ for all $j \in I$ together with linear independence of the sets $\{\beta_{j}\}_{j \in I}$ in $\mathfrak{h}_{\mathscr{G}}^{*}$ and $\{\alpha_{j}\}_{j \in I}$ in $\mathfrak{E}^{*}$ gives us an evident one-to-one correspondence between elements in the orbits $\mathcal{W}\lambda \subset \mathcal{P}$ and $\mathcal{W}\lambda' \subset \Lambda$. 
So $\mathcal{W}\lambda'$ is finite. 
By \OrbitProposition, $\mathscr{G}$ must be a Coxeter--Dynkin flower. 

For the converse, if $\mathscr{G}$ is a Coxeter-Dynkin flower, then (for example) any irreducible representation whose highest weight is the quasi-minuscule dominant weight (cf.\ \S \QuasiExampleSection\ below, which does not depend on this result) is a nontrivial and finite-dimensional weight module.\hfill\QED

\begin{center}
\underline{\hspace*{4in}}
\end{center}

\vspace*{0.5cm} 
\noindent
{\bf \S \LFKSection. {\em La Florado Klasado} (The Flowering Classification).} 
The abundance of Dynkin-diagram classification results from preceding sections deserves special notice. 
In this section, we will define a number of properties related to connected integral embryophytic graphs (IEG's) that are clearly intended to evoke connections with some of our preceding results. 
All of these properties are finitistic in some sense, and their equivalence will be recorded here as just a few amongst many other such equivalences. 
Taken all together, these equivalences comprise what we believe is the quintessential classification result in all of mathematics: We call it {\em La Florado Klasado}, which is Esperantese for `The Flowering Classification'. 
These equivalences have many different flavors -- some are elementary, some are obscure, some are technical, some are fundamental.  
This rich profusion of equivalences has been observed by many (e.g.\ \cite{HHSV}, \cite{Arnold}, \cite{Trapa}), but the extent to which these equivalences continue to proliferate and, indeed, flourish, seems to be a phenomenon that ought to be much more widely known and even celebrated, both within and outside of the mathematics community. 

{\bf [\S \LFKSection.1:\! Some reflections on the nature of mathematical classifications.]} In the December 2014 issue of the {\em Notices of the AMS}, Princeton philosopher J.\ P.\ Burgess reviews I.\ Hacking's book {\em Why is There a Philosophy of Mathematics at All?} \cite{Burgess}. 
Burgess notes that, broadly speaking, there are two reasons why the philosophy of mathematics persists as an engaging area of interest. 
One of these reasons is its often unexpected applicability, an idea perhaps best expressed by E.\ Wigner as ``the unreasonable effectiveness of mathematics.'' 
The second reason is the often compelling nature of mathematical proofs, which can make mathematical results seem as if they were ever casting about our universe, languishing, and waiting upon us, to be apprehended. 

In the review, Burgess briefly considers Hacking's discussion of classification theorems, and specifically their relationship to the philosophy and popular understanding of mathematics. 
Burgess and Hacking agree that the classification of platonic solids is too prosaic to serve as an effective illustration of this kind of result. 
But Burgess disagrees with Hacking's subsequent use of ``The Enormous Theorem'' -- the classification of the finite simple groups -- as a replacement illustration, as the latter is simply too complicated to serve as an effective example of the phenomenon of mathematical taxonomies. 
Burgess then momentarily speculates upon the following: 
\begin{center}
\parbox{5.5in}{\footnotesize 
Surely there must be examples---two-manifolds, perhaps, or nonplanar graphs---less hackneyed than the Platonic solids but less ferociously technical than finite simple groups.}
\end{center}
Burgess' musing could be framed as a challenge to the mathematical community: What classification result might effectively highlight the beauty and utility of such taxonomies while being at least partly accessible and also deeply rooted within the main currents of contemporary mathematical practice? 
{\em La Florado Klasado} is, it seems, ideally suited to these purposes and is, in addition, beautifully interconnected with the aforementioned classifications of regular polytopes and also of the finite simple groups. 

{\em La Florado Klasado} ({\em LFK}) first entered human consciousness in the 1880's through the mind of the German mathematician W.\ Killing and was first recorded in his groundbreaking paper \cite{Killing}. 
Aspects of Killing's life, and of the brilliance of Killing's paper \cite{Killing}, are wonderfully recounted in A.\ J.\ Coleman's 1989 {\em Mathematical Intelligencer} article entitled ``The greatest mathematical paper of all time'' \cite{Coleman}. 
Some thirty years on now, Coleman's assessment of Killing's paper is arresting and yet reasonable. 
Now, many aspects of {\em LFK} that we have considered in preceding sections are somewhat elementary and combinatorial, but Killing's encounter with {\em LFK} was occasioned by his study of what is perhaps the most esoteric entry point imaginable: {\sl The classification of finite-dimensional simple Lie algebras over the complex numbers}. 
Killing's classification result was formally completed by E.\ Cartan in the 1890's \cite{Cartan}, and efforts to further simplify this classification continued well into the twentieth century. 
In the decades since, many other occurrences of this classification have been discovered, so that now the proliferation of these results, and their tendril-like reach into many areas of mathematics, is an efflorescence of exquisite and luminous beauty. 

{\bf [\S \LFKSection.2:\! Our statement of {\em LFK}, with proofs.]} For our statement of {\em La Florado Klasado} and our descriptions of various named properties related by equivalences of {\em LFK}, our focal object will be a connected integral embryophytic graph (IEG) $\mathscr{G}$. 
That said, the named properties below apply both to $\mathscr{G}$ and to $M(\mathscr{G})$, since any IEG is completely determined by its pulsation matrix.  
Say $\mathscr{G}$ (or $M(\mathscr{G})$) has property \fbox{\scriptsize F} if it is an integral Coxeter-Dynkin flower. 
(These are the `{\em floroj}', or flowers, of The Flowering Classification.) 
Say $\mathscr{G}$ is {\em ranked poset propitious}, or has property \fbox{\scriptsize P}, if there exists a finite and $\mathscr{G}$-structured poset with at least one edge. 
Say $\mathscr{G}$ is {\em orbitally obliging}, or has property \fbox{\scriptsize O}, if the orbit of some nonzero dominant weight under the action of the Weyl group $\mathcal{W}(\mathscr{G})$ is finite. 
Say $\mathscr{G}$ is {\em Coxeter group copacetic}, or has property \fbox{\scriptsize C}, if its Weyl group $\mathcal{W}(\mathscr{G})$ is finite. 
Say $\mathscr{G}$ is {\em root system refreshing}, or has property \fbox{\scriptsize R}, if its root system ${\Phi}(\mathscr{G})$ is finite. 
Say $\mathscr{G}$ is {\em index idiosyncratic}, or has property \fbox{\scriptsize I}, if some proper parabolic subgroup of the Weyl group $\mathcal{W}(\mathscr{G})$ has finite index. 
Say $\mathscr{G}$ is {\em symmetric function salutary}, or has property \fbox{\scriptsize S}, if there exists a non-constant $\mathcal{W}(\mathscr{G})$-symmetric function. 
Say $\mathscr{G}$ is {\em usefully Euclidean}, or has property \fbox{\scriptsize U}, if there exists a positive definite diagonal matrix $D$ such that $DM$ is symmetric and positive definite, where $M=M(\mathscr{G})$ is our pulsation matrix.
Say $\mathscr{G}$ is {\em game gratifying}, or has property \fbox{\scriptsize G}, if there exists a terminating Networked-numbers Game played from a nonzero dominant starting position. 
Say $\mathscr{G}$ is {\em cordially Kac--Moody}, or has property \fbox{\scriptsize K}, if the Kac--Moody Lie algebra $\mathfrak{g}(\mathscr{G})$ is finite-dimensional. 
Say $\mathscr{G}$ is {\em linear representation laudable}, or has property \fbox{\scriptsize L}, if the Kac--Moody Lie algebra $\mathfrak{g}(\mathscr{G})$ has a nontrivial and finite-dimensional weight module. 

What follows are just some of the equivalences comprising {\em LFK}. 

\noindent 
{\bf \LFK\ {\em La Florado Klasado}, or The Flowering Classification (originated by W.\ Killing and E.\ Cartan)}\ \ {\sl Let $\mathscr{G}$ be a connected integral embryophyte. 
Then the following are equivalent:\\ 
\hspace*{0.75in} (F) $\mathscr{G}$ is an integral Coxeter--Dynkin flower, i.e.\ has property \fbox{\scriptsize F}.\\ 
\hspace*{0.75in} (P) $\mathscr{G}$ is ranked poset propitious, i.e.\ has property \fbox{\scriptsize P}.\\ 
\hspace*{0.75in} (O) $\mathscr{G}$ is orbitally obliging, i.e.\ has property \fbox{\scriptsize O}.\\ 
\hspace*{0.75in} (C) $\mathscr{G}$ is Coxeter group copacetic, i.e.\ has property \fbox{\scriptsize C}.\\ 
\hspace*{0.75in} (R) $\mathscr{G}$ is root system refreshing, i.e.\ has property \fbox{\scriptsize R}.\\ 
\hspace*{0.75in} (I) $\mathscr{G}$ is index idiosyncratic, i.e.\ has property \fbox{\scriptsize I}.\\ 
\hspace*{0.75in} (S) $\mathscr{G}$ is symmetric function salutary, i.e.\ has property \fbox{\scriptsize S}.\\ 
\hspace*{0.75in} (U) $\mathscr{G}$ is usefully Euclidean, i.e.\ has property \fbox{\scriptsize U}.\\ 
\hspace*{0.75in} (G) $\mathscr{G}$ is game gratifying, i.e.\ has property \fbox{\scriptsize G}.\\ 
\hspace*{0.75in} (K) $\mathscr{G}$ is cordially Kac--Moody, i.e.\ has property \fbox{\scriptsize K}.\\
\hspace*{0.75in} (L) $\mathscr{G}$ is linear representation laudable, i.e.\ has property \fbox{\scriptsize L}.} 

Some comments on references for these equivalences are contained within the exposition of the previous section. 
Our proof in \S \WeylSection\ of the previously-known equivalence of \fbox{\scriptsize F} and \fbox{\scriptsize I} appears to be new. 
The equivalence of \fbox{\scriptsize F} and \fbox{\scriptsize O} also seems to be new, as is the equivalence of \fbox{\scriptsize F} and \fbox{\scriptsize S}.

{\em Proof of LFK.} The equivalence of \fbox{\scriptsize F} and \fbox{\scriptsize P} was offered in \S \GStructureIntroSection\ as \GCMTheorem, borrowing from one of this author's contributions to \cite{DE}. 
The remaining equivalences of the above version of {\em LFK} were established in the preceding section using \OrbitProposition. 
Specifically, the equivalence \fbox{\scriptsize F} and \fbox{\scriptsize O} was established in \OrbitProposition\ itself; 
the equivalence of \fbox{\scriptsize F} with \fbox{\scriptsize C} and with \fbox{\scriptsize R} in \FiniteWeylGroupResult; 
of \fbox{\scriptsize F} and \fbox{\scriptsize I} in \FiniteIndexCorollary; 
of \fbox{\scriptsize F} and \fbox{\scriptsize S} in \NonconstantTheorem; 
of \fbox{\scriptsize F} and \fbox{\scriptsize U} in \InnerProductTheorem; 
of \fbox{\scriptsize F} and \fbox{\scriptsize G} in \NGCorollary; 
of \fbox{\scriptsize F} and \fbox{\scriptsize K} in \KacMoodyCorollary;
and of \fbox{\scriptsize F} and \fbox{\scriptsize L} in \KMRepCorollary.\hfill\QED 

{\bf [\S \LFKSection.3:\! Some further remarks on our version of {\em LFK}.]} 
Each part of the above version of {\em LFK} generalizes naturally to IEG's that are not necessarily connected. 
Of course, our term `Coxeter--Dynkin posy' refers to a disjoint union of Coxeter--Dynkin flowers.  
\GCMTheorem\ of \S \GStructureIntroSection\ explicitly addresses ranked poset sturdiness in the case of Coxeter--Dynkin posies. 
When $\mathscr{G} = \mathscr{G}_{1} \oplus \cdots \oplus \mathscr{G}_{m}$ is a disjoint union of Coxeter--Dynkin flowers, we can view the pulsation matrix $M(\mathscr{G})$ as block diagonal with blocks $M(\mathscr{G}_{i})$ for $i \in [1,m]_{\mathbb{Z}}$ in which case the Euclidean space $\mathfrak{E}^{*}$ naturally decomposes as a direct sum $\mathfrak{E}_{1}^{*} \oplus \cdots \oplus \mathfrak{E}_{m}^{*}$ of pairwise orthogonal subspaces. 
Then $\mathcal{W}(\mathscr{G})$ is naturally isomorphic to $\mathcal{W}(\mathscr{G}_{1}) \times \cdots \times \mathcal{W}(\mathscr{G}_{m})$ and acts stably on the corresponding inner product subspaces. 
Further, the associated Lie algebra $\mathfrak{g}(\mathscr{G})$ decomposes as the direct sum $\mathfrak{g}(\mathscr{G})_{1} \oplus \cdots \oplus \mathfrak{g}(\mathscr{G})_{m}$ of simple Lie algebras. 
See \S 2 of \cite{DonPosetModels} for further details on how the associated root system and weight lattice decompose when the Coxeter--Dynkin posy is not connected. 

A version of La Florado Klasado can be given for connected \underline{unital} embryophytic graphs, i.e.\ connected UEG's, wherein all pulsation factors are unity. 
In fact, the connected UEG versions of these properties and/or equivalences are sometimes easier to state and/or motivate than in the case of IEG's. 
Moreover, there is a version of La Florado Klasado for connected `dihedral' embryophytic graphs, i.e.\ connected DEG's. 
Pulsation matrices for DEG's may have real number entries as long as any product of pulsation factors $M_{ij}M_{ji}$ is a {\em dihedral number}, meaning that if the said product is smaller than $4$, it must satisfy $M_{ij}M_{ji} = 4\cos^{2}(\pi/m_{ij})$ for some integer $m_{ij} > 1$. 
This results in some additional Coxeter--Dynkin flowers besides those encountered in \IEGGraphFigure\ above.  
However not all IEG properties extend to DEG's (for example, property \fbox{\scriptsize P}). 
See \cite{DonEur} for more discussion. 


\begin{center}
\underline{\hspace*{4in}}
\end{center}

\vspace*{0.5cm} 
\noindent 
{\bf \S \FlowerSection.\ Consequences of LFK for our algebraic-combinatorial environment.} 
Throughout this section, we take $\mathscr{G} = (\Gamma_{I},M_{I \times I})$ to be a Coxeter--Dynkin flower or posy, with connectedness explicitly invoked when needed, and we let $n := |I|$. 
By {\em La Florado Klasado} (\LFK), this hypothesis has many consequences for structures related to $\mathscr{G}$, many of which appeared in \S \WeylSection. 
In particular, the Weyl group $\mathcal{W}(\mathscr{G})$ and root system ${\Phi}(\mathscr{G})$ are finite, the $n$-dimensional real vector space $\mathfrak{E}^{*}$ has a natural and $\mathcal{W}$-invariant inner product, and there exist non-constant $\mathcal{W}$-symmetric functions. 
The more specific aspects we develop in this section can be lightly browsed on a first reading and then consulted later as needed. 
What we present here mostly borrows from standard sources (\cite{Hum}, \cite{HumCoxeter}, \cite{Kac}) as well as the monograph  \cite{DonPosetModels}. 
The latter provides a treatment of many of these ideas in language and notation mostly consistent with this paper.  

{\bf [\S \FlowerSection.1:\! A $\mathcal{W}$-invariant inner product space.]} 
Our $n$-dimensional real vector space $\mathfrak{E}^{*}$ with fundamental basis $\mathscr{B} = \{\omega_{i}\}_{i \in I}$ and simple root basis $\{\alpha_{i}\}_{i \in I}$ is now endowed with an inner product $\langle \cdot,\cdot \rangle$, as in \InnerProductTheorem. 
Regard any $\mu = \sum_{i \in I}\mu_{i}\omega_{i}$ in $\mathfrak{E}^{*}$ as a row vector $[\mu]_{\mathscr{B}} = [\mu_{i}]_{i \in I}$. 
Then $\langle \mu,\nu \rangle$ can be computed as $[\mu]M^{-1}D^{-1}[\nu]^{\myT}$, where $D_{ii} = 2/\langle \alpha_{i},\alpha_{i} \rangle$ for each $i \in I$. 
We have $\alpha_{i} = \sum_{j \in I}M_{ij}\omega_{j}$ and $\omega_{j} = \sum_{i \in I}(M^{-1})_{ij}\alpha_{i}$. 
For any nonzero $\mu \in \mathfrak{E}^{*}$, we set $\mu^{\vee} := \frac{2}{\langle \mu,\mu \rangle}v$. 
Then $\langle \omega_{i},\alpha_{j}^{\vee} \rangle = \delta_{ij}$, where $\delta_{ij}$ is Kronecker's delta, so $M = (M_{ij})_{i,j \in I} = \left(\langle \alpha_{i},\alpha_{j}^{\vee} \rangle\right)_{i,j \in I}$.  
 We have $\mathcal{W} = \left\langle \rule[-0.75mm]{0mm}{4mm} \mygens_{i} \right\rangle_{i \in I}$ acting on $\mathfrak{E}^{*}$ by the rule $\mygens_{i}.\mu = \mu - \langle \mu,\alpha_{i}^{\vee} \rangle\alpha_{i}$. 
The inner product is preserved by $\mathcal{W}$ in the sense that $\langle \sigma.\mu,\nu \rangle = \langle \mu,\sigma^{-1}.\nu \rangle$ for all $\sigma \in \mathcal{W}$ and $\mu,\nu \in \mathfrak{E}^{*}$. 

{\bf [\S \FlowerSection.2:\! The root system and co-root system.]} 
In \S \WeylSection, the root system ${\Phi} = \Phi(\mathscr{G})$ was realized within $\mathfrak{E}^{*}$ as $\{\sigma.\alpha_{i}\}_{\sigma \in \mathcal{W}, i \in I}$. 
We define the co-root system within $\mathfrak{E}^{*}$ as $\Phi^{\vee} = \Phi^{\vee}(\mathscr{G}) := \{\sigma.\alpha^{\vee}_{i}\}_{\sigma \in \mathcal{W}, i \in I}$. 
Observe that $\alpha \longmapsto \frac{2}{\langle \alpha,\alpha \rangle}\alpha$ gives us a one-to-one correspondence between $\Phi$ and $\Phi^{\vee}$
Some routine facts about root/co-root systems are that $\Phi$ is partitioned into a set $\Phi^{+}$ of {\em positive roots} and a set $\Phi^{-}$ of {\em negative roots} and that for any $\alpha \in \Phi$ we have $k\alpha \in \Phi$ if and only if $k \in \{\pm 1\}$. 
Recall from \S \WeylSection\ that: 
For any $\sigma \in \mathcal{W}$ and $i \in I$, we have $\ell(\sigma s_{i}) > \ell(\sigma)$ if and only if $\sigma.\alpha_{i} \in \Phi^{+}$, and $\ell(\sigma s_{i}) < \ell(\sigma)$ if and only if $\sigma.\alpha_{i} \in \Phi^{-}$. 
Moreover, $\ell(\sigma) = \rule[-1.5mm]{0.2mm}{5mm}\{\alpha \in \Phi^{+}\, |\, \sigma.\alpha \in \Phi^{-}\}\rule[-1.5mm]{0.2mm}{5mm}$. 
These statements still hold when we replace roots with co-roots and $\Phi$ with $\Phi^{\vee}$. 
The following elements from $\mathfrak{E}^{*}$ play distinguished roles in the theory to follow and so are given special names: 
$\varrho := \sum_{i \in I}\omega_{i} = \frac{1}{2}\sum_{\alpha \in \Phi^{+}}\alpha$ and $\varrho^{\vee} := \sum_{i \in I}\frac{2}{\langle \alpha_{i},\alpha_{i} \rangle}\omega_{i} = \frac{1}{2}\sum_{\alpha \in \Phi^{+}}\alpha^{\vee}$.

{\bf [\S \FlowerSection.3:\! The longest Weyl group element.]} 
Since our Weyl group is finite, there is some least upper bound on the lengths of its members.  
In fact, such a longest element is unique, and we denote it by $w_{0}$.  
The following lemma recapitulates this and other useful and well-known facts about $w_{0}$.  

\noindent
{\bf \LongestLemma}\ \ {\sl Let $\mathscr{G}$ be a Coxeter--Dynkin posy. 
(1) There exists a unique longest element $w_{0}$ of $\mathcal{W}$, and $w_{0}^{2} = \varepsilon$; (2) $w_{0}(\Phi^{+}) = \Phi^{-}$ and therefore $\ell(w_{0}) = |\Phi^{+}|$; (3) the length of any Networked-numbers Game played on $\mathscr{G}$ from initial position $\varrho$, or from any strongly dominant weight, is $l = \ell(w_{0})$, and $s_{i_{l}}\cdots s_{i_{1}}$ is a reduced expression for $w_{0}$ if and only if $(\gamma_{i_{1}},\ldots,\gamma_{i_{l}})$ is a game sequence from said initial position; (4) $w_{0}.\varrho = -\varrho$ and that $w_{0}.\varrho^{\vee} = -\varrho^{\vee}$; and (5) there exists an involution $\sigma_{0}:I \longrightarrow I$ of our color palette such that $w_{0}.\alpha_{i} = -\alpha_{\sigma_{0}(i)}$ and $w_{0}.\omega_{i} = -\omega_{\sigma_{0}(i)}$ for all $i \in I$.} 

{\em Proof.} 
For any longest element $w$ and for each $i \in I$, we have $\ell(w\mygens_{i}) < \ell(w)$ and therefore $w.\alpha_{i} \in \Phi^{-}$. 
It follows that $w.\alpha \in \Phi^{-}$ for all $\alpha \in \Phi^{+}$, so $\ell(w) = |\Phi^{+}|$, proving {\sl (2)}. 
If $v$ is another longest element, then for all $i \in I$ we have $(vw).\alpha_{i} = v.(w.\alpha_{i}) \in \Phi^{+}$ and therefore $\ell(vw\mygens_{i}) > \ell(vw)$; this can only be true if $vw = \varepsilon$.  
Since the latter is true for any pair of longest elements, there can be only one such element, demonstrating the uniqueness claim of part {\sl (1)} and the fact that $w_{0}$ is its own inverse. 
For {\sl (3)}, that every game sequence from a strongly dominant initial position has length $l = \ell(w)$ follows from Corollary 3.3 of \cite{DonEur}. 
The fact that $(\gamma_{i_{1}},\ldots,\gamma_{i_{l}})$ is a game sequence from any strongly dominant initial position if and only if $s_{i_{l}}\cdots s_{i_{1}}$ is a reduced expression follows from Eriksson's Reduced Word Result, stated as Theorem 2.8 in \cite{DonEur}. 
Part {\sl (4)} follows from the facts that $\varrho$ (resp.\ $\varrho^{\vee}$) is half the sum of the positive roots (resp.\ positive coroots) and that $-w_{0}$ permutes the positive roots (resp.\ positive coroots).
For {\sl (5)}, observe that $\langle \varrho,-w_{0}.\alpha_{i}^{\vee} \rangle = \langle -w_{0}.\varrho,\alpha_{i}^{\vee} \rangle = \langle \varrho,\alpha_{i}^{\vee} \rangle = 1$, which means that $-w_{0}.\alpha_{i}^{\vee}$ to be a simple coroot. 
Therefore, $-w_{0}$ permutes the simple coroots, i.e.\ there is a permutation $\sigma_{0}: I \longrightarrow I$ such that $-w_{0}.\alpha_{i}^{\vee} = \alpha_{\sigma_{0}(i)}^{\vee}$. 
Since $w_{0}^{2} = \varepsilon$, then $\sigma_{0}$ is an involution. 
Also, $\delta_{ij} = \langle \omega_{i},\alpha_{j}^{\vee} \rangle = \langle -w_{0}.\omega_{i},-w_{0}.\alpha_{j}^{\vee} \rangle = \langle -w_{0}.\omega_{i},\alpha_{\sigma_{0}(j)}^{\vee} \rangle$ for all $i,j \in I$, which forces $-w_{0}.\omega_{i} = \omega_{\sigma_{0}(i)}$. 
Finally, $M_{ij} = \langle \alpha_{i},\alpha_{j}^{\vee} \rangle = \langle -w_{0}.\alpha_{i},-w_{0}.\alpha_{j}^{\vee} \rangle = \langle -w_{0}.\alpha_{i},\alpha_{\sigma_{0}(j)}^{\vee} \rangle$ for all $i, j \in I$. 
Since $M$ is invertible in our finitistic context, we must therefore have $-w_{0}.\alpha_{i} = \alpha_{\sigma_{0}(i)}$.\hfill\QED

{\bf [\S \FlowerSection.4:\! Root lattice cosets as connected components of the weight lattice.]} 
The root lattice $\{\sum_{i \in I}k_{i}\alpha_{i}\, |\, k_{i} \in \mathbb{Z} \mbox{ for each } i \in I\} = \mathbb{Z}\{\alpha_{i}| i \in I\}$ and the weight lattice $\Lambda = \{\sum_{i \in I}k_{i}\omega_{i}\, |\, k_{i} \in \mathbb{Z} \mbox{ for each } i \in I\} = \mathbb{Z}\{\omega_{i}| i \in I\}$ interact in interesting and useful ways. 
Of course $\Phi \subset \mathbb{Z}\{\alpha_{i}| i \in I\} \subset \Lambda$. 
Partially order $\Lambda$ by the rule $\nu \leq \mu$ if and only if $\mu - \nu = \sum_{i \in I}k_{i}\alpha_{i}$ for some set of nonnegative integers $\{k_{i}\}_{i \in I}$. 
It follows from these definitions that the connected components of $\Lambda$ are just the cosets of $\mathbb{Z}\{\alpha_{i}| i \in I\}$. 
In Proposition 4.11 of \cite{DonPosetModels} we recorded the following well-known facts about these connected components. 

\noindent 
{\bf \MinimalWeightsLemma}\ \ {\sl Let $\mathscr{G}$ be a Coxeter--Dynkin flower. 
The connected components of $\Lambda$ are finite in number; 
amongst the dominant weights within some such connected component, there is a unique minimal dominant weight; and 
amongst the nonzero dominant weights within the connected component of the zero weight, there is a unique minimal nonzero dominant weight.}   

{\em Proof.} As noted above, this follows from Proposition 4.11 of \cite{DonPosetModels}.\hfill\QED

Amongst the `minimal' dominant weights of the preceding Lemma, the nonzero weights are called {\em minuscule}. 
The unique minimal nonzero dominant weight within the connected component of the zero weight is called {\em quasi-minuscule}. 
So, associated to any Coxeter--Dynkin flower is a unique quasi-minuscule dominant weight. 
However, the Coxeter--Dynkin flowers $\myE_{8}$, $\myF_{4}$, and $\myG_{2}$ have no associated minuscule dominant weights, while $\myA_{n}$ ($n \geq 2$), $\myD_{n}$, and $\myE_{6}$ each have more than one associated minuscule dominant weight. 

{\bf [\S \FlowerSection.5:\! $\mathscr{G}$-structured posets and saturated sets of weights.]} 
Let $R$ be a $\mathscr{G}$-structured poset. 
In \S \GStructureIntroSection, we defined the weight of an element $\xelt$ from $R$ to be the vector $wt(\xelt) = \left(\mym_{i}(\xelt)\right)_{i \in I}$, where $\mym_{i}(\xelt) = \rho_{i}(\xelt)-\delta_{i}(\xelt)$ is the rank of $\xelt$ within its $i$-component less the depth of $\xelt$ within this same component. 
We now express weights in terms of our fundamental basis: $wt(\xelt) = \sum_{i \in I}\mym_{i}(\xelt)\omega_{i}$. 
Define the {\em weight-generating function} of $R$ to be the Laurent polynomial $\WGF(R;\myvarZ) := \sum_{\xelt \in R}\myvarZ^{wt(\xelt)}$. 
Moreover, we declare $\Pi(R)$ to be the set of weights $\{wt(\xelt)\, |\, \xelt \in R\} \subset \Lambda$, and for any $\mu \in \Lambda$, let $R_{\mu} := \{\xelt \in R\, |\, wt(\xelt) = \mu\}$. 
So, $\WGF(R;\myvarZ) := \sum_{\mu \in \Lambda}\myabs R_{\mu} \myabs\myvarZ^{\mu} = \sum_{\mu \in \Pi(R)}\myabs R_{\mu} \myabs\myvarZ^{\mu}$. 

Other subsets of $\Lambda$ of interest are any so-called {\em saturated} set of weights $\mathscr{W}$ having the property that for any $\mu \in \mathscr{W}$, any root $\alpha$, and any integer $i$ between $0$ and $\langle \mu,\alpha^{\vee} \rangle$, it is the case that $\mu+i\alpha \in \mathscr{W}$. 
We give any saturated set the induced ordering from $\Lambda$. 
Corresponding to each dominant weight $\lambda$, there is a distinguished saturated set of weights:  
Let $\Pi(\lambda)$ be the set consisting of the weights in all orbits $\mathcal{W}\nu$ for all $\nu \in \Lambda^{+}$ such that $\nu \leq \lambda$, and give $\Pi(\lambda)$ this induced ordering from $\Lambda$. 
Moreover, if $R$ is a $\mathscr{G}$-structured poset, give $\Pi(R)$ the induced ordering from $\Lambda$. 
For any such $R$, say an element $\pelt \in R$ is {\em prominent} if $\delta_{i}(\pelt) = 0$ for all $i \in I$. 
What follows are some key facts about saturated sets of weights, which include $\Pi(\lambda)$ and $\Pi(R)$.  
These are mainly borrowed from Proposition 2.1/Proposition 3.6 and Theorem 3.8/Corollaries 3.9 and 3.11 of \cite{DonPosetModels}. 

\noindent 
{\bf \SaturatedLemma}\ \ {\sl (1) A finite saturated set of weights $\mathscr{W}$ is ranked, and moreover $\nu \rightarrow \mu$ in $\mathscr{W}$ if and only if $\mu - \nu = \alpha_{i}$ for some $i \in I$, in which case we assign color $i$ to this edge. 
(2) Now regarding our finite saturated set $\mathscr{W}$ to be a ranked poset with edges colored by $I$, then all monochromatic components of $\mathscr{W}$ are chains, $\mathscr{W}$ is diamond-colored, and, moreover, for any $\mu \in \mathscr{W}$, the quantity $wt(\mu) := \sum_{i \in I}(\mym_{i}(\mu))\omega_{i}$ is exactly $\mu$; it follows that $\mathscr{W}$ is $\mathscr{G}$-structured. 
(3) For any dominant $\lambda$, $\Pi(\lambda)$ is finite and topographically balanced with unique maximal element $\lambda$ and unique minimal element $w_{0}.\lambda$ and is the unique saturated set of weights in $\Lambda$ that has $\lambda$ as its only maximal element; moreover, $\Pi(\lambda)^{*} \cong \Pi(-w_{0}.\lambda)$ and $\Pi(\lambda)^{\bowtie} \cong \Pi(\lambda)$ are isomorphisms of edge-colored posets. 
(4) Let $\mathscr{M}$ be the set of maximal elements in some finite saturated set of weights $\mathscr{W}$. 
Observe that any $\lambda \in \mathscr{M}$ is dominant. 
Then $\displaystyle \mathscr{W} = \bigcup_{\lambda \in \mathscr{M}}\Pi(\lambda)$, an equality of edge-colored directed graphs. 
(5) If $R$ is $\mathscr{G}$-structured,  let $\mathscr{D} := \{\xelt \in R\, |\, \delta_{i}(\xelt)=0 \mbox{ for all } i \in I\}$ be the set of prominent elements in $R$. 
Then $\displaystyle \Pi(R) = \bigcup_{\delt \in \mathscr{D}} \Pi(wt(\delt))$, an equality of edge-colored directed graphs. 
(6) If $R$ is $\mathscr{G}$-structured and connected, we let $\mathscr{L} := \{\lambda \in \Pi(R)\, |\, \langle \lambda,\varrho^{\vee} \rangle \geq \langle wt(\xelt),\varrho^{\vee} \rangle \mbox{ for all } \xelt \in R\}$. 
Then $\belt \in R$ has rank zero if and only if $w_{0}.wt(\belt) \in \mathscr{L}$, and, moreover, $\telt \in R$ has maximal rank if and only if $wt(\telt) \in \mathscr{L}$. 
Now let $\lambda$ be any element of $\mathscr{L}$. 
Then} $\posetlength(R) = 2\langle \lambda,\varrho^{\vee} \rangle =  \langle \lambda-w_{0}.\lambda,\varrho^{\vee} \rangle = -2\langle w_{0}.\lambda,\varrho^{\vee} \rangle$. 
{\sl Also, for any $\xelt \in R$, we have $\rho(\xelt) = \langle wt(\xelt)+\lambda,\varrho^{\vee} \rangle = \langle wt(\xelt)-w_{0}.\lambda,\varrho^{\vee} \rangle$ and $\delta(\xelt) = \langle \lambda-wt(\xelt),\varrho^{\vee} \rangle = \langle -w_{0}.\lambda-wt(\xelt),\varrho^{\vee} \rangle$.} 

{\em Proof.} Parts {\sl (1)} and {\sl (2)} are from Proposition 3.6 of \cite{DonPosetModels}. 
All claims of part {\sl (3)} are from Proposition 2.1 of \cite{DonPosetModels} except the claims that $\Pi(\lambda)$ is topographically balanced and that $\Pi(\lambda)^{*} \cong \Pi(-w_{0}.\lambda)$ and $\Pi(\lambda)^{\bowtie} \cong \Pi(\lambda)$. 
For $\Pi(\lambda)^{*} \cong \Pi(-w_{0}.\lambda)$, observe that the bijection $\phi: \Pi(\lambda)^{*} \longrightarrow \Pi(-w_{0}.\lambda)$ given by $\phi(\mu^{*}) = -\mu$ preserves edges and edge colors: $\mu^{*} \myarrow{i} \nu^{*}$ if and only if $\nu \myarrow{i} \mu$ if and only if $-\mu \myarrow{i} -\nu$ if and only if $\phi(\mu^{*}) \myarrow{i} \phi(\nu^{*})$. 
For $\Pi(\lambda)^{\bowtie} \cong \Pi(\lambda)$, observe that the bijection $\psi: \Pi(\lambda)^{*} \longrightarrow \Pi(\lambda)^{\sigma_{0}}$ given by $\psi(\mu^{*}) = -w_{0}.\mu$ preserves edges and edge colors: $\mu^{*} \myarrow{i} \nu^{*}$ if and only if $\nu \myarrow{i} \mu$ if and only if $-w_{0}.\mu \mylongarrow{\sigma_{0}(i)} -w_{0}.\nu$ if and only if $\psi(\mu^{*}) \mylongarrow{\sigma_{0}(i)} \phi(\nu^{*})$. 
Part {\sl (4)} follows from Corollary 3.9 of \cite{DonPosetModels}, part {\sl (5)} follows from Theorem 3.8 of \cite{DonPosetModels}, and part {\sl (6)} follows from Corollary 3.11 of \cite{DonPosetModels}.\hfill\QED

{\em Proof of \GCMTheorem.3.} 
From parts {\sl (1)}--{\sl (3)} of the above lemma, it follows that for any dominant weight $\lambda$, the saturated set of weights $\Pi(\lambda)$ is $\mathscr{G}$-structured and its unique maximal element has weight $\lambda$.\hfill\QED

{\bf [\S \FlowerSection.6:\! Monomial symmetric functions and Weyl bialternants.]} 
Let us now take a closer look at $\mathcal{W}$-symmetric functions in the current context. 
(See \S \WeylSection.9 for the definition of a $\mathcal{W}$-symmetric function and some other related initial concepts.)  
In this subsection, $\lambda$ represents a generic dominant weight. 
Since $\mathcal{W}$ is finite, the orbit $\mathcal{W}\lambda$ and the stabilizer $\mbox{\sffamily stabilizer($\lambda$)}$ are both finite.  
We see, then, that the quantity $\zeta_{\lambda} := \frac{1}{|\mbox{\sffamily \tiny stabilizer($\lambda$)}|}\sum_{\sigma \in \mathcal{W}}\myvarZ^{\sigma.\lambda} = \sum_{\mu \in \mathcal{W}\lambda}\myvarZ^{\mu}$ is a $\mathcal{W}$-symmetric function.  
We call $\zeta_{\lambda}$ a {\em monomial symmetric function}. 

Among the most important $\mathcal{W}$-symmetric functions are the {\em Weyl bialternants}, which, as their name indicates, are a combination of two so-called `alternants'.  
Say a weight $\mu$ is {\em strongly dominant} and write $\mu \in \Lambda^{++}$ if $\mu_{i} >0$ for all $i \in I$, where $\mu = \sum_{i \in I}\mu_{i}\omega_{i}$. 
Let $\varrho$ denote the strongly dominant weight $\sum_{i \in I}\omega_{i}$, so $\mu$ is strongly dominant if and only if $\mu = \varrho + \mu'$ for some dominant weight $\mu'$. 
For any strongly dominant $\mu$, we define the {\em alternant} $\mathcal{A}(\myvarZ^{\mu})$ to be the alternating sum $\sum_{\sigma \in \mathcal{W}}\mysign(\sigma)\myvarZ^{\sigma.\mu}$, where $\mysign(\sigma) = (-1)^{\ell(\sigma)} = \det(\sigma)$. 

The next theorem is a synthesis of several foundational results from the theory of Weyl symmetric functions. 
All of these are proved in \cite{DonPosetModels} without any recourse to representation theory. 

\noindent 
{\bf \FTWSF\ (The Fundamental Theorem of Weyl Symmetric Functions)}\ \ 
{\sl Let $\mathscr{G}$ be a Coxeter--Dynkin posy, and let $\lambda$ be a dominant weight. Then: (1) There exists a unique Laurent polynomial $\chi = \chi_{_{\lambda}}$ satisfying} $\mathcal{A}(\myvarZ^{\varrho})\chi = \mathcal{A}(\myvarZ^{\varrho+\lambda})${\sl , and moreover $\chi_{_{\lambda}}$ is a $\mathcal{W}$-symmetric function.}

\noindent 
\hspace*{0.1in}\rule[-1.5mm]{0.5mm}{10mm}\hspace*{0.1in}\parbox[b]{6in}{\small We call each $\chi_{_{\lambda}}$ a {\em Weyl bialternant}. When we write $\chi_{_{\lambda}} = \sum_{\mu \in \Lambda}d_{\lambda,\mu}\myvarZ^{\mu}$, the integers $\{d_{\lambda,\mu}\}_{\mu \in \Lambda}$, at most finitely many of which are nonzero, are called $\mathscr{G}${\em -Kostka numbers}.}

\noindent 
{\sl (2) The set $\{\chi_{_{\lambda}}\}_{\lambda \in \Lambda^{+}}$ of Weyl bialternants is a $\mathbb{Z}$-basis and the set $\{\chi_{_{\omega_{i}}}\}_{i \in I}$ of fundamental $\mathcal{W}$-symmetric functions is an algebraic basis for the ring of $\mathcal{W}$-symmetric functions.}\\  
{\sl (3)  The $\mathscr{G}$-Kostka number $d_{\lambda,\mu}$ is nonzero if and only if $\mu \in \Pi(\lambda)$, in which case $d_{\lambda,\mu}$ is positive. 
Moreover, $\chi_{_{\lambda}} = \sum_{\nu \in \Pi(\lambda) \cap \Lambda^{+}} d_{\lambda,\nu}\zeta_{\nu}$; we have $d_{\lambda,\mu} = d_{\lambda,\sigma.\mu}$ for all $\sigma \in \mathcal{W}$ and $\mu \in \Lambda$; and $d_{\lambda,\lambda} = 1$.} 

{\em Proof.} Part {\sl (1)} follows from Theorems 2.6 and 2.9 of \cite{DonPosetModels}. 
Part {\sl (2)} follows from Theorem 2.14 of \cite{DonPosetModels}. 
The claims in {\sl (3)} that `$d_{\lambda,\mu}$ is nonzero if and only if $\mu \in \Pi(\lambda)$, in which case $d_{\lambda,\mu}$ is positive' and that `$\chi_{_{\lambda}} = \sum_{\nu \in \Pi(\lambda) \cap \Lambda^{+}} d_{\lambda,\nu}\zeta_{\nu}$' follow from Corollary 2.12 of \cite{DonPosetModels}. 
From Corollary 2.10 of \cite{DonPosetModels} we get that $d_{\lambda,\mu} = d_{\lambda,\sigma.\mu}$ for all $\sigma \in \mathcal{W}$ and $\mu \in \Lambda$. 
That $d_{\lambda,\lambda} = 1$ follows from Lemma 2.8 of \cite{DonPosetModels}.\hfill\QED

{\bf [\S \FlowerSection.7:\! Specializations.]} Given a $\mathcal{W}$-symmetric function $\chi = \sum_{\mu \in \Lambda}c_{\mu}\myvarZ^{\mu}$, where the monomial $\myvarZ^{\mu}$ is $\prod_{i \in I}z_{i}^{\langle \mu,\alpha_{i}^{\vee} \rangle}$, let $\chi\, \rule[-2mm]{0.2mm}{5mm}_{\left\{z_{i} := q^{\langle \omega_{i},\varrho^{\vee} \rangle}\right\}}$ be the Laurent polynomial in $q$ obtained by replacing each $z_{i}$ with $q^{\langle \omega_{i},\varrho^{\vee} \rangle}$. 
That is, $\chi\, \rule[-2mm]{0.2mm}{5mm}_{\left\{z_{i} := q^{\langle \omega_{i},\varrho^{\vee} \rangle}\right\}} = \sum_{\mu \in \Lambda}c_{\mu}q^{\langle \mu,\varrho^{\vee} \rangle}$. 
Notice that if we evaluate $\chi\, \rule[-2mm]{0.2mm}{5mm}_{\left\{z_{i} := q^{\langle \omega_{i},\varrho^{\vee} \rangle}\right\}}$ at $q=1$, we get $\sum_{\mu \in \Lambda}c_{\mu}$. 
Now, on the one hand, $w_{0}.\chi = \chi$, so $w_{0}.\chi\, \rule[-2mm]{0.2mm}{5mm}_{\left\{z_{i} := q^{\langle \omega_{i},\varrho^{\vee} \rangle}\right\}} = \chi\, \, \rule[-2mm]{0.2mm}{5mm}_{\left\{z_{i} := q^{\langle \omega_{i},\varrho^{\vee} \rangle}\right\}}$. 
But on the other hand, $w_{0}.\chi = \sum_{\mu \in \Lambda}c_{\mu}\myvarZ^{w_{0}.\mu}$, so $w_{0}.\chi\, \rule[-2mm]{0.2mm}{5mm}_{\left\{z_{i} := q^{\langle \omega_{i},\varrho^{\vee} \rangle}\right\}} = \sum_{\mu \in \Lambda}c_{\mu}q^{\langle w_{0}.\mu,\varrho^{\vee} \rangle} = \sum_{\mu \in \Lambda}c_{\mu}q^{\langle \mu,w_{0}.\varrho^{\vee} \rangle} = \sum_{\mu \in \Lambda}c_{\mu}q^{- \langle \mu,\varrho^{\vee} \rangle}$. 
Therefore the Laurent polynomial $\chi\, \rule[-2mm]{0.2mm}{5mm}_{\left\{z_{i} := q^{\langle \omega_{i},\varrho^{\vee} \rangle}\right\}} = \sum_{\mu \in \Lambda}c_{\mu}q^{\langle \mu,\varrho^{\vee} \rangle} = \sum_{\mu \in \Lambda}c_{\mu}q^{- \langle \mu,\varrho^{\vee} \rangle}$ is symmetric in that the coefficient of $q^{k}$ is the same as the coefficient of $q^{-k}$ for any integer $k$. 

There is a somewhat nicer polynomial version of this quantity which is developed in the next result. 
Part {\sl (1)} of the following statement seems to be due to Jacobson \cite{Jacobson}; {\sl (2)} and {\sl (3)} are immediate consequences of {\sl (1)}. 
Some terminology: A degree $l$ polynomial $a_{0} + a_{1}q + \cdots +  a_{l}q^{l}$ is {\em symmetric} if $a_{i} = a_{l-i}$ for all $i \in [0,l]_{\mathbb{Z}}$ and {\em unimodal} if there is some $u \in [0,l]_{\mathbb{Z}}$ such that $a_{0} \leq \cdots \leq a_{u} \geq \cdots \geq a_{l}$. 

\noindent 
{\bf \SpecializeTheorem\ (Specialization)}\ \ 
{\sl Let $\mathscr{G}$ be a Coxeter--Dynkin posy. Let $\lambda$ be any dominant weight. 
(1) The following is an identity of symmetric and unimodal polynomials of degree $2\langle \lambda,\varrho^{\vee} \rangle$ in the variable $q$:}
\[q^{\langle \lambda,\varrho^{\vee} \rangle}\chi_{\lambda}\, \rule[-2mm]{0.2mm}{5mm}_{\left\{z_{i} := q^{\langle \omega_{i},\varrho^{\vee} \rangle}\right\}} = \sum_{\mu \in \Lambda}d_{\lambda,\mu}q^{\langle \mu+\lambda,\varrho^{\vee} \rangle} = \frac{\prod_{\alpha \in \Phi^{+}}\left(1-q^{\langle \lambda+\varrho,\alpha^{\vee} \rangle}\right)}{\prod_{\alpha \in \Phi^{+}}\left(1-q^{\langle \varrho,\alpha^{\vee} \rangle}\right)}.\]
{\sl (2) Moreover, we have:} $\displaystyle \sum_{\mu \in \Lambda}d_{\lambda,\mu} = \frac{\prod_{\alpha \in \Phi^{+}}\langle \lambda+\varrho,\alpha^{\vee} \rangle}{\prod_{\alpha \in \Phi^{+}}\langle \varrho,\alpha^{\vee} \rangle}$, {\sl which is the number of terms of the polynomial of (1), counting multiplicities.}\\ 
{\sl (3) Further, the degree of the polynomial of (1) is given by} $\displaystyle \sum_{\alpha \in \Phi^{+}}\langle \lambda,\alpha^{\vee} \rangle$. 

{\em Proof.} See the proof of Theorem 2.17 of \cite{DonPosetModels} for a proof of part {\sl (1)}, from which parts {\sl (2)} and {\sl (3)} obviously follow.\hfill\QED
 
Curiously, the unimodality claim in part {\sl (1)} of the above statement is the only result in this section so far that depends on the representation theory of the semisimple Lie algebra associated with $\mathscr{G}$. 
We call the polynomial in part {\sl (1)} the {\em quantum dimension polynomial (QDP)} for $\chi_{_{\lambda}}$, the quantity in part {\sl (2)} the {\em dimension formula (DF)} for $\chi_{_{\lambda}}$, and the quantity in part {\sl (3)} the {\em height formula (HF)} for $\chi_{_{\lambda}}$ for reasons relating to representation theory. 
The DF counts the number of terms (including multiplicities) in the Laurent polynomial $\chi_{_{\lambda}}$, and the QDP  organizes these terms according to height within the set of weights $\Pi(\lambda)$ in such a way that the Weyl bialternant term $d_{\lambda,\mu}\myvarZ^{\mu}$ contributes the monomial $d_{\lambda,\mu}q^{\langle \mu+\lambda,\varrho^{\vee}\rangle}$ term to the QDP. 
The up-coming notion of `splitting posets' offers a more overtly combinatorial interpretation of the QDP and the DF. 

{\bf [\S \FlowerSection.8:\! Splitting posets.]} 
A {\em splitting poset} for a given $\mathcal{W}$-symmetric function $\chi$ is a $\mathscr{G}$-structured poset $R$ for which $\WGF(R;\myvarZ) = \chi$. 
The combinatorics of splitting posets for Weyl bialternants is rich and connects to many enumerative and extremal problems. 
We use the following notation in our statement of the next theorem.  
For any positive integer $m$, we define the $q$-integer $[m]_{q} := q^{m-1} + \cdots + q^{1} + 1$. 
It will be helpful to identify fundamental weights by type, so $\{\omega_{k}^{\mytinyX}\}_{k \in \{1,2,\ldots,n\}}$ are the fundamental weights in type $\myX \in \{\myA,\myB,\myC,\myD\}$. 
When $\displaystyle \lambda = \sum_{k=1}^{n}\lambda_{k}\omega_{k}^{\mytinyX}$ is a dominant weight, we let $\displaystyle \lambda_{i}^{j} := \sum_{k=i}^{j}\lambda_{k}$. 
We thank undergraduate student Nicholas Gaubatz for his substantive help in finding the product-of-quotients expression for the $\myD_{n}$ case in part {\sl (3)} of the following theorem; see \cite{Gaubatz} for a proof of these formulas that uses reasoning related to the Networked-numbers Game. 

\noindent 
{\bf \SplittingPosetTheorem\ (The Weyl Bialternant Splitting Poset Theorem)}\ \ {\sl Let $\mathscr{G}$ be a Coxeter--Dynkin posy. 
Let $\lambda$ be a dominant weight and$R$ a splitting poset for $\chi_{_{\lambda}}$.  
Regard $R$ to be ranked by $\rho: R \longrightarrow [0,\langle 2\lambda,\varrho^{\vee} \rangle]_{\mathbb{Z}}$ given by $\rho(\xelt) := \langle wt(\xelt)+\lambda,\varrho^{\vee} \rangle$ for all $\xelt \in R$, cf.\ \SaturatedLemma.}\\ 
{\sl (1) The rank-generating function for $R$ satisfies}  
\[\RGF(R;q) = \sum_{\xelt \in R}q^{\langle wt(\xelt)+\lambda,\varrho^{\vee} \rangle} = \frac{\prod_{\alpha \in \Phi^{+}}\left(1-q^{\langle \lambda+\varrho,\alpha^{\vee} \rangle}\right)}{\prod_{\alpha \in \Phi^{+}}\left(1-q^{\langle \varrho,\alpha^{\vee} \rangle}\right)},\] 
{\sl agrees with the QDP for $\chi_{_{\lambda}}$, and is a symmetric and unimodal polynomial of degree $2\langle \lambda,\varrho^{\vee} \rangle$. Therefore the cardinality and length of $R$ can also be expressed in terms of the DF and HF for $\chi_{_{\lambda}}$:} 
\[\CARD(R) = \frac{\prod_{\alpha \in \Phi^{+}}\langle \lambda+\varrho,\alpha^{\vee} \rangle}{\prod_{\alpha \in \Phi^{+}}\langle \varrho,\alpha^{\vee} \rangle} \hspace*{0.2in}\mbox{\sl and}\hspace*{0.2in} \posetlength(R) = \sum_{\alpha \in \Phi^{+}}\langle \lambda,\alpha^{\vee} \rangle.\] 
{\sl (2) Let $l$ be the number of positive roots in $\Phi$ and let $s_{i_{l}}\cdots s_{i_{2}}s_{i_{1}}$ be a reduced expression for $w_{0}$, the longest element of $\mathcal{W}$. 
Then $(\gamma_{i_{1}},\gamma_{i_{2}},\ldots,\gamma_{i_{l}})$ is a game sequence for a Networked-numbers Game played on $\mathscr{G}$ from initial position $\lambda+\varrho$ or from initial position $\varrho$. 
Also, if $c_{k}$ (respectively, $d_{k}$) is the number at node $\gamma_{i_{k}}$ when this node is fired in play originating from $\lambda+\varrho$ (resp.\ $\varrho$), then $c_{k} = \langle \lambda+\varrho,s_{i_{1}}s_{i_{2}} \cdots s_{i_{k-1}}.\alpha_{i_{k}}^{\vee} \rangle$ (resp.\ $d_{k} = \langle \varrho,s_{i_{1}}s_{i_{2}} \cdots s_{i_{k-1}}.\alpha_{i_{k}}^{\vee} \rangle$). Moreover,} 
\[\RGF(R;q) = \prod_{k=1}^{l}\frac{1-q^{c_{k}}}{1-q^{d_{k}}} \hspace*{0.2in}\mbox{\sl and}\hspace*{0.2in} \CARD(R) = \prod_{k=1}^{l}\frac{c_{k}}{d_{k}} \hspace*{0.2in}\mbox{\sl and}\hspace*{0.2in} \posetlength(R) = \sum_{k=1}^{l}(c_{k}-d_{k}).\]
{\sl (3) Now suppose} $\myX_{n} \in \{\myA_{n},\myB_{n},\myC_{n},\myD_{n}\}$, {\sl and write} $\displaystyle \lambda := \sum_{k=1}^{n}\lambda_{k}\omega_{k}^{\mytinyX}$. {\sl Then:}

\noindent 
\hspace*{1.5 em}
\fbox{$\myA_{n}$} $\displaystyle \RGF(R;q) = \prod_{i=1}^{n}\prod_{j=i}^{n}\frac{[\lambda_{i}^{j}+j+1-i]_{q}}{[j+1-i]_{q}}$

\noindent 
\hspace*{1.5 em}
\fbox{$\myB_{n}$} $\displaystyle \RGF(R;q) = \prod_{i=1}^{n-1}\prod_{j=i}^{n-1}\frac{[\lambda_{i}^{j}+j+1-i]_{q}}{[j+1-i]_{q}}\, \prod_{i=1}^{n}\prod_{j=i}^{n}\frac{[\lambda_{i}^{n}+\lambda_{j}^{n-1}+2n+1-i-j]_{q}}{[2n+1-i-j]_{q}}$

\noindent 
\hspace*{1.5 em}
\fbox{$\myC_{n}$} $\displaystyle \RGF(R;q) = \prod_{i=1}^{n}\prod_{j=i}^{n}\frac{[\lambda_{i}^{j}+j+1-i]_{q}}{[j+1-i]_{q}}\, \prod_{i=1}^{n-1}\prod_{j=i+1}^{n}\frac{[\lambda_{i}^{n}+\lambda_{j}^{n}+2n+2-i-j]_{q}}{[2n+2-i-j]_{q}}$

\noindent 
\hspace*{1.5 em}
\fbox{$\myD_{n}$} $\displaystyle \RGF(R;q) = \prod_{i=1}^{n-1}\prod_{j=i}^{n-1}\frac{[\lambda_{i}^{j}+j+1-i]_{q}}{[j+1-i]_{q}}\, \prod_{i=1}^{n-1}\prod_{j=i+1}^{n}\frac{[\lambda_{i}^{n-2}+\lambda_{j}^{n}+2n-i-j]_{q}}{[2n-i-j]_{q}}$

The proof of part {\sl (3)} is deferred until the end of this section. 

{\em Proof.} 
Part {\sl (1)} follows from Proposition 4.7 of \cite{DonPosetModels}. 
For {\sl (2)}, the facts that $\ell(w_{0}) = l = |\Phi^{+}|$ and that $(\gamma_{i_{1}},\gamma_{i_{2}},\ldots,\gamma_{i_{l}})$ is a game sequence when played from from $\lambda+\varrho$ when $w_{0} = s_{i_{l}}\cdots s_{i_{2}}s_{i_{1}}$ follow from \LongestLemma\ parts (1) and (2). 
Here is how we calculate the number $c_{k}$ at node $\gamma_{i_{k}}$ when that node is fired in the given game sequence: $c_{k} := \langle \mygens_{i_{k-1}}\cdots\mygens_{i_{1}}.(\lambda+\varrho),\alpha_{i_{k}}^{\vee} \rangle = \langle \lambda+\varrho,\mygens_{i_{1}}\cdots\mygens_{i_{k-1}}.\alpha_{i_{k}}^{\vee} \rangle$, where the latter equality follows from the $\mathcal{W}$-invariance of $\langle \cdot,\cdot \rangle$. 
By Lemma 5.1 of \cite{DonEur}, the set $\{\mygens_{i_{1}}\cdots\mygens_{i_{k-1}}.\alpha_{i_{k}}^{\vee}\}_{k=1}^{l}$ is exactly the set $\{\alpha^{\vee} \in (\Phi^{\vee})^{+}\, |\, w_{0}.\alpha^{\vee} \in (\Phi^{\vee})^{-}\}$, and by \LongestLemma.1 this latter set is just $(\Phi^{\vee})^{+}$. 
That is, $\{c_{k}\}_{k=1}^{l} \eqmulti \{\langle \lambda+\varrho,\alpha^{\vee} \rangle\}_{\alpha \in \Phi^{+}}$. 
The expressions for $\RGF(R;q)$, $\CARD(R)$, and $\posetlength(R)$ given in part {\sl (2)} of the theorem statement now follow from part {\sl (1)}.
The proof of {\sl (3)} can be found at the end of this section 
\hfill\QED

Some other (elementary) properties of splitting posets for Weyl bialternants can be found in \SaturatedLemma.6. 
The fact that every Weyl bialternant has a splitting poset is covered, for example, by Proposition 4.12 of \cite{DonPosetModels}. 

{\bf [\S \FlowerSection.9:\! Vertex coloring.]} 
There are many ways to confirm that a given $\mathscr{G}$-structured poset $R$ is splitting. 
If $R$ is shown to realize a representation of the associated complex semisimple Lie algebra $\mathfrak{g}(\mathscr{G})$, then $R$ is splitting. 
If $R'$ is $\mathscr{G'}$-structured for some Dynkin--Coxeter posy $\mathscr{G}'$ that contains $\mathscr{G}$ as a sub-embryophyte and if $R'$ is splitting, then we might be able to conclude that $R$ is splitting if $R$ is realized as a $\mathscr{G}$-component of $R'$ or as some felicitous recoloring of $R'$. 
Another approach, which we call the vertex-coloring method, is developed in \cite{DonPosetModels}, and we will utilize a special version of this method to support our hands-on work with examples in the next section. 

Let $\mathcal{C} := \mathcal{C}_{1} \times \cdots \times \mathcal{C}_{p}$ be a product of chains $\mathcal{C}_{1}, \ldots, \mathcal{C}_{p}$, and let $l_{q}$ be the length of chain $\mathcal{C}_{q}$. 
Fix $q \in [1,p]_{\mathbb{Z}}$. 
Any subset of $\mathcal{C}$ of the form $\{(\xelt_{1},\ldots,\xelt_{p}) \in \mathcal{C}\, |\, \mbox{rank of $\xelt_{q}$ in $\mathcal{C}_{q}$ is less than $l_{q}$}\}$ is a {\em sub-face of} $\mathcal{C}$ when given the induced partial order from $\mathcal{C}$. 
Say $\phi: \mathcal{C} \longrightarrow P$ is an isomorphism from $\mathcal{C}$ to some poset $P$. 
A subset $Q$ of $P$, with partial order induced from $P$, is a {\em sub-face of} $P$ if $Q$ is the image under $\phi$ of some sub-face of $\mathcal{C}$.  

Now suppose that $R$ is a $\mathscr{G}$-structured poset and that $\melt$ is a prominent element of $R$. 
Since $\delta_{i}(\melt)=0$ for all $i \in I$, then $wt(\melt)$ is a dominant weight.
Call any $\kappa: R \setminus \{\melt\} \longrightarrow I$ a {\em vertex-coloring function} on $R \setminus \{\melt\}$. 
For any $\xelt \in R \setminus \{\melt\}$, define the {\em kindred set} $\myK(\xelt)$ to be the subset $\{\yelt \in \comp_{\kappa(\xelt)}(\xelt)\, |\, \yelt \not= \melt \mbox{ and } \kappa(\yelt) = \kappa(\xelt)\}$ with partial order induced by $\comp_{\kappa(\xelt)}(\xelt)$. 
We say $\kappa$ is a {\em sub-face friendly} vertex-coloring function if, for all $\xelt \in R \setminus \{\melt\}$, the component $\comp_{\kappa(\xelt)}(\xelt)$ is isomorphic to a chain product and the kindred set $\myK(\xelt)$ is a sub-face of $\comp_{\kappa(\xelt)}(\xelt)$. 

\noindent 
{\bf \VertexColoringTheorem\ (Vertex Coloring)}\ \ {\sl Let $\mathscr{G}$ be a Coxeter--Dynkin posy.  
Let $R$ be a $\mathscr{G}$-structured poset with some prominent element $\melt$, and let $\lambda$ be the dominant weight $wt(\melt)$. 
Let $\kappa: R \setminus \{\melt\} \longrightarrow I$ be a vertex-coloring function and $\tau: R \setminus \{\melt\} \longrightarrow R \setminus \{\melt\}$ a bijection. 
Suppose at least one of the following two criteria is true: 
(1) All monochromatic components are rank symmetric and for all $\xelt \in R \setminus \{\melt\}$ we have $wt(\tau(\xelt)) = wt(\xelt) - (1+\mym_{\kappa(\xelt)})\alpha_{\kappa(\xelt)}$;  or 
(2) All monochromatic components of $R$ are (isomorphic to) chain products and $\kappa$ is sub-face friendly. 
Then $R$ is a splitting poset for $\chi_{_{\lambda}}$, so in particular} $\WGF(R) = \chi_{_{\lambda}}$. 

{\em Proof.} That criterion {\sl (1)} implies $R$ is a splitting poset for $\chi_{_{\lambda}}$ is Theorem 4.1 of \cite{DonPosetModels}; that {\sl (2)} implies this same result is Corollary 8.2.B of \cite{DonPosetModels}.\hfill\QED

The intricacy of the hypotheses of \VertexColoringTheorem\ can make it tricky to apply. 
However, as we will see in \S \QuasiExampleSection, it is actually fairly easy to apply when the monochromatic components of some given $\mathscr{G}$-structured poset are small and simply described. 

{\bf [\S \FlowerSection.10:\! Kac--Moody Lie algebras and representations in finite dimensions.]} 
Let us write our integral Coxeter--Dynkin posy $\mathscr{G}$ as the disjoint sum $\mathscr{G}_{1} \oplus \cdots \oplus \mathscr{G}_{m}$, where each $\mathscr{G}_{i}$ is an integral Coxeter--Dynkin flower with connected vascular graph $\Gamma(\mathscr{G}_{i})$. 
The Kac--Moody Lie algebra $\mathfrak{g}(\mathscr{G})$, defined in \S \WeylSection.12, is therefore isomorphic to $\mathfrak{g}(\mathscr{G}_{1}) \oplus \cdots \oplus \mathfrak{g}(\mathscr{G}_{m})$, and each $\mathfrak{g}(\mathscr{G}_{i})$ is finite-dimensional and simple by \KacMoodyCorollary. 
So, $\mathfrak{g} = \mathfrak{g}(\mathscr{G})$ is a finite-dimensional semisimple Lie algebra. 

Our next aim is to produce a realization $(\mathfrak{h}_{\mathscr{G}},\Delta_{\mathscr{G}},\widetilde{\Delta}_{\mathscr{G}})$ of $M = M(\mathscr{G})$. 
To create $\mathfrak{h}_{\mathscr{G}}$, we extend scalars in $\mathfrak{E}$ and take $\mathfrak{h}_{\mathscr{G}} := \mathfrak{E} \otimes_{\mathbb{R}} \mathbb{C}$ with $\widetilde{\Delta}_{\mathscr{G}} := \{\myroot_{i} \otimes 1\}_{i \in I} \subset \mathfrak{h}_{\mathscr{G}}$.  
For convenience, we set $\myqh_{i} := \myroot_{i} \otimes 1$ for $i \in I$. 
We also identify $\mathfrak{h}_{\mathscr{G}}^{*}$ with the complexification $\mathfrak{E}^{*} \otimes_{\mathbb{R}} \mathbb{C}$. 
By \InnerProductTheorem, we have a $\mathcal{W}$-invariant inner product $\langle \cdot,\cdot \rangle: \mathfrak{E}^{*} \times \mathfrak{E}^{*} \longrightarrow \mathbb{R}$ given by $\langle \mu,\nu \rangle := [\mu]S^{-1}[\nu]^{\mytinyT}$, where $[\xi]$ is the vector $\xi = \sum_{i \in I}\xi_{i}\omega_{i}$ thought of as a row vector $[\xi_{i}]_{i \in I}$ and where $D$ is a diagonal matrix with positive main diagonal entries such that $S = DM$ is a symmetric and positive-definite matrix. 
For any ${\boldsymbol \mu} \in \mathfrak{h}_{\mathscr{G}}^{*}$, let $[{\boldsymbol \mu}]$ denote the row vector associated with ${\boldsymbol \mu}$ when the latter is expressed in the basis $\{\omega_{i} \otimes 1\}_{i \in I}$, and let $[\overline{\boldsymbol \mu}]^{\mytinyT}$ be its conjugate transpose. 
Now define $B: \mathfrak{h}_{\mathscr{G}}^{*} \times \mathfrak{h}_{\mathscr{G}}^{*} \longrightarrow \mathbb{C}$ by $B({\boldsymbol \mu},{\boldsymbol \nu}) := [{\boldsymbol \mu}]S^{-1}[\overline{\boldsymbol \nu}]^{\mytinyT}$. 
It is easy to see that $B$ is a Hermitian form\footnote{It was not necessary to impose sesquilinearity here; we could accomplish  the same eventual purposes by using the bilinear form $B'$ given by $B'({\boldsymbol \mu},{\boldsymbol \nu}) := [{\boldsymbol \mu}]S^{-1}[{\boldsymbol \nu}]^{\mytinyT}$.} on $\mathfrak{h}_{\mathscr{G}}^{*}$. 
Let ${\boldsymbol \alpha}_{i} \in \mathfrak{h}_{\mathscr{G}}^{*}$ be the linear functional defined by ${\boldsymbol \mu} \mapsto B({\boldsymbol \mu},\alpha_{i}^{\vee} \otimes 1)$, and let $\widetilde{\Delta}_{\mathscr{G}} := \{{\boldsymbol \alpha}_{i}\}_{i \in I} \subset \mathfrak{h}_{\mathscr{G}}^{*}$. 
At this point, to see that $(\mathfrak{h}_{\mathscr{G}},\Delta_{\mathscr{G}},\widetilde{\Delta}_{\mathscr{G}})$ is a realization of $M = M(\mathscr{G})$, we must check two things. (1) We require that $\dim_{\mathbb{C}}(\mathfrak{h}_{\mathscr{G}}) = n$, which follows from the fact that co-rank of $M$ is $n - \mbox{\sffamily rank}(M) = 0$. 
(2) We require that ${\boldsymbol \alpha}_{j}(\alpha_{i}^{\vee} \otimes 1) = M_{ij}$ for all $i,j \in I$, which follows from a straightforward application of the definitions. 

Within this context, the defining relations $\mathcal{R}_{1}$ through $\mathcal{R}_{3}$ from \S \WeylSection.12 for our given Kac--Moody (and finite-dimensional semisimple) Lie algebra $\mathfrak{g}$ can be simplified as follows: 
\begin{description}
\item[($\mathcal{R}'_{1}$)] $[\myqh_{i}, \myqh_{j}] = 0$ for all $i,j \in I$ 
\item[($\mathcal{R}'_{2}$)] $[\myqh_{i}, \myqx_{j}] = M_{ji}\myqx_{i}$, $[\myqh_{i}, \myqy_{j}] = -M_{ji}\myqy_{i}$ for all $i,j \in I$
\item[($\mathcal{R}'_{3}$)] $[\myqx_{i}, \myqy_{j}] = \delta_{ij}\myqh_{i}$ for all $i, j \in I$
\item[($\mathcal{R}_{4}^{\pm}$)] For all $i, j \in I$ with $i \not= j$ we have ($+$) $(\mbox{\sffamily ad}\myqx_{i})^{1-M_{ji}}(\myqx_{j}) = 0$ and ($-$) $(\mbox{\sffamily ad}\myqy_{i})^{1-M_{ji}}(\myqy_{j}) = 0$ 
\end{description}
Note that the relations $\mathcal{R}_{4}^{\pm}$ remain unchanged. 
Let $\mathcal{R}'_{\mathscr{G}}$ be the collection of all the above relations.  
Take as our generating set $\mathcal{S}'_{\mathscr{G}} := \{\myqx_{i},\myqy_{i},\myqh_{i}\}_{i \in I(\mathscr{G})}$. 
So our finite-dimensional semisimple Lie algebra $\mathfrak{g}(\mathscr{G})$ is isomorphic to the Lie algebra $\langle \mathcal{S}'_{\mathscr{G}}\, |\, \mathcal{R}'_{\mathscr{G}} \rangle$ generated by $\{\myqx_{i},\myqy_{i},\myqh_{i}\}_{i \in I(\mathscr{G})}$ subject to the above relations.  

By restricting scalars from the $n$-dimensional complex vector spaces $\mathfrak{h}_{\mathscr{G}}$ and $\mathfrak{h}_{\mathscr{G}}^{*}$ to the $n$-dimensional (real) Euclidean spaces $\mathfrak{E}$ and $\mathfrak{E}^{*}$, we can view the set of roots $\mbox{\Fontauri \small R}(\mathscr{G})$, the root lattice $\mathcal{Q}_{\mathscr{G}} = \mathbb{Z}\Delta_{\mathscr{G}}$, and the weight lattice $\mathcal{P}_{\mathscr{G}}$ of \S \WeylSection.12 \& \S \WeylSection.13 as, respectively, our root system $\Phi(\mathscr{G})$, root lattice $\mathbb{Z}\{\alpha_{i}\, |\, i \in I\}$, and weight lattice $\Lambda$ from \S \FlowerSection.1, \S \FlowerSection.2, and \S \FlowerSection.4 above. 

Also, non-trivial finite-dimensional representations of $\mathfrak{g}$, which exist by \KMRepCorollary, are obviously integrable members of category $\mathcal{O}$. 
Say $V$ is a finite-dimensional $\mathfrak{g}$-module. 
Then there is a weight basis for $V$. 
We know from \S \WeylSection.13 that, amongst the category $\mathcal{O}$ representations of $\mathfrak{g}$, the irreducible representations are precisely the Verma module quotients $V(\lambda)$ for all $\lambda \in \Lambda^{+}$. 
In the language of \cite{Kac} and \cite{Kumar}, our pulsation matrix is, by \InnerProductTheorem, symmetrizable. 
This observation effects the next two results: the `Complete Reducibility' result below will follow from Theorem 10.7 of \cite{Kac}, and the `Weyl Character Formula' result below will follow from the more general Weyl--Kac Character Formula stated as Theorem 10.4 in \cite{Kac}. 

\noindent
{\bf \WeylKacCharacter\ (Weyl Character Formula)}\ \ {\sl Let $\mathscr{G}$ be an integral Coxeter--Dynkin posy with associated finite-dimensional semisimple Lie algebra $\mathfrak{g}(\mathscr{G})$. 
Let $\lambda \in \Lambda^{+}$ be a dominant weight with associated irreducible Verma module quotient $V(\lambda)$. 
Then} \[\mathcal{A}(\myvarZ^{\varrho})\, \mychar\left(V(\lambda)\right) = \mathcal{A}(\myvarZ^{\varrho+\lambda}),\] 
{\sl and so} $\mychar\left(V(\lambda)\right)$ {\sl is exactly the Weyl bialternant $\chi_{_{\lambda}}^{\mathscr{G}}$ of The Fundamental Theorem of Weyl Symmetric Functions (\FTWSF\ above). 
Moreover, $V(\lambda)$ is finite-dimensional.}

{\em Proof.} This is a special case of Weyl--Kac Character Formula; see Theorem 10.4 in \cite{Kac}.\hfill\QED

\noindent 
{\bf \CompletelyReducible\ (Complete Reducibility)}\ \ {\sl Let $\mathscr{G}$ be an integral Coxeter--Dynkin posy. 
Any finite-dimensional $\mathfrak{g}(\mathscr{G})$-module $V$ is completely reducible in that there exists a unique multiset of dominant weights $\{\lambda_{1},\ldots,\lambda_{p}\}$ such that}
\[V \cong V(\lambda_{1}) \oplus \cdots \oplus V(\lambda_{p}).\]
{\sl Also, if $W$ is a finite-dimensional $\mathfrak{g}(\mathscr{G})$-module, then $V \cong W$ if and only if} $\mychar(V) = \mychar(W)$.

{\em Proof.} The fact that $V \cong V(\lambda_{1}) \oplus \cdots \oplus V(\lambda_{p})$ for some multiset of dominant weights $\{\lambda_{1},\ldots,\lambda_{p}\}$ follows directly from Theorem 10.7 of \cite{Kac}. 
Uniqueness of this multiset follows from the uniqueness of the Verma module quotients as irreducible representations: $V(\lambda) \cong V(\nu)$ for dominant weights $\lambda$ and $\nu$ if and only if $\lambda = \nu$. 

Now suppose $\mychar(V) = \mychar(W)$, where $V \cong V(\lambda_{1}) \oplus \cdots \oplus V(\lambda_{p})$ and $W \cong V(\nu_{1}) \oplus \cdots \oplus V(\nu_{s})$. 
Fix some total order `$\leq_{I}$' on the colors in $I$. 
Let $\mymult(\lambda_{q})$ be the multiplicity of $\lambda_{q}$ within the multiset $\{\lambda_{1},\ldots,\lambda_{p}\}$, and similarly define $\mymult(\nu_{t})$ for $\nu_{t}$ in $\{\nu_{1},\ldots,\nu_{s}\}$. 
Temporarily impose a total order `$>_{\mbox{\tiny temp}}$' on the {\em set} $\{\lambda_{1},\ldots,\lambda_{p}\}$ by requiring $\lambda_{q} >_{\mbox{\tiny temp}} \lambda_{r}$ exactly when either (1) $\langle \lambda_{q},\varrho^{\vee} \rangle > \langle \lambda_{r},\varrho^{\vee} \rangle$ or (2) $\langle \lambda_{q},\varrho^{\vee} \rangle = \langle \lambda_{r},\varrho^{\vee} \rangle$ and there is some $i \in I$ such that $\langle \lambda_{q},\alpha_{j}^{\vee} \rangle = \langle \lambda_{q},\alpha_{j}^{\vee} \rangle$ for $j <_{I} i$ but $\langle \lambda_{q},\alpha_{i}^{\vee} \rangle > \langle \lambda_{q},\alpha_{i}^{\vee} \rangle$. 
Re-index elements of the set $\{\lambda_{1},\ldots,\lambda_{p}\}$ so that $\lambda_{q} >_{\mbox{\tiny temp}} \lambda_{r}$ if and only if $q < r$. 
Impose the same indexing/total order on the set $\{\nu_{1},\ldots,\nu_{s}\}$. 
Note that the monomial `$\myvarZ^{\lambda_{1}}$' has coefficient $\mymult(\lambda_{1})$ in $\mychar(V)$, as this monomial can only be contributed by the $\mymult(\lambda_{1})$ copies of $V(\lambda_{1})$ in the unique decomposition of $V$. 
Similarly, `$\myvarZ^{\nu_{1}}$' has coefficient $\mymult(\nu_{1})$ in $\mychar(W)$. 
Equality of the characters for $V$ and $W$ forces $\myvarZ^{\lambda_{1}} = \myvarZ^{\nu_{1}}$ (and hence $\lambda_{1} = \nu_{1}$) and $\mymult(\lambda_{1}) = \mymult(\nu_{1})$. 
Then $\mychar(V) - \mymult(\lambda_{1})\mychar(V(\lambda_{1})) = \mychar(W) - \mymult(\nu_{1})\mychar(V(\nu_{1}))$, and we can repeat the process to get $\myvarZ^{\lambda_{2}} = \myvarZ^{\nu_{2}}$ (and hence $\lambda_{2} = \nu_{2}$) and $\mymult(\lambda_{2}) = \mymult(\nu_{2})$, etc. 
This shows that the mutlisets $\{\lambda_{1},\ldots,\lambda_{p}\}$ and $\{\nu_{1},\ldots,\nu_{s}\}$ are identical, so $V \cong W$.\hfill\QED

{\bf [\S \FlowerSection.11:\! DCML's as supporting graphs.]} 
Here we will say precisely what it means for a DCML to serve as a model for a semisimple Lie algebra representation. 
We take $\mathscr{G} = (\Gamma,M)$ to be a Coxeter--Dynkin posy and $\mathfrak{g} = \mathfrak{g}(\mathscr{G})$ to be its associated finite-dimensional semisimple Lie algebra. 
All vector spaces here (including $\mathfrak{g}$) are over $\mathbb{C}$. 
If $V$ is a finite-dimensional $\mathfrak{g}$-module, then $V$ is an integrable representation in category $\mathcal{O}$. 
In particular, $V$ is a direct sum of its finite-dimensional weight spaces: $V = \bigoplus_{\mu \in \Lambda} V_{\mu}$. 
A {\em weight basis} for $V$ is any basis that respects the foregoing weight space decomposition. 
For any nonzero $v \in V_{\mu}$, we call $v$ a {\em weight vector} and say the {\em weight of} $v$ is $wt(v) := \mu$.  
Let $R$ be an indexing set of size $d:=\dim(V)$ for some weight basis $\{v_{\relt}\}_{\relt \in R}$ of $V$. 
Fix $\selt, \telt \in R$. 
Expand $\myqx_{i}.v_{\selt}$ in the basis $\{v_{\xelt}\}_{\xelt \in R}$ by writing $\myqx_{i}.v_{\selt} = \sum_{\uelt \in R}\myqX^{(i)}_{\uelt,\selt}v_{\uelt}$ for some scalars $\{\myqX^{(i)}_{\uelt,\selt}\}_{\uelt \in R}$, and similarly write $\myqy_{i}.v_{\telt} = \sum_{\relt \in R}\myqY^{(i)}_{\relt,\telt}v_{\relt}$ for some scalars $\{\myqY^{(i)}_{\relt,\telt}\}_{\relt \in R}$. 
For a weight basis vector $v_{\selt}$ with weight $\mu$, i.e.\ $v_{\selt} \in V_{\mu}$, recall from \S \WeylSection.13 that $\myqx_{i}.v_{\selt} \in V_{\mu + \alpha_{i}}$. 
This means that $\myqX^{(i)}_{\uelt,\selt} = 0$ if $v_{\uelt} \not\in V_{\mu + \alpha_{i}}$. 
Similarly, if $v_{\telt} \in V_{\nu}$, then $\myqy_{i}.v_{\telt} \in V_{\nu - \alpha_{i}}$, so $\myqY^{(i)}_{\relt,\telt} = 0$ if $v_{\relt} \not\in V_{\nu-\alpha_{i}}$. 

Continuing with the set-up from the previous paragraph, create an edge-colored directed graph as follows.  
Use elements of our indexing set $R$ as vertices, and create a colored and directed edge $\selt \myarrow{i} \telt$ for any $\selt, \telt \in R$ if and only if at least one of $\myqX^{(i)}_{\telt,\selt}$ or $\myqY^{(i)}_{\selt,\telt}$ is nonzero. 
This edge-colored directed graph, also referred to as $R$, is the {\em supporting graph} for the given weight basis $\{v_{\relt}\}_{\relt \in R}$ of $V$. 
The {\em representation diagram} for the given weight basis $\{v_{\relt}\}_{\relt \in R}$ of $V$ is the supporting graph $R$ together with the collection of scalar pairs $\{(\myqX^{(i)}_{\telt,\selt},\myqY^{(i)}_{\selt,\telt})\}_{i \in I\, \mbox{\large ,}\ \selt \myarrow{i} \telt \in \SmallEdgeSet_{i}(R)}$.  
As an abuse of notation, this representation diagram is also referred to as $R$. 
A consequence of the observations at the end of the previous paragraph is that if $\selt \myarrow{i} \telt$ in $R$, then $wt(v_{\selt}) + \alpha_{i} = wt(v_{\telt})$. 
In particular, if $\selt \myarrow{i} \telt$ and $\selt \myarrow{j} \telt$, then $wt(v_{\telt}) = wt(v_{\selt}) + \alpha_{i} = wt(v_{\selt}) + \alpha_{j}$, and independence of our simple root vectors forces $i=j$. 
Also, if $\selt \myarrow{i} \telt$ and $\telt \myarrow{j} \selt$, then $wt(v_{\telt}) = wt(v_{\selt}) + \alpha_{i} = wt(v_{\telt}) + \alpha_{i} + \alpha_{j}$, which is impossible due to independence of our simple root vectors. 
We conclude that there is at most one edge between any two elements of the supporting graph $R$.  

Our supporting graph $R$ is actually the covering digraph for a poset. 
To see this, we begin with an argument that $R$ is an acyclic digraph. 
Suppose $(\relt_{0},\relt_{1},\ldots,\relt_{p})$ is a path in $R$ such that $(\relt_{0},\relt_{1},\ldots,\relt_{p-1})$ is simple, $\relt_{0} = \relt_{p}$, and $\relt_{q} \myarrow{i_{q}} \relt_{q+1}$ for some $i_{q} \in I$ whenever $q \in [0,p-1]_{\mathbb{Z}}$.  Now say $wt(v_{\relt_{0}}) = \mu$. 
Then $\mu + \alpha_{i_{0}} + \alpha_{i_{1}} + \cdots + \alpha_{i_{p-1}} = \mu$. 
Now linear independence of the $\alpha_{i}$'s forces us to have $p=0$, and hence $R$ is an acyclic graph. 
So, we may impose the natural partial order on $R$ whereby $\selt \leq_{R} \telt$ if and only if there is some path $(\selt = \relt_{0},\ldots,\relt_{p}=\telt)$ from $\selt$ to $\telt$ wherein $\relt_{q} \myarrow{i_{q}} \relt_{q+1}$ for some $i_{q} \in I$ whenever $q \in [0,p-1]_{\mathbb{Z}}$. 
Notice that some $\telt$ covers some $\selt$ in $R$ only if we have $\selt \myarrow{i} \telt$ for some color $i$. 
On the other hand, suppose that $\selt \myarrow{i} \telt$ and that there is a path $(\selt = \relt_{0},\ldots,\relt_{p}=\telt)$ in the supporting graph $R$ from $\selt$ to $\telt$ wherein $\relt_{q} \myarrow{i_{q}} \relt_{q+1}$ for some $i_{q} \in I$ whenever $q \in [0,p-1]_{\mathbb{Z}}$. 
Then $wt(v_{\selt})+\alpha_{i} = wt(v_{\telt}) = wt(v_{\selt}) + \sum_{q=1}^{p}\alpha_{i_{q}}$. 
Again, independence of our simple roots means that $p=1$, that $i_{1} = i$, and that $\telt$ covers $\selt$ in our partial order. 
That is, the colored and directed edges of the supporting graph $R$ are exactly the covering relations for the poset $(R,\leq_{R})$. 

Supporting graphs have many nice combinatorial properties. 
Some elementary combinatorial properties were explored in Lemmas 3.1 and 3.2 of \cite{DonSupp}, and here we recapitulate some statements from those lemmas.  
Every connected component of $R$ is a ranked poset. 
Moreover, $wt(v_{\relt}) = \sum_{i \in I} \mym_{i}(\relt)\omega_{i} = \sum_{i \in I}\big(\rho_{i}(\relt)-\delta_{i}(\relt)\big)\omega_{i} = wt(\relt)$. 
That is, the weight of an element of $R$ is the weight of the corresponding weight basis vector. 
Also, within each connected component of $R$, elements with the same weight must have the same rank. 
The supporting graph $R$ is $\mathscr{G}$-structured. 
Also, suppose weights $\mu$ and $\mu+\alpha_{i}$ are weights from the saturated set of weights $\mathscr{W}(V)$. 
Then there is at least one pair of weight basis vectors $v_{\selt} \in V_{\mu}$ and $v_{\telt} \in V_{\mu+\alpha_{i}}$ such that $\selt \myarrow{i} \telt$. 

These latter facts have consequences for diamonds in our $\mathscr{G}$-structured representation diagram $R$. 
Consider a diamond $\,$ \parbox{1.4cm}{\begin{center}
\setlength{\unitlength}{0.2cm}
\begin{picture}(6.5,3.5)
\put(3,0){\circle*{0.5}} 
\put(1,2){\circle*{0.5}}
\put(3,4){\circle*{0.5}} 
\put(5,2){\circle*{0.5}}
\put(1,2){\line(1,1){2}} 
\put(3,0){\line(-1,1){2}}
\put(5,2){\line(-1,1){2}} 
\put(3,0){\line(1,1){2}}
\put(1.75,0.55){\em \small k} 
\put(3.5,0.5){\em \small l}
\put(1.7,2.7){\em \small i} 
\put(3.75,2.55){\em \small j}
\put(3.5,-0.75){\footnotesize $\relt$} 
\put(5.75,1.75){\footnotesize $\telt$}
\put(3.5,4.5){\footnotesize $\uelt$} 
\put(-0.5,1.75){\footnotesize $\selt$}
\end{picture} \end{center}} in $R$. 
Since $wt(\selt) + \alpha_{i}-\alpha_{j} = wt(\telt) = wt(\selt) -\alpha_{k}+\alpha_{l}$, we get $\alpha_{i}+\alpha_{k} = \alpha_{j}+\alpha_{l}$. 
Linear independence of these simple roots now forces (1) $i=j$ and $k=l$ or (2) $i=l$ and $j=k$. 
Assume for the moment that all scalars on all edges of this representation diagram are nonzero and that $R$ is topographically balanced. 
Then, in our given diamond, when we expand $\myqy_{j}.(\myqx_{i}.v_{\selt})$ in the basis $\{v_{\xelt}\}_{\xelt \in R}$ for our representing space $V$, the coefficient of $v_{\telt}$ is $\myqX_{\uelt,\selt}^{(i)}\myqY_{\telt,\uelt}^{(j)}$. 
Similarly, the coefficient of $v_{\telt}$ in $\myqx_{l}.(\myqx_{k}.v_{\selt})$ is $\myqY_{\relt,\selt}^{(k)}\myqX_{\telt,\relt}^{(l)}$. 
But, considering the relations enjoyed by these generators of $\mathfrak{g}(\mathscr{G})$, the coefficient of $v_{\telt}$ when we expand $[\myqx_{i},\myqy_{j}].v_{\selt}$ must be zero. 
In particular, there must be some $\relt'$ and edges $\selt \mybackarrow{j} \relt' \myarrow{i} \telt$ in our supporting graph $R$. 
But topographical balance forces $\relt'=\relt$ and therefore $i=l$ and $j=k$. 
That is, diamond-coloring of topographically balanced supporting graphs is a natural consequence of having nonzero edge coefficients in the companion representation diagram. 

Next we describe how representation diagrams can be obtained synthetically. 
In what follows, all edge colors are from our color palette $I$. 
An {\em edge-tagged poset} is an edge-colored ranked poset $R$ together with the assignment of a scalar pair $(\myqX_{\telt,\selt},\myqY_{\selt,\telt})$ to each edge $\selt \myarrow{i} \telt$. 
The quantity $\myqP_{\selt,\telt} := \myqX_{\telt,\selt}\myqY_{\selt,\telt}$ is the {\em edge product} for edge $\selt \myarrow{i} \telt$. 
Now suppose $Q$ is another edge-tagged poset. 
We say $R$ is {\em edge-product similar to} $Q$ if there is an edge- and edge-color-preserving isomorphism $\phi: R \longrightarrow Q$ such that $\myqP_{\selt,\telt} = \myqP_{\phi(\selt),\phi(\telt)}$ whenever there is an edge from $\selt$ to $\telt$ in $R$. 
Call $R$ a {\em properly} edge-tagged poset if, whenever $\selt \myarrow{i} \telt$ is an edge in $R$, then at least one of the scalars from the scalar pair $(\myqX_{\telt,\selt},\myqY_{\selt,\telt})$ is nonzero. 
For such a properly edge-tagged poset $R$, let $V[R]$ be the complex vector space freely generated by the symbols $\{v_{\relt}\}_{\relt \in R}$. 
For each $i \in I$ and $\selt, \telt \in R$, define actions of generators $\myqx_{i}$ and $\myqy_{i}$ on $V[R]$ as follows: 
\begin{eqnarray*}
\myqx_{i}.v_{\selt} & := & \sum_{\telt' \in R\, \mbox{\large ,}\, \selt \myarrow{i} \telt'}\myqX_{\telt',\selt}v_{\telt'}\\
\myqy_{i}.v_{\telt} & := & \sum_{\selt' \in R\, \mbox{\large ,}\, \selt' \myarrow{i} \telt}\myqY_{\selt',\telt}v_{\selt'}
\end{eqnarray*}
Now suppose $R$ is topographically balanced. 
The {\em crossing relation} of color $i \in I$ at vertex $\selt \in R$ is the equation 
\[\sum_{\relt \in R\, \mbox{\large ,}\, \relt \myarrow{i} \selt}\myqP_{\relt,\selt} - \sum_{\telt \in R\, \mbox{\large ,}\, \selt \myarrow{i} \telt}\myqP_{\selt,\telt} = \mym_{i}(\selt) = \rho_{i}(\selt) -\delta_{i}(\selt).\]
The {\em diamond relations} for a given diamond of edges$\,$ 
\parbox{1.4cm}{\begin{center}
\setlength{\unitlength}{0.2cm}
\begin{picture}(6.5,3.5)
\put(3,0){\circle*{0.5}} 
\put(1,2){\circle*{0.5}}
\put(3,4){\circle*{0.5}} 
\put(5,2){\circle*{0.5}}
\put(1,2){\line(1,1){2}} 
\put(3,0){\line(-1,1){2}}
\put(5,2){\line(-1,1){2}} 
\put(3,0){\line(1,1){2}}
\put(1.75,0.55){\em \small j} 
\put(3.5,0.5){\em \small i}
\put(1.7,2.7){\em \small i} 
\put(3.75,2.55){\em \small j}
\put(3.5,-0.75){\footnotesize $\relt$} 
\put(5.75,1.75){\footnotesize $\telt$}
\put(3.5,4.5){\footnotesize $\uelt$} 
\put(-0.5,1.75){\footnotesize $\selt$}
\end{picture} \end{center}} is the pair of equations $\myqY_{\telt,\uelt}\myqX_{\uelt,\selt} = \myqX_{\selt,\relt}\myqY_{\relt,\telt}$ and $\myqY_{\selt,\uelt}\myqX_{\uelt,\telt} = \myqX_{\telt,\relt}\myqY_{\relt,\selt}$. 
The set of all possible crossing and diamond relations is called the {\em CD relations}, and our topographically balanced and edge-tagged poset $R$ {\em satisfies the CD relations} if the scalar pairs assigned to the edges of $R$ together comprise a solution to the CD relations. 
The {\em diamond relation for edge products} on a given diamond of edges$\,$ 
\parbox{1.4cm}{\begin{center}
\setlength{\unitlength}{0.2cm}
\begin{picture}(6.5,3.5)
\put(3,0){\circle*{0.5}} 
\put(1,2){\circle*{0.5}}
\put(3,4){\circle*{0.5}} 
\put(5,2){\circle*{0.5}}
\put(1,2){\line(1,1){2}} 
\put(3,0){\line(-1,1){2}}
\put(5,2){\line(-1,1){2}} 
\put(3,0){\line(1,1){2}}
\put(1.75,0.55){\em \small j} 
\put(3.5,0.5){\em \small i}
\put(1.7,2.7){\em \small i} 
\put(3.75,2.55){\em \small j}
\put(3.5,-0.75){\footnotesize $\relt$} 
\put(5.75,1.75){\footnotesize $\telt$}
\put(3.5,4.5){\footnotesize $\uelt$} 
\put(-0.5,1.75){\footnotesize $\selt$}
\end{picture} \end{center}} is the equation $\myqP_{\selt,\uelt}\myqP_{\telt,\uelt} = \myqP_{\relt,\selt}\myqP_{\relt,\telt}$. 
The {\em CD relations for edge products} are the usual crossing relations together with the diamond relations for edge products. 

Let $R$ be the supporting graph for a weight basis $\{v_{\relt}\}_{\relt \in R}$ of a $\mathfrak{g}(\mathscr{G})$-module $V$. 
We say $R$ is {\em solitary} if, whenever $\{w_{\relt}\}_{\relt \in R}$ is another weight basis for $V$ that has supporting graph $R$, there exist scalars $\{c_{\relt}\}_{\relt \in R}$ such that $w_{\relt} = c_{\relt}v_{\relt}$ for all $\relt \in R$. 
(In general, when two bases differ by scalars in this way, we say the bases are {\em diagonally equivalent}, since, relative to appropriate total orderings of the bases, the change-of-basis matrix is diagonal; if this diagonal matrix is a multiple of the identity matrix, we say the bases are {\em scalar equivalent}.)
It is easy to see that if two weight bases share the same solitary supporting graph, then the bases are diagonally equivalent and their respective representation diagrams are edge-product similar. 
We say $R$ is {\em edge-minimizing} if no other supporting graph for $V$ has fewer edges than $R$, and we say $R$ is {\em edge-minimal} if no proper subgraph of $R$ is a supporting graph for $V$. 
Intuitively, we think of an edge-minimal supporting graph as being `locally edge-minimizing'. 
Of course, an edge-minimizing supporting graph is edge-minimal.
Finally, we say $R$ is a {\em modular (resp.\ distributive) lattice} supporting graph if the poset $R$ is a modular (resp.\ distributive) lattice. 
Each of the terms {\em modular / distributive lattice}, {\em edge-minimal}, {\em edge-minimizing}, and/or {\em solitary} is also used as an adjective describing any associated representation diagram $R$ and affiliated weight basis $\{w_{\relt}\}_{\relt \in R}$, if said term also describes the supporting graph $R$.  

The next result can be applied to DCML's but also to connected, topographically balanced, and diamond-colored posets that are not lattices. 
Parts {\sl (1)--(3)} of this theorem have appeared in some form elsewhere, but part {\sl (4)} is new. 

\noindent 
{\bf \MainLieRepTheorem}\ \ {\sl Suppose $R$ is a connected, topographically balanced, diamond-colored, and properly edge-tagged poset. 
(1) Then $R$ is $\mathscr{G}$-structured and satisfies the CD relations if and only if $R$ is a representation diagram for the weight basis $\{v_{\relt}\}_{\relt \in R}$ of the $\mathfrak{g}(\mathscr{G})$-module $V[R]$, where generator actions are as defined in the preceding paragraphs. 
(2) Assume that the equivalent conditions stated in part (1) are true for $R$ and that the edge product $\myqP_{\selt,\telt} \not= 0$  whenever $\selt \myarrow{i} \telt$ is an edge in $R$. 
Further, suppose that there is an isomorphism $\phi: R \longrightarrow Q$ of edge-colored posets where $Q$ is edge-tagged and $\myqP_{\phi(\selt),\phi(\telt)} = \myqP_{\selt,\telt}$. 
Then there exist scalars $\{c_{\relt}\}_{\relt \in R}$ such that the weight basis $\{c_{\relt}v_{\relt}\}_{\relt \in R}$ for the $\mathfrak{g}(\mathscr{G})$-module $V[R]$ has representation diagram $Q$. 
(3) Keep the hypotheses from the first sentence of (2), and assume that the $\mathfrak{g}(\mathscr{G})$-module $V[R]$ is irreducible. 
Further, assume that for every edge-tagged poset $Q$ such that $Q$ and $R$ are isomorphic edge-colored posets, it is the case that $Q$ is edge-product similar to $R$. 
Then $R$ is solitary. 
(4) Again use the hypotheses of the first sentence of (2), and suppose that the CD relations for edge products on $R$ have a unique solution. Then $R$ is edge-minimal.}

{\em Proof.} The equivalence stated as {\sl (1)} follows easily from Proposition 3.4 of \cite{DonSupp}; see Lemma 3.1 of \cite{DLP2} for a statement that uses the language of crossing and diamond relations. 
Part {\sl (2)} is Lemma 4.2 of \cite{DLP2}. 
Part {\sl (3)} is Lemma 4.3 of \cite{DLP2}.

For part {\sl (4)}, we need the following claim concerning uncolored ranked posets. \underline{Claim:} Say $P$ is a finite poset expressed as a disjoint sum of ranked posets. 
For any integer $\mu$, let $P_{\mu} := \{\xelt \in P\, |\, \mym(\xelt) = \mu\}$, where as usual, $\mym(\xelt) = \rho(\xelt)-\delta(\xelt)$ when $\rho$ and $\delta$ are the unique rank and depth functions in the connected component of $\xelt$. 
Next, suppose $P'$ is a disjoint sum of ranked posets obtained by removing some edges from $P$ but keeping all vertices. 
The elements of $P'$ are $\{\xelt'\, |\, \xelt \in P\}$, and define $P'_{\mu}$, $\mym'$, $\rho'$, and $\delta'$ as for $P$. 
If $|P_{\mu}| = |P'_{\mu}|$ for all $\mu \in \mathbb{Z}$, then $\mym(\xelt) = \mym'(\xelt')$ for all $\xelt \in P$. 

\underline{Proof of Claim:} Let $M := \max\{\mym(\xelt)\, |\, \xelt \in P\}$. 
We induct on $M-\mu$ to prove that for any integer $\mu \in [0,M]_{\mathbb{Z}}$ and any $\xelt \in P$, we have $\mym'(\xelt') = \mu \Longrightarrow \mym(\xelt) = \mu$ and $\mym'(\xelt') = -\mu \Longrightarrow \mym(\xelt) = -\mu$. 
If so, then the mapping $P'_{\mu} \longrightarrow P_{\mu}$ given by $\xelt' \mapsto \xelt$ is injective and, since $|P'_{\mu}| = |P_{\mu}|$, also surjective; similarly the mapping $P'_{-\mu} \longrightarrow P_{-\mu}$ given by $\xelt' \mapsto \xelt$ is a bijection. 
It will follow, then, that $\mym'(\xelt') = \mym(\xelt)$ for all $\xelt \in P$. 

For the basis step of the induction, say $M - \mu = 0$, i.e.\ $\mu = M$, and suppose $\mym'(\xelt') = M$. 
If $\xelt'$ is not rank maximal within its connected component of $P'$, then there is some rank-maximal $\yelt'$ with $\rho'(\yelt') > \rho'(\xelt')$. 
Moreover, $\delta'(\xelt') = \rho'(\yelt')-\rho'(\xelt')$. 
So, $M= \rho'(\xelt')-\delta'(\xelt) = \rho'(\yelt')-2\delta'(\xelt') < \rho'(\yelt')$. 
However, $\rho'(\yelt') \leq \rho(\yelt)$, and in $\comp(\yelt) \subseteq P$ we can find a rank-maximal $\zelt$ with $\rho(\yelt) \leq \rho(\zelt)$. 
Then, $M < \rho(\zelt)$, which contradicts the maximality of $M$. 
So, when $\mym'(\xelt') = M$, then $\xelt'$ is rank maximal within its connected component. 
Similarly, when $\mym'(\xelt') = -M$, then $\xelt'$ is rank minimal in its connected component. 

Continuing with our hypothesis that $\mym'(\xelt') = M$, then rank maximality of $\xelt'$ in its connected component means that $\rho'(\xelt') = M$. 
Since $\rho'(\xelt') \leq \rho(\xelt)$, then $M \leq \rho(\xelt)$. 
Let $\yelt$ be rank maximal in $\comp(\xelt) \subseteq P$, so that $M \leq \rho(\xelt) \leq \rho(\yelt) = \mym(\yelt)$, and maximality of $M$ means that $M = \mym(\yelt)$ and that the preceding inequalities are actually equalities. 
So $\xelt$ is rank maximal in $\comp(\xelt) \subseteq P$, hence $M = \mym(\xelt)$. 
In a similar way, we see that $\mym'(\xelt') = -M \Longrightarrow \mym(\xelt) = -M$. 

For the induction step, suppose that for some nonnegative integer $K$, it is the case that for all $\kappa \in [K+1,M]_{\mathbb{Z}}$, we have $\mym'(\xelt') = \kappa \Longrightarrow \mym(\xelt) = \kappa$ and $\mym'(\xelt') = -\kappa \Longrightarrow \mym(\xelt) = -\kappa$. 
We wish to show that $\mym'(\xelt') = K \Longrightarrow \mym(\xelt) = K$ and $\mym'(\xelt') = -K \Longrightarrow \mym(\xelt) = -K$. 
Suppose now that $\mym'(\xelt') = K$. 
If $\xelt'$ is not rank maximal in $\comp(\xelt') \subseteq P'$, then let $\yelt'$ be rank maximal within this component. 
So, $K < \mym'(\yelt)$, and, by our inductive hypothesis, $\mym'(\yelt') = \mym(\yelt)$. 
A path in $P'$ from $\xelt'$ to $\yelt'$ is also a path in $P$ from $\xelt$ to $\yelt$. 
Because $\mym(\qelt)-\mym(\pelt) = 2$ whenever $\pelt \rightarrow \qelt$ in $P$ and $\mym(\qelt')-\mym(\pelt') = 2$ whenever $\pelt' \rightarrow \qelt'$ in $P'$, we conclude that $\mym'(\zelt') = \mym(\zelt)$ for any $\zelt'$ in our chosen path from $\xelt'$ to $\yelt'$. 
It follows that $\mym(\xelt')=\mym(\xelt)$. 

So, now suppose $\xelt'$ is rank maximal within its connected component in $P'$. 
We intend to rule out $\mym(\xelt) < K$ and $\mym(\xelt) > K$. 
Let $\kappa := \mym(\xelt)$. 
First, consider $\kappa > K$. 
By our inductive hypothesis, whenever $\mym'(\yelt') = \kappa$, then $\mym(\yelt) = \kappa$, so the mapping $\yelt' \mapsto \yelt$ injects $P'_{\kappa}$ into $P_{\kappa}$, and, since $|P'_{\kappa}| = |P_{\kappa}|$, then this mapping is a bijection. 
So, $\xelt \in P_{\kappa}$ means $\xelt' \in P'_{\kappa}$, contradicting the fact that $\mym'(\xelt') = K < \kappa$. 
We thusly rule out the possibility that $\kappa > K$. 
Second, consider $\kappa < K$. 
If $\xelt$ is rank maximal in its connected component, we have $K = \mym'(\xelt') = \rho'(\xelt') \leq \rho(\xelt) = \mym(\xelt) = \kappa$, violating the hypothesis that $\kappa < K$. 
So, $\xelt$ is not rank maximal in $\comp(\xelt) \subseteq P$. 
Let $\yelt$ (resp.\ $\belt$) be rank maximal (resp.\ rank minimal) in this connected component, and set $N := \mym(\yelt) = \rho(\yelt)$ (so $-N = \mym(\belt) = -\delta(\belt)$). 
Also, let $\aelt'$ be rank minimal in $\comp(\xelt') \subseteq P'$, so $\mym(\aelt') = -K$. 
Since $\mym'(\xelt')-\mym'(\aelt') = 2K$, then $\mym(\xelt)-\mym(\aelt)=2K$ also, so $\melt(\aelt) = \kappa-2K$. 
Since $K > \kappa$, then $-K > \kappa-2K = \mym(\aelt)$. 
Again, by our inductive hypothesis, whenever $\mym'(\yelt') = \kappa-2K$, then $\mym(\yelt) = \kappa-2K$, so the mapping $\yelt' \mapsto \yelt$ injects $P'_{\kappa-2K}$ into $P_{\kappa-2K}$, and, since $|P'_{\kappa-2K}| = |P_{\kappa-2K}|$, then this mapping is a bijection. 
So, $\aelt \in P_{\kappa-2K}$ means $\aelt' \in P'_{\kappa-2K}$, contradicting the fact that $\mym'(\aelt') = -K > \kappa-2K$. 
This reasoning rules out the possibility that $\kappa < K$. 
We conclude that $\mym(\xelt) = \kappa = K$, as desired. 
This completes the induction proof of our Claim.

We return now to the setting of part {\sl (4)} of the theorem statement. 
Regard $R$ to be a representation diagram for a $\mathfrak{g}(\mathscr{G})$-module $V = V[R]$, and assume that $R$ is connected, topographically balanced, and diamond-colored with nonzero edge products. 
Suppose further that $R'$ is some subgraph of $R$ that is also a supporting graph for $V$. 
We will show that $R'$ and $R$ are identical as edge-colored directed graphs. 
Let $\myqP'_{\selt',\telt'}$ be the edge product on an edge $\selt' \myarrow{i} \telt'$ associated with a solution to the CD relations for edge products on the supporting graph $R'$. 
Define a new set of edge products on $R$ by the rule 
\[\myqQ_{\selt,\telt} := \left\{\begin{array}{cl} \myqP'_{\selt',\telt'} & \hspace*{0.0in}\mbox{if $\selt' \myarrow{i} \telt'$ is an edge in $R'$}\\ 0 & \hspace*{0.0in}\mbox{otherwise}\end{array}\right.\]
whenever $\selt \myarrow{i} \telt$ is an edge in the supporting graph $R$. 
We claim that $\{\myqQ_{\pelt,\qelt}\}_{\pelt \rightarrow \qelt \mbox{\scriptsize\, in } R}$ comprises a solution to the CD relations for edge products on $R$. 

Indeed, for $i \in I$ and $\selt' \in R'$, we have 
\[\sum_{\relt' \in R'\, \mbox{\large ,}\, \relt' \myarrow{i} \selt'}\myqP'_{\relt',\selt'} - \sum_{\telt' \in R'\, \mbox{\large ,}\, \selt' \myarrow{i} \telt'}\myqP'_{\selt',\telt'} = \mym_{i}(\selt')\] 
and, by our Claim, $\mym_{i}(\selt) = \mym_{i}(\selt')$. 
It follows that 
\[\sum_{\relt \in R\, \mbox{\large ,}\, \relt \myarrow{i} \selt}\myqQ_{\relt,\selt} - \sum_{\telt \in R\, \mbox{\large ,}\, \selt \myarrow{i} \telt}\myqQ_{\selt,\telt} = \mym_{i}(\selt).\] 
So, the numbers $\{\myqQ_{\pelt,\qelt}\}_{\pelt \rightarrow \qelt \mbox{\scriptsize\, in } R}$ satisfies all crossing relations in $R$. 
Now consider a diamond$\,$ \parbox{1.4cm}{\begin{center}
\setlength{\unitlength}{0.2cm}
\begin{picture}(6.5,3.5)
\put(3,0){\circle*{0.5}} 
\put(1,2){\circle*{0.5}}
\put(3,4){\circle*{0.5}} 
\put(5,2){\circle*{0.5}}
\put(1,2){\line(1,1){2}} 
\put(3,0){\line(-1,1){2}}
\put(5,2){\line(-1,1){2}} 
\put(3,0){\line(1,1){2}}
\put(1.75,0.55){\em \small j} 
\put(3.5,0.5){\em \small i}
\put(1.7,2.7){\em \small i} 
\put(3.75,2.55){\em \small j}
\put(3.5,-0.75){\footnotesize $\relt$} 
\put(5.75,1.75){\footnotesize $\telt$}
\put(3.5,4.5){\footnotesize $\uelt$} 
\put(-0.5,1.75){\footnotesize $\selt$}
\end{picture} \end{center}} in $R$, but suppose$\,$ \parbox{1.4cm}{\begin{center}
\setlength{\unitlength}{0.2cm}
\begin{picture}(6.5,3.5)
\put(3,0){\circle*{0.5}} 
\put(1,2){\circle*{0.5}}
\put(3,4){\circle*{0.5}} 
\put(5,2){\circle*{0.5}}
\put(3,0){\line(-1,1){2}}
\put(5,2){\line(-1,1){2}} 
\put(3,0){\line(1,1){2}}
\put(1.75,0.55){\em \small j} 
\put(3.5,0.5){\em \small i}
\put(3.75,2.55){\em \small j}
\put(3.5,-0.75){\footnotesize $\relt'$} 
\put(5.75,1.75){\footnotesize $\telt'$}
\put(3.5,4.5){\footnotesize $\uelt'$} 
\put(-0.5,1.75){\footnotesize $\selt'$}
\end{picture} \end{center}} is the corresponding set of edges in $R'$. 
Then the coefficient $c_{1}$ of the $V$-basis vector $v_{\telt'}$ when we expand $\myqx_{i}.(\myqy_{j}.v_{\selt'})$ in the weight basis $\{v_{\xelt}\}_{\xelt \in R}$ is $\myqX_{\telt',\relt'}\myqY_{\relt',\selt'}$. 
Suppose for the moment that the coefficient $c_{2}$ of $v_{\telt'}$ in $\myqy_{j}.(\myqx_{i}.v_{\selt'})$ is nonzero. 
Then, there must be some $\tilde{\uelt}' \in R'$ such that $\selt' \myarrow{i} \tilde{\uelt}' \mybackarrow{j} \telt'$. 
If so, these latter edges in $R'$ are also edges in $R$, and the topographical balance property that $R$ possesses then forces $\tilde{\uelt} = \uelt$. 
In particular, $\selt' \myarrow{i} \uelt'$ is an edge in $R'$, contrary to our hypothesis. 
We conclude that the coefficient $c_{2} = 0$. 
According to the relations satisfied by our defining generators for $\mathfrak{g}(\mathscr{G})$, the coefficient of $v_{\telt}$ in $[\myqx_{i},\myqy_{j}].v_{\selt}$ must be zero, so $\myqX_{\telt',\relt'}\myqY_{\relt',\selt'} = c_{1} = c_{2} = 0$. 
Therefore $\myqP'_{\relt',\selt'}\myqP'_{\relt',\telt'} = 0$, and hence $\myqQ_{\relt,\selt}\myqQ_{\relt,\telt} = 0$. 
Since $\myqQ_{\selt,\uelt} = 0$, then $\myqQ_{\selt,\uelt}\myqQ_{\telt,\uelt} = 0$, which means that $\myqQ_{\relt,\selt}\myqQ_{\relt,\telt} = \myqQ_{\selt,\uelt}\myqQ_{\telt,\uelt}$. 
We can similarly analyze other configurations of edges amongst $\relt'$, $\selt'$, $\telt'$, and $\uelt'$ in $R'$ to see that on our given diamond in $R$, the numbers $\myqQ_{\relt,\selt}$, $\myqQ_{\relt,\telt}$, $\myqQ_{\selt,\uelt}$, and $\myqQ_{\telt,\uelt}$ satisfy the associated diamond relation for edge products. 
Since we chose our diamond in $R$ generically, it follows that the numbers $\{\myqQ_{\pelt,\qelt}\}_{\pelt \rightarrow \qelt \mbox{\scriptsize\, in } R}$ satisfy all diamond relations for edge products in $R$. 

However, by hypothesis $\{\myqP_{\pelt,\qelt}\}_{\pelt \rightarrow \qelt \mbox{\scriptsize\, in } R}$ is the unique set of numbers that satisfy the CD relations for edge products in $R$. 
It follows that each $\myqQ_{\pelt,\qelt}$ coincides with the nonzero quantity $\myqP_{\pelt,\qelt}$. 
So, each edge in $R$ must also be an edge in $R'$, meaning $R'$ and $R$ are the same edge-colored digraph.\hfill\QED 

The next exercise is intended to help the reader suss out some of the notational/linguistic technicalities. 

\noindent 
{\bf \TopoNotLatticeExercise}\ \ Regard the connected and topographically balanced (and therefore ranked) poset $R$ of \TopoBalancedRankedFigure\ to be edge-colored by a one-element set, so $R$ is monochromatic. 
(a) Verify that $R$ is $\myA_{1}$-structured. 
(a) Find a set of positive and rational edge coefficients so that the resulting edge-tagged poset, which we name `$R_{1}$', satisfies the CD relations. 
(c) Find another set of positive and rational edge coefficients so that the resulting edge-tagged poset, which we name `$R_{2}$', satisfies the CD relations and is {\sc not} edge-product similar to $R_{1}$. 

We close this section with a proof of part deferred from \S \FlowerSection.8.

{\em Proof of part {\sl (3)} of \SplittingPosetTheorem.} 
The formulas of {\sl (3)} can be derived as follows.  
For dominant $\lambda = \sum_{i \in I}\lambda_{i}\omega_{i}$ and a positive root $\alpha$ such that $\alpha^{\vee} = \sum_{i \in I}k_{i}\alpha_{i}^{\vee}$, the quantity $\langle \lambda+\varrho,\alpha^{\vee} \rangle$ can be re-expressed as $\sum_{i \in I}k_{i}(\lambda_{i}+1)$. 
So, to concretely evaluate the numerator terms in our part {\sl (1)} formula for $\RGF(R;q)$, it suffices to have concrete descriptions of each positive co-root $\alpha^{\vee}$ in terms of the simple co-roots $\alpha_{i}^{\vee}$. 
For $\myA_{n}$,\ $\myB_{n}$ (with co-roots from\ $\myC_{n}$),\ $\myC_{n}$ (with co-roots from\ $\myB_{n}$), and\ $\myD_{n}$, we can describe the set of positive roots in terms of the simple roots by analyzing the `top half' of any

\newcommand{\ASixFlower}{
\setlength{\unitlength}{0.5cm}
\begin{picture}(17,1)
\put(-2.5,-0.25){\LARGE $\myA_{6}$}
\put(-0.9,0){\vector(1,0){1.2}}
{\color{Cyan}\put(1,0){\circle*{0.7}}}
\put(1,0){\circle{0.7}}
{\color{Red}\put(4,0){\circle*{0.7}}}
\put(4,0){\circle{0.7}}
{\color{Purple}\put(7,0){\circle*{0.7}}}
\put(7,0){\circle{0.7}}
{\color{OliveGreen}\put(10,0){\circle*{0.7}}}
\put(10,0){\circle{0.7}}
{\color{BurntOrange}\put(13,0){\circle*{0.7}}}
\put(13,0){\circle{0.7}}
{\color{SpringGreen}\put(16,0){\circle*{0.7}}}
\put(16,0){\circle{0.7}}
\thicklines
\put(1.35,0){\line(1,0){2.3}}
\put(4.35,0){\line(1,0){2.3}}
\put(7.35,0){\line(1,0){2.3}}
\put(10.35,0){\line(1,0){2.3}}
\put(13.35,0){\line(1,0){2.3}}
\put(0.8,0.6){\textcolor{Cyan}{\em 1}}
\put(3.8,0.6){\textcolor{Red}{\em 2}}
\put(6.8,0.6){\textcolor{Purple}{\em 3}}
\put(9.8,0.6){\textcolor{OliveGreen}{\em 4}}
\put(12.8,0.6){\textcolor{BurntOrange}{\em 5}}
\put(15.8,0.6){\textcolor{SpringGreen}{\em 6}}
\end{picture}
}

\newcommand{\DSixFlower}{
\setlength{\unitlength}{0.5cm}
\begin{picture}(14,1)
\put(16,-0.35){\LARGE $\myD_{6}$}
\put(15.6,0.05){\vector(-1,0){1.2}}
{\color{Cyan}\put(1,0){\circle*{0.7}}}
\put(1,0){\circle{0.7}}
{\color{Red}\put(4,0){\circle*{0.7}}}
\put(4,0){\circle{0.7}}
{\color{Purple}\put(7,0){\circle*{0.7}}}
\put(7,0){\circle{0.7}}
{\color{OliveGreen}\put(10,0){\circle*{0.7}}}
\put(10,0){\circle{0.7}}
{\color{SpringGreen}\put(13,-1){\circle*{0.7}}}
\put(13,-1){\circle{0.7}}
{\color{BurntOrange}\put(13,1){\circle*{0.7}}}
\put(13,1){\circle{0.7}}
\thicklines
\put(1.35,0){\line(1,0){2.3}}
\put(4.35,0){\line(1,0){2.3}}
\put(7.35,0){\line(1,0){2.3}}
\put(10.36,0.1){\line(3,1){2.3}}
\put(10.36,-0.1){\line(3,-1){2.3}}
\put(0.8,0.6){\textcolor{Cyan}{\em 1}}
\put(3.8,0.6){\textcolor{Red}{\em 2}}
\put(6.8,0.6){\textcolor{Purple}{\em 3}}
\put(9.8,0.6){\textcolor{OliveGreen}{\em 4}}
\put(13.4,0.5){\textcolor{BurntOrange}{\em 5}}
\put(13.4,-1){\textcolor{SpringGreen}{\em 6}}
\end{picture}
}

\newcommand{\ASixRoots}{
\setlength{\unitlength}{1cm}
\begin{picture}(7,5)
\put(-2,5){\large $\myA_{6}$ positive roots}
\put(-2,0){\TypeEboxDot{Gray}}
\put(-1.9,-0.175){\tiny 000001}
\put(0,0){\TypeEboxDot{Gray}}
\put(0.1,-0.175){\tiny 000010}
\put(2,0){\TypeEboxDot{Gray}}
\put(2.1,-0.175){\tiny 000100}
\put(4,0){\TypeEboxDot{Gray}}
\put(4.1,-0.175){\tiny 001000}
\put(6,0){\TypeEboxDot{Gray}}
\put(6.1,-0.175){\tiny 010000}
\put(8,0){\TypeEboxDot{Gray}}
\put(8.1,-0.175){\tiny 100000}
\put(-1,1){\TypeEboxDot{Gray}}
\put(-1.3,1.5){\tiny 000011}
\put(1,1){\TypeEboxDot{Gray}}
\put(1.75,1.175){\tiny 000110}
\put(3,1){\TypeEboxDot{Gray}}
\put(3.75,1.175){\tiny 001100}
\put(5,1){\TypeEboxDot{Gray}}
\put(5.75,1.175){\tiny 011000}
\put(7,1){\TypeEboxDot{Gray}}
\put(7.75,1.175){\tiny 110000}
\put(0,2){\TypeEboxDot{Gray}}
\put(-0.3,2.5){\tiny 000111}
\put(2,2){\TypeEboxDot{Gray}}
\put(2.75,2.175){\tiny 001110}
\put(4,2){\TypeEboxDot{Gray}}
\put(4.75,2.175){\tiny 011100}
\put(6,2){\TypeEboxDot{Gray}}
\put(6.75,2.175){\tiny 111000}
\put(1,3){\TypeEboxDot{Gray}}
\put(0.7,3.5){\tiny 001111}
\put(3,3){\TypeEboxDot{Gray}}
\put(3.75,3.175){\tiny 011110}
\put(5,3){\TypeEboxDot{Gray}}
\put(5.75,3.175){\tiny 111100}
\put(2,4){\TypeEboxDot{Gray}}
\put(1.7,4.5){\tiny 011111}
\put(4,4){\TypeEboxDot{Gray}}
\put(4.75,4.175){\tiny 111110}
\put(3,5){\TypeEboxDot{Gray}}
\put(2.7,5.5){\tiny 111111}
\thicklines
\put(6.625,0.375){\color{Cyan}\qbezier(0,0)(0.375,0.375)(0.75,0.75)}
\put(5.625,1.375){\color{Cyan}\qbezier(0,0)(0.375,0.375)(0.75,0.75)}
\put(4.625,2.375){\color{Cyan}\qbezier(0,0)(0.375,0.375)(0.75,0.75)}
\put(3.625,3.375){\color{Cyan}\qbezier(0,0)(0.375,0.375)(0.75,0.75)}
\put(2.625,4.375){\color{Cyan}\qbezier(0,0)(0.375,0.375)(0.75,0.75)}
\put(4.625,0.375){\color{Red}\qbezier(0,0)(0.375,0.375)(0.75,0.75)}
\put(3.625,1.375){\color{Red}\qbezier(0,0)(0.375,0.375)(0.75,0.75)}
\put(2.625,2.375){\color{Red}\qbezier(0,0)(0.375,0.375)(0.75,0.75)}
\put(1.625,3.375){\color{Red}\qbezier(0,0)(0.375,0.375)(0.75,0.75)}
\put(2.625,0.375){\color{Purple}\qbezier(0,0)(0.375,0.375)(0.75,0.75)}
\put(1.625,1.375){\color{Purple}\qbezier(0,0)(0.375,0.375)(0.75,0.75)}
\put(0.625,2.375){\color{Purple}\qbezier(0,0)(0.375,0.375)(0.75,0.75)}
\put(0.625,0.375){\color{OliveGreen}\qbezier(0,0)(0.375,0.375)(0.75,0.75)}
\put(-0.375,1.375){\color{OliveGreen}\qbezier(0,0)(0.375,0.375)(0.75,0.75)}
\put(-1.375,0.375){\color{BurntOrange}\qbezier(0,0)(0.375,0.375)(0.75,0.75)}
\put(4.39,4.375){\color{SpringGreen}\qbezier(0,0)(-0.375,0.375)(-0.75,0.75)}
\put(5.39,3.375){\color{BurntOrange}\qbezier(0,0)(-0.375,0.375)(-0.75,0.75)}
\put(6.39,2.375){\color{OliveGreen}\qbezier(0,0)(-0.375,0.375)(-0.75,0.75)}
\put(7.39,1.375){\color{Purple}\qbezier(0,0)(-0.375,0.375)(-0.75,0.75)}
\put(8.39,0.375){\color{Red}\qbezier(0,0)(-0.375,0.375)(-0.75,0.75)}
\put(3.39,3.375){\color{SpringGreen}\qbezier(0,0)(-0.375,0.375)(-0.75,0.75)}
\put(4.39,2.375){\color{BurntOrange}\qbezier(0,0)(-0.375,0.375)(-0.75,0.75)}
\put(5.39,1.375){\color{OliveGreen}\qbezier(0,0)(-0.375,0.375)(-0.75,0.75)}
\put(6.39,0.375){\color{Purple}\qbezier(0,0)(-0.375,0.375)(-0.75,0.75)}
\put(2.39,2.375){\color{SpringGreen}\qbezier(0,0)(-0.375,0.375)(-0.75,0.75)}
\put(3.39,1.375){\color{BurntOrange}\qbezier(0,0)(-0.375,0.375)(-0.75,0.75)}
\put(4.39,0.375){\color{OliveGreen}\qbezier(0,0)(-0.375,0.375)(-0.75,0.75)}
\put(1.39,1.375){\color{SpringGreen}\qbezier(0,0)(-0.375,0.375)(-0.75,0.75)}
\put(2.39,0.375){\color{BurntOrange}\qbezier(0,0)(-0.375,0.375)(-0.75,0.75)}
\put(0.39,0.375){\color{SpringGreen}\qbezier(0,0)(-0.375,0.375)(-0.75,0.75)}
\put(-1.1,0.65){\color{BurntOrange}{\em 5}}
\put(-0.1,1.65){\color{OliveGreen}{\em 4}}
\put(0.9,0.65){\color{OliveGreen}{\em 4}}
\put(0.9,2.65){\color{Purple}{\em 3}}
\put(1.9,1.65){\color{Purple}{\em 3}}
\put(2.9,0.65){\color{Purple}{\em 3}}
\put(1.9,3.65){\color{Red}{\em 2}}
\put(2.9,2.65){\color{Red}{\em 2}}
\put(3.9,1.65){\color{Red}{\em 2}}
\put(4.9,0.65){\color{Red}{\em 2}}
\put(2.9,4.65){\color{Cyan}{\em 1}}
\put(3.9,3.65){\color{Cyan}{\em 1}}
\put(4.9,2.65){\color{Cyan}{\em 1}}
\put(5.9,1.65){\color{Cyan}{\em 1}}
\put(6.9,0.65){\color{Cyan}{\em 1}}
\put(7.9,0.65){\color{Red}{\em 2}}
\put(6.9,1.65){\color{Purple}{\em 3}}
\put(5.9,0.65){\color{Purple}{\em 3}}
\put(5.9,2.65){\color{OliveGreen}{\em 4}}
\put(4.9,1.65){\color{OliveGreen}{\em 4}}
\put(3.9,0.65){\color{OliveGreen}{\em 4}}
\put(4.9,3.65){\color{BurntOrange}{\em 5}}
\put(3.9,2.65){\color{BurntOrange}{\em 5}}
\put(2.9,1.65){\color{BurntOrange}{\em 5}}
\put(1.9,0.65){\color{BurntOrange}{\em 5}}
\put(3.9,4.65){\color{SpringGreen}{\em 6}}
\put(2.9,3.65){\color{SpringGreen}{\em 6}}
\put(1.9,2.65){\color{SpringGreen}{\em 6}}
\put(0.9,1.65){\color{SpringGreen}{\em 6}}
\put(-0.1,0.65){\color{SpringGreen}{\em 6}}
\end{picture}}

\newcommand{\DSixRoots}{
\setlength{\unitlength}{1cm}
\begin{picture}(7,8)
\put(-0.3,9){\large $\myD_{6}$ positive roots}
\put(-1,0){\TypeEboxDot{Gray}}
\put(-0.9,-0.175){\tiny 000010}
\put(0,0){\TypeEboxDot{Gray}}
\put(0.1,-0.175){\tiny 000001}
\put(1,0){\TypeEboxDot{Gray}}
\put(1.1,-0.175){\tiny 000100}
\put(3,0){\TypeEboxDot{Gray}}
\put(3.1,-0.175){\tiny 001000}
\put(5,0){\TypeEboxDot{Gray}}
\put(5.1,-0.175){\tiny 010000}
\put(7,0){\TypeEboxDot{Gray}}
\put(7.1,-0.175){\tiny 100000}
\put(0,1){\TypeEboxDot{Gray}}
\put(-0.5,1.175){\tiny 000110}
\put(1,1){\TypeEboxDot{Gray}}
\put(1.75,1.175){\scriptsize $\pelt$}
\put(6.75,1.7){\tiny $\pelt = 000101$}
\put(2,1){\TypeEboxDot{Gray}}
\put(2.75,1.175){\tiny 001100}
\put(4,1){\TypeEboxDot{Gray}}
\put(4.75,1.175){\tiny 011000}
\put(6,1){\TypeEboxDot{Gray}}
\put(6.75,1.175){\tiny 110000}
\put(0,2){\TypeEboxDot{Gray}}
\put(-0.5,2.175){\tiny 000111}
\put(1,2){\TypeEboxDot{Gray}}
\put(1.75,2.175){\scriptsize $\qelt$}
\put(6.75,2){\tiny $\qelt = 011100$}
\put(2,2){\TypeEboxDot{Gray}}
\put(2.75,2.175){\scriptsize $\relt$}
\put(6.75,2.3){\tiny $\relt = 001101$}
\put(3,2){\TypeEboxDot{Gray}}
\put(3.75,2.175){\tiny 011100}
\put(5,2){\TypeEboxDot{Gray}}
\put(5.75,2.175){\tiny 111000}
\put(1,3){\TypeEboxDot{Gray}}
\put(0.5,3.175){\tiny 001111}
\put(2,3){\TypeEboxDot{Gray}}
\put(2.75,3.175){\scriptsize $\selt$}
\put(5.75,3){\tiny $\selt = 011110$}
\put(3,3){\TypeEboxDot{Gray}}
\put(3.75,3.175){\scriptsize $\telt$}
\put(5.75,3.3){\tiny $\telt = 011101$}
\put(4,3){\TypeEboxDot{Gray}}
\put(4.75,3.175){\tiny 111100}
\put(0,4){\TypeEboxDot{Gray}}
\put(-0.5,4.175){\tiny 001211}
\put(2,4){\TypeEboxDot{Gray}}
\put(1.5,4.175){\tiny 011111}
\put(3,4){\TypeEboxDot{Gray}}
\put(3.75,4.175){\scriptsize $\uelt$}
\put(5.75,4.175){\tiny $\uelt = 111111$}
\put(4,4){\TypeEboxDot{Gray}}
\put(4.75,4.175){\tiny 111101}
\put(1,5){\TypeEboxDot{Gray}}
\put(0.5,5.175){\tiny 011211}
\put(3,5){\TypeEboxDot{Gray}}
\put(3.75,5.175){\tiny 111111}
\put(0,6){\TypeEboxDot{Gray}}
\put(-0.5,6.175){\tiny 012211}
\put(2,6){\TypeEboxDot{Gray}}
\put(2.75,6.175){\tiny 111211}
\put(1,7){\TypeEboxDot{Gray}}
\put(1.75,7.175){\tiny 112211}
\put(0,8){\TypeEboxDot{Gray}}
\put(0.75,8.175){\tiny 122211}
\thicklines
\put(5.625,0.375){\color{Cyan}\qbezier(0,0)(0.375,0.375)(0.75,0.75)}
\put(4.625,1.375){\color{Cyan}\qbezier(0,0)(0.375,0.375)(0.75,0.75)}
\put(3.625,2.375){\color{Cyan}\qbezier(0,0)(0.375,0.375)(0.75,0.75)}
\put(3.625,3.375){\color{Cyan}\qbezier(0,0)(0.375,0.375)(0.75,0.75)}
\put(2.625,3.375){\color{Cyan}\qbezier(0,0)(0.375,0.375)(0.75,0.75)}
\put(2.625,4.375){\color{Cyan}\qbezier(0,0)(0.375,0.375)(0.75,0.75)}
\put(1.625,5.375){\color{Cyan}\qbezier(0,0)(0.375,0.375)(0.75,0.75)}
\put(0.625,6.375){\color{Cyan}\qbezier(0,0)(0.375,0.375)(0.75,0.75)}
\put(3.625,0.375){\color{Red}\qbezier(0,0)(0.375,0.375)(0.75,0.75)}
\put(2.625,1.375){\color{Red}\qbezier(0,0)(0.375,0.375)(0.75,0.75)}
\put(2.625,2.375){\color{Red}\qbezier(0,0)(0.375,0.375)(0.75,0.75)}
\put(1.625,2.375){\color{Red}\qbezier(0,0)(0.375,0.375)(0.75,0.75)}
\put(1.625,3.375){\color{Red}\qbezier(0,0)(0.375,0.375)(0.75,0.75)}
\put(0.625,4.375){\color{Red}\qbezier(0,0)(0.375,0.375)(0.75,0.75)}
\put(1.625,0.375){\color{Purple}\qbezier(0,0)(0.375,0.375)(0.75,0.75)}
\put(1.625,1.375){\color{Purple}\qbezier(0,0)(0.375,0.375)(0.75,0.75)}
\put(0.625,1.375){\color{Purple}\qbezier(0,0)(0.375,0.375)(0.75,0.75)}
\put(0.625,2.375){\color{Purple}\qbezier(0,0)(0.375,0.375)(0.75,0.75)}
\put(0.625,0.375){\color{OliveGreen}\qbezier(0,0)(0.375,0.375)(0.75,0.75)}
\put(-0.375,0.375){\color{OliveGreen}\qbezier(0,0)(0.375,0.375)(0.75,0.75)}
\put(1.39,7.375){\color{Red}\qbezier(0,0)(-0.375,0.375)(-0.75,0.75)}
\put(2.39,6.375){\color{Purple}\qbezier(0,0)(-0.375,0.375)(-0.75,0.75)}
\put(3.39,5.375){\color{OliveGreen}\qbezier(0,0)(-0.375,0.375)(-0.75,0.75)}
\put(4.39,4.375){\color{BurntOrange}\qbezier(0,0)(-0.375,0.375)(-0.75,0.75)}
\put(4.39,3.375){\color{BurntOrange}\qbezier(0,0)(-0.375,0.375)(-0.75,0.75)}
\put(5.39,2.375){\color{OliveGreen}\qbezier(0,0)(-0.375,0.375)(-0.75,0.75)}
\put(6.39,1.375){\color{Purple}\qbezier(0,0)(-0.375,0.375)(-0.75,0.75)}
\put(7.39,0.375){\color{Red}\qbezier(0,0)(-0.375,0.375)(-0.75,0.75)}
\put(1.39,5.375){\color{Purple}\qbezier(0,0)(-0.375,0.375)(-0.75,0.75)}
\put(2.39,4.375){\color{OliveGreen}\qbezier(0,0)(-0.375,0.375)(-0.75,0.75)}
\put(3.39,3.375){\color{BurntOrange}\qbezier(0,0)(-0.375,0.375)(-0.75,0.75)}
\put(3.39,2.375){\color{BurntOrange}\qbezier(0,0)(-0.375,0.375)(-0.75,0.75)}
\put(4.39,1.375){\color{OliveGreen}\qbezier(0,0)(-0.375,0.375)(-0.75,0.75)}
\put(5.39,0.375){\color{Purple}\qbezier(0,0)(-0.375,0.375)(-0.75,0.75)}
\put(1.39,3.375){\color{OliveGreen}\qbezier(0,0)(-0.375,0.375)(-0.75,0.75)}
\put(2.39,2.375){\color{BurntOrange}\qbezier(0,0)(-0.375,0.375)(-0.75,0.75)}
\put(2.39,1.375){\color{BurntOrange}\qbezier(0,0)(-0.375,0.375)(-0.75,0.75)}
\put(3.39,0.375){\color{OliveGreen}\qbezier(0,0)(-0.375,0.375)(-0.75,0.75)}
\put(1.39,1.375){\color{BurntOrange}\qbezier(0,0)(-0.375,0.375)(-0.75,0.75)}
\put(1.39,0.375){\color{BurntOrange}\qbezier(0,0)(-0.375,0.375)(-0.75,0.75)}
\put(-0.1,0.65){\color{OliveGreen}{\em 4}}
\put(0.7,0.45){\color{OliveGreen}{\em 4}}
\put(0.9,2.65){\color{Purple}{\em 3}}
\put(1.1,1.85){\color{Purple}{\em 3}}
\put(1.7,1.45){\color{Purple}{\em 3}}
\put(1.9,0.65){\color{Purple}{\em 3}}
\put(0.9,4.65){\color{Red}{\em 2}}
\put(1.9,3.65){\color{Red}{\em 2}}
\put(2.1,2.85){\color{Red}{\em 2}}
\put(2.7,2.45){\color{Red}{\em 2}}
\put(2.9,1.65){\color{Red}{\em 2}}
\put(3.9,0.65){\color{Red}{\em 2}}
\put(0.9,6.65){\color{Cyan}{\em 1}}
\put(1.9,5.65){\color{Cyan}{\em 1}}
\put(2.9,4.65){\color{Cyan}{\em 1}}
\put(3.1,3.85){\color{Cyan}{\em 1}}
\put(3.7,3.45){\color{Cyan}{\em 1}}
\put(3.9,2.65){\color{Cyan}{\em 1}}
\put(4.9,1.65){\color{Cyan}{\em 1}}
\put(5.9,0.65){\color{Cyan}{\em 1}}
\put(0.9,7.65){\color{Red}{\em 2}}
\put(6.9,0.65){\color{Red}{\em 2}}
\put(0.9,5.65){\color{Purple}{\em 3}}
\put(1.9,6.65){\color{Purple}{\em 3}}
\put(4.9,0.65){\color{Purple}{\em 3}}
\put(5.9,1.65){\color{Purple}{\em 3}}
\put(0.9,3.65){\color{OliveGreen}{\em 4}}
\put(1.9,4.65){\color{OliveGreen}{\em 4}}
\put(2.9,5.65){\color{OliveGreen}{\em 4}}
\put(2.9,0.65){\color{OliveGreen}{\em 4}}
\put(3.9,1.65){\color{OliveGreen}{\em 4}}
\put(4.9,2.65){\color{OliveGreen}{\em 4}}
\put(3.7,3.85){\color{BurntOrange}{\em 5}}
\put(2.7,2.85){\color{BurntOrange}{\em 5}}
\put(1.7,1.85){\color{BurntOrange}{\em 5}}
\put(0.7,0.85){\color{BurntOrange}{\em 5}}
\put(1.1,1.45){\color{BurntOrange}{\em 5}}
\put(2.1,2.45){\color{BurntOrange}{\em 5}}
\put(3.1,3.45){\color{BurntOrange}{\em 5}}
\put(3.9,4.65){\color{BurntOrange}{\em 5}}
%
\put(3.505,4.41){\color{SpringGreen}\line(0,1){0.68}}
\put(3.4,4.61){\color{SpringGreen}{\em 6}}
\put(2.505,3.41){\color{SpringGreen}\line(0,1){0.68}}
\put(2.4,3.61){\color{SpringGreen}{\em 6}}
\put(1.505,2.41){\color{SpringGreen}\line(0,1){0.68}}
\put(1.4,2.61){\color{SpringGreen}{\em 6}}
\put(0.505,1.41){\color{SpringGreen}\line(0,1){0.68}}
\put(0.4,1.61){\color{SpringGreen}{\em 6}}
\put(4.505,3.41){\color{SpringGreen}\line(0,1){0.68}}
\put(4.4,3.61){\color{SpringGreen}{\em 6}}
\put(3.505,2.41){\color{SpringGreen}\line(0,1){0.68}}
\put(3.4,2.61){\color{SpringGreen}{\em 6}}
\put(2.505,1.41){\color{SpringGreen}\line(0,1){0.68}}
\put(2.4,1.61){\color{SpringGreen}{\em 6}}
\put(1.505,0.41){\color{SpringGreen}\line(0,1){0.68}}
\put(1.4,0.61){\color{SpringGreen}{\em 6}}
\end{picture}}

\begin{figure}
\begin{center}
{\bf \RootsFigure.1}\ \ Depictions of positive roots for $\myA_{6}$ and $\myD_{6}$\\
{\small (The notation `$a_{1}a_{2} \cdots a_{n}$' beside a given vertex means $\sum_{i \in I}a_{i}\alpha_{i}$.)}

\vspace*{0.2in} 
\ASixFlower

\vspace*{0.5in} 
\DSixFlower

\setlength{\unitlength}{1cm}
\begin{picture}(14,16)
\put(6,8.5){\ASixRoots}
\put(-0.5,1){\DSixRoots}
\end{picture}
\end{center}
\end{figure}

\newcommand{\BFiveFlower}{
\setlength{\unitlength}{0.5cm}
\begin{picture}(15,1)
\put(15.3,-0.25){\LARGE $\myB_{5}$}
\put(15.1,0){\vector(-1,0){1.2}}
{\color{Cyan}\put(1,0){\circle*{0.7}}}
\put(1,0){\circle{0.7}}
{\color{Red}\put(4,0){\circle*{0.7}}}
\put(4,0){\circle{0.7}}
{\color{Purple}\put(7,0){\circle*{0.7}}}
\put(7,0){\circle{0.7}}
{\color{OliveGreen}\put(10,0){\circle*{0.7}}}
\put(10,0){\circle{0.7}}
{\color{BurntOrange}\put(13,0){\circle*{0.7}}}
\put(13,0){\circle{0.7}}
\thicklines
\put(11,0){\vector(1,0){0.2}}
\put(11.35,0){\vector(1,0){0.2}}
\put(12,0){\vector(-1,0){0.2}}
\put(1.35,0){\line(1,0){2.3}}
\put(4.35,0){\line(1,0){2.3}}
\put(7.35,0){\line(1,0){2.3}}
\put(10.35,0){\line(1,0){2.3}}
\put(0.8,0.6){\textcolor{Cyan}{\em 1}}
\put(3.8,0.6){\textcolor{Red}{\em 2}}
\put(6.8,0.6){\textcolor{Purple}{\em 3}}
\put(9.8,0.6){\textcolor{OliveGreen}{\em 4}}
\put(12.8,0.6){\textcolor{BurntOrange}{\em 5}}
\end{picture}
}

\newcommand{\CFiveFlower}{
\setlength{\unitlength}{0.5cm}
\begin{picture}(15,1)
\put(-2.5,-0.25){\LARGE $\myC_{5}$}
\put(-0.9,0){\vector(1,0){1.2}}
{\color{Cyan}\put(1,0){\circle*{0.7}}}
\put(1,0){\circle{0.7}}
{\color{Red}\put(4,0){\circle*{0.7}}}
\put(4,0){\circle{0.7}}
{\color{Purple}\put(7,0){\circle*{0.7}}}
\put(7,0){\circle{0.7}}
{\color{OliveGreen}\put(10,0){\circle*{0.7}}}
\put(10,0){\circle{0.7}}
{\color{BurntOrange}\put(13,0){\circle*{0.7}}}
\put(13,0){\circle{0.7}}
\thicklines
\put(11,0){\vector(1,0){0.2}}
\put(11.65,0){\vector(-1,0){0.2}}
\put(12,0){\vector(-1,0){0.2}}
\put(1.35,0){\line(1,0){2.3}}
\put(4.35,0){\line(1,0){2.3}}
\put(7.35,0){\line(1,0){2.3}}
\put(10.35,0){\line(1,0){2.3}}
\put(0.8,0.6){\textcolor{Cyan}{\em 1}}
\put(3.8,0.6){\textcolor{Red}{\em 2}}
\put(6.8,0.6){\textcolor{Purple}{\em 3}}
\put(9.8,0.6){\textcolor{OliveGreen}{\em 4}}
\put(12.8,0.6){\textcolor{BurntOrange}{\em 5}}
\end{picture}
}

\newcommand{\BFiveRoots}{
\setlength{\unitlength}{1cm}
\begin{picture}(7,8)
\put(1,8.15){\large $\myB_{5}$ positive roots}
\put(-1,0){\TypeEboxDot{Gray}}
\put(-0.9,-0.175){\tiny 00001}
\put(1,0){\TypeEboxDot{Gray}}
\put(1.1,-0.175){\tiny 00010}
\put(3,0){\TypeEboxDot{Gray}}
\put(3.1,-0.175){\tiny 00100}
\put(5,0){\TypeEboxDot{Gray}}
\put(5.1,-0.175){\tiny 01000}
\put(7,0){\TypeEboxDot{Gray}}
\put(7.1,-0.175){\tiny 10000}
\put(0,1){\TypeEboxDot{Gray}}
\put(-0.5,1.175){\tiny 00011}
\put(2,1){\TypeEboxDot{Gray}}
\put(2.75,1.175){\tiny 00110}
\put(4,1){\TypeEboxDot{Gray}}
\put(4.75,1.175){\tiny 01100}
\put(6,1){\TypeEboxDot{Gray}}
\put(6.75,1.175){\tiny 11000}
\put(-1,2){\TypeEboxDot{Gray}}
\put(-0.25,2.175){\tiny 00012}
\put(1,2){\TypeEboxDot{Gray}}
\put(1.75,2.175){\tiny 00111}
\put(3,2){\TypeEboxDot{Gray}}
\put(3.75,2.175){\tiny 01110}
\put(5,2){\TypeEboxDot{Gray}}
\put(5.75,2.175){\tiny 11100}
\put(0,3){\TypeEboxDot{Gray}}
\put(-0.5,3.175){\tiny 00112}
\put(2,3){\TypeEboxDot{Gray}}
\put(2.75,3.175){\tiny 01111}
\put(4,3){\TypeEboxDot{Gray}}
\put(4.75,3.175){\tiny 11110}
\put(-1,4){\TypeEboxDot{Gray}}
\put(-0.25,4.175){\tiny 00122}
\put(1,4){\TypeEboxDot{Gray}}
\put(1.75,4.175){\tiny 01112}
\put(3,4){\TypeEboxDot{Gray}}
\put(3.75,4.175){\tiny 11111}
\put(0,5){\TypeEboxDot{Gray}}
\put(-0.5,5.175){\tiny 01122}
\put(2,5){\TypeEboxDot{Gray}}
\put(2.75,5.175){\tiny 11112}
\put(-1,6){\TypeEboxDot{Gray}}
\put(-0.25,6.175){\tiny 01222}
\put(1,6){\TypeEboxDot{Gray}}
\put(1.75,6.175){\tiny 11122}
\put(0,7){\TypeEboxDot{Gray}}
\put(0.75,7.175){\tiny 11222}
\put(-1,8){\TypeEboxDot{Gray}}
\put(-0.25,8.175){\tiny 12222}
\thicklines
\put(5.625,0.375){\color{Cyan}\qbezier(0,0)(0.375,0.375)(0.75,0.75)}
\put(4.625,1.375){\color{Cyan}\qbezier(0,0)(0.375,0.375)(0.75,0.75)}
\put(3.625,2.375){\color{Cyan}\qbezier(0,0)(0.375,0.375)(0.75,0.75)}
\put(2.625,3.375){\color{Cyan}\qbezier(0,0)(0.375,0.375)(0.75,0.75)}
\put(1.625,4.375){\color{Cyan}\qbezier(0,0)(0.375,0.375)(0.75,0.75)}
\put(0.625,5.375){\color{Cyan}\qbezier(0,0)(0.375,0.375)(0.75,0.75)}
\put(-0.375,6.375){\color{Cyan}\qbezier(0,0)(0.375,0.375)(0.75,0.75)}
\put(3.625,0.375){\color{Red}\qbezier(0,0)(0.375,0.375)(0.75,0.75)}
\put(2.625,1.375){\color{Red}\qbezier(0,0)(0.375,0.375)(0.75,0.75)}
\put(1.625,2.375){\color{Red}\qbezier(0,0)(0.375,0.375)(0.75,0.75)}
\put(0.625,3.375){\color{Red}\qbezier(0,0)(0.375,0.375)(0.75,0.75)}
\put(-0.375,4.375){\color{Red}\qbezier(0,0)(0.375,0.375)(0.75,0.75)}
\put(1.625,0.375){\color{Purple}\qbezier(0,0)(0.375,0.375)(0.75,0.75)}
\put(0.625,1.375){\color{Purple}\qbezier(0,0)(0.375,0.375)(0.75,0.75)}
\put(-0.375,2.375){\color{Purple}\qbezier(0,0)(0.375,0.375)(0.75,0.75)}
\put(-0.375,0.375){\color{OliveGreen}\qbezier(0,0)(0.375,0.375)(0.75,0.75)}
\put(0.39,7.375){\color{Red}\qbezier(0,0)(-0.375,0.375)(-0.75,0.75)}
\put(1.39,6.375){\color{Purple}\qbezier(0,0)(-0.375,0.375)(-0.75,0.75)}
\put(2.39,5.375){\color{OliveGreen}\qbezier(0,0)(-0.375,0.375)(-0.75,0.75)}
\put(3.39,4.375){\color{BurntOrange}\qbezier(0,0)(-0.375,0.375)(-0.75,0.75)}
\put(4.39,3.375){\color{BurntOrange}\qbezier(0,0)(-0.375,0.375)(-0.75,0.75)}
\put(5.39,2.375){\color{OliveGreen}\qbezier(0,0)(-0.375,0.375)(-0.75,0.75)}
\put(6.39,1.375){\color{Purple}\qbezier(0,0)(-0.375,0.375)(-0.75,0.75)}
\put(7.39,0.375){\color{Red}\qbezier(0,0)(-0.375,0.375)(-0.75,0.75)}
\put(0.39,5.375){\color{Purple}\qbezier(0,0)(-0.375,0.375)(-0.75,0.75)}
\put(1.39,4.375){\color{OliveGreen}\qbezier(0,0)(-0.375,0.375)(-0.75,0.75)}
\put(2.39,3.375){\color{BurntOrange}\qbezier(0,0)(-0.375,0.375)(-0.75,0.75)}
\put(3.39,2.375){\color{BurntOrange}\qbezier(0,0)(-0.375,0.375)(-0.75,0.75)}
\put(4.39,1.375){\color{OliveGreen}\qbezier(0,0)(-0.375,0.375)(-0.75,0.75)}
\put(5.39,0.375){\color{Purple}\qbezier(0,0)(-0.375,0.375)(-0.75,0.75)}
\put(0.39,3.375){\color{OliveGreen}\qbezier(0,0)(-0.375,0.375)(-0.75,0.75)}
\put(1.39,2.375){\color{BurntOrange}\qbezier(0,0)(-0.375,0.375)(-0.75,0.75)}
\put(2.39,1.375){\color{BurntOrange}\qbezier(0,0)(-0.375,0.375)(-0.75,0.75)}
\put(3.39,0.375){\color{OliveGreen}\qbezier(0,0)(-0.375,0.375)(-0.75,0.75)}
\put(0.39,1.375){\color{BurntOrange}\qbezier(0,0)(-0.375,0.375)(-0.75,0.75)}
\put(1.39,0.375){\color{BurntOrange}\qbezier(0,0)(-0.375,0.375)(-0.75,0.75)}
\put(-0.1,0.65){\color{OliveGreen}{\em 4}}
\put(-0.1,2.65){\color{Purple}{\em 3}}
\put(0.9,1.65){\color{Purple}{\em 3}}
\put(1.9,0.65){\color{Purple}{\em 3}}
\put(-0.1,4.65){\color{Red}{\em 2}}
\put(0.9,3.65){\color{Red}{\em 2}}
\put(1.9,2.65){\color{Red}{\em 2}}
\put(2.9,1.65){\color{Red}{\em 2}}
\put(3.9,0.65){\color{Red}{\em 2}}
\put(-0.1,6.65){\color{Cyan}{\em 1}}
\put(0.9,5.65){\color{Cyan}{\em 1}}
\put(1.9,4.65){\color{Cyan}{\em 1}}
\put(2.9,3.65){\color{Cyan}{\em 1}}
\put(3.9,2.65){\color{Cyan}{\em 1}}
\put(4.9,1.65){\color{Cyan}{\em 1}}
\put(5.9,0.65){\color{Cyan}{\em 1}}
\put(-0.1,7.65){\color{Red}{\em 2}}
\put(6.9,0.65){\color{Red}{\em 2}}
\put(-0.1,5.65){\color{Purple}{\em 3}}
\put(0.9,6.65){\color{Purple}{\em 3}}
\put(4.9,0.65){\color{Purple}{\em 3}}
\put(5.9,1.65){\color{Purple}{\em 3}}
\put(-0.1,3.65){\color{OliveGreen}{\em 4}}
\put(0.9,4.65){\color{OliveGreen}{\em 4}}
\put(1.9,5.65){\color{OliveGreen}{\em 4}}
\put(2.9,0.65){\color{OliveGreen}{\em 4}}
\put(3.9,1.65){\color{OliveGreen}{\em 4}}
\put(4.9,2.65){\color{OliveGreen}{\em 4}}
\put(3.9,3.65){\color{BurntOrange}{\em 5}}
\put(2.9,2.65){\color{BurntOrange}{\em 5}}
\put(1.9,1.65){\color{BurntOrange}{\em 5}}
\put(0.9,0.65){\color{BurntOrange}{\em 5}}
\put(2.9,4.65){\color{BurntOrange}{\em 5}}
\put(1.9,3.65){\color{BurntOrange}{\em 5}}
\put(0.9,2.65){\color{BurntOrange}{\em 5}}
\put(-0.1,1.65){\color{BurntOrange}{\em 5}}
\end{picture}}

\newcommand{\CFiveRoots}{
\setlength{\unitlength}{1cm}
\begin{picture}(7,8)
\put(1,8.15){\large $\myC_{5}$ positive roots}
\put(-1,0){\TypeEboxDot{Gray}}
\put(-0.9,-0.175){\tiny 00001}
\put(1,0){\TypeEboxDot{Gray}}
\put(1.1,-0.175){\tiny 00010}
\put(3,0){\TypeEboxDot{Gray}}
\put(3.1,-0.175){\tiny 00100}
\put(5,0){\TypeEboxDot{Gray}}
\put(5.1,-0.175){\tiny 01000}
\put(7,0){\TypeEboxDot{Gray}}
\put(7.1,-0.175){\tiny 10000}
\put(0,1){\TypeEboxDot{Gray}}
\put(-0.5,1.175){\tiny 00011}
\put(2,1){\TypeEboxDot{Gray}}
\put(2.75,1.175){\tiny 00110}
\put(4,1){\TypeEboxDot{Gray}}
\put(4.75,1.175){\tiny 01100}
\put(6,1){\TypeEboxDot{Gray}}
\put(6.75,1.175){\tiny 11000}
\put(-1,2){\TypeEboxDot{Gray}}
\put(-0.25,2.175){\tiny 00021}
\put(1,2){\TypeEboxDot{Gray}}
\put(1.75,2.175){\tiny 00111}
\put(3,2){\TypeEboxDot{Gray}}
\put(3.75,2.175){\tiny 01110}
\put(5,2){\TypeEboxDot{Gray}}
\put(5.75,2.175){\tiny 11100}
\put(0,3){\TypeEboxDot{Gray}}
\put(-0.5,3.175){\tiny 00121}
\put(2,3){\TypeEboxDot{Gray}}
\put(2.75,3.175){\tiny 01111}
\put(4,3){\TypeEboxDot{Gray}}
\put(4.75,3.175){\tiny 11110}
\put(-1,4){\TypeEboxDot{Gray}}
\put(-0.25,4.175){\tiny 00221}
\put(1,4){\TypeEboxDot{Gray}}
\put(1.75,4.175){\tiny 01121}
\put(3,4){\TypeEboxDot{Gray}}
\put(3.75,4.175){\tiny 11111}
\put(0,5){\TypeEboxDot{Gray}}
\put(-0.5,5.175){\tiny 01221}
\put(2,5){\TypeEboxDot{Gray}}
\put(2.75,5.175){\tiny 11121}
\put(-1,6){\TypeEboxDot{Gray}}
\put(-0.25,6.175){\tiny 02221}
\put(1,6){\TypeEboxDot{Gray}}
\put(1.75,6.175){\tiny 11221}
\put(0,7){\TypeEboxDot{Gray}}
\put(0.75,7.175){\tiny 12221}
\put(-1,8){\TypeEboxDot{Gray}}
\put(-0.25,8.175){\tiny 22221}
\thicklines
\put(5.625,0.375){\color{Cyan}\qbezier(0,0)(0.375,0.375)(0.75,0.75)}
\put(4.625,1.375){\color{Cyan}\qbezier(0,0)(0.375,0.375)(0.75,0.75)}
\put(3.625,2.375){\color{Cyan}\qbezier(0,0)(0.375,0.375)(0.75,0.75)}
\put(2.625,3.375){\color{Cyan}\qbezier(0,0)(0.375,0.375)(0.75,0.75)}
\put(1.625,4.375){\color{Cyan}\qbezier(0,0)(0.375,0.375)(0.75,0.75)}
\put(0.625,5.375){\color{Cyan}\qbezier(0,0)(0.375,0.375)(0.75,0.75)}
\put(-0.375,6.375){\color{Cyan}\qbezier(0,0)(0.375,0.375)(0.75,0.75)}
\put(3.625,0.375){\color{Red}\qbezier(0,0)(0.375,0.375)(0.75,0.75)}
\put(2.625,1.375){\color{Red}\qbezier(0,0)(0.375,0.375)(0.75,0.75)}
\put(1.625,2.375){\color{Red}\qbezier(0,0)(0.375,0.375)(0.75,0.75)}
\put(0.625,3.375){\color{Red}\qbezier(0,0)(0.375,0.375)(0.75,0.75)}
\put(-0.375,4.375){\color{Red}\qbezier(0,0)(0.375,0.375)(0.75,0.75)}
\put(1.625,0.375){\color{Purple}\qbezier(0,0)(0.375,0.375)(0.75,0.75)}
\put(0.625,1.375){\color{Purple}\qbezier(0,0)(0.375,0.375)(0.75,0.75)}
\put(-0.375,2.375){\color{Purple}\qbezier(0,0)(0.375,0.375)(0.75,0.75)}
\put(-0.375,0.375){\color{OliveGreen}\qbezier(0,0)(0.375,0.375)(0.75,0.75)}
\put(0.39,7.375){\color{Cyan}\qbezier(0,0)(-0.375,0.375)(-0.75,0.75)}
\put(1.39,6.375){\color{Red}\qbezier(0,0)(-0.375,0.375)(-0.75,0.75)}
\put(2.39,5.375){\color{Purple}\qbezier(0,0)(-0.375,0.375)(-0.75,0.75)}
\put(3.39,4.375){\color{OliveGreen}\qbezier(0,0)(-0.375,0.375)(-0.75,0.75)}
\put(4.39,3.375){\color{BurntOrange}\qbezier(0,0)(-0.375,0.375)(-0.75,0.75)}
\put(5.39,2.375){\color{OliveGreen}\qbezier(0,0)(-0.375,0.375)(-0.75,0.75)}
\put(6.39,1.375){\color{Purple}\qbezier(0,0)(-0.375,0.375)(-0.75,0.75)}
\put(7.39,0.375){\color{Red}\qbezier(0,0)(-0.375,0.375)(-0.75,0.75)}
\put(0.39,5.375){\color{Red}\qbezier(0,0)(-0.375,0.375)(-0.75,0.75)}
\put(1.39,4.375){\color{Purple}\qbezier(0,0)(-0.375,0.375)(-0.75,0.75)}
\put(2.39,3.375){\color{OliveGreen}\qbezier(0,0)(-0.375,0.375)(-0.75,0.75)}
\put(3.39,2.375){\color{BurntOrange}\qbezier(0,0)(-0.375,0.375)(-0.75,0.75)}
\put(4.39,1.375){\color{OliveGreen}\qbezier(0,0)(-0.375,0.375)(-0.75,0.75)}
\put(5.39,0.375){\color{Purple}\qbezier(0,0)(-0.375,0.375)(-0.75,0.75)}
\put(0.39,3.375){\color{Purple}\qbezier(0,0)(-0.375,0.375)(-0.75,0.75)}
\put(1.39,2.375){\color{OliveGreen}\qbezier(0,0)(-0.375,0.375)(-0.75,0.75)}
\put(2.39,1.375){\color{BurntOrange}\qbezier(0,0)(-0.375,0.375)(-0.75,0.75)}
\put(3.39,0.375){\color{OliveGreen}\qbezier(0,0)(-0.375,0.375)(-0.75,0.75)}
\put(0.39,1.375){\color{OliveGreen}\qbezier(0,0)(-0.375,0.375)(-0.75,0.75)}
\put(1.39,0.375){\color{BurntOrange}\qbezier(0,0)(-0.375,0.375)(-0.75,0.75)}
\put(-0.1,0.65){\color{OliveGreen}{\em 4}}
\put(-0.1,2.65){\color{Purple}{\em 3}}
\put(0.9,1.65){\color{Purple}{\em 3}}
\put(1.9,0.65){\color{Purple}{\em 3}}
\put(-0.1,4.65){\color{Red}{\em 2}}
\put(0.9,3.65){\color{Red}{\em 2}}
\put(1.9,2.65){\color{Red}{\em 2}}
\put(2.9,1.65){\color{Red}{\em 2}}
\put(3.9,0.65){\color{Red}{\em 2}}
\put(-0.1,6.65){\color{Cyan}{\em 1}}
\put(0.9,5.65){\color{Cyan}{\em 1}}
\put(1.9,4.65){\color{Cyan}{\em 1}}
\put(2.9,3.65){\color{Cyan}{\em 1}}
\put(3.9,2.65){\color{Cyan}{\em 1}}
\put(4.9,1.65){\color{Cyan}{\em 1}}
\put(5.9,0.65){\color{Cyan}{\em 1}}
\put(-0.1,7.65){\color{Cyan}{\em 1}}
\put(6.9,0.65){\color{Red}{\em 2}}
\put(-0.1,5.65){\color{Red}{\em 2}}
\put(0.9,6.65){\color{Red}{\em 2}}
\put(4.9,0.65){\color{Purple}{\em 3}}
\put(5.9,1.65){\color{Purple}{\em 3}}
\put(-0.1,3.65){\color{Purple}{\em 3}}
\put(0.9,4.65){\color{Purple}{\em 3}}
\put(1.9,5.65){\color{Purple}{\em 3}}
\put(2.9,0.65){\color{OliveGreen}{\em 4}}
\put(3.9,1.65){\color{OliveGreen}{\em 4}}
\put(4.9,2.65){\color{OliveGreen}{\em 4}}
\put(3.9,3.65){\color{BurntOrange}{\em 5}}
\put(2.9,2.65){\color{BurntOrange}{\em 5}}
\put(1.9,1.65){\color{BurntOrange}{\em 5}}
\put(0.9,0.65){\color{BurntOrange}{\em 5}}
\put(2.9,4.65){\color{OliveGreen}{\em 4}}
\put(1.9,3.65){\color{OliveGreen}{\em 4}}
\put(0.9,2.65){\color{OliveGreen}{\em 4}}
\put(-0.1,1.65){\color{OliveGreen}{\em 4}}
\end{picture}}

\begin{figure}
\begin{center}
{\bf \RootsFigure.2}\ \ Depictions of positive roots for $\myB_{5}$ and $\myC_{5}$\\
{\small (The notation `$a_{1}a_{2} \cdots a_{n}$' beside a given vertex means $\sum_{i \in I}a_{i}\alpha_{i}$.)}

\vspace*{0.2in} 
\BFiveFlower

\vspace*{0.5in} 
\CFiveFlower

\setlength{\unitlength}{1cm}
\begin{picture}(14,16)
\put(6.75,6.5){\BFiveRoots}
\put(0.25,0){\CFiveRoots}
\end{picture}
\end{center}
\end{figure}

\clearpage
\noindent 
supporting graph for the appropriate adjoint representation as in \cite{DonAdjoint}; for examples, see \RootsFigure.  
For $\myA_{n}$, these are seen to be $\displaystyle \sum_{k \in [i,j]_{\mathbb{Z}}}\alpha_{k}^{\vee}$ for any integer pair $(i,j)$ with $1 \leq i \leq j \leq n$. 
For $\myB_{n}$, whose co-roots are from $\myC_{n}$, these are seen to be either $\displaystyle \sum_{k \in [i,j]_{\mathbb{Z}}}\alpha_{k}^{\vee}$ for any integer pair $(i,j)$ with $1 \leq i \leq j \leq n-1$ or\ \ $\displaystyle \sum_{k \in [p,n]_{\mathbb{Z}}}\alpha_{k}^{\vee} + \sum_{k \in [q,n-1]_{\mathbb{Z}}}\alpha_{k}^{\vee}$\hfill for any integer pair $(p,q)$ with $1 \leq p \leq n$ and $p \leq q \leq n$. 
For $\myC_{n}$, whose co-roots are from $\myB_{n}$, these are seen to be either $\displaystyle \sum_{k \in [i,j]_{\mathbb{Z}}}\alpha_{k}^{\vee}$ for any integer pair $(i,j)$ with 
$1 \leq i \leq j \leq n$ or $\displaystyle \sum_{k \in [p,n]_{\mathbb{Z}}}\alpha_{k}^{\vee} + \sum_{k \in [q,n]_{\mathbb{Z}}}\alpha_{k}^{\vee}$ for any integer pair $(p,q)$ with $1 \leq p \leq n-1$ and $p+1 \leq q \leq n$. 
For $\myD_{n}$, these are seen to be either $\displaystyle \sum_{k \in [i,j]_{\mathbb{Z}}}\alpha_{k}^{\vee}$ for any integer pair $(i,j)$ with $1 \leq i \leq j \leq n-1$ or $\displaystyle \sum_{k \in [p,n-2]_{\mathbb{Z}}}\alpha_{k}^{\vee} + \sum_{k \in [q,n]_{\mathbb{Z}}}\alpha_{k}^{\vee}$ for any integer pair $(p,q)$ with $1 \leq p \leq n-1$ and $p+1 \leq q \leq n$.\hfill\QED

\vspace*{0.5cm}
\begin{center}
\fbox{\Large \bf Part IV}\\ 
\underline{\large \bf Examples and applications}
\end{center}
In the next several sections, we put the preceding ideas to work. 
Specifically, in \S \QuasiExampleSection, for each integral Coxeter--Dynkin flower $\mathscr{G}$ we construct a family of DCML's that are realized as splitting posets and supporting graphs associated with the `quasi-miniscule' or `short adjoint' $\mathfrak{g}(\mathscr{G})$-module$\,$/$\,$$\mathcal{W}(\mathscr{G})$-symmetric function. 
These DCML's can be easily described uniformly (i.e.\ across `$\myA$--$\myG$' type) and possess many beautiful properties. 
Our work in \S \MinusculeExampleSection\ with `minuscule' $\mathfrak{g}(\mathscr{G})$-modules$\,$/$\,$$\mathcal{W}(\mathscr{G})$-symmetric functions follows a similar pattern, but the bulk of our effort is towards \underline{uniform} demonstrations that, in each minuscule case, (a) the unique splitting poset/supporting graph is a DCDL and that (b) its vertex-colored compression poset is itself a distributive lattice that can be realized as a full-length sublattice of a product of two chains. 
In \S \MinusculeLatticePosetSection, we use this understanding of minuscule compression posets together with theory from some prior sections to demonstrate sturdiness of DCDL's obtained from certain skew-stacks of minuscule compression posets.

\vspace*{0.5cm}
\noindent 
{\bf \S \QuasiExampleSection. Examples of DCML's that are splitting posets: The quasi-minuscule case.} 
In \S \FlowerSection, we exhibited some of the combinatorial virtues of a splitting poset for a Weyl bialternant. 
Our production of examples here and in the next section will utilize or refer to most of the structures developed in \S \FlowerSection, as those are intrinsic to the general landscape of a given Coxeter--Dynkin flower or posy. 
For this section, we assume $\mathscr{G} = (\Gamma_{I},M_{I \times I})$ is a Coxeter--Dynkin flower, and we take $n = |I|$. 
Our work in this section with the irreducible $\mathfrak{g}(\mathscr{G})$-module and $\mathcal{W}(\mathscr{G})$-Weyl bialternant associated with the so-called quasi-minuscule dominant weight -- which is known to be the highest short root in the root system $\Phi(\mathscr{G})$ -- is patterned after our work in \cite{DonAdjoint} with the adjoint representation of $\mathfrak{g}(\mathscr{G})$ -- whose associated dominant weight is the highest root. 

{\bf [\S \QuasiExampleSection.1:\! Root lengths.]} 
By inspection of all Coxeter--Dynkin flowers, we see that there are at most two root lengths amongst the simple roots.  
Two different simple root lengths only occur in the type $\myB$, $\myC$, $\myF$, and $\myG$ cases. 
Then, $\mathcal{W}$-invariance of the inner product implies that there are at most two root lengths amongst all roots. 
We use the adjectives {\em short} and {\em long} generally to describe the lengths of roots when there are exactly two root lengths available ($\myB$, $\myC$, $\myF$, and $\myG$). 
In the $\myA$--$\myD$--$\myE$ cases, we consider all roots to be short. 
We observe that $\Phi^{\vee}$ has exactly two root lengths available if and only if $\Phi$ does as well and that $\alpha^{\vee} \in \Phi^{\vee}$ is long if and only if $\alpha \in \Phi$ is short. 
Let $\Phi_{\mbox{\tiny short}}$ (respectively, $\Phi_{\mbox{\tiny long}}$) be the set of all short (resp.\ long) roots in $\Phi$, and use analogous notation for the co-roots $\Phi^{\vee}$. 
Let $\#_{\mbox{\tiny short}}(\mathscr{G})$ denote the number of short simple roots in $\Phi$ and $\#_{\mbox{\tiny long}}(\mathscr{G}) = n-\#_{\mbox{\tiny short}}(\mathscr{G})$ denote the number of long simple roots. 

{\bf [\S \QuasiExampleSection.2:\! Highest short root.]} 
According to \MinimalWeightsLemma, our partially ordered weight lattice $\Lambda$ has a finite number of connected components, one of which contains the zero weight.  
Give the induced order from $\Lambda$ to the set of nonzero dominant weights within the connected component of the zero weight. 
Amongst these weights, there is a unique minimal nonzero dominant weight, which we hereby denote $\varpi$. 
It follows from the definitions that $\Pi(\varpi) = \mathcal{W}\varpi \disjointunion \{\zeroweight\}$ is a set equality. 
Several equivalent characterizations of $\varpi$ are given in Proposition 4.14 of \cite{DonPosetModels} and are largely drawn from \cite{StemDom}. 
In the next lemma we re-iterate these and add a new characterization. 

\noindent 
{\bf \HighestShawtyLemma}\ \ {\sl Let $\lambda$ be a nonzero dominant weight.  The following are equivalent:}\\ 
{\sl (1) $\Pi(\lambda) = \mathcal{W}\lambda \disjointunion \{\zeroweight\}$ is an equality of sets.}\\
{\sl (2) For all $\nu \in \Lambda^{+}$, if $\nu \in \Pi(\lambda)$ then $\nu=\lambda$ or $\nu=\zeroweight$.}\\
{\sl (3) We have $\langle \lambda,\alpha^{\vee} \rangle \in \{0,\pm 1, \pm 2\}$ for all roots $\alpha$ and there is a unique $\beta \in \Phi$ for which $\langle \lambda,\beta^{\vee} \rangle = 2$.}\\
{\sl (4) $\lambda = \varpi$, in which case the orbit $\mathcal{W}\varpi$ is the set} $\Phi_{\mbox{\tiny short}}$.\\
{\sl (5) In any Networked-numbers Game played to completion from initial position $\lambda$, exactly once the number at a fired node is `$+2$' and otherwise this number is `$+1$'.}

{\em Proof.} 
See \cite{DonPosetModels} Proposition 4.14 for {\sl (1)} $\Leftrightarrow$ {\sl (2)} $\Leftrightarrow$ {\sl (3)} $\Leftrightarrow$ {\sl (4)}. 
So, to complete the proof it suffices to demonstrate {\sl (5)} $\Leftrightarrow$ {\sl (3)}. 
It will be useful to establish two general claims that depend upon the finiteness of $\mathcal{W}$, but first some general set-up. 
Suppose, in the language of \S \WeylSection.4, that $\nu$ is $K^{c}$-dominant for some subset $K \subseteq I$. 
In particular, the stabilizer of $\nu$ under the action of the Weyl group $\mathcal{W}$ is the parabolic subgroup $\mathcal{W}_{K}$, and any $\sigma \in \mathcal{W}$ can be written uniquely as $\sigma^{K}\sigma_{K}$ where $\sigma^{K} \in \mathcal{W}^{K}$ is the minimal length representative of the coset $\sigma\mathcal{W}_{K}$. 

\noindent 
\underline{Claim 1:} For any intermediate position $\mu$ in any Networked-numbers Game played from initial position $\nu$, an illegal firing move played from $\mu$ results in a position $\mu'$ that can be obtained from $\nu$ by some legal firing sequence. 
\underline{\em Proof.} Say $\mu = \mygens_{i_{p}}\cdots\mygens_{i_{1}}.\nu$ is the position resulting from the legal firing sequence $(\gamma_{i_{1}},\ldots,\gamma_{i_{p}})$. 
But suppose firing $\gamma_{i_{p+1}}$ from $\mu$ is illegal, resulting in position $\mu' = \mygens_{i_{p+1}}\mygens_{i_{p}}\cdots\mygens_{i_{1}}.\nu$. 
Then by Proposition 3.2 of \cite{DonEur}, the expression $\sigma := \mygens_{i_{p+1}}\mygens_{i_{p}}\cdots\mygens_{i_{1}}$ is not reduced. 
Then let $\mygens_{j_{q}}\cdots\mygens_{j_{1}}$ be a reduced expression for $\sigma^{K}$. 
Again by Proposition 3.2 of \cite{DonEur}, $(\gamma_{j_{i}},\ldots,\gamma_{j_{q}})$ is legal to play from $\nu$ and yields the position $\mu' = \mygens_{j_{q}}\cdots\mygens_{j_{1}}.\nu$. 

\noindent 
\underline{Claim 2:} For any game sequence $(\gamma_{i_{1}},\ldots,\gamma_{i_{k}})$ played from $\nu$, we have $\{\mygens_{i_{1}}\cdots\mygens_{i_{p}}.\alpha^{\vee}_{i_{p+1}}\}_{p \in [0,k-1]_{\mathbb{Z}}} = \{\eta^{\vee} \in (\Phi^{\vee})^{+}\, |\, \langle \nu,\eta^{\vee} \rangle > 0\}$. 
\underline{\em Proof.} 
Well, we know that $(w_{0})^{K} = \mygens_{i_{k}}\cdots\mygens_{i_{1}}$ and $k = \ell((w_{0})^{K})$ by Corollary 3.4 of \cite{DonEur}. 
Borrowing from the notation of that paper, let $N^{\vee}\big((w_{0})^{K}\big) := \{\eta^{\vee} \in (\Phi^{\vee})^{+}|(w_{0})^{K}.\eta^{\vee} \in (\Phi^{\vee})^{-}\}$. 
In this context, it is well known that $k = \myabs N^{\vee}\big((w_{0})^{K}\big)\myabs$, see e.g.\ Proposition 2.12 of \cite{DonEur}. 
Now, by Lemma 5.1 of that paper, $N^{\vee}\big((w_{0})^{K}\big) = \{\mygens_{i_{1}}\cdots\mygens_{i_{p}}.\alpha_{i_{p+1}}^{\vee}\, |\, p \in [0,k-1]_{\mathbb{Z}} \}$. 
If $\langle \nu,\eta^{\vee} \rangle > 0$ for some $\eta^{\vee} \in \Phi^{\vee}$, then necessarily $\eta^{\vee} \in (\Phi^{\vee})^{+}$, and the fact that $0 < \langle \nu,\eta^{\vee} \rangle = \langle (w_{0})^{K}.\nu,(w_{0})^{K}.\eta^{\vee} \rangle$ means that we must have $(w_{0})^{K}.\eta^{\vee} \in (\Phi^{-})^{\vee}$; hence $\eta^{\vee} \in N^{\vee}\big((w_{0})^{K}\big)$. 
On the other hand, if $\eta^{\vee} \in N^{\vee}\big((w_{0})^{K}\big)$, then $\eta^{\vee} = \mygens_{i_{1}}\cdots\mygens_{i_{p}}.\alpha_{i_{p+1}}^{\vee}$ for some $p \in [0,k-1]_{\mathbb{Z}}$, and moreover $\langle \nu,\eta^{\vee} \rangle = \langle \mygens_{i_{p}}\cdots\mygens_{i_{1}}.\nu,\alpha_{i_{p+1}}^{\vee} \rangle$ is positive since it represents the number at node $\gamma_{i_{p+1}}$ after the legal firing sequence $(\gamma_{i_{1}},\ldots,\gamma_{i_{p}})$ has been played. 
That is, $N^{\vee}\big((w_{0})^{K}\big) = \{\eta^{\vee} \in (\Phi^{\vee})^{+}\, |\, \langle \nu,\eta^{\vee} \rangle > 0\}$. 

Now we return our attention to the yet-to-be-proved equivalence {\sl (5)} $\Leftrightarrow$ {\sl (3)}. 
Assuming {\sl (3)}, we suppose that $(\gamma_{i_{1}},\ldots,\gamma_{i_{l}})$ is a game sequence from initial position $\lambda$.
Since $\langle \lambda,\beta^{\vee} \rangle = 2$, then, by Claim 2, there exists some $q \in [0,l-1]_{\mathbb{Z}}$ such that $\beta^{\vee} = \mygens_{i_{1}}\cdots\mygens_{i_{q}}.\alpha_{i_{q+1}}^{\vee}$. 
Then, for the intermediate position resulting from playing the legal firing sequence $(\gamma_{i_{1}},\ldots,\gamma_{i_{q}})$, the number at node $\gamma_{i_{q+1}}$ is $\langle \mygens_{i_{q}}\cdots\mygens_{i_{1}}.\lambda,\alpha_{i_{q+1}}^{\vee} \rangle = \langle \lambda,\mygens_{i_{1}}\cdots\mygens_{i_{q}}.\alpha_{i_{q+1}}^{\vee} \rangle = \langle \lambda,\beta^{\vee} \rangle = 2$.  
Now say $\alpha^{\vee} = \mygens_{i_{1}}\cdots\mygens_{i_{p}}.\alpha_{i_{p+1}}^{\vee}$ for some $p \in [0,l-1]_{\mathbb{Z}} \setminus \{q\}$. 
Then the number at node $\gamma_{i_{p+1}}$ after the legal firing sequence $(\gamma_{i_{1}},\ldots,\gamma_{i_{p}})$ has been played from $\nu$ is $\langle \mygens_{i_{p}}\cdots\mygens_{i_{1}}.\lambda,\alpha_{i_{p+1}}^{\vee} \rangle = \langle \lambda,\mygens_{i_{1}}\cdots\mygens_{i_{p}}.\alpha_{i_{p+1}}^{\vee} \rangle$. 
Since this latter quantity must be positive (because firing $\gamma_{i_{p+1}}$ at this point is legal), must be in $\{0,\pm1,\pm2\}$ (by the hypothesis of the statement of {\sl (3)}), and cannot be $2$ (by statement {\sl (3)} this can only be 2 if $\alpha^{\vee} = \beta^{\vee}$ and by Claim 2 the latter can be true only if $p=q$), then the number $\langle \mygens_{i_{p}}\cdots\mygens_{i_{1}}.\lambda,\alpha_{i_{p+1}}^{\vee} \rangle$ at node $\gamma_{i_{p+1}}$ when $\gamma_{i_{p+1}}$ is fired must be $1$.   
This establishes that {\sl (3)} implies {\sl (5)}. 

Conversely, suppose {\sl (5)} holds. 
Say $\lambda$ is $J^{c}$-dominant. 
Any given co-root $\alpha^{\vee}$ can be written as $\sigma^{-1}.\alpha_{j}^{\vee}$ for some simple co-root $\alpha_{j}^{\vee}$ and some $\sigma \in \mathcal{W}$. 
Now, $\langle \lambda,\alpha^{\vee} \rangle = \langle \lambda,\sigma^{-1}.\alpha_{j}^{\vee} \rangle = \langle \sigma.\lambda,\alpha_{j}^{\vee} \rangle$. 
Let $\mygens_{i_{p}}\cdots\mygens_{i_{1}}$ be a reduced expression for $\sigma^{J}$. 
By Proposition 3.2 of \cite{DonEur}, it follows that $(\gamma_{i_{1}},\ldots,\gamma_{i_{p}})$ is a legal firing sequence from initial position $\lambda$, resulting in an intermediate position $\mu$. 
The number at node $\gamma_{j}$ of $\mu$ is $\langle \mu,\alpha_{j}^{\vee} \rangle$. 
When $\langle \mu,\alpha_{j}^{\vee} \rangle \not= 0$, there are exactly two scenarios. 
($+$) If $\langle \mu,\alpha_{j}^{\vee} \rangle$ is positive, then let $i_{p+1} := j$ and extend the legal firing sequence $(\gamma_{i_{1}},\ldots,\gamma_{i_{p}},\gamma_{i_{p+1}})$ to a game sequence $(\gamma_{i_{1}},\ldots,\gamma_{i_{l}})$ simply by continuing game play. 
By the hypothesis expressed in statement {\sl (5)}, then, we must have $\langle \lambda,\alpha^{\vee} \rangle = \langle \mu,\alpha_{i_{p+1}}^{\vee} \rangle \in \{1,2\}$. 
($-$) On the other hand, if the number $\langle \mu,\alpha_{j}^{\vee} \rangle$ at node $\gamma_{j}$ of the intermediate position $\mu$ is negative, the firing node $\gamma_{j}$ illegally would give us a position $\mu'$ that could, by Claim 1, also be reached by a legal firing sequence. 
From $\mu'$, firing node $\gamma_{j}$ would be legal, since $-\langle \mu,\alpha_{j}^{\vee} \rangle$ would be positive. 
By our reasoning in the ($+$) case, $-\langle \mu,\alpha_{j}^{\vee} \rangle \in \{1,2\}$, so $\langle \lambda,\alpha^{\vee} \rangle \in \{-1,-2\}$. 
So, we have shown that $\langle \lambda,\alpha^{\vee} \rangle \in \{0,\pm 1,\pm 2\}$ for all $\alpha \in \Phi$. 

Now suppose $(\gamma_{i_{1}},\ldots,\gamma_{i_{l}})$ is any game sequence played from $\lambda$. 
Suppose $\langle \lambda,\beta^{\vee} \rangle = 2$; positivity of this quantity implies $\beta^{\vee} \in (\Phi^{\vee})^{+}$. 
So, $\beta^{\vee} \in \{\eta^{\vee} \in (\Phi^{\vee})^{+}\, |\, \langle \lambda,\eta^{\vee} \rangle > 0\}$, and therefore by Claim 2 $\beta^{\vee} \in \{\mygens_{i_{1}}\cdots\mygens_{i_{p}}.\alpha^{\vee}_{i_{p+1}}\}_{p \in [0,l-1]_{\mathbb{Z}}}$. 
By the hypothesis of {\sl (5)}, there is a unique $q \in [0,l-1]_{\mathbb{Z}}$ such that $\langle \lambda,\mygens_{i_{1}}\cdots\mygens_{i_{q}}.\alpha^{\vee}_{i_{q+1}} \rangle = 2$ while $\langle \lambda,\mygens_{i_{1}}\cdots\mygens_{i_{p}}.\alpha^{\vee}_{i_{p+1}} \rangle =1$ for $p \not= q$. 
So, $\beta^{\vee}$ must be the co-root $\mygens_{i_{1}}\cdots\mygens_{i_{q}}.\alpha^{\vee}_{i_{q+1}}$. 
This completes the proof that {\sl (5)} implies {\sl (3)}.\hfill\QED 

The nonzero dominant weight $\varpi$ is therefore the {\em highest short root} in the senses that $\alpha \leq \varpi$ and $\langle \alpha,\varrho^{\vee} \rangle \leq \langle \varpi,\varrho^{\vee} \rangle$ for all short roots $\alpha$. 
The highest short root is also sometimes called the {\em quasi-minuscule dominant weight} uniquely associated to any given Coxeter--Dynkin flower, and we call the associated $\mathcal{W}$-symmetric function the {\em quasi-minuscule Weyl bialternant}. 
Here are the highest short roots for the root system $\Phi(\mathscr{G})$ for each Coxeter--Dynkin flower, written as $(\mathscr{G},\varpi)$ pairs and with node numbering from \IEGGraphFigure: $(\myA_{n},\omega_{1}+\omega_{n})$, $(\myB_{n},\omega_{1})$, $(\myC_{n},\omega_{2})$, $(\myD_{n},\omega_{2})$, $(\myE_{6},\omega_{2})$, $(\myE_{7},\omega_{1})$, $(\myE_{8},\omega_{8})$, $(\myF_{4},\omega_{4})$, and $(\myG_{2},\omega_{1})$. 

{\bf [\S \QuasiExampleSection.3:\! DCML's associated with quasi-minuscule dominant weights.]}   
Our aim is to describe all of the diamond-colored modular lattice splitting posets (henceforth, `splitting DCML's') for the quasi-minuscule Weyl bialternant $\chi^{\mathscr{G}}_{_{\varpi}}$. 
For reference, some examples are depicted in \ShawtyFigureList. 
What follows is a variation on constructions from \cite{DonAdjoint} of splitting DCML's for the Weyl bialternant $\chi^{\mathscr{G}}_{_{\lambda}}$ when $\lambda$ is the highest long root of $\Phi$. 
Let $k$ be the index affiliated with some short simple root $\alpha_{k}$. 
We form an edge-colored directed graph $\myfancyQ^{(k)} = \myfancyQ^{(k)}(\varpi)$, called, for now, the $k^{\mbox{\tiny th}}$ ``quasi-minuscule poset'', as follows. 
As a set, $\myfancyQ^{(k)} = \Phi_{\mbox{\tiny short}} \disjointunion \{\lsem k,k \rsem\} \disjointunion \{\lsem i,j \rsem\, |\, \alpha_{i},\alpha_{j} \in \Phi_{\mbox{\tiny short}} \mbox{ with } M_{ij} < 0\}$, a disjoint union of three sets and where, in the latter set, we regard $\lsem i,j \rsem$ and $\lsem j,i \rsem$ to be equivalent. 
For $\xelt, \yelt \in \myfancyQ^{(k)}(\varpi)$, let $\xelt \myarrow{i} \yelt$ if and only if exactly one of the following is true: (a) $\xelt = \alpha = \Phi_{\mbox{\tiny short}}$ and $\yelt = \beta = \Phi_{\mbox{\tiny short}}$ with $\beta-\alpha=\alpha_{i}$; (b) $\xelt = \lsem i,j \rsem$ and $\yelt = \alpha_{i}$; or (c) $\xelt = -\alpha_{i}$ and $\yelt = \lsem i,j \rsem$. 

Since $\alpha \leq 0 \leq \beta$ for every $\alpha \in \Phi_{\mbox{\tiny short}}^{-}$ and $\beta \in \Phi_{\mbox{\tiny short}}^{+}$, it follows that $\myfancyQ^{(k)}(\varpi)$ is the edge-colored covering digraph for a connected poset with rank function $\rho: \myfancyQ^{(k)}(\varpi) \longrightarrow [0,2\langle \varpi,\varrho^{\vee} \rangle]_{\mathbb{Z}}$ given by $\rho(\xelt) = \langle \xelt+\varpi,\varrho^{\vee} \rangle$ when $\xelt \in \Phi_{\mbox{\tiny short}}$ and $\rho(\xelt) = \langle \varpi,\varrho^{\vee} \rangle$ when $\xelt = \lsem i,j \rsem$. 
It is not hard to see that, since $\Pi(\myfancyQ^{(k)}) = \Pi(\varpi)$ is $\mathscr{G}$-structured (cf.\ \SaturatedLemma), then $\myfancyQ^{(k)}$ is $\mathscr{G}$-structured as well. 
The set $\mathcal{Z}^{(k)} := \{\lsem k,k \rsem\} \disjointunion \{\lsem i,j \rsem\, |\, \alpha_{i},\alpha_{j} \in \Phi_{\mbox{\tiny short}} \mbox{ with } M_{ij} < 0\}$ is exactly comprised of the elements of $\myfancyQ^{(k)}$ with weight equal to $\zeroweight$. 

Before we establish that quasi-minuscule posets are $\mathscr{G}$-structured DCML's, we make one more observation about their structure. 
Now, in the sub-embryophyte of $\mathscr{H}$ of $\mathscr{G}$ whose nodes correspond to short simple roots, the vascular graph $\Gamma(\mathscr{H})$ is necessarily a tree which therefore has $\#_{\mbox{\tiny short}}(\mathscr{G})-1$ adjacencies. 
So, $\rule[-1.5mm]{0.2mm}{5mm}\{\lsem k,k \rsem\} \disjointunion \{\lsem i,j \rsem\, |\, \alpha_{i},\alpha_{j} \in \Phi_{\mbox{\tiny short}} \mbox{ with } M_{ij} < 0\}\rule[-1.5mm]{0.2mm}{5mm} = 1 + (\#_{\mbox{\tiny short}}(\mathscr{G})-1) = \#_{\mbox{\tiny short}}(\mathscr{G})$. 
It follows that $\CARD\left(\myfancyQ^{(k)}(\varpi)\right) = \rule[-1.5mm]{0.2mm}{5mm}\, \Phi_{\mbox{\tiny short}}\rule[-1.5mm]{0.2mm}{5mm} + \#_{\mbox{\tiny short}}(\mathscr{G})$. 

We record some of the preceding observations about $\myfancyQ^{(k)}$ in the next result. 
The main content of this result that has yet to be addressed is the DCML claim of the second sentence. 
Once this is established, we are justified in calling any $\myfancyQ^{(k)}$ a `quasi-minuscule DCML'.

\noindent 
{\bf \ShortAdjointDCML}\ \ {\sl Let $\mathscr{G}$ be a Coxeter--Dynkin flower. 
Then the} $k^{\mbox{\tiny th}}$ {\sl quasi-minuscule poset $\myfancyQ^{(k)} = \myfancyQ^{(k)}(\varpi)$ is a $\mathscr{G}$-structured and diamond-colored modular lattice. 
Moreover, $\Pi\left(\myfancyQ^{(k)}\right) = \Pi(\varpi)$; we have} $\rule[-1.5mm]{0.2mm}{5mm}\, \mathcal{Z}^{(k)}\rule[-1.5mm]{0.2mm}{5mm} = \#_{\mbox{\tiny short}}(\mathscr{G})${\sl ; we have $wt(\xelt)=\xelt$ whenever $\xelt$ is from the} `$\Phi_{\mbox{\tiny short}}$' part of $\myfancyQ^{(k)}${\sl ; we have} $wt(\lsem i,j \rsem)=\zeroweight$ {\sl for any $\lsem i,j \rsem$ from the `$\mathcal{Z}^{(k)}$' part of $\myfancyQ^{(k)}$; and the unique rank function on $\myfancyQ^{(k)}$ is $\rho: \myfancyQ^{(k)}(\varpi) \longrightarrow [0,2\langle \varpi,\varrho^{\vee} \rangle]_{\mathbb{Z}}$ given by $\rho(\xelt) = \langle wt(\xelt)+\varpi,\varrho^{\vee} \rangle$.}

{\em Proof.} The only part of the statement not addressed in paragraphs preceding the lemma statement is the claim the quasi-minuscule posets are DCML's. 
We prove this using case analysis.  
In types $\myA$--$\myD$--$\myE$, all roots have the same length and may be considered short.  
So, the associated quasi-minuscule pre-DCML's are exactly the extremal supporting graphs constructed in \cite{DonAdjoint} for the adjoint representations. 
In \cite{DonAdjoint}, it is shown that these are modular lattices. 
That there is only one possible splitting poset in each of the $\myB_{n}(\omega_{1})$ and $\myG_{2}(\omega_{1})$ cases is a consequence of Theorem 4.20 of \cite{DonPosetModels}. 
It is easy to see (e.g.\ \cite{DonSupp}) that the unique splitting poset in each case is a chain (and therefore a $\myB_{n}$- or $\myG_{2}$-structured DCML) and that it coincides with our quasi-minuscule pre-DCML. 
See \ShawtyFFourFigure\ for depictions of the two quasi-minuscule posets in the $\myF_{4}(\omega_{1})$ case.  

In fact both of these are DCDL's, as can be seen by constructing $\Jcolor\left(\jcolor(L)\right)$ for each $L$. 
These two splitting DCDL's are more thoroughly discussed in [Gilliland]. 
The $\myC_{n}(\omega_{2})$ case is studied in [DonSymplectic], where the quasi-minuscule posets $\myfancyQ^{(1)}$ and $\myfancyQ^{(n-1)}$ correspond respectively to the `KN' and `De Concini' symplectic lattices, which are both $\myC_{n}$-structured and splitting since both are shown in that paper to serve as supporting graphs for the associated fundamental representation of the symplectic Lie algebra $\mathfrak{g}(\myC_{n})$. 
Let $\widetilde{K_{\mytinyC_{n}}(\omega_{2})}$ be the $\myC_{n}$-structured DCDL obtained by removing the unique join irreducible element whose weight is $\zeroweight$ from either of $\myfancyQ^{(1)}$ or $\myfancyQ^{(n-1)}$. 
Then, for $1 < k < n-1$, the quasi-minuscule poset $\myfancyQ^{(k)}$ is obtained by modifying the unique color $k$ diamond \parbox[c]{1cm}{\begin{center}
\setlength{\unitlength}{0.2cm}
\begin{picture}(4.5,4.5)
\put(2.2,0){\circle*{0.6}} \put(0.2,3){\circle*{0.6}} 
\put(4.2,3){\circle*{0.6}} \put(2.2,6){\circle*{0.6}}
\put(0.2,3){\line(2,3){2}} \put(2.2,0){\line(-2,3){2}} \put(4.2,3){\line(-2,3){2}}
\put(2.2,0){\line(2,3){2}}
\put(0.95,1.5){\tiny $k$} \put(2.95,1.5){\tiny $k$}
\put(0.95,4.5){\tiny $k$} \put(2.95,4.5){\tiny $k$}
\end{picture} \vspace*{-0.3in} 
\end{center}} of $\widetilde{K_{\mytinyC_{n}}(\omega_{2})}$ to be \parbox[c]{1cm}{\begin{center}
\setlength{\unitlength}{0.2cm}
\begin{picture}(4.5,4.5)
\put(2.2,0){\circle*{0.6}} \put(0.2,3){\circle*{0.6}} \put(2.2,3){\circle*{0.6}}
\put(4.2,3){\circle*{0.6}} \put(2.2,6){\circle*{0.6}}
\put(0.2,3){\line(2,3){2}} \put(2.2,0){\line(-2,3){2}} \put(4.2,3){\line(-2,3){2}}
\put(2.2,0){\line(0,1){3}} \put(2.2,3){\line(0,1){3}} \put(2.2,0){\line(2,3){2}}
\put(0.95,1.5){\tiny $k$} \put(2.95,1.5){\tiny $k$}
\put(0.95,4.5){\tiny $k$} \put(2.95,4.5){\tiny $k$}
\put(1.95,1.5){\tiny $k$} \put(1.95,4.5){\tiny $k$}
\end{picture} \vspace*{-0.3in} 
\end{center}}. 
It is not hard to see that this poset, obtained from a $\myC_{n}$-structured DCDL, is a $\myC_{n}$-structured DCML. 
See \ShawtyCFourFigure\ for a depiction of the $n=4$ case. 
\hfill\QED 

{\bf [\S \QuasiExampleSection.4:\! Comments and questions about quasi-minuscule DCML's.]} 
We note that, by inspection of cases, a quasi-minuscule $\myfancyQ^{(k)}$ is distributive if and only if the the node $\gamma_{k}$ of $\mathscr{G}$ is adjacent to at most one other node $\gamma_{j}$ whose associated simple root $\alpha_{j}$ is short. 
So, the quasi-miniscule DCML's that are distributive are exactly these: $\myfancyQ^{(1)}$ and $\myfancyQ^{(n)}$ in the $\myA_{n}$ case; $\myfancyQ^{(1)}$ in the $\myB_{n}$ case; $\myfancyQ^{(1)}$ and $\myfancyQ^{(n-1)}$ in the $\myC_{n}$ case; $\myfancyQ^{(3)}$ and $\myfancyQ^{(4)}$ in the $\myF_{4}$ case; and $\myfancyQ^{(1)}$ in the $\myG_{2}$ case. 
That is, amongst the quasi-minuscule DCML's in the $\myD$ and $\myE$ cases, none are distributive. 

We have long believed that Stembridge's `Smash' Theorem \cite{StemQuasi} can be used to argue uniformly (i.e.\ without depending on the $\myA$--$\myG$ type of the given Coxeter--Dynkin flower) that our $\myfancyQ^{(k)}$'s are, indeed, modular lattices. 
We claimed as much in \cite{DonAdjoint}, but in a recent re-examination of that claim, we realized that our putative proof is missing some details. 
So, this is an open question of interest to us.  
We also speculate that there is a version of Stembridge's Smash Theorem for highest long roots that could lead to a uniform proof that the extremal supporting graphs of the adjoint representations found in \cite{DonAdjoint} are DCML's. 
To us, this is another interesting open question. 

Another interesting open question is to find, and uniformly describe, splitting DCDL's for the Weyl bialternants $\chi^{\mathscr{G}}_{_{k\varpi}}$ when $k \in [1,\infty)_{\mathbb{Z}}$. 
According to our comments in the first paragraph of this subsection, this is probably only plausible in the cases $\mathscr{G} \in \{\myA_{n},\myB_{n},\myC_{n},\myF_{4},\myG_{2}\}$. 
In the $\myA_{n}$ cases, two such examples are the `GT-left' and `GT-right' DCDL's described in \cite{DonSupp} for highest weight $k\varpi = k\omega_{1}+k\omega_{n}$. 
In the $\myB_{n}$ cases, two such examples are the `RS' and `Molev' DCDL's described in \cite{DLP1} for highest weight $k\varpi = k\omega_{1}$ 
In the $\myG_{2}$ cases, two such examples are the `RS' and `LM' DCDL's described in \cite{DLP1} and \cite{DLP2} for highest weight $k\varpi = k\omega_{1}$. 
In the $\myF_{4}$ case, the highest weight of interest is $k\varpi = k\omega_{4}$. 
A splitting DCDL for $\chi^{\mytinyF_{4}}_{_{2\omega_{4}}}$ was found in [Gilliland]; moreover, the quasi-minuscule splitting DCDL's depicted in \ShawtyFFourFigure\ (with compression posets also depicted there) are splitting DCDL's for $\chi^{\mytinyF_{4}}_{_{\omega_{4}}}$. 
In the $\myC_{n}$ cases, the highest weight of interest is $k\varpi = k\omega_{2}$. 
The quasi-minuscule splitting DCDL's depicted in \ShawtyCFourDistributiveFigure\ (with compression posets also depicted there) are splitting DCDL's for $\chi^{\mytinyC_{4}}_{_{\omega_{2}}}$.
In the $\myA_{n}$--$\myB_{n}$--$\myG_{2}$ cases, compression posets for some of our $\chi^{\mathscr{G}}_{_{k\varpi}}$-splitting DCDL's can be viewed nicely as a kind of stacking of $k$ copies a compression poset for a $\chi^{\mathscr{G}}_{_{\varpi}}$-splitting DCDL (which is necessarily one of our $\mathscr{G}$-quasi-minuscule splitting DCML's). 
We believe some similar viewpoint might be useful in discovering $\chi^{\mathscr{G}}_{_{k\varpi}}$-splitting DCDL's in the $\myC_{n}$--$\myF_{4}$ cases. 

If we relax the diamond-colored requirement for our study of splitting modular/distributive lattices for quasi-minuscule Weyl bialternants, then some other possibilities emerge. 
For example, the following $\myA_{2}$-structured distributive lattice is not diamond-colored but is a splitting poset for the quasi-minuscule Weyl bialternant $\chi_{_{\omega_{1}+\omega_{2}}}^{\mytinyA_{2}}$:
\parbox{2cm}{\begin{center}
\setlength{\unitlength}{0.2cm}
\begin{picture}(9,6.5)
\put(3,0){\circle*{0.5}} 
\put(1,2){\circle*{0.5}} \put(5,2){\circle*{0.5}}
\put(3,4){\circle*{0.5}} \put(7,4){\circle*{0.5}} 
\put(5,6){\circle*{0.5}} \put(9,6){\circle*{0.5}} 
\put(7,8){\circle*{0.5}}
\put(3,0){\line(-1,1){2}} \put(3,0){\line(1,1){2}}
\put(1,2){\line(1,1){2}} 
\put(5,2){\line(-1,1){2}} \put(5,2){\line(1,1){2}} 
\put(3,4){\line(1,1){2}} 
\put(7,4){\line(-1,1){2}} \put(7,4){\line(1,1){2}} 
\put(5,6){\line(1,1){2}} 
\put(9,6){\line(-1,1){2}}  
\put(1.7,0.65){\em \scriptsize 1} 
\put(3.7,2.65){\em \scriptsize 1}
\put(5.7,2.65){\em \scriptsize 1}
\put(7.7,4.65){\em \scriptsize 1}
\put(5.7,6.65){\em \scriptsize 1}
\put(3.7,0.7){\em \scriptsize 2}
\put(1.7,2.7){\em \scriptsize 2} 
\put(3.7,4.7){\em \scriptsize 2} 
\put(5.7,4.7){\em \scriptsize 2} 
\put(7.7,6.7){\em \scriptsize 2} 
\end{picture} \end{center}}. 
However, the latter cannot serve as a supporting graph for a representation of $\mathfrak{g}(\myA_{2})$, due to the imposition of Lemma 2.2 from \cite{DonAdjoint} on length two monochromatic components. 
The latter lemma requires that if a monochromatic component of some supporting graph has length two, has a unique element $\xelt$ of rank zero, and has a unique element $\yelt$ of rank two, then $\xelt$ must be the unique minimal element of this component and $\yelt$ its unique maximal element.  

{\bf [\S \QuasiExampleSection.5:\! Quasi-minuscule DCML's are splitting posets and supporting graphs.]} 
Within the present development, we have yet to demonstrate that quasi-minuscule DCML's are splitting posets for $\chi_{_{\varpi}}^{\mathscr{G}}$. 
In fact, we can show much more. 
We will see that no other diamond-colored modular lattices can serve as splitting posets for $\chi_{_{\varpi}}^{\mathscr{G}}$. 
We will also see that quasi-minuscule DCML's are distinguished supporting graphs for an irreducible $\mathfrak{g}(\mathscr{G})$-module with highest weight $\varpi$; this irreducible representation is sometimes called the `short adjoint representation' since its highest weight is the highest short root, while the (irreducible) adjoint representation has highest weight equal to the highest long root.
These claims are part of the content of the next result, which we prove in parts. 
Only the proofs of the supporting graph claims use case-analysis based on the type classification (i.e.\ $\myA, \myB, \ldots, \myG$) of Coxeter--Dynkin flowers. 
In view of \QMSplitting, from here on we call any $\myfancyQ^{(k)}$ a {\em quasi-minuscule splitting DCML}. 

\noindent 
{\bf \QMSplitting}\ \ {\sl Let $\mathscr{G}$ be a Coxeter--Dynkin flower with highest short root $\varpi$. 
A splitting poset for the Weyl bialternant $\chi_{_{\varpi}}^{\mathscr{G}}$ is a DCML if and only if it is a quasi-minuscule DCML  $\myfancyQ^{(k)}(\varpi)$ for some $k \in I$ such that the corresponding simple root $\alpha_{k}$ is short. 
In particular,} $\WGF\left(\myfancyQ^{(k)}(\varpi);\myvarZ\right) = \chi_{_{\varpi}}^{\mathscr{G}}$. 
{\sl Moreover, each $\myfancyQ^{(k)}(\varpi)$ is a positive-rational, edge-minimizing, and solitary supporting graph for a highest-weight-$\varpi$ irreducible representation of the simple Lie algebra $\mathfrak{g}(\mathscr{G})$, and any edge-minimizing supporting graph for this representation must be a quasi-minuscule DCML.}  

{\em Proof \#1 that quasi-minuscule DCML's are splitting.} 
Since $\myfancyQ^{(k)}(\varpi)$ is $\mathscr{G}$-structured by \ShortAdjointDCML, it suffices to show that $\WGF\left(\myfancyQ^{(k)}(\varpi);\myvarZ\right)$. 
We argue this directly, relying crucially on Freudenthal's Multiplicity Formula (FMF). 
A version of FMF that requires no representation theory was developed in Theorem 2.11 of \cite{DonPosetModels}. 

By \HighestShawtyLemma, $\Pi(\varpi) \cap \Lambda^{+} = \{\zeroweight,\varpi\}$. 
Since $\chi_{_{\varpi}}^{\mathscr{G}} = \sum_{\mu \in \Pi(\varpi)\cap\Lambda^{+}}d_{\varpi,\mu}\zeta_{\mu}$, where `$\zeta_{\mu}$' is a monomial symmetric function, it suffices to determine $d_{\varpi,\varpi}$ and $d_{\varpi,\tinyzeroweight}$. 
By \FTWSF.3, $d_{\varpi,\varpi} = 1$. 
It is an easy exercise to use FMF to get $d_{\varpi,\tinyzeroweight} = \#_{\mbox{\tiny short}}(\mathscr{G})$.
So, $\chi_{_{\varpi}}^{\mathscr{G}} = \zeta_{\varpi} + (\#_{\mbox{\tiny short}}(\mathscr{G}))\zeta_{\tinyzeroweight}$.

Now we show that $\WGF\left(\myfancyQ^{(k)}(\varpi);\myvarZ\right) = \zeta_{\varpi} + (\#_{\mbox{\tiny short}}(\mathscr{G}))\zeta_{\tinyzeroweight}$. 
This follows from the facts, recorded in \ShortAdjointDCML, that $wt(\xelt) = \xelt$ for any $\xelt \in \myfancyQ^{(k)}(\varpi) \cap \Phi_{short}$, where the latter set is just $\Phi_{short}$, and that $wt(\xelt)=\zeroweight$ for any $\xelt \in \myfancyQ^{(k)}(\varpi) \setminus \Phi_{short}$, where the latter set is just $\mathcal{Z}^{(k)}$.\hfill\QED 

{\em Proof \#2 that quasi-minuscule DCML's are splitting.} 
This proof uses \VertexColoringTheorem. 
Most often, the hypotheses of criteria {\sl (2)} of \VertexColoringTheorem\ will not apply to $\myfancyQ^{(k)}(\varpi)$, as there will be $\mbox{\sffamily M}_{5}$-type sublattices. 
So, to apply \VertexColoringTheorem, we need to check that, for all $i \in I$, the $i$-components of $\myfancyQ^{(k)}$ are rank symmetric and we must produce a vertex-coloring function $\kappa: \myfancyQ^{(k)}(\varpi) \setminus \{\varpi\} \longrightarrow I$ and a bijection $\tau: \myfancyQ^{(k)}(\varpi) \setminus \{\varpi\} \longrightarrow \myfancyQ^{(k)}(\varpi) \setminus \{\varpi\}$ satisfying criteria {\sl (1)} of that result. 

We start with a brief analysis of the monochromatic components of $\myfancyQ^{(k)}$. 
Observe that for any $i \in I$ and $\xelt \in \myfancyQ^{(k)}$, the $i$-component $\comp_{i}(\xelt)$ is a two-element chain unless $\xelt \in \{-\alpha_{i},\lsem i,j \rsem,\alpha_{i}\}$ wherein $\alpha_{i}$ and $\alpha_{j}$ are short simple roots and $\langle \alpha_{j},\alpha_{i}^{\vee} \rangle < 0$. 
In the latter case, $\comp_{i}(\xelt)$ has the structure \parbox[c]{1cm}{\begin{center}
\setlength{\unitlength}{0.2cm}
\begin{picture}(4.5,4.5)
\put(2.2,0){\circle*{0.6}} \put(0.2,3){\circle*{0.6}} 
\put(4.2,3){\circle*{0.6}} \put(2.2,6){\circle*{0.6}}
\put(0.2,3){\line(2,3){2}} \put(2.2,0){\line(-2,3){2}} 
\put(4.2,3){\line(-2,3){2}} \put(2.2,0){\line(2,3){2}}
\multiput(2.2,0)(-0.1,0.3){8}{\qbezier(0,0)(-0.025,0.075)(-0.05,0.15)}
\multiput(2.2,0)(0.1,0.3){8}{\qbezier(0,0)(0.025,0.075)(0.05,0.15)}
\multiput(2.2,6)(-0.1,-0.3){8}{\qbezier(0,0)(-0.025,-0.075)(-0.05,-0.15)}
\multiput(2.2,6)(0.1,-0.3){8}{\qbezier(0,0)(0.025,-0.075)(0.05,-0.15)}
\put(0.9,1.5){\scriptsize $i$} \put(3.00,1.5){\scriptsize $i$}
\put(0.9,4.5){\scriptsize $i$} \put(3.00,4.5){\scriptsize $i$}
\end{picture} \vspace*{-0.3in} 
\end{center}}, where the dotted lines indicate the possible presence of other length-two color $i$ chains within this $i$-component. 
So, the $i$-components of $\myfancyQ^{(k)}$ are rank symmetric. 

Next, we define $\kappa$ and $\tau$. 
To this end, we utilize a partial ordering of the color palette $I(\mathscr{H})$ of the sub-embryophyte of $\mathscr{G}$ corresponding to the short simple roots. 
In fact, we view $I(\mathscr{H})$ as a rooted tree with index $k$ at the top and with $i \leq j$ if there is a shortest path in $\Gamma(\mathscr{H})$ from $\gamma_{i}$ to $\gamma_{k}$ that includes $\gamma_{j}$. 
We refer to this poset as $\left(I(\mathscr{H}),\leq\right)$. 
Note that any two adjacent nodes from $\Gamma(\mathscr{G})$ have comparable colors in $\left(I(\mathscr{H}),\leq\right)$. 

We first define $\kappa$ and $\tau$ for those $\xelt$ residing in $\myfancyQ^{(k)}(\varpi) \setminus \{\varpi\}$ but not in $\{-\alpha_{i}\}_{i \in I(\mathscr{H})} \cup \mathcal{Z}^{(k)}$. 
For any such $\xelt$, we set $\tau(\xelt) := \xelt$.  
For $\kappa(\xelt)$, we freely pick any ascendant $\yelt$ of $\xelt$ where, say, $\xelt \myarrow{i} \yelt$ and then assign $\kappa(\xelt) := i$. 
Now, $wt(\tau(\xelt)) = wt(\xelt) = wt(\xelt) - (1 - 1)\alpha_{i} = wt(\xelt) - (1 + \mym_{i}(\xelt))\alpha_{i}$, so criteria {\sl (1)} of \VertexColoringTheorem\ holds so far. 
Next we define $\kappa$ and $\tau$ on $\{-\alpha_{i}\}_{i \in I(\mathscr{H})} \cup \mathcal{Z}^{(k)}$. 
Pick $i \in I(\mathscr{G})$. 
Set $\kappa(-\alpha_{i}) := i$ and $\tau(-\alpha_{i}) := \lsem i,j' \rsem$ where $j'$ is largest in $\left(I(\mathscr{H}),\leq\right)$ such that $\lsem i,j' \rsem \in \mathcal{Z}^{(k)}$. 
Moreover, set $\kappa(\lsem i,j \rsem) := i$ if $i$ is no larger that $j$ in $\mathcal{Z}^{(k)}$, in which case we set $\tau(\lsem i,j \rsem) := -\alpha_{i}$. 
It is easy to see that $\tau^{2} = \tau \circ \tau$ is the identity on $\{-\alpha_{i}\}_{i \in I(\mathscr{H})} \cup \mathcal{Z}^{(k)}$. 
In the preceding notation, $wt(\tau(-\alpha_{i})) = wt(\lsem i,j' \rsem) = \zeroweight = -\alpha_{i} + \alpha_{i} = -\alpha_{i} - (1-2)\alpha_{i} = wt(-\alpha_{i}) - (1-\mym_{i}(-\alpha_{i}))\alpha_{i}$, and $wt(\tau(\lsem i,j \rsem)) = wt(-\alpha_{i}) = -\alpha_{i} = \zeroweight - \alpha_{i} = wt(\lsem i,j \rsem) - (1+0)\alpha_{i} = wt(\lsem i,j \rsem) - (1+\mym_{i}(\lsem i,j \rsem))\alpha_{i}$. 
This confirms all aspects of criteria {\sl (1)} of \VertexColoringTheorem\ and completes the proof.\hfill\QED

{\em Proof that splitting DCML's for $\chi_{_{\varpi}}^{\mathscr{G}}$ must be quasi-minuscule DCML's.} 
Suppose $L$ is a splitting DCML for $\chi_{_{\varpi}}^{\mathscr{G}}$. 
Now $\Pi(L) = \Pi(\varpi)$, and the latter is the set $\Phi_{\mbox{\tiny short}} \disjointunion \{\zeroweight\}$ with the induced order from $\Lambda$. 
Since $d_{\varpi,\alpha} = 1$ for all $\alpha \in \Phi_{\mbox{\tiny short}} = \mathcal{W}\varpi$, it follows that $\xelt \myarrow{i} \yelt$ in $L$ if $wt(\xelt)$ and $wt(\yelt)$ are short roots with $wt(\xelt) \myarrow{i} wt(\yelt)$ in $\Pi(\varpi)$. 
Note that the set of elements of $L$ with weight $\zeroweight$ are exactly those with rank $\langle \varpi,\varrho^{\vee} \rangle = \frac{1}{2}\posetlength(L)$, which is the same for any quasi-minuscule DCML. 
From here on, we identify any non-middle-rank element $\xelt \in L$ with its associated short root $wt(\xelt) \in \Phi_{\mbox{\tiny short}}$. 
So, with possible exception of edges incident with the middle rank of $L$, the edges (and edge colors) of $L$ agree with those of any quasi-minuscule DCML $\myfancyQ^{(k)}(\varpi)$. 

Now, $\alpha_{i}+\alpha_{j}$ is a short root if and only if $\alpha_{i}$ and $\alpha_{j}$ are short simple roots and nodes $\gamma_{i}$ and $\gamma_{j}$ are adjacent $\Gamma(\mathscr{H})$. 
In this case, we have $\alpha_{i} \myarrow{j}\  \alpha_{i}+\alpha_{j}\  \mybackarrow{i} \alpha_{j}$ in $L$. 
Therefore, since $L$ is topographically balanced, there must exist a unique middle rank element, notated `$\zelt_{ij} = \zelt_{ji}$', that is covered by both $\alpha_{i}$ and $\alpha_{j}$. 
Since $L$ is diamond-colored, we have $\alpha_{i} \mybackarrow{i} \zelt_{ij} \myarrow{j} \alpha_{j}$. 
Suppose $k, l \in I(\mathscr{H})$ with $\gamma_{k}$ adjacent to $\gamma_{l}$ in $\Gamma(\mathscr{H})$. 
If $\zelt_{ij} = \zelt_{kl}$, then, without loss of generality, we may assume $i \not= k$. 
So, $\alpha_{i} \mybackarrow{i} \zelt_{ij}=\zelt_{kl} \myarrow{k} \alpha_{k}$ for distinct short simple roots $\alpha_{i}$ and $\alpha_{k}$. 
The topographical balance, diamond coloring, and sturdiness of $L$ imply that $\alpha_{i}+\alpha_{k}$ is also a short root and that 
$\alpha_{i} \myarrow{k}\  \alpha_{i}+\alpha_{k}\  \mybackarrow{i} \alpha_{k}$. 
In particular, $\gamma_{i}$ and $\gamma_{k}$ are adjacent $\Gamma(\mathscr{H})$. 
If $i$ and $l$ are also distinct, then we see that $\gamma_{i}$ and $\gamma_{l}$ are adjacent in $\Gamma(\mathscr{H})$, thus creating a 3-cycle in amongst the nodes $\gamma_{i}$, $\gamma_{k}$, and $\gamma_{l}$. 
This contradicts the fact that $\Gamma(\mathscr{G})$ and hence $\Gamma(\mathscr{H})$ is a tree. 
So, we get $i = l$. 
Similar reasoning forces $j$ to be one of $k$ or $l$, and since $i=l$ then we must have $j=k$. 
We have therefore shown that $\zelt_{ij} = \zelt_{kl}$ if and only if $\{i,j\}=\{k,l\}$. 

Let $\mathcal{N} := \left\{\{i,j\} \subseteq I(\mathscr{H})\, \rule[-1.7mm]{0.2mm}{5mm}\, \gamma_{i} \mbox{ is adjacent to } \gamma_{j} \mbox{ in } \Gamma(\mathscr{H})\right\}$, which is the set of neighboring pairs in $\Gamma(\mathscr{H})$. 
Since there are $\#_{\mbox{\tiny short}}(\mathscr{G})$ short simple roots whose indices comprise the set $I(\mathscr{H})$, and since $\Gamma(\mathscr{H})$ is a tree, then $\mathcal{N}$ consists of exactly $\#_{\mbox{\tiny short}}(\mathscr{G}) - 1$ pairs. 
Therefore the elements $\{\zelt_{ij}\}_{\{i,j\} \in \mathcal{N}}$ comprise all but one of the middle rank elements of $L$. 
Note that any such $\zelt_{ij}$ has exactly two ascendants, one along an edge of color $i$ and the other along an edge of color $j$. 
Now, $\mym_{i}(\alpha_{i}) = \langle \alpha_{i},\alpha_{i}^{\vee} \rangle = 2$; the only short simple root that can have any descendants along edges of color $i$ is $\alpha_{i}$; and $2\alpha_{i}$ is not a root (short or long). 
So $\alpha_{i}$ is the unique maximal element of a length two $i$-component $\comp_{i}(\alpha_{i})$, and, for similar reasons, $-\alpha_{i}$ is its unique minimal element. 
We now see that $\comp_{i}(\alpha_{i})$ is the set $\{\pm \alpha_{i}\} \disjointunion \{\zelt_{ij}\}_{\{i,j\} \in \mathcal{N}}$ with $-\alpha_{i} \myarrow{i} \zelt_{ij} \alpha_{i}$ for each such $\zelt_{ij}$. 
In particular, we see that $\zelt_{ij}$ has exactly two descendants, one along an edge of color $i$ and the other along an edge of color $j$. 

We have now described the adjacencies of all but one of the middle rank elements. 
Let us suppose that this remaining middle rank element is incident with at least one of $\alpha_{k}$ and $-\alpha_{k}$ for some $k \in I(\mathscr{H})$. 
We will argue that this element, for now named `$\zelt_{kk}$', is adjacent to each of $\alpha_{k}$ and $-\alpha_{k}$ and that there are its only two adjacencies. 
Since $\comp_{k}(\alpha_{k}) = \comp_{k}(-\alpha_{k})$ is a length two modular lattice, then, since at least one of $\zelt_{kk} \myarrow{k} \alpha_{k}$ or $-\alpha_{k} \myarrow{k} \zelt_{kk}$ is an edge in $\comp_{k}(\alpha_{k})$, then both are. 
Suppose now that $\zelt_{kk} \myarrow{i} \alpha_{i}$ for some $i \not= k$. 
In view of the topographical balance, diamond-coloring, and sturdiness of $L$, the fact that $\alpha_{k} \mybackarrow{k} \zelt_{kk} \myarrow{i} \alpha_{i}$ implies that $\alpha_{i} + \alpha_{k}$ is a root with $\alpha_{i} \myarrow{k}\  \alpha_{i}+\alpha_{k}\  \mybackarrow{i} \alpha_{k}$ and that $\{i,k\} \in \mathcal{N}$. 
So, $\zelt_{kk} = \alpha_{i} \wedge \alpha_{k} = \zelt_{ik}$, which contradicts the fact that $\zelt_{kk}$ is not one of the middle rank elements $\{\zelt_{pq}\}_{\{p,q\} \in \mathcal{N}}$. 
So, $\alpha_{k}$ is the only simple short root that is incident with $\zelt_{kk}$ in $L$. 
Similarly see that $-\alpha_{k}$ is the only negative simple short root incident with $\zelt_{kk}$. 
This completes our analysis of edges incident with $\zelt_{kk}$ and accounts for all edges incident with the middle rank elements of $L$.  
Therefore the mapping $\phi: L \longrightarrow \myfancyQ^{(k)}$ given by $\alpha \stackrel{\phi}{\longmapsto} \alpha$ and $\zelt_{ij} \stackrel{\phi}{\longmapsto} \lsem i,j \rsem$ is an isomorphism of edge-colored directed graphs.\hfill\QED

{\em Proof that quasi-minuscule splitting DCML's are positive-rational and solitary supporting graphs for the short adjoint representation.} 
For the $\myA$, $\myD$, and $\myE$ cases, the claims that our quasi-minuscule splitting DCML's are positive-rational and solitary supporting graphs for the short adjoint representation follow from Propositions 4.1 and 5.1 of \cite{DonAdjoint}. 
In that paper, edge coefficients for these splitting DCML's are explicitly prescribed. 
In the $\myB$ and $\myG$ cases, there is exactly one quasi-minuscule splitting DCML, and in each case it is a chain. 
In general, in any representation diagram $R$, if an $i$-component $\comp_{i}(\selt)$ of some $\selt \in R$ is a chain with an edge $\selt \myarrow{i} \telt$, then $\myqX_{\telt,\selt}\myqY_{\selt,\telt} = \rho_{i}(\telt)\delta_{i}(\selt)$, cf.\ Lemma 3.9 of \cite{DonSupp}.   
In fact, the latter lemma says that in such an $i$-component, if we assign to each edge $\pelt \myarrow{i} \qelt$ coefficients $\myqX_{\qelt,\pelt}$ and $\myqY_{\pelt,\qelt}$ whose product is $\rho_{i}(\qelt)\delta_{i}(\pelt)$, then the color $i$ crossing relations at the vertices of this $i$-component are satisfied. 
Applying the preceding facts to the quasi-minuscule splitting DCML's in the $\myB$ and $\myG$ cases, we see that these chains satisfy the hypotheses of \MainLieRepTheorem.3 and hence are solitary and positive-rational supporting graphs for the short adjoint representation. 

For case $\myF$, the two quasi-minuscule splitting DCML's are depicted in \ShawtyFFourFigure. 
In each lattice, the $\{3,4\}$-component containing the middle rank vertices is isomorphic to a quasi-minuscule splitting DCML for $\myA_{2}$. 
Assign coefficients to the edges of these $\{3,4\}$-components according to \cite{DonAdjoint}, and then assign $\myqX_{\qelt,\pelt} := 1$ and $\myqY_{\pelt,\qelt} := 1$ to all other splitting DCML edges $\pelt \myarrow{i} \qelt$. 
It is easy to check that all CD relations hold, so these quasi-minuscule splitting DCML's are positive rational supporting graphs for the short adjoint representation. 
Moreover, the edge products are uniquely determined on edges of the $\{3,4\}$-component, since each is a solitary quasi-miniscule splitting DCDL for $\myA_{2}$; edge products are uniquely determined to be unity on any other edge of (say) color $i$ since its $i$-component is a length one chain. 
So, on all edges of each quasi-minuscule DCML, the edge products are uniquely determined. 
We can apply \MainLieRepTheorem.3 to see that each of the quasi-minuscule splitting DCML's from \ShawtyFFourFigure\ are solitary. 

For case $\myC$, the argument is entirely similar to case $\myF$ (see \ShawtyCFourFigure). 
Fix a quasi-minuscule splitting DCML $\myfancyQ^{(k)}$ for $\myC_{n}$. 
The $\{1,2,\ldots,n-1\}$-component of $\myfancyQ^{(k)}$ that contains the middle rank vertices is isomorphic to a quasi-minuscule splitting DCML for $\myA_{n-1}$. 
The latter has edge coefficients explicitly prescribed in \cite{DonAdjoint}, and the product of the coefficients on any such edge is uniquely determined. 
Any edge of (say) color $i$ that is not part of this $\myA_{n-1}$-component comprises the entirety of its $i$-component and therefore has edge product uniquely determined to be unity, and we take each of the $\myqX$- and $\myqY$-coefficients to be `$1$'. 
Given the structure of $\myfancyQ^{(k)}$, it is easy to see that our choices for edge coefficients satisfy the CD relations, and hence $\myfancyQ^{(k)}$ is a positive rational supporting graph for the short adjoint representation. 
Again, we can apply \MainLieRepTheorem.3 to see that our given $\myfancyQ^{(k)}$ is solitary.\hfill\QED

{\em Proof that quasi-minuscule splitting DCML's are exactly the edge-minimizing supporting graphs for the short adjoint representation.} 
That our quasi-minuscule splitting DCML's are positive rational supporting graphs for the short adjoint representation of $\mathfrak{g}(\mathscr{G})$ was just established above. 
It is easy to see that the colors $J$ of the edges incident with the middle rank vertices of a given quasi-minuscule splitting DCML are precisely the colors of the short simple roots. 
Let $\mathscr{H}$ be the corresponding $J$-sub-embryophyte of $\mathscr{G}$. 
Observe that $\mathscr{H}$ is a Coxeter-Dynkin flower of type $\myA_{m}$, $\myD_{m}$, or $\myE_{m}$, where $m := |J|$. 
Now, $\mathscr{H}$ is all of $\mathscr{G}$ precisely in the $\myA$--$\myD$--$\myE$ cases, in which case  Proposition 7.1 of \cite{DonAdjoint} applies to prove that the quasi-minuscule splitting DCML's are exactly the edge-minimizing supporting graphs for the short adjoint representation. 
In the $\myB$, $\myC$, $\myF$, and $\myG$ cases, we observe that Proposition 7.1 of \cite{DonAdjoint} applies to the $J$-component of the middle rank vertices of any quasi-minuscule splitting poset, which forces our quasi-minuscule splitting DCML's to be exactly the edge-minimizing supporting graphs for the short adjoint representation.\hfill\QED

\begin{center}
\underline{\hspace*{4in}}
\end{center}

\vspace*{0.5cm}
\noindent 
{\bf \S \MinusculeExampleSection. Examples of DCDL's that are splitting posets: The minuscule case.} 
In this section, we continue with our development of fundamental examples and take $\mathscr{G} = (\Gamma_{I},M_{I \times I})$ to be a Coxeter--Dynkin flower with $n=|I|$.  
In the previous section, we learned that associated to every such $\mathscr{G}$ is a unique nonzero weight, called the quasi-minuscule dominant weight $\varpi$.  
The quasi-minuscule Weyl bialternant $\chi_{_{\varpi}}^{\mathscr{G}}$ has exactly $\#_{\mbox{\tiny short}}(\mathscr{G})$ splitting diamond-colored modular lattices, where $\#_{\mbox{\tiny short}}(\mathscr{G})$ is the number of short simple roots, cf.\ \QMSplitting. 
However, most of these splitting DCML's are not distributive, and for non-distributive splitting DCML's, companion compression posets are not available. 

The main objects of interest in this section are splitting posets for the minuscule Weyl bialternants. 
It is not hard to establish (see \MinusculeWeylBialt.2 below) that any given minuscule Weyl bialternant has a unique splitting poset. 
One of the main contributions of this section is to establish, using general (i.e.\ type-independent) reasoning, that this unique splitting poset is a diamond-colored distributive lattice. 
This result is not new, but our proof -- which uses some lattice theory together with explicit descriptions of meets and joins -- might be. 
As a follow-up, we use general reasoning to prove that the companion vertex-colored compression poset is itself a distributive lattice that can be realized as a full-length sublattice of a product of two chains. 
It seems there is no such type-independent proof of the latter result in the literature. 
Many of the more well-known results concerning these objects can be found for example in R.\ A.\ Proctor's seminal paper \cite{PrEur} and also in \cite{StrayerThesis}, \cite{StrayerArxiv}, and \cite{Green}. 
For examples, see \ESixFigures.

Here is a useful reminder of some facts about integral embryophytic graphs and their associated Weyl groups: Given distinct colors $i$ and $j$ from the color palette, the nodes $\gamma_{i}$ and $\gamma_{j}$ are nonadjacent in $\Gamma$ if and only if the pulsation matrix entries $M_{ij}$ and $M_{ji}$ are identically zero if and only if the Weyl group generators $\mygens_{i}$ and $\mygens_{j}$ commute. 

{\bf [\S \MinusculeExampleSection.1:\! Some characterizations of minuscule dominant weights.]} 
A minuscule dominant weight is any weight $\omega$ that meets one of the equivalent criteria from the next lemma. 
The equivalences of characterizations {\sl (1)} through {\sl (4)} are well known and demonstrated (for example) in Proposition 4.14 of \cite{DonPosetModels}. 
The equivalence of {\sl (5)} with each of the characterizations {\sl (1)} through {\sl (4)} is easily discerned, but the characterization itself might be new.

\noindent 
{\bf \MinusculeCharacterizationLemma}\ \ {\sl Let $\omega$ be a nonzero dominant weight.  The following are equivalent:}\\ 
{\sl (1) $\Pi(\omega) = \mathcal{W}\omega$ is an equality of sets.}\\
{\sl (2) For all $\nu \in \Lambda^{+}$, if $\nu \in \Pi(\omega)$ then $\nu=\omega$.}\\
{\sl (3) We have $\langle \omega,\alpha^{\vee} \rangle \in \{0,\pm 1\}$ for all roots $\alpha$.}\\
{\sl (4) When we view the weight lattice $\Lambda$ as a poset under the partial ordering given \S \FlowerSection, $\omega$ is the minimal element amongst the dominant weights in some connected component of $\Lambda$.}\\
{\sl (5) In any intermediate position in any Networked-numbers Game played from the initial position $\omega$, the number at any node has absolute value at most $1$. In this case, the numbers at any two distinct nodes of such an intermediate position are both positive only if the two nodes are non-adjacent.} 

{\em Proof.} See \cite{DonPosetModels} Proposition 4.14 for {\sl (1)} $\Leftrightarrow$ {\sl (2)} $\Leftrightarrow$ {\sl (3)} $\Leftrightarrow$ {\sl (4)}. 
Our goal now is to demonstrate {\sl (5)} $\Leftrightarrow$ {\sl (3)}. 
But first, we note that the second sentence of {\sl (5)} follows from the first, because if an intermediate position has positive numbers at a pair of adjacent nodes, then firing either of these nodes is legal and produces a number greater than one at the other node. 

Now, assuming {\sl (3)}, we suppose that $(\gamma_{i_{1}},\ldots,\gamma_{i_{p}})$ is a legal firing sequence from initial position $\omega$. 
Then, for the resulting intermediate position $\mu$, the number at node $\gamma_{j}$ is $\langle \mu,\alpha_{j}^{\vee} \rangle = \langle \mygens_{i_{p}}\cdots\mygens_{i_{1}}.\omega,\alpha_{j}^{\vee} \rangle = \langle \omega,\mygens_{i_{1}}\cdots\mygens_{i_{p}}.\alpha_{j}^{\vee} \rangle = \langle \omega,\alpha^{\vee} \rangle \in \{0,\pm 1\}$, where $\alpha^{\vee}$ is the co-root $\mygens_{i_{1}}\cdots\mygens_{i_{p}}.\alpha_{j}^{\vee}$. 
This establishes that {\sl (3)} implies the first sentence of {\sl (5)}. 

Conversely, suppose the first sentence of {\sl (5)} holds. 
Any given co-root $\alpha^{\vee}$ can be written as $\sigma^{-1}.\alpha_{j}^{\vee}$ for some simple co-root $\alpha_{j}^{\vee}$ and some $\sigma \in \mathcal{W}$. 
Now, $\langle \omega,\alpha^{\vee} \rangle = \langle \omega,\sigma^{-1}.\alpha_{j}^{\vee} \rangle = \langle \sigma.\omega,\alpha_{j}^{\vee} \rangle$. 
Let $J^{c}$ be the set of colors $i \in I$ such that $\langle \omega,\alpha_{i}^{\vee} \rangle > 0$. 
Following \S \WeylSection.3/4, we can write $\sigma$ uniquely as $\sigma^{J} \sigma_{J}$, so that $\sigma.\omega = \sigma^{J}.\omega$. 
Let $\mygens_{i_{p}}\cdots\mygens_{i_{1}}$ be a reduced expression for $\sigma^{J}$. 
By Proposition 3.2 of \cite{DonEur}, it follows that $(\gamma_{i_{1}},\ldots,\gamma_{i_{p}})$ is a legal firing sequence from initial position $\omega$, resulting in an intermediate position $\mu$. 
The number at node $\gamma_{j}$ of $\mu$ is $\langle \mu,\alpha_{j}^{\vee} \rangle$, which, by {\sl (5)}, is in the set $\{0,\pm 1\}$. 
Since $\langle \mu,\alpha_{j}^{\vee} \rangle = \langle \mygens_{i_{p}}\cdots\mygens_{i_{1}}.\omega,\alpha_{j}^{\vee} \rangle = \langle \sigma^{J}.\omega,\alpha_{j}^{\vee} \rangle = \langle \sigma^{J}\sigma_{J}.\omega,\alpha_{j}^{\vee} \rangle = \langle \sigma.\omega,\alpha_{j}^{\vee} \rangle = \langle \omega,\sigma^{-1}.\alpha_{j}^{\vee} \rangle = \langle \omega,\alpha^{\vee} \rangle$, then $\langle \omega,\alpha^{\vee} \rangle \in \{0,\pm 1\}$. 
So, the first sentence of the statement of {\sl (5)} implies {\sl (3)}.\hfill\QED 

That $\omega$ must be some fundamental weight follows nicely from part {\sl (5)}:  
Write $\omega = \sum_{k \in I}a_{k}\omega_{k}$, so by {\sl (5)} each $a_{k} \in \{0,1\}$.  
Suppose for some distinct $i, j \in I$ we have $a_{i} = 1 = a_{j}$. 
Without loss of generality, assume there is a simple path $(\gamma_{i} = \gamma_{i_0},\gamma_{i_1},\ldots,\gamma_{i_{p-1}},\gamma_{i_{p}}=\gamma_{j})$ in $\Gamma$ such that $a_{i_{k}} = 0$ for $0 < k < p$. 
Regarding $\omega$ as our initial position for some Networked-numbers Game, we see that $(\gamma_{i} = \gamma_{i_0},\gamma_{i_1},\ldots,\gamma_{i_{p-1}})$ is a legal firing sequence when played from $\omega$ and, contrary to {\sl (5)}, results in a number larger than $1$ at node $\gamma_{i_{p}}=\gamma_{j}$. 
So, $\omega$ is fundamental, which we generically identify from here on as $\omega_{f}$ and call a {\em minuscule fundamental weight}. 

Following the discussion of \S \WeylSection.3/4, the stabilizer of $\omega_{f}$ under the action of the Weyl group $\mathcal{W}$ is $\mathcal{W}_{\mysmallsetF}$, where $\mysetF := I \setminus \{f\}$, and, moreover, there is a one-to-one correspondence between elements of the $\mathcal{W}$-orbit of $\omega_{f}$ and the minimal coset representatives from the set $\mathcal{W}^{\mysmallsetF}$ given by $\mu \leftrightarrow \sigma$ if and only if $\mu = \sigma.\omega_{f}$. 
In the language of \cite{DonEur}, part {\sl (5)} implies that the weight $\omega_{f}$ is `adjacency-free'. 
Proposition 4.2 of \cite{DonEur} implies that every $\sigma$ from the set $\mathcal{W}^{\mysmallsetF}$ of minimal coset representatives is `fully commutative', which we define as follows:  
Suppose we are given a product $\mygens_{i_{p}} \mygens_{i_{p-1}} \cdots \mygens_{i_{2}}\mygens_{i_{1}}$ of Weyl group generators.    
To {\em apply a commutative switch} to this expression is to find some adjacent pair of factors that commute, say $\mygens_{i_{q}}\mygens_{i_{q-1}}$, and rewrite it as $\mygens_{i_{p}} \mygens_{i_{p-1}} \cdots \mygens_{i_{q-1}}\mygens_{i_{q}} \cdots \mygens_{i_{2}}\mygens_{i_{1}}$. 
A Weyl group element $\sigma$ is {\em fully commutative} if and only if each reduced expression for $\sigma$ can be obtained from some given reduced expression by applying some sequence of commutative switches. 

{\bf [\S \MinusculeExampleSection.2:\! Classification of minuscule fundamental weights.]} 
The minuscule fundamental weights for each Coxeter--Dynkin flower, written as $(\mathscr{G},\omega_{f})$ pairs and with node numbering from \IEGGraphFigure, are $(\myA_{n},\omega_{f})$ for $1 \leq f \leq n$; $(\myB_{n},\omega_{n})$; $(\myC_{n},\omega_{1})$; $(\myD_{n},\omega_{f})$ for $f \in \{1,n-1,n\}$; $(\myE_{6},\omega_{f})$ for $f \in \{1,6\}$; and $(\myE_{7},\omega_{7})$. 
See, for example, Exercise  13 from \S 13 of \cite{Hum}, which applies to minuscule fundamental weights in view of \MinusculeCharacterizationLemma.4. 

Here is an alternative approach: The adjacency-free fundamental weights are classified in Theorem 4.3 of \cite{DonEur}.  
These are identified as those in the above list as well as $(\myB_{n},\omega_{1})$; $(\myC_{n},\omega_{n})$; and $(\myG_{2},\omega_{f})$ for $f \in \{1,2\}$. 
It is easy to see that, in the latter cases, Networked-numbers Game play from the indicated fundamental initial positions violates the main requirement of \MinusculeCharacterizationLemma.5. 
So, the only possible $(\mathscr{G},\omega_{f})$ pairs are those from the first sentence of the previous paragraph; that NG play from these initial positions satisfies the main requirement of \MinusculeCharacterizationLemma.5 is easy to check case-by-case. 

{\bf [\S \MinusculeExampleSection.3:\! Some poset-structural aspects of the saturated set of weights $\Pi(\omega_{f})$.]} 
Most of the following poset-structural claims follow from \SaturatedLemma.2. 
The remaining claims are consequences of \MinusculeCharacterizationLemma. 

\noindent 
{\bf \MinusculeStructureLemma}\ \ {\sl For our given minuscule fundamental weight $\omega_{f}$, any monochromatic component of the edge-colored poset $\Pi(\omega_{f})$ is a chain with no more than two elements
Also, $\Pi(\omega_{f})$ is connected, ranked, $\mathscr{G}$-structured, diamond-colored, and topographically balanced. 
Moreover, if an element of $\Pi(\omega_{f})$ is above an edge of color $i$ and an edge of color $j$, or is below two such colored edges, then $\gamma_{i}$ and $\gamma_{j}$ are distinct and nonadjacent nodes of $\Gamma$, and if $\gamma_{i}$ and $\gamma_{j}$ are distinct and nonadjacent nodes, then any $\{i,j\}$-component of $\Pi(\omega_{f})$ that contains a length two chain is a diamond.} 

{\em Proof.} The lemma claims not covered by \SaturatedLemma.2 are: $\Pi(\omega_{f})$ is topographically balanced, has $i$-components of length at most one for each $i \in I$, only has colors from nonadjacent nodes appearing as edge colors within its diamonds, and has diamonds as its length two bichromatic components when the two colors correspond to nonadjacent nodes from the vascular graph of our given Coxeter--Dynkin flower $\mathscr{G} = (\Gamma,M)$.  
Suppose a weight $\lambda \in \Pi(\omega_{f})$ is $i$-prominent, so that $\delta_{i}(\lambda) = 0$ and the quantity $\rho_{i}(\lambda) = \langle \lambda,\alpha_{i}^{\vee} \rangle$ is nonnegative and also measures the length of $\comp_{i}(\lambda)$. 
Since $\lambda = \sigma.\omega_{f}$ and $\langle \cdot,\cdot \rangle$ is $\mathcal{W}$-invariant, then $0 \leq \langle \sigma.\omega_{f},\alpha_{i}^{\vee} \rangle = \langle \omega_{f},\sigma^{-1}.\alpha_{i}^{\vee} \rangle = \langle \omega_{f},\alpha^{\vee} \rangle$. 
This fact, together with \MinusculeCharacterizationLemma.3, means that $\langle \omega_{f},\alpha^{\vee} \rangle \in \{0,1\}$. 
So, an $i$-component of $\Pi(\omega_{f})$ is a chain with at most two elements. 

Now suppose$\,$ 
\parbox{1.4cm}{\begin{center}
\setlength{\unitlength}{0.2cm}
\begin{picture}(6.5,2)
\put(3,0){\circle*{0.5}} 
\put(1,2){\circle*{0.5}}
\put(5,2){\circle*{0.5}}
\put(3,0){\line(-1,1){2}}
\put(3,0){\line(1,1){2}}
\put(1.75,0.55){\em \small j} 
\put(3.5,0.5){\em \small i}
\put(3.5,-0.75){\footnotesize $\xi$} 
\put(5.75,1.75){\footnotesize $\nu$}
\put(-0.5,1.75){\footnotesize $\mu$}
\end{picture} \end{center}} depicts a pair of covering relations in $\Pi(\omega_{f})$. 
We have $i \not= j$ since $i$-components are chains, so both $\langle \alpha_{j},\alpha_{i}^{\vee} \rangle$ and $\langle \alpha_{i},\alpha_{j}^{\vee} \rangle$ are nonpositive. 
Now, 
$\langle \xi,\alpha_{j}^{\vee} \rangle = \langle \xi,\alpha_{i}^{\vee} \rangle = -1$. 
So, $-1 \leq \langle \mu,\alpha_{i}^{\vee} \rangle = \langle \xi+\alpha_{j},\alpha_{i}^{\vee} \rangle = -1 + \langle \alpha_{j},\alpha_{i}^{\vee} \rangle \leq -1$. 
Then, the inequalities within the previous expression are equalities, which implies that  $\langle \alpha_{j},\alpha_{i}^{\vee} \rangle = 0$ and that $-1 = \langle \mu,\alpha_{i}^{\vee} \rangle$. 
In particular, the distinct nodes $\gamma_{i}$ and $\gamma_{j}$ in $\Gamma$ are nonadjacent, and $\mu \myarrow{i} \mu+\alpha_{i}$ is an edge in $\Pi(\omega_{f})$. 
Similarly, we can show that $\nu \myarrow{i} \nu+\alpha_{j}$ is an edge in $\Pi(\omega_{f})$. 
Of course, $\nu+\alpha_{j} = \mu+\alpha_{i}$, giving us$\,$ 
\parbox{1.4cm}{\begin{center}
\setlength{\unitlength}{0.2cm}
\begin{picture}(6.5,3.5)
\put(3,0){\circle*{0.5}} 
\put(1,2){\circle*{0.5}}
\put(3,4){\circle*{0.5}} 
\put(5,2){\circle*{0.5}}
\put(1,2){\line(1,1){2}} 
\put(3,0){\line(-1,1){2}}
\put(5,2){\line(-1,1){2}} 
\put(3,0){\line(1,1){2}}
\put(1.75,0.55){\em \small j} 
\put(3.5,0.5){\em \small i}
\put(1.7,2.7){\em \small i} 
\put(3.75,2.55){\em \small j}
\put(3.5,-0.75){\footnotesize $\xi$} 
\put(5.75,1.75){\footnotesize $\nu$}
\put(-1,4.5){\scriptsize $\mu+\mysmallestroot_{i} = \nu+\mysmallestroot_{j}$} 
\put(-0.5,1.75){\footnotesize $\mu$}
\end{picture} \end{center}}. 
Since $\Pi(\omega_{f})$ is diamond colored, then any other diamond that includes the edges $\xi \myarrow{j} \mu$ and $\xi \myarrow{i} \nu$ must have parallel edges of colors $i$ and $j$ as depicted here:$\,$ 
\parbox{1.4cm}{\begin{center}
\setlength{\unitlength}{0.2cm}
\begin{picture}(6.5,3.5)
\put(3,0){\circle*{0.5}} 
\put(1,2){\circle*{0.5}}
\put(3,4){\circle*{0.5}} 
\put(5,2){\circle*{0.5}}
\put(1,2){\line(1,1){2}} 
\put(3,0){\line(-1,1){2}}
\put(5,2){\line(-1,1){2}} 
\put(3,0){\line(1,1){2}}
\put(1.75,0.55){\em \small j} 
\put(3.5,0.5){\em \small i}
\put(1.7,2.7){\em \small i} 
\put(3.75,2.55){\em \small j}
\put(3.5,-0.75){\footnotesize $\xi$} 
\put(5.75,1.75){\footnotesize $\nu$}
\put(3.5,4.5){\footnotesize $\lambda$} 
\put(-0.5,1.75){\footnotesize $\mu$}
\end{picture} \end{center}}. 
Since the $i$-component (respectively, $j$-component) of $\mu$ (resp., of $\nu$) is a chain, then we must have $\lambda = \mu+\alpha_{i} = \nu+\alpha_{j}$. 
In particular, there exists a unique $\lambda$ in $\Pi(\omega_{f})$ such that $\mu \rightarrow \lambda$ and $\nu \rightarrow \lambda$. 
Similarly, one can see that if, for some $\lambda', \mu', \nu' \in \Pi(\omega_{f})$, we have $\mu' \myarrow{i'} \lambda'$ and $\nu' \myarrow{j'} \lambda'$, then nodes $\gamma_{i'}$ and $\gamma_{j'}$ in $\Gamma$ are distinct and nonadjacent and there exists a unique $\xi' \in \Pi(\omega_{f})$ such that $\xi' \rightarrow \mu'$ and $\xi' \rightarrow \mu'$.

Last, suppose $\gamma_{i}$ and $\gamma_{j}$ are distinct and nonadjacent nodes in $\Gamma$, and say $\xi \myarrow{j} \mu \myarrow{i} \lambda$ is part of an $\{i,j\}$-component in $\Pi(\omega_{f})$. 
So, $\mym_{i}(\xi) = \langle \xi,\alpha_{i}^{\vee} \rangle = \langle \mu-\alpha_{j},\alpha_{i}^{\vee} \rangle = \langle \mu,\alpha_{i}^{\vee} \rangle = -1$, so $\xi \myarrow{i} \nu$ for some $\nu \in \Pi(\omega_{f})$. 
Moreover, $\mym_{j}(\nu) = \langle \nu,\alpha_{j}^{\vee} \rangle = \langle \xi+\alpha_{i},\alpha_{j}^{\vee} \rangle = \langle \xi,\alpha_{j}^{\vee} \rangle = -1$, so $\nu \myarrow{j} \nu+\alpha_{j}$ is an edge in $\Pi(\omega_{f})$. 
Of course, $\nu+\alpha_{j} = \xi + \alpha_{i} + \alpha_{j} = \lambda$. 
Since no monochromatic component of $\Pi(\omega_{f})$ can have length exceeding one, the edges described in this paragraph must account for all of our given $\{i,j\}$-component, which we now see is$\,$ 
\parbox{1.4cm}{\begin{center}
\setlength{\unitlength}{0.2cm}
\begin{picture}(6.5,3.5)
\put(3,0){\circle*{0.5}} 
\put(1,2){\circle*{0.5}}
\put(3,4){\circle*{0.5}} 
\put(5,2){\circle*{0.5}}
\put(1,2){\line(1,1){2}} 
\put(3,0){\line(-1,1){2}}
\put(5,2){\line(-1,1){2}} 
\put(3,0){\line(1,1){2}}
\put(1.75,0.55){\em \small j} 
\put(3.5,0.5){\em \small i}
\put(1.7,2.7){\em \small i} 
\put(3.75,2.55){\em \small j}
\put(3.5,-0.75){\footnotesize $\xi$} 
\put(5.75,1.75){\footnotesize $\nu$}
\put(3.5,4.5){\footnotesize $\lambda$} 
\put(-0.5,1.75){\footnotesize $\mu$}
\end{picture} \end{center}}.\hfill\QED

{\bf [\S \MinusculeExampleSection.4:\! Some properties and distinctions of minuscule Weyl bialternants.]} 
We call the Weyl symmetric function $\chi_{_{\omega_{f}}}$ a {\em minuscule Weyl bialternant}. 
By \WeylKacCharacter, it is equal to the formal character of the irreducible representation of $\mathfrak{g}(\mathscr{G})$ with highest weight $\omega_{f}$. 
The latter is called a {\em minuscule representation} of $\mathfrak{g}(\mathscr{G})$.  

In \MinusculeWeylBialt\ we record some aspects of minuscule Weyl bialternants and the associated minuscule representations that are well known and easy to prove. 
We re-establish them here as a warm-up for the lattice-theoretic study of their splitting posets. 
The uniqueness result in part {\sl (2)} is equivalent to the fact that all weight spaces of the minuscule representation are one-dimensional, cf.\ Theorem 4.20 from \cite{DonPosetModels}. 
The classification of irreducible representations with one-dimensional weight spaces was re-derived in Theorem 6.7 of the 2003 paper \cite{DonSupp}, having also appeared in \cite{Howe} and \cite{StemMultFree}. 
The uniform and explicit construction of such representations\footnote{Written as $(\mathscr{G},\lambda)$ pairs and with node numbering from \IEGGraphFigure, the irreducible nonminuscule representations with one-dimensional weight spaces are: $(\myA_{n},k\omega_{f})$ with $f \in \{1,n\}$ and $k \in [2,\infty)_{\mathbb{Z}}$; $(\myB_{n},\omega_{1})$; $(\myC_{n},\omega_{n})$ with $n \in \{2,3\}$; and $(\myG_{2};\omega_{1})$. 
In \S 6.7 note 7 of the 2013 book \cite{Green}, Green cites Corollary 1.4.1 from the 2008 paper \cite{Ginzburg} as a source for examples of nonminuscule representations with one-dimensional weight spaces. 
However, in the paragraph of \cite{Ginzburg} that follows Corollary 1.4.1, the classification of nonminuscule representations errantly omits $(\myC_{3},\omega_{3})$.} was first obtained in Theorem 6.4 of the 2003 paper \cite{DonSupp}.

\noindent 
{\bf \MinusculeWeylBialt}\ \ {\sl Take as given our minuscule fundamental weight $\omega_{f}$.  
(1) The minuscule Weyl bialternant $\chi_{_{\omega_{f}}}$ coincides with the monomial symmetric function $\zeta_{\omega_{f}}$. 
(2) The edge-colored poset of weights $\Pi(\omega_{f})$ is, up to edge-colored directed graph isomorphism, the unique splitting poset for $\chi_{_{\omega_{f}}}$.
(3) The edge-colored poset of weights $\Pi(\omega_{f})$ becomes a representation diagram for an irreducible $\mathfrak{g}(\mathscr{G})$-module with highest weight $\omega_{f}$ when we assign the scalars} $\myqX^{(i)}_{\nu,\mu} := 1 =: \myqY^{(i)}_{\mu,\nu}$ {\sl to any edge $\mu \myarrow{i} \nu$ in $\Pi(\omega_{f})$.}

{\em Proof.} For {\sl (1)}, since $\Pi(\omega_{f}) = \mathcal{W}\omega_{f}$, then for every $\mu \in \Pi(\omega_{f})$, the $\mathscr{G}$-Kostka number $d_{\omega_{f},\mu}$ is unity (\FTWSF.3), and therefore $\chi_{_{\omega_{f}}} = \sum_{\mu \in \Pi(\omega_{f})}\myvarZ^{\mu}$.
The latter quantity is $\zeta_{\omega_{f}}$. 
For {\sl (2)}, suppose now that $R$ is any splitting poset for $\chi_{_{\omega_{f}}}$. 
Let $\melt$ be any prominent element of $R$, so $\delta_{i}(\melt) = 0$ for all $i \in I$. 
Since $wt(\melt)$ must be dominant and a member of $\Pi(\omega_{f})$, then by \MinusculeCharacterizationLemma.2, $wt(\melt) = \omega_{f}$. 
So, $R$ has only one prominent element, which must be its unique maximal element. 
By \SaturatedLemma.5, $\Pi(R) = \Pi(\omega_{f})$ is an equality of edge-colored directed graphs. 
Since $d_{\omega_{f},\mu} = 1$ for all $\mu \in \Pi(\omega_{f})$, then $R$ can have only one element $\xelt$ such that $wt(\xelt) = \mu$. 
We see now that $wt: R \longrightarrow \Pi(\omega_{f})$ is an isomorphism of edge-colored directed graphs. 
In particular, $\Pi(\omega_{f})$ is the unique splitting poset for $\chi_{_{\omega_{f}}}$. 
For {\sl (3)}, note that by \MinusculeStructureLemma, any edge $\mu \myarrow{i} \nu$ constitutes an $i$-component that is a chain of length one, and any diamond $\,$ 
\parbox{1.4cm}{\begin{center}
\setlength{\unitlength}{0.2cm}
\begin{picture}(6.5,3.5)
\put(3,0){\circle*{0.5}} 
\put(1,2){\circle*{0.5}}
\put(3,4){\circle*{0.5}} 
\put(5,2){\circle*{0.5}}
\put(1,2){\line(1,1){2}} 
\put(3,0){\line(-1,1){2}}
\put(5,2){\line(-1,1){2}} 
\put(3,0){\line(1,1){2}}
\put(1.75,0.55){\em \small j} 
\put(3.5,0.5){\em \small i}
\put(1.7,2.7){\em \small i} 
\put(3.75,2.55){\em \small j}
\put(3.5,-0.75){\footnotesize $\xi$} 
\put(5.75,1.75){\footnotesize $\nu$}
\put(3.5,4.5){\footnotesize $\lambda$} 
\put(-0.5,1.75){\footnotesize $\mu$}
\end{picture} \end{center}} requires that $\gamma_{i}$ and $\gamma_{j}$ be distinct and nonadjacent, i.e.\ that $\langle \alpha_{i},\alpha_{j}^{\vee} \rangle = 0 = \langle \alpha_{j},\alpha_{i}^{\vee} \rangle$. 
Given these structural facts, it is trivial now to confirm that the assignment of scalars to the edges of $\Pi(\omega_{f})$ as prescribed in the statement of part {\sl (3)} satisfies the CD relations.\hfill\QED 

{\bf [\S \MinusculeExampleSection.5:\! Meets and joins in $\Pi(\omega_{f})$.]} 
The following result is somewhat technical but demonstrates that the saturated set of weights $\Pi(\omega_{f})$ -- which is a connected, ranked, $\mathscr{G}$-structured, diamond-colored, and topographically balanced poset -- is a lattice by explicitly identifying the meet and join of any two weights in $\Pi(\omega_{f})$. 

\noindent
{\bf \MeetJoinMinusculeProp}\ \ {\sl For our given minuscule fundamental weight $\omega_{f}$, let $\mu, \nu \in \Pi(\omega_{f})$. 
Write $\mu = \sigma.\omega_{f}$ with $\sigma \in \mathcal{W}^{\mysmallsetF}$ and $\ell(\sigma) = s$. 
Also, write $\nu = \tau.\omega_{f}$ with $\tau \in \mathcal{W}^{\mysmallsetF}$ and $\ell(\tau) = t$. 
If each of $\sigma$ and $\tau$ has positive length, then we can locate, by algorithm, a number $r \in [1,\min(s,t)]_{\mathbb{Z}}$ and reduced expressions $\sigma = \mygens_{i_{s}}\cdots\mygens_{i_{1}}$ and $\tau = \mygens_{j_{t}}\cdots\mygens_{j_{1}}$ such that (a) $i_{1} = j_{1}, \ldots, i_{r}=j_{r}$; (b) $\mygens_{j_{t}}\cdots\mygens_{j_{r+1}}\mygens_{i_{r+1}}\cdots\mygens_{i_{s}}$ is a reduced expression for $\tau\sigma^{-1}$; (c) the sets $\{i_{r+1},\ldots,i_{s}\} \subseteq I$ and $\{j_{r+1},\ldots,j_{t}\} \subseteq I$ have no colors in common; and (d) $\mygens_{i_{x}}$ commutes with $\mygens_{j_{y}}$ whenever $x \in [r+1,s]_{\mathbb{Z}}$ and $y \in [r+1,t]_{\mathbb{Z}}$.  
Moreover, these outcomes do not depend on any initial choices in the sense that if $\sigma = \mygens_{i'_{s}}\cdots\mygens_{i'_{1}}$ and $\tau = \mygens_{j'_{t}}\cdots\mygens_{j'_{1}}$ such that $i'_{1} = j'_{1}, \ldots, i'_{q}=j'_{q}$ for some $q \geq 1$ and $\mygens_{j'_{t}}\cdots\mygens_{j'_{q+1}}\mygens_{i'_{q+1}}\cdots\mygens_{i'_{s}}$ is a reduced expression for $\tau\sigma^{-1}$, then $q=r$,} $\{i'_{r+1},\ldots,i'_{s}\} \eqmulti \{i_{r+1},\ldots,i_{s}\}$, $\{j'_{r+1},\ldots,j'_{t}\} \eqmulti \{j_{r+1},\ldots,j_{t}\}$,  {\sl and} $\{i'_{1},\ldots,i'_{r}\} \eqmulti \{i_{1},\ldots,i_{r}\}$. 

{\em Proof.} 
Within our finitistic context, we note at the outset the general principle that $(\gamma_{k_{1}},\ldots,\gamma_{k_{q}})$ is a legal firing sequence from initial position $\omega_{f}$ if and only if $\mygens_{k_{q}} \cdots \mygens_{k_{1}}$ is a reduced expression for an element of $\mathcal{W}^{\mysmallsetF}$, cf.\ Proposition 3.2 and Corollary 3.4 of \cite{DonEur}. 
The preceding equivalence can be rephrased as: $(\gamma_{k_{1}},\ldots,\gamma_{k_{q}})$ is a legal firing sequence from initial position $\omega_{f}$ if and only if for each $p \in [1,q-1]_{\mathbb{Z}}$ and weight $\lambda^{(p)} := \mygens_{k_{p}} \cdots \mygens_{k_{1}}.\omega_{f}$ there is a color $k_{p+1}$ edge below $\lambda^{(p)}$ in the weight diagram $\Pi(\omega_{f})$. 
Our argument begins with some fixed reduced expressions for $\sigma$ and $\tau$, say $\sigma = \mygens_{i_{s}}\cdots\mygens_{i_{1}}$ and $\tau = \mygens_{j_{t}}\cdots\mygens_{j_{1}}$. 
Since $\sigma, \tau \in \mathcal{W}^{\mysmallsetF}$, then $i_{1} = j_{1} = f$. 
So, we can write $\tau \sigma^{-1} = \mygens_{j_{t}} \cdots \mygens_{j_{2}} \mygens_{i_{2}} \cdots \mygens_{i_{s}}$. 

\begin{figure}[ht]
\begin{center}
{\small {\bf \MinusculeMeetJoinFigure}\ \  Depiction of part of $\Pi(\omega_{f})$ under the hypothesis that the expression $\mygens_{j_{t}} \cdots \mygens_{j_{2}} \mygens_{i_{2}} \cdots \mygens_{i_{s}}$\\
for $\tau \sigma^{-1}$ is reduced when $\sigma = \mygens_{i_{s}}\cdots\mygens_{i_{1}} \in \mathcal{W}^{\mysmallsetF}$ has length $s$ and $\tau = \mygens_{j_{t}}\cdots\mygens_{j_{1}} \in \mathcal{W}^{\mysmallsetF}$ has length $t$.}

\setlength{\unitlength}{0.4cm}
\begin{picture}(18,18.5)
\put(2,15.9){\large $\omega_{f}$}
\put(-1.15,5.8){\large $\mu$}
\put(16.4,9.8){\large $\nu$}
\put(3,16.5){\line(3,2){2.25}}
\put(5.25,18){\vector(1,0){4.2}}
\put(10,16){\line(0,1){2}}
\put(10,18){\circle*{0.35}}
\put(8.2,16.8){\small $i_{1} = f = j_{1}$}
\multiput(2,8)(0.4,0.4){10}{\qbezier(0.1,0.1)(0.2,0.2)(0.3,0.3)}
\multiput(4,6)(0.4,0.4){10}{\qbezier(0.1,0.1)(0.2,0.2)(0.3,0.3)}
\multiput(6,4)(0.4,0.4){10}{\qbezier(0.1,0.1)(0.2,0.2)(0.3,0.3)}
\multiput(8,2)(0.4,0.4){10}{\qbezier(0.1,0.1)(0.2,0.2)(0.3,0.3)}
\multiput(4,2)(-0.4,0.4){5}{\qbezier(-0.1,0.1)(-0.2,0.2)(-0.3,0.3)}
\multiput(6,4)(-0.4,0.4){5}{\qbezier(-0.1,0.1)(-0.2,0.2)(-0.3,0.3)}
\multiput(10,8)(-0.4,0.4){5}{\qbezier(-0.1,0.1)(-0.2,0.2)(-0.3,0.3)}
\multiput(12,10)(-0.4,0.4){5}{\qbezier(-0.1,0.1)(-0.2,0.2)(-0.3,0.3)}
\multiput(14,12)(-0.4,0.4){5}{\qbezier(-0.1,0.1)(-0.2,0.2)(-0.3,0.3)}
\put(2,4){\circle*{0.35}} 
\put(0,6){\circle*{0.35}}
\put(2,8){\circle*{0.35}} 
\put(4,6){\circle*{0.35}}
\put(0,6){\line(1,1){2}} 
\put(2,4){\line(-1,1){2}}
\put(4,6){\line(-1,1){2}} 
\put(2,4){\line(1,1){2}}
\put(0.7,4.65){\em \small $j_{2}$} 
\put(2.7,4.65){\em \small $i_{s}$}
\put(0.7,6.65){\em \small $i_{s}$} 
\put(2.7,6.65){\em \small $j_{2}$}
\put(8,10){\circle*{0.35}} 
\put(6,12){\circle*{0.35}}
\put(6,12){\line(1,1){2}} 
\put(8,10){\line(-1,1){2}}
\put(8,10){\line(1,1){2}}
\put(6.7,10.65){\em \small $j_{2}$} 
\put(8.7,10.65){\em \small $i_{3}$}
\put(6.7,12.65){\em \small $i_{3}$} 
\put(10,12){\circle*{0.35}} 
\put(8,14){\circle*{0.35}}
\put(10,16){\circle*{0.35}} 
\put(12,14){\circle*{0.35}}
\put(8,14){\line(1,1){2}} 
\put(10,12){\line(-1,1){2}}
\put(12,14){\line(-1,1){2}} 
\put(10,12){\line(1,1){2}}
\put(8.7,12.65){\em \small $j_{2}$} 
\put(10.7,12.65){\em \small $i_{2}$}
\put(8.7,14.65){\em \small $i_{2}$} 
\put(10.7,14.65){\em \small $j_{2}$}
\put(14,8){\circle*{0.35}} 
\put(12,10){\circle*{0.35}}
\put(14,12){\circle*{0.35}} 
\put(16,10){\circle*{0.35}}
\put(12,10){\line(1,1){2}} 
\put(14,8){\line(-1,1){2}}
\put(16,10){\line(-1,1){2}} 
\put(14,8){\line(1,1){2}}
\put(12.7,8.65){\em \small $j_{t}$} 
\put(14.7,8.65){\em \small $i_{2}$}
\put(12.7,10.65){\em \small $i_{2}$} 
\put(14.7,10.65){\em \small $j_{t}$}
\put(12,6){\circle*{0.35}} 
\put(10,8){\circle*{0.35}}
\put(10,8){\line(1,1){2}} 
\put(12,6){\line(-1,1){2}}
\put(12,6){\line(1,1){2}}
\put(10.7,6.65){\em \small $j_{t}$} 
\put(12.7,6.65){\em \small $i_{3}$}
\put(10.7,8.65){\em \small $i_{3}$} 
\put(6,0){\circle*{0.35}} 
\put(4,2){\circle*{0.35}}
\put(6,4){\circle*{0.35}} 
\put(8,2){\circle*{0.35}}
\put(4,2){\line(1,1){2}} 
\put(6,0){\line(-1,1){2}}
\put(8,2){\line(-1,1){2}} 
\put(6,0){\line(1,1){2}}
\put(4.7,0.65){\em \small $j_{t}$} 
\put(6.7,0.65){\em \small $i_{s}$}
\put(4.7,2.65){\em \small $i_{s}$} 
\put(6.7,2.65){\em \small $j_{t}$}
\end{picture} 
\end{center}
\end{figure}

Suppose the expression $\mygens_{j_{t}} \cdots \mygens_{j_{2}} \mygens_{i_{2}} \cdots \mygens_{i_{s}}$ for $\tau \sigma^{-1}$ is reduced. 
We intend to argue that for all $x \in [2,s]_{\mathbb{Z}}$ and $y \in [2,t]_{\mathbb{Z}}$, we have $i_{x} \not= j_{y}$ and $\mygens_{i_{x}}\mygens_{j_{y}} = \mygens_{j_{y}}\mygens_{i_{x}}$. 
Observe that $i_{2} \not= j_{2}$, else our expression for $\tau \sigma^{-1}$ is not reduced. 
The firing sequences $(\gamma_{f},\gamma_{i_{2}})$ and $(\gamma_{f},\gamma_{j_{2}})$ are both legal, and therefore there is an edge of color $i_{2}$ and an edge of color $j_{2}$ below $\omega_{f}-\alpha_{f}$ in $\Pi(\omega_{f})$. 
It follows from \MinusculeStructureLemma\ that this $\{i_{2},j_{2}\}$ component of $\Pi(\omega_{f})$ is a diamond and that the generators $\mygens_{i_{2}}$ and $\mygens_{j_{2}}$ commute. 
In particular, the legal firing sequence $(\gamma_{f},\gamma_{i_{2}})$ can be legally followed by firing $\gamma_{i_{3}}$ (which we already knew from our hypothesis that $\mygens_{i_{s}} \cdots \mygens_{i_{1}}$ is a reduced expression for $\sigma$) or $\gamma_{j_{2}}$ (which we now know since $\omega_{f}-\alpha_{f}-\alpha_{i_{2}}-\alpha_{j_{2}} = \mygens_{j_{2}}\mygens_{i_{2}}\mygens_{f}.\omega_{f}$ is a member of $\Pi(\omega_{f})$). 
Now, if we rewrite $\tau \sigma^{-1}$ as $\mygens_{j_{t}} \cdots \mygens_{j_{3}} \mygens_{i_{2}} \mygens_{j_{2}} \mygens_{i_{3}} \cdots \mygens_{i_{s}}$ by switching `$\mygens_{i_{2}}$' and `$\mygens_{j_{2}}$' in the original expression, then the latter expression remains reduced, and therefore $j_{2} \not= i_{3}$. 
It follows that there is an edge of color $i_{3}$ and an edge of color $j_{2}$ below $\omega_{f}-\alpha_{f}-\alpha_{i_{2}}$ in $\Pi(\omega_{f})$. 
It follows from \MinusculeStructureLemma\ that this $\{i_{3},j_{2}\}$ component of $\Pi(\omega_{f})$ is a diamond and that the generators $\mygens_{i_{3}}$ and $\mygens_{j_{2}}$ commute. 
In particular, the legal firing sequence $(\gamma_{f},\gamma_{i_{2}},\gamma_{i_{3}})$ can be legally followed by firing $\gamma_{i_{4}}$ (which we already knew from our hypothesis that $\mygens_{i_{s}} \cdots \mygens_{i_{1}}$ is a reduced expression for $\sigma$) or $\gamma_{j_{2}}$ (which we now know since $\omega_{f}-\alpha_{f}-\alpha_{i_{2}}-\alpha_{i_{3}}-\alpha_{j_{2}} = \mygens_{j_{2}}\mygens_{i_{3}}\mygens_{i_{2}}\mygens_{f}.\omega_{f}$ is a member of $\Pi(\omega_{f})$). 
Continue reasoning in this way to see that for all $x \in [2,s]_{\mathbb{Z}}$, we have $i_{x} \not= j_{2}$ and $\mygens_{i_{x}}\mygens_{j_{2}} = \mygens_{j_{2}}\mygens_{i_{x}}$. 
Extends this same argument to see that for all $x \in [2,s]_{\mathbb{Z}}$ and $y \in [2,t]_{\mathbb{Z}}$, we have $i_{x} \not= j_{y}$ and $\mygens_{i_{x}}\mygens_{j_{y}} = \mygens_{j_{y}}\mygens_{i_{x}}$. 
So, under the hypotheses at the top of this paragraph, $\Pi(\omega_{f})$ possesses the subposet structure depicted in \MinusculeMeetJoinFigure. 

\begin{figure}[t]
\begin{center}
{\small {\bf \MinusculeMeetJoinFigureTwo}\ \  Depiction of part of $\Pi(\omega_{f})$ under certain hypotheses\\ from the third paragraph of our proof of \MeetJoinMinusculeProp.}

\setlength{\unitlength}{0.4cm}
\begin{picture}(18,22.5)
\put(2,19.9){\large $\omega_{f}$}
\put(3,20.5){\line(3,2){2.25}}
\put(5.25,22){\vector(1,0){4.2}}
\put(-1.15,5.8){\large $\mu$}
\put(18.4,7.8){\large $\kappa$}
\put(8.2,20.8){\small $i_{1} = f = j_{1}$}
\put(10,20){\circle*{0.35}}
\put(10,22){\circle*{0.35}}
\put(10,20){\line(0,1){2}}
\multiput(10,18)(0,0.4){5}{\qbezier(0,0.1)(0,0.2)(0,0.3)}
\put(9,16.8){\small $i_{p} = j_{p}$}
\put(10,16){\line(0,1){2}}
\put(10,18){\circle*{0.35}}
\multiput(2,8)(0.4,0.4){10}{\qbezier(0.1,0.1)(0.2,0.2)(0.3,0.3)}
\multiput(4,6)(0.4,0.4){10}{\qbezier(0.1,0.1)(0.2,0.2)(0.3,0.3)}
\multiput(6,4)(0.4,0.4){10}{\qbezier(0.1,0.1)(0.2,0.2)(0.3,0.3)}
\multiput(8,2)(0.4,0.4){10}{\qbezier(0.1,0.1)(0.2,0.2)(0.3,0.3)}
\multiput(10,0)(0.4,0.4){10}{\qbezier(0.1,0.1)(0.2,0.2)(0.3,0.3)}
\multiput(4,2)(-0.4,0.4){5}{\qbezier(-0.1,0.1)(-0.2,0.2)(-0.3,0.3)}
\multiput(6,4)(-0.4,0.4){5}{\qbezier(-0.1,0.1)(-0.2,0.2)(-0.3,0.3)}
\multiput(10,8)(-0.4,0.4){5}{\qbezier(-0.1,0.1)(-0.2,0.2)(-0.3,0.3)}
\multiput(12,10)(-0.4,0.4){5}{\qbezier(-0.1,0.1)(-0.2,0.2)(-0.3,0.3)}
\multiput(14,12)(-0.4,0.4){5}{\qbezier(-0.1,0.1)(-0.2,0.2)(-0.3,0.3)}
\put(2,4){\circle*{0.35}} 
\put(0,6){\circle*{0.35}}
\put(2,8){\circle*{0.35}} 
\put(4,6){\circle*{0.35}}
\put(0,6){\line(1,1){2}} 
\put(2,4){\line(-1,1){2}}
\put(4,6){\line(-1,1){2}} 
\put(2,4){\line(1,1){2}}
\put(0.7,4.65){\em \small $j_{p+1}$} 
\put(2.7,4.65){\em \small $i_{s}$}
\put(0.7,6.65){\em \small $i_{s}$} 
\put(2.7,6.65){\em \small $j_{p+1}$}
\put(8,10){\circle*{0.35}} 
\put(6,12){\circle*{0.35}}
\put(6,12){\line(1,1){2}} 
\put(8,10){\line(-1,1){2}}
\put(8,10){\line(1,1){2}}
\put(6.7,10.65){\em \small $j_{p+1}$} 
\put(8.7,10.65){\em \small $i_{p+2}$}
\put(6.7,12.65){\em \small $i_{p+2}$} 
\put(10,12){\circle*{0.35}} 
\put(8,14){\circle*{0.35}}
\put(10,16){\circle*{0.35}} 
\put(12,14){\circle*{0.35}}
\put(8,14){\line(1,1){2}} 
\put(10,12){\line(-1,1){2}}
\put(12,14){\line(-1,1){2}} 
\put(10,12){\line(1,1){2}}
\put(8.7,12.65){\em \small $j_{p+1}$} 
\put(10.7,12.65){\em \small $i_{p+1}$}
\put(8.7,14.65){\em \small $i_{p+1}$} 
\put(10.7,14.65){\em \small $j_{p+1}$}
\put(14,8){\circle*{0.35}} 
\put(12,10){\circle*{0.35}}
\put(14,12){\circle*{0.35}} 
\put(16,10){\circle*{0.35}}
\put(12,10){\line(1,1){2}} 
\put(14,8){\line(-1,1){2}}
\put(16,10){\line(-1,1){2}} 
\put(14,8){\line(1,1){2}}
\put(12.7,8.65){\em \small $j_{q}$} 
\put(14.7,8.65){\em \small $i_{p+1}$}
\put(12.7,10.65){\em \small $i_{p+1}$} 
\put(14.7,10.65){\em \small $j_{q}$}
\put(12,6){\circle*{0.35}} 
\put(10,8){\circle*{0.35}}
\put(10,8){\line(1,1){2}} 
\put(12,6){\line(-1,1){2}}
\put(12,6){\line(1,1){2}}
\put(10.7,6.65){\em \small $j_{q}$} 
\put(12.7,6.65){\em \small $i_{p+2}$}
\put(10.7,8.65){\em \small $i_{p+2}$} 
\put(6,0){\circle*{0.35}} 
\put(4,2){\circle*{0.35}}
\put(6,4){\circle*{0.35}} 
\put(8,2){\circle*{0.35}}
\put(4,2){\line(1,1){2}} 
\put(6,0){\line(-1,1){2}}
\put(8,2){\line(-1,1){2}} 
\put(6,0){\line(1,1){2}}
\put(4.7,0.65){\em \small $j_{q}$} 
\put(6.7,0.65){\em \small $i_{s}$}
\put(4.7,2.65){\em \small $i_{s}$} 
\put(6.7,2.65){\em \small $j_{q}$}
\put(16,6){\circle*{0.35}} 
\put(18,8){\circle*{0.35}}
\put(14,8){\line(1,1){2}} 
\put(16,6){\line(-1,1){2}}
\put(18,8){\line(-1,1){2}} 
\put(16,6){\line(1,1){2}}
\put(14.7,6.65){\em \small $j_{q+1}$} 
\put(16.7,6.65){\em \small $i_{p+1}$}
\put(16.7,8.65){\em \small $j_{q+1}$}
\put(14,4){\circle*{0.35}} 
\put(12,6){\circle*{0.35}}
\put(14,4){\line(-1,1){2}}
\put(14,4){\line(1,1){2}}
\put(12.7,4.65){\em \small $j_{q+1}$} 
\put(14.7,4.65){\em \small $i_{p+2}$}
\put(8,-2){\circle*{0.35}} 
\put(10,0){\circle*{0.35}}
\put(8,-2){\line(-1,1){2}}
\put(10,0){\line(-1,1){2}} 
\put(8,-2){\line(1,1){2}}
\put(6.7,-1.35){\em \small $j_{q+1}$} 
\put(8.7,-1.35){\em \small $i_{s}$}
\put(8.7,0.65){\em \small $j_{q+1}$}
\end{picture} 

\vspace*{-0.1in}
\hspace*{0.1in}
\end{center}
\end{figure}

The main part of our argument concerns the following algorithm. 
Suppose that for some $p \in [1,\min(s,t)]_{\mathbb{Z}}$, we have $i_{1}=j_{1},\ldots,i_{p}=j_{p}$ and $\tau \sigma^{-1} = \mygens_{j_{t}}\cdots\mygens_{j_{p+1}}\mygens_{i_{p+1}}\cdots\mygens_{i_{s}}$. 
If this latter expression is reduced, then skip the remainder of this paragraph as well as the two subsequent paragraphs. 
If the expression is not reduced, we proceed as follows.
Then, let $q$ be the unique number in $[p+1,t-1]_{\mathbb{Z}}$ such that the expression $\mygens_{j_{q}}\cdots\mygens_{j_{p+1}}\mygens_{i_{p+1}}\cdots\mygens_{i_{s}}$ is reduced but $\mygens_{j_{q+1}}\mygens_{j_{q}}\cdots\mygens_{j_{p+1}}\mygens_{i_{p+1}}\cdots\mygens_{i_{s}}$ is not. 
Borrowing the reasoning of the preceding paragraph, we see that for all $x \in [p+1,s]_{\mathbb{Z}}$ and $y \in [p+1,q]_{\mathbb{Z}}$, the colors $i_{x}$ and $j_{y}$ are distinct and the generators $\mygens_{i_{x}}$ and $\mygens_{j_{y}}$ commute. 
Let us suppose for the moment that $\mygens_{j_{q+1}}$ is distinct from and commutes with each such $\mygens_{i_{x}}$, and let $\kappa := \mygens_{j_{q+1}}\mygens_{j_{q}}\cdots\mygens_{j_{1}}.\omega_{f}$.  
As argued in the previous paragraph, we obtain a subposet structure within $\Pi(\omega_{f})$ as depicted in \MinusculeMeetJoinFigureTwo. 
Since $\Pi(\omega_{f})$ is $\mathscr{G}$-structured, $\mu-\kappa$ is the integer linear combination of simple roots $\sum_{x =p+1}^{r}\alpha_{i_{x}} - \sum_{y=p+1}^{q+1}\alpha_{j_{y}}$, where there are exactly $(q+1-p)+(r-p)$ simple root vectors on the right-hand side, counting their multiplicities by the absolute value of their integer coefficients. 
Now, 
$\kappa  =  \mygens_{j_{q+1}}\mygens_{j_{q}}\cdots\mygens_{j_{1}}.\omega_{f}  
= (\mygens_{j_{q+1}}\mygens_{j_{q}}\cdots\mygens_{j_{1}})(\mygens_{i_{1}}\cdots\mygens_{i_{s}}).\mu 
= (\mygens_{j_{q+1}}\mygens_{j_{q}}\cdots\mygens_{j_{p+1}}\mygens_{i_{p+1}}\cdots\mygens_{i_{s}}).\mu$, 
and since $\mygens_{j_{q+1}}\mygens_{j_{q}}\cdots\mygens_{j_{p+1}}\mygens_{i_{p+1}}\cdots\mygens_{i_{s}}$ can be reduced, then $\kappa-\mu$ can be expressed as an integer-linear combination of fewer than $(q+1-p)+(r-p)$ simple root vectors (counting multiplicities as before). 
However, this contradicts the independence of simple root vectors within our finitistic context. 
So, it cannot be the case that $\mygens_{j_{q+1}}$ is distinct from and commutes with each such $\mygens_{i_{x}}$, where $x \in [p+1,s]_{\mathbb{Z}}$. 

Therefore, there is a smallest integer $z \in [p+1,s]_{\mathbb{Z}}$ such that one of the following is true: (1) $\mygens_{j_{q+1}} = \mygens_{i_{z}}$ or (2) $j_{q+1} \not\in \{i_{p+1},\ldots,i_{s}\}$ but $\mygens_{j_{q+1}}$ does not commute with $\mygens_{i_{z}}$. 
Let us suppose that we are in case (2). 
So, $\mygens_{i_{x}}$ is distinct from and commutes with each $\mygens_{j_{y}}$ for $x \in [p+1,z]_{\mathbb{Z}}$ and $y \in [p+1,q]_{\mathbb{Z}}$. 
Reasoning as in the second paragraph of the proof we see that the following equalities hold and that all six of the generator product expressions in these equalities are reduced:  
\begin{eqnarray*}
\mygens_{i_{z-1}}\cdots\mygens_{i_{p+1}}\mygens_{j_{q}}\cdots\mygens_{j_{1}}.\omega_{f} & = & \mygens_{j_{q}}\cdots\mygens_{j_{p+1}}\mygens_{i_{z-1}}\cdots\mygens_{i_{1}}.\omega_{f}\\
\mygens_{i_{z}}\cdots\mygens_{i_{p+1}}\mygens_{j_{q}}\cdots\mygens_{j_{1}}.\omega_{f} & = & \mygens_{j_{q}}\cdots\mygens_{j_{p+1}}\mygens_{i_{z}}\cdots\mygens_{i_{1}}.\omega_{f}\\
\mygens_{i_{z-1}}\cdots\mygens_{i_{p+1}}\mygens_{j_{q+1}}\cdots\mygens_{j_{1}}.\omega_{f} & = & \mygens_{j_{q+1}}\cdots\mygens_{j_{p+1}}\mygens_{i_{z-1}}\cdots\mygens_{i_{1}}.\omega_{f}. 
\end{eqnarray*}
With $\lambda := \mygens_{i_{z-1}}\cdots\mygens_{i_{p+1}}\mygens_{j_{q}}\cdots\mygens_{j_{1}}.\omega_{f}$, i.e.\ the first of the three different weights described by the three above equalities, we see that $\lambda - \alpha_{i_{z}}$ is the weight described by the second equality and that $\lambda - \alpha_{j_{q+1}}$ is the weight described by the third equality. 
Then, in $\Pi(\omega_{f})$, the weight $\lambda$ is above an edge of color $i_{z}$ and is also above a second distinct edge of color $j_{q+1}$. 
But this means that $\mygens_{j_{q+1}}$ commutes with $\mygens_{i_{z}}$, violating the hypothesis of case (2).

We have ruled out the second of the two cases described at the top of the preceding paragraph.  
So, case (1) must hold, with $\mygens_{j_{q+1}} = \mygens_{i_{z}}$. 
Directly from the definition of $q$, we discerned in the third paragraph of the proof that  $\mygens_{i_{z}}$ is distinct from and commutes with each of $\mygens_{j_{p+1}},\ldots,\mygens_{i_{q}}$. 
From the definition of $z$, it follows that $\mygens_{j_{q+1}}$ is distinct from and commutes with each of $\mygens_{i_{p+1}},\ldots,\mygens_{i_{z-1}}$.   
So, $\tau\sigma^{-1}$ is 
{\small \begin{eqnarray*}
(\mygens_{j_{t}}\cdots\mygens_{j_{q+1}}\mygens_{j_{q}}\cdots\mygens_{j_{p+1}})(\mygens_{i_{p+1}}\cdots\mygens_{i_{z-1}}\mygens_{i_{z}}\cdots\mygens_{i_{s}}) & = & (\mygens_{j_{t}}\cdots\mygens_{j_{q+2}}\mygens_{j_{q}}\cdots\mygens_{j_{p+1}})(\mygens_{j_{q+1}}\mygens_{i_{z}})(\mygens_{i_{p+1}}\cdots\mygens_{i_{z-1}}\mygens_{i_{z+1}}\cdots\mygens_{i_{s}}),
\end{eqnarray*}}which reduces to $(\mygens_{j_{t}}\cdots\mygens_{j_{q+2}}\mygens_{j_{q}}\cdots\mygens_{j_{p+1}})(\mygens_{i_{p+1}}\cdots\mygens_{i_{z-1}}\mygens_{i_{z+1}}\cdots\mygens_{i_{s}})$. 
Take $p' := p+1$ and reindex the expressions for $\sigma$ and $\tau$ as follows: 
Write $\sigma = \mygens_{i'_{s}}\cdots\mygens_{i'_{p'+1}}\mygens_{i'_{p'}}\cdots\mygens_{i'_{1}}$, with $i'_{p'} = i_{z}$, $i'_{x} = i_{x-1}$ for $x \in [p'+1,z]_{\mathbb{Z}}$, and $i'_{x} = i_{x}$ otherwise; and write $\tau = \mygens_{j'_{t}}\cdots\mygens_{j'_{p'+1}}\mygens_{j'_{p'}}\cdots\mygens_{j'_{1}}$, with $j'_{p'}=j_{q+1}$, $j'_{y} = j_{y-1}$ for $y \in [p'+1,q+1]_{\mathbb{Z}}$, and $j'_{y} = j_{y}$ otherwise.
So, $\tau\sigma^{-1} = (\mygens_{j'_{t}}\cdots\mygens_{j'_{p'+1}})(\mygens_{i'_{p'+1}}\cdots\mygens_{i'_{s}})$ with $i'_{1}=j'_{1},\ldots,i'_{p'}=j'_{p'}$. 
Now return to the beginning of the third paragraph of the proof and apply the prescribed steps to these revised expressions for $\sigma$, $\tau$, and $\tau\sigma^{-1}$. 

\begin{figure}[t]
\begin{center}
{\small {\bf \MinusculeMeetJoinFigureThree}\ \  Depiction of part of $\Pi(\omega_{f})$ representing the conclusions expressed in the sixth\\ paragraph of our proof of \MeetJoinMinusculeProp\ resulting from the algorithm of paragraphs three through five.}

\setlength{\unitlength}{0.4cm}
\begin{picture}(18,22.5)
\put(2,19.9){\large $\omega_{f}$}
\put(3,20.5){\line(3,2){2.25}}
\put(5.25,22){\vector(1,0){4.2}}
\put(-1.15,5.8){\large $\mu$}
\put(18.4,7.8){\large $\nu$}
\put(8.2,20.8){\small $i_{1} = f = j_{1}$}
\put(10,20){\circle*{0.35}}
\put(10,22){\circle*{0.35}}
\put(10,20){\line(0,1){2}}
\multiput(10,18)(0,0.4){5}{\qbezier(0,0.1)(0,0.2)(0,0.3)}
\put(9,16.8){\small $i_{r} = j_{r}$}
\put(10,16){\line(0,1){2}}
\put(10,18){\circle*{0.35}}
\multiput(2,8)(0.4,0.4){10}{\qbezier(0.1,0.1)(0.2,0.2)(0.3,0.3)}
\multiput(4,6)(0.4,0.4){10}{\qbezier(0.1,0.1)(0.2,0.2)(0.3,0.3)}
\multiput(6,4)(0.4,0.4){10}{\qbezier(0.1,0.1)(0.2,0.2)(0.3,0.3)}
\multiput(8,2)(0.4,0.4){10}{\qbezier(0.1,0.1)(0.2,0.2)(0.3,0.3)}
\multiput(10,0)(0.4,0.4){10}{\qbezier(0.1,0.1)(0.2,0.2)(0.3,0.3)}
\multiput(4,2)(-0.4,0.4){5}{\qbezier(-0.1,0.1)(-0.2,0.2)(-0.3,0.3)}
\multiput(6,4)(-0.4,0.4){5}{\qbezier(-0.1,0.1)(-0.2,0.2)(-0.3,0.3)}
\multiput(10,8)(-0.4,0.4){5}{\qbezier(-0.1,0.1)(-0.2,0.2)(-0.3,0.3)}
\multiput(12,10)(-0.4,0.4){5}{\qbezier(-0.1,0.1)(-0.2,0.2)(-0.3,0.3)}
\multiput(14,12)(-0.4,0.4){5}{\qbezier(-0.1,0.1)(-0.2,0.2)(-0.3,0.3)}
\put(2,4){\circle*{0.35}} 
\put(0,6){\circle*{0.35}}
\put(2,8){\circle*{0.35}} 
\put(4,6){\circle*{0.35}}
\put(0,6){\line(1,1){2}} 
\put(2,4){\line(-1,1){2}}
\put(4,6){\line(-1,1){2}} 
\put(2,4){\line(1,1){2}}
\put(0.7,4.65){\em \small $j_{r+1}$} 
\put(2.7,4.65){\em \small $i_{s}$}
\put(0.7,6.65){\em \small $i_{s}$} 
\put(2.7,6.65){\em \small $j_{r+1}$}
\put(8,10){\circle*{0.35}} 
\put(6,12){\circle*{0.35}}
\put(6,12){\line(1,1){2}} 
\put(8,10){\line(-1,1){2}}
\put(8,10){\line(1,1){2}}
\put(6.7,10.65){\em \small $j_{r+1}$} 
\put(8.7,10.65){\em \small $i_{r+2}$}
\put(6.7,12.65){\em \small $i_{r+2}$} 
\put(10,12){\circle*{0.35}} 
\put(8,14){\circle*{0.35}}
\put(10,16){\circle*{0.35}} 
\put(12,14){\circle*{0.35}}
\put(8,14){\line(1,1){2}} 
\put(10,12){\line(-1,1){2}}
\put(12,14){\line(-1,1){2}} 
\put(10,12){\line(1,1){2}}
\put(8.7,12.65){\em \small $j_{r+1}$} 
\put(10.7,12.65){\em \small $i_{r+1}$}
\put(8.7,14.65){\em \small $i_{r+1}$} 
\put(10.7,14.65){\em \small $j_{r+1}$}
\put(14,8){\circle*{0.35}} 
\put(12,10){\circle*{0.35}}
\put(14,12){\circle*{0.35}} 
\put(16,10){\circle*{0.35}}
\put(12,10){\line(1,1){2}} 
\put(14,8){\line(-1,1){2}}
\put(16,10){\line(-1,1){2}} 
\put(14,8){\line(1,1){2}}
\put(12.7,8.65){\em \small $j_{t-1}$} 
\put(14.7,8.65){\em \small $i_{r+1}$}
\put(12.7,10.65){\em \small $i_{r+1}$} 
\put(14.7,10.65){\em \small $j_{t-1}$}
\put(12,6){\circle*{0.35}} 
\put(10,8){\circle*{0.35}}
\put(10,8){\line(1,1){2}} 
\put(12,6){\line(-1,1){2}}
\put(12,6){\line(1,1){2}}
\put(10.7,6.65){\em \small $j_{t-1}$} 
\put(12.7,6.65){\em \small $i_{r+2}$}
\put(10.7,8.65){\em \small $i_{r+2}$} 
\put(6,0){\circle*{0.35}} 
\put(4,2){\circle*{0.35}}
\put(6,4){\circle*{0.35}} 
\put(8,2){\circle*{0.35}}
\put(4,2){\line(1,1){2}} 
\put(6,0){\line(-1,1){2}}
\put(8,2){\line(-1,1){2}} 
\put(6,0){\line(1,1){2}}
\put(4.7,0.65){\em \small $j_{t-1}$} 
\put(6.7,0.65){\em \small $i_{s}$}
\put(4.7,2.65){\em \small $i_{s}$} 
\put(6.7,2.65){\em \small $j_{t-1}$}
\put(16,6){\circle*{0.35}} 
\put(18,8){\circle*{0.35}}
\put(14,8){\line(1,1){2}} 
\put(16,6){\line(-1,1){2}}
\put(18,8){\line(-1,1){2}} 
\put(16,6){\line(1,1){2}}
\put(14.7,6.65){\em \small $j_{t}$} 
\put(16.7,6.65){\em \small $i_{r+1}$}
\put(16.7,8.65){\em \small $j_{t}$}
\put(14,4){\circle*{0.35}} 
\put(12,6){\circle*{0.35}}
\put(14,4){\line(-1,1){2}}
\put(14,4){\line(1,1){2}}
\put(12.7,4.65){\em \small $j_{t}$} 
\put(14.7,4.65){\em \small $i_{r+2}$}
\put(8,-2){\circle*{0.35}} 
\put(10,0){\circle*{0.35}}
\put(8,-2){\line(-1,1){2}}
\put(10,0){\line(-1,1){2}} 
\put(8,-2){\line(1,1){2}}
\put(6.7,-1.35){\em \small $j_{t}$} 
\put(8.7,-1.35){\em \small $i_{s}$}
\put(8.7,0.65){\em \small $j_{t}$}
\end{picture} 

\vspace*{-0.1in}
\hspace*{0.1in}
\end{center}
\end{figure}

The result of the preceding process is that we obtain an integer $r \in [1,\min(s,t)]_{\mathbb{Z}}$ and reduced expressions $\sigma = \mygens_{i_{s}}\cdots\mygens_{i_{r+1}}\mygens_{i_{r}}\cdots\mygens_{i_{1}}$, $\tau = \mygens_{j_{t}}\cdots\mygens_{j_{r+1}}\mygens_{j_{r}}\cdots\mygens_{j_{1}}$, and $\tau\sigma^{-1} = \mygens_{j_{t}}\cdots\mygens_{j_{r+1}}\mygens_{i_{r+1}}\cdots\mygens_{i_{s}}$, where $i_{1}=j_{1},\ldots,i_{r}=j_{r}$ and where $\mygens_{i_{x}}$ is distinct from and commutes with $\mygens_{j_{y}}$ for all $x \in [r+1,s]_{\mathbb{Z}}$ and $y \in [r+1,t]_{\mathbb{Z}}$. 
We claim that $(\mygens_{j_{t}}\cdots\mygens_{j_{r+1}})(\mygens_{i_{s}}\cdots\mygens_{i_{r+1}}\mygens_{i_{r}}\cdots\mygens_{i_{1}})$ is a reduced expression for an element $\upsilon$ of the Weyl group $\mathcal{W}(\mathscr{G})$. 
Of course, the `$\mygens_{i_{s}}\cdots\mygens_{i_{r+1}}\mygens_{i_{r}}\cdots\mygens_{i_{1}}$' part of this expression is a reduced expression for $\sigma$, so the given expression for $\upsilon$ can only be reduced by interactions between the generators $\mygens_{i_{x}}$ (for $x \in [1,s]_{\mathbb{Z}}$) and $\mygens_{j_{y}}$ (for $y \in [r+1,t]_{\mathbb{Z}}$). 
Since each $\mygens_{i_{x}}$ is distinct from and commutes with each $\mygens_{j_{y}}$ for all $x \in [1,s]_{\mathbb{Z}}$ and $y \in [r+1,t]_{\mathbb{Z}}$, then our given expression for $\upsilon$ can only be reduced by interactions between the generators $\mygens_{i_{x}}$ (for $x \in [1,r]_{\mathbb{Z}}$) and $\mygens_{j_{y}}$ (for $y \in [r+1,t]_{\mathbb{Z}}$). 
But $\mygens_{i_{x}} = \mygens_{j_{x}}$ for $x \in [1,r]_{\mathbb{Z}}$, hence $\mygens_{j_{t}}\cdots\mygens_{j_{r+1}}\mygens_{i_{r}}\cdots\mygens_{i_{1}} = \mygens_{j_{t}}\cdots\mygens_{j_{r+1}}\mygens_{j_{r}}\cdots\mygens_{j_{1}}$, where the latter is a reduced expression for $\tau$. 
So the given expression $(\mygens_{j_{t}}\cdots\mygens_{j_{r+1}})(\mygens_{i_{s}}\cdots\mygens_{i_{r+1}}\mygens_{i_{r}}\cdots\mygens_{i_{1}})$ is reduced. 
Moreover, reasoning as in the second paragraph of this proof, we see that $\Pi(\omega_{f})$ possesses the subposet structure depicted in \MinusculeMeetJoinFigureThree. 

We aim to demonstrate a kind of uniqueness of the output of the foregoing algorithm. 
We have, in hand, reduced expressions for $\sigma$, $\tau$, and $\tau \sigma^{-1}$ as in the preceding paragraph. 
Now suppose there exists some $q \in [1,\min(s,t)]_{\mathbb{Z}}$ such that $\sigma = \mygens_{i'_{s}}\cdots\mygens_{i'_{q+1}}\mygens_{i'_{q}}\cdots\mygens_{i'_{1}}$; $\tau = \mygens_{j'_{t}}\cdots\mygens_{j'_{q+1}}\mygens_{j'_{q}}\cdots\mygens_{j'_{1}}$; $i'_{1}=j'_{1},\ldots,i'_{q}=j'_{q}$; and $\mygens_{j'_{t}}\cdots\mygens_{j'_{q+1}}\mygens_{i'_{q+1}}\cdots\mygens_{i'_{s}}$ is a reduced expression for $\tau\sigma^{-1}$. 
The expressions `$\mygens_{j_{t}}\cdots\mygens_{j_{r+1}}\mygens_{i_{r+1}}\cdots\mygens_{i_{s}}$' and `$\mygens_{j'_{t}}\cdots\mygens_{j'_{q+1}}\mygens_{i'_{q+1}}\cdots\mygens_{i'_{s}}$' for $\tau\sigma^{-1}$ are both reduced, so $q=r$. 
Since $\sigma$ is fully commutative, then the multisets $\{i_{1},\ldots,i_{s}\}$ and $\{i'_{1},\ldots,i'_{s}\}$, which represent our two reduced expressions for $\sigma$, must coincide. 
Similarly, $\{j_{1},\ldots,j_{t}\} \eqmulti \{j'_{1},\ldots,j'_{r}\}$. 
Therefore, the hoped-for multiset equalities $\{i_{r+1},\ldots,i_{s}\} \eqmulti \{i'_{r+1},\ldots,i'_{s}\}$ and $\{j_{r+1},\ldots,j_{t}\} \eqmulti \{j'_{r+1},\ldots,j'_{t}\}$ will follow if we can establish the equality $\{i_{1},\ldots,i_{r}\} \eqmulti \{i'_{1},\ldots,i'_{r}\}$. 

We will establish that $\{i_{1},\ldots,i_{r}\} \eqmulti \{i'_{1},\ldots,i'_{r}\}$ by induction on $\min(s,t)$. 
Take as our base case $\min(s,t)=1$; then $r = 1$. 
In this case, we have $i_{1} = f = i'_{1}$. 
For our inductive hypothesis, suppose for some positive integer $m$, it is the case that whenever $\sigma$ and $\tau$ as above satisfy $\min(s,t) \leq m$ with $r \in [1,\min(s,t)]_{\mathbb{Z}}$, we get $\{i_{1},\ldots,i_{r}\} \eqmulti \{i'_{1},\ldots,i'_{r}\}$. 
Now suppose $\min(s,t) = m+1$. 
Without loss of generality, assume $s = m+1$. 
Set $\lambda := \mygens_{i_{r+1}}\cdots\mygens_{i_{s}}.\mu = \mygens_{j_{r+1}}\cdots\mygens_{j_{t}}.\nu$, and let $\lambda' := \mygens_{i'_{r+1}}\cdots\mygens_{i'_{s}}.\mu = \mygens_{j'_{r+1}}\cdots\mygens_{j'_{t}}.\nu$. 
If $r = m+1 = s$, then the equality of multisets $\{i_{1},\ldots,i_{s}\}$ and $\{i'_{1},\ldots,i'_{s}\}$ means that $\{i_{1},\ldots,i_{r}\} \eqmulti \{i'_{1},\ldots,i'_{r}\}$. 
From here on, then, assume that $r < m+1$. 
We consider the following exhaustive and mutually exclusive cases: (1) $i_{s} = i'_{x}$ for some largest $x \in [r+1,s]_{\mathbb{Z}}$, (2) $i_{s} \not\in \{i'_{r+1},\ldots,i'_{s}\}$ but $i_{s} = j'_{y}$ for some smallest $y \in [r+1,t]_{\mathbb{Z}}$, or (3) $i_{s} \not\in \{i'_{r+1},\ldots,i'_{s}\}$ and $i_{s} \not\in \{j'_{r+1},\ldots,j'_{t}\}$. 

For case (1), consider $\mygens_{j'_{t}}\cdots\mygens_{j'_{r+1}}\mygens_{i'_{r+1}}\cdots\mygens_{i'_{x}}\cdots\mygens_{i'_{s}}$. 
If $x=s$, then take $\sigma' := \mygens_{i'_{s-1}}\cdots\mygens_{i'_{r+1}}\mygens_{i'_{r}}\cdots\mygens_{i'_{1}} = \mygens_{i_{s-1}}\cdots\mygens_{i_{r+1}}\mygens_{i_{r}}\cdots\mygens_{i_{1}}$ and $\mu' := \sigma'.\omega_{f} = \mu+\alpha_{i_{s}}$. 
Notice that $\tau\sigma^{-1} = (\tau(\sigma')^{-1})(\mygens_{i_{s}})$, so the expression $\mygens_{j'_{t}}\cdots\mygens_{j'_{r+1}}\mygens_{i'_{r+1}}\cdots\mygens_{i'_{s}}$ for $\tau(\sigma')^{-1}$ is reduced. 
We can now apply the inductive hypothesis to the pair $\sigma'$ and $\tau$ to get $\{i_{1},\ldots,i_{r}\} \eqmulti \{i'_{1},\ldots,i'_{r}\}$. 
If $x < s$, then $i_{s} \not\in \{i'_{x+1},\ldots,i'_{s}\}$. 
Then, $\mu$ is below edges of distinct colors $i_{s}$ and $i'_{s}$, so $\mygens_{i_{s}}$ and $\mygens_{i'_{s}}$ commute. 
Observe, then, that $\mu+\alpha_{i'_{s}}$ is below edges of distinct colors $i_{s}$ and $i'_{s-1}$, so $\mygens_{i_{s}}$ and $\mygens_{i'_{s-1}}$ commute. 
Continue in this way to see that $\mygens_{i_{s}}$ also commutes with $\mygens_{i'_{s-2}},\ldots,\mygens_{i'_{x+1}}$. 
So, 
\begin{eqnarray*}
\sigma & = & \mygens_{i'_{s}}\mygens_{i'_{s-1}}\cdots\mygens_{i'_{x+1}}\mygens_{i'_{x}}\mygens_{i'_{x-1}}\cdots\mygens_{i'_{r+1}}\mygens_{i'_{r}}\cdots\mygens_{i'_{1}}\\
& = & (\mygens_{i_{s}})(\mygens_{i'_{s}}\mygens_{i'_{s-1}}\cdots\mygens_{i'_{x+1}}\mygens_{i'_{x-1}}\cdots\mygens_{i'_{r+1}}\mygens_{i'_{r}}\cdots\mygens_{i'_{1}}).
\end{eqnarray*}
Take $\sigma' := \mygens_{i'_{s}}\mygens_{i'_{s-1}}\cdots\mygens_{i'_{x+1}}\mygens_{i'_{x-1}}\cdots\mygens_{i'_{r+1}}\mygens_{i'_{r}}\cdots\mygens_{i'_{1}}$. 
Then, as before, we have $\tau\sigma^{-1} = (\tau(\sigma')^{-1})(\mygens_{i_{s}})$, so $\mygens_{j'_{t}}\cdots\mygens_{j'_{r+1}}\mygens_{i'_{r+1}}\cdots\mygens_{i'_{x-1}}\mygens_{i'_{x+1}}\cdots\mygens_{i'_{s}}$ is a reduced expression for $\tau(\sigma')^{-1}$. 
Again, we apply the inductive hypothesis to the pair $\sigma'$ and $\tau$ to get $\{i_{1},\ldots,i_{r}\} \eqmulti \{i'_{1},\ldots,i'_{r}\}$. 

For case (2), we can use the reasoning of the previous paragraph to see that the hypothesis $i_{s} \not\in \{i'_{r+1},\ldots,i'_{s}\}$ implies that $\mygens_{i_{s}}$ commutes with each $\mygens_{i'_{x}}$ for $x \in [r+1,s]_{\mathbb{Z}}$.  
We are given that $y \in [r+1,t]_{\mathbb{Z}}$ is smallest such that $i_{s} = j'_{y}$, so we similarly see that $\mygens_{i_{s}}$ commutes with each of $\mygens_{j'_{r+1}},\ldots,\mygens_{j'_{y-1}}$. 
Let $\lambda' := \mygens_{i'_{r}}\cdots\mygens_{i'_{1}}.\omega_{f}$, so $\mu = \mygens_{i'_{s}}\cdots\mygens_{i'_{r+1}}.\lambda'$. 
Check that $\mygens_{i_{s}}.\lambda' = \lambda'+\alpha_{i_{s}}$, so in particular there is a color $i_{s}$ edge above $\lambda'$. 
If $i_{s} \not= i'_{r}$, then $\lambda'$ is below an edge of color $i'_{r}$, so $\mygens_{i_{s}}$ commutes with $\mygens_{i'_{r}}$ and $\lambda'+\alpha_{i'_{r}}$ is below an edge of color $i_{s}$. 
Further, if $i_{s} \not\in \{i'_{r},i'_{r-1}\}$, then, since $\lambda'+\alpha_{i'_{r}}$ is also below an edge of color $i'_{r-1}$, we see that $\mygens_{i_{s}}$ commutes with $\mygens_{i'_{r-1}}$ and $\lambda'+\alpha_{i'_{r}}+\alpha_{i'_{r-1}}$ is below an edge of color $i_{s}$. 
Continuing in this way, we see that if $i_{s} \not\in \{i'_{r},i'_{r-1},\ldots,i'_{1}\}$ with $\lambda'+\alpha_{i'_{r}}+\alpha_{i'_{r-1}}+\cdots+\alpha_{i'_{2}}$ below an edge of color $i_{s}$, then $\mygens_{i_{s}}$ commutes with $\mygens_{i'_{1}}$ and $\lambda'+\alpha_{i'_{r}}+\alpha_{i'_{r-1}}+\cdots+\alpha_{i'_{2}}+\alpha_{i'_{1}}$ is below an edge of color $i_{s}$. 
But $\lambda'+\alpha_{i'_{r}}+\alpha_{i'_{r-1}}+\cdots+\alpha_{i'_{2}}+\alpha_{i'_{1}} = \mygens_{i'_{1}}\mygens_{i'_{2}}\cdots\mygens_{i'_{r-1}}\mygens_{i'_{r}}.\lambda' = \omega_{f}$, which is maximal in $\Pi(\omega_{f})$ and therefore cannot be below an edge of color $i_{s}$. 
We conclude that there is some largest $z \in [1,r]_{\mathbb{Z}}$ such that $i_{s} = i'_{z}$ wherein $\mygens_{i_{s}}$ commutes each $\mygens_{i'_{x}}$ for each $x \in [z+1,r]_{\mathbb{Z}}$. 
Then, 
\begin{eqnarray*}
\tau & = & \mygens_{j'_{t}}\cdots\mygens_{j'_{y+1}}\mygens_{j'_{y}}\mygens_{j'_{y-1}}\cdots\mygens_{j'_{r+1}}\mygens_{j'_{r}}\cdots\mygens_{j'_{z+1}}\mygens_{j'_{z}}\mygens_{j'_{z-1}}\cdots\mygens_{j'_{1}}\\
& = & \mygens_{j'_{t}}\cdots\mygens_{j'_{y+1}}\mygens_{j'_{y-1}}\cdots\mygens_{j'_{r+1}}\mygens_{j'_{y}}\mygens_{j'_{z}}\mygens_{j'_{r}}\cdots\mygens_{j'_{z+1}}\mygens_{j'_{z-1}}\cdots\mygens_{j'_{1}}\\
& = & \mygens_{j'_{t}}\cdots\mygens_{j'_{y+1}}\mygens_{j'_{y-1}}\cdots\mygens_{j'_{r+1}}\mygens_{i_{s}}\mygens_{i_{s}}\mygens_{j'_{r}}\cdots\mygens_{j'_{z+1}}\mygens_{j'_{z-1}}\cdots\mygens_{j'_{1}}\\
& = & \mygens_{j'_{t}}\cdots\mygens_{j'_{y+1}}\mygens_{j'_{y-1}}\cdots\mygens_{j'_{r+1}}\mygens_{j'_{r}}\cdots\mygens_{j'_{z+1}}\mygens_{j'_{z-1}}\cdots\mygens_{j'_{1}}, 
\end{eqnarray*}
where the latter is an expression for $\tau$ with fewer than $t$ factors. 
This contradicts the fact that $\ell(\tau) = t$. 
We conclude that case (2) cannot occur. 

Our hypothesis for case (3) is that $i_{s} \not\in \{i'_{r+1},\ldots,i'_{s}\}$ and $i_{s} \not\in \{j'_{r+1},\ldots,j'_{t}\}$. 
Using reasoning from the previous two paragraphs, we see that $\mygens_{i_{s}}$ commutes with each $\mygens_{i'_{x}}$ and $\mygens_{j'_{y}}$ for $x \in [r+1,s]_{\mathbb{Z}}$ and $y \in [r+1,t]_{\mathbb{Z}}$ and that each of $\mu$, $\nu$, and $\lambda' := \mygens_{i'_{r}}\cdots\mygens_{i'_{1}}.\omega_{f}$ is below an edge of color $i_{s}$. 
As in the previous paragraph, we are forced to conclude that $i_{s} \in \{i'_{r},i'_{r-1},\ldots,i'_{1}\}$, so let $z \in [1,r]_{\mathbb{Z}}$ be largest such that $i_{s} = i'_{z}$, in which case $\mygens_{i_{s}}$ commutes with $\mygens_{i'_{r}},\ldots,\mygens_{i'_{z+1}}$. 
Then $\sigma$ can be expressed as $\mygens_{i_{s}}\mygens_{i'_{s}}\cdots\mygens_{i'_{r+1}}\mygens_{i'_{r}}\cdots\mygens_{i'_{z+1}}\mygens_{i'_{z-1}}\cdots\mygens_{i'_{1}}$ and $\tau$ can be expressed as $\mygens_{i_{s}}\mygens_{j'_{t}}\cdots\mygens_{j'_{r+1}}\mygens_{j'_{r}}\cdots\mygens_{j'_{z+1}}\mygens_{j'_{z-1}}\cdots\mygens_{j'_{1}}$. 
Let $\sigma' := \mygens_{i'_{s}}\cdots\mygens_{i'_{r+1}}\mygens_{i'_{r}}\cdots\mygens_{i'_{z+1}}\mygens_{i'_{z-1}}\cdots\mygens_{i'_{1}}$ and $\tau' := \mygens_{j'_{t}}\cdots\mygens_{j'_{r+1}}\mygens_{j'_{r}}\cdots\mygens_{j'_{z+1}}\mygens_{j'_{z-1}}\cdots\mygens_{j'_{1}}$, which have lengths $s-1$ and $t-1$ respectively. 
Then, $\tau'(\sigma')^{-1} = \mygens_{j'_{t}}\cdots\mygens_{j'_{r+1}}\mygens_{i'_{r+1}}\cdots\mygens_{i'_{s}}$, which we know is reduced. 
Moreover, we also know that $i'_{1}=j'_{1},\ldots,i'_{z-1}=j'_{z-1},i'_{z+1}=j'_{z+1},\ldots,i'_{r}=j'_{r}$. 
We also know that $\sigma' = \mygens_{i_{s}}\sigma = \mygens_{i_{s-1}}\cdots\mygens_{i_{r+1}}\mygens_{i_{r}}\cdots\mygens_{i_{1}}$ and $\tau' = \mygens_{i_{s}}\tau = \mygens_{j_{t}}\cdots\mygens_{j_{r+1}}\mygens_{j_{r}}\cdots\mygens_{j_{1}} = \mygens_{j_{t}}\cdots\mygens_{j_{r+1}}\mygens_{i_{s}}\mygens_{j_{r}}\cdots\mygens_{j_{1}}$. 

Continuing our analysis of case (3), take $\lambda := \mygens_{i_{r}}\cdots\mygens_{i_{1}}.\omega_{f} = \mygens_{j_{r+1}}\cdots\mygens_{j_{t}}.\nu$. 
Then $\langle \lambda,\alpha_{i_{s}}^{\vee} \rangle = \langle \nu+\sum_{y=r+1}^{t}\alpha_{j_{y}},\alpha_{i_{s}}^{\vee} \rangle = \langle \nu,\alpha_{i_{s}}^{\vee} \rangle + \sum_{y=r+1}^{t}\langle \alpha_{j_{y}},\alpha_{i_{s}}^{\vee} \rangle = -1$, since $\nu$ is below an edge of color $i_{s}$ and each $\gamma_{j_{y}}$ ($y \in [r+1,t]_{\mathbb{Z}}$) is distinct from and not adjacent to $\gamma_{i_{s}}$. 
That is, $\lambda$ is below an edge of color $i_{s}$. 
Then, $\lambda+\alpha_{i_{s}} = \mygens_{i_{s}}.\lambda = \mygens_{i_{s}}\mygens_{i_{r}}\cdots\mygens_{i_{1}}.\omega_{f}$. 
Again, as in the previous paragraph, we are forced to conclude that $i_{s} \in \{i_{r},\ldots,i_{1}\}$, so let $w \in [1,r]_{\mathbb{Z}}$ be largest such that $i_{s} = i_{w} = j_{w}$, in which case $\mygens_{i_{s}}$ commutes with $\mygens_{i_{r}},\ldots,\mygens_{i_{w+1}}$. 
Then we can see that $\tau'  =  \mygens_{j_{t}}\cdots\mygens_{j_{r+1}}\mygens_{j_{r}}\cdots\mygens_{j_{w+1}}\mygens_{j_{w-1}}\cdots\mygens_{j_{1}}$ and $\sigma' = \mygens_{i_{t}}\cdots\mygens_{i_{r+1}}\mygens_{i_{s}}\mygens_{i_{r}}\cdots\mygens_{i_{w+1}}\mygens_{i_{w-1}}\cdots\mygens_{i_{1}}$. 
So, $\mygens_{j_{t}}\cdots\mygens_{j_{r+1}}\mygens_{i_{s}}\mygens_{i_{r+1}}\cdots\mygens_{i_{s-1}}$ is another reduced expression for $\tau'(\sigma')^{-1}$. 
We may therefore invoke the inductive hypothesis to conclude that if we remove one instance of $i'_{z}$ from $\{i'_{1},\ldots,i'_{r}\}$ and remove one instance of $i_{w}$ from the multiset $\{i_{1},\ldots,i_{r}\}$, then we get an exact coincidence of multisets. 
Since $i_{w} = i_{s} = i'_{z}$, we therefore obtain the multiset equality $\{i_{1},\ldots,i_{r}\} \eqmulti \{i'_{1},\ldots,i'_{r}\}$. 
This completes our analysis of case (3) and the proof.\hfill\QED

\noindent 
{\bf \MeetJoinFollowUp}\ \ {\sl Take as given our minuscule fundamental weight $\omega_{f}$. 
Let $\sigma$ and  $\tau$ be in $\mathcal{W}^{\mysmallsetF}$ with lengths $s = \ell(\sigma)$ and $t = \ell(\tau)$, and let $\mu := \sigma.\omega_{f}$ and $\nu := \tau.\omega_{f}$.  
If $s$ and $t$ are both positive, then apply the algorithm from the proof of \MeetJoinMinusculeProp\ to get a number $r \in [1,\min(s,t)]_{\mathbb{Z}}$ and reduced expressions $\sigma = \mygens_{i_{s}}\cdots\mygens_{i_{1}}$ and $\tau = \mygens_{j_{t}}\cdots\mygens_{j_{1}}$ possessing properties (a) through (d) from the statement of that proposition, and let $\lambda^{\mathbf{join}}(\mu,\nu) = \lambda^{\mathbf{join}} := \mygens_{i_{r}}\cdots\mygens_{i_{1}}.\omega_{f} = \mygens_{j_{r}}\cdots\mygens_{j_{1}}.\omega_{f}$ and $\lambda^{\mathbf{meet}}(\mu,\nu) = \lambda^{\mathbf{meet}} := \mygens_{j_{t}}\cdots\mygens_{j_{r+1}}.\mu = \mygens_{i_{s}}\cdots\mygens_{i_{r+1}}.\nu$. 
If at least one of $s$ or $t$ is zero, then let $\lambda^{\mathbf{join}}(\mu,\nu) = \lambda^{\mathbf{join}} := \omega_{f}$ and let $\lambda^{\mathbf{meet}}(\mu,\nu) = \lambda^{\mathbf{meet}} := \tau.\mu = \sigma.\nu$. 
Then $\mu, \nu \in [\lambda^{\mathbf{meet}},\lambda^{\mathbf{join}}]_{\Pi(\omega_{f})} \subseteq [\lambda^{(0)},\lambda^{(1)}]_{\Pi(\omega_{f})}$ whenever $\lambda^{(0)}$ and $\lambda^{(1)}$ are weights in $\Pi(\omega_{f})$ such that $\mu, \nu \in [\lambda^{(0)},\lambda^{(1)}]_{\Pi(\omega_{f})}$.}  

{\em Proof.} We first address the case that at least one of $s$ or $t$ is zero. 
Then, at least one of $\mu$ or $\nu$ is $\omega_{f}$; without loss of generality, we take $\mu := \omega_{f}$. 
Under these assumptions, $\lambda^{\mathbf{join}} = \mu$, $\lambda^{\mathbf{meet}} = \nu$, and $\lambda^{\mathbf{meet}} \leq \lambda^{\mathbf{join}}$. 
Clearly $\mu, \nu \in [\lambda^{\mathbf{meet}},\lambda^{\mathbf{join}}]_{\Pi(\omega_{f})}$. 
Now suppose $\mu, \nu \in [\lambda^{(0)},\lambda^{(1)}]_{\Pi(\omega_{f})}$ for some weights $\lambda^{(0)}$ and $\lambda^{(1)}$ in $\Pi(\omega_{f})$. 
In particular, $\lambda^{(0)} \leq \nu$ and $\mu \leq \lambda^{(1)}$. 
That is, $\lambda^{(0)} \leq \lambda^{\mathbf{meet}}$ and $\lambda^{\mathbf{join}} \leq \lambda^{(1)}$, i.e.\ $[\lambda^{\mathbf{meet}},\lambda^{\mathbf{join}}]_{\Pi(\omega_{f})} \subseteq [\lambda^{(0)},\lambda^{(1)}]_{\Pi(\omega_{f})}$. 

Now assume $s$ and $t$ are both positive. 
From our given reduced expressions for $\sigma$ and $\tau$, we discern that for a Networked-numbers Game played from initial position $\omega_{f}$, each of the firing sequences $(\gamma_{i_{1}},\ldots,\gamma_{i_{r}},\gamma_{i_{r+1}},\ldots,\gamma_{i_{s}})$ and $(\gamma_{j_{1}},\ldots,\gamma_{j_{r}},\gamma_{j_{r+1}},\ldots,\gamma_{j_{t}})$ is legal. 
The firing sequence $(\gamma_{i_{1}},\ldots,\gamma_{i_{r}}) = (\gamma_{j_{1}},\ldots,\gamma_{j_{r}})$ takes us to the intermediate position $\lambda^{\mathbf{join}}$. 
From there, the firing sequence $(\gamma_{i_{r+1}},\ldots,\gamma_{i_{s}})$ takes us to $\mu$ and $(\gamma_{j_{r+1}},\ldots,\gamma_{j_{t}})$ takes us to $\nu$. 
So, $\mu = \lambda^{\mathbf{join}} - \sum_{x=r+1}^{s}\alpha_{i_{x}}$ and $\nu = \lambda^{\mathbf{join}} - \sum_{y=r+1}^{t}\alpha_{j_{y}}$, hence $\mu \leq \lambda^{\mathbf{join}} \geq \nu$. 
From the sixth paragraph of the proof of \MeetJoinMinusculeProp, we know that $\mygens_{j_{t}}\cdots\mygens_{j_{r+1}}\mygens_{i_{s}}\cdots\mygens_{i_{1}}$ and $\mygens_{i_{s}}\cdots\mygens_{i_{r+1}}\mygens_{j_{t}}\cdots\mygens_{j_{1}}$ are reduced expressions for the same element $\upsilon$ of $\mathcal{W}(\mathscr{G})$. 
Then, for a Networked-numbers Game played from $\omega_{f}$, the firing sequences $(\gamma_{i_{1}},\ldots,\gamma_{i_{s}},\gamma_{j_{r+1}},\ldots,\gamma_{j_{t}})$ and $(\gamma_{i_{1}},\ldots,\gamma_{i_{s}},\gamma_{j_{r+1}},\ldots,\gamma_{j_{t}})$ are both legal, and both take us to $\lambda^{\mathbf{meet}} := \upsilon.\omega_{f}$. 
So, $\lambda^{\mathbf{meet}} = \mu - \sum_{y=r+1}^{t}\alpha_{j_{y}} = \nu - \sum_{x=r+1}^{s}$, and hence $\mu \geq \lambda^{\mathbf{meet}} \leq \nu$. 
So, $\mu, \nu \in [\lambda^{\mathbf{meet}},\lambda^{\mathbf{join}}]_{\Pi(\omega_{f})}$. 

Now suppose $\mu$ and $\nu$ are in some interval $[\lambda^{(0)},\lambda^{(1)}]_{\Pi(\omega_{f})}$. 
So, we have paths $\mu \mylongarrow{i'_{s}} \cdots \mylongarrow{i'_{q+1}} \lambda^{(1)} \mylongbackarrow{j'_{p+1}} \cdots \mylongbackarrow {j'_{t}} \nu$ from $\mu$ and $\nu$ up to $\lambda^{(1)}$ together with a path $\lambda^{(1)} \mylongarrow{i'_{q}} \cdots \mylongarrow{i'_{1}} \omega_{f}$ from $\lambda^{(1)}$ up to $\omega_{f}$. 
We may therefore write $\sigma = \mygens_{i'_{s}}\cdots\mygens_{i'_{q+1}}\mygens_{i'_{q}}\cdots\mygens_{i'_{1}}$ and $\tau = \mygens_{j'_{t}}\cdots\mygens_{j'_{q+1}}\mygens_{j'_{q}}\cdots\mygens_{j'_{1}}$ with $i'_{1}=j'_{1},\ldots,i'_{q}=j'_{q}$. 
If we apply the algorithm from the third paragraph of the proof of \MeetJoinMinusculeProp\ to the expressions in the preceding sentence, we not only obtain expressions for $\sigma$ and $\tau$ matching the hypotheses of the proposition statement, but we also see that the multiset $\{i'_{1},\ldots,i'_{q}\}$ is a multisubset of $\{i_{1},\ldots,i_{r}\}$. 
It follows that $\lambda^{\mathbf{join}} \leq \lambda^{(1)}$. 

Using \MinusculeCharacterizationLemma, we can see that $-w_{0}.\omega_{f}$ is also a minuscule dominant weight, where $w_{0}$ is the unique longest element of the finite Weyl group $\mathcal{W}(\mathscr{G})$, cf.\ \S \FlowerSection.3. 
Observe that $-w_{0}.\mu \leq -w_{0}.\lambda^{(0)} \geq -w_{0}.\nu$ in $\Pi(-w_{0}.\omega_{f})$. 
Within this context, it is easy to see that the preceding paragraph applies to $-w_{0}.\lambda^{\mathbf{meet}}$ to yield $-w_{0}.\lambda^{\mathbf{meet}} \leq -w_{0}.\lambda^{(0)}$ in $\Pi(-w_{0}.\omega_{f})$ and therefore $\lambda^{\mathbf{meet}} \geq \lambda^{(0)}$ in $\Pi(\omega_{f})$.\hfill\QED

{\bf [\S \MinusculeExampleSection.6:\! The saturated set of weights $\Pi(\omega_{f})$ as a DCDL.]} 
\MeetJoinMinusculeProp\ demonstrates that $\Pi(\omega_{f})$ is a lattice by precisely locating meets and joins. 
We now show that $\Pi(\omega_{f})$ is a DCDL by proving that, as a lattice, $\Pi(\omega_{f})$ is cancellative (cf.\ \CancellationTheorem). 

{\bf \MinusculeDCDL}\ \ {\sl For our given minuscule fundamental weight $\omega_{f}$, $\Pi(\omega_{f})$ is a cancellative lattice with meets and joins given by $\mu \wedge \nu = \lambda^{\mathbf{meet}}(\mu,\nu)$ and $ \mu \vee \nu = \lambda^{\mathbf{join}}(\mu,\nu)$, for all $\mu, \nu \in \Pi(\omega_{f})$, as in \MeetJoinMinusculeProp. 
In particular, $\Pi(\omega_{f})$ is a $\mathscr{G}$-structured diamond-colored distributive lattice.}

{\em Proof.} That $\Pi(\omega_{f})$ is $\mathscr{G}$-structured is established in \MinusculeStructureLemma. 
Assume that $\mu = \sigma.\omega_{f}$ and $\nu = \tau.\omega_{f}$ with $\sigma$ and $\tau$ as in \MeetJoinMinusculeProp. 
By part {\sl (5)} of that proposition, $\lambda^{\mathbf{join}}(\mu,\nu)$ is the join of $\mu$ and $\nu$ and $\lambda^{\mathbf{meet}}(\mu,\nu)$ is their meet. 
So, once we show that the lattice $\Pi(\omega_{f})$ is cancellative, then, since it is diamond-colored, we can conclude it is a $\mathscr{G}$-structured DCDL. 

Suppose $\xi \in \Pi(\omega_{f})$ has the property that $\mu \vee \xi = \nu \vee \xi$ and $\mu \wedge \xi = \nu \wedge \xi$. 
We wish to show that $\mu = \nu$. 
Write $\xi = \mygreeku.\omega_{f}$ where $\mygreeku = \mygens_{k_{u}}\cdots\mygens_{k_{1}}$ is a length $u$ element of $\mathcal{W}^{\mysmallsetF}$ such that $\mu \vee \xi = \mygens_{i_{p}}\cdots\mygens_{i_{1}}.\omega_{f}$ with $i_{1}=k_{1},\ldots,i_{p}=k_{p}$ and $\ell(\mygens_{k_{u}}\cdots\mygens_{k_{p+1}}\mygens_{i_{p+1}}\cdots\mygens_{i_{s}}) = (u-p)+(s-p)$. 
We can also write $\mygreeku = \mygens_{k'_{u}}\cdots\mygens_{k'_{1}}$ such that $\nu \vee \xi = \mygens_{j_{q}}\cdots\mygens_{j_{1}}.\omega_{f}$ with $j_{1}=k'_{1},\ldots,j_{q}=k'_{q}$ and $\ell(\mygens_{k'_{u}}\cdots\mygens_{k'_{q+1}}\mygens_{j_{q+1}}\cdots\mygens_{j_{t}}) = (u-q)+(t-q)$. 
Now, the expression $\mygens_{i_{p}}\cdots\mygens_{i_{1}}$ is reduced since it is a suffix of this reduced expression $\mygens_{i_{s}}\cdots\mygens_{i_{p+1}}\mygens_{i_{p}}\cdots\mygens_{i_{1}}$ for $\sigma$; since $\sigma \in \mathcal{W}^{\mysmallsetF}$, it follows from the definitions that the suffix $\mygens_{i_{p}}\cdots\mygens_{i_{1}}$ is also in $\mathcal{W}^{\mysmallsetF}$. 
Similarly see that $\mygens_{j_{q}}\cdots\mygens_{j_{1}}$ is a reduced expression for an element of $\mathcal{W}^{\mysmallsetF}$. 
These latter facts, together with the equality $\mygens_{i_{p}}\cdots\mygens_{i_{1}}.\omega_{f} = \mu \vee \xi = \nu \vee \xi = \mygens_{j_{q}}\cdots\mygens_{j_{1}}.\omega_{f}$, imply that $\mygens_{i_{p}}\cdots\mygens_{i_{1}} = \mygens_{j_{q}}\cdots\mygens_{j_{1}}$ and that $p = q$. 
Without loss of generality, we may assume that $k_{1}=i_{1}=j_{1}=k'_{1},\ldots,k_{p}=i_{p}=j_{p}=k'_{p}$. 
That is, $\xi = \mygens_{k_{u}}\cdots\mygens_{k_{p+1}}\mygens_{k_{p}}\cdots\mygens_{k_{1}}.\omega_{f} = \mygens_{k_{u}}\cdots\mygens_{k_{p+1}}.(\mu \vee \xi) = \mu \vee \xi - (\alpha_{k_{p+1}}+\cdots+\alpha_{k_{u}})$ while  $\xi = \mygens_{k'_{u}}\cdots\mygens_{k'_{p+1}}\mygens_{k_{p}}\cdots\mygens_{k_{1}}.\omega_{f} = \mygens_{k'_{u}}\cdots\mygens_{k'_{p+1}}.(\nu \vee \xi) = \nu \vee \xi - (\alpha_{k'_{p+1}}+\cdots+\alpha_{k'_{u}})$. 
Therefore, $\{k_{p+1},\ldots,k_{u}\} \eqmulti \{k'_{p+1},\ldots,k'_{u}\}$. 
Again, without losing generality, we can take $k'_{p+1}=k_{p+1},\ldots,k'_{u}=k_{u}$. 
Now, $\mu - (\alpha_{k_{p+1}} + \cdots + \alpha_{k_{u}}) = \mygens_{k_{u}}\cdots\mygens_{k_{p+1}}.\mu = \mu \wedge \xi = \nu \wedge \xi = \mygens_{k_{u}}\cdots\mygens_{k_{p+1}}.\nu = \nu - (\alpha_{k_{p+1}} + \cdots + \alpha_{k_{u}})$. 
Therefore $\mu = \nu$.\hfill\QED 

From here on, we refer to the $\mathscr{G}$-structured DCDL $\Pi(\omega_{f})$ by the notation $L_{\mathscr{G}}(\omega_{f})$ and call it a {\em $\mathscr{G}$-minuscule splitting DCDL}, or just a {\em minuscule splitting DCDL} when the Coxeter--Dynkin flower $\mathscr{G}$ is clear from context. 
Its compression poset $P_{\mathscr{G}}(\omega_{f}) := \jcolor(L_{\mathscr{G}}(\omega_{f}))$, to be called the associated {\em $\mathscr{G}$-minuscule compression poset} (or simply {\em minuscule compression poset}), is therefore $\mathscr{G}$-structured in the vertex-colored sense.  
Now, the finite Weyl group $\mathcal{W}(\mathscr{G})$ for our given Coxeter--Dynkin flower $\mathscr{G}$ has longest element denoted $w_{0}$ wherein $-w_{0}.\omega_{k} = \omega_{\sigma_{0}(k)}$ for the involution $\sigma_{0}: I \longrightarrow I$ of the color palette, cf.\ \S \FlowerSection.3.
Then, by \MinusculeCharacterizationLemma.1/2, the fundamental weight $-w_{0}.\omega_{f} = \omega_{\sigma_{0}(f)}$ is also a minuscule fundamental weight. 

\noindent
{\bf \StarBowtieLemma}\ \ {\sl Take as given our minuscule fundamental weight $\omega_{f}$. 
Then $L_{\mathscr{G}}(\omega_{f})^{*}$ is isomorphic to the minuscule splitting DCDL $L_{\mathscr{G}}(\omega_{\sigma_{0}(f)}) = L_{\mathscr{G}}(-w_{0}.\omega_{f})$, and $L_{\mathscr{G}}(\omega_{f})^{\bowtie}$ is isomorphic to $L_{\mathscr{G}}(\omega_{f})$. 
Moreover, $P_{\mathscr{G}}(\omega_{f})^{*} \cong P_{\mathscr{G}}(\omega_{\sigma_{0}(f)})$ is an isomorphism of vertex-colored posets, as is $P_{\mathscr{G}}(\omega_{f})^{\bowtie} \cong P_{\mathscr{G}}(\omega_{f})$. 
The minuscule compression poset $P_{\mathscr{G}}(\omega_{f})$ has a unique maximal element and a unique minimal element, and the colors of these elements are $f$ and $\sigma_{0}(f)$ respectively.} 

{\em Proof.} The claims of the second sentence follow from \SaturatedLemma.3. 
Then the claims of the third sentence follow from \JMCorollary.

Now, by \SaturatedLemma.3, $\omega_{f}$ is the unique maximal element of $L_{\mathscr{G}}(\omega_{f}) = \Pi(\omega_{f}) = \mathcal{W}\omega_{f}$. 
If $\mu \myarrow{i} \omega_{f}$ in $L_{\mathscr{G}}(\omega_{f})$, then within the saturated set of weights $\Pi(\omega_{f})$, the $i$-component of $\omega_{f}$ is a saturated chain of color $i$ edges $\mu_{0} \myarrow{i} \mu_{1} \myarrow{i} \cdots \myarrow{i} \mu_{p} = \mu \myarrow{i} \omega_{f}$, where $\mu_{0} = \mygens_{i}.\omega_{f}$. 
But by \MinusculeStructureLemma, such a chain can have at most two elements, so $p = 0$; we must also have $i = f$ since otherwise $\mygens_{i}$ stabilizes $\omega_{f}$. 
Therefore $\mu = \mygens_{f}.\omega_{f} = \omega_{f}-\alpha_{f}$. 
In particular, the only edge below $\omega_{f}$ is $(\omega_{f}-\alpha_{f}) \myarrow{f} \omega_{f}$, making $\omega_{f}$ a join irreducible element of $L_{\mathscr{G}}(\omega_{f})$. 
Within $P_{\mathscr{G}}(\omega_{f})$ this join irreducible has color $f$. 

Set $g := \sigma_{0}(f)$. 
By the same reasoning of the previous paragraph, in $P_{\mathscr{G}}(\omega_{g}) \cong P_{\mathscr{G}}(\omega_{f})^{*}$, the unique maximal element has color $g$, so in $P_{\mathscr{G}}(\omega_{f}) \cong \big(P_{\mathscr{G}}(\omega_{f})^{*}\big)^{*}$ the unique minimal element has color $g$. 
Alternatively, we can argue directly that $-\omega_{g} \mylongarrow{g} -\omega_{g}+\alpha_{g}$ is the only edge above $-\omega_{g}$ in $L_{\mathscr{G}}(\omega_{f})$, making $-\omega_{g}+\alpha_{g}$ a color $g$ join irreducible element and also the unique minimal element of $P_{\mathscr{G}}(\omega_{f})$.\hfill\QED 

{\bf [\S \MinusculeExampleSection.7:\! Properties of certain subposets of minuscule compression posets.]} 
Not surprisingly, $\mathscr{G}$-minuscule compression posets have a number of distinctive properties. 
Our list of these begins with the next (well-known) result, which demonstrates that for two colors of neighboring elements in the vascular graph $\Gamma$, the induced-order subposet of $P_{\mathscr{G}}(\omega_{f})$ comprised of all elements of these two colors is totally ordered and also demonstrates that two elements with the same color cannot be neighbors in $P_{\mathscr{G}}(\omega_{f})$. 

Before that, some bookkeeping: From here on, we use Greek letters such as $\lambda, \mu, \nu$ (representing weights) as well as boldface letters like $\uelt, \velt, \welt$ (representing order ideals from $P_{\mathscr{G}}(\omega_{f})$) to refer to elements of the minuscule splitting DCDL $\Pi(\omega_{f}) \cong L_{\mathscr{G}}(\omega_{f}) \cong \Jcolor(P_{\mathscr{G}}(\omega_{f}))$. 
Similarly, elements of $P_{\mathscr{G}}(\omega_{f})$ might be denoted as boldface letters like $\uelt, \velt, \welt$ (when thought of as join- or meet-irreducible elements in $L_{\mathscr{G}}(\omega_{f})$) or denoted as lower-case letters like $u, v, w$ (when we are thinking of $P_{\mathscr{G}}(\omega_{f})$ independently of its associated minuscule splitting DCDL). 
Also, take $\mathscr{E}: \EdgeSet(L_{\mathscr{G}}(\omega_{f})) \longrightarrow I$ to be the edge-coloring function for our generically given minuscule splitting DCDL and $\mathscr{V}: \VertexSet(P_{\mathscr{G}}(\omega_{f})) \longrightarrow I$ to be the vertex-coloring function for the associated minuscule compression poset. 

\noindent 
{\bf \ijChainLemma}\ \ {\sl Let $P = P_{\mathscr{G}}(\omega_{f})$ be our given minuscule compression poset. 
Let $u$ and $v$ be elements of $P$ with vertex colors $i := \vcolor(u)$ and $j := \vcolor(v)$. 
(1) If $u \rightarrow v$ in $P$, then $\gamma_{i}$ and $\gamma_{j}$ are distinct and adjacent nodes in $\Gamma$. 
(2) If $u$ and $v$ are incomparable elements of $P$, then nodes $\gamma_{i}$ and $\gamma_{j}$ of $\Gamma$ are distinct and nonadjacent; in particular, when $k$ and $l$ are distinct colors from our color palette $I$ such that $\gamma_{k}$ and $\gamma_{l}$ are adjacent in $\Gamma$, then $\vcolor^{-1}(\{k,l\})$ is a chain in $P$.} 

{\em Proof.} The ``in particular'' claim in part {\em (2)} is the contrapositive of the conditional statement that precedes it. 
To prove said conditional, let $L := L_{\mathscr{G}}(\omega_{f})$. 
Let $\uelt := (-\infty,u]_{P}$ and $\velt := (-\infty,v]_{P}$ be the down-sets from $P$ generated by elements $u$ and $v$ respectively. 
Since $u$ and $v$ are incomparable in $P$, then we have distinct edges $\uelt \vee \velt \setminus \{u\} \myarrow{i} \uelt \vee \velt$ and $\uelt \vee \velt \setminus \{v\} \myarrow{j} \uelt \vee \velt$ in $L$ below $\uelt \vee \velt$. 
By \MinusculeStructureLemma, $\gamma_{i}$ and $\gamma_{j}$ are distinct and nonadjacent in $\Gamma$. 

For {\sl (1)}, say we are given an edge $u \rightarrow v$ in $P$ with $i := \vcolor(u)$ and $j := \vcolor(v)$. 
We claim that $(-\infty,v] \setminus \{v\}$ and $(-\infty,v] \setminus \{v,u\}$ are both down-sets from $P$. 
First, suppose $y \in (-\infty,v] \setminus \{v\}$ with $x \leq y$. 
Of course, $y < v$, and since $y$ is in the down-set $(-\infty,v]$, then $x \in (-\infty,v]$. 
Clearly $x \not= v$, else we get $x \leq y < x$, which is absurd. 
So $x \in (-\infty,v] \setminus \{v\}$ as well. 
Second, suppose $y \in (-\infty,v] \setminus \{u,v\}$ with $x \leq y$. 
As before, we must have $y < v$, $x \in (-\infty,v]$, and $x \not= v$. 
Consider the possibility that $x = u$. 
Since $y \not= u$, then $x < y < v$, meaning $u < y < v$. 
But this violates the fact that $v$ covers $u$. 
So we must have $x \not= u$ as well. 
Therefore, $x \in (-\infty,v] \setminus \{u,v\}$ as well. 
Let $\selt := (-\infty,v] \setminus \{u,v\}$, $\telt := (-\infty,v] \setminus \{v\}$, and $\velt := (-\infty,v]$, which we now know are three down-sets from $P$. 
Then $\selt \myarrow{i} \telt \myarrow{j} \velt$ in $L$. 
Clearly $i$ and $j$ must be distinct, else the $i$-component of $\velt$ in $L$ is a chain with at least three elements, in violation of \MinusculeStructureLemma. 
Moreover, since $\velt \setminus \{u\}$ is not a down-set from $P$, then there can be no $\uelt \in L$ such that $\selt \myarrow{j} \uelt \myarrow{i} \velt$. 
By \MinusculeStructureLemma, it cannot be the case that the distinct nodes $\gamma_{i}$ and $\gamma_{j}$ are nonadjacent in $\Gamma$. 
\hfill\QED

The next result states how we can naturally match certain lower order ideals from one $\mathscr{G}$-minuscule compression poset to a certain upper order ideals from a second $\mathscr{G}$-minuscule compression poset. 
The result appears to be new. 

\noindent 
{\bf \IdealMatchingProp}\ \ {\sl With respect to our given Coxeter--Dynkin flower $\mathscr{G}$, let $\omega_{f}$ and $\omega_{g}$ be any two $\mathscr{G}$-minuscule fundamental weights. 
For brevity, set $P_{f} := P_{\mathscr{G}}(\omega_{f})$, $L_{f} := L_{\mathscr{G}}(\omega_{f})$, $P_{g} := P_{\mathscr{G}}(\omega_{g})$, and $L_{g} := L_{\mathscr{G}}(\omega_{g})$. 
Then $P_{f}$ has a color $g$ vertex and $P_{g}$ has a color $f$ vertex.} 
{\sl Next, let $\tilde{f} := \sigma_{0}(f)$, so $\omega_{\tilde{f}}$ is the minuscule fundamental weight $-w_{0}.\omega_{f}$, cf.\ \StarBowtieLemma., and similarly let $\tilde{g} := \sigma_{0}(g)$.}
{\sl Suppose $t \in P_{f}$ has color $g$, and let $\telt := (-\infty,t]_{P_{f}}$ be the down-set from $P_{f}$ generated by $t$.} 
{\sl Then there exists a unique color $\tilde{f}$ element $b \in P_{g}$ with associated up-set $\belt := [b,\infty)_{P_{g}}$ such that the intervals} $[\mymin(L_{f}),\telt]_{L_{f}}$ {\sl and} $[\belt,\mymax(L_{g})]_{L_{g}}$ {\sl are isomorphic as edge-colored posets and hence as DCDL's. 
Moreover, $(-\infty,t]_{P_{f}}$ and $[b,\infty)_{P_{g}}$ are isomorphic vertex-colored posets. 
Similarly, for every color $\tilde{f}$ element $b' \in P_{g}$ with associated up-set $\belt' := [b',\infty)_{P_{g}}$, there exists a unique color $g$ element $t' \in P_{f}$ with associated down-set $\telt' := (-\infty,t']_{P_{f}}$ such that} $[\mymin(L_{f}),\telt']_{L_{f}} \cong [\belt',\mymax(L_{g})]_{L_{g}}$ {\sl and} $(-\infty,t']_{P_{f}} \cong [b',\infty)_{P_{g}}$. 

{\em Proof.} 
By \NewLFKTheorem.2, connectedness of the vascular graph $\Gamma(\mathscr{G})$ associated with our Coxeter--Dynkin flower $\mathscr{G}$ implies that $L_{f}$ has an edge of color $g$ and therefore its associated compression poset $P_{f}$ must have a vertex of color $g$; similarly see that $P_{g}$ has a vertex of color $f$. 
We will argue that, given $t \in P_{f}$ with color $g$, there exists a unique color $\tilde{f} := \sigma_{0}(f)$ element $b \in P_{g}$ such that $[\mymin(L_{f}),\telt]_{L_{f}} \cong [\belt,\mymax(L_{g})]_{L_{g}}$. 
In this case, we automatically get the isomorphisms $(-\infty,t]_{P_{f}} \cong \jcolor([\mymin(L_{f}),\telt]) \cong \mcolor([\belt,\mymax(L_{g})]) \cong [b,\infty)_{P_{g}}$. 
The uniqueness of $b$ follows from the fact that the color $\tilde{f}$ elements of $P_{g}$ are totally ordered; in particular, for any $b' \not= b$ with color $\tilde{f}$, the up-set $\belt' := [b',\infty)_{P_{g}}$ will not have the same distance from $\mymax(L_{g})$ as $\belt$ and therefore cannot serve as the minimal element for the interval we seek in $L_{g}$. 
Moreover, the claim ``for every color $\tilde{f}$ element $b' \in P_{g}$ there exists a unique color $g$ element $t' \in P_{f}$ such that $[\mymin(L_{f}),\telt']_{L_{f}} \cong [\belt',\mymax(L_{g})]_{L_{g}}$'' will follow by entirely similar reasoning. 

So, take $t \in P_{f}$ with color $g$ and $\telt := (-\infty,t]_{P_{f}}$ as the down-set from $P_{f}$ generated by $t$. 
We propose a mapping $\psi: [\mymin(L_{f}),\telt] \longrightarrow L_{g}$ as follows. 
Suppose $\xelt \in [\mymin(L_{f}),\telt]$ such that $\xelt = \xelt_{p} \myarrow{i_{p}} \xelt_{p-1} \myarrow{i_{p-1}} \cdots \myarrow{i_{1}} \xelt_{0} = \telt$ is a path in $L_{f}$ from $\xelt$ up to $\telt$. 
Of course, $i_{1} = g$. 
We wish to show that there is a unique $\yelt \in L_{g}$ such that $\mymax(L_{g})=\yelt_{0} \mybackarrow{i_{1}} \cdots \mybackarrow{i_{p}} \yelt_{p} = \yelt$ is a path in $L_{g}$ from $\yelt$ up to $\mymax(L_{g})$. 

For the moment, assume we can argue that for any such path from $\xelt$ up to $\telt$ in $L_{f}$, there is a corresponding path from some $\yelt \in L_{g}$ up to $\mymax(L_{g})$. 
We argue first that $\yelt$ is uniquely identified by this correspondence. 
Since $L_{f}$ is ranked, any other path in $L_{f}$ from $\xelt$ up to $\telt$ must have length $p$, so suppose $\xelt = \xelt'_{p} \myarrow{j_{p}} \xelt'_{p-1} \myarrow{j_{p-1}} \cdots \myarrow{j_{1}} \xelt'_{0} = \telt$ is another such path, and let $\mymax(L_{g})=\yelt'_{0} \mybackarrow{j_{1}} \cdots \mybackarrow{j_{p}} \yelt'_{p} = \yelt'$ be the corresponding path from some $\yelt'$ up to $\mymax(L_{g})$ in $L_{g}$. 
Now $wt(\telt)-wt(\xelt) = \sum_{q=1}^{p}\alpha_{i_{q}} = \sum_{q=1}^{p}\alpha_{j_{q}}$, forcing the equality $\{i_{1},\ldots,i_{p}\} \eqmulti \{j_{1},\ldots,j_{p}\}$. 
In particular, $wt(\yelt) = \omega_{g} - \sum_{q=1}^{p}\alpha_{i_{q}} = \omega_{g} - \sum_{q=1}^{p}\alpha_{j_{q}} = wt(\yelt')$. 
Since any element of $L_{g} = \Pi(\omega_{g})$ is uniquely identified by its weight, then we must have $\yelt'=\yelt$. 

Next, we confirm for any such path from $\xelt$ up to $\telt$ in $L_{f}$, there indeed exists a corresponding path from some $\yelt \in L_{g}$ up to $\mymax(L_{g})$. 
To do so, we induct on $p$. 
When $p=0$, then $\xelt = \telt$ and the corresponding path in $L_{g}$ has $\yelt := \mymax(L_{g})$. 
In what follows next, we use the following notation: If $\pelt$ is an intermediate position in some Networked-numbers Game played on $\mathscr{G}$, then let $\gamma_{i}(\pelt) := \langle \pelt,\alpha_{i}^{\vee} \rangle$ denote the number at node $\gamma_{i}$; by \MinusculeCharacterizationLemma, if $\pelt$ is an intermediate position for a game played from an initial position that is a minuscule fundamental weight, we must have $|\gamma_{i}(\pelt)| \leq 1$. 
Suppose now that for some $p \geq 0$, it is the case that if $0 \leq q \leq p$ and $\xelt = \xelt_{q} \myarrow{i_{q}} \cdots \myarrow{i_{1}} \xelt_{0}=\telt$ is a path from some $\xelt$ up to $\telt$ in $L_{f}$, then there exists some $\yelt \in L_{g}$ such that $\mymax(L_{g})=\yelt_{0} \mybackarrow{i_{1}} \cdots \mybackarrow{i_{q}} \yelt_{q} = \yelt$ is a path in $L_{g}$ from $\yelt$ up to $\mymax(L_{g})$ with $\gamma_{i}(\xelt) \leq \gamma_{i}(\yelt)$ for all $i \in I$. 
Now, say $\xelt = \xelt_{p+1} \mylongarrow{i_{p+1}} \xelt_{p} \myarrow{i_{p}} \cdots \myarrow{i_{1}} \xelt_{0}=\telt$ is a path in $L_{f}$ from some $\xelt$ up to $\telt$. 
So, there is some $\yelt_{p} \in L_{g}$ such that $\mymax(L_{g})=\yelt_{0} \mybackarrow{i_{1}} \cdots \mybackarrow{i_{p}} \yelt_{p}$ is a path in $L_{g}$ from $\yelt_{p}$ up to $\mymax(L_{g})$ with $\gamma_{i}(\xelt_{p}) \leq \gamma_{i}(\yelt_{p})$ for all $i \in I$. 
In particular, $\gamma_{i_{p+1}}(\xelt_{p}) \leq \gamma_{i_{p+1}}(\yelt_{p})$. 
Since $\xelt_{p+1} \mylongarrow{i_{p+1}} \xelt_{p}$ in $L_{f}$, then $\gamma_{i_{p+1}}(\xelt_{p}) = 1$ which forces $\gamma_{i_{p+1}}(\yelt_{p}) = 1$. 
Therefore we can legally fire node $\gamma_{i_{p+1}}$ from position $\yelt_{p}$ to obtain a position $\yelt =: \yelt_{p+1}$. 
We have $\yelt = \mygens_{i_{p+1}}.\yelt_{p} \in \mathcal{W}\omega_{g} = \Pi(\omega_{g}) = L_{g}$ with $\yelt  \mylongarrow{i_{p+1}} \yelt_{p}$. 
So, $\yelt$ is the unique element of $L_{g}$ for which there is a path $\mymax(L_{g})=\yelt_{0} \mybackarrow{i_{1}} \cdots \mybackarrow{i_{p}} \yelt_{p} \mylongbackarrow{i_{p+1}} \yelt_{p+1} = \yelt$. 
Now, if $\gamma_{j}$ is distinct from and nonadjacent to $\gamma_{i_{p+1}}$, then $\gamma_{j}(\xelt) = \gamma_{j}(\xelt_{p}) \leq \gamma_{j}(\yelt_{p}) = \gamma_{j}(\yelt)$. 
On the other hand, if $\gamma_{j}$ is distinct from and adjacent to $\gamma_{i_{p+1}}$, then $\gamma_{j}(\xelt) = \gamma_{j}(\xelt_{p})+|M_{i_{p+1},j}| \cdot 1 \leq \gamma_{j}(\yelt_{p})+|M_{i_{p+1},j}| \cdot 1 = \gamma_{j}(\yelt)$. 
And, $\gamma_{i_{p+1}}(\xelt) = -1 = \gamma_{i_{p+1}}(\yelt)$. 
This confirms that $\gamma_{i}(\xelt) \leq \gamma_{i}(\yelt)$ for all $i \in I$ and completes the induction argument.

Therefore, we have a well-defined function $\psi: [\mymin(L_{f}),\telt] \longrightarrow L_{g}$. 
If $\yelt_{1} := \psi(\xelt_{1}) = \psi(\xelt_{2}) =: \yelt_{2}$, then $wt(\yelt_{1}) = \omega_{g} - \sum_{q=1}^{p}\alpha_{i_{q}} = \omega_{g} - \sum_{q=1}^{r}\alpha_{j_{q}} = wt(\yelt_{2})$, where $\xelt_{1} \myarrow{i_{p}} \cdots \myarrow{i_{1}} \telt$ and $\xelt_{2} \myarrow{j_{r}} \cdots \myarrow{j_{1}} \telt$. 
So,  $\{i_{1},\ldots,i_{p}\} \eqmulti \{j_{1},\ldots,j_{r}\}$, forcing $r = p$. 
Then, $wt(\xelt_{1}) = wt(\telt) - \sum_{q=1}^{p}\alpha_{i_{q}} = wt(\telt) - \sum_{q=1}^{p}\alpha_{j_{q}} = wt(\xelt_{2})$.  Since any element of $L_{f} = \Pi(\omega_{f})$ is uniquely identified by its weight, then $\xelt_{1} = \xelt_{2}$. 
That is, $\psi$ is injective. 
Observe that by our definition of $\psi$, we have $\xelt_{1} \myarrow{i} \xelt_{2}$ in $[\mymin(L_{f}),\telt]$ if and only if $\psi(\xelt_{1}) \myarrow{i} \psi(\xelt_{2})$ in $L_{g}$. 

Let $\belt := \psi(\mymin(L_{f}))$, and fix corresponding paths $\mymin(L_{f}) \myarrow{i_{r}} \cdots \myarrow{i_{1}} \telt$ and $\belt \myarrow{i_{r}} \cdots \myarrow{i_{1}} \mymax(L_{g})$. 
We aim to show that $\psi([\mymin(L_{f}),\telt]) = [\belt,\mymax(L_{g})]$. 
For any $\xelt \in [\mymin(L_{f}),\telt]$, there is a path $\mymin(L_{f}) \myarrow{j_{r}} \cdots \mylongarrow{j_{p+1}} \xelt \myarrow{j_{p}} \cdots \myarrow{j_{1}} \telt$ and hence a path $\belt \myarrow{j_{r}} \cdots \mylongarrow{j_{p+1}} \psi(\xelt) \myarrow{j_{p}} \cdots \myarrow{j_{1}} \mymax(L_{g})$. 
In particular, $\psi(\xelt) \in [\belt,\mymax(L_{g})]$. 

Suppose the containment $\psi([\mymin(L_{f}),\telt]) \subseteq [\belt,\mymax(L_{g})]$ is proper. 
Then there is some positive number $p$ and some $\yelt_{p} \in [\belt,\mymax(L_{g})]$ with a path $\belt = \yelt_{r} \myarrow{j_{r}} \cdots \mylongarrow{j_{p+2}} \yelt_{p+1} \mylongarrow{j_{p+1}} \yelt_{p}$ such that $\yelt_{q} \in \psi([\mymin(L_{f}),\telt])$ for $q \in [p+1,r]_{\mathbb{Z}}$ but $\yelt_{p} \not\in \psi([\mymin(L_{f}),\telt])$. 
(For such a path, if $p=0$, then necessarily $\yelt_{0} = \mymax(L_{g})$, which is why the smallest such number $p$ in the preceding sentence must be positive.)
So, there exists a path $\mymin(L_{f}) = \xelt_{r} \myarrow{j_{r}} \cdots \mylongarrow{j_{p+2}} \xelt_{p+1}$ in $L_{f}$ with $\yelt_{q} = \psi(\xelt_{q})$ for each $q \in [p+1,r]_{\mathbb{Z}}$. 
We know from our earlier induction argument that $\gamma_{i}(\xelt_{q}) \leq \gamma_{i}(\yelt_{q})$ for all appropriate $q$. 
Also, since $\yelt_{p+1} \mylongarrow{j_{p+1}} \yelt_{p}$, then $\gamma_{j_{p+1}}(\yelt_{p+1}) = -1$, which forces $\gamma_{j_{p+1}}(\xelt_{p+1}) = -1$. 
Therefore, there is some $\xelt_{p} \in L_{f}$ such that $\xelt_{p+1} \mylongarrow{j_{p+1}} \xelt_{p}$.
 Now, if $\xelt_{p} \in [\mymin(L_{f}),\telt]$, then we would have $\psi(\xelt_{p}) = \yelt_{p}$, but this would violate the fact that $\yelt_{p} \not\in \psi([\mymin(L_{f}),\telt])$. 
So, $\xelt_{p} \not\in [\mymin(L_{f}),\telt]$. 
That said, we must have $\xelt_{p+1} \mylongarrow{j'_{p+1}} \xelt'_{p}$ for some $\xelt'_{p} \in [\mymin(L_{f}),\telt]$, otherwise we have $\xelt_{p+1} = \telt$ and $\yelt_{p+1} = \psi(\xelt_{p+1}) = \psi(\telt) = \mymax(L_{g})$, contradicting the fact that $\yelt_{p+1}$ is below $\yelt_{p}$. 
By \MinusculeStructureLemma, $\gamma_{j_{p+1}}$ and $\gamma_{j'_{p+1}}$ are distinct and nonadjacent in $\Gamma$. 
Let $\xelt_{p-1} := \xelt_{p} \vee \xelt'_{p}$, which cannot be in $[\mymin(L_{f}),\telt]$ else $\xelt_{p}$ is also. 
We have $\xelt_{p} \mylongarrow{j_{p+1}} \xelt_{p-1} \mylongbackarrow{j'_{p+1}} \xelt'_{p}$. 
Let $\yelt'_{p} := \psi(\xelt'_{p})$, so $\yelt_{p+1} \mylongarrow{j'_{p+1}} \yelt'_{p}$; let $\yelt_{p-1} := \yelt_{p} \vee \yelt'_{p}$, so $\yelt_{p} \mylongarrow{j_{p+1}} \yelt_{p-1} \mylongbackarrow{j'_{p+1}} \yelt'_{p}$. 
There must be some $\xelt'_{p-1} \in [\mymin(L_{f}),\telt]$ with $\xelt'_{p} \myarrow{j'_{p}} \xelt'_{p-1}$, otherwise we have $\xelt'_{p} = \telt$ and $\yelt'_{p} = \psi(\xelt'_{p}) = \psi(\telt) = \mymax(L_{g})$, contradicting the fact that $\yelt'_{p}$ is below $\yelt_{p-1}$. 
Continuing in this way, we obtain, for any $p \geq q \geq 1$, the following diamonds in $L_{f}$ and $L_{g}$ respectively: 

\hspace*{1in}
\parbox{2.7cm}{\begin{center}
\setlength{\unitlength}{0.4cm}
\begin{picture}(6.5,5.5)
\put(3,1){\circle*{0.35}} 
\put(1,3){\circle*{0.35}}
\put(3,5){\circle*{0.35}} 
\put(5,3){\circle*{0.35}}
\put(1,3){\line(1,1){2}} 
\put(3,1){\line(-1,1){2}}
\put(5,3){\line(-1,1){2}} 
\put(3,1){\line(1,1){2}}
\put(1,3){\line(-1,-1){1}}
\multiput(0.4,2.7)(-0.2,-0.2){6}{\color{white}{\rule{3mm}{0.4mm}}}
\put(3,1){\line(-1,-1){1}}
\multiput(2.4,0.7)(-0.2,-0.2){6}{\color{white}{\rule{3mm}{0.4mm}}}
\put(5,3){\line(1,1){1}}
\multiput(5.1,3.2)(0.2,0.2){6}{\color{white}{\rule{3mm}{0.4mm}}}
\put(3,5){\line(1,1){1}}
\multiput(3.1,5.2)(0.2,0.2){6}{\color{white}{\rule{3mm}{0.4mm}}}
\put(1.7,1.75){\scriptsize $j_{p+1}$} 
\put(3.55,1.75){\scriptsize $j'_{q+1}$}
\put(1.55,3.75){\scriptsize $j'_{q+1}$} 
\put(3.7,3.75){\scriptsize $j_{p+1}$}
\put(3.5,0.5){\footnotesize $\xelt'_{q+1}$} 
\put(5.5,2.75){\footnotesize $\xelt'_{q}$}
\put(1.25,5.25){\footnotesize $\xelt_{q-1}$} 
\put(-0.1,3.25){\footnotesize $\xelt_{q}$}
\end{picture} \end{center}} 
\hspace*{1.5in}
\parbox{2.7cm}{\begin{center}
\setlength{\unitlength}{0.4cm}
\begin{picture}(6.5,5.5)
\put(3,1){\circle*{0.35}} 
\put(1,3){\circle*{0.35}}
\put(3,5){\circle*{0.35}} 
\put(5,3){\circle*{0.35}}
\put(1,3){\line(1,1){2}} 
\put(3,1){\line(-1,1){2}}
\put(5,3){\line(-1,1){2}} 
\put(3,1){\line(1,1){2}}
\put(1,3){\line(-1,-1){1}}
\multiput(0.4,2.7)(-0.2,-0.2){6}{\color{white}{\rule{3mm}{0.4mm}}}
\put(3,1){\line(-1,-1){1}}
\multiput(2.4,0.7)(-0.2,-0.2){6}{\color{white}{\rule{3mm}{0.4mm}}}
\put(5,3){\line(1,1){1}}
\multiput(5.1,3.2)(0.2,0.2){6}{\color{white}{\rule{3mm}{0.4mm}}}
\put(3,5){\line(1,1){1}}
\multiput(3.1,5.2)(0.2,0.2){6}{\color{white}{\rule{3mm}{0.4mm}}}
\put(1.7,1.75){\scriptsize $j_{p+1}$} 
\put(3.55,1.75){\scriptsize $j'_{q+1}$}
\put(1.55,3.75){\scriptsize $j'_{q+1}$} 
\put(3.7,3.75){\scriptsize $j_{p+1}$}
\put(3.5,0.5){\footnotesize $\yelt'_{q+1}$} 
\put(5.5,2.75){\footnotesize $\yelt'_{q}$}
\put(1.25,5.25){\footnotesize $\yelt_{q-1}$} 
\put(-0.1,3.25){\footnotesize $\yelt_{q}$}
\end{picture} \end{center}} 

\noindent
In the above figures, $\xelt'_{q+1}$ and $\xelt'_{q}$ are in $[\mymin(L_{f}),\telt]$, and $\yelt'_{q+1} = \psi(\xelt'_{q+1})$ and $\yelt'_{q} = \psi(\xelt'_{q})$ are in $\psi([\mymin(L_{f}),\telt]) \subset [\belt,\mymax(L_{g})]$, while $\xelt_{q}$ and $\xelt_{q-1}$ are not in $[\mymin(L_{f}),\telt]$ and $\yelt_{q}$ and $\yelt_{q-1}$ are not in $\psi([\mymin(L_{f}),\telt])$.
When $q=1$, there must be some $\xelt'_{0} \in [\mymin(L_{f}),\telt]$ with $\xelt'_{1} \myarrow{j'_{1}} \xelt'_{0}$, otherwise we have $\xelt'_{1} = \telt$ and $\yelt'_{1} = \psi(\xelt'_{1}) = \psi(\telt) = \mymax(L_{g})$, contradicting the fact that $\yelt'_{1}$ is below $\yelt_{0}$. 
We let $\yelt'_{0} := \psi(\xelt'_{0})$.

\hspace*{1in}
\parbox{2.7cm}{\begin{center}
\setlength{\unitlength}{0.4cm}
\begin{picture}(6.5,5)
\put(3,1){\circle*{0.35}} 
\put(1,3){\circle*{0.35}}
\put(3,5){\circle*{0.35}} 
\put(5,3){\circle*{0.35}}
\put(7,5){\circle*{0.35}}
\put(1,3){\line(1,1){2}} 
\put(3,1){\line(-1,1){2}}
\put(5,3){\line(-1,1){2}} 
\put(3,1){\line(1,1){2}}
\put(1,3){\line(-1,-1){1}}
\multiput(0.4,2.7)(-0.2,-0.2){6}{\color{white}{\rule{3mm}{0.4mm}}}
\put(3,1){\line(-1,-1){1}}
\multiput(2.4,0.7)(-0.2,-0.2){6}{\color{white}{\rule{3mm}{0.4mm}}}
\put(5,3){\line(1,1){2}}
\put(1.7,1.75){\scriptsize $j_{p+1}$} 
\put(3.55,1.75){\scriptsize $j'_{2}$}
\put(1.55,3.75){\scriptsize $j'_{2}$} 
\put(3.7,3.75){\scriptsize $j_{p+1}$}
\put(5.55,3.75){\scriptsize $j'_{1}$}
\put(3.5,0.5){\footnotesize $\xelt'_{2}$} 
\put(5.5,2.5){\footnotesize $\xelt'_{1}$}
\put(1.9,5.25){\footnotesize $\xelt_{0}$} 
\put(-0.1,3.25){\footnotesize $\xelt_{1}$}
\put(7.4,5.2){\footnotesize $\xelt'_{0}$}
\end{picture} \end{center}} 
\hspace*{1.5in}
\parbox{2.7cm}{\begin{center}
\setlength{\unitlength}{0.4cm}
\begin{picture}(6.5,5)
\put(3,1){\circle*{0.35}} 
\put(1,3){\circle*{0.35}}
\put(3,5){\circle*{0.35}} 
\put(5,3){\circle*{0.35}}
\put(7,5){\circle*{0.35}}
\put(1,3){\line(1,1){2}} 
\put(3,1){\line(-1,1){2}}
\put(5,3){\line(-1,1){2}} 
\put(3,1){\line(1,1){2}}
\put(1,3){\line(-1,-1){1}}
\multiput(0.4,2.7)(-0.2,-0.2){6}{\color{white}{\rule{3mm}{0.4mm}}}
\put(3,1){\line(-1,-1){1}}
\multiput(2.4,0.7)(-0.2,-0.2){6}{\color{white}{\rule{3mm}{0.4mm}}}
\put(5,3){\line(1,1){2}}
\put(1.7,1.75){\scriptsize $j_{p+1}$} 
\put(3.55,1.75){\scriptsize $j'_{2}$}
\put(1.55,3.75){\scriptsize $j'_{2}$} 
\put(3.7,3.75){\scriptsize $j_{p+1}$}
\put(5.55,3.75){\scriptsize $j'_{1}$}
\put(3.5,0.5){\footnotesize $\yelt'_{2}$} 
\put(5.5,2.5){\footnotesize $\yelt'_{1}$}
\put(1.9,5.25){\footnotesize $\yelt_{0}$} 
\put(-0.1,3.25){\footnotesize $\yelt_{1}$}
\put(7.4,5.2){\footnotesize $\yelt'_{0}$}
\end{picture} \end{center}} 

\noindent 
Since $\xelt'_{0}$ is now $r$ steps above $\mymin(L_{f})$ and $\xelt'_{0} \in [\mymin(L_{f}),\telt]$, then we must have $\xelt'_{0} = \telt$. 
This means that $\yelt'_{0} = \psi(\xelt'_{0}) = \psi(\telt) = \mymax(L_{g})$. 
However, $\yelt_{0}$ is not comparable to $\yelt'_{0}$, which means $\yelt'_{0}$ cannot be the unique maximal element of $L_{g}$. 
From this contradiction, we see that the hypothesis that the containment $\psi([\mymin(L_{f}),\telt]) \subset [\belt,\mymax(L_{g})]$ is proper is false, so $\psi([\mymin(L_{f}),\telt]) = [\belt,\mymax(L_{g})]$. 

We have therefore shown that $\psi: [\mymin(L_{f}),\telt] \longrightarrow [\belt,\mymax(L_{g})]$ is an isomorphism of DCDL's, which completes the proof.\hfill\QED

{\bf [\S \MinusculeExampleSection.8:\! Distributivity of minuscule compression posets.]} 
In Theorem 8.3.10$\!$ / $\!$Corollary 8.3.11 of the 2013 book \cite{Green}, Green proves distributivity of minuscule compression posets by identifying them as certain induced order subposets of $\Pi(\varpi_{\mbox{\tiny long}})$ and then analyzing cases according to the `$\myA$--$\myG$'-type-based classification of minuscule fundamental weights. 
However, to our knowledge, the first demonstration of distributivity of minuscule compression posets appeared in 1984 in Proposition 4.2 of \cite{PrEur}, where Proctor showed, by invoking the classification of minuscule fundamental weights, that each such $P$ is the distributive lattice of order ideals from some explicitly identified poset.  

Our main goal for this penultimate subsection of \S \MinusculeExampleSection\ is to demonstrate, using general principles, that each minuscule compression poset is a distributive lattice. 
In fact, we intend to demonstrate that each minuscule compression poset is a full-length sublattice of some bi-fold product of chains by showing how \TwoChainLemma\ applies. 
To begin this work, we apply \CompressionPosetLattice\ as part of a general proof demonstrating that each minuscule compression poset is a lattice. 

\noindent 
{\bf \MinusculePosetIsLattice}\ \ {\sl Take as given our minuscule fundamental weight $\omega_{f}$. 
The join of any two meet irreducible elements of the minuscule splitting DCDL $L_{\mathscr{G}}(\omega_{f})$ is meet irreducible, and the meet of any two join irreducible elements of $L_{\mathscr{G}}(\omega_{f})$ is join irreducible. 
In particular, the associated minuscule compression poset $P_{\mathscr{G}}(\omega_{f})$ is a lattice.}

{\em Proof.} 
The claim in the last sentence of the result statement follows from \CompressionPosetLattice. 
For the claim that the join of any two meet irreducible elements is join irreducible, we proceed as follows. 
Let $\mathscr{T}$ be a total ordering of $L = L_{\mathscr{G}}(\omega_{f})$ such that $\xelt <_{\mathscr{T}} \yelt$ whenever $\delta(\xelt) < \delta(\yelt)$, where the depth of any element of $L$ is measured, say, by the length of any shortest path to $\mymax(L)$ or by the cardinality of the element of $L$ when viewed as an upper order ideal from $P = P_{\mathscr{G}}(\omega_{f})$ or {\sl etc}. 
(To obtain such a total ordering $\mathscr{T}$, one could read the elements of $L$ across its ranks, starting at the top.)  
Let $\uelt = [u,\infty)$ and $\velt = [v,\infty)$ be two meet irreducible elements of $L$ identified as up-sets from $P$ and generated by elements $u, v \in P$ respectively. 
When $\uelt$ and $\velt$ are incomparable with $\uelt <_{\mathscr{T}} \velt$, call the pair $(\uelt,\velt)$ a $\mathscr{T}$-{\em pairing} of meet irreducibles. 
Order all such $\mathscr{T}$-pairings lexicographically according to $\mathscr{T}$, which means that $(\uelt',\telt')$ precedes $(\uelt,\velt)$ when either $\uelt' <_{\mathscr{T}} \uelt$ or $\uelt'=\uelt$ with $\velt' <_{\mathscr{T}} \velt$. 
We wish to prove that for any $\mathscr{T}$-pairing $(\uelt,\velt)$, the join $\uelt \vee \velt$ in $L$ is also meet irreducible. 
As a contradiction hypothesis, suppose the desired principle fails, so there is a first $\mathscr{T}$-pairing $(\uelt,\velt)$ of meet irreducibles such that $\uelt \vee \velt$ is not meet irreducible. 

Let $\mu := wt(\uelt) = \sigma.\omega_{f}$ and $\nu := wt(\velt) = \tau.\mu$, where $\sigma = \mygens_{i_{s}} \cdots \mygens_{i_{1}}$ and $\tau = \mygens_{j_{t}} \cdots \mygens_{j_{1}}$ enjoy all the properties of the statement of  \MeetJoinMinusculeProp. 
We can view $\lambda = wt(\uelt \vee \velt)$ as $\mygens_{i_{r}} \cdots \mygens_{i_{1}}.\omega_{f}$ where $i_{r} = j_{r},\ldots,i_{1}=j_{1}$. 
For $p \in [r+1,s]_{\mathbb{Z}}$, let $\mathcal{I}_{p}$ be the multiset $\{i_{r+1},\ldots,i_{p}\}$, and take $\mathcal{I} := \mathcal{I}_{s}$; for $q \in [r+1,t]_{\mathbb{Z}}$, let $\mathcal{J}_{q}$ be the multiset $\{j_{r+1},\ldots,j_{q}\}$, and take $\mathcal{J} := \mathcal{J}_{t}$. 
When we write $\uelt = [u,\infty) \subseteq P$, then $\vcolor(u) = i_{s}$, and with $\velt = [v,\infty) \subseteq P$ we have $\vcolor(v) = j_{t}$. 
Write $\uelt = \uelt_{s} \myarrow{i_{s}} \cdots \mylongarrow{i_{r+2}} \uelt_{r+1} \mylongarrow{i_{r+1}} \uelt_{r} := \uelt \vee \velt$ and $\velt = \velt_{t} \myarrow{j_{t}} \cdots \mylongarrow{j_{r+2}} \velt_{r+1} \mylongarrow{j_{r+1}} \velt_{r} := \uelt \vee \velt$. 
For such an edge $\uelt_{p+1} \mylongarrow{i_{p+1}} \uelt_{p}$, let $x_{p+1} \in P$ be the unique element of color $i_{p+1}$ such that $\uelt_{p} = \uelt_{p+1} \setminus \{x_{p+1}\}$, and for edge $\velt_{q+1} \mylongarrow{j_{q+1}} \velt_{q}$, let $y_{q+1} \in P$ be the unique element of color $j_{q+1}$ such that $\velt_{q} = \velt_{q+1} \setminus \{y_{q+1}\}$. 
Since $\welt := \uelt \vee \velt$ is not meet irreducible, then it must be the case that $\welt$ is below an edge of color $i_{r-1}$ and an edge of color $i_{r}$, where necessarily the distinct Weyl group generators $\mygens_{i_{r-1}}$ and $\mygens_{i_{r}}$ commute and moreover (by \ijChainLemma) nodes $\gamma_{i_{r-1}}$ and $\gamma_{i_{r}}$ are distinct and nonadjacent in the vascular graph $\Gamma$. 
So, there exists a unique color $i_{r-1}$ element $x \in P$ that is minimal in the up-set $\welt$, and there exists a unique color $i_{r}$ element $y \in P$ that is minimal in the up-set $\welt$. 

If $i_{r-1} \in \mathcal{J}$, then $\mygens_{i_{r-1}}$ commutes with $\mygens_{i_{p}}$ and $\gamma_{i_{r-1}}$ and $\gamma_{i_{p}}$ are distinct and nonadjacent in $\Gamma$ for each $p \in [r+1,s]_{\mathbb{Z}}$. 
But then, by \ijChainLemma, we must have a diamond$\,$ 
\parbox{2.75cm}{\begin{center}
\setlength{\unitlength}{0.4cm}
\begin{picture}(5.5,3.5)
\put(2,0){\circle*{0.35}} 
\put(0,2){\circle*{0.35}}
\put(2,4){\circle*{0.35}} 
\put(4,2){\circle*{0.35}}
\put(0,2){\line(1,1){2}} 
\put(2,0){\line(-1,1){2}}
\put(4,2){\line(-1,1){2}} 
\put(2,0){\line(1,1){2}}
\put(0.7,0.65){\em \small $i_{r-1}$} 
\put(2.5,0.65){\em \small $i_{p}$}
\put(0.7,2.65){\em \small $i_{p}$} 
\put(2.7,2.65){\em \small $i_{r-1}$}
\put(2.4,-0.75){\footnotesize $\uelt_{p}$} 
\put(4.4,1.75){\footnotesize $\uelt_{p-1}$}
\end{picture} \end{center}} in $L$. 
But when $p = s$, this implies that $\uelt = \uelt_{s}$ is not meet irreducible. 
So, $i_{r-1} \not\in \mathcal{J}$. 
Similarly see that $i_{r-1} \not\in \mathcal{I}$, $i_{r} \not\in \mathcal{J}$, and $i_{r} \not\in \mathcal{I}$. 
Moreover, our argument shows that $\mygens_{i_{r-1}}$ cannot commute with $\mygens_{i_{p}}$  for all $p \in [r+1,s]_{\mathbb{Z}}$. 
So, let $a \in [r+1,s]_{\mathbb{Z}}$ be largest so that $\mygens_{i_{r-1}}$ commutes with $\mygens_{i_{r+1}},\ldots,\mygens_{i_{a-1}}$ but not $\mygens_{i_{a}}$; in particular, $\gamma_{i_{r-1}}$ and $\gamma_{i_{a}}$ are distinct and adjacent in $\Gamma$. 
Also, let $b \in [r+1,s]_{\mathbb{Z}}$ be largest so that $\mygens_{i_{r}}$ commutes with $\mygens_{i_{r+1}},\ldots,\mygens_{i_{b-1}}$ but not $\mygens_{i_{b}}$, let $a' \in [r+1,t]_{\mathbb{Z}}$ be largest so that $\mygens_{i_{r-1}}$ commutes with $\mygens_{j_{r+1}},\ldots,\mygens_{j_{a'-1}}$ but not $\mygens_{j_{a'}}$, and let $b' \in [r+1,t]_{\mathbb{Z}}$ be largest so that $\mygens_{i_{r}}$ commutes with $\mygens_{j_{r+1}},\ldots,\mygens_{j_{b'-1}}$ but not $\mygens_{j_{b'}}$. 
Then, $\gamma_{i_{r-1}}$ and $\gamma_{i_{b}}$ are distinct and adjacent, as are the node pairs $(\gamma_{i_{r-1}},\gamma_{j_{a'}})$ and $(\gamma_{i_{r}},\gamma_{j_{b'}})$. 

Now, since $x_{a}$ is a minimal element of $\uelt_{a}$, there is a path in $L$ from $\uelt_{a}$ up to $\uelt^{(a)} := [x_{a},\infty)$. 
Notate this as $\uelt_{a} \myarrow{e_{1}} \cdots \myarrow{e_{c}} \uelt^{(a)} := [x_{a},\infty)$ formed by successively removing from the up-set $\uelt_{a}$ the sequence of elements $(x^{(a)}_{1},\ldots,x^{(a)}_{c})$. 
Similarly there are paths $\uelt_{b} \myarrow{f_{1}} \cdots \myarrow{f_{d}} \uelt^{(b)} := [x_{b},\infty)$ formed by successively removing from $\uelt_{b}$ the sequence of elements $(x^{(b)}_{1},\ldots,x^{(b)}_{d})$; $\velt_{a'} \myarrow{e'_{1}} \cdots \myarrow{e'_{c'}} \velt^{(a')} := [y_{a'},\infty)$ formed by successively removing from $\velt_{a'}$ the sequence of elements $(y^{(a')}_{1},\ldots,y^{(a')}_{c'})$; and $\velt_{b'} \myarrow{f'_{1}} \cdots \myarrow{f'_{d'}} \velt^{(b')} := [y_{b'},\infty)$ formed by successively removing from $\velt_{b'}$ the sequence of elements $(y^{(b')}_{1},\ldots,y^{(b')}_{d'})$. 
It is routine to see that each of $\velt$ and $\welt$ is incomparable to each of $\uelt^{(a)}$ and $\uelt^{(b)}$, that $\velt \vee \uelt^{(a)} = \welt \vee \uelt^{(a)}$, and that $\velt \vee \uelt^{(b)} = \welt \vee \uelt^{(b)}$. 
Similarly see that each of $\uelt$ and $\welt$ is incomparable to each of $\velt^{(a')}$ and $\velt^{(b')}$, that $\uelt \vee \velt^{(a')} = \welt \vee \velt^{(a')}$, and that $\uelt \vee \velt^{(b')} = \welt \vee \velt^{(b')}$. 
Moreover, $\uelt_{a} = \welt \wedge \uelt^{(a)}$, $\uelt_{b} = \welt \wedge \uelt^{(b)}$, $\velt_{a'} = \welt \wedge \velt^{(a')}$, and $\velt_{b'} = \welt \wedge \velt^{(b')}$. 
Since the $\mathscr{T}$-pairing $(\uelt,\velt)$ serves as a minimal counterexample, then each of $\velt \vee \uelt^{(a)}$ and $\uelt \vee \velt^{(a')}$ is a meet irreducible, as are each of $\velt \vee \uelt^{(b)}$ and $\uelt \vee \velt^{(b')}$. 
Elementary set-theoretic reasoning shows that the meet irreducible $\velt \vee \uelt^{(a)} = \uelt \vee \velt^{(a')}$ is the same as $\uelt^{(a)} \vee \velt^{(a')} = [x,\infty)$ and is below an edge of color $i_{r-1}$ and that the meet irreducible $\velt \vee \uelt^{(b)} = \uelt \vee \velt^{(b')}$ coincides with $\uelt^{(b)} \vee \velt^{(b')} = [y,\infty)$ and sits below an edge of color $i_{r}$. 

Let $\aelt := \uelt^{(a)} \vee \velt^{(a')}$, $\belt := \uelt^{(b)} \vee \velt^{(b')}$, $\zelt^{-} := \aelt \wedge \belt$, and $\zelt^{+} := \aelt \vee \belt$. 
Let $\mathcal{W}$ (respectively, $\mathcal{E}$, $\mathcal{F}$) be the multiset of edge colors in any path from $\welt$ up to $\zelt^{-}$ (respectively, $\zelt^{-}$ up to $\aelt$, $\zelt^{-}$ up to $\belt$). 
So, $\mathcal{E}$ (respectively, $\mathcal{F}$) is the multiset of edge colors in any path from $\belt$ up to $\zelt^{+}$ (respectively, $\aelt$ up to $\zelt^{+}$). 
Then, $\mathcal{E} \cap \mathcal{F} = \emptyset$, and for every $e \in \mathcal{E}$ and $f \in \mathcal{F}$ the nodes $\gamma_{e}$ and $\gamma_{f}$ are nonadjacent. 
Since $\aelt = \uelt^{(a)} \vee \welt$ and $\uelt_{a} = \uelt^{(a)} \wedge \welt$, then the multiset of edge colors in any path from $\uelt_{a}$ up to $\uelt^{(a)}$ is the multiset union $\mathcal{E} \cup \mathcal{W}$ and the multiset of edge colors in any path from $\uelt^{(a)}$ up to $\aelt$ is $\mathcal{I}_{a}$. 
Similarly see that the multiset of edge colors in any path from $\velt_{a'}$ up to $\velt^{(a')}$ is the multiset union $\mathcal{E} \cup \mathcal{W}$ and that the multiset of edge colors in any path from $\velt^{(a)}$ up to $\aelt$ is $\mathcal{J}_{a'}$. 
Analogously, the multiset of edge colors in any path from $\uelt_{b}$ up to $\uelt^{(b)}$ is the multiset union $\mathcal{F} \cup \mathcal{W}$, and the multiset of edge colors in any path from $\uelt^{(b)}$ up to $\belt$ is $\mathcal{I}_{b}$; also, the multiset of edge colors in any path from $\velt_{b'}$ up to $\velt^{(b')}$ is the multiset union $\mathcal{F} \cup \mathcal{W}$,  and the multiset of edge colors in any path from $\velt^{(b)}$ up to $\belt$ is $\mathcal{J}_{b'}$. 

Our `minimal counterexample' hypothesis for $(\uelt,\velt)$ means that none of $\uelt_{p}$ can be meet irreducible for $p \in [r+1,s-1]_{\mathbb{Z}}$, and for all $q \in [r+1,t-1]_{\mathbb{Z}}$ no $\velt_{q}$ can be meet irreducible. 
Therefore, any meet irreducible elements in $L$ that are strictly between $\uelt_{s}$ and $\uelt^{(a)}$ must have colors from the set $\mathcal{E} \cup \mathcal{W}$, and the same is true for any meet irreducible elements in $L$ strictly between $\velt_{t}$ and $\velt^{(a')}$. 
Similarly, any meet irreducible elements in $L$ that are strictly between $\uelt_{s}$ and $\uelt^{(b)}$ must have colors from the set $\mathcal{F} \cup \mathcal{W}$, and the same is true for any meet irreducible elements in $L$ strictly between $\velt_{t}$ and $\velt^{(b')}$. 
We take the preceding information and interpret it within the context of our compression poset $P$. 
There is a path in $P$ from $u$ (of color $i_{s}$) up to $x_{a}$ (with color $i_{a}$) that only passes through vertices with colors from the set $\mathcal{E} \cup \mathcal{W}$. 
This means there is path of nodes with colors from the set $\mathcal{E} \cup \mathcal{W}$ connecting $\gamma_{i_{s}}$ to $\gamma_{i_{a}}$ in $\Gamma$. 
Similarly, there is a path of nodes with colors from the set $\mathcal{E} \cup \mathcal{W}$ connecting $\gamma_{j_{t}}$ to $\gamma_{j_{a'}}$ in $\Gamma$; there is a path of nodes with colors from the set $\mathcal{F} \cup \mathcal{W}$ connecting $\gamma_{i_{s}}$ to $\gamma_{i_{b}}$; and there is a path of nodes with colors from the set $\mathcal{F} \cup \mathcal{W}$ connecting $\gamma_{j_{t}}$ to $\gamma_{j_{b'}}$. 
Here is a schematic of this information as it appears within $\Gamma$: 

\begin{center}
\setlength{\unitlength}{1cm}
\begin{picture}(8,4)
\put(4,3.5){\circle{0.7}}
\put(3.66,3.45){\small $\gamma_{i_{r-1}}$}
\put(2,3.5){\circle{0.7}}
\put(1.75,3.45){$\gamma_{i_{a}}$}
\put(6,3.5){\circle{0.7}}
\put(5.75,3.45){$\gamma_{j_{a'}}$}
\put(2.35,3.5){\line(1,0){1.3}}
\put(4.35,3.5){\line(1,0){1.3}}
\put(4,0.5){\circle{0.7}}
\put(3.8,0.45){\small $\gamma_{i_{r}}$}
\put(2,0.5){\circle{0.7}}
\put(1.75,0.45){$\gamma_{i_{b}}$}
\put(6,0.5){\circle{0.7}}
\put(5.75,0.45){$\gamma_{j_{b'}}$}
\put(2.35,0.5){\line(1,0){1.3}}
\put(4.35,0.5){\line(1,0){1.3}}
\put(0,2){\circle{0.7}}
\put(-0.25,1.95){$\gamma_{i_{s}}$}
\put(8,2){\circle{0.7}}
\put(7.75,1.95){$\gamma_{j_{t}}$}
\multiput(0.3,2.2)(0.2,0.15){7}{\qbezier(0,0)(0.05,0.0375)(0.1,0.075)}
\put(1.7,3.25){\qbezier(0,0)(0.016,0.012)(0.032,0.024)}
\put(-1.95,3.3){\scriptsize ($1$)}
\put(-1.5,3.3){\scriptsize Nodes with colors}
\put(-1.5,3){\scriptsize from $\mathcal{E} \cup \mathcal{W}$}
\put(0.2,3.05){\vector(2,-1){0.6}}
\multiput(0.3,1.8)(0.2,-0.15){7}{\qbezier(0,0)(0.05,-0.0375)(0.1,-0.075)}
\put(1.7,0.75){\qbezier(0,0)(0.016,-0.012)(0.032,-0.024)}
\put(-1.95,1){\scriptsize ($2$)}
\put(-1.5,1){\scriptsize Nodes with}
\put(-1.5,0.7){\scriptsize colors from $\mathcal{F} \cup \mathcal{W}$}
\put(0.2,1.05){\vector(2,1){0.6}}
\multiput(7.7,2.2)(-0.2,0.15){7}{\qbezier(0,0)(-0.05,0.0375)(-0.1,0.075)}
\put(6.3,3.25){\qbezier(0,0)(-0.016,0.012)(-0.032,0.024)}
\put(7.35,3.3){\scriptsize ($3$)}
\put(7.8,3.3){\scriptsize Nodes with colors}
\put(7.8,3){\scriptsize from $\mathcal{E} \cup \mathcal{W}$}
\put(7.7,3.05){\vector(-2,-1){0.6}}
\multiput(7.7,1.8)(-0.2,-0.15){7}{\qbezier(0,0)(-0.05,-0.0375)(-0.1,-0.075)}
\put(6.3,0.75){\qbezier(0,0)(-0.016,-0.012)(-0.032,-0.024)}
\put(7.35,1){\scriptsize ($4$)}
\put(7.8,1){\scriptsize Nodes with colors}
\put(7.8,0.7){\scriptsize from $\mathcal{F} \cup \mathcal{W}$}
\put(7.75,0.7){\vector(-2,1){0.75}}
\end{picture} 
\end{center}

\noindent 
In the above schematic, the (possibly empty) set of nodes indicated by the dashed line (1) are on a path in $\Gamma$ from $\gamma_{i_{s}}$ and $\gamma_{i_{a}}$ and occur strictly between those two nodes, and similarly for nodes in the sets indicated by (2), (3), and (4). 
The following nodes are necessarily pairwise distinct: $\gamma_{i_{r-1}}$, $\gamma_{i_{r}}$, $\gamma_{i_{s}}$, $\gamma_{j_{t}}$. 
Moreover, $\{i_{a},i_{b},i_{s}\} \cap \{j_{a'},j_{b'},j_{t}\} = \emptyset$, and the nodes $\gamma_{i_{r-1}}$ and $\gamma_{i_{r}}$ are distinct from all other nodes in the schematic. 
There are also many possible coincidences of nodes: Any two of $\gamma_{i_{a}}$, $\gamma_{i_{b}}$, and $\gamma_{i_{s}}$ can be the same;  any two of $\gamma_{j_{a'}}$, $\gamma_{j_{b'}}$, and $\gamma_{j_{t}}$ can be the same; and any two of the sets (1), (2), (3), and/or (4) might have a node in common. 
In all of these circumstances, we obtain a cycle in $\Gamma$. 
This is contrary to the assumption that $\mathscr{G} = (\Gamma,M)$ is a Coxeter--Dynkin flower. 
Therefore our hypothesis that the join of two meet irreducibles in $L$ may fail to be meet irreducible is false. 
That is, the join of any two meet irreducibles in $L$ is also meet irreducible. 
Since the preceding statement also applies to the minuscule splitting DCDL $L^{*}$, we conclude that the meet of any two join irreducibles in $L$ is also join irreducible.\hfill\QED

The previous result showed that the minuscule compression poset of a given minuscule splitting DCDL is itself a lattice. 
The next result (\MinusculePosetIsTopoBalanced) further investigates the lattice structure of a minuscule compression poset. 
Once we show that a given minuscule compression poset is topographically balanced, then we will know it is a modular lattice and can consider its max cover number, cf.\ \S \SubstructureSection.5. 

\noindent 
{\bf \MinusculePosetIsTopoBalanced}\ \ {\sl Take as given our minuscule fundamental weight $\omega_{f}$.  
The minuscule compression poset $P_{\mathscr{G}}(\omega_{f})$ is a topographically balanced lattice, i.e.\ a modular lattice, whose max cover number does not exceed two.}

{\em Proof.} The compression poset $P := P_{\mathscr{G}}(\omega_{f})$ of $L = L_{\mathscr{G}}(\omega_{f})$ is a lattice by \MinusculePosetIsLattice. 
Since $P^{*}$ is also a minuscule compression poset, it suffices to show that for any three distinct elements $u,v,w \in P$ such that $u \rightarrow w \leftarrow v$, there exists a unique $z \in P$ such that $u \leftarrow z \rightarrow v$. 
Set $j := \vcolor(u)$, $k := \vcolor(v)$, and $l := \vcolor(w)$. 
Since $u$ and $v$ are incomparable in $P$, then $j$ and $k$ are distinct and nodes $\gamma_{j}$ and $\gamma_{k}$ of $\Gamma$ are nonadjacent. 
Moreover, since $u$ and $w$ are neighbors in $P$, then $j$ and $l$ are distinct and nodes $\gamma_{j}$ and $\gamma_{l}$ are adjacent; similarly see that $k \not= l$ and that $\gamma_{k}$ and $\gamma_{l}$ are adjacent. 
The $\{j,k,l\}$-submatrix of the pulsation matrix $M$ therefore looks like
\[\begin{array}{rc}
 & \parbox{1.1in}{\hspace*{0.1mm} $\displaystyle \begin{array}{ccc} \hspace*{1.2mm}j\hspace*{1.2mm} & \hspace*{1.2mm}k\hspace*{1.2mm} & \hspace*{1.2mm}l\hspace*{1.2mm}\\ \downarrow & \downarrow & \downarrow\end{array}$}\\
\parbox{0.55in}{\hfill $\displaystyle \begin{array}{r} j\hspace*{1.5mm} \longrightarrow\\ k\hspace*{1.5mm} \longrightarrow\\ l\hspace*{1.5mm} \longrightarrow\end{array}$} & \parbox{1.1in}{$\displaystyle \left(\begin{array}{ccc} 2 & 0 & -a\\ 0 & 2 & -b\\ -d & -c & 2\end{array}\right)$}\mbox{\huge ,}
\end{array}\]
where $a, b, c, d$ are \underline{positive} integers. 

Let $\welt := (-\infty,w]$ be the down-set in $L$ corresponding to $w \in P$. 
As depicted below, in $L$ we have $\telt := \welt \setminus \{w\} \myarrow{l} \welt$, $\uelt = \welt \setminus \{v,w\} \myarrow{k} \telt$, and $\velt = \welt \setminus \{u,w\} \myarrow{j} \telt$; we also have $\selt \myarrow{j} \uelt$ and $\selt \myarrow{k} \velt$, where $\selt := \uelt \cap \velt = \welt \setminus \{u,v,w\}$. 
Let $wt(\welt) =: \lambda = (\lambda_{j},\lambda_{k},1)$. 
Then $wt(\telt) = (\lambda_{j}+d,\lambda_{k}+c,-1) = (1,1,-1)$, $wt(\uelt) = (-1,1,-1+a)$, $wt(\velt) = (1,-1,-1+b)$, and $wt(\selt) = (-1,-1,-1+a+b)$. 
We see that $\lambda_{j}+d = 1$, i.e.\ $\lambda_{j} = 1-d$, hence $\lambda_{j} \leq 0$ and $d \in \{1,2\}$; similarly, $\lambda_{k} \leq 0$ and $c \in \{1,2\}$. 
Also, since $-1+a \in \{-1,0,1\}$ forces $a \in \{1,2\}$; similarly, $b \in \{1,2\}$. 
Further, $-1+a+b \in \{-1,0,1\}$, which forces $a+b=2$ and therefore $a = b = 1$. 
Therefore, $wt(\selt) = (-1,-1,1)$, which means $\selt$ is above an edge of color $j$. 
This implies the existence of a color $j$ element $s \in P$ that is maximal in $\selt$. 
Since $\gamma_{l}$ is adjacent both to $\gamma_{j}$ and $\gamma_{k}$, then neither $\uelt$ nor $\velt$ is below an edge of color $l$. 
Since $\mym_{l}(\uelt) = -1+a = 0 = -1+b = \mym_{l}(\velt)$, then neither $\uelt$ nor $\velt$ is above an edge of color $l$. 
This means that $s$ is not maximal in $\uelt$ and is not maximal in $\velt$. 
But, since $s$ is maximal in $\uelt \setminus \{u\} = \selt = \velt \setminus \{v\}$, then we must have $s \rightarrow u$ and $s \rightarrow v$ in $P$. 
Since $P$ is a lattice, there can be no other $s' \not= s$ such that $s' \rightarrow u$ and $s' \rightarrow v$. 
We have shown there exists a unique $s \in P$ such that $u \leftarrow s \rightarrow v$ and, moreover, that $\vcolor(s) = \vcolor(w)$. 

(The work of the previous paragraph can be extended, as described in the footnote\footnote{We claim here that, in addition to knowing $a=b=1$, we also have $(c,d) \in \{(1,1),(1,2),(2,1)\}$. 
Based on our work in the preceding paragraph, we only need to rule out $c=d=2$. 
Since $\selt$ is above an edge of color $l$ and since $wt(\selt) = (-1,-1,1)$, we conclude that $\telt' := \selt \setminus \{s\}$ has weight $wt(\telt') = (1,1,-1)$. 
This implies $\telt'$ is above an edge of color $j$ (say $\uelt' \myarrow{j} \telt'$) and above an edge of color $k$ (say $\velt' \myarrow{k} \telt'$).  
Since $L$ is a DCDL, there is some $\selt'$ with $\selt' \myarrow{k} \uelt'$ and $\selt' \myarrow{j} \velt'$. 
We have $wt(\uelt') = (-1,1,0)$, $wt(\velt') = (1,-1,0)$, and hence $wt(\selt') = (-1,-1,1)$. 
Therefore we can repeat the process from the top of the paragraph to obtain elements $\telt''$, $\uelt''$, $\velt''$, $\selt''$, {\sl ad infinitum}. 
This contradicts the finiteness of $L$. 
Now using the facts that $a=b=1$ and $(c,d) \in \{(1,1),(1,2),(2,1)\}$, one can easily describe the $\{j,k,l\}$-components of $L$ as a diamond, with parallel edges of colors $j$ and $k$, that is above and below chains of length one each with edge color $l$ when $(c,d)=(1,1)$; above and below chains of length two, each with edges of colors $k$ and $l$ when $(c,d)=(2,1)$; above and below chains of length two, each with edges of colors $j$ and $l$ when $(c,d)=(1,2)$.} below, to a classification of the $\{j,k,l\}$-components of $L$.)
 
Now suppose some $v \in P$ covers exactly $m$ elements, say $u_{1},\ldots,u_{m}$, with $m$ as large as possible. 
Let $j := \vcolor(v)$, and set $i_{p} := \vcolor(u_{p})$ for $p \in [1,m]_{\mathbb{Z}}$. 
We can see that $\gamma_{j}$ is distinct from and adjacent to each $\gamma_{i_{p}}$ and that $\gamma_{i_{p}}$ and $\gamma_{i_{q}}$ are distinct and nonadjacent when $p \not= q$. 
Let $\uelt := \langle u_{1},\ldots,u_{m} \rangle^{\downarrow}$, the down-set from $P$ generated by $u_{1},\ldots,u_{m}$, and let $\velt := (-\infty,v]$. 
Of course, $\velt$ is join irreducible in $L$, $\uelt \myarrow{j} \velt$, and $\mym_{j}(\uelt) = -1$. 
Since each $u_{p}$ is maximal in $\uelt$, then $\uelt_{p} := \uelt \setminus \{u_{p}\}$ is a down-set from $P$. 
Indeed, $\uelt_{1},\ldots,\uelt_{m}$ comprise exactly the elements of $L$ that are covered by $\uelt$. 
Let $\telt := \uelt \setminus \{u_{1},\ldots,u_{m}\}$. 
By \IntervalProp, $[\telt,\uelt]$ is isomorphic to a diamond-colored Boolean lattice of size $2^{m}$ with edges colored by $\{i_{1},\ldots,i_{m}\}$. 
One can also see that the subset $\{u_{1},\ldots,u_{m}\}$ of $P$, which is an antichain consisting of $m$ distinct elements, is an $\{i_{1},\ldots,i_{m}\}$-subordinate of $P$. 
Therefore the $\{i_{1},\ldots,i_{m}\}$-component of $\uelt$ in $L$ is precisely $[\telt,\uelt]$.  
Since $\telt$ is only below the down-sets $\telt \cup \{u_{p}\}$, for $p \in [1,m]_{\mathbb{Z}}$, then $\telt$ is not below an edge of color $j$, and hence $\mym_{j}(\telt) \geq 0$. 
Within our minuscule context, of course, we have $\mym_{j}(\telt) \leq 1$. 
So, $-1 = \mym_{j}(\uelt) = \mym_{j}(\telt) + \sum_{p=1}^{m}M_{i_{p},j}$, i.e.\ $\sum_{p=1}^{m}M_{i_{p},j} = -1-\mym_{j}(\telt)$. 
Since each $M_{i_{p},j}$ is a negative integer, then $-1-\mym_{j}(\telt) \leq -m$, i.e.\ $m \leq 1+\mym_{j}(\telt) \leq 2$.\hfill\QED

In view of the preceding result, we can invoke \TwoChainLemma\ to obtain the following theorem. 
Note that the proof of this result is completely general in that it does not depend upon the type classification of embryophytic graphs nor upon the classification of minuscule fundamental weights.  

\noindent 
{\bf \MinusculePosetIsDistributive}\ \ {\sl Take as given our minuscule fundamental weight $\omega_{f}$. 
The minuscule compression poset $P_{\mathscr{G}}(\omega_{f})$ is a distributive lattice and is isomorphic to some full-length sublattice of a  product of two (nonempty) chains.} 

{\em Proof.} In view of \MinusculePosetIsTopoBalanced, \TwoChainLemma\ applies.\hfill\QED

{\bf [\S \MinusculeExampleSection.9:\! Vascular graph doppelg\"{a}ngers in minuscule compression posets.]} 
In this subsection we observe that, at the top and at the bottom of a given minuscule compression poset, there are mirror copies of the vascular graph of our Coxeter--Dynkin flower.  
This idea was discerned and utilized in a more general setting by Proctor in \cite{PrJACO}, where he referred to the topmost copy of the vascular graph as a ``top tree''. 
Here, we call it a ``vascular graph doppelg\"{a}nger''. 
In particular, with reference to the next theorem, the {\em lower (respectively, upper) vascular graph doppelg\"{a}nger} the unique lower (respectively, upper) order ideal from a given minuscule poset that is isomorphic, as an undirected graph, to the vascular graph for our given Coxeter--Dynkin flower. 

\noindent 
{\bf \MinusculeTopTree}\ \ {\sl Take as given our minuscule fundamental weight $\omega_{f}$. 
There is a unique lower order ideal from the minuscule compression poset $P_{\mathscr{G}}(\omega_{f})$ that is isomorphic as a vertex-colored undirected graph to our vascular graph $\Gamma$. 
Likewise, there is a unique upper order ideal from $P_{\mathscr{G}}(\omega_{f})$ that is isomorphic as a vertex-colored undirected graph to our vascular graph $\Gamma$.}

{\em Proof.} 
Let $\mymin = \mymin\big(P_{\mathscr{G}}(\omega_{f})\big)$ be the unique minimal element of $P := P_{\mathscr{G}}(\omega_{f})$, which has color $g := \vcolor(\mymin) = \sigma_{0}(f)$. 
Throughout the proof, we think of elements of $L := L_{\mathscr{G}}(\omega_{f})$ as weights from $\Pi(\omega_{f})$. 
So, $-\omega_{g} \myarrow{g} -\omega_{g}+\alpha_{g}$ in $L := L_{\mathscr{G}}(\omega_{f})$. 
We know, by \S \MinusculeExampleSection.5 and \StarBowtieLemma, that $\omega_{g}$ is also minuscule. 
Say $\gamma_{h}$ is a neighbor of $\gamma_{g}$ in $\Gamma$, so $\langle \alpha_{g},\alpha_{h}^{\vee} \rangle < 0$. 
Then, $\langle -\omega_{g}+\alpha_{g},\alpha_{h}^{\vee} \rangle = \langle -\omega_{g},\alpha_{h}^{\vee} \rangle + \langle \alpha_{g},\alpha_{h}^{\vee} \rangle = \langle \alpha_{g},\alpha_{h}^{\vee} \rangle < 0$, so $-\omega_{g}+\alpha_{g}$ is below an edge of color $h$. 
Since, by \MinusculePosetIsTopoBalanced, the max cover number of $L$ is at most two, it follows that $\gamma_{g}$ has at most two neighbors in $\Gamma$. 

Now suppose $(\gamma_{g},\gamma_{h_{1}},\ldots,\gamma_{h_{p}})$ is a simple path in $\Gamma$ from $\gamma_{g}$ to any leaf $\gamma_{h_{p}}$; call this a {\em leaf path from} $\gamma_{g}$. 
So, in $L$ we have $-\omega_{g} \mylongarrow{g} (-\omega_{g}+\alpha_{g}) \mylongarrow{h_{1}} (-\omega_{g}+\alpha_{g}+\alpha_{h_{1}}) \mylongarrow{h_{2}} \cdots \mylongarrow{h_{p}} \big(-\omega_{g}+\sum_{i=1}^{p}\alpha_{h_{i}}\big)$; call this the $L$-{\em companion leaf path from} $\omega_{g}$ and call $-\omega_{g}+\sum_{i=1}^{p}\alpha_{h_{i}}$ its {\em leaf vertex}. 
We now show that for each $q \in [1,p]_{\mathbb{Z}}$, the weight $-\omega_{g}+\sum_{i=1}^{q}\alpha_{h_{i}}$ is a join irreducible in $L$. 
For convenience, set $h_{0} := g$.  
Let $q \in [1,p]_{\mathbb{Z}}$, and consider $\nu := -\omega_{g}+\sum_{i=1}^{q}\alpha_{h_{i}}$. 
Of course, we have $\mu \mylongarrow{h_{q}} \nu$ when $\mu := -\omega_{g}+\sum_{i=1}^{q-1}\alpha_{h_{i}}$ because $\langle \mu,\alpha_{h_{q}}^{\vee} \rangle = \langle -\omega_{g},\alpha_{h_{q}}^{\vee} \rangle + \sum_{i=1}^{q-1}\langle \alpha_{h_{i}},\alpha_{h_{q}}^{\vee} \rangle = \langle \alpha_{h_{q-1}},\alpha_{h_{q}}^{\vee} \rangle < 0$. 
Suppose $\xi \myarrow{j} \nu$ for some $\xi \in L$. 
Then, $\langle \nu,\alpha_{j}^{\vee} \rangle = 1$. 
Now, the facts that $\langle -\omega_{g}+\alpha_{g},\alpha_{g}^{\vee} \rangle = 1$, $\langle \alpha_{h_{1}},\alpha_{g}^{\vee} \rangle < 0$, and $\sum_{i=2}^{q}\langle \alpha_{h_{i}},\alpha_{g}^{\vee} \rangle = 0$ imply that $\langle \nu,\alpha_{g}^{\vee} \rangle \leq 0$, so $j \not= g$. 
Now suppose $j \in \{h_{1},\ldots,h_{q-1}\}$, say $j = h_{r}$. 
Then, $\gamma_{j} = \gamma_{h_{r}}$ is adjacent to $\gamma_{h_{r-1}}$ and $\gamma_{h_{r+1}}$, meaning that $\langle \alpha_{h_{r-1}},\alpha_{j}^{\vee} \rangle$ and $\langle \alpha_{h_{r+1}},\alpha_{j}^{\vee} \rangle$ are both negative, while $\langle \alpha_{h_{i}},\alpha_{j}^{\vee} \rangle = 0$ for $i \in [0,q]_{\mathbb{Z}} \setminus \{r-1,r,r+1\}$. 
But then $\langle \nu,\alpha_{j}^{\vee} \rangle = \langle -\omega_{g},\alpha_{j}^{\vee} \rangle + \langle \alpha_{j},\alpha_{j}^{\vee} \rangle + \sum_{h_{i}\not=j}\langle \alpha_{h_{i}},\alpha_{j}^{\vee} \rangle = 2 + \langle \alpha_{h_{r-1}},\alpha_{j}^{\vee} \rangle + \langle \alpha_{h_{r+1}},\alpha_{j}^{\vee} \rangle \leq 0 \not= 1$.  
Next, suppose $j \in \{h_{q+1},\ldots,h_{p}\}$, say $j = h_{r}$. 
Then, $\langle \nu,\alpha_{j}^{\vee} \rangle = \langle -\omega_{g},\alpha_{j}^{\vee} \rangle + \sum_{i=0}^{q}\langle \alpha_{h_{i}},\alpha_{j}^{\vee} \rangle = \langle \alpha_{h_{q}},\alpha_{j}^{\vee} \rangle$, and the latter cannot equal $1$. 
If $\gamma_{j}$ is distinct from and not adjacent to any of $\{\gamma_{h_{i}}\}_{i=1}^{q-1}$, then computations similar to the above show that $\langle \nu,\alpha_{j}^{\vee} \rangle$ cannot be $1$. 
Now suppose $j \not \in \{h_{0},h_{1},\ldots,h_{q}\}$ but that $\gamma_{j}$ is adjacent to some $\gamma_{h_{r}}$ with $r \in [0,q]_{\mathbb{Z}}$. 
Since $\Gamma$ is acyclic, then $\gamma_{j}$ cannot be adjacent to $\gamma_{h_{i}}$ if $i \in [0,q]_{\mathbb{Z}} \setminus \{r\}$. 
Then $\langle \nu,\alpha_{j}^{\vee} \rangle = \langle -\omega_{g},\alpha_{j}^{\vee} \rangle + \sum_{i=0}^{q}\langle \alpha_{h_{i}},\alpha_{j}^{\vee} \rangle = \langle \alpha_{h_{r}},\alpha_{j}^{\vee} \rangle$, and the latter cannot equal $1$. 
The foregoing analysis leaves $j=h_{q}$ as the only possibility. 
 
Now take the union $\mathcal{I}$ of all vertices from all $L$-companion leaf paths from $-\omega_{g}$, but discard $-\omega_{g}$. 
So, all vertices in $\mathcal{I}$ are join irreducibles in $L$ and therefore elements of $P$, and these are in an evident one-to-one color-preserving correspondence with the nodes of $\Gamma$. 
Since all elements of $\mathcal{I}$ are comparable to some leaf vertex from some $L$-companion leaf path, then $\mathcal{I}$ is a lower order ideal from $P$ that is therefore isomorphic to $\Gamma$ as an undirected vertex-colored graph. 
Suppose $u \rightarrow v$ in $\mathcal{I}$. 
If $u \rightarrow w$ for some $w \not= v$, then $w$ must have a color distinct from the color of $v$, by \ijChainLemma.2. 
Uniqueness of $\mathcal{I}$ follows from the preceding fact, as we are prevented in branching away from $u$ to find another copy of our vascular graph.  
An entirely similar argument will produce a unique upper order ideal $\mathcal{F}$ from $P$ that is isomorphic to $\Gamma$ as an undirected vertex-colored graph.\hfill\QED

\begin{center}
\underline{\hspace*{4in}}
\end{center}

\vspace*{0.5cm} 
\noindent
{\bf \S \MinusculeLatticePosetSection. Skew-stacks of minuscule compression posets.} 
Minuscule splitting DCDL's are among the most compelling and naturally occurring diamond-colored distributive lattices.  
The combinatorial distinctiveness of minuscule splitting DCDL's and their companion minuscule compression posets seems first to have been recognized by Stanley (e.g.\ \cite{StanLefschetz}) and Proctor (e.g\ \cite{PrMonthly}, \cite{PrEur}). 
Moreover, as recorded in \cite{PrEur}, Proctor and Stanley also recognized the algebraic-combinatorial significance of the distributive lattice of order ideals taken from the product of a chain and a minuscule compression poset. 

As in \S \FrameSection, we now view the product $P_{\mathscr{G}}(\omega_{f}) \times [1,m]_{\mathbb{Z}}$ of a minuscule compression poset $P_{\mathscr{G}}(\omega_{f})$ by a chain with $m$ elements as the {\em minuscule prism}  $P_{\mathscr{G}}(\omega_{f})^{(\mysmallposetstack m)}$, which is isomorphic to the skew-stack $P_{1} \myposetstack_{\phi} P_{2} \myposetstack_{\phi} \cdots \myposetstack_{\phi} P_{m}$ when we take $P_{i} := P_{\mathscr{G}}(\omega_{f})$ for each $i \in [1,m]_{\mathbb{Z}}$ and take $\phi: P_{\mathscr{G}}(\omega_{f}) \longrightarrow P_{\mathscr{G}}(\omega_{f})$ to be the identity function. 
When $\mystackedP$ is a minuscule prism, we call $\Jcolor(\mystackedP)$ its associated {\em prismatic-minuscule DCDL}. 
More generally, we call any skew-stack $\mystackedP$ of minuscule compression posets a {\em polyminuscule skew-stack}, and we call $\Jcolor(\mystackedP)$ its associated {\em polymimuscule DCDL}. 

The main purpose of this closing section is to present, in \GStructureDCDLCorollary, the statement and proof of a sturdiness result for polyminuscule skew-stacks, including minuscule prisms. 
This result is used in \cite{DDW} to demonstrate sturdiness of some $\myE_{6}$-polyminuscule DCDL's and of some $\myE_{7}$-prismatic minuscule DCDL's. 
We also hope to apply it in \cite{DDgeneral}, a paper whose aim will be to extend and generalize the results of \cite{DD}. 
Before the main result of this section, we need the following.

\noindent 
{\bf \DoppelLemma}\ \ {\sl Suppose that $P_{1}$ and $P_{2}$ are $\mathscr{G}$-minuscule compression posets and that $\phi: \mathcal{I}_{1} \longrightarrow \mathcal{F}_{2}$ is, as in \IdealMatchingProp, a vertex-color- and directed-edge-preserving isomorphism from a nonempty lower order ideal $\mathcal{I}_{1}$ from $P_{1}$ to a nonempty upper order ideal $\mathcal{F}_{2}$ from $P_{2}$. 
Further suppose that nodes $\gamma_{i}$ and $\gamma_{j}$ are distinct and adjacent nodes in the vascular graph $\Gamma(\mathscr{G})$ of our Coxeter--Dynkin flower $\mathscr{G}$. 
Let $u \in \mathcal{I}_{1}$ and $\phi(u) \in \mathcal{F}_{2}$ be color $i$ vertices. 
Then, $u$ is a vertex of the lower vascular graph doppelg\"{a}nger of $P_{1}$ whose neighboring color $j$ doppelg\"{a}nger vertex is not in $\mathcal{I}_{1}$ if and only if $\phi(u)$ is a vertex of the upper vascular graph doppelg\"{a}nger of $P_{2}$ whose neighboring color $j$ doppelg\"{a}nger vertex is not in $\mathcal{F}_{2}$.}

{\em Proof.} Suppose $u$ is a vertex of the lower vascular graph doppelg\"{a}nger of $P_{1}$ whose neighboring color $j$ doppelg\"{a}nger vertex is not in $\mathcal{I}_{1}$. 
Now, $\mathcal{I}_{1}$ can have no color $j$ vertices, for the following reason: By \ijChainLemma.2, the color $j$ vertices within $P_{1}$ are totally ordered, and, amongst these, the minimal color $j$ vertex is a lower doppelg\"{a}nger vertex which, by hypothesis, is not in $\mathcal{I}_{1}$. 
Therefore $\mathcal{F}_{2} = \phi(\mathcal{I}_{1})$ has no color $j$ vertices. 
By using the fact that all color $i$ vertices in $P_{2}$ are also totally ordered, then the color $i$ upper doppelg\"{a}nger vertex of $P_{2}$ must be a member of $\mathcal{F}_{2}$ since $\phi(u)$ is a color $i$ member of the upper order ideal $\mathcal{F}_{2}$. 
It only remains to be seen that $\phi(u)$ is, in fact, an upper doppelg\"{a}nger vertex of $P_{2}$. 
Suppose otherwise, so $\phi(u) < \phi(u')$, where $\phi(u')$ is the upper doppelg\"{a}nger vertex of $P_{2}$. 
Let $v$ be the upper doppelg\"{a}nger vertex of color $j$ that is a neighbor of $\phi(u')$. 
Of course, $v \not\in \mathcal{F}_{2}$, so therefore $v \rightarrow \phi(u')$ in $P_{2}$. 
Since the set of all color $i$ or $j$ vertices in $P_{2}$ is totally ordered (\ijChainLemma.2, again), then the facts $\phi(u) < \phi(u')$ and that $v \rightarrow \phi(u')$ forces us to have $\phi(u) < v$. 
But this contradicts the fact that $v \not\in \mathcal{F}_{2}$. 
So, $\phi(u)$ is an upper doppelg\"{a}nger vertex of $P_{2}$, and the upper doppelg\"{a}nger vertex that is the color $j$ neighbor of $\phi(u)$ is not in $\mathcal{F}_{2}$. 
The argument for the converse is entirely similar.\hfill\QED

\noindent 
{\bf \GStructureDCDLCorollary}\ \ {\sl Let $\mathscr{G} = (\Gamma_{I},M_{I \times I})$ be a Coxeter--Dynkin flower with palette of colors $I$.  
Consider the connected polyminuscule skew-stack} $\mystackedP := P_{1} \myposetstack_{\phi_{1}} P_{2} \myposetstack_{\phi_{2}} \cdots \myposetstack_{\phi_{m-1}} P_{m}$ {\sl and its corresponding birds-eye scaffold} $\scaffoldS^{\mybirdseye} = (\mystackedP^{\mybirdseye},\mathcal{A}^{\mybirdseye},\mathcal{B}^{\mybirdseye},\vcolor^{\mybirdseye}),$ {\sl where each tier $P_{t}$ of} $\mystackedP$ {\sl is a $\mathscr{G}$-minuscule compression poset $P_{\mathscr{G}}(\omega_{f_{t}})$ and each bond $\phi_{t}$ matches a nonempty lower order ideal from $P_{t}$ with a nonempty upper order ideal from $P_{t+1}$ as in \IdealMatchingProp. 
Then the vertex-colored birds-eye poset} $\mystackedP^{\mybirdseye}$ {\sl enjoys properties (1) and (2) from the statement of \ijChainLemma, is a distributive lattice that is a full-length sublattice of a product of two (nonempty) chains, and has maximal element} $\overline{\mymax(P_{1})}$ {\sl of color $f_{1}$ and minimal element} $\overline{\mymin(P_{m})}$ {\sl of color $\sigma_{0}(f_{m})$. 
Moreover, the polyminuscule DCDL identified by the isomorphic equivalence} $\Jcolor(\mystackedP) \cong \Lcolor(\scaffoldS^{\mybirdseye})$ {\sl is $\mathscr{G}$-structured and its monochromatic components are chain products.}
 
{\em Proof.} 
We begin by demonstrating the claims about $\mystackedP^{\mybirdseye}$. 
It follows from \MinusculePosetIsDistributive\ together with \StackLemma.5 that $\mystackedP^{\mybirdseye}$ is a distributive lattice that is a full-length sublattice of a product of two (nonempty) chains. 
That $\mystackedP^{\mybirdseye}$ has maximal element $\overline{\mymax(P_{1})}$ of color $f_{1}$ and minimal element $\overline{\mymin(P_{m})}$ of color $\sigma_{0}(f_{m})$ follows from \StarBowtieLemma\ together with \StackLemma.4. 
Next we confirm that $\mystackedP^{\mybirdseye}$ enjoys property (1) from the statement of \ijChainLemma. 
So, suppose $\overline{u} \rightarrow \overline{v}$ is a covering relation in $\mystackedP^{\mybirdseye}$, and set $i := \vcolor^{\mybirdseye}(\overline{u})$ and $j := \vcolor^{\mybirdseye}(\overline{v})$. 
By \StackLemma.3, we must have $u' \rightarrow v'$ as a plank within some tier $P_{t}$, where $\mytier(u') = t = \mytier(v')$ for some $u' \in \overline{u}$ and $v' \in \overline{v}$. 
With $\vcolor^{(t)}: P_{t} \longrightarrow I$ as the vertex-coloring function for this tier and with $i = \vcolor^{(t)}(u')$ and $j = \vcolor^{(t)}(v')$, we have, by applying \ijChainLemma\ to the minuscule compression poset $P_{t}$, that $\gamma_{i}$ and $\gamma_{j}$ are distinct and adjacent in the vascular graph $\Gamma$. 

To confirm property (2) of \ijChainLemma\ for $\mystackedP^{\mybirdseye}$, we will say why $(\vcolor^{\mybirdseye})^{-1}(\{i,j\})$ is a totally ordered subset of $\mystackedP^{\mybirdseye}$ for distinct colors $i,j \in I$ such that $\gamma_{i}$ and $\gamma_{j}$ are adjacent in $\Gamma$. 
Say $\overline{u}$ and $\overline{v}$ are distinct elements of $(\vcolor^{\mybirdseye})^{-1}(\{i,j\})$. 
We will demonstrate that $\overline{u}$ and $\overline{v}$ are comparable in $\mystackedP^{\mybirdseye}$. 
Suppose there exist $u' \in \overline{u}$ and $v' \in \overline{v}$ such that $\mytier(u') = \mytier(v')$. 
By \ijChainLemma.2, the color $i$ or $j$ vertices of this tier are totally ordered. 
So, $u'$ and $v'$ are comparable within this tier, and therefore by \StackLemma.2, $\overline{u}$ and $\overline{v}$ are comparable in $\mystackedP^{\mybirdseye}$. 
Now suppose there is no tier which contains elements from each of $\overline{u}$ and $\overline{v}$. 
We may assume that $r := \mytier(v^{\uparrow}) < \mytier(u^{\downarrow}) =: t$. 

We claim that if $w_{t-1}$ is the least element of $P_{t-1} \setminus \mathcal{I}_{t-1}$ of color $i$ or $j$, then $\overline{u} \leq_{\mysmallbirdseye} \overline{w_{t-1}}$. 
Note that since $t = \mytier(u^{\downarrow})$, then $u^{\downarrow} \not\in \mathcal{F}_{t}$. 
Let us consider the case that an upper doppelg\"{a}nger vertex of $P_{t}$ of color $i$ or $j$ -- momentarily call it $u'$ -- is within $\mathcal{F}_{t}$. 
Then $u^{\downarrow} < u'$ in $P_{t}$. 
Since the color $i$ or $j$ vertices of $P_{t-1}$ are totally ordered, then $\phi_{t-1}^{-1}(u') < w_{t-1}$ in $P_{t-1}$. 
Now by \StackLemma.2, $\overline{u} \leq_{\mysmallbirdseye} \overline{w_{t-1}}$. 
Now consider the case that no upper doppelg\"{a}nger vertex of $P_{t}$ of color $i$ or $j$ is within $\mathcal{F}_{t}$. 
It follows that $\mathcal{F}_{t}$ has no vertices of color $i$ or $j$. 
Now let $u'$ be the maximal element of color $i$ or $j$ in $P_{t}$, so $u^{\downarrow} \leq u'$ in $P_{t}$. 
Within the upper doppelg\"{a}nger of $P_{t}$, walk from $u'$ up to the first upper doppelg\"{a}nger vertex $u''$ that is in $\mathcal{F}_{t}$. 
By \DoppelLemma, $\phi_{t-1}^{-1}(u'')$ is the vertex of the lower doppelg\"{a}nger of $P_{t-1}$ with the same color as $u''$. 
Therefore, within the lower doppelg\"{a}nger of $P_{t-1}$, we can walk from $\phi_{t-1}^{-1}(u'')$ up to the first lower doppelg\"{a}nger vertex of color $i$ or $j$, which is $w_{t-1}$. 
Again by \StackLemma.2, $\overline{u} \leq_{\mysmallbirdseye} \overline{w_{t-1}}$. 

If $r = t-1$, then we have $\overline{u} \leq_{\mysmallbirdseye} \overline{w_{t-1}} \leq_{\mysmallbirdseye} \overline{v}$, and we are done with our confirmation of property (2) of \ijChainLemma\ for $\mystackedP^{\mybirdseye}$. 
If $r < t-1$, then repeatedly re-apply the reasoning of the preceding paragraph to see that $\overline{u} \leq_{\mysmallbirdseye} \overline{w_{t-1}} \leq_{\mysmallbirdseye} \overline{w_{t-2}} \leq_{\mysmallbirdseye} \cdots \leq_{\mysmallbirdseye} \overline{w_{r}} \leq_{\mysmallbirdseye} \overline{v}$. 

Having verified that the birds-eye poset $\mystackedP^{\mybirdseye}$ enjoys properties (1) and (2) from \ijChainLemma, then the claim that the monochromatic components of $\Lcolor(\scaffoldS^{\mybirdseye})$ are chain products now follows from \ChainProp. 
To prove the remaining part of the concluding statement of the theorem, we will apply \GStructureDCDLTheorem\ by verifying that its hypotheses hold within our present setting. 
Let $L_{t} := \Jcolor(P_{t}) = L_{\mathscr{G}}(\omega_{f_{t}})$. 
That the $\mathscr{G}$-minuscule splitting DCDL $L_{t}$ is $\mathscr{G}$-structured follows from \MinusculeDCDL. 
By \NewLFKTheorem.2, $\ecolor\! :\! \EdgeSet(L_{t}) \rightarrow I$ is surjective. 
That $\vcolor^{(t)}(x) \not= \vcolor^{(t)}(y)$ whenever $x \rightarrow y$ in $P_{t}$ follows from \ijChainLemma. 

Take two colors $\{i,j\} =: J \subseteq I$ with $i \not= j$, and let $\mathscr{G}' = (\Gamma',M')$ be the $J$-sub-embryophyte of $\mathscr{G}$. 
Also, take two indices $r,s \in \{1,2,\ldots,m\}$ and two connected and nonempty $J$-subordinates $\mathcal{C}$ and $\mathcal{C}'$ from tiers $P_{r}$ and $P_{s}$ respectively. 
To complete the proof, we argue that $\mathcal{C}$ and $\mathcal{C}'$ are saturated chains of the same length. 
This can be done by classifying the possible $J$-components of any $\mathscr{G}$-minuscule splitting DCDL. 
But here we will use more general reasoning. 
Since any nonempty $J$-subordinate is the result of removing a down-set $\mathcal{I}'$ from a down-set $\mathcal{I}$ that properly contains $\mathcal{I}'$, it follows that our given $J$-subordinates $\mathcal{C}$ and $\mathcal{C}'$ are saturated chains if they are chains at all. 
Our focus for the moment will be on the $J$-subordinate $\mathcal{C}$. 
Let $K = \Jcolor(\mathcal{C})$ be a $J$-component of $L_{r}$ whose corresponding $J$-subordinate in $P_{r}$ is $\mathcal{C}$. 
Observe that $K$ is a chain if and only if $\mathcal{C}$ is a chain. 
Suppose for the moment that $K$ is not a chain.  
Then some $\xelt \in K$ is above two elements of $K$ or is below two elements of $K$. 
In the former case, $\xelt$ is above an edge of color $i$ and an edge of color $j$, and, since any monochromatic component of $L_{r}$ is either a singleton or a two-element chain, then $\xelt$ must be a maximal element in $K$. 
Rename this maximal element to be `$\uelt$'. 
Diamond-coloring of $L_{r}$ implies we have the following substructure within $K$:$\,$ \parbox{1.4cm}{\begin{center}
\setlength{\unitlength}{0.2cm}
\begin{picture}(6.5,3.5)
\put(3,0){\circle*{0.5}} \put(1,2){\circle*{0.5}}
\put(3,4){\circle*{0.5}} \put(5,2){\circle*{0.5}}
\put(1,2){\line(1,1){2}} \put(3,0){\line(-1,1){2}}
\put(5,2){\line(-1,1){2}} \put(3,0){\line(1,1){2}}
\put(1.75,0.55){\em \small j} \put(3.7,0.7){\em \small i}
\put(1.7,2.7){\em \small i} \put(3.75,2.55){\em \small j}
\put(3.5,-0.75){\footnotesize $\relt$} \put(5.75,1.75){\footnotesize $\telt$}
\put(3.5,4){\footnotesize $\uelt$} \put(-0.5,1.75){\footnotesize $\selt$}
\end{picture} \end{center}}. 
There can be no other incident edges of color $i$ or $j$ without violating the constraint that all monochromatic components are either singleton elements or two-element chains. 
That is, $K$, in its entirety, is the diamond just depicted. 
If $\xelt$ is below two other elements of $K$, then by similar reasoning we see that $K$ is the preceding diamond with $\xelt = \relt$. 
In both cases, it follows that $M' = \left(\begin{array}{cc}2 & 0\\ 0 & 2\end{array}\right)$. 
But then our $J$-subordinate of $P_{r}$ is the disjoint sum of an element of color $i$ and an incomparable element of color $j$ and is therefore not connected. 
Therefore $K$ must be a chain in $L_{r}$, and therefore $\mathcal{C}$ is a saturated chain in $P_{r}$.  
We can also conclude from our work in this paragraph that if some $J$-component of $L_{r}$ is not a chain, then every connected and nonempty $J$-subordinate of every $P_{t}$ (where $t \in [1,m]_{\mathbb{Z}}$) is a singleton vertex of color $i$ or of color $j$.  

For any $\xelt \in K$, let the quantity $\wt_{\mathcal{C}}(\xelt) := \mym_{i}(\xelt)\omega_{i}^{J}+\mym_{j}(\xelt)\omega_{j}^{J}$ denote the weight of $\xelt$ as an element of the $\mathscr{G}'$-structured DCDL $K$. 
Let $\alpha_{i}' := 2\omega_{i}^{J}+M_{ij}\omega_{j}^{J}$ and $\alpha_{j}' := M_{ji}\omega_{i}^{J}+2\omega_{j}^{J}$. 
For the maximal element $\uelt$ in $K$, we have $\lambda := \wt_{\mathcal{C}}(\uelt) \in \{\omega_{i}^{J}\, ,\, \omega_{j}^{J}\}$. 
Also, for each element $\xelt$ of $K$, we have $\wt_{\mathcal{C}}(\xelt) = \lambda - \#_{i}(\xelt)\alpha_{i}' - \#_{j}(\xelt)\alpha_{j}'$, where $\#_{k}(\xelt)$ is the number of color $k$ edges in any path in $K$ from $\xelt$ up to $\uelt$, for $k \in \{i,j\}$. 
That is, the elements of $K$ are in one-to-one correspondence with the set $\wt_{\mathcal{C}}(K)$. 
Now, generators for the parabolic subgroup $\mathcal{W}_{J}$ of the Weyl group $\mathcal{W}(\mathscr{G})$ act on $\mbox{\sffamily span}_{\mathbb{R}}\{\omega_{i}^{J},\omega_{j}^{J}\}$ by the rules $\mygens_{i}.(a\omega_{i}^{J}+b\omega_{j}^{J}) = 2a\omega_{i}^{J}-b\omega_{j}^{J}$ and $\mygens_{j}.(a\omega_{i}^{J}+b\omega_{j}^{J}) = -a\omega_{i}^{J}+2b\omega_{j}^{J}$.  
That is, we can view these generator actions as moves of the $(1,1)$-seeded two-node networked-numbers game as studied in \S 4 of \cite{DDHK}, with $p := M_{ij}$ and $q = M_{ji}$. 
By Theorem 4.2 of that paper, it follows that the length of the game is the same when played from starting position $(a,b) = (1,0)$ as when played from starting position $(a,b) = (0,1)$. 
Then, the number of elements in the set of weights $\wt_{\mathcal{C}}(K)$ is the same whether $\lambda = \omega_{i}^{J}$ or $\lambda = \omega_{j}^{J}$. 
So, if $K' = \Jcolor(\mathcal{C}')$ is a $J$-component of $L_{s}$ associated with the $J$-subordinate $\mathcal{C}'$ from $P_{s}$, then $K'$ is a chain with the same length as $K$. 
It follows that the saturated chain $\mathcal{C}'$ must have the same length as the saturated chain $\mathcal{C}$.\hfill\QED

\begin{center}
\underline{\hspace*{4in}}
\end{center}

%
\renewcommand{\refname}{\bf \large
References}
\renewcommand{\baselinestretch}{1.1}

\end{document}